\documentclass[11pt]{book}
\usepackage{stmaryrd}
\usepackage{amssymb}

\usepackage{mathrsfs}

\usepackage{amsmath}

\usepackage{verbatim}

\usepackage{txfonts}
\usepackage{cases}
\usepackage{cases}

\usepackage{makeidx}
\makeindex
\usepackage{amsfonts}
\usepackage{amssymb}
\usepackage{amssymb}
\usepackage{amssymb}
\usepackage[numbers,sort&compress]{natbib}
\usepackage{amsfonts}

\newtheorem{Theorem}{Theorem}[section]
\newtheorem{Lemma}[Theorem]{Lemma}
\newtheorem{Definition}[Theorem]{Definition}
\newtheorem{Remark}{Remark}[section]

\newtheorem{Proposition}{Proposition}[section]
\newtheorem{Corollary}{Corollary}[section]
\newcommand{\benu}{\begin{enumerate}}
\newcommand{\beqa}{\begin{eqnarray}}
\newcommand{\beqan}{\begin{eqnarray*}}
\newcommand{\eay}{\end{array}}
\newcommand{\edm}{\end{displaymath}}
\newcommand{\eenu}{\end{enumerate}}
\newcommand{\eeq}{\end{equation}}
\newcommand{\eeqa}{\end{eqnarray}}
\newcommand{\eeqan}{\end{eqnarray*}}

\newcommand{\bqa}{\begin{eqnarray}}
\newcommand{\eqa}{\end{eqnarray}}
\newcommand{\bqw}{\begin{eqnarray*}}
\newcommand{\eqw}{\end{eqnarray*}}

\newcommand{\non}{\nonumber}
\newcommand{\bea}{\begin{array}{cc}}
\newcommand{\ena}{\end{array}}

\newcommand{\beq}{\begin{equation}}
\newcommand{\enq}{\end{equation}}

 \numberwithin{equation}{section}
\usepackage{ifpdf}
\ifpdf
\usepackage{CJK}
\usepackage[CJKbookmarks=true,
        hyperindex=true,
      pdfstartview=FitH,
     bookmarksnumbered=true,
         bookmarksopen=true,
         citecolor=blue,
         linkcolor=blue,
        unicode,
         colorlinks=true,
         pdfborder=001,
      ]{hyperref}
\else
\usepackage[
         hypertex,
         hyperindex,
      linkcolor=blue,
         unicode,
         citecolor=blue%
         ]{hyperref}
\fi

\usepackage{hypernat} 

\usepackage{makeidx}
\makeindex

\begin{document}
\pagenumbering{roman}
\title{Prandtl Equations and Related Boundary Layer Equations}
\author{Yuming Qin\\
Department of Mathematics, College of Science\\
Institute for Nonlinear Science\\
Donghua University, Shanghai 201620, P. R. China\\
yuming\_qin@hotmail.com,yuming@dhu.edu.cn\\
\\Xiaolei Dong\\
School of Mathematics and Statistics, \\
Zhoukou Normal University\\
Zhoukou Henan,  466001, P. R. China\\
E-mail: xld0908@163.com\\
\\Xiuqing Wang\\
Department of Mathematics, Kunming University of Science and Technology\\
Kunming University of Science and Technology\\
Kunming
Yunnan, 650500, P. R. China\\
E-mail: daqingwang@kust.edu.cn\\}
\date{2024.6.27}
\maketitle
\newpage

\thispagestyle{empty}
\begin{center}
{\bf In memory of Yuming's father, Zhenrong QIN, and Yuming's mother,  Xilan XIA\ \\

Dedicated to Xiaolei's parents, Anmin DONG and Xiurong SUN\\

Dedicated to Xiuqing's parents, Dianqin WANG and Jun LI\\
Dedicated to Yuming's wife and son, Yu YIN and Jia QIN\\
and Xiuqiung's husband  Hongjun Liu}
\end{center}

\newpage

\thispagestyle{empty}

\setcounter{page}{5}
\tableofcontents
\cleardoublepage

\section*{Preface}
\addcontentsline{toc}{chapter}{Preface}
\markboth{Preface}{Preface}\thispagestyle{empty}

At the International Mathematical Congress held in Heidelberg in 1904 Prandtl \cite{prandtl}, in his lecture ``Fluid Motion with Very Small Friction" suggested an explanation of this phenomenon and described the main principles underlying his new theory, currently called the theory of boundary layer. He showed that the flow about a solid body can be divided into two regions: a very thin layer in the neighborhood of the body (the boundary) where viscous friction plays an essential part, and the remaining
region outside this layer where friction may be neglected (the outer flow). Thus, for fluids whose viscosity is small, its influence is perceptible only in a very thin region adjacent to the walls of a body in the flow; the said region, according to Prandtl, is called the boundary layer. Prandtl derived the system of equations for the first approximation of the flow velocity in the boundary layer. This system served as a basis for the development of the boundary layer theory, which has now become one of the fundamental parts of fluid dynamics.

Normally, the boundary layer equation of Navier-Stokes equations is generally called the Prandtl equations.
According to the Prandtl theory, the Prandtl equations can be derived from the Navier-Stokes equations.
In the more than 100 years since the Prandtl equation was put forward, there have been many research results about the Prandtl equation. These results mainly focus on the well-posedness of local and global solutions.
Since the Prandtl equations and the Navier-Stokes equations are closely related, therefore, it is also very important to study the relationship between the Navier-Stokes equations and the Prandtl equations, that is the zero-viscosity limit of Navier-Stokes equations.
The construction of the Navier-Stokes solution is performed as a composite asymptotic expansion involving an Euler solution, a Prandtl boundary layer solution and a correction term. The zero-viscosity limit for the incompressible Navier-Stokes equations  is a challenging problem due to the formation of a boundary layer whose thickness is proportional to the square root of the viscosity. The Navier-Stokes solution  goes to an Euler solution outside a boundary layer, and that it is close to a solution of the Prandtl equations within the boundary layer.

Of course, many mathematicians also study the boundary layer of other fluids such as boundary layer in Magnetohydrodynamics, boundary layer in non-Newtonian flows, Maxwell fluid and other boundary layer equations.

This book aims to present some recent results on the Prandtl equations and MHD boundary layer equations. This book is essentially divided into two parts. The first part includes Chapter 1 in which we shall systematically survey the results till 2020 on the Prandtl equations and MHD boundary layer equations. The second part includes Chapter 2-Chapter 6 where the local  and the global wellposedness  of solutions to the  Prandtl equations  and MHD boundary layer equations. In detail, Chapter 2 is concerned with global well-posedness of solutions to the 2D Prandtl-Hartmann equations in analytic framework. Chapter 3 investigates the local existence of solutions to the 2D Prandtl equations in a weighted Sobolev space. Chapter 4 studies the local well-posedness of solutions to the 2D mixed Prandtl equations in a Sobolev space without monotonicity and lower bound. Chapter 5 is concerned with local existence of solutions to the 2D magnetic Prandtl equations in the Prandtl-Hartmann regime. Chapter 6 proves the local existence of solutions to the 3D Prandtl equations with a special structure. The results in Chapters 2-6 have been obtained recently by the authors and  have never been published before, i.e., which are updated new results on the Prandtl equations and MHD boundary layer equations.

We sincerely wish that the reader will know the main ideas and essence of the basic theories and methods of the Prandtl equations and the MHD boundary layer equations in deriving the local and global well-posedness of solutions. We also wish that the reader can be stimulated by some ideas from this book and continue to undertake the further study and research after having read this book.

We also want to take this opportunity to thank all the people who concern about us including our teachers, colleagues and collaborators, etc.

Yuming Qin also acknowledges the NNSF of China for its support. Currently, this book project is being supported by the NNSF of China with contract number 12171082, the Fundamental Research Funds for the Central Universities with contract numbers 2232022G-13, 2232023G-13. Xiuqing  Wang is being supported by the Yunnan Fundamental Research Projects of China with  contract number 202401AT070411. Xiaolei Dong is being supported by the Natural Science Foundation of Henan with contract number 242300420673 and the Zhoukou Normal University high level talents research start funding project with contract number ZKNUC2022008.\\

\newpage
\pagenumbering{arabic}

\chapter{Survey on the Prandtl Equations and Related Boundary Layer Equations}
\setcounter{equation}{0}
The incompressible Navier-Stokes equations  are the equations of motion describing the momentum conservation of viscous incompressible fluid. The motion equation of viscous fluid was first proposed by Navier in 1827. He only considered the motion of incompressible fluids. Poisson proposed the motion equation of compressible fluids in 1831, and Stokes proposed that the viscosity coefficient is a positive constant in 1845. Fluid equations, such equations are collectively called  Navier-Stokes  equations, and the 2D incompressible Navier-Stokes equations are as follows:
\begin{equation}\left\{
\begin{array}{ll}
{\bf u}_{t}+ ({\bf u}\cdot\nabla){\bf u} -\nu\triangle {\bf u}+\nabla p=0, \\
 \text{div }  {\bf u}= 0 .
\end{array}
 \label{1.1.1}         \right.\end{equation}
Here ${\bf u}=(u,v)$ is the velocity, $p(t,x,y)$ is the pressure, $\nu>0$ is the coefficient of kinematic viscosity.

At the Heidelberg International Mathematical Conference in 1904, Prandtl put forward a basic criterion for describing a phenomenon in his lecture ``Fluid motion with very small friction". He pointed out that the flow motion of a solid can be divided into two areas: in a thin layer near the object, viscous friction plays an important role in it; in the area outside this thin layer, the friction effect can be ignored. Prandtl called this thin layer a boundary layer, and can use the so-called Prandtl equations to describe.
The Prandtl equations are the theoretical model of the Navier-Stokes equations in the thin layer, and the boundary layer equations are the theoretical model of the fluid equations in the thin layer. According to the Prandtl's lecture,
  let $$ t'=t, \ \ x'=x,\ \ y' =\varepsilon^{-1}y, $$
and let $\overline{u}(t',x',y')=u (t',x',\varepsilon y')$, $\overline{v}(t',x',y')=\varepsilon^{-1} v(t',x',\varepsilon y')$, from (\ref{1.1.1}),
 making $\varepsilon$ tends to zero, he obtained
\begin{equation}\left\{
\begin{array}{ll}
u_{t}+ u u_{x}+  v u_{y}=-p_{x}+\nu u_{yy}, \\
 u_{x}+  v_{y}= 0 ,
\end{array}
 \label{1.1.2}         \right.\end{equation}
and $p_{x}$ satisfies the Bernoulli's law:
\begin{eqnarray}
U_{t}+UU_{x}+ p_{x}=0.
\label{1.1.3}
\end{eqnarray}
Before the 20th century, there were very few results about the Prandtl equation, but until 1999,
Oleinik (\cite{O,O1}),  Oleinik and  Samokhin (\cite{OS}) gave the analysis framework of boundary layer theory,
for which there have been many results, for example, wellposedness, viscosity limit, blow up, etc.

Euler equations have a closed relationship with the Navier-Stokes equations, and the 2D incompressible Euler equations are as follows ($ \nu=0 $ in \eqref{1.1.1}):
\begin{equation}\left\{
\begin{array}{ll}
{\bf u}_{t}+ ({\bf u}\cdot\nabla){\bf u} +\nabla p=0, \\
 \text{div }  {\bf u}= 0 .
\end{array}
 \label{1.1.4}         \right.\end{equation}
In recent decades, the relationship among the solutions of Navier-Stokes equations, Euler equations and Prandtl equations has always been a hot research topic of many scholars. There have been some results, for example, wellposedness, viscosity limit, blow up, etc.

K$\acute{a}$rm$\acute{a}$n and  Tsien (\cite{KT}) in 1938 concerned with the theory of the laminar boundary layer in compressible fluid:
\begin{eqnarray*}
\rho u u_{x}+\rho v u_{y}=\nu u_{yy}.
\label{w.1}
\end{eqnarray*}
 Narasimha and  Vasantha  (\cite{NV}) used numerical solution to boundary layer flows of fluids with a high Prandtl number in heat transfer to liquids
and vapours.
Nickel  (\cite{N1})  reviewed the results about the two-dimensional, steady, and incompressible Prandtl's   boundary layer theory.

There are currently four methods  used for boundary layer equations. The first method is von-Mises transformation and Coroco transformation. And the second method is to separate out the shear flow method, that is, $u(t,x,y)=u^{s}(t,y)+\widetilde{u}(t,x,y)$, where the shear flow satisfies the heat equation. The third method is to study the Prandtl equations in Gevrey space.  The fourth method is the energy estimation method, such as. In this chapter, we list some results and main theorems about the boundary layer equations.

\section{2D Prandtl Equations}

Serrin (\cite{S1}) in 1968  illustrated approximation results by studying in detail a specific example of fluid flow where a boundary layer occurs,
to some extent, obtained in this way a miniature version of the full theory.

\subsection{  2D Steady  Prandtl Equations  }
 Fife (\cite{F2,F1}) in 1965, Serrin (\cite{S}) in 1967  considered  the approximation between Navier-Stokes equations and the  steady two-dimensional  Prandtl equations:
\begin{equation}\left\{
\begin{array}{ll}
 u u_{x}+  v u_{y}=UU_{x}+\nu u_{yy}, \\
 u_{x}+  v_{y}= 0 .
\end{array}
 \label{1.1.5}         \right.\end{equation}
 Walter (\cite{W})  considered  the steady 2D (two-dimensional)  Prandtl equations (\ref{1.1.5}),
 and gave various estimates on the rate of convergence of $u(x,y)\rightarrow U(x)$ as $y\rightarrow \infty$ if $u(0,y)\rightarrow U(0)$.
Let dimensionless parameters $\lambda$, $\xi$, and $z$ be determined by the relations
$$\lambda=\frac{\delta}{2aP},\ \delta=\inf_{x\in \mathbb{R}} U^2(x),\ W=\sup_{x\in \mathbb{R}} U^2(x),$$
$$P=\sup_{x\in \mathbb{R}}|2\frac{dp(x)}{dx}|,\ \xi=\frac{\delta}{W},\ z=\frac{\psi^2P}{\nu W^{3/2}},\ \psi=\int_0^y u_0(\tau)d\tau.$$

  Kuznetsov (\cite{K1}) in 1994 considered the boundary layer equations   (\ref{1.1.5})  in the domain $\Omega=\{(x,y)| x\in (0,a), y\in (0,\infty)  \}$, and used von Mises Transform to prove that there exists a unique global solution, i.e., the following theorem.
 \begin{Theorem}(\cite{K1})\label{t.2}
Assume the following inequalities holds,
$$\frac{(2+\lambda)^{2}}{\lambda^{2}(1+\lambda)}<z<\frac{3}{32}\Big(\frac{\xi(2+\lambda)}{2(1+\lambda)}\Big)^{3/2}.$$
Then there exists a unique solution $u,v$ to problem  (\ref{1.1.5}) in $\Omega$ such that $u>0$, the derivatives $\partial u/\partial y $ and $\partial^{2} u/\partial y ^{2}$ are continuous and bounded in $\Omega$,$ \partial u/\partial y\geq m>0$  for $y=0$, and
the functions $\partial u/\partial x,v $,  and $\partial v/\partial y $ are continuous and bounded in any finite part of $\overline{\Omega}$.
\end{Theorem}

  Liu and  Wan (\cite{LW}) in 1985 considered a steady, two-dimensional, laminar, incompressible flow  (\ref{1.1.5}), and found an exact solution if
the transverse component of the external flow velocity is allowed to be deflected upward while the longitudinal component is assumed to be linearly retarded.
It is that this solution is regular at separation and exhibits reverse flow thereafter.
\begin{Theorem} (\cite{LW}) \label{t.1}
For the given positive pressure gradient (\ref{w.3}),   let $g_{0}$, $g_{1}$ be the solution of  (\ref{w.4}), with constants
 $\alpha\approx 9.45, \beta\approx 1.255$. Then,
\begin{equation*}
u=g_{0}'+ xg_{1}',\ \ \ \ v= - g_{1}
 \label{w.5}
 \end{equation*}
in an analytic solution of the problem  (\ref{1.1.5}) in $D_{B}=\{-1\leq x\leq 1, 0\leq y\leq \delta\}$,  the boundary layer thickness to be $\delta=B /R^{\frac{1}{2}}$, $R$ denotes  the Reynolds number, $B=\beta/k^{1/2}$, possessing
the following properties in the region
$$ \{ (x,y)\in D_{B}| y\leq \beta^{*}/k^{\frac{1}{2}} R^{\frac{1}{2}}\}    $$
for some $\beta^{*}$, $0<\beta^{*}\leq \beta$: \\
(i) $v\geq 0$, $v=O(R^{-\frac{1}{2}})$; moreover $u_{y}>0$ on $y=0$ if $x>0$, and $u_{y}>0$ on $y=0$ if $x>0$, and $u_{y}<0$ on $y=0$ if $x<0$;\\
(ii) the separating streamline $\gamma: \varphi=0\ (i.e., x=-g_{0}(y)/g_{1}(y), y>0)$, has slope
$$\frac{3\alpha k^{3/2}}{p_0R^{1/2}}$$
at $(0,0)$ and $x>0$ on $\gamma$;\\
(iii) the curve $T;u=0\ (i.e., x=-g_{0}'(y)/g_{1}'(y), y>0)$ starting from $(0,0)$
lies to the right of $\gamma$. Thus $u<0 ~(u>0)$ to the right (left) of $T$.

Here
\begin{equation}
p_{x}=p_{0}-k^{2}x,\ \ \ \  p_{y}= -k^{2}y,\ \ \ \   -1\leq x\leq 1,\ \ \ \ 0\leq y<+\infty,
 \label{w.3}
 \end{equation}
 \begin{equation}\left\{
\begin{array}{ll}
\nu g_{0}'''+g_{1}g_{0}''-g_{1}'g_{0}''=p_{0}, \\
\nu g_{1}'''+g_{1}g_{1}''-g_{1}'^{2} =-k^{2}.
\end{array}
 \label{w.4}         \right.\end{equation}

\end{Theorem}

Oleinik and  Samokhin (\cite{OS}) in 1999 considered 2D stationary boundary layer  based on the von Mises transformation which reduces the boundary layer system to a single quasilinear parabolic equation. Consider
\begin{equation}\left\{\begin{array}{ll}
uu_{x}+vu_{y}+p_{x}=u_{yy},\\
u_{x}+v_{y}=0,\\
u(0,y)=u_{0}(y),
\end{array}\right.\label{ww1}\end{equation}
in the domain $D=\{0<x<X,0<y<+\infty\}$  with the conditions
\begin{equation}\left\{
\begin{array}{ll}
u(0,y)=u_{0}(y),\ \ u( x,0) =0, \ \  v( x,0) =v_{0}(x),  \\
\lim\limits_{y\rightarrow+\infty}u( x,y)=U( x).
\end{array}
 \label{ww2}         \right.\end{equation}

Introduce von Mises transformation of the variables by
 $$x=x, \ \ \psi =\psi(x,y), \ \ w(x,\psi)=u^{2}(x,y), $$
where
$$u=\frac{\partial \psi}{\partial y}, \ \  v-v_{0}=-\frac{\partial \psi}{\partial x}, \ \ \psi(x,0)=0, $$
the domain $D$ turns into $G=\{0<x<X, 0<\psi<\infty\}$. Then the system  (\ref{ww1}) reduces to the equation
\begin{equation}
\frac{\partial w}{\partial x}+v_{0}\frac{\partial w}{\partial \psi} =\sqrt{w}\frac{\partial^{2} w}{\partial \psi^{2}}-2\frac{\partial p(x)}{\partial x}.
 \label{ww6}
\end{equation}

\begin{Theorem} \label{OS.2.1} (\cite{OS})
Assume that $u_{0}(y)$ for $y>0$; $u_{0}(0)=0$, $u_{0}'(0)>0$, $u_{0}(y)\rightarrow U(0)\neq 0$ as $y\rightarrow \infty$;
$p_{x}$ and $v_{0}(x)$ are continuously differentiable on [0,X];
$u_{0}(y), u_{0}'(y), u_{0}''(y)$ are bounded for $0\leq y\leq \infty$ and satisfy the H$\ddot{o}$lder condition.
Moreover, assume that for small $y$ the following compatibility condition is satisfied at the point $(0,0)$:
\begin{eqnarray}
\frac{d^{2}u_{0}(y)}{d y^{2}}-\frac{dp(0)}{dx}-v_{0}\frac{du_{0}(y)}{dy}=\mathcal{O}(y^{2}).
\label{ww5}
\end{eqnarray}
Then, for some $X>0$ there exists a solution $u(x,y),v(x,y)$ of problem (\ref{ww1})-(\ref{ww2}) in $D$, which has the following properties: $u(x,y)$ is bounded and continuous in $\overline{D}$, $u>0$ for $y>0$;  $u_{y}>m>0$ for $0\leq y\leq y_{0}$, where $m$ and $y_{0}$ are constants; $u_{y}$ and $u_{yy}$ are bounded and continuous in $D$; $u_{x}, v, v_{y}$ are bounded and continuous in any finite portion of $\overline{D}$. \\
Moreover, if $|u'_{0}(y)|\leq m_{1}\exp(-m_{2}y), m_{1}, m_{2}=$const.$>0$, then $u_{x}$ and  $v_{y}$ are bounded in $D$. If either $p_{x}\leq 0$ and $v_{0}(x)\leq 0$ or $p_{x}\leq 0$, then such a solution of problem (\ref{ww1})-(\ref{ww2}) exists in $D$ for any $X>0$.

\end{Theorem}

\begin{Theorem} \label{OS.2.2} (\cite{OS})
Let $u$ and $v$ be two functions satisfying system (\ref{ww1}) in  $D$ , continuous in $\overline{D}$ and satisfying the conditions  (\ref{ww2}). Moreover, assume the following estimates hold: $0<u<c_{1}$ for $\psi>0$,
\begin{eqnarray}
k_{1}y\leq u\leq k_{2}y, \ \  \text{for}\ \  0<y<y_{0}, \label{ww3}\\
u_{yy}\leq k_{3} \ \  \text{in }\ \  D,
\label{ww4}
\end{eqnarray}
where $k_{1},k_{2},k_{3}, y_{0}$,  are positive constants. Then $u,v$ is the only solution of problem (\ref{ww1})-(\ref{ww2}) with these properties.
If $p_{x}\leq 0$ and $v_{0}(x)\leq 0$, then the conditions (\ref{ww3})-(\ref{ww4}) can be dropped. If $p_{x}\leq 0$ and $v_{0}(x)$ is allowed to change sign, then it is only the condition (\ref{ww4}) that is required for the uniqueness.
\end{Theorem}

\begin{Theorem} \label{OS.2.3} (\cite{OS})
Let all the assumptions of Theorem \ref{OS.2.1}, except the compatibility condition (\ref{ww5}), hold for  $u_{0}(y), v_{0}(x),p_{x}$. Then, for some $X>0$, problem (\ref{ww1})-(\ref{ww2}) in the domain $D$ admits a solution $u,v$ with the following properties: $u(x,y)$ is continuous and bounded in $\overline{D}$;
$u>0$ for $y>0$;  $u_{y}\geq m>0$ for $0\leq y\leq y_{0}$ and $x\geq x_{0}>0$, where $m,y_{0}$ are constants depending on $x_{0}$; $u_{y}$ is continuous and bounded in $D$; $u_{yy}$ is continuous in $D$ and bounded in $D \backslash \Omega_{\rho}$, where $\Omega_{\rho}$ is the neighborhood of $(0,0)$ given by the inequality $x^{2}+y^{2}\leq \rho$; $u_{x}, v, v_{y}$ are continuous in $D$ and bounded in any finite portion of the domain $D\backslash\Omega_{\rho}$. If
$|u'_{0}(y)|\leq k_{1}\exp(-k_{2}y)$ for $ y\rightarrow \infty$ , then $u_{x},v_{y}$ are bounded in $D\ \Omega_{\rho}$. If $p_{x}\leq 0$ and $v_{0}\leq 0$ or $p_{x}\leq 0$, then problem (\ref{ww1})-(\ref{ww2}) admits a solution in $D$ for any $X>0$. If $p_{x}\leq 0$ in a neighborhood of $x=0$, this solution is unique.
\end{Theorem}

\begin{Theorem} \label{OS.2.4} (\cite{OS})
Let $U(x)=U_{1}(x+d)^{(\alpha+1)/2}$, $U_{1}=$const.$>0,\alpha\geq -1$, $d=$const.$\geq 0$, $v_{0}(x)\equiv 0$. Let $u(x,y)$ be a solution of problem (\ref{ww1})- (\ref{ww2}) corresponding to the initial velocity profile $u_{0}(y)$, and let $\overline{u}(x,y), \overline{v}(x,y)$ be a self-similar solution of system (\ref{ww1}), i.e., a solution of the form
$$ \overline{u}(x,y)=U(x) f'(\eta), \ \ \overline{v}(x,y)=\int_{0}^{y}\frac{\partial\overline{u} }{\partial x}dy +v_{0}(x).$$
 Then
$$\left|\frac{u(x,y)}{U(x)} - \frac{\overline{u}(x,y)}{U(x)} \right|=o\left(\left[ 1+(\frac{\alpha+1}{2}\ln 2) \right]x^{-\frac{\alpha+1}{2}}\right)\ \ \text{as}\ \  x\rightarrow \infty, $$
uniformly with respect to $y$.\\
Further, consider two solutions $u_{1}(x,y)$ and $u_{2}(x,y)$ of problem (\ref{ww1})-(\ref{ww2})  corresponding to two distinct initial velocity profiles $u_{10}(y)$ and $u_{20}(y)$, let $v_{0}(x)\equiv 0$ and assume that
$$ C_{1}(x+d)^{2m-1}\leq U U_{x} \leq  C_{2}(x+d)^{2n-1}, \ \  C_{1}, C_{2}=\text{const.}>0, $$
where $m\leq n<5m/3$. Then
$$\left|\frac{u_{1}}{U } - \frac{u_{2}}{U } \right|=o\left(1\right)\ \ \text{as}\ \  x\rightarrow \infty, $$
uniformly with respect to $y$.
\end{Theorem}

\begin{Theorem} \label{OS.2.5} (\cite{OS})
Assume that
\begin{eqnarray*}
&& \lim\limits_{y\rightarrow+\infty} U( x)=U_{\infty}= \text{const.}>0, \ \  v_{0}(x)\equiv 0, \\
&& 0\leq u_{0}(y)\leq U(0), \ \  u_{0}(0)=0, \ \   u'_{0}(0)>0,  \\
&& 0\leq \frac{dU}{dx} \leq \frac{M_{0}}{(x+1)^{\gamma_{0}+1}}, \ \ \gamma_{0}>0.
\end{eqnarray*}
Then
$$|u(x,y)- u_{1}(x,y)|\rightarrow 0, \ \   \text{as}\ \  x\rightarrow \infty,$$
uniformly in $y\in [0,\infty)$, where $u(x,y)$ is the solution of problem (\ref{ww1})-(\ref{ww2}) corresponding to the initial velocity profile $u_{0}(y)$;
$$u_{1}(x,y)=U_{\infty}f'(\eta), \ \   \text{when}\ \  \eta=y\sqrt{\frac{U_{\infty}}{2(x+1)}} $$
and  $f(\eta)$ is the solution of the following boundary value problem
\begin{eqnarray*}
&& f'''+ff''=0, \ \  0<\eta<\infty ,  \\
&& f(0)=0, \ \ f'(0)=0, \ \  f'(\infty)=1.
\end{eqnarray*}
Under the additional conditions
\begin{eqnarray*}
&& U(0)f'(y-N)\leq u_{0}(y) , \ \ y\in [N,\infty),  \\
&& |u_{0}(y)-U(0)|\leq M_{1}\exp(-\gamma_{1}y^{2}), \ \ y\in[0,\infty),
\end{eqnarray*}
with some constants $N,M,\gamma_{1}>0$, the following estimate holds:
$$|u(x,y)- u_{1}(x,y)|\leq \frac{M_{2}}{(x+1)^{\gamma_{2}}}$$
where $M_{2}, 0\leq \gamma_{2},\ \gamma_{0}$ are constants depending only on the data of the
problem.
\end{Theorem}

\begin{Theorem} \label{OS.2.6} (\cite{OS})
Assume that problem (\ref{ww1})-(\ref{ww2}) admits a solution in a domain $D$. Then $X<x_{0}$, where $x_{0}$ is determined by the conditions
$$\max u^{2}_{0}(y)-2\int_{0}^{x_{0}} \frac{dp}{dx} dx=0, \ \  \frac{dp(x_{0})}{dx}>0.  $$
\end{Theorem}

\begin{Theorem} \label{OS.2.7} (\cite{OS})
Assume that $dp/dx\geq 0$ and there exist $x_{1}, x_{2}$ such that $x_{1}\leq x_{2}$ and
\begin{eqnarray*}
 \int_{x_{1}}^{x_{2}} v_{0}(t)dt> \sqrt{6M_{0}(x_{2}-x_{1})}.
\end{eqnarray*}
Then separation of the boundary layer occurs at some points $x_{1}<x_{2}$.

\end{Theorem}

\begin{Theorem} \label{OS.2.8} (\cite{OS})
Assume the following: $v_{0}(x)\leq 0$, $v'_{0}(x)\leq 0$ in the interval $0\leq x\leq X$; $dp/dx$ and $v_{0}(x)$ belong to the class
 $C^{1}([0,X])$; the initial velocity
 profile $u_{0}(y)$ is such that
\begin{eqnarray*}
&& u_{0}(0)=0, \ \ u_{0}(y)>0 \ \  \text{for}\ \  y>0, \ \   \lim\limits_{y\rightarrow \infty}u_{0}(y)= U(0), \\
&& |u_{0}(y)|, |u'_{0}(y)|, |u''_{0}(y)|\leq M, \ \   y\in [0,\infty );
\end{eqnarray*}
the H$\ddot{o}$lder inequality with exponent $\alpha=\alpha(N)\in (0,1)$ holds for $u''_{0}(y)$ on the segment $N^{-1}\leq y\leq N$,  for any $N>1$;  moreover,
$$l_{0}(u_{0})=u''_{0}(y)- u'_{0}(y)- \frac{dp(0)}{dx}=\mathcal{O}(y^{2}), \ \ \text{as}\ \  y\rightarrow 0. $$
Assume also that
$$\frac{dp }{dx}\leq (Cx+1)^{-3}\frac{dp(0)}{dx} \ \  \text{for}\ \  0\leq x<\infty, \ \ C=\sup\limits_{y\geq 0}\left|\frac{l_{0}(u_{0})}{u_{0}^{2}} \right|. $$
Then for any finite $X>0$, problem (\ref{ww1})-(\ref{ww2}) in $D$ admits a unique solution $u,v$ such that: $u,u_{y},u_{yy}$ are bounded and
continuous in $D$; $v,u_{x},v_{y}$ are continuous in $D$, and the velocity profile $u(x_{0},y)$, for any $x_{0}\in [0,X]$ has the properties of
the initial profile $u_{0}(y)=u(0,y)$.

\end{Theorem}

Assume that $dp/dx\leq 0$ and $v_{0}(x)\equiv 0$. Applying the von Mises transformation,  then we obtain the equation
\begin{equation}
\frac{\partial w}{\partial x} =\sqrt{w}\frac{\partial^{2} w}{\partial \psi^{2}}-2\frac{\partial p(x)}{\partial x}
 \label{ww7}
\end{equation}
in  $G=\{0<x<X, 0<\psi<\infty\}$, with the boundary conditions
\begin{equation}\left\{\begin{array}{ll}
w(x,0)=0,  \ \ w(0,\psi)=w_{0}(\psi), \\
\lim\limits_{y \rightarrow \infty} w(x,\psi)=U^{2}(x)  \ \ \text{uniformly on } \ \ [0,X].
\end{array}\right.\label{ww8}\end{equation}

\begin{Theorem} \label{OS.2.9} (\cite{OS})
Assume that: $U(x)$ has a bounded second derivative on the segment $[0,X]$, $U(x)>0, U'(x)\geq 0$; $u_{0}(y)\in C^{2}(0\leq y\leq \infty)$,
$u_{0}(0)=0, u'_{0}(0)>0, u''_{0}(y)\leq 0$, and $u''_{0}(y)$ satisfies the local H$\ddot{o}$lder condition for $y>0$. Moreover, let
\begin{eqnarray*}
&&u_{0}(y)\rightarrow U(0), \ \  \text{as } \ \ y \rightarrow 0,\\
&& u''_{0}(y)+U(0)U'(0)=\mathcal{O}(y^{2})\  \text{as } \ \ y \rightarrow 0.
\end{eqnarray*}
Then problem (\ref{ww7})-(\ref{ww8}) has a unique solution w with the following properties: $w, \partial w/\partial x, \partial w/\partial \psi$ are
 bounded and continuous in $\overline{G}$, and satisfy the inequalities
\begin{eqnarray*}
&& w>0, \ \ \left|\frac{\partial w}{\partial x}\right|\leq C\psi,\\
&& \frac{\partial w}{\partial \psi}(x,0)\geq \alpha>0, \ \  \frac{\partial w}{\partial x}\leq 2 UU' \ \ \text{in} \ \ G, \ \  C,\alpha= \text{const}.
\end{eqnarray*}

\end{Theorem}

\begin{Theorem} \label{OS.2.10} (\cite{OS})
 Let $f(x,z)=\left(f_{0}(x,z_{0},z_{1},...), f_{1}(x,z_{0},z_{1},...),...\right)$ be a mapping $[0,X]\times B\rightarrow B$  which is continuous in $x$
 and satisfies the Lipschitz condition in $z$:
$$\|f(x,z)-f(x,\overline{z});B\|\leq L\|z-\overline{z};B\|, \ \  z,\overline{z}\in B; \ \  x\in[0,X],$$
where $L=$const.$>0$. Then the Cauchy problem
\begin{equation}
\frac{du}{d x} =f(x,u), \ \ x\in [0,X], \ \ u(0)\in B,
 \label{ww9}
\end{equation}
has a unique solution $u(x)\in C^{1}(0,X; B)$. Moreover, if $u_{0}(x)$ belongs to the space $ C^{0}(0,X; B)$, then the successive approximations $u_{n}(x)$
given by
$$u_{n+1}=u(0)+\int_{0}^{x}f(\tau, u_{n}(\tau))d\tau, \ \ n=0,1,2,...,$$
are uniformly convergent to $u(x)$ in the norm of $B$, as $n\rightarrow \infty$.

\end{Theorem}

\begin{Theorem} \label{OS.2.11} (\cite{OS})
Assume that the conditions of Theorem \ref{OS.2.10} hold for $f(x,z)$, and $f(x,z)$ increases quasimonotonically in $z$. Let $u\in C^{1}(0,X; B)$
be a solution of problem $(\ref{ww9})$, and let $v\in C^{1}(0,X; B)$ be a function such that
$$\frac{dv}{dx}\leq f(x,v), \ \ x\in [0,X], \ \ v(0)\leq  u(0). $$
Then, $v(x)\leq u(x)$ for $x\in [0,X]$.

\end{Theorem}

\begin{Theorem} \label{OS.2.12} (\cite{OS})
Problem (\ref{ww7})-(\ref{ww8}) with $p_{x}(x)\leq 0$ admits a unique solution such that $w_{x}, w_{\psi},w_{\psi\psi}$ are continuous in $G$,
 and $w\geq0, w_{\psi\psi}\leq 0$ in $G$.

\end{Theorem}



 Oleinik and  Samokhin (\cite{OS}) in  1999 considered
the boundary layer system for an axially symmetric  incompressible stationary flow has the form
\begin{equation}\left\{
\begin{array}{ll}
 u\frac{\partial u}{\partial x}+v\frac{\partial u}{\partial y}=  \nu\frac{\partial^{2} u}{\partial y^{2}} +U \frac{dU}{d x},  \\
\frac{\partial (ru)}{\partial x}+\frac{\partial (rv)}{\partial y} =0,
\end{array}
 \label{ww12}         \right.\end{equation}
in a domain $D=\{0<x<X, 0<y<\infty\}$, with the boundary conditions
\begin{equation}\left\{
\begin{array}{ll}
u(0,y)=0,\ \ u(x,0)=0, \ \  v( x,0) =v_{0}(x) ,  \\
u(x,y)\rightarrow U(x)  \ \  \text{as}\ \  y\rightarrow \infty.
\end{array}
 \label{ww13}         \right.\end{equation}
The function $r(x)$ determines the surface of the body past which the fluid flows: $r(0)=0,r_{x}(0)\neq 0, r(x)>0$ for $x>0$; $U(x)$ is a given
longitudinal velocity component of the outer flow: $U(0)=0, U(x)>0$ for $x>0$, $U_{x}>0$ for $x\geq 0$; $\nu$ is the viscosity coefficient.

Introduce new independent Crocco variables by
$$\xi=x, \ \ \eta=\frac{u(x,y)}{U(x)},$$
and a new unknown function
$$w(\xi,\eta)= \frac{u_{y}(x,y)}{U(x)}.$$
From  (\ref{ww12}), they obtained the following equation for $w(\xi,\eta)$
\begin{eqnarray}
\nu w^{2}w_{\eta\eta}-\eta Uw_{\xi}+Aw_{\eta}+Bw=0
 \label{ww14}
\end{eqnarray}
in the domain $\Omega=\{0<\xi<X, 0<\eta <1 \}$, with the boundary conditions
\begin{eqnarray}
w|_{\eta=1}=0, \ \  (\nu w^{2}w_{\eta\eta}-v_{0}w+C)|_{\eta=0}=0,
 \label{ww15}
\end{eqnarray}
where
$$A=(\eta^{2}-1)U_{x}(\xi), \ \  B=-\eta U_{x}(\xi)+\frac{\eta r_{x}(\xi)U(\xi)}{r(\xi)}, \ \  C=U_{x}(\xi).$$
Note that from the above assumptions on $U(x)$ and $r(x)$, one has
$$\lim\limits_{\xi\rightarrow 0}B(\xi, \eta)=0. $$

\begin{Theorem} \label{OS.2.16} (\cite{OS})
Assume that the functions $r(x), U(x)$ are twice continuously differentiable,$U_{x}$, $U(0)=0,|U_{xx}|\leq N_{1}x$, $v_{0}\leq N_{2}x$, $N_{i}=$const.$>0$,
  and $v_{0}$ has its first derivative bounded.

Then problem  (\ref{ww14})-(\ref{ww15}) has a solution $w$ in the domain $\Omega=\{0<\xi<X ,0<\eta<1\}$, where $X>0$ depends only on $U,r,v_{0}$. This solution has the following properties: $w( \xi,\eta)$ is continuous in $\overline{\Omega}$; $M_{9}(1-\eta)\sigma\leq w \leq  M_{2}(1-\eta)\sigma$; $w_{\eta}$ is continuous in $\eta$ for $\eta<1$; $-M_{22} \sigma\leq w_{\eta} \leq - M_{23} \sigma$ in $\Omega$, $|w_{\xi}|\leq M_{24}(1-\eta)\sigma$, $ww_{\eta\eta}\leq -M_{25}$,  and $ww_{\eta\eta}$ is bounded in $\Omega$, $\sigma=(-\ln \mu(1-\eta))^{1/2}, 0<\mu<1$. In any closed
interior subdomain of $\Omega$, the function $w$ and its derivatives in (\ref{ww14}) satisfy the H$\ddot{o}$lder condition.
Furthermore, $w$ satisfies
 equation (\ref{ww14}) in $\Omega$ and  also the conditions (\ref{ww15}) for $0\leq \xi\leq X$.

\end{Theorem}

\begin{Theorem} \label{OS.2.17} (\cite{OS})
Assume the assumptions of Theorem \ref{OS.2.16} hold for $U(x), r(x), v_{0}(x)$. Then problem (\ref{ww12})-(\ref{ww13}) in the domain $D$ admits a solution with
 the following properties: $u/U, u_{y}/U$ are continuous and bounded in $\overline{D}$; $u>0$ for $y>0$ and $x>0$; $u\rightarrow U$ as $y\rightarrow \infty$,
$u(x,0)=u(0,y)=0$; $ u_{y}/U>0$ for $y>0$; $ u_{y}/U\rightarrow 0$  as $y\rightarrow \infty$; $u_{yy}, u_{x}, u_{y}$  are bounded and continuous with respect
 to $y$ in $\overline{D}$;  $v$ is bounded for bounded $y$; $u_{yy}$ is bounded in $\overline{D}$;   $u_{xy}$ is bounded in $D$ for bounded $y$;  $u_{xy}$
and  $u_{yy}$ are continuous in $\overline{D}$; $u_{yy}/u_{y}$ is continuous in $\overline{D}$ with respect to $y$. Moreover,
\begin{eqnarray*}
&& M_{9}(U-u)\sigma\leq u_{y} \leq  M_{2}(U-u)\sigma, \\
&& -M_{22} \sigma\leq \frac{u_{yy}}{u_{y}} \leq - M_{23} \sigma, \ \ -M_{26} \leq  \frac{u_{yyy}u_{y}-u_{yy}^{2}}{u_{y}^{2}} \leq - M_{25},  \\
&& \left| \frac{1}{u_{y}}(u_{xy}u_{y}-u_{x}u_{yy} )+ \frac{U_{x}}{Uu_{y}} (uu_{yy}- u_{y }^{2})\right| \leq M_{24}(U-u)\sigma, \\
&& M_{32} U\exp (-M_{33}y^{2}) \leq U-u\leq M_{34} U\exp (-M_{35}y^{2}),
\end{eqnarray*}
where $\sigma =(-\ln \mu(1-u/U))^{1/2}$, $M_{i}, \mu$ are positive constants, $\mu<e^{-1/2}$.

\end{Theorem}

\begin{Theorem} \label{OS.2.18} (\cite{OS})
Assume $(u,v)$ is a solution of problem (\ref{ww12})-(\ref{ww13}) such that: the derivatives $u_{x},u_{y},v_{y},u_{yy}$, $u_{yyy},u_{xy}$ are continuous in $D$;
 $ u/U$ and $ u_{y}/U $ are continuous in $\overline{D}$;  $ u_{y} >0$ for $y\geq 0, x>0$; $ u_{y}/U>0$ for $y=0$;  $u_{y}/U\rightarrow 0$ as
 $y\rightarrow \infty$; $u_{yy}/u_{y}, u_{x}$ are continuous with respect to $y$ at  $y=0$; $(u_{yyy}u_{y}-u^{2}_{y})/u^{2}_{y}\leq 0$.
  Then $u,v$ is the only solution of problem (\ref{ww12})-(\ref{ww13}) with these properties.

\end{Theorem}

Oleinik and  Samokhin (\cite{OS}) in 1999 considered
the boundary layer system for the symmetric  incompressible stationary flow:
\begin{equation}\left\{
\begin{array}{ll}
 u\frac{\partial u}{\partial x}+v\frac{\partial u}{\partial y}=  \nu\frac{\partial^{2} u}{\partial y^{2}} +U \frac{dU}{d x},  \\
\frac{\partial u}{\partial x}+\frac{\partial v}{\partial y} =0,
\end{array}
 \label{ww16}         \right.\end{equation}
in a domain $D=\{0<x<X, 0<y<\infty\}$, with the boundary conditions
\begin{equation}\left\{
\begin{array}{ll}
u(0,y)=0,\ \ u(x,0)=0, \ \  v( x,0) =v_{0}(x) ,  \\
u(x,y)\rightarrow U(x)  \ \  \text{as}\ \  y\rightarrow \infty.
\end{array}
 \label{ww17}         \right.\end{equation}
Let
\begin{eqnarray}
U(x)=\sum\limits_{m=0}^{q} a_{m}x^{2m+1}+a_{q+1}(x),
 \label{ww18}
\end{eqnarray}
where $a_{0}=a, a_{j}=$const., $j=1,...,q$, and
\begin{eqnarray*}
|a_{q+1}(x)| \leq C_{1}x^{2(q+1)+1}, \ \  |a_{(q+1)x}(x)| \leq C_{2}x^{2(q+1) }, \ \  |a_{(q+1)xx}(x)| \leq C_{3}x^{2 q+1  };
\end{eqnarray*}
and let
\begin{eqnarray}
v_{0}(x)=\sum\limits_{m=0}^{q} b_{m}x^{2m}+b_{q+1}(x),
 \label{ww19}
\end{eqnarray}
where $b_{0}=b, b_{j}=$const., $j=1,...,q$,
\begin{eqnarray*}
|b_{q+1}(x)| \leq C_{4}x^{2(q+1) }, \ \  |b_{(q+1)x}(x)| \leq C_{5}x^{2 q+1 },
\end{eqnarray*}
$$ C_{j}= \text{const.}>0, \ \ j=1,...,5 , \ \  q\geq 0. $$

\begin{Theorem} \label{OS.2.19} (\cite{OS})
Suppose that conditions (\ref{ww18})-(\ref{ww19}) with some $q\geq 0$ hold for $U(x), v_{0}(x)$. Then problem (\ref{ww16})-(\ref{ww17}) in the domain
$D=\{0<x<X,0<y<\infty\}$ ( $X$ depends on $U, v_{0}$) admits a solution with the following properties: $u/U, u_{y}/U$ are continuous and bounded in
$\overline{D}$; $u>0$ for $y>0,x>0$; $u\rightarrow U(x)$ as  $y\rightarrow \infty$ so that
$$1-\frac{u}{U}=\exp  \left\{-\frac{U_{x}(0)}{2\nu}[y^{2}(1+\mathcal{O}(x^{2}))+ \mathcal{O}(y^{1+\varepsilon}))]  \right\}, \ \  \text{as} \ \
 y\rightarrow \infty, $$
where $\varepsilon>0$ is arbitrarily small; $ u_{y}/U>0$ for $y\geq 0$;  $ u_{y}/U\rightarrow 0$  as $y\rightarrow \infty$; $u_{yy},u_{x},u_{y}$ are
bounded and continuous in $\overline{D}$ with respect to $y$, and are continuous in $D$ with respect to $x,y$; $v$ is continuous in $\overline{D}$ with
respect to $y$, and in $D$ it is continuous with respect to  $x,y$ and bounded for bounded $y$; $u_{xy}$ and $u_{yyy}$ are continuous in $D$; $u_{xy}$ is
bounded in $D$ for bounded $y$; $u_{yy}/u_{y}$ is continuous in $\overline{D}$ with respect to $y$. Moreover, the following estimates hold:
\begin{eqnarray*}
&& U(x)Y(\frac{u}{U})(1-M_{22}x^{2}) \leq u_{y}\leq U(x)Y(\frac{u}{U})(1+M_{22}x^{2}), \\
&& Y_{\eta}(\frac{u}{U})(1+M_{14}x^{2}) \leq \frac{u_{yy}}{u_{y}}\leq  Y_{\eta}(\frac{u}{U})(1-M_{15}x^{2}), \\
&& -M_{25} <  \frac{u_{yyy}u_{y}-u_{yy}^{2}}{u_{y}^{2}}< - M_{26}.
\end{eqnarray*}

\end{Theorem}

\begin{Theorem} \label{OS.2.20} (\cite{OS})
A solution $u,v$ of problem (\ref{ww16})-(\ref{ww17}) in $D$ is unique if it possesses the following properties: the derivatives
$u_{x}, u_{y},v_{y},u_{yy},u_{yyy}, u_{xy}$ are continuous in $D$; $u/U,u_{y}/U$ are continuous in  $\overline{D}$, $v$ is continuous
in $\overline{D}$ with respect to $y$; $u_{y}>0$ for $y\geq 0$, $x>0$; $ u_{y}/U>0$ for $y=0$;
$ u_{y}/U\rightarrow  0$ os $y\rightarrow \infty$; $u_{yy}/u_{y}, u_{x}$ are continuous with respect to $y$ at $y=0$.

\end{Theorem}

\begin{Theorem} \label{OS.2.21} (\cite{OS})
 Assume that conditions (\ref{ww18})-(\ref{ww19}) with some $q\geq 0$ hold for $U(x), v_{0}(x)$. Then the following estimate holds for
the solution $u,v$ of problem (\ref{ww16})-(\ref{ww17}):
$$\left| \frac{u_{y}}{U}-\sum\limits_{m=0}^{q} Y_{m}\frac{u }{U}x^{2m}  \right|\leq M_{24}Y_{0}\frac{u }{U}x^{2q+2} , $$
where $M_{24} $ is a positive constant.

\end{Theorem}

\begin{Theorem} \label{OS.2.22} (\cite{OS})
 Let $U(x)=x^{m}V(x)$, $v_{0}(x)=x^{(m-1)/2}v_{1}(x)$,  where $V(x)$ and $v_{1}(x)$ satisfy
\begin{eqnarray*}
&& V(x)=a+a_{1}(x), \ \ a>0, \ \  v_{1}(x)=b+b_{1}(x), \\
&& |a_{1}(x)|\leq M_{32}x^{3}, \ \  |a_{1x}(x)|\leq M_{33}x^{2}, \ \ |a_{1xx}(x)|\leq M_{34}x , \\
&&  |b_{1}(x)|\leq M_{35}x^{2}, \ \  |b_{1x}(x)|\leq M_{36}x .
\end{eqnarray*}
  Then problem (\ref{ww16})-(\ref{ww17}) in the domain $D$, for some $X$ depending on $U$ and $v_{0}$, admits a solution $u(x,y), v(x,y)$ with the following
properties: $u_{y}>0$ for $y\geq 0$, $x>0$; $u/U,u_{y}/ (x^{(m-1)/2}U(x) )$ are bounded and continuous in $\overline{D}$; $u>0$ for $y>0$ and $x>0$;
$u(x,y)\rightarrow U(x), u_{y}\rightarrow 0$ as $y\rightarrow \infty$; $u_{x},u_{y}, u_{yy}$ are bounded and continuous in $\overline{D}$ with respect to $y$,
and continuous with respect to $x,y$ in $D$; $v$ is continuous in  $\overline{D}$ with respect to $y$, and continuous with respect to $x,y$ in $D$;
$u_{yy}/(x^{(m-1)/2}u_{y})$ is continuous in $\overline{D}$ with respect to $y$; $u_{xy}$ and $u_{yyy}$ are continuous in $D$. Moreover,
the following inequalities hold:
\begin{eqnarray*}
&& M_{44}\exp (-M_{45}x^{m-1}y^{2})\leq 1-\frac{u}{U}\leq M_{44}\exp (-M_{46}x^{m-1}y^{2}), \\
&& x^{\frac{m-1}{2}} Y(\frac{u}{U})(1-M_{37}x^{2}) \leq \frac{u_{y}}{U}\leq x^{\frac{m-1}{2}} Y(\frac{u}{U})(1+M_{38}x^{2}), \\
&& Y_{\eta}(\frac{u}{U})(1+M_{39}x^{2}) \leq \frac{u_{yy}}{u_{y}}\leq Y_{\eta}(\frac{u}{U})(1-M_{40}x^{2})  .
\end{eqnarray*}
The solution with the above properties is unique.

\end{Theorem}

 Oleinik and  Samokhin (\cite{OS}) in 1999 considered
the boundary layer system (\ref{ww16}) in a domain $D=\{0<x<X, 0<y<\infty\}$, with the boundary conditions
\begin{equation}\left\{
\begin{array}{ll}
u(0,y)=u_{0}(y),\ \ u(x,0)=0, \ \  v( x,0) =v_{0}(x) ,  \\
u(x,y)\rightarrow U(x)  \ \  \text{as}\ \  y\rightarrow \infty.
\end{array}
 \label{ww20}         \right.\end{equation}

\begin{Theorem} \label{OS.2.23} (\cite{OS})
Assume that $U(x)>0$ for $0\leq x\leq X$; $U_{x}(x)$ and $v_{0}(x)$ are continuously differentiable on the segment $0\leq x\leq X$;
$u_{0}(y)$ is continuous for $0\leq y<\infty, u_{0}(0)=0, u_{0}(y)>0$ for $y>0$; $u_{0}(y)\rightarrow U(0)$ as $y\rightarrow \infty$;
 the derivatives $u_{0y}(y), u_{0yy}(y), u_{0yyy}(y)$ exist and satisfy the inequalities
\begin{eqnarray*}
&& M_{1}\left(1-\frac{u_{0}(y)}{U(0)} \right)\sqrt{-\ln \mu\left(1-\frac{u_{0 }}{U(0)} \right)} \leq \frac{u_{0y}(y)}{U(0)}
\leq  M_{2}\left(1-\frac{u_{0}(y)}{U(0)} \right)\sqrt{-\ln \mu\left(1-\frac{u_{0 }}{U(0)} \right)}, \\
&&- M_{3} \sqrt{-\ln \mu\left(1-\frac{u_{0 }}{U(0)} \right)} \leq \frac{u_{0yy}}{u_{0y}}
\leq  -M_{4} \sqrt{-\ln \mu\left(1-\frac{u_{0 }}{U(0)} \right)},  \\
&& \left|\frac{u_{0yyy}u_{0y}-u_{0yy}^{2}}{u_{0y}^{2}}\right|\leq M_{5},
\end{eqnarray*}
where $\mu, M_{i}(i=1,2,\cdots,5)$ are positive constants and $0<\mu<1$. Moreover, let the following compatibility conditions hold:
\begin{eqnarray*}
&& v_{0}(0)u_{0y}(0) =U(0)U_{x}(0)+\nu u_{0yy}(0), \\
&&v_{0}(0)u_{0yy}(y)-\nu u_{0yyy}(y)=\mathcal{O}(y), \ \  \text{as}\ \ y\rightarrow 0.
\end{eqnarray*}
Then problem (\ref{ww16}), (\ref{ww20}) in the domain $D$, with $X$ depending on $U, v_{0}, u_{0}$, admits a uniqtie solution with the following properties:
$u, u_{y}$ are bounded and continuous in $\overline{D}$; $u>0$ for $y>0$, $u_{y}>0$ for $y\geq 0$; $u_{yy}, u_{x}, v_{y} $ are bounded and continuous
in $\overline{D}$ with respect to $y$;  $v$ is bounded for  bounded $y$; $u_{xy}, u_{yyy}$ are continuous in $\overline{D}$, $u_{yyy}$ is bounded in $D$,
$u_{xy}$ is bounded in $D$ for bounded $y$; $u_{yy}/u_{y}$ is continuous in $D$ with respect to $y$.
Moreover,
\begin{eqnarray*}
&& M_{6}\left(1-\frac{u  }{U } \right)\sqrt{-\ln \mu\left(1-\frac{u }{U } \right)} \leq \frac{u_{ y}(y)}{U }
\leq  M_{7}\left(1-\frac{u }{U } \right)\sqrt{-\ln \mu\left(1-\frac{u }{U } \right)}, \\
&& M_{8}\exp (-M_{9}y^{2})\leq 1-\frac{u  }{U }\leq M_{10}\exp (-M_{11}y^{2}), \\
&&- M_{12} \sqrt{-\ln \mu\left(1-\frac{u }{U } \right)} \leq \frac{u_{ yy}}{u_{ y}}
\leq  -M_{13} \sqrt{-\ln \mu\left(1-\frac{u }{U } \right)},  \\
&& \left|\frac{u_{ yyy}u_{ y}-u_{ yy}^{2}}{u_{ y}^{2}}\right|\leq M_{14},
\end{eqnarray*}
where $M_{i} (i=1,2,...,5), \mu$ are positive constants, $\mu<1$.

\end{Theorem}

\begin{Theorem} \label{OS.2.24} (\cite{OS})
Assume that $U_{x}(x), v_{0}(x)$ have piecewise continuous first derivatives, and $u_{0}(y)$ is continuous for $0<y<\infty$, $u_{0}(0)=0 , u_{0}(y)>0$
for $y>0$, $u_{0}(y)\rightarrow U(0)$ as $y\rightarrow \infty$, $U(x)>0$ for $x\geq 0$,
$$M_{15}(U(0)-u_{0}(y))\leq u_{0y}(y) \leq M_{16}(U(0)-u_{0}(y))$$
where $M_{i}$ are positive constants. Assume also that $u_{0}$ has bounded measurable derivatives up to the third order such that
$$\frac{u_{0yy}}{u_{0y}}, \ \  \ \   \frac{u_{0yyy}u_{0y}-u_{0yy}^{2}}{u_{0y}^{2}}$$
are bounded. Moreover, assume the following compatibility conditions hold:
\begin{eqnarray*}
&&v_{0}(0)u_{0y}(0)=-p_{x}(0)+\nu u_{0yy}(0), \\
&& v_{0}(0)u_{0yy}(y)-\nu u_{0yyy}(y)=\mathcal{O}(y), \ \  y\rightarrow 0.
\end{eqnarray*}
Then problem (\ref{ww16}), (\ref{ww20})  in $D$, with $X$ depending on the given functions, admits a unique solution with the following properties:
$u,u_{y}$ are continuous and bounded in $\overline{D}$;  $u_{y}>0$ for $y\geq 0$; $v$ is continuous in $y$ and
bounded for bounded $y$; the weak derivatives $u_{x}, u_{yx}, u_{yy}, u_{yyy}, v_{y}$ exist
and are bounded measurable functions in $D$. Moreover,
$$\frac{u_{ yy}}{u_{ y}}, \ \  \ \   \frac{u_{ yyy}u_{ y}-u_{ yy}^{2}}{u_{ y}^{2}}$$
are bounded, and
\begin{eqnarray*}
&& M_{17}\exp (-M_{18}y )\leq 1-\frac{u  }{U }\leq M_{19}\exp (-M_{20}y ), \\
&& M_{21}(U(x)-u(x,y))\leq u_{y}(x,y)\leq M_{22}   (U(x)-u(x,y))\leq u_{y}(x,y), \\
&& \left| \frac{1}{u_{y}}(u_{yx} u_{y}-u_{x}u_{yy}) +\frac{U_{x}}{u_{y}U} (uu_{yy}-u_{y}^{2})\right|\leq M_{23}(U-u).
\end{eqnarray*}

\end{Theorem}

\begin{Theorem} \label{OS.2.25} (\cite{OS})
Assume that the function $U(x)$ is continuous in the interval $0\leq x\leq X$, and $U_{x}(x), v_{0}(x)$ are piecewise continuous, $U_{x}(x)\geq$const.$>0$.
 Assume also that $u_{0}(y)$ is continuous in the interval $0\leq y<\infty$, $u_{0}(0)=0, u_{0}(y)>0$ for $y>0$, $u_{0}(y)\rightarrow U(0)$ as
 $y\rightarrow \infty$,
  the derivative $u_{0y}(y)$ is continuous, $u_{0y}>0$; there exists the weak derivative $u_{0yy}(y)$ and the following inequalities hold:
\begin{eqnarray*}
&& M_{21}\left(1-\frac{u_{0}(y) }{U(0)} \right)\left(-\ln \mu\left(1-\frac{u_{0}(y) }{U(0)} \right)\right)^{\frac{1}{2}} \leq \frac{u_{0}(y)}{U(0) }
\leq  M_{22}\left(1-\frac{u_{0}(y) }{U(0)} \right)\left(-\ln \mu\left(1-\frac{u_{0}(y) }{U(0)} \right)\right)^{\frac{1}{2}}, \\
&& \int_{0}^{\infty} \frac{u_{0yy}^{2}(y)}{u_{0y}(y)} \leq \infty, \ \  \mu=\text{const.}, \ \  0<\mu<1.
\end{eqnarray*}
Then problem (\ref{ww16}), (\ref{ww20})  has a weak solution $u(x,y) ,v(x,y)$ with the following properties: $u(x,y)$ is a bounded measurable function continuous
in $\overline{D}$ with respect to $y$; $u>0$ for $y>0$, $u=0$ for $y=0$; $u_{y}(x,y)$ is bounded in $D$ and
\begin{eqnarray*}
&& M_{23}\left(1-\frac{u  }{U } \right)\left(-\ln \mu\left(1-\frac{u }{U } \right)\right)^{\frac{1}{2}} \leq \frac{u_{ y} }{U }
\leq  M_{24}\left(1-\frac{u  }{U } \right)\left(-\ln \mu\left(1-\frac{u }{U } \right)\right)^{\frac{1}{2}} , \\
&& M_{25}\exp (-M_{26}y^{2})\leq 1-\frac{u  }{U }\leq M_{27}\exp (-M_{28}y^{2});
\end{eqnarray*}
there exist weak derivatives $u_{yy}$ and $u_{yyy}$ such that
$$\int_{D}\frac{u_{yy}^{2}}{u_{y}}dxdy  < \infty, \ \  \int_{D}\frac{(u_{yyy}u_{y}-u_{y}^{2})^{2}}{u_{y}^{4}}u^{2}dxdy  < \infty; $$
and $v(x,y)$ is a measurable function in $D$ such that
$$\int_{D\cap\{(x,y):0\leq x\leq X-\varepsilon\} } v(x,y)\exp(-My^{2}) dxdy  < \infty, \ \  \forall M,\varepsilon =\text{const.}>0.$$
Equations (\ref{ww16}) and the boundary conditions (\ref{ww20}) at $x=0, y=0$ hold
in the sense of the following integral identities:
\begin{eqnarray*}
&& \int_{D}\left(\nu \frac{u_{y}\phi_{yy}-u_{yy}\phi_{y}}{u_{y}} -u\phi_{x}+(u^{2}-U^{2})\frac{U_{x}}{U}\frac{\phi_{y}}{u_{y}} \right)dxdy
 -\int_{0}^{X}v_{0}(x)\phi(x,0)dx-\int_{0}^{\infty}u_{0}(y)\phi(0,y)dy =0, \\
&& \int_{D}\left(\nu u_{yy}\phi +\frac{u^{2}}{2}\phi_{x} -v u_{y}\phi+UU_{x}\phi \right)dxdy+ \int_{0}^{\infty}\frac{u_{0}(y)^{2}}{2}\phi(0,y)dy=0,
\end{eqnarray*}
for any $\phi(x,y)$ such that $\phi_{x}, \phi_{y},\phi_{yy}$ continuous in  $\overline{D}$ and $\phi(X,y)=0$,
$\phi_{y}(x,0)=0, |\phi|\leq M_{29}\exp(-M_{30}y^{2})$ as $y\rightarrow \infty$.
A solution $(u, v)$ of problem (\ref{ww16}), (\ref{ww20}) with the above properties is unique.

\end{Theorem}

 Oleinik and  Samokhin (\cite{OS}) in 1999
reached  boundary layer with unknown border between two media, i.e., they considered a stationary plane-parallel incompressible flow past a porous surface.
 Assume that another incompressible fluid, whose properties differ from those of the main stream, is injected into the boundary layer through the
porous wall. \\

The velocity components of the injected fluid $u_{1}(x,y), v_{1}(x,y)$ satisfy the Prandtl system
 \begin{equation}\left\{
\begin{array}{ll}
  u_{1}\frac{\partial u_{1}}{\partial x}+v_{1}\frac{\partial u_{1}}{\partial y}
=  \alpha_{1}\frac{\partial^{2} u_{1}}{\partial y^{2}} -\frac{1}{\rho_{1}}\frac{\partial p }{\partial x} ,  \\
\frac{\partial u_{1}}{\partial x}+\frac{\partial v_{1}}{\partial y} =0
\end{array}
 \label{ww51}         \right.\end{equation}
in the domain $D_{1}=\{0<x<X, 0<y<y_{*}(x)\}$ , where $\alpha_{1}$ is the viscosity of the injected fluid and $\rho_{1}$ is its density.

In the domain  $D_{2}=\{0<x<X, y_{*}(x)<y<\infty\}$ , another system is given, namely,
 \begin{equation}\left\{
\begin{array}{ll}
  u_{2}\frac{\partial u_{2}}{\partial x}+v_{2}\frac{\partial u_{2}}{\partial y}
=  \alpha_{2}\frac{\partial^{2} u_{2}}{\partial y^{2}} -\frac{1}{\rho_{2}}\frac{\partial p }{\partial x} ,  \\
\frac{\partial u_{1}}{\partial x}+\frac{\partial v_{1}}{\partial y} =0
\end{array}
 \label{ww52}         \right.\end{equation}
where $u_{2}(x,y), v_{2}(x,y)$  are the velocity components of the main fluid, $\alpha_{2}$ is its viscosity, and $\rho_{2}$ is its density.

The boundary conditions are given in the form
 \begin{equation}
\begin{array}{ll}
u_{1}(0,y)=u_{10}(y), \ \  \text{for} \ \  0\leq y\leq y_{*}(0), \\
u_{2}(0,y)=u_{20}(y), \ \  \text{for} \ \ y_{*}(0)\leq y\leq \infty, \\
u_{1}(x,0)=0, \ \  v_{1}(x,0)=v_{0}(x)>0, \\
u_{2}(x,y)\rightarrow U(x),  \ \ \text{as}\ \   y\rightarrow\infty,
\end{array}
 \label{ww53}         \end{equation}
where $y_{*}(0)$ is regarded as given, and $2p(x)+\rho_{2}U^{2}(x)=$const.
Assume that the fluids do not penetrate the interface between the two media. In this case, the curve $y=y_{*}(x)$ coincides with a streamline and,
therefore,
\begin{eqnarray}
\frac{d y_{*}(x)}{dx}=\frac{v_{1}(x,y_{*}(x) )}{u_{1}(x,y_{*}(x) )}=\frac{v_{2}(x,y_{*}(x) )}{u_{2}(x,y_{*}(x) )} .
 \label{ww54}
\end{eqnarray}
On the line $y=y_{*}(x)$, transmission conditions are imposed, which amount to the continuity of the velocity and of the viscous forces, namely,
\begin{eqnarray}
u_{1}(x,y_{*}(x) )=u_{2}(x,y_{*}(x) ), \label{ww55} \\
\nu_{1}\rho_{1}\frac{\partial u_{1}(x,y_{*}(x) )}{\partial y} =\nu_{2}\rho_{2}\frac{\partial u_{2}(x,y_{*}(x) )}{\partial y}.
 \label{ww56}
\end{eqnarray}
It follows from (\ref{ww54}) and (\ref{ww55}) that $v_{1}(x,y_{*}(x) )=v_{2}(x,y_{*}(x) )$.

\begin{Theorem} \label{OS.7.40} (\cite{OS})
For some $X>0$, problem (\ref{ww51})-(\ref{ww56})  has a solution $u_{1}, v_{1},u_{2},v_{2},y_{*}(x)$ with the following properties: $u_{1}$ is continuous in
$\overline{D}_{1}$, $u_{1}>0$ for $y>0$; $\partial u_{1}/\partial y>m_{1}>0$ for $0<y\leq y_{0}$ with constant $m_{1}$ and $y_{0}$;
$\partial u_{1}/\partial y, \partial^{2} u_{1}/\partial y^{2}$ are continuous and bounded in $D_{1}$; $\partial u_{1}/\partial x, v_{1}, \partial v_{1}/\partial y$
 are bounded in $D_{1}$ and continuous at $y=0$; $u_{2}$ is continuous and bounded in $\overline{D}_{2}$; $u_{2}>0$ in  $\overline{D}_{2}$;
$\partial u_{2}/\partial y, \partial^{2} u_{2}/\partial y^{2}$  are continuous and bounded in $D_{2}$;
 $\partial u_{2}/\partial x, v_{2}, \partial v_{2}/\partial y$  are continuous and bounded in any finite subdomain of $D_{2}$;
$\partial u_{i}/\partial y, v_{i}$ are continuous in $D_{i}\ (i = 1,2)$ up to the curve $y=y_{*}(x)$; if $u'_{20}(y)\rightarrow 0$ as $y\rightarrow\infty$
and this convergence is fast enough, so that $u'_{20}(y)\leq m_{2}\exp(-m_{3}y^{2}), m_{2}, m_{3}>0$, then $\partial u_{2}/\partial x, \partial v_{2}/\partial y$
are bounded in $D_{2}$, the function $ y_{*}(x)$ is continuously differentiable on $[0,X]$. If $dp/dx<0$, then such a solution of problem
(\ref{ww51})-(\ref{ww56})  exists in $D$ for any $X>0$.

\end{Theorem}

\begin{Theorem} \label{OS.7.41} (\cite{OS})
The solution of problem (\ref{ww51})-(\ref{ww56})  with the properties specified in Theorem \ref{OS.7.40} is unique.

\end{Theorem}

 Oleinik and  Samokhin (\cite{OS}) in 1999 considered
 the stationary mixing layer at the interface between two fluids with distinct properties. Assume that the interface coincides
with the straight line $y=0$.
 In this case, the fluid motion in the domain $D_{1}=\{0<x<X, 0<y<+\infty \}$ is described by the system (\ref{ww51}) with the boundary conditions
 \begin{equation}
\begin{array}{ll}
u_{1}(0,y)=u_{10}(y), \ \ v_{1}(x,0)=0, \\
u_{1}(x,y)\rightarrow U_{1}(x),  \ \ \text{as}\ \   y\rightarrow+\infty,
\end{array}
 \label{ww57}         \end{equation}
where $u_{1}(x,y), v_{1}(x,y)$ are the velocity components; $\nu_{1}$ is the viscosity of the fluid, and $\rho_{1}$ is its density.

In the domain  $D_{2}=\{0<x<X, -\infty<y<0 \}$, the velocity field is described by the system (\ref{ww52}) with the boundary conditions
 \begin{equation}\left\{
\begin{array}{ll}
u_{2}(0,y)=u_{20}(y), \ \ v_{2}(x,0)=0, \\
u_{2}(x,y)\rightarrow U_{2}(x),  \ \ \text{as}\ \   y\rightarrow-\infty,
\end{array}
 \label{ww58}    \right.     \end{equation}
where $u_{2}(x,y), v_{2}(x,y)$ are the velocity components; $\nu_{2}$ is the viscosity of the fluid, and $\rho_{2}$ is its density.

For $y=0$, the following transmission conditions are imposed:
 \begin{equation}\left\{
\begin{array}{ll}
u_{1}(x,0)=u_{2 }(x, 0),  \\
\nu_{1}\rho_{1} u_{1y}(x,0)= \nu_{2}\rho_{2} u_{2y}(x,0).
\end{array}
 \label{ww59}      \right.   \end{equation}

 \begin{Theorem} \label{OS.7.42} (\cite{OS})
Under the compatibility assumptions, problem (\ref{ww51})-(\ref{ww52}), (\ref{ww57})-(\ref{ww59}),
for some $X>0$, admits a solution $u_{1}, v_{1},u_{2},v_{2} $ with the following properties: $u_{i}$ is continuous and bounded in $\overline{D}_{i}$;
$u_{i}>0$; $\partial u_{i}/\partial y, \partial^{2} u_{i}/\partial y^{2}$ are continuous and bounded in $D_{i}$;
$\partial u_{i}/\partial x, v_{i}, \partial v_{i}/\partial y$ are continuous and bounded in any finite subdomain of $D_{i}$;
$\partial u_{i}/\partial y, v_{i}$ are continuous in $D_{i}$ up to the line $y=0$; $i=1,2$. If
$$|u'_{0}(y)|\leq k_{1}\exp(-k_{2}|y|) \ \  \text{as}\ \  |y|\rightarrow \infty$$
with constant $k_{1} ,k_{2}>0$, then $\partial u_{i}/\partial x, \partial v_{i}/\partial y$ are bounded in $D_{i}$.
If $dp/dx<0$, the solution of problem (\ref{ww51})-(\ref{ww52}), (\ref{ww57})-(\ref{ww59}) exists for any $X>0$.
The solution of problem (\ref{ww51})-(\ref{ww52}), (\ref{ww57})-(\ref{ww59}) with these properties is unique.

\end{Theorem}


 Oleinik and  Samokhin (\cite{OS}) in 1999 considered
  the stationary Prandtl system for a planar incompressible flow by the von Mises, the stationary Prandtl system
\begin{equation}\left\{
\begin{array}{ll}
   u \frac{\partial u }{\partial x}+v \frac{\partial u }{\partial y}
=  \nu  \frac{\partial^{2} u }{\partial y^{2}}  +U(x)\frac{d U }{d x}, \\
\frac{\partial u }{\partial x}+\frac{\partial v }{\partial y} =0,
\end{array}
 \label{ww92}         \right.\end{equation}
 in the domain $D =\{  0<x<X, 0<y< \infty \}$, with the boundary conditions
 \begin{equation}\left\{
\begin{array}{ll}
  u ( 0,y)=u_{1}(y), \ \ u (x,0)=0, \ \ v (x,0)=v_{0}( x), \\
u ( x,y)\rightarrow U ( x),  \ \ \text{as}\ \   y\rightarrow \infty.
\end{array}
 \label{ww93}    \right.     \end{equation}

\begin{Theorem} \label{OS.10.59} (\cite{OS})
Suppose that in (\ref{ww92})-(\ref{ww93}), we have $U_{x}(x)>0$ for $0\leq x\leq X$, and $U(x), u_{1}(y), v_{0}(x)\equiv v_{0}(\varepsilon^{-1}x)$
 (for any $\varepsilon\in (0,\varepsilon_{0}]$) satisfy the assumptions that $u_{0}(y)$ for $y>0$; $u_{0}(0)=0$, $u_{0}'(0)>0$,
 $u_{0}(y)\rightarrow U(0)\neq 0$ as $y\rightarrow \infty$;
$p_{x}$ and $v_{0}(x)$ are continuously differentiable on [0,X];
$u_{0}(y), u_{0}'(y), u_{0}''(y)$ are bounded for $0\leq y\leq \infty$ and satisfy the H$\ddot{o}$lder condition.
Moreover, assume that for small $y$ the following compatibility condition is satisfied at the point $(0,0)$:
\begin{eqnarray*}
\frac{d^{2}u_{0}(y)}{d y^{2}}-\frac{dp(0)}{dx}-v_{0}\frac{du_{0}(y)}{dy}=\mathcal{O}(y^{2}).
\end{eqnarray*}
Then the solution $u_{\varepsilon}, v_{\varepsilon}$ of problem (\ref{ww92})-(\ref{ww93})  converges, as $\varepsilon\rightarrow 0$, to a solution of the system
\begin{equation*}\left\{
\begin{array}{ll}
   \overline{u} \frac{\partial \overline{u} }{\partial x}+\overline{v} \frac{\partial \overline{u} }{\partial y}
=  \nu  \frac{\partial^{2} \overline{u} }{\partial y^{2}}  +U(x)\frac{d U }{d x}, \\
\frac{\partial u }{\partial x}+\frac{\partial v }{\partial y} =0,
\end{array}\right.
 \end{equation*}
in the domain $D$, with the conditions
 \begin{equation*}\left\{
\begin{array}{ll}
 \overline{ u} ( 0,y)=u_{1}(y), \ \ \overline{u} (x,0)=0, \ \ \overline{v} (x,0)=v_{0}( x), \\
\overline{u} ( x,y)\rightarrow U ( x),  \ \ \text{as}\ \   y\rightarrow \infty,
\end{array}\right.
       \end{equation*}
so that
\begin{eqnarray*}\left\{
\begin{array}{ll}
 u_{\varepsilon}\rightarrow \overline{u} \ \ \text{uniformly for} \ \ 0\leq x\leq X, \ \  0\leq y\leq y_{1}<\infty, \\
 v_{\varepsilon}\rightarrow \overline{v} \ \ \text{uniformly for} \ \ x_{0}\leq x\leq X, \ \  0\leq y\leq y_{1} , \\
 \text{grad}~u_{\varepsilon}\rightarrow \text{grad}~\overline{u} \ \ \text{weakly in} \ \  L_{2}(\{0\leq x\leq X, 0\leq y\leq y_{1}\}) ,
\end{array}\right.
\end{eqnarray*}
for arbitrary fixed $y_{1}\in (0,\infty), x_{0}\in (0,X)$.

\end{Theorem}

In the case of favorable pressure gradient, Oleinik proved the global existence
of classical solution for the 2D steady Prandtl equations for a class of positive data. In
the case of adverse pressure gradient, an important physical phenomena is the boundary
layer separation. Shen, Wang and Zhang   (\cite{SWZ}) in 2019 proved the boundary layer separation for a large class
of Oleinik's data and confirmed Goldstein's hypothesis concerning the local behavior of the
solution near the separation.
The authors studied the 2D steady Prandtl equations
\begin{equation}
\left\{\begin{array}{l}
u \partial_{x} u+v \partial_{y} u-\partial_{y}^{2} u+\partial_{x} p=0, \quad x \geq 0, y \geq 0 ,\\
\partial_{x} u+\partial_{y} v=0, \\
\left.u\right|_{y=0}=\left.v\right|_{y=0}=0 ,\quad \text { and } \quad \lim\limits _{y \rightarrow+\infty} u(x, y)=U(x),
\end{array}\right.
\end{equation}
where the outer flow $(U(x), p(x))$ satisfies
$$
U(x) U^{\prime}(x)+p^{\prime}(x)=0.
$$
This system derived by Prandtl is used to describe the behavior of the solution near $y=0$ for the steady Navier-Stokes equations with non slip boundary condition when the viscosity coefficient $\nu$ tends to zero:
$$
\left\{\begin{array}{l}
u^{\nu} \cdot \nabla u^{\nu}-\nu \Delta u^{\nu}+\nabla p=f^{\nu}, \\
\operatorname{div} u^{\nu}=0 ,\\
\left.u^{\nu}\right|_{y=0}=0.
\end{array}\right.
$$
Roughly speaking, away from the boundary, the solution $u^{\nu}$ can be described by the Euler equations; near the boundary $y=0, u^{\nu}$ behaves as
$$
u^{\nu}(x, y)=(u(x, y / \sqrt{\nu}), \sqrt{\nu} v(x, y / \sqrt{\nu})),
$$
where $(u, v)$ satisfies the Prandtl type equation.

The authors of \cite{SWZ} first introduced a class of data denoted by $\mathcal{K},$ which satisfies
$$
u \in C_{b}^{3, \alpha}\ ([0,+\infty))~(\alpha>0), \quad u(0)=0, u_{y}(0)>0, u_{y}(y) \geq 0 \text { for } y \in[0,+\infty),
$$

$$\lim _{y \rightarrow+\infty} u(0, y)=U(0)>0, \quad u_{y y}(y)-\partial_{x} p(0)=O\left(y^{2}\right).$$
The data $u_{0} \in \mathcal{K}$.
\begin{Proposition}{\label{lemma1}}
If $u_{0} \in \mathcal{K},$ then there exists $X>0$ such that the steady Prandtl equation admits a solution $u \in C^{1}\left([0, X) \times \mathbb{R}_{+}\right)$ with the following properties:
1. Regularity: $u$ is bounded and continuous in $[0, X] \times[0,+\infty) ; u_{y}, u_{y y}$ are bounded and continuous in $[0, X) \times \mathbb{R}_{+} ;$ and $v, v_{y}, u_{x}$ are locally bounded and continuous in $[0, X) \times \mathbb{R}_{+}$.

2. Non-degeneracy: $u(x, y)>0$ in $[0, X) \times(0,+\infty)$ and for all $\bar{x}<X,$ there exists $y_{0}>0, m>0$ so that
$\partial_{y} u(x, y) \geq m  $ in $[0, \bar{x}] \times\left[0, y_{0}\right]$.

3. Global existence: if $p^{\prime}(x) \leq 0,$ then the solution is global in $x$,
together with
$$
w(0, \psi)=w_{0}(\psi)=u_{0}(y)^{2}, \quad w(x, 0)=0, \quad \lim _{\psi \rightarrow+\infty} w(x, \psi)=U(x)^{2}.
$$
\end{Proposition}
\begin{Theorem}(\cite{SWZ})
Fix any $\mu \in(0,1) .$ Let $u$ be a solution   constructed in  Proposition \ref{lemma1} with $u_{0} \in \mathcal{K}$ satisfying
$$
\left\|\partial_{y} u_{0}\right\|_{L^{\infty}\left(0, y_{0}\right)} \leq \frac{1}{2} \epsilon_{0} x_{0}^{\frac{1}{4}},
$$
where $y_{0}$ is determined by
$$
B x_{0}^{\frac{3}{4}}=\psi_{0}=\int_{0}^{y_{0}} u_{0}(z) d z,
$$
and $\epsilon_{0}, B$ are positive constants depending only on $\mu .$ Then there exists a separation point $x_{s}=X^{*}$ with $X^{*}<\frac{\mu}{2} x_{0}$.

\end{Theorem}

Wang and Zhang (\cite{WZ}) in  2019 proved the global $C^{\infty}$ regularity of the Oleinik's solution for the steady Prandtl equation with favorable pressure.
The authors  studied the steady Prandtl equations:
\begin{equation}
\left\{\begin{array}{l}
u \partial_{x} u+v \partial_{y} u-\partial_{y}^{2} u=-\frac{d p}{d x}(x), \quad x \geq 0, y \geq 0, \\
\partial_{x} u+\partial_{y} v=0, \\
\left.u\right|_{y=0}=\left.v\right|_{y=0}=0, \quad \text { and } \quad \lim _{y \rightarrow+\infty} u(x, y)=U(x),
\end{array}\right.
\end{equation}
where the outer flow $(U(x), p(x))$ satisfies
$$
U(x) U^{\prime}(x)+p^{\prime}(x)=0.
$$
This system derived by Prandtl could be used to describe the behavior of the solution near $y=0$ for the steady Navier-Stokes equations when the viscosity coefficient is small.

Oleinik ( see Theorem 2.1.1 in \cite{OS}) proved the following classical result:

If $u_{0} \in \mathcal{K}$ and $\frac{d p}{d x}(x)$ is smooth, then there exists $X>0$ such that the steady Prandtl equation (1.1) admits a solution $u \in C^{1}\left([0, X) \times \mathbb{R}_{+}\right)$ with the following properties:

1. Regularity: $u$ is bounded and continuous in $[0, X] \times \mathbb{R}_{+} ; u_{y}, u_{y y}$ are bounded and continuous in $[0, X) \times \mathbb{R}_{+} ;$ and $v, v_{y}, u_{x}$ are locally bounded and continuous in $[0, X) \times \mathbb{R}_{+} .$

2. Non-degeneracy: $u(x, y)>0$ in $[0, X) \times(0,+\infty)$ and for all $\bar{x}<X$, there exists $y_{0}>0, m>0$ so that
$$
\partial_{y} u(x, y) \geq m \quad \text { in }[0, \bar{x}] \times\left[0, y_{0}\right].
$$

3. Global existence: if $p^{\prime}(x) \leq 0,$ then the solution is global in $x$.

\begin{Theorem}(\cite{WZ})
Let $u$ be a global solution to the equation constructed by Oleinik with $u_{0} \in \mathcal{K}$ and $\frac{d p(x)}{d x} \leq 0$ smooth. For any positive integers $m, k$ and any positive constants $X, Y$ with $\epsilon<X,$ there exists a positive constant $C$ depending only on $\epsilon, X, Y, u_{0}, p, k, m$ so that
$$
\left|\partial_{x}^{k} \partial_{y}^{m} u(x, y)\right| \leq C \quad \text { for } \quad(x, y) \in[\epsilon, X] \times[0, Y].
$$
\end{Theorem}

Dalibard and Masmoudi  (\cite{DM}) in 2019 were interested  in the stationary version of the Prandtl equation, namely
\begin{equation}
\left\{\begin{array}{c}
u u_{x}+v u_{y}-u_{y y}=-\frac{d p_{E}(x)}{d x}, \quad x>0, y>0, \\
u_{x}+v_{y}=0, \quad x>0, y>0, \\
u|_{  x=0}=u_{0}, \quad u|_{  y=0}=0, \lim\limits _{y \rightarrow \infty} u(x, y)=u_{E}(x),
\end{array}\right.\label{y2.11}
\end{equation}
where $y=0$ stands for the rigid wall, $x$ (resp. $y$ ) is the tangential (resp. normal) variable to the wall. The functions $u_{E}, p_{E}$ are given by the outer flow: more precisely $u_{E}$ (resp. $\left.p_{E}\right)$ is the trace at the boundary of the tangential velocity (resp. of the pressure) of a flow satisfying the Euler equations. The functions $u_{E}, p_{E}$ are linked by the relation
$$
u_{E} u_{E}^{\prime}=-\frac{d p_{E}(x)}{d x}.
$$

Let $\alpha>0, X \in ( 0,+\infty] .$ Let $u_{0} \in \mathcal{C}_{b}^{2, \alpha}(\mathbb{R})$ such that $u_{0}(0)=0$, $u_{0}^{\prime}(0)>0, \lim _{y \rightarrow \infty} u_{0}(y)=u_{E}(0)>0,$ and such that $u_{0}(y)>0$ for $y>0$. Assume that
$d p_{E} / d x \in \mathcal{C}^{1}([0, X]),$ and that for $y \ll 1$ the following compatibility condition is satisfied
$$
u_{0}^{\prime \prime}(y)-\frac{d p_{E}(0)}{d x}=O\left(y^{2}\right).
$$

The authors in \cite{DM} are  interested in the case where the solution of (\ref{y2.11}) is not global: more precisely, they considered the equation with $d p_{E} / d x=1,$ i.e.,
$$
\left\{\begin{array}{r}
u u_{x}+v u_{y}-u_{y y}=-1, \quad x \in\left(0, x_{0}\right), y>0, \\
u_{x}+v_{y}=0, \quad x \in\left(0, x_{0}\right), y>0, \\
u_{\mid x=0}=u_{0}, \quad u_{\mid y=0}=0, \lim _{y \rightarrow \infty} u(x, y)=u_{E}(x),
\end{array}\right.
$$
with $u_{E}(x)=\sqrt{2\left(x_{0}-x\right)},$ for some $x_{0}>0,$ and $u_{0}$ satisfies the assumptions.

The authors's assumptions on the initial data $u_{0}$ are the following:

(H1) $u_{0} \in \mathcal{C}^{7}\left(\mathbb{R}_{+}\right), u_{0}$ is increasing in $y$ and $\lambda_{0}:=u_{0}^{\prime}(0)>0$,

(H2) There exists a constant $C_{0}>0$ such that
$$
\begin{aligned}
\forall y \geq 0, \quad-C_{0} \inf \left(y^{2}, 1\right) \leq u_{0}^{\prime \prime}(y)-1 \leq 0, \\
C_{0}^{-1} \leq-u_{0}^{(4)}(0) \leq C_{0}, \\
\left\|u_{0}\right\|_{W^{7, \infty}} \leq C_{0} .
\end{aligned}
$$

(H3) $u_{0}=u_{0}^{\text {app }}+v_{0},$ where
$$
\begin{aligned}
u_{0}^{\mathrm{app}}=& \lambda_{0} y+\frac{y^{2}}{2}+u_{0}^{(4)}(0) \frac{y^{4}}{4 !} \\
&-c_{7}\left(u_{0}^{(4)}(0)\right)^{2} \frac{y^{7}}{\lambda_{0}}+c_{10}\left(u_{0}^{(4)}(0)\right)^{3} \frac{y^{10}}{\lambda_{0}^{2}}+c_{11}\left(u_{0}^{(4)}(0)\right)^{3} \frac{y^{11}}{\lambda_{0}^{3}} \quad \text { for } y \leq \lambda_{0}^{3 / 7} ,\\
\left|u_{0}^{\mathrm{app}}\right| \leq & C_{0} \text { for } y \geq \lambda_{0}^{3 / 7},
\end{aligned}
$$
and
$$
\left|v_{0}\right| \leq C_{0}\left(\lambda_{0}^{-\frac{3}{2}}\left(\lambda_{0} y^{7}+c_{8} y^{8}\right)+\lambda_{0}^{-2} y^{10}+\lambda_{0}^{-3} y^{11}\right) \quad \text { for } y \leq \lambda_{0}^{3 / 7}.
$$
\begin{Theorem}(\cite{DM})
Consider the Prandtl equation with adverse pressure gradient (P) and with an initial data $u_{0} \in \mathcal{C}^{7}\left(\mathbb{R}_{+}\right)$ satisfying (H1)-(H3). Then for any $\eta>0, C_{0}>0,$ there exists $\epsilon_{0}>0$ such that if $\lambda_{0}<\epsilon_{0},$ the existence time $x^{*}$ is finite, and $x^{*}=O\left(\lambda_{0}^{2}\right) .$ Furthermore, setting $\lambda(x):=\partial_{y} u|_{y=0}(x),$ there exists a constant $C>0,$ depending on $u_{0},$ such that
$$
\lambda(x) \sim C \sqrt{x^{*}-x} \quad \text { as } x \rightarrow x^{*}.
$$
\end{Theorem}

Gao and Zhang  (\cite{GZ}) in 2020 considered the zero-viscosity limit of the $2 \mathrm{D}$ steady NavierStokes equations in $(0, L) \times \mathbb{R}^{+}$ with non-slip boundary conditions. By estimating the stream-function of the remainder, we justify the validity of the Prandtl boundary layer expansions.
The authors considered the vanishing viscosity limit of steady Navier-Stokes equations
\begin{equation}
\left\{\begin{array}{l}
U^{\varepsilon} U_{X}^{\varepsilon}+V^{\varepsilon} U_{Y}^{\varepsilon}-\varepsilon \Delta U^{\varepsilon}+P_{X}^{\varepsilon}=0, \\
U^{\varepsilon} V_{X}^{\varepsilon}+V^{\varepsilon} V_{Y}^{\varepsilon}-\varepsilon \Delta V^{\varepsilon}+P_{Y}^{\varepsilon}=0 ,\\
U_{X}^{\varepsilon}+V_{Y}^{\varepsilon}=0, \\
\left.U^{\varepsilon}\right|_{Y=0}=\left.V^{\varepsilon}\right|_{Y=0}=0,
\end{array}\right.\label{+5}
\end{equation}
in a two dimensional domain $\Omega=\{(X, Y): 0 \leqslant X \leqslant L, Y \geqslant 0\} .$ A formal limit $\varepsilon \rightarrow 0$ should lead to the Euler flow $\left[U^{0}, V^{0}\right]$ inside $\Omega$ :
\begin{equation}
\left\{\begin{array}{l}
U^{0} U_{X}^{0}+V^{0} U_{Y}^{0}+P_{X}^{0}=0, \\
U^{0} V_{X}^{0}+V^{0} V_{Y}^{0}+P_{Y}^{0}=0 ,\\
U_{X}^{0}+V_{Y}^{0}=0, \\
\left.V^{0}\right|_{Y=0}=0.
\end{array}\right.
\end{equation}

In their first result, they assumed that the outside Euler flow $\left[U^{0}, V^{0}\right] \equiv\left(u_{e}^{0}(X, Y), v_{e}^{0}(X, Y)\right)$ satisfying the following hypothesis:
$$
\begin{aligned}
&0<c_{0} \leqslant u_{e}^{0} \leqslant C_{0}<\infty,\\
&\left\|v_{e}^{0}\right\|_{L^{\infty}} \ll 1,\\
&\left\|\langle Y\rangle^{k} \nabla^{m} v_{e}^{0}\right\|_{L^{\infty}}<\infty \text { for sufficiently large } k, m \geqslant 0,\\
&\left\|\langle Y\rangle^{k} \nabla^{m} u_{e}^{0}\right\|_{L^{\infty}}<\infty \text { for sufficiently large } k, m \geqslant 1.
\end{aligned}
$$
Here $\langle Y\rangle=Y+1$.
The authors  in  \cite{GZ} considered   the Prandtl equations with the positive data
$$
\left\{\begin{array}{l}
u_{p}^{0} u_{p x}^{0}+v_{p}^{0} u_{p y}^{0}-u_{p y y}^{0}+p_{p x}^{0}=0, \quad p_{p y}^{0}=0, \quad u_{p x}^{0}+v_{p y}^{0}=0, \quad(x, y) \in(0, L) \times \mathbb{R}_{+}, \\
\left.u_{p}^{0}\right|_{x=0}=U_{P}^{0}(y),\left.\quad u_{p}^{0}\right|_{y=0}=\left.v_{p}^{0}\right|_{y=0}=0,\left.\quad u_{p}^{0}\right|_{y \uparrow \infty}=\left.u_{e}^{0}\right|_{Y=0},
\end{array}\right.
$$
where $U_{P}^{0}$ is a prescribing smooth function such that
$$
U_{P}^{0}>0 \text { for } y>0, \quad \partial_{y} U_{P}^{0}(0)>0, \quad \partial_{y}^{2} U_{P}^{0}-u_{e}^{0} u_{e x}^{0}(0) \sim y^{2} \text { near } y=0,
$$
$\partial_{y}^{m}\left\{U_{P}^{0}-u_{e}^{0}(0)\right\}$ decay fast for any $m \geqslant 0$.
\begin{Theorem}(\cite{GZ})
Assume the Euler flows $\left[u_{e}^{0}, v_{e}^{0}\right]$ satisfy satisfy the former assumptions, $U_{P}^{0}$ is a smooth function satisfying and high order compatibility conditions, $L$ is a constant small enough, then there exist $C(L), \varepsilon_{0}(L)>0$ depending on $L,$ such that for $0<\varepsilon \leqslant \varepsilon_{0},$ equations \eqref{+5} admits a solution $\left[U^{\varepsilon}, V^{\varepsilon}\right] \in W^{2,2}(\Omega),$ satisfying:
\begin{equation}\left\{
\begin{array}{ll}
\|U^{\varepsilon}-u_{e}^{0}+u_{e}^{0}|_{Y=0}-u_{p}^{0}\|_{L^{\infty}} \leqslant C \sqrt{\varepsilon}, \\
\|V^{\varepsilon}-v_{e}^{0}\|_{L^{\infty}} \leqslant C \sqrt{\varepsilon},
\end{array}\right.\end{equation}
with the following boundary conditions:
\begin{equation}\left\{
\begin{array}{ll}
(U^{\varepsilon}, V^{\varepsilon})|_{Y=0}=0, \\
(U^{\varepsilon}, V^{\varepsilon})|_{X=0}=\left(u_e^0(0,Y)-u_e^0(0,0)+u_p^0(0,\frac{Y}{\sqrt\varepsilon})+\sqrt\varepsilon a_0,v_e^0(0,Y)+\sqrt\varepsilon b_0\right), \\
(U^{\varepsilon}, V^{\varepsilon})|_{X=L}=\left(u_{e}^{0}(L, Y)-u_{e}^{0}(L, 0)+u_{p}^{0}(L, \frac{Y}{\sqrt{\varepsilon}})+\sqrt{\varepsilon} a_{L}, v_{e}^{0}(L, Y)+\sqrt{\varepsilon} b_{L}\right).
\end{array}\right.\end{equation}
Here
\begin{equation}\left\{
\begin{array}{ll}
a_{0}(Y)=u_{e}^{1}(0, Y)+u_{b}^{1}\left(0, \frac{Y}{\sqrt{\varepsilon}}\right)+\sqrt{\varepsilon} u_{e}^{2}(0, Y)+\sqrt{\varepsilon} \hat{u}_{b}^{2}\left(0, \frac{Y}{\sqrt{\varepsilon}}\right),\\
a_{L}(Y)=u_{e}^{1}(L, Y)+u_{b}^{1}\left(L, \frac{Y}{\sqrt{\varepsilon}}\right)+\sqrt{\varepsilon} u_{e}^{2}(L, Y)+\sqrt{\varepsilon} \hat{u}_{b}^{2}\left(L, \frac{Y}{\sqrt{\varepsilon}}\right), \\
b_{0}(Y)=v_{b}^{0}\left(0, \frac{Y}{\sqrt{\varepsilon}}\right)+v_{e}^{1}(0, Y)+\sqrt{\varepsilon} v_{b}^{1}\left(0, \frac{Y}{\sqrt{\varepsilon}}\right)+\sqrt{\varepsilon} v_{e}^{2}(0, Y)+\varepsilon \hat{v}_{b}^{2}\left(0, \frac{Y}{\sqrt{\varepsilon}}\right) ,\\
b_{L}(Y)=v_{b}^{0}\left(L, \frac{Y}{\sqrt{\varepsilon}}\right)+v_{e}^{1}(L, Y)+\sqrt{\varepsilon} v_{b}^{1}\left(L, \frac{Y}{\sqrt{\varepsilon}}\right)+\sqrt{\varepsilon} v_{e}^{2}(L, Y)+\varepsilon \hat{v}_{b}^{2}\left(L, \frac{Y}{\sqrt{\varepsilon}}\right),
\end{array}\right.\end{equation}
are smooth functions.
\end{Theorem}

For the second result, they considered $L$ is any given positive constant. The authors  assumed the Euler flow $\left[U^{0}, V^{0}\right] \equiv\left[u_{e}^{0}(Y), 0\right]$ is a shear flow, that is, it satisfies the following hypothesis:
$0<c_{0} \leqslant u_{e}^{0} \leqslant C_{0}<\infty$,
$\left\|\langle Y\rangle^{k} \nabla^{m} u_{e}^{0}\right\|_{L^{\infty}}<\infty$ , for sufficiently large $k, m \geqslant 1$.
Here $\langle Y\rangle=Y+1$.

While the authors assumed $\left[u_{p}^{0}, v_{p}^{0}\right]$ is a smooth solution of Prandtl equations satisfying the following hypothesis:
$$
\begin{array}{l}
u_{p}^{0}>0,-u_{p y y}^{0} \geqslant 0, \text { for } y>0, \\
u_{p y}^{0}>0, \text { for } y \geqslant 0,
\end{array}
$$
and $\nabla^{m}\left\{u_{p}^{0}-u_{e}^{0}(0)\right\}$ decays fast for any $m \geqslant 0$.

\begin{Theorem}(\cite{GZ})
Assume the Euler flows $\left[u_{e}^{0}, v_{e}^{0}\right]$ satisfy the assumptions, the Prandtl profiles $\left[u_{p}^{0}, v_{p}^{0}\right]$ satisfy the assumptions, and $L>0$ is any given constant, then there exist $C(L), \varepsilon_{0}(L)>0$ depending on $L,$ such that for $0<\varepsilon \leqslant \varepsilon_{0},$ equations \eqref{+5} admits a solution $\left[U^{\varepsilon}, V^{\varepsilon}\right] \in W^{2,2}(\Omega)$, satisfying:
$$
\begin{array}{c}
\left\|U^{\varepsilon}-u_{e}^{0}+\left.u_{e}^{0}\right|_{Y=0}-u_{p}^{0}\right\|_{L^{\infty}} \leqslant C \sqrt{\varepsilon}, \\
\left\|V^{\varepsilon}\right\|_{L^{\infty}} \leqslant C \sqrt{\varepsilon},
\end{array}
$$
with the boundary conditions:
$$
\begin{array}{l}
{\left.\left[U^{\varepsilon}, V^{\varepsilon}\right]\right|_{Y=0}=0}, \\
{\left.\left[U^{\varepsilon}, V^{\varepsilon}\right]\right|_{X=0}=\left[u_{e}^{0}(Y)-u_{e}^{0}(0)+u_{p}^{0}\left(0, \frac{Y}{\sqrt{\varepsilon}}\right)+\sqrt{\varepsilon} a_{0}, \sqrt{\varepsilon} b_{0}\right]}, \\
{\left[U^{\varepsilon}, V^{\varepsilon}\right]_{X=L}=\left[u_{e}^{0}(Y)-u_{e}^{0}(0)+u_{p}^{0}\left(L, \frac{Y}{\sqrt{\varepsilon}}\right)+\sqrt{\varepsilon} a_{L}, \sqrt{\varepsilon} b_{L}\right]}.
\end{array}
$$
Here
$$
\begin{array}{l}
a_{0}(Y)=u_{e}^{1}(0, Y)+u_{b}^{1}\left(0, \frac{Y}{\sqrt{\varepsilon}}\right)+\sqrt{\varepsilon} u_{e}^{2}(0, Y)+\sqrt{\varepsilon} \hat{u}_{b}^{2}\left(0, \frac{Y}{\sqrt{\varepsilon}}\right),\\
a_{L}(Y)=u_{e}^{1}(L, Y)+u_{b}^{1}\left(L, \frac{Y}{\sqrt{\varepsilon}}\right)+\sqrt{\varepsilon} u_{e}^{2}(L, Y)+\sqrt{\varepsilon} \hat{u}_{b}^{2}\left(L, \frac{Y}{\sqrt{\varepsilon}}\right) ,\\
b_{0}(Y)=v_{b}^{0}\left(0, \frac{Y}{\sqrt{\varepsilon}}\right)+v_{e}^{1}(0, Y)+\sqrt{\varepsilon} v_{b}^{1}\left(0, \frac{Y}{\sqrt{\varepsilon}}\right)+\sqrt{\varepsilon} v_{e}^{2}(0, Y)+\varepsilon \hat{v}_{b}^{2}\left(0, \frac{Y}{\sqrt{\varepsilon}}\right), \\
b_{L}(Y)=v_{b}^{0}\left(L, \frac{Y}{\sqrt{\varepsilon}}\right)+v_{e}^{1}(L, Y)+\sqrt{\varepsilon} v_{b}^{1}\left(L, \frac{Y}{\sqrt{\varepsilon}}\right)+\sqrt{\varepsilon} v_{e}^{2}(L, Y)+\varepsilon \hat{v}_{b}^{2}\left(L, \frac{Y}{\sqrt{\varepsilon}}\right),
\end{array}
$$
are smooth functions constructed.
\end{Theorem}

Iyer and Masmoudi  (\cite{IM}) in 2020 established the convergence of $2 \mathrm{D},$ stationary Navier-Stokes flows, $\left(u^{\varepsilon}, v^{\varepsilon}\right)$ to the classical Prandtl boundary layer, $\left(\bar{u}_{p}, \bar{v}_{p}\right) $ posed on the domain $(0, \infty) \times(0, \infty):$
$$
\left\|u^{\varepsilon}-\bar{u}_{p}\right\|_{L_{y}^{\infty}} \lesssim \sqrt{\varepsilon}\langle x\rangle^{-\frac{1}{4}+\delta}, \quad\left\|v^{\varepsilon}-\sqrt{\varepsilon} \bar{v}_{p}\right\|_{L_{y}^{\infty}} \lesssim \sqrt{\varepsilon}\langle x\rangle^{-\frac{1}{2}}.
$$
This validates Prandtl's boundary layer theory globally in the $x$ -variable for a large class of boundary layers, including the entire one parameter family of the classical Blasius profiles, with sharp decay rates.

Consider the Navier-Stokes (NS) equations posed on the domain $\mathcal{Q}:=(0, \infty) \times(0, \infty)$ :
\begin{equation}
\left\{\begin{array}{l}
u^{\varepsilon} u_{x}^{\varepsilon}+v^{\varepsilon} u_{Y}^{\varepsilon}+P_{x}^{\varepsilon}=\varepsilon \Delta u^{\varepsilon} ,\\
u^{\varepsilon} v_{x}^{\varepsilon}+v^{\varepsilon} v_{Y}^{\varepsilon}+P_{Y}^{\varepsilon}=\varepsilon \Delta v^{\varepsilon}, \\
u_{x}^{\varepsilon}+v_{Y}^{\varepsilon}=0.
\end{array}\right.
\end{equation}
The authors in \cite{IM} took the following boundary conditions in the vertical direction
$$
\left.\left[u^{\varepsilon}, v^{\varepsilon}\right]\right|_{Y=0}=[0,0], \quad\left[u^{\varepsilon}(x, Y), v^{\varepsilon}(x, Y)\right] \stackrel{Y \rightarrow \infty}{\longrightarrow}\left[u_{E}(x, \infty), v_{E}(x, \infty)\right],
$$
which coincide with the classical no-slip boundary condition at $\{Y=0\}$ and the Euler matching condition as $Y \rightarrow \infty$. The authors in (\cite{IM}) fixed the vector field
$$
\left[u_{E}, v_{E}\right]:=[1,0], \quad P_{E}=0,
$$
as a solution to the steady, Euler equations condition above reads $\left[u^{\varepsilon}, v^{\varepsilon}\right] \stackrel{Y \rightarrow \infty}{\longrightarrow}[1,0]$.
They expanded the rescaled solution as
\begin{equation*}
\left\{\begin{array}{l}
U^{\varepsilon}:=1+u_{p}^{0}+\sum_{i=1}^{N_{1}} \varepsilon^{\frac{i}{2}}\left(u_{E}^{i}+u_{P}^{i}\right)+\varepsilon^{\frac{N_{2}}{2}} u=: \bar{u}+\varepsilon^{\frac{N_{2}}{2}} u, \\
V^{\varepsilon}:=v_{p}^{0}+v_{E}^{1}+\sum_{i=1}^{N_{1}-1} \varepsilon^{\frac{i}{2}}\left(v_{P}^{i}+v_{E}^{i+1}\right)+\varepsilon^{\frac{N_{1}}{2}} v_{p}^{N_{1}}+\varepsilon^{\frac{N_{2}}{2}} v=: \bar{v}+\varepsilon^{\frac{N_{2}}{2}} v, \\
P^{\varepsilon}:=\sum_{i=0}^{N_{1}+1} \varepsilon^{\frac{i}{2}} P_{p}^{i}+\sum_{i=1}^{N_{1}} \varepsilon^{\frac{i}{2}} P_{E}^{i}+\varepsilon^{\frac{N_{2}}{2}} P.
\end{array}\right.
\end{equation*}
Above,
$$
\left[u_{E}^{0}, v_{E}^{0}\right]:=[1,0], \quad\left[u_{p}^{i}, v_{p}^{i}\right]=\left[u_{p}^{i}(x, y), v_{p}^{i}(x, y)\right], \quad\left[u_{E}^{i}, v_{E}^{i}\right]=\left[u_{E}^{i}(x, Y), v_{E}^{i}(x, Y)\right],
$$
and $N_{1}, N_{2}$ are the  the expansion parameters for the sake of precision.

Introduce the  Blasius boundary layer, which takes the following form
\begin{equation}
\left\{\begin{array}{l}
{\left[\bar{u}_{*}^{x_{0}}, \bar{v}_{*}^{x_{0}}\right]=\left[f^{\prime}(z), \frac{1}{\sqrt{x+x_{0}}}\left(z f^{\prime}(z)-f(z)\right)\right]} ,\\
z:=\frac{y}{\sqrt{x+x_{0}}}, \\
f f^{\prime \prime}+f^{\prime \prime \prime}=0, \quad f^{\prime}(0)=0, \quad f^{\prime}(\infty)=1, \quad \frac{f(z)}{z} \stackrel{\eta \rightarrow \infty}{\longrightarrow} 1,
\end{array}\right.
\end{equation}
where above, $f^{\prime}=\partial_{z} f(z)$ and $x_{0}$ is a free parameter. Physically, $x_{0}$ has the meaning that at $x=-x_{0},$ the fluid interacts with the leading edge of, say, a plate (hence the singularity at $\left.x=-x_{0}\right)$.
The authors's complete norm will be
$$
\|U, V\|_{\mathcal{X}}:=\sum_{n=0}^{10}\left(\|U, V\|_{X_{n}}+\|U, V\|_{X_{n+\frac{1}{2}}}+\|U, V\|_{Y_{n+\frac{1}{2}}}\right)+\|U, V\|_{X_{11}}+\|U, V\|_{E}.
$$
\begin{Theorem}(\cite{IM})
Fix $N_{1}=400$ and $N_{2}=200$. Fix the leading order Euler flow to be
$$
\left[u_{E}^{0}, v_{E}^{0}, P_{E}^{0}\right]:=[1,0,0].
$$
Assume the following pieces of initial data at $\{x=0\}$ are given for $i=0, \ldots, N_{1},$ and $j=1, \ldots N_{1},$
$$
\left.u_{p}^{i}\right|_{x=0}=: U_{p}^{i}(y),\left.\quad v_{E}^{j}\right|_{x=0}=: V_{E}^{j}(Y),
$$
where we make the following assumptions on the initial datum:

(1) For $i=0,$ the boundary layer datum $\bar{U}_{p}^{0}(y)$ is in a neighborhood of Blasius. More precisely, we will assume
$$
\left\|\left(\bar{U}_{p}^{0}(y)-\bar{u}_{*}(0, y)\right)\langle y\rangle^{m_{0}}\right\|_{C^{\ell} 0} \leq \delta_{*},
$$
where $0<\delta_{*}<<1$ is small relative to universal constants, where $m_{0}, \ell_{0},$ are large, explicitly computable numbers. Assume also the difference $\bar{U}_{p}^{0}(y)-\bar{u}_{*}(0, y)$ satisfies generic parabolic compatibility conditions at $y=0$.

(2) For $i=1, . ., N_{1},$ the boundary layer datum, $U_{p}^{i}(\cdot)$ is sufficiently smooth and decays rapidly:
$$
\left\|U_{p}^{i}\langle y\rangle^{m_{i}}\right\|_{C^{\ell} i} \lesssim 1,
$$
where $m_{i}, \ell_{i}$ are large, explicitly computable constants (for instance, we can take $m_{0}=10,000,$ $\ell_{0}=10,000$ and $\left.m_{i+1}=m_{i}-5, \ell_{i+1}=\ell_{i}-5\right),$ and satisfies generic parabolic compatibility conditions at $y=0$.

(3) The Euler datum $V_{E}^{i}(Y)$ satisfies generic elliptic compatibility conditions,

(4) Assume Dirichlet datum for the remainders, that is
$$
\left.[u, v]\right|_{x=0}=\left.[u, v]\right|_{x=\infty}=0.
$$
Then there exists an $\varepsilon_{0}<<1$ small relative to universal constants such that for every $0 \leq \varepsilon \leq \varepsilon_{0},$ there exists a unique solution $\left(u^{\varepsilon}, v^{\varepsilon}\right)$ to system which satisfies the expansion in
the quadrant, Q. Each of the intermediate quantities in the expansion satisfies the following estimates for $i=1, \ldots, N_{1}$
$$
\begin{array}{l}
{\left\|\partial_{x}^{k} \partial_{y}^{j} u_{p}^{i}\langle z\rangle^{M}\right\|_{L_{y}^{\infty}} \leq C_{M, k, j}\langle x\rangle^{-\frac{1}{4}-k-\frac{j}{2}+\sigma_{*}}, \quad\left\|v_{p}^{i}\right\|_{L_{y}^{\infty}} \leq C_{M, k, j}\langle x\rangle^{-\frac{3}{4}-k-\frac{j}{2}+\sigma_{*}}} ,\\
\left\|\partial_{x}^{k} \partial_{Y}^{j} u_{E}^{i}\right\|_{L_{y}^{\infty}} \leq C_{k, j}\langle x\rangle^{-\frac{1}{2}-k-j} , \quad\left\|v_{E}^{i}\right\|_{L_{y}^{\infty}} \leq C_{k, j}\langle x\rangle^{-\frac{1}{2}-k-j},
\end{array}
$$
 where $\sigma_{*}:=\frac{1}{10,000}$. Finally, the remainder $(u, v)$ exists globally in the quadrant, $\mathcal{Q},$ and satisfies the following estimates
$$
\|u, v\|_{\mathcal{X}} \lesssim 1.
$$
\end{Theorem}

Li and Ding  (\cite{LD}) in  2020 concerned with the validity of the Prandtl boundary layer theory in the inviscid limit of the steady incompressible Navier-Stokes equations. Under the symmetry assumption, they established the validity of the Prandtl boundary layer expansions and the error estimates. The convergence rate as $\varepsilon \rightarrow 0$ is also given.
They considered the following steady incompressible Navier-Stokes equations
\begin{equation}
\left\{\begin{array}{l}
U U_{X}+V U_{Y}+P_{X}=\varepsilon U_{X X}+\varepsilon U_{Y Y}, \\
U V_{X}+V V_{Y}+P_{Y}=\varepsilon V_{X X}+\varepsilon V_{Y Y}, \\
U_{X}+V_{Y}=0,
\end{array}\right.\label{LD.1}
\end{equation}
in the domain
$$
\Omega:=\{(X, Y) \mid 0 \leq X \leq L, 0 \leq Y \leq 2\},
$$
with moving boundary conditions
$$
U(X, 0)=U(X, 2)=u_{b}>0, V(X, 0)=V(X, 2)=0.
$$

They further  focused on the problem when $\varepsilon \rightarrow 0 .$ As $\varepsilon \rightarrow 0$, a formal limit of the solution of $(\ref{LD.1})$ should be the shear flow $\left[U_{0}, V_{0}\right]=\left[u_{e}^{0}(Y), 0\right]$, which satisfies the corresponding Euler equations. We assume that this smooth positive function $u_{e}^{0}(\cdot)$ satisfies $u_{e}^{0}(1-Y)=u_{e}^{0}(1+Y)$, for any $Y \in[0,1]$ and $u_{e}^{0}(0)=u_{e}^{0}(2)=u_{e} \neq u_{b} .$ Accordingly, they assumed that the solution $[U, V]$ to $(\ref{LD.1})$  satisfies the following symmetrical conditions with respect to $Y=1$,
$$
U(X, 1-Y)=U(X, 1+Y), V(X, 1-Y)=-V(X, 1+Y), Y \in[0,1],
$$
It should be noted that, due to this assumption, the pair $[U, V]_{1 \leq Y \leq 2}$ satisfies equations $(\ref{LD.1})$  as long as $[U, V]_{0 \leq Y \leq 1}$ does. Then their discussion can be restricted to the domain
$$
\Omega_{0}:=\{(X, Y) \mid 0 \leq X \leq L, 0 \leq Y \leq 1\},
$$
and the boundary conditions turn to
$$
[U, V](X, 0)=\left[u_{b}, 0\right], \quad\left[U_{Y}, V\right](X, 1)=[0,0],
$$
  introduce the Prandtl's scaling
$$
x=X, y=\frac{Y}{\sqrt{\varepsilon}},
$$
and the new unknown functions
$$
U^{\varepsilon}(x, y)=U(X, Y), V^{\varepsilon}(x, y)=\frac{1}{\sqrt{\varepsilon}} V(X, Y).
$$
Under this transformation, system $(\ref{LD.1})$ can be rewritten as
\begin{equation}
\left\{\begin{array}{l}
U^{\varepsilon} U_{x}^{\varepsilon}+V^{\varepsilon} U_{y}^{\varepsilon}+P_{x}^{\varepsilon}=U_{y y}^{\varepsilon}+\varepsilon U_{x x}^{\varepsilon}, \\
U^{\varepsilon} V_{x}^{\varepsilon}+V^{\varepsilon} V_{y}^{\varepsilon}+P_{y}^{\varepsilon} / \varepsilon=V_{y y}^{\varepsilon}+\varepsilon V_{x x}^{\varepsilon}, \\
U_{x}^{\varepsilon}+V_{y}^{\varepsilon}=0
\end{array}\right.\label{LD.2},
\end{equation}
in the domain
$$
\Omega_{\varepsilon}:=\left\{(x, y) \mid 0 \leq x \leq L, 0 \leq y \leq \frac{1}{\sqrt{\varepsilon}}\right\},
$$
with the boundary conditions
$$
\left[U^{\varepsilon}, V^{\varepsilon}\right](x, 0)=\left[u_{b}, 0\right], \quad\left[U_{y}^{\varepsilon}, V^{\varepsilon}\right]\left(x, \frac{1}{\sqrt{\varepsilon}}\right)=[0,0] .
$$
In what follows, the authors in \cite{IM} intend to find the exact solutions $\left[U^{\varepsilon}, V^{\varepsilon}, P^{\varepsilon}\right]$ in form of
\begin{equation}
\left\{\begin{array}{l}
U^{\varepsilon}(x, y)=u_{a p p}(x, y)+\varepsilon^{\gamma+\frac{1}{2}} u^{\varepsilon}(x, y) \\
V^{\varepsilon}(x, y)=v_{a p p}(x, y)+\varepsilon^{\gamma+\frac{1}{2}} v^{\varepsilon}(x, y) \\
P^{\varepsilon}(x, y)=p_{a p p}(x, y)+\varepsilon^{\gamma+\frac{1}{2}} p^{\varepsilon}(x, y)
\end{array}\right.\label{LD.3}
\end{equation}
where
\begin{equation}
\left\{\begin{array}{l}
u_{a p p}(x, y)=u_{e}^{0}(\sqrt{\varepsilon} y)+u_{p}^{0}(x, y)+\sqrt{\varepsilon} u_{e}^{1}(x, \sqrt{\varepsilon} y)+\sqrt{\varepsilon} u_{p}^{1}(x, y), \\
v_{a p p}(x, y)=v_{p}^{0}(x, y)+v_{e}^{1}(x, \sqrt{\varepsilon} y)+\sqrt{\varepsilon} v_{p}^{1}(x, y), \\
p_{a p p}(x, y)=\sqrt{\varepsilon} p_{e}^{1}(x, \sqrt{\varepsilon} y)+\sqrt{\varepsilon} p_{p}^{1}(x, y)+\varepsilon p_{p}^{2}(x, y).
\end{array}\right.\label{LD.4}
\end{equation}

For convenience, denote $z:=\sqrt{\varepsilon} y .$ Boundary conditions on $\{y=0\}:$
\begin{equation}
u_{e}^{0}(0)+u_{p}^{0}(x, 0)=u_{b}, \quad u_{e}^{1}(x, 0)+u_{p}^{1}(x, 0)=0, \quad u^{\varepsilon}(x, 0)=0,\label{LD.5}
\end{equation}
\begin{equation}
v_{p}^{0}(x, 0)+v_{e}^{1}(x, 0)=0, \quad v_{p}^{1}(x, 0)=0, \quad v^{\varepsilon}(x, 0)=0.\label{LD.6}
\end{equation}
Boundary conditions on $\left\{y=\frac{1}{\sqrt{\varepsilon}}\right\}$ :
\begin{equation}
u_{p y}^{0}\left(x, \frac{1}{\sqrt{\varepsilon}}\right)=0, \quad u_{e z}^{1}(x, 1)=0, u_{p y}^{1}\left(x, \frac{1}{\sqrt{\varepsilon}}\right)=0, \quad u_{y}^{\varepsilon}\left(x, \frac{1}{\sqrt{\varepsilon}}\right)=0;\label{LD.7}
\end{equation}
\begin{equation}
v_{p}^{0}\left(x, \frac{1}{\sqrt{\varepsilon}}\right)=0, \quad v_{e}^{1}(x, 1)=0, v_{p}^{1}\left(x, \frac{1}{\sqrt{\varepsilon}}\right)=0, \quad v^{\varepsilon}\left(x, \frac{1}{\sqrt{\varepsilon}}\right)=0.\label{LD.8}
\end{equation}
Boundary conditions on $\{x=0\}$ :
\begin{equation}
u_{p}^{0}(0, y)=\bar{u}_{0}(y), \quad u_{e}^{1}(0, z)=u_{b}^{1}(z), \quad u_{p}^{1}(0, y)=\bar{u}_{1}(y), \quad u^{\varepsilon}(0, y)=0;\label{LD.9}
\end{equation}
\begin{equation}
v_{e}^{1}(0, z)=V_{b 0}(z), \quad v^{\varepsilon}(0, y)=0.\label{LD.10}
\end{equation}
Boundary conditions on $\{x=L\}$ :
\begin{equation}
v_{e}^{1}(L, z)=V_{b L}(z), \quad\left[p^{\varepsilon}-2 \varepsilon u_{x}^{\varepsilon}, u_{y}^{\varepsilon}+\varepsilon v_{x}^{\varepsilon}\right](L, y)=0.\label{LD.11}
\end{equation}

\begin{Theorem}(\cite{LD})
Let $u_{b}>0$ be a constant tangential velocity of the Navier-Stokes flow on the boundary $\{Y=0\}$, and let $u_{e}^{0}(Y)$ be a smooth positive Euler flow satisfies $u_{e z}^{0}(1)=0 .$ Suppose that the boundary conditions prescribed in $(\ref{LD.5})$-$(\ref{LD.11})$ hold. Suppose further that the positive condition $\min _{y}\left\{u_{e}^{0}(\sqrt{\varepsilon} y)+\bar{u}_{0}(y)\right\}>0$ holds. Then there exists a constant $L_{0}>0$, which depends only on the prescribed data, such that for $0<L \leq L_{0}$ and $\gamma \in\left(0, \frac{1}{5}\right)$, the asymptotic expansions stated in $(\ref{LD.3})$-$(\ref{LD.4})$ is a solution to equations $(\ref{LD.2})$ on $\Omega_{\varepsilon}$ together with the corresponding boundary conditions. The remainder solutions $\left[u^{\varepsilon}, v^{\varepsilon}\right]$ satisfies the estimate
$$
\left\|\nabla_{\varepsilon} u^{\varepsilon}\right\|_{L^{2}\left(\Omega_{\varepsilon}\right)}+\left\|\nabla_{\varepsilon} v^{\varepsilon}\right\|_{L^{2}\left(\Omega_{\varepsilon}\right)}+\left\|u^{\varepsilon}\right\|_{L^{\infty}\left(\Omega_{\varepsilon}\right)}+\sqrt{\varepsilon}\left\|v^{\varepsilon}\right\|_{L^{\infty}\left(\Omega_{\varepsilon}\right)} \leq C_{0}.
$$
\end{Theorem}

\subsection{  2D    Prandtl Equations-Local Existence  }

Usually von-Mises transformation and Coroco transformation are used to prove the existence of the global solutions. Therefore, in this section, we list the papers that used the shear flow method to study of the Prandtl equations in Gevrey space, and the energy estimation method to prove the existence of local solutions.

The authors in \cite{CLS1,CLS2,AWXY,W1,DDGM,cwz1,cwz,WZ2,WW} studied the local well-posedness of the Prandtl equations around a shear flow.

Cannone, Lombardo and Sammartino   (\cite{CLS1}) in  2013  concerned with Prandtl  equations with incompatible
data, i.e., with initial data that, in general, which do not fulfil the boundary conditions
imposed on the solution. Under the hypothesis of analyticity in the streamwise
variable, the authors proved that Prandtl equations, on the half-plane or on the
half-space, are well-posed for a short time.
Prandtl  equations on the half-plane are written as
\begin{equation}
\left\{\begin{array}{l}
\partial_{t} u^{P}+u^{P} \partial_{x} u^{P}+v^{P} \partial_{Y} u^{P}+\partial_{x} p^{E}=\partial_{Y Y} u^{P}, \quad(x, Y) \in \mathbb{R} \times \mathbb{R}^{+}, \\
\partial_{x} u^{P}+\partial_{Y} v^{P}=0, \\
u^{P}(x, Y=0, t)=0, \\
u^{P}(x, Y, t=0)=u_{0}(x, Y), \\
u^{P}(x, Y \rightarrow \infty, t)=U(x, t), \\
v^{P}(x, Y=0, t)=0
\end{array}\right.\label{+1}
\end{equation}
where the Euler matching datum $U(x, t)$ and the pressure $p^{E}=p^{E}(x, t)$ satisfy the Bernoulli's law:
$$
\partial_{t} U+U \partial_{x} U+\left.\partial_{x} p^{E}\right|_{y=0}=0.
$$

We first introduce the norms of the spaces $\mathcal{K}^{l, \rho}$, $\mathcal{K}^{l, \rho, \alpha}$ and $\mathcal{K}_{\beta,T}^{l, \rho, \alpha}$, respectively,
 \begin{equation}
|f|_{l, \rho}=\sum_{j=0}^l\sup\limits_{\lambda\in(-\rho,\rho)}\|\partial_x^j f(\cdot+i\lambda)\|_{L_x}<+\infty,\label{-1}
\end{equation}
 \begin{equation}
|f|_{l, \rho ,\alpha} =\sum\limits_{0\leq j\leq 2}\sum\limits_{i\leq l-j}\sup\limits_{Y\in \mathbb{R}^+}\langle \label{-2} Y\rangle|^\alpha\partial_{Y}^{j}\partial_{x}^{i}f(\cdot,Y)|_{0,\rho}<+\infty,
\end{equation}
\begin{equation}
|f|_{l, \rho ,\alpha,\beta,T} =\sum\limits_{0\leq j\leq 2}\sum\limits_{i\leq l-j}\sup\limits_{0\leq t\leq T}|\partial_{Y}^{j}\partial_{x}^{i}f(\cdot,\cdot,t)|_{0,\rho-\beta t,\alpha}
+\sum\limits_{i\leq l-2}\sup\limits_{0\leq t\leq T}|\partial_{t}\partial_{x}^{i} f(\cdot,\cdot,t)|_{0,\rho-\beta t,\alpha}<+\infty.\label{-3}
\end{equation}

\begin{Theorem}(\cite{CLS1}) Suppose that $U \in \mathcal{K}_{\beta_{0}, T_{0}}^{l, \rho_{0}}$ and $u_{0}-U \in K^{l, \rho_{0}, \alpha} .$ Then there exist $0<\rho<\rho_{0}$, $\beta>\beta_{0}>0$ and $0<T<T_{0}$ such that Prandtl equations \eqref{+1} admit a unique solution $u^{P} .$ This solution can be written as
$$
u^{P}(x, Y, t)=u^{S}+u^{R},
$$
where
(i) $u^{S}$ is an initial layer corrector and has the following form:
$$
u^{S}=-2 u_{0}(x, Y=0) \operatorname{erfc}\left(\frac{Y}{2 \sqrt{t}}\right),
$$
where $\operatorname{erfc}(y)=\frac{2}{\sqrt\pi}\int_0^ye^{-z^2}dz$.
(ii) $u^{R}$ is the solution of the Prandtl equation with compatible data and with a source term that keeps into account the interaction with $u^{S},$ and has following form:
$$
u^{R}(x, Y, t)=\tilde{u}(x, Y, t)+U(x, t),
$$
where $\tilde{u} \in \mathcal{K}_{\beta, T}^{l, \rho, \alpha}$;\\

(iii) The solution $u^{P}(x, Y, t),$ for $0<t<T$, is analytic both in $x$ and $Y$.
\end{Theorem}

 Cannone,  Lombardo and  Sammartino (\cite{CLS2}) in 2014 proved the well-posedness of the Prandtl boundary layer equations   (\ref{1.1.6})
on a periodic strip when the initial and the boundary data are not assigned to be compatible in the domain $(x,y)\in[0,2\pi]\times\mathbb{R}_T{}$,
\begin{equation}\left\{
\begin{array}{ll}
u_{t}+ u u_{x}+  v u_{y}= \nu u_{yy}+p_{x}, \\
 u_{x}+  v_{y}= 0,\\
u(0,x,y)=u_{0},\ \ u(t,x,0)=v(t,x,0)=0,\\
\lim\limits_{y\rightarrow+\infty}u(t,x,y)=U(t,x).
\end{array}
 \label{1.1.6}         \right.\end{equation}
But, initial and boundary data are not compatible, i.e., $$u_{0}(x,0)\neq 0=u^{P}(0,x, 0).$$
 The space $\mathcal{H}^{\sigma}$ is the space of the $2\pi$ periodic real functions $f(x)$ such that
$$|f|_{\sigma}:=\sum\limits_{k\in\mathbb{Z}}|\widehat{f}(k)|e^{|k|\sigma}< \infty, $$
where $\widehat{f}(k)$ are the Fourier coefficients of $f$. \\
The space  $\mathcal{H}^{\sigma,\alpha}_{\beta,T}$ is the space of the functions $f(t,x,y)$, $2\pi$ periodic with respect to $x$, such that
\begin{equation*}
f,\partial_{t}f,\partial_{y}^{j}f\in \mathcal{H}^{\sigma-\beta t,\alpha}, \forall 0\leq t\leq T, \text{and} ~j\leq 2.
 \end{equation*}
Moreover,
\begin{equation*}
|f|_{\sigma ,\alpha,\beta,T} =\sum\limits_{0\leq j\leq 2}\sup\limits_{0\leq t\leq T}|\partial_{y}^{j}f(\cdot,\cdot,t)|_{\sigma-\beta t,\alpha}
+\sup\limits_{0\leq t\leq T}|\partial_{t} f(\cdot,\cdot,t)|_{\sigma-\beta t,\alpha} <\infty.
 \end{equation*}

\begin{Theorem}(\cite{CLS2})\label{ }
Suppose that the Euler datum $U$ and the initial datum $u_{0}$ satisfy the following hypotheses: \\
(i) $U$ satisfies the Bernoulli's law equation (\ref{1.1.3});\\
(ii)$U(x,t)\in \mathcal{H}^{\sigma_{0}}_{\beta_{0},T_{0}}$;\\
(iii) $u_{0}(x,y)-U(0,x)\in \mathcal{H}^{\sigma_{0},\alpha}$.\\
Then there exist $0<\sigma<\sigma_{0}$, $\beta>\beta_{0}>0$ and $0<T<T_{0}$ such that Prandtl equations (\ref{1.1.6}) admit, in [0,T ], a unique
 solution $u^{P}$. This solution can be written as
$$u^{P}(t,x,y)=u^{s}+u^{R},$$
where\\
1. $u^{s}$ is an initial layer corrector and has the following form:
$$u^{s}=-2u_{0}(x,0)\operatorname{erfc}\left(\frac{y}{2\sqrt{t}}\right).$$
2. $u^{R}$ is the solution of the Prandtl equation with compatible data and with a source term that keeps into account the interaction with $u^{s}$;  $u^{R}$ can be decomposed as follows:
$$u^{R}(t,x,y)=\widetilde{u}(t,x,y)+U$$
where
 $$\widetilde{u}\in \mathcal{H}^{\sigma,\alpha}_{\beta,T}. $$
Moreover, the solution $u^{P}(t,x,y)$, for $t>0$, is analytic both in $x$ and $Y$.

\end{Theorem}

  Alexander,  Wang,  Xu and  Yang (\cite{AWXY}) in 2015   studied the local well-posedness of the Prandtl equations  (\ref{1.1.6}) around a monotonic
 shear flow in the domain $   \mathbb{R}_{+} ^{2}=\{(x,y)| x\in  \mathbb{R} , y\in[0, +\infty)  \}$.   In this paper, assume that
$$u_{0}(x,y)=u_{0}^{s}(y)+\widehat{u}_{0}(x,y),$$
where $u_{0}^{s}$ is monotonic in $y$,
\begin{equation}
\partial_{y} u_{0}^{s}(y)>0,\ \ \forall y\geq 0.
\label{1.14}
\end{equation}
 Set $\Omega_{T}=[0,T]\times  \mathbb{R}_{+} ^{2}$. For any non-negative integer $k$ and real number $\ell$, define
$$\|u\|_{\mathcal{A}^{k}_{\ell}(\Omega_{T} ) }=\left(\sum\limits_{k_{1}+[\frac{k_{2}+1}{2}]\leq k}
\|\langle y \rangle^{\ell} \partial_{t,x}^{k_{1}} \partial_{y}^{k_{2}} u\|^{2}_{L^{2}([0,T]\times  \mathbb{R}_{+} ^{2} )} \right)^{\frac{1}{2}},$$
and
$$\|u\|_{\mathcal{D}^{k}_{\ell}(\Omega_{T} ) }= \sum\limits_{k_{1}+[\frac{k_{2}+1}{2}]\leq k}
\|\langle y \rangle^{\ell} \partial_{t,x}^{k_{1}} \partial_{y}^{k_{2}} u\| _{L^{\infty}(L^{2}_{t,x} )} ,$$
with  $\langle y \rangle=(1+y^{2}) ^{\frac{1}{2}}$.
 \begin{Theorem}(\cite{AWXY})
 Concerning the problem  (\ref{1.1.6}) with $U$ being a constant, we have the following existence and stability results. \\
(1) Given any integer $k\geq 5$ and real number $\ell> \frac{1}{2}$, let the initial data $u_{0}(x,y)=u_{0}^{s}(y)+\tilde{u}_{0}(x,y)$ satisfy the compatibility
 conditions of the initial boundary value problem (\ref{1.1.6})  up to order $k+4$. Assume the following two conditions: \\
(i) The monotonicity condition (\ref{1.14}) holds for $ u_{0}^{s}$, and
\begin{equation}\left\{
\begin{array}{ll}
 \partial_{y}^{2j}u_{0}^{s}(0)= 0, \ \ \ \   \forall 0\leq j\leq k+4, \\
\|\langle y\rangle^{\ell} (u_{0}^{s}-1)\|_{H^{2k+9}(\mathbb{R}_{+} ^{2})}+\|\frac{\partial_{y}^{2 }u_{0}^{s}}{\partial_{y} u_{0}^{s}}  \|_{H^{2k+7}(\mathbb{R}_{+} ^{2})}\leq C,
\end{array}
 \label{w.20}         \right.\end{equation}
for a fixed constant $C>0$.\\
(ii) There exists a small constant $\epsilon>0$ depending only on $ u_{0}^{s}$, such that
\begin{equation}
\|\tilde{u}_{0}\|_{\mathcal{A}_{\ell}^{2k+9}(\mathbb{R}_{+} ^{2})}  +\|\frac{\partial_{y}\tilde{u}_{0} }{\partial_{y} u_{0}^{s}}  \|_{\mathcal{A}_{\ell}^{2k+9}(\mathbb{R}_{+} ^{2})}\leq \epsilon.
\label{w.21}
\end{equation}
Then there is a $T>0$, such that the problem (\ref{1.1.6}) admits a classical solution $u$ satisfying
\begin{equation}
u- u ^{s}\in \mathcal{A}_{\ell}^{k}(\Omega_{T}),     \frac{\partial_{y}(u-u ^{s} ) }{\partial_{y} u_{0}^{s}}  \in \mathcal{A}_{\ell}^{k}(\Omega_{T}),
v\in  \mathcal{D}_{0}^{k-1}(\Omega_{T})  .
\label{w.22}
\end{equation}
Here, note that $( u ^{s}(t,y) , 0)$ is the shear flow of the Prandtl equation defined by the initial data $( u _{0}^{s}(y) , 0)$. \\
(2) The solution is unique in the function space described in (\ref{w.22}). Moreover, we have the stability with respect to the initial data in the following sense. For any given two initial data
$$u_{0}^{1}=  u ^{s}_{0}+\tilde{u}_{0}^{1},\ \ u_{0}^{2}=  u ^{s}_{0}+\tilde{u}_{0}^{2}, $$
if $u ^{s}_{0}$ satisfy  (\ref{1.14})-(\ref{w.20}), and $\tilde{u}_{0}^{1},\tilde{u}_{0}^{2} $ satisfy (\ref{w.21}), then the corresponding
solutions $(u^{1}, v^{1})$ and  $(u^{2}, v^{2})$ of (\ref{1.1.6})  satisfy
$$\|u^{1}-u^{2}\|_{ \mathcal{A}_{\ell}^{p}(\Omega_{T}) }+\|v^{1}-v^{2}\|_{\mathcal{D}_{0}^{p-1}(\Omega_{T}) }\leq C
\|\frac{\partial}{\partial y}\left(\frac{u_{0}^{1}-u_{0}^{2}}{\partial y u _{0}^{s}}  \right)\|_{\mathcal{A}_{\ell}^{p}(\Omega_{T}) }$$
for all $p\leq k-1$, where the constant $C>0$ depends only on $T$ and the upper bounds of the norms of $u_{0}^{1}, u_{0}^{2}$.

\end{Theorem}

  Wu (\cite{W1}) in 2015 studied the local well-posedness of classical solutions to the unsteady Prandtl equations with Robin boundary condition in half space $\mathbb{R}_{+}$ in weighted Sobolev spaces with shear flow. His results are also valid for the Dirichlet boundary case.
He studied the following Prandtl equations with Robin boundary condition
\begin{equation}\left\{
\begin{array}{ll}
u_{t}+uu_{x}+vu_{y} =u_{yy}, \\
u_{x}+v_{y}=0,\\
(u_{y}-\beta u)|_{y=0}=0,\ \ v_{y=0}=0,\\
u(0,x,y)=u_{0}, \\
\lim\limits_{y\rightarrow+\infty}u(t,x,y)=U=1
\end{array}
 \label{1.1.8}         \right.\end{equation}
and defined the following functional spaces:
$$\|u\|_{\mathcal{A}^{k}_{\ell}(\Omega_{T} ) }=\left(\sum\limits_{k_{1}+[\frac{k_{2}+1}{2}]\leq k}
\|\langle y \rangle^{\ell} \partial_{t,x}^{k_{1}} \partial_{y}^{k_{2}} u\|^{2}_{L^{2}([0,T]\times  \mathbb{R}_{+} ^{2} )} \right)^{\frac{1}{2}},$$
$$\|u\|_{\mathcal{C}^{k}_{\ell}(\Omega_{T} ) }= \sum\limits_{k_{1}+[\frac{k_{2}+1}{2}]\leq k}
\|\langle y \rangle^{\ell} \partial_{t,x}^{k_{1}} \partial_{y}^{k_{2}} u\| _{L^{2}_{y}(L^{\infty}_{t,x} )} ,$$
$$\|u\|_{\mathcal{D}^{k}_{\ell}(\Omega_{T} ) }= \sum\limits_{k_{1}+[\frac{k_{2}+1}{2}]\leq k}
\|\langle y \rangle^{\ell} \partial_{t,x}^{k_{1}} \partial_{y}^{k_{2}} u\| _{L^{\infty}_{y}(L^{2}_{t,x} )} ,$$
with  $\langle y \rangle=(1+y^{2}) ^{\frac{1}{2}}$.

\begin{Theorem}(\cite{W1})
 Concerning the nonlinear unsteady Prandtl equations with Robin boundary condition (\ref{1.1.8}), giving any integer $k\geq 5$ and real number $\ell>\frac{1}{2}$, we have the following existence, uniqueness and stability results.

(1) For any $\delta_{\beta}>0$, $\delta_{\beta}\leq \beta<+\infty$, assume the initial data $u_{0}(x,y)=u_{0}^{s}(y)+\tilde{u}_{0}(x,y)$ satisfies the following two conditions: \\
(i) $u^{s}_{0}$ satisfies
\begin{equation}\left\{
\begin{array}{ll}
u_{0}^{s}(y)>0, \ \  \partial_{y}u_{0}^{s}(y)>0, \ \ \beta-\frac{\partial_{yy}u_{0}^{s}(y)}{\partial_{y}u_{0}^{s}(y)}\geq\delta_{s,0}>0 , \forall y\in [0,+\infty), \\
\partial_{y}^{2j}(\partial_{y}u_{0}^{s}(y)-\beta  u_{0}^{s}(y))|_{y=0}=0 , \ \ \forall 0\leq j\leq 2k+10, \\
\lim\limits_{y\rightarrow+\infty}u^{s}_{0}( y) =1, \\
\|u^{s}_{0}-1 \|_{L^{2}}+\|u^{s}_{0} \|_{\infty}+\|u^{s}_{0} \|_{\mathcal{C}_{\ell}^{2k+11} }+\|\frac{\partial_{yy}u_{0}^{s}(y)}{\partial_{y}u_{0}^{s}(y)} \|_{\mathcal{C}_{\ell}^{2k+10} } \leq C,
\end{array}
 \label{w.24}         \right.\end{equation}
for a fixed constant $C>0$.\\
 (ii) There exists a small constant $\epsilon >0$ such that $\tilde{u}_{0}=u_{0}-u_{0}^{s}$ satisfies
\begin{equation}\left\{
\begin{array}{ll}
\partial_{y}^{2j}(\partial_{y}\tilde{u}_{0} (x,y)-\beta  \tilde{u}_{0} (x,y))|_{y=0}=0 , \ \ \forall 0\leq j\leq 2k+8, \\
\lim\limits_{y\rightarrow+\infty}\tilde{u} _{0}( x,y) =0, \\
\|\tilde{u} _{0}  \| _{\mathcal{A}_{\ell}^{2k+9}(\mathbb{R}_{+}^{2},t=0) }+\|\frac{\partial_{y }\tilde{u}_{0} }{\partial_{y}u_{0}^{s} } \|_{\mathcal{A}_{\ell}^{2k+9}(\mathbb{R}_{+}^{2},t=0) } \leq \epsilon.
\end{array}
 \label{w.25}         \right.\end{equation}
Then there exists a time $T\in [0,+\infty)$, such that the Prandtl system (\ref{1.1.8}) admits a unique classical solution $(u,v)$ satisfying
\begin{equation}
u>0,\ \ \partial_{y} u>0, \ \  \beta-\frac{\partial_{yy}u }{\partial_{y}u }\geq\delta >0 ,
\label{w.26}
\end{equation}
and
\begin{equation}\left\{
\begin{array}{ll}
 u-u^{s}\in \mathcal{A}_{\ell}^{ k }([0,T]\times\mathbb{R}_{+}^{2} ) ,\ \
\partial_{y}(u-u^{s}),\frac{\partial_{y }(u-u^{s}) }{\partial_{y}u ^{s}} \in \mathcal{A}_{\ell}^{ k }([0,T]\times\mathbb{R}_{+}^{2} ), \\
 v\in \mathcal{D}_{0}^{ k-1 }([0,T]\times\mathbb{R}_{+}^{2} ) ,\ \
\partial_{y} v, \partial_{yy}v\in \mathcal{A}_{\ell}^{ k-1 }([0,T]\times\mathbb{R}_{+}^{2} ) , \\
 \partial_{y}^{j}u|_{y=0}-\partial_{y}^{j}u^{s}|_{y=0}\in \mathcal{A}_{\ell}^{ k -[\frac{j+1}{2}] }([0,T]\times\mathbb{R}  ) ,\ \ 0\leq j\leq 2k, \\
 \partial_{y}^{j+1}v|_{y=0}\in \mathcal{A}_{\ell}^{ k-1 -[\frac{j+1}{2}] }([0,T]\times\mathbb{R}  ) ,\ \ 0\leq j\leq 2k-2.
\end{array}
 \label{w.27}         \right.\end{equation}

(2)  For any $\delta_{\beta}>0$, $\delta_{\beta}\leq \beta<+\infty$, the classical solution to the Prandtl system (\ref{1.1.8}) is stable with respect to the initial data
in the following sense: for any given two initial data
$$ u_{0}^{1}=u_{0}^{s}+  \tilde{u}_{0}^{1}, \ \  u_{0}^{2}=u_{0}^{s}+  \tilde{u}_{0}^{2},  $$
if $u_{0}^{s}$ satisfies  (\ref{w.24})  and $ u_{0}^{1}, u_{0}^{2}$ satisfy  (\ref{w.25}) , then for all $p\leq k-2$, the corresponding
solutions $( u^{1},v^{1})$ and $( u^{2},v^{2})$ of the Prandtl system (\ref{1.1.8}) satisfy
\begin{eqnarray*}
&& \|u^{1}-u^{2}\|_{ \mathcal{A}_{\ell}^{ p}([0,T]\times\mathbb{R}_{+}^{2} )}
+ \| \frac{\partial_{y}(u^{1}-u^{2})}{\partial_{y} u^{s}}\| _{ \mathcal{A}_{\ell}^{ p}([0,T]\times\mathbb{R}_{+}^{2} )}\nonumber \\
&& \quad +\|v^{1}-v^{2}\|_{ \mathcal{D}_{\ell}^{ p-1}([0,T]\times\mathbb{R}_{+}^{2} )}
+ \sum\limits_{j=0}^{2p}\|  \partial_{y}^{j}u^{1}|_{y=0}-\partial_{y}^{j}u^{2}|_{y=0}\| _{ \mathcal{A}_{\ell}^{ p- [\frac{j+1}{2}] }([0,T]\times\mathbb{R} )}\nonumber \\
&& \leq C(T,\epsilon, u_{0}^{s} ) \|\partial_{y}(\frac{u_{0}^{1}-u_{0}^{2}}{\partial_{y}(u_{0}^{1}+u_{0}^{2} )}) \| _{ \mathcal{A}_{\ell}^{ p}( \mathbb{R}_{+}^{2},t=0 )}
+\frac{C(T,\epsilon, u_{0}^{s} )}{\max\{\sqrt{\beta-C_{\eta}} \} , \sqrt{\delta}} \|\partial_{y}(\frac{u_{0}^{1}-u_{0}^{2}}{\partial_{y}(u_{0}^{1}+u_{0}^{2} )})|_{y=0} \| _{ \mathcal{A}_{\ell}^{ p}( \mathbb{R} )} .
\label{ }
\end{eqnarray*}

(3) As $\beta\rightarrow+\infty$, $\partial_{t}^{j}u^{s}|_{y=0}=\mathcal{O}(\frac{1}{\beta})$, $0\leq j\leq k$,
$\|u|_{y=0}-u^{s}|_{y=0}\|_{\mathcal{A}_{\ell}^{ k}} =\mathcal{O}(\frac{1}{\beta})$ and $(u,v)$ satisfies (\ref{w.27}) uniformly.
 When $\beta=+\infty$, $(u,v)$ satisfies (\ref{w.27}) and for $p\leq k-2$,
\begin{eqnarray*}
&& \|u^{1}-u^{2}\|_{ \mathcal{A}_{\ell}^{ p}([0,T]\times\mathbb{R}_{+}^{2} )}
+ \| \frac{\partial_{y}(u^{1}-u^{2})}{\partial_{y} u^{s}}\| _{ \mathcal{A}_{\ell}^{ p}([0,T]\times\mathbb{R}_{+}^{2} )}
+\|v^{1}-v^{2}\|_{ \mathcal{D}_{\ell}^{ p-1}([0,T]\times\mathbb{R}_{+}^{2} )} \nonumber \\
&&\quad + \sum\limits_{j=0}^{2p}\|  \partial_{y}^{j}u^{1}|_{y=0}-\partial_{y}^{j}u^{2}|_{y=0}\| _{ \mathcal{A}_{\ell}^{ p- [\frac{j+1}{2}] }([0,T]\times\mathbb{R} )}
\leq C  \|\partial_{y}(\frac{u_{0}^{1}-u_{0}^{2}}{\partial_{y}(u_{0}^{1}+u_{0}^{2} )}) \| _{ \mathcal{A}_{\ell}^{ p}( \mathbb{R}_{+}^{2},t=0 )}
  .
\label{ }
\end{eqnarray*}

\end{Theorem}

Authors of \cite{DDGM} in 2017 studied the linear stability of shear flow solutions to the interactive boundary layer with high frequency perturbations. The general concern of \cite{DDGM} is the boundary layer behavior of high Reynolds number flows. The authors of \cite{DDGM} restricted to the two-dimensional case near a flat boundary, and
considered the Navier-Stokes system for $t>0$, $\textbf{x}=(x,Y)\in \mathbb{R}\times \mathbb{R}_+$,
\begin{eqnarray}
 \partial_tu^{\nu}-\nu\Delta u^{\nu}+u^{\nu}\cdot \nabla u^{\nu}+\nabla P^{\nu}=0,\ \ divu^{\nu}=0,\ u^{\nu}(x,Y=0)=0.   \label{d6}
\end{eqnarray}
The parameter $\nu<<1$ refers to the inverse Reynolds number.
They wrote the linearized system as the following system, with unknowns $U$, $V$ and $u_e$,
\begin{eqnarray}\left\{\begin{array}{ll}
 \partial_tU+U_s\partial_xU+VU'_s-\partial_y^2U=\partial_tu_e+\partial_xu_e,\\
  \partial_xU+\partial_yV=0,
\end{array}\label{d7}\right.\end{eqnarray}
$$U|_{y=0}=V|_{y=0}=0,\ \lim\limits_{y\rightarrow+\infty}=u_e,\ u_e=\sqrt\nu DN\int_0^{+\infty}(u_e-U)dy,$$
where we  use the convention $DN f=-\partial_Y\psi|_{y=0}$, $\psi$ is the solution of $\Delta \psi=0$.

After linearization around the shear flow, we find
\begin{eqnarray}
 \partial_tU+U_s\partial_xU+VU'_s-\partial_y^2U=\partial_tu_e+\partial_xu_e,\ \partial_xU+\partial_yV=0,\label{d8}
\end{eqnarray}
$$U|_{y=0}=V|_{y=0}=0,\ \lim\limits_{y\rightarrow+\infty}=u_e,\ \Delta_su_e=\int_0^{+\infty}(u_e-U)dy,$$
with $\Delta_su_e$  the displacement thickness associated to the shear flow.

Criterion 1. Assume that $\xi(\Delta_s)=0$. Assume, moreover, that  either $U_s^{''}(0)>0$ and $\xi_+(\Delta_s)\geq\xi_-(\Delta_s)$ or $U_s^{''}(0)<0$ and $\xi_+(\Delta_s)>\xi_-(\Delta_s)$.

Criterion 2. Assume that $\xi(0)=0$. Assume, moreover, that  either $U_s^{''}(0)>0$ and $\xi_+(0)\geq\xi_-(0)$ or $U_s^{''}(0)<0$ and $\xi_+(0)>\xi_-(0)$.

Criterion 3. Assume that there exists an arbitrary abscissa $\gamma>0$ such that $\xi(\gamma)=0$, and moreover, either $U_s^{''}(0)>0$ and $\xi_+(\gamma)\geq\xi_-(\gamma)$ or $U_s^{''}(0)<0$ and $\xi_+(\gamma)>\xi_-(\gamma)$.

\begin{Theorem}(\cite{DDGM})
 Let $U_s$ be a shear flow satisfying Criterion 1. There exist $K,\eta > 0$ such that, for $k\geq K$, there exists
$\lambda_k\in \mathbb{C}$ with $Re\lambda_k\geq \eta k$ such that the system \eqref{d8} displays $(\lambda_k,k)$ instabilities.
\end{Theorem}
\begin{Theorem}(\cite{DDGM})\\
$\cdot$ Let $U_s$ be a shear flow satisfying Criterion 1. There exist $K,\eta,\gamma_+ > 0$ such that, for $k\geq K$ and $\nu>0$ satisfying $(\sqrt\nu k)^{-1}\leq\gamma_+$, there exists
$\lambda_{k,\nu}\in \mathbb{C}$ with $Re\lambda_{k,\nu}\geq \eta k$ such that the  system \eqref{d7} displays $(\lambda_k,k)$ instabilities.\\
$\cdot$ Let $U_s$ be a shear flow satisfying Criterion 1. There exist $K,\eta,\gamma_-,\gamma_+ > 0$ such that, for $k\geq K$ and $\nu>0$ satisfying $\gamma_-\leq(\sqrt\nu k)^{-1}\leq\gamma_+$, there exists
$\lambda_{k,\nu}\in \mathbb{C}$ with $Re\lambda_{k,\nu}\geq \eta k$ such that the  system \eqref{d7} displays $(\lambda_k,k)$ instabilities.\\
$\cdot$ If $U_s$ satisfies $U_s^{''}(0)>0$, there exist $\eta,\gamma_-,> 0$ such that, for $\nu k^3\geq 1$, there exists
$\lambda_{k,\nu}$  such that the  system \eqref{d7} displays $(\lambda_{k,\nu},k)$ instabilities with $Re\lambda_{k,\nu}\geq \eta \nu k^3$ .
\end{Theorem}
\begin{Theorem}(\cite{DDGM})
Let $U_s$ be an arbitrary monotone shear flow. There exist positive constants $\bar{\nu},K,S,\eta$ such that for $\nu\leq\bar{\nu},|\bar{k}|\geq K$
and $k\nu^{3/4}\geq S$, there exists $\lambda_{k,\nu}$ such that the  system \eqref{d7} displays $(\lambda_{k,\nu},k)$ instabilities with $Re\lambda_{k,\nu}\geq \eta \nu^{3/4} k^2$.
\end{Theorem}

\begin{Theorem}(\cite{DDGM})
Both for the \eqref{d7} and \eqref{d8} systems, we can build shear flows displaying instabilities with chosen behavior within
a given spectral domain. More precisely,\\
$\cdot$System \eqref{d8}: let $\mu\in \Gamma_1=\Big\{a+ib\in\mathbb{C},
 0<a<1,0<b<\sqrt{\frac{3a(1-a)^2}{(4-3a)}}\Big\}$. There exist a shear flow $U_s$  and $K>0$ such that, for $k\geq K$, there exists a sequence $\lambda_k$ with $\lambda_k\sim-ik\mu$ such
that the system \eqref{d8} displays $(\lambda_k,k)$ instabilities.\\
$\cdot$System \eqref{d7}: let $\mu\in \Gamma_2=\Big\{a+ib\in\mathbb{C},
 0<a<\frac{2}{3},0<b<\sqrt{\frac{a(2-3a)}{(3)}}\Big\}$. For any $\varepsilon>0$, there exist a shear flow $U_s$  and $K>0$ and $\gamma_+>0$ such that $k\geq K$ and $(\sqrt\nu k)^{-1}\leq\gamma_+$ there exists a sequence $\lambda_{\nu,k}$ with $|\lambda_{\nu,k}/(k+i\mu)|\leq\varepsilon$ such
that the system \eqref{d7} displays $(\lambda_{\nu,k},k)$ instabilities.\\
$\cdot$ System\eqref{d7}: let $\mu\in \Gamma_3=\Big\{a+ib\in\mathbb{C},
 0<a<1,0<b<\sqrt{a(1-a)}\Big\}$. For any $\gamma> 0$ such that, for any $\varepsilon>0$, there exist  a shear flow $U_s$ and  $K>0$ and $\gamma_+>\gamma>\gamma_-$ such that $k\geq K$ and $\gamma_-<(\sqrt\nu k)^{-1}\leq\gamma_+$, there exists a sequence $\lambda_{\nu,k}$ with $|\lambda_{\nu,k}/(k+i\mu)|\leq\varepsilon$ such
that the system \eqref{d7} displays $(\lambda_{\nu,k},k)$ instabilities.\\
\end{Theorem}

Using the paralinearization technique, Chen, Wang and Zhang (\cite{cwz1}) in 2018 proved  the local well-posedness of the Prandtl equation  for monotonic data in anisotropic Sobolev space with an exponential weighted and low regularity and provided a new  way for the zero-viscosity limit problem of the Navier-Stokes equations with the non-slip boundary condition. They obtained the local  well-posedness of Prandtl equations  when weighted function $\mu$ is an exponential function in $H^{3,1}_\mu(\mathbb{R}_+^2)\cap H^{1,2}_\mu(\mathbb{R}_+^2)$. They studied the Prandtl equation in $\mathbb{R}_+\times \mathbb{R}_+^2$,
\begin{eqnarray}\left\{\begin{array}{ll}
u_t+u\partial_xu+v\partial_yu-\partial_{y}^2u+\partial_xp=0,\\
\partial_xu+\partial_yv=0,\\
(u,v)|_{y=0}=0, \lim\limits_{y\rightarrow+\infty}u=U,\\
u|_{t=0}=u_0.
\end{array}\right.\label{9.1}\end{eqnarray}
For the simplicity, they considered  in \cite{cwz1} the case of a uniform outflow $U(t, x) =1$, which implies that $p$ is a constant. Then the Prandtl equation \eqref{9.1}  reduces to
\begin{eqnarray}\left\{\begin{array}{ll}
u_t+u\partial_xu+v\partial_yu-\partial_{y}^2u=0,\\
\partial_xu+\partial_yv=0,\\
(u,v)|_{y=0}=0, \lim\limits_{y\rightarrow+\infty}u=1,\\
u|_{t=0}=u_0.
\end{array}\right.\label{9.2}\end{eqnarray}
They considered the data around a shear flow, i.e.,
\begin{eqnarray*}\left\{\begin{array}{ll}
\partial_t u_s-\partial_{y}^2u_s=0,\\
u_s|_{y=0}=0, \lim_{y\rightarrow+\infty}u=1,\\
u_s|_{t=0}=u_{0s}(y).
\end{array}\right. \end{eqnarray*}
Obviously, $(u^s (t, y), 0)$ is a shear flow solution of the Prandtl equation \eqref{9.1}.
Let $u(t, x, y) = u^s (t, y) +\tilde{u}(t, x, y)$. Then $ \tilde{u}$ satisfies the following equations
\begin{eqnarray*}\left\{\begin{array}{ll}
\partial_t \tilde{u}+u\partial_x\tilde{u}+v\partial_yu-\partial_{y}^2\tilde{u}=0,\\
\partial_x\tilde{u}+\partial_yv=0,\\
\tilde{u}|_{y=0}=0, \lim_{y\rightarrow+\infty}u=0,\\
\tilde{u}|_{t=0}=\tilde{u}_{0}.
\end{array}\right. \end{eqnarray*}
\begin{Theorem}( \cite{cwz1})
Let $\mu(y)=e^{\frac{y}{2}}$. Assume that $u_0^s-1\in H_{y,\mu}^2,\ \ \tilde{u_0}\in H_\mu^{3,1}\cap H^{1,2}_\mu$ with $u_0^s(0)=\tilde{u_0}|_{y=0}=0$. Moreover, for some $c>0$ and $k=0,1,2$
\begin{eqnarray*}
 \partial_yu_0^s(y)\geq ce^{-y}, |\partial_y^k(u_0^s(y)-1)|\leq c^{-1}e^{-y}, \ \text{for}\ y\in[0,+\infty), \\
|\partial_y\tilde {u}_0^s(x,y)|+|\partial_x\partial_y\tilde {u}_0^s(x,y)|\leq \epsilon e^{-y}\ \text{for}\ (x,y)\in \mathbb{R}_+^2, \\
\|\tilde{u}_0\|_{{H_\mu}^{3,1}}+\|\tilde{u}_0\|_{{H_\mu}^{1,2}}\leq\epsilon.
\end{eqnarray*}
There exist $\epsilon_0> 0$ and $T > 0$ such that for any $\epsilon \in (0, \epsilon_0 )$, the Prandtl equation \eqref{9.1} has a unique solution in $[0, T]$,
 which satisfies
\begin{eqnarray*}
 c_1e^{-y}\leq \partial_yu(t,x,y)\leq C_1e^{-y} \ \text{for}\ (t,x,y)\in [0,T]\times \mathbb{R}_+^2. \\
|\partial_y^2 u(t,x,y)|+|\partial_x\partial_y u(t,x,y)|\leq C_1e^{-y}, \\
u-u^s\in L^\infty(0,T;H^{3,1}_\mu\cap H^{1,2}_\mu)\cap L^2(0,T;H^{1,3}_\mu).
\end{eqnarray*}
for some $c_1 > 0$ and $C_1 > 0$. Here $H^{s,\sigma}_
\mu$ is the weighted anisotropic Sobolev space.
\end{Theorem}

 When shear flow is non-monotonic, Chen, Wang and Zhang (\cite{cwz}) in 2018  proved the well-posedness of the linearized Prandtl equation  in Gevrey class $2-\theta$ for any $\theta> 0$ and constructed a class of solution with the growth like $e^{\sqrt kt}$ for the linearized Prandtl equation by energy method.
The authors in \cite{cwz} studied the Prandtl equation in $\mathbb{R}_+\times \mathbb{R}_+^2$,
\begin{eqnarray}\left\{\begin{array}{ll}
u_t+u\partial_xu+v\partial_yu-\partial_{y}^2u+\partial_xp=0,\\
\partial_xu+\partial_yv=0,\\
(u,v)|_{y=0}=0, \lim\limits_{y\rightarrow+\infty}u=U,\\
u|_{t=0}=u_0.
\end{array}\right.\label{8.1}\end{eqnarray}
where $(u, v)$ denotes the tangential and normal velocity of the boundary layer flow, and $(U(t, x), p(t, x))$ are the values on the boundary of the tangential velocity and pressure of the outflow, which satisfies the Bernoulli's law
 $$\partial_tU+U\partial_xU+\partial_xp=0.$$
The linearized Prandtl equation around $(u^s , 0)$ takes as follows
\begin{eqnarray}\left\{\begin{array}{ll}
\partial_t u+u^s\partial_xu+v\partial_yu^s-\partial_{y}^2u_s=0,\\
u|_{y=0}=0, \lim\limits_{y\rightarrow+\infty}u=0,\\
u|_{t=0}=u_{0}
\end{array}\right.\label{8.2}\end{eqnarray}
The main result of \cite{cwz} is stated as follows.
\begin{Theorem}( \cite{cwz})
Let $\theta\in(0,\frac{1}{2}]$. Assume that $e^{{\langle D_x\rangle}^{\frac{1}{2}+2\theta}}u_0\in H_\mu^{\frac{1}{2},1}$ with $\partial_y^ku_0|_{y=0}=0$ for $k=0,2$. Then there exists $T > 0$,
so that \eqref{8.1} has a unique solution $u$ in $[0, T]$.  In particular, we have
$$u_\Phi\in L^\infty(0,T;H_\mu^{\frac{1}{4}+\theta,1}).$$
Here we denote
$$f_\Phi=\mathcal{F}^{-1}(e^{\Phi(t,\xi)}\hat{f}(\xi)),\Phi(t,\xi)=(1-\lambda t)\langle\xi\rangle^{\frac{1}{2}+2\theta},$$
and $H_\mu^{s,\sigma}$ is the weighted Sobolev space with $\mu=e^{\frac{y}{2}}$.
\end{Theorem}

Yang and Zhu  (\cite{WZ2}) in 2020  considered the 2D Prandtl equation with constant outer flow and monotonic data, and proved that if the curvature of the velocity distribution (i.e., $\left.\partial_{y}^{2} u\right)$ is bounded near the boundary, then the solution can not develop the singularity.
They studied the Prandtl equation in $\mathbb{R}_{+} \times \mathbb{R}_{+}^{2}$ :
\begin{equation}
\left\{\begin{array}{l}
\partial_{t} u+u \partial_{x} u+v \partial_{y} u-\partial_{y}^{2} u+\partial_{x} p=0, \\
\partial_{x} u+\partial_{y} v=0, \\
\left.u\right|_{y=0}=\left.v\right|_{y=0}=0 ,\quad   \quad \lim\limits _{y \rightarrow+\infty} u(t, x, y)=U(t, x), \\
\left.u\right|_{t=0}=u_{0},
\end{array}\right.
\end{equation}
where $(u, v)$ denotes the tangential and normal velocity of the boundary layer flow, and $(U(t, x), p(t, x))$ is the values on the boundary of the tangential velocity and pressure of the outer flow, which satisfies the Bernoulli's law
$$
\partial_{t} U+U \partial_{x} U+\partial_{x} p=0.
$$

They considered a class of monotonic data $u_{0}(x, y),$ which holds that there exist positive constants $c, C>0$ such that
$$
\partial_{y} u_{0}(x, y) \geq c e^{-y} \quad \text { for } \quad y \in \mathbb{R}_{+},
$$
and for $k=0,1,2,$
$$
\left|\partial_{y}^{k}\left(u_{0}(x, y)-u_{0}^{s}(y)\right)\right|+\left|\partial_{x} \partial_{y} u_{0}(x, y)\right| \leq C e^{-y} \quad \text { for }(x, y) \in \mathbb{R}_{+}^{2}
$$
where $u_{0}^{s}(y) \in C^{2}\left(\mathbb{R}_{+}\right)$ satisfies
$$
u_{0}^{s}(0)=0, \quad\left|\partial^{k}\left(u_{0}^{s}(y)-1\right)\right| \leq C e^{-y} \quad \text { for } y \in \mathbb{R}_{+}\ \text {and } k=0,1,2.
$$

Let $\omega(y)$ be a nonnegative function in $\mathbb{R}_{+}$. They introduced the weighted $L^{p}$ norm:
$$
\|f\|_{L_{\omega}^{p}}=\|\omega(y) f(x, y)\|_{L^{p}}, \quad\|f\|_{L_{y, \omega}^{p}}=\|\omega(y) f(y)\|_{L^{p}},
$$
Let $k, \ell \in \mathbb{N},$ the weighted anisotropic Sobolev space $H_{\omega}^{k, \ell}$ consists of all functions $f \in L_{\omega}^{2}$ satisfying
$$
\|f\|_{H_{\omega}^{k, \ell}}^{2}=\sum_{\alpha \leq k} \sum_{\beta \leq \ell}\left\|\partial_{x}^{\alpha} \partial_{y}^{\beta} f\right\|_{L_{\omega}^{2}}^{2}<+\infty.
$$
They denoted by $H_{y, \omega}^{\ell}$ the weighted Sobolev space in $\mathbb{R}_{+},$ which consists of all functions $f \in L_{y, \omega}^{2}$ satisfying
$$
\|f\|_{H_{y, \omega}^{\ell}}^{2} \stackrel{\text { def }}{=} \sum_{\beta \leq \ell}\left\|\partial_{y}^{\beta} f\right\|_{L_{y, \omega}^{2}}^{2}<+\infty.
$$
When $\omega=1$, they denoted $H_{\omega}^{k, \ell}$ by $H^{k, \ell}$, and $H_{y, \omega}^{\ell}$ by $H_{y}^{\ell}$ for the simplicity.

\begin{Theorem}(\cite{WZ2})
Assume that the initial data satisfies the assumptions and $\widetilde{u}_{0} \in H_{\mu}^{3,1} \cap H_{\mu}^{1,2}$ with $\widetilde{u}_{0}(x, 0)=0 .$ There exists a constant $T>0$ and a unique solution of the Prandtl equation in $[0, T],$ which satisfies
$$
c_{1} e^{-y} \leq \partial_{y} u(t, x, y) \leq C_{1} e^{-y} \quad \text { for }(t, x, y) \in[0, T] \times \mathbb{R}_{+}^{2},
$$
$$
\left|\partial_{y}^{2} u(t, x, y)\right|+\left|\partial_{x} \partial_{y} u(t, x, y)\right| \leq C_{1} e^{-y} \quad \text { for }(t, x, y) \in[0, T] \times \mathbb{R}_{+}^{2},
$$
$\widetilde{u} \in L^{\infty}\left(0, T ; H_{\mu}^{3,1} \cap H_{\mu}^{1,2}\right) \cap L^{2}\left(0, T ; H_{\mu}^{1,3}\right)$
 for some $c_{1}, C_{1}>0 .$ Let $T^{*}$ be the maximal existence time of the solution. If the solution satisfies the following conditions
$$
\begin{array}{c}
c e^{-y} \leq \partial_{y} u(t, x, y) \leq C e^{-y} \quad \text { for } \quad(t, x, y) \in\left[0, T^{*}\right) \times \mathbb{R}_{+}^{2}, \\
\sup _{t \in\left[0, T^{*}\right)} A(t) \leq C,
\end{array}
$$
for some constants  $c, C>0,$ then it can be extended after $t=T^{*}$. Here
$$
A(t)=\sup _{(s, x, y) \in[0, t] \times \mathbb{R}_{+}^{2}} e^{y}\left(\sum_{k=0}^{2}\left|\partial_{y}^{k}(u-1)\right|
+(1+y)^{-1}\left(\left|\partial_{x} u\right|+\left|\partial_{x} \partial_{y} u\right|\right)\right),
$$
and $H_{\mu}^{k, \ell}$ is the weighted anisotropic Sobolev space.
\end{Theorem}
\begin{Theorem}(\cite{WZ2})
Let $u$ be a solution in $\left[0, T^{*}\right)$ constructed in the Theorem 1.1.46. Then there exists $c_{0}, C_{0}$ such that
$$
c_{0} e^{-y} \leq \partial_{y} u(t, x, y) \leq C_{0} e^{-y} \quad \text { for }(t, x, y) \in\left[0, T^{*}\right)
\times \mathbb{R}_{+}^{2}.
$$
Moreover, if
$$
\left|\partial_{y}^{2} u(t, x, y)\right| \leq C \quad \text { for }(t, x) \in\left[0, T^{*}\right)
 \times \mathbb{R}, y \in[0, \delta],
$$
for some constants $C>0, \delta>0,$ then
$$
A(t) \leq C_{*} \quad \text { for any } t \in\left[0, T^{*}\right),
$$
where $C_{*}>0$ is a constant depending only on $C_{0}, c_{0}, T^{*}, \delta .$
\end{Theorem}

Wang and Wang  (\cite{WW}) in 2020 under the assumption that the initial velocity and outer flow velocity are analytic in the horizontal variable, obtained the local well-posedness of the geophysical boundary layer problem by using energy method in the weighted Chemin-Lerner spaces. Moreover, when the initial velocity and outer flow velocity satisfy certain condition on a transversal plane, for any smooth solution decaying exponentially in the normal variable to the geophysical boundary layer problem, it is proved that its $W^{1, \infty}$-norm blows up in a finite time.

The authors in \cite{WW} considered the following initial boundary value problem in the domain $Q_{T}=\{0<t<$ $T, x \in \mathbb{R}, y>0\},$
\begin{equation}
\left\{\begin{array}{l}
\partial_{t} u+u \partial_{x} u+v \partial_{y} u-\int_{y}^{+\infty}(u-U) d y^{\prime}-\partial_{y}^{2} u=\partial_{t} U+U \partial_{x} U, \\
\partial_{x} u+\partial_{y} v=0, \\
\left.u\right|_{t=0}=u_{0}(x, y), \\
\left.(u, v)\right|_{y=0}=(0,0), \quad \lim\limits _{y \rightarrow+\infty} u(t, x, y)=U(t, x),
\end{array}\right. \label{y3.18}
\end{equation}
where $(u, v)$ is the velocity field, and $U(t, x)$ is the tangential velocity of the outer flow.

Let $u^{s}(t, x, y)=U(t, x) \phi(t, y)$ and $w=u-u^{s}$. From (\ref{y3.18}), we know that $w$ satisfies the following problem in $\{t>0, x \in \mathbb{R}, y>0\},$
\begin{equation}
\left\{\begin{array}{l}
\partial_{t} w+\left(w+u^{s}\right) \partial_{x} w+w \partial_{x} u^{s}-\int_{0}^{y} \partial_{x}\left(w+u^{s}\right) d y^{\prime} \partial_{y}\left(w+u^{s}\right)-\int_{y}^{+\infty} w d y^{\prime}-\partial_{y}^{2} w, \\
=(1-\phi)\left(\partial_{t} U+(1+\phi) U \partial_{x} U\right)-\int_{y}^{+\infty} U(1-\phi) d y^{\prime} \\
\left.w\right|_{y=0}=0, \lim \limits_{y \rightarrow+\infty} w=0 ,\\
\left.w\right|_{t=0}=w_{0}(x, y) \triangleq u_{0}(x, y)-U(0, x) \operatorname{erfc}\left(\frac{y}{2}\right).
\end{array}\right.
\end{equation}
Introduce the following function spaces with parameters $s>0, l \in \mathbb{N}_{+}$ and $p \in[1,+\infty)$.\\
(i) The space $B^{s}$ is the set of functions $u \in \mathcal{S}^{\prime}(\mathbb{R})$ such that
$$
\|u\|_{B^{s}}:=\sum_{k \in \mathbb{Z}} 2^{k s}\left\|\Delta_{k} u\right\|_{L^{2}(\mathbb{R})}<+\infty.
$$
(ii) The space $B_{\psi}^{s, l}$, with a positive function $\psi(y)$, is the set of functions $u \in \mathcal{S}^{\prime}\left(\mathbb{R}_{+}^{2}\right)$ such that
$$
\|u\|_{B_{\psi}^{s, l}}:=\sum_{j=0}^{l} \sum_{k \in \mathbb{Z}} 2^{k s}\left\|e^{\psi(y)} \Delta_{k} \partial_{y}^{j} u\right\|_{L^{2}\left(\mathbb{R}_{x} \times \mathbb{R}^{+}\right)}<+\infty.
$$
(iii) The space $\tilde{L}_{t}^{p}\left(B^{s}\right)$ is defined as the completion of $C([0, t] ; \mathcal{S}(\mathbb{R}))$ with the norm
$$
\|u\|_{\tilde{L}_{t}^{p}\left(B^{s}\right)}:=\sum_{k \in \mathbb{Z}} 2^{k s}\left(\int_{0}^{t}\left\|\Delta_{k} u\left(t^{\prime}, \cdot\right)\right\|_{L^{2}(\mathbb{R})}^{p} d t^{\prime}\right)^{\frac{1}{p}}.
$$
(iv) For any positive function $\psi\left(t^{\prime}, y\right)$ and nonnegative $f\left(t^{\prime}\right) \in L_{l o c}^{1}\left(\mathbb{R}^{+}\right),$ the space $\tilde{L}_{t, f}^{p}\left(B_{\psi}^{s, l}\right)$ is defined as the completion of $C\left([0, t] ; \mathcal{S}\left(\mathbb{R}_{+}^{2}\right)\right)$ with the norm
$$
\|u\|_{\tilde{L}_{t, f}^{p}\left(B_{\psi}^{s, l}\right)}:=\sum_{j=0}^{l} \sum_{k \in \mathbb{Z}} 2^{k s}\left(\int_{0}^{t} f\left(t^{\prime}\right)\left\|e^{\psi\left(t^{\prime}, y\right)} \Delta_{k} \partial_{y}^{j} u\left(t^{\prime}, \cdot\right)\right\|_{L^{2}\left(\mathbb{R}_{x} \times \mathbb{R}^{+}\right)}^{p} d t^{\prime}\right)^{\frac{1}{p}}.
$$
Denote $\tilde{L}_{t, 1}^{p}\left(B_{\psi}^{s, l}\right)$ by $\tilde{L}_{t}^{p}\left(B_{\psi}^{s, l}\right)$ for simplicity when $f\left(t^{\prime}\right) \equiv 1 .$ The above notations can be properly modified to the case $p=+\infty$.

To obtain energy estimates for solutions, we introduce the weights
$$
\psi(t, y)=\frac{1+y^{2}}{16(1+t)^{\gamma}} \quad \text { and } \quad \psi_{0}(y)=\frac{1+y^{2}}{16},
$$
with $\gamma \geq 2$, then there holds
$$
\sqrt{-\left(\psi_{t}+2 \psi_{y}^{2}\right)} \geq(1+y) \sqrt{\frac{\gamma-1}{32(1+t)^{\gamma+1}}}.
$$

\begin{Theorem}(\cite{WW})
For a given $T_{0}>0$, assume that the initial velocity $w_{0}(x, y)$ and the outer flow velocity $U(t, x)$ are analytic in $x \in \mathbb{R},$ and
$$
e^{\langle D\rangle} w_{0} \in B_{\psi_{0}}^{\frac{1}{2}, 0},
$$
and
$$
e^{\langle D\rangle} U \in \tilde{L}_{T_{0}}^{\infty}\left(B^{\frac{1}{2}}\right) \cap \tilde{L}_{T_{0}}^{\infty}\left(B^{1}\right) \cap \tilde{L}_{T_{0}}^{\infty}\left(B^{\frac{3}{2}}\right), \quad e^{\langle D\rangle} U_{t} \in \tilde{L}_{T_{0}}^{\infty}\left(B^{\frac{1}{2}}\right).
$$
Then, there exists $0<T^{*} \leq T_{0}$ such that the problem has a unique solution $e^{\Phi(t, D)} w \in \tilde{L}_{T^{*}}^{\infty}\left(B_{\psi}^{\frac{1}{2}, 0}\right) $ and $\Phi(t, D)$ is the operator
associated with the symbol $\Phi(t, \xi)$. $D$ denotes the Fourier multiplier.
\end{Theorem}

\begin{Theorem}(\cite{WW})
For a given $T>0$, assume that the initial data and outer flow satisfy
$$
u_{0}(0, y)=U(t, 0)=0, U_{x}(t, 0) \geq 0, u_{0 x}(0, y) \leq U_{x}(0,0),
$$
and there is $M \geq 0$ depending on $T$ and $\left\|U_{x}(\cdot, 0)\right\|_{L^{\infty}([0, T])}$ such that
$$
\int_{0}^{+\infty} \rho(y)\left(U_{x}(0,0)-u_{0 x}(0, y)\right) d y>M,
$$
for a weight function $\rho(y)$ if the smooth solution $u \in C^{3}\left(Q_{T}\right)$ of the first equation \eqref{y3.18} satisfies
$$
\sup _{0 \leq t \leq T} \int_{0}^{+\infty}(u(t, 0, y))^{2} e^{2 \mu y} d y<+\infty,
$$
and
$$
\lim _{y \rightarrow+\infty} \sup _{0 \leq t \leq T}\left(\left|u_{x}(t, 0, y)-U_{x}(t, 0)\right| e^{\mu y}\right)=0,
$$
for some $\mu>0,$ then the $W^{1, \infty}$- norm of $u$ blows up in a finite time.
\end{Theorem}

The authors in \cite{KMVW,mw,GHY} studied the local well-posedness of the Prandtl equations by the energy   method.
In the 2D setting, it is convenient to study the evolution of the vorticity $w=\partial_y u$. The nonlinear transport
equation is as follows,

\bqa \partial_t w+u\partial_xw+v\partial_y w=0,\label{+6}\eqa
where, one may compute $(u,v)$ from $w$ via

\bqa u=-\partial_y\mathcal{A}(w)\ \text{and}\ v=\partial_x\mathcal{A}(w),\label{+7}\eqa

where the stream function $\mathcal{A}(w)$ solves $$-\partial_{yy}\mathcal{A}(w)=w$$
with the boundary condition $\mathcal{A}(w)|_{y=0,1}=0$.

Kukavica,  Masmoudi,  Vicol and  Wong (\cite{KMVW}) in 2014 found a new class of data for which the Prandtl boundary layer equations and the
hydrostatic Euler equations are locally in time well-posed.   They showed that the local existence and uniqueness of solution hold.
 \begin{Theorem}(\cite{KMVW})\label{ }
Let $I\subset\mathbb{R}$ be an open interval, $s\geq 4$ be an even integer, $\gamma\geq 1$, $\sigma>\gamma+\frac{1}{2}$, and $\delta\in(0,\frac{1}{4})$.
Assume that the initial velocity obeys $u_{0}-U_{0}\in H^{s,\gamma-1}(I)$, the initial vorticity satisfies $w_{0}\in H^{s,\gamma}_{\sigma,2\delta}(I)$,
and that the outer Euler flow $U$ is sufficiently smooth (for example, the following equation will suffice)
$$\|U\|_{W^{s+9,\infty}_{t,x}}=\sup\limits_{t}\sum\limits_{0\leq 2l\leq s+9}\|\partial_{t}^{l}U\|_{W^{s-2l+9,\infty}(I)}< \infty .$$
Then there exists a $T>0$ and a smooth solution $u$ of (\ref{1.1.6}) such that the velocity obeys
$u-U\in L^{\infty}([0,T];H^{s,\gamma-1}(I))\cap C_{w}([0,T];H^{s}_{loc}(I))$ and the vorticity obeys
$w\in L^{\infty}([0,T];H^{s,\gamma}_{\sigma,\delta}(I))\cap C_{w}([0,T],H^{s}_{loc}(I))$. When $I=\mathbb{R}$ or $I=[0,T]$, the solution constructed above is
the unique solution in this regularity class. When $\partial I\neq\emptyset$, e.g., if $I=(a,b)$, there exists a positive $M\leq \infty$ such that $u(t)$ is
the unique solution with
$u(t)-U(t)\in H^{s,\gamma-1}(I_{t})$ and $w(t)\in H^{s,\gamma}_{\sigma,\delta}(I_{t})$, where $I_{t}=\{x\in\mathbb{R}:(x-Mt,x+Mt)\subseteq I\}$.

\end{Theorem}

\begin{Theorem}(\cite{KMVW})
Let $\sigma,\tau_{0}>0$ and $s\geq 4$. Assume that the initial vorticity satisfies
a uniform Rayleigh condition on $\mathfrak{D}_{m}$, we have $\partial_{y} u_{0}\in H^{s}_{2\sigma,loc}(\mathfrak{D}_{m})$, and $\partial_{y} u_{0}$ is uniformly
real-analytic on $\mathfrak{D}_{a}$ with radius of analyticity at least $\tau_{0}$, i.e., $u_{0}\in Y_{ \tau_{0}}(\mathfrak{D}_{a})$. Also
assume that the initial velocity $u_{0}$ satisfies the compatibility condition. Then there exist $T>0$ and a unique smooth solution
 $\partial_{y} u\in C(0,T;X_{\tau_{0}/2}(\mathfrak{D}_{a}))\cap H^{s}_{\sigma,loc}(\mathfrak{D}_{m})$
of (\ref{+6})-\eqref{+7} on $[0,T]$, which has zero vertical mean.\\
Here
\begin{equation*}
H^{s,\gamma}_{\sigma,\delta}(I)=\left\{ w\in H^{s,\gamma}(I):(1+y)^{\sigma}w(x,y)\geq \delta,
 \sum\limits_{|\alpha|\leq 2}|(1+y)^{\sigma+\alpha_{2}}D^{\alpha}w|\leq \delta^{-2} \text{for all } (x,y)\in  I\times \mathbb{R}_{+}  \right\}
 \label{ }
 \end{equation*}
$$  \mathfrak{D}_{m}=(-\infty, 1) \times (0,1), \ \   \mathfrak{D}_{a}=(0,  \infty ) \times (0,1). $$
\end{Theorem}

  Masmoudi and  Wong (\cite{mw}) in 2015 proved local existence and uniqueness for the 2D Prandtl system in weighted Sobolev spaces under
  the Oleiniks monotonicity assumption.  This new energy estimate is based on a cancellation
property which is valid under the monotonicity assumption, but not using the Crocco transform. \\

They considered the   Prandtl equations  (\ref{1.1.6}) in a periodic domain
$ \mathbb{T}\times \mathbb{R}_{+} =\{(x,y)| x\in  \mathbb{R}/\mathbb{Z}, y\in[0, +\infty)  \}$,
$U(t,x)$  satisfy Bernoulli's law (\ref{1.1.3}). In \cite{mw}, let $w=\partial_{y} u$, and denote
$$H^{s,\gamma}_{\sigma,\delta}:=\left\{w:\mathbb{T}\times\mathbb{R_{+}}\rightarrow\mathbb{R},
\|w\|_{H^{s,\gamma}}\leq +\infty, (1+y)^{\sigma}|w|\geq\delta,
\sum \limits_{|\alpha|\leq 2}  |(1+y)^{\sigma+\alpha_{2}}D^{\alpha}w | ^{2}\leq\frac{1}{\delta^{2}} \right\}.$$

\begin{Theorem} (\cite{mw})
 (Local $H^{s,\gamma}_{\sigma,\delta}$ Existence and Uniqueness to the Prandtl Equations ). Let $s\geq 4$ be an even integer,
 $\gamma\geq 1,\sigma>\gamma+\frac{1}{2}$ and $\delta\in(0,\frac{1}{2})$. For simplicity, we suppose that the outer flow $U$ satisfies
$$\sup\limits_{t}\| |U| \|_{s+9,\infty}:=\sup\limits_{t} \sup\limits_{l=0}^{[\frac{s+9}{2}]} \|\partial_{t}^{l}U\|_{W^{s-2l+9,\infty}(\mathbb{T})}< +\infty .$$
Assume that $u_{0}-U\in H^{s,\gamma-1}$ and the initial vorticity $w_{0}=\partial_{y}u_{0}\in H^{s,\gamma}_{\sigma,2\delta}$. In addition, when $s=4$,
 we further assume that $\delta>0$ is chosen small enough such that
$$\|w_{0}\|_{H^{s,\gamma}_{g} }\leq C\delta^{-1},$$
where the norm $\|\cdot\|_{H^{s,\gamma}_{g} }$ will be defined by (\ref{1.13}) and $C$ is a universal constant. Then there exist a time
$T:=T(s,\gamma,\sigma,\delta, \|w_{0}\|_{H^{s,\gamma}  },U)>0$ and a unique classical solution $(u,v)$ to the Prandtl equations (\ref{1.1.6}) such that
$u-U\in L^{\infty}([0,T],H^{s,\gamma-1})\cap C([0,T]; H^{s}-w)$
and the vorticity $w:=\partial_{y} u\in L^{\infty}([0,T]; H^{s,\gamma}_{\sigma,\delta} )\cap C([0,T]; H^{s}-w)$, where $ H^{s}-w$ is the space $H^{s}$
endowed with its weak topology.
\end{Theorem}
Here
\begin{equation}
\|w\|_{H^{s,\gamma}_{g}}^{2}:=\|(1+y)^{\gamma}g_{s}\|^{2}_{L^{2}}+
\sum\limits_{\substack{ |\alpha|\leq s\\ \alpha_{1}\leq s-1}}  \|(1+y)^{\gamma+\alpha_{2}}D^{\alpha}w\|_{L^{2}}^{2}.
\label{1.13}
\end{equation}

Gao, Huang and Yao (\cite{GHY}) in 2018 investigated the local-in-time well-posedness for the two-dimensional Prandtl equations in weighted Sobolev spaces under the Oleinik
 monotonicity condition.
Throughout this paper, they are concerned with the two-dimensional Prandtl equations  in a periodic domain $\mathbb{T}\times\mathbb{R}_+$,
\begin{eqnarray}\left\{\begin{array}{ll}
u_t+u\partial_xu+v\partial_yu-\partial_{y}^2u+\partial_xp=0,\\
\partial_xu+\partial_yv=0,\\
(u,v)|_{y=0}=0, \lim\limits_{y\rightarrow+\infty}u=U,\\
u|_{t=0}=u_0.
\end{array}\right.\label{30.1}\end{eqnarray}
Define the weighted Sobolev space
$$H^{s,\gamma}_{\sigma,\delta}:=\Big\{w:\|w\|_{H^{s,\gamma}}<+\infty, (1+y)^{\sigma}w\geq \delta, \sum\limits_{|\alpha|\leq 2}|(1+y)^{\sigma+\alpha_{2}}D^{\alpha }w|^{2}\leq \frac{1}{\delta^{2}}\Big\},$$
where $s\geq 4,\gamma\geq 1,\sigma>\gamma+\frac{1}{2}, \delta\in (0,\frac{1}{2})$.

\begin{Theorem}(\cite{GHY})
Let $s \geq 4$ be an even integer, $\gamma\geq 1,\sigma>\gamma+\frac{1}{2}, \delta\in (0,\frac{1}{2})$. Suppose the outer flow $U$ satisfies,
\begin{eqnarray*}
M_U:=\sup\limits_{t}\sum\limits_{t=0}^{\frac{s+1}{2}}\|\partial^{l}_{t}U\|_{W^{s-2l+2}(\mathbb{R})}<+\infty.
\end{eqnarray*}
Assume that the initial tangential velocity,$u_{0}-U|_{t=0}\in H^{s,\gamma-1}$ and the initial vorticity $w_{0}:=\partial_{y}u_{0}\in H^{s,\gamma}_{\sigma, 2\delta}$. Then there exist a time  $T=T(s,\gamma,\delta,\delta_0,M_U)>0$  and a unique classical solution $(u,v)$ to the Prandtl equations \eqref{30.1} such that
\begin{eqnarray*}
\sup\limits_{0\leq t\leq T}\|w\|^2_{H^{s,\gamma}_{\sigma,\delta}}\leq C_{s,\gamma,\delta}(1+\|w_0\|^8_{H^{s,\gamma}_{\sigma,\delta}}+M_U^4)<+\infty.
\end{eqnarray*}
\end{Theorem}
\begin{Theorem}(\cite{GHY})
Let $s \geq 4$ be an even integer, $\gamma\geq 1,\sigma>\gamma+\frac{1}{2}, \delta\in (0,\frac{1}{2})$, and $\varepsilon\in(0,1)$,  the smooth solution $(u^\varepsilon,v^\varepsilon,w^\varepsilon)$ to the regularized Prandtl equations \eqref{30.2}-\eqref{30.4}. Under the assumptions of the above theorem,  there exists a time $T_a=T_a(s,\gamma,\delta,\delta_0,M_U)>0$, independent of $\varepsilon$, such the following estimates hold on
\begin{eqnarray*}
\sup\limits_{0\leq t\leq T}\|w^\varepsilon\|^2_{H^{s,\gamma}_{\sigma,\delta}}
\leq C_{s,\gamma,\delta}(1+\|w^\varepsilon_0\|^8_{H^{s,\gamma}_{\sigma,\delta}}+M_U^4)<+\infty.
\end{eqnarray*}
\end{Theorem}

Denote vorticity $w^\varepsilon=\partial_y u^\varepsilon$,
using the the regularized Prandtl equations, we find that this vorticity satisfies the following evolution equations:
\begin{eqnarray}\left\{\begin{array}{ll}
w^\varepsilon_t+u^\varepsilon\partial_xw^\varepsilon+v\partial_yw^\varepsilon-\partial_{y}^2w^\varepsilon=0,\\
w^\varepsilon|_{t=0}=w_0=\partial_yu_0,\\
\partial_yw^\varepsilon=\partial_xp,
\end{array}\right.\label{30.2}\end{eqnarray}
where the velocity field $(u^\varepsilon,v^\varepsilon)$ is given by
\begin{eqnarray}
 u^\varepsilon(t,x,y)=U(t,x)-\int_y^{+\infty}w^\varepsilon dy\label{30.3}
\end{eqnarray}
and
\begin{eqnarray}
 v^\varepsilon(t,x,y)=-\int_0^y\partial_xu^\varepsilon dy.  \label{30.4}
\end{eqnarray}

The authors in \cite{KV,gvm,gmm,lwx,dgv} studied the local well-posedness of  the Prandtl equations in Gevrey space.

Kukavica and Vicol   (\cite{KV}) in 2013 addressed the local well-posedness of the Prandtl boundary layer equations. Using a new change
of variables, the authors in \cite{KV} allowed for more general data than previously considered, that is, they required the matching at the
top of the boundary layer to be at a polynomial rather than exponential rate. The proof is direct, via analytic
energy estimates in the tangential variables.

 The foundations for the boundary layer theory were laid by Prandtl, who made the ansatz
$u^{N S}(x, \tilde{y}, t)=(u(x, \tilde{y} / \sqrt{\nu}, t), \sqrt{\nu} w(x, \tilde{y} / \sqrt{\nu}, t))$. Inserting this velocity field in the Navier-Stokes equations and sending the kinematic viscosity $\nu$ to zero, one formally obtains the Prandtl boundary layer equations for the unknown velocity field $(u, \sqrt{\nu} w)$
\begin{equation}
\left\{\begin{array}{l}
\partial_{t} u-\partial_{Y Y} u+u \partial_{x} u+w \partial_{Y} u+\partial_{x} P=0, \\
\partial_{x} u+\partial_{Y} w=0, \\
\partial_{Y} P=0,
\end{array}\label{y1.7}  \right.
\end{equation}
in $\mathbb{H}=\left\{(x, Y) \in \mathbb{R}^{2}: Y>0\right\},$ where $Y=\tilde{y} / \sqrt{\nu}$ is the normal variable in the boundary layer.  For simplicity of the presentation, they in \cite{KV} considered the two-dimensional setting, but all the methods and results presented here extend to the three-dimensional case as well. The system $(\ref{y1.7})$ is supplemented with the no-slip and the no-influx boundary conditions
\begin{equation}\left\{\begin{array}{ll}
u(x, Y, t)|_{Y=0}=0, \\
w(x, Y, t)|_{Y=0}=0,
\end{array}\non\right.\end{equation}
for $t>0,$ and the matching conditions with the Euler flow as $Y \rightarrow \infty$, via the Bernoulli's law
\begin{equation}\left\{\begin{array}{ll}
\lim \limits_{Y \rightarrow \infty} u(x, Y, t)=U(x, t),\\
\partial_{x} P(x, t)=-\left(\partial_{t}+U(x, t) \partial_{x}\right) U(x, t),
\end{array}\non\right.\end{equation}
for $x \in \mathbb{R}^{2}, t>0,$ where $U(x, t)$ is the trace at $\tilde{y}=0$ of the tangential component of the Euler flow $u^{E}$.

Under the change of variables, the Prandtl system reads
$$
\partial_{t} v-A^{2} \partial_{y y} v+N(v)+L(v)=F
$$
where
\begin{eqnarray*}
N(v)&=&v \partial_{x} v-\partial_{x} W(v) \partial_{y} v+\partial_{x} a W(v) \partial_{y} v, \quad
W(v)(x, y)=\int_{0}^{y} v(x, \zeta) d \zeta, \\
L(v)&=&\partial_{x} W(v) \partial_{y} \phi U+\partial_{x} v(1-\phi) U+\partial_{y} v\left(\Phi \partial_{x} U-\partial_{x} a \Phi U\right)-W(v) \partial_{x} a \partial_{y} \phi U+v(1-\phi) \partial_{x} U, \\
F&=&\left(\phi(1-\phi)+\Phi \partial_{y} \phi\right) U \partial_{x} U-\partial_{x} a \partial_{y} \phi \Phi U^{2}-A^{2} \partial_{y y} \phi U-\phi \partial_{x} P.
\end{eqnarray*}

The system \eqref{y1.7} is supplemented with the boundary conditions
$$
\begin{array}{l}
\left.v(x, y, t)\right|_{y=0}=\left.u(x, Y, t)\right|_{Y=0}-(1-\phi(0)) U(x, t)=0, \\
\lim\limits _{y \rightarrow \infty} v(x, y, t)=\lim _{Y \rightarrow \infty} u(x, Y, t)-U(x, t)=0,
\end{array}
$$
for all $(x, t) \in \mathbb{R} \times[0, \infty),$ and initial condition
$$
\left.v(x, y, t)\right|_{t=0}=v_{0}(x, y)=u_{0}(x, Y)-(1-\phi(y)) U_{0}(x) .
$$
The authors of \cite{KV}  considered the $x$ -analytic norm with $y$ -weight given by
$$
\rho(y)=\langle y\rangle^{\alpha}
$$
for some $\alpha>0$ to be fixed later. Namely, for a function $V(x, y)$ and a number $\tau_{0}>0$, denote
$$
\|V\|_{X_{\tau_{0}}}^{2}=\sum_{m \geq 0}\left\|\rho(y) \partial_{x}^{m} V(x, y)\right\|_{L^{2}(\mathbb{H})}^{2} \tau_{0}^{2 m} M_{m}^{2},
$$
where the analytic weights $M_{m}$ are defined as
$$
M_{m}=\frac{(m+1)^{r}}{m !},
$$
for some $r>0$ to be determined. If $v(x, y, t)$ and $\tau(t)$ have $t$ -dependence, similarly denote
$$
\|v(t)\|_{X_{\tau(t)}}^{2}=\sum_{m \geq 0}\left\|\rho(y) \partial_{x}^{m} v(x, y, t)\right\|_{L^{2}(\mathbb{H})}^{2} \tau(t)^{2 m} M_{m}^{2}.
$$
If the $t$ dependence is clear from the context, we will omit it. Since the weight $\rho(y)$ does not depend on $x$, the analytic norm may also be written as
$$
\|v\|_{X_{\tau}}^{2}=\sum_{m \geq 0}\|\rho v\|_{\dot{H}_{x}^{m}}^{2} \tau^{2 m} M_{m}^{2}.
$$
For a positive number $\tau>0$,  we write $v \in X_{\tau}$ if $\|v\|_{X_{\tau}}<\infty$.

\begin{Theorem}(\cite{KV})
Fix real numbers $\alpha>1 / 2, \theta>\alpha+1 / 2,$ and $r>1 .$ Assume that the initial data for the underlying Euler flow is uniformly real analytic, with radius of analyticity at least $\tau_{E}>0$ and analytic norm bounded by $G_{E}>0 .$ There exists $\tau_{0}=\tau_{0}\left(r, \tau_{E}, G_{E}\right)>0$ such that for all $v_{0} \in X_{\tau_{0}}$ there is $T_{*}=T_{*}\left(r, \alpha, \theta, \tau_{E}, G_{E}, \tau_{0},\left\|v_{0}\right\|_{X_{\tau_{0}}}\right)>0$ such that the initial value problem has a unique real-analytic solution on $\left[0, T_{*}\right]$.

\end{Theorem}

 G$\acute{e}$rard-Varet and  Masmoudi (\cite{gvm}) in 2015 showed that the Prandtl system (\ref{1.1.6})  is actually locally well-posed for data that
  belong to the Gevrey class $7/4$ in the horizontal variable $x$. Their result improved the classical local well-posedness result for data that are
  analytic in $x$ (that is Gevrey class 1). \\
In \cite{gvm},    let $m\leq 1$, they  defined  the Gevrey space $G^{m}(\mathbb{T})$ which is the set of functions $f$ satisfying:
there exist $C,\tau>0$ such that
$$|f^{j}(x)|\leq  C\tau^{-j}(j!)^{m},\ \  \forall j\in \mathbb{N},\ \ x\in \mathbb{T}, $$
$G^{m}(\mathbb{T})= \cup_{\tau>0}G^{m}_{\tau}(\mathbb{T}) $,
where
$$ G^{m}_{\tau}(\mathbb{T}):=\left\{f\in C^{\Omega}(\mathbb{T}), \sup\limits_{j} \tau^{j}(j!)^{-m}(j+1)^{10}\|f ^{(j)}\|_{L^{2}}<+ \infty  \right\}$$
is a Banach space. For $s\in \mathbb{N},\gamma\geq 0$,  define the spaces
$$H^{s}_{\gamma}:= \left\{ g\in H^{s}(\mathbb{R}_{+}),\ \  (1+y)^{\gamma+k}g^{(k)}\leq L^{2}(\mathbb{R}_{+}), k=0,...s \right\},
\|g\|_{H^{s}_{\gamma}}:= \sum\limits_{k=0}^{s} \|(1+y)^{\gamma+k}g^{(k)}\|_{L^{2}}^{2} .$$
They assumed the following  assumptions:\\
$(H) $ \ \  $w_{0}(x,y)=0$ if $y=a_{0}(x)$, for some curve $a_{0}>0$, with $\partial_{y}w_{0}(x,a_{0}(x))>0$, $\forall x$.\\
$ (H')$ there exist  $ \sigma,\delta>0$ such that for all $y>3$, for all $x$, for all $\alpha\in\mathbb{N}^{1}$, $|\alpha|\leq 2$,
$$|w_{0}(x,y)|\geq \frac{2\delta}{(1+y)^{\sigma}},  \ \  |\partial^{\alpha}w_{0}(x,y)|\leq \frac{1}{2\delta (1+y)^{\sigma+\alpha_{2}}} .$$

\begin{Theorem}(\cite{gvm})
Let  $\tau_{0}>0,s>8$ even, $\gamma\geq 1,\sigma\geq \gamma+\frac{1}{2},\delta>0$. Let $U$ be a constant and $u_{0}$ satisfy
$$u_{0}\in G^{7/4}_{\tau_{0}}(\mathbb{T}, H^{s+1}_{\gamma-1} ),\ \  w_{0}:=\partial_{y}u_{0}\in   G^{7/4}_{\tau_{0}}(\mathbb{T}, H^{s }_{\gamma } ). $$
Assume the compatibility condition: $u_{0}|_{y=0}=0$, and assumptions $(H),(H')$ above hold. Then there exists $T>0,0<\tau<\tau_{0}$ and a unique solution
$$ u\in L^{\infty}(0,T; G^{7/4}_{\tau }(\mathbb{T}, H^{s+1}_{\gamma-1} ), w\in L^{\infty}(0,T; G^{7/4}_{\tau }(\mathbb{T}, H^{s }_{\gamma } ),  $$
of the Prandtl equations  (\ref{1.1.6})  with initial data $u_{0}$.

\end{Theorem}

 G$\acute{e}$rard-Varet,  Maekawa and Masmoudi (\cite{gmm}) in 2016   investigated the stability of boundary layer solutions of the
  two-dimensional incompressible Navier-Stokes equations (\ref{1.1.1}). They considered shear flow solutions of Prandtl type :
$$u^{\nu}(t,x,y)=(U^{E}(t,y)+U^{BL}(t,\frac{y}{\sqrt{\nu}}), 0), \ \  0\leq \nu\ll 1.$$
They showed that if $U^{BL}$ is monotonic, then $u^{\nu}$ is stable over some time interval $(0,T)$, $T$ independent of $\nu$, under perturbations with Gevrey
regularity in $x$ and Sobolev regularity in $y$. They improved  the classical stability results of Sammartino and Caflisch (\cite {SC1,SC})  in analytic class
 (both in $x$ and $y$),
where $u^{\nu}(t,x,y)$ is the solution of (\ref{1.1.1}),  $(U^{E},V^{E})$ is the solution of (\ref{1.1.4}),  $(U^{BL},V^{BL})$ is a corrector localized near
the boundary,   $(U^{P},V^{P})$ is the solution of   Prandtl equations,
they emphasized that the condition is always satisfied when $U^{P}$ is the solution to the problem
\begin{equation} \left\{
\begin{array}{ll}
\partial_{t} U^{P}-\partial_{Y}^{2}U^{P}=0, \ \ \ \ t>0, Y>0,\\
U^{P}|_{Y=0}=0, \ \  \lim\limits_{Y\rightarrow\infty} U^{P}=U^{E}|_{y=0}, \ \ \ \  t\geq 0,\\
U^{P}|_{t=0}=U.
\end{array}
 \label{1.9}         \right.\end{equation}
They assumed that
 \begin{equation}
\begin{array}{ll}
(i)~ U|_{Y=0}=0, \lim\limits_{Y\rightarrow\infty} U =U^{E}|_{y=0}, \ \  U^{E},U\in BC^{2}(\mathbb{R}_{+}), \\
(ii)~ \|U\|:=\sum\limits_{k=0,1,2}\sup\limits_{Y\geq 0} (1+Y)^{k}|\partial_{Y}^{k}U(Y)|< \infty, \\
(iii)~ \text{For each}~\sigma\in(0,1] ~\text{there exists}~ M_{\sigma}>0 ~\text{such that }~-M_{\sigma}\partial_{Y}^{2}U\geq (\partial_{Y} U)^{2} ~\text{for}~Y\geq \sigma ,
\end{array}
 \label{1.10}          \end{equation}
and  assumed
 \begin{equation}
\begin{array}{ll}
(i)\ U|_{Y=0}=0, \lim\limits_{Y\rightarrow\infty} U =U^{E}|_{y=0}, \ \  U^{E},U\in BC^{2}(\mathbb{R}_{+}), \\
(ii)\ \|U\|:=\sum\limits_{k=0,1,2}\sup\limits_{Y\geq 0} (1+Y)^{k}|\partial_{Y}^{k}U(Y)|< \infty, \\
(iii)\   \text{There exists an}~ M >0 ~\text{such that }~-M \partial_{Y}^{2}U\geq (\partial_{Y} U)^{2} ~\text{for}~Y\geq 0 .
\end{array}
 \label{1.11}        \end{equation}
 To state our main theorem, we need to introduce a few notations. Let

 $$(\mathcal{P}_nf)(y)=f_n(y)e^{inx},\ f_n(y)=\frac{1}{2\pi}\int_0^{2\pi}f(x,y)e^{-inx}dx,\ n\in\mathbb{Z},$$
 the projection on the Fourier mode $n$ in $x$. We then introduce, for $\gamma\in (0,1]$, $d\geq0$ and $K\geq0$,
  the Banach space $X_{d,\gamma,K}$ as
  \bqa X_{d,\gamma,K}=\{f\in L^2_\sigma(\Omega)|\|f\|_{X_{d,\gamma,K}}
  =\sup\limits_{n\in\mathbb{Z}}(1+|n|^d)e^{K\theta_{\gamma,n}}\|\mathcal{P}_nf\|_{L^2(\Omega)}<\infty\},\label{+8}\eqa
 $$\theta_{\gamma,n}=|n|^\gamma(1+(1-\gamma)log(1+|n|)).$$

\begin{Theorem} (\cite{gmm})
Assume that (\ref{1.9}) holds for some $T>0$. For any $\gamma\in[\frac{7}{9},1)$, $d>\frac{9}{2}-\frac{7}{2}\gamma$, and $K>0$,
there exist $C,T', K'>0$ such that the following statement holds for any sufficiently small $\nu>0$.
If $\|a\|_{X_{d,\gamma,K}}\leq \nu^{\frac{1}{2}+\beta}$ with $\beta=\frac{2(1-\gamma)}{\gamma}$, then the system (\ref{1.1.1}) admits a unique solution
$u\in C([0,T'];L^{2}_{\sigma}(\Omega))\cap L^{2}([0,T'];W^{1,2}_{0}(\Omega))$ satisfying the estimate
$$\sup\limits_{0< t\leq T'}(\|u(t)\|_{X_{d,\gamma,K}}+ (\nu t)^{\frac{1}{4}}\|u(t)\|_{L^{\infty}(\Omega)} + (\nu t)^{\frac{1}{2}}\|\nabla u(t)\|_{L^{2}(\Omega)}  )\leq C \|a\|_{X_{d,\gamma,K}}.$$
If $U^{E}(t,y)=U^{E}(y)$ and $U^{P}(t,Y)=U(Y)$ are steady and satisfy (\ref{1.10}) instead of (\ref{1.9}),
the above result holds for any $\gamma\in[ \frac{5}{7} ,1]$. \\
If $U^{E}(t,y)=U^{E}(y)$ and $U^{P}(t,Y)=U(Y)$ are steady and satisfy (\ref{1.11}) instead of (\ref{1.9}),
the  above result holds for any $\gamma\in[ \frac{2}{3} ,1]$.
\end{Theorem}

  Li,  Wu and Xu (\cite{lwx}) in 2016 studied the Gevrey smoothing effects of the local solution for the Prandtl boundary layer equation (\ref{1.1.6}) with
 $U$ being a constant.

\begin{Theorem}(\cite{lwx})
Let $u(t,x,y)$ be a classical local in time solution to Prandtl equation (\ref{1.1.6}) which $U$ is a constant on $[0,T]$ with the properties subsequently
listed below: \\
(i) There exist two constants $C_\ast>1,\sigma>\frac{1}{2}$  such that for any  $(t,x,y)\in [0,T]\times \mathbb{R}_+^2$,
\begin{eqnarray*}
 && C_\ast\langle y\rangle^{-\sigma}\leq u_y(t,x,y)\leq C_\ast \langle y\rangle^{-\sigma}, \\
&&|u_{yy}(t,x,y)|+|u_{yyy}(t,x,y)|\leq C_\ast \langle y\rangle^{-\sigma},
\end{eqnarray*}
where $\langle y\rangle=(1+|y|^2)^{\frac{1}{2}}$.\\
(ii) There exist constants $c > 0,C_0 > 0$ and integer  $N_0\geq7$  such that
$$\|e^{2cy}u_x\|_{L^\infty([0,T];H^{N_0})}+\|e^{2cy}u_{xy}\|_{L^\infty([0,T];H^{N_0})}\leq C_0. $$
Then for any  $0 < T_1 < T$, there exists a constant $L>0$, such that for any $0 < t\in(0 ,T_1]$,
$$\forall m>1+N_0, \ \|e^{\bar{c}y}\partial_x^mu\|_{L^2(\mathbb{R}_+^2)}\leq t^{-3(m-N_0-1)L^m(m!)^{3(1+\sigma)}},$$
where $0<\bar{c}<c$. The constant $L$ depends only on $C_0 ,T_1 ,C_\ast ,c,\bar{c}$, and $\sigma$. Therefore,
the solution $u$ belongs to the Gevrey class of index  $3(1 + \sigma)$ with respect to $x \in \mathbb{R}$ for
any $0 < t \leq T_1$.
\end{Theorem}

Dietert and  G$\acute{e}$rard-Varet (\cite{dgv}) in 2019 showed the local in time well-posedness of the following Prandtl equations for data with Gevrey
2 regularity in $x$ and Sobolev regularity in $y$,
\begin{equation}
\left\{\begin{array}{c}
\partial_{t} U^{P}+U^{P} \partial_{x} U^{P}+V^{P} \partial_{y} U^{P}-\partial_{y}^{2} U^{P}=\partial_{t} U^{E}+U^{E} \partial_{x} U^{E}, \\
\partial_{x} U^{P}+\partial_{y} V^{P}=0,
\end{array}\right.\label{y2.13}
\end{equation}
 in the domain $\Omega=\mathbb{T} \times \mathbb{R}_{+},$ completed with boundary conditions
\begin{equation}
\left.U^{P}\right|_{y=0}=\left.V^{P}\right|_{y=0}=0, \quad \lim _{y \rightarrow+\infty} U^{P}=U^{E}.\label{y2.14}
\end{equation}

Let $\gamma \geq 1, \tau>0, r \in \mathbb{R}$. For functions $f=f(x)$ of one variable, the authors in \cite{dgv} defined the Gevrey norm
$$
|f|_{\gamma, \tau, r}^{2}=\sum_{j \in \mathbb{N}}\left(\frac{\tau^{j+1}(j+1)^{r}}{(j !)^{\gamma}}\right)^{2}\left\|f^{(j)}\right\|_{L^{2}(\mathbb{T})}^{2},
$$
and for functions $f=f(x, y)$ of two variables, the norm
$$
\|f\|_{\gamma, \tau, r}^{2}=\sum_{j \in \mathbb{N}}\left(\frac{\tau^{j+1}(j+1)^{r}}{(j !)^{\gamma}}\right)^{2}\left\|\partial_{x}^{j} f\right\|_{j}^{2},
$$
where $\|\cdot\|_{j}, j \geq 0,$ denotes a family of weighted $L^{2}$ norms. Namely,
$$
\|f\|_{j}^{2}=\int_{\mathbb{T} \times \mathbb{R}^{+}}|f(x, y)|^{2} \rho_{j}(y) \mathrm{d} x \mathrm{~d} y,
$$
where $\rho_{j}, j \geq 0$, is the family of weights given by
$$
\rho_{0}(y)=(1+y)^{2 m}, \quad \rho_{j}(y)=\frac{\rho_{j-1}(y)}{\left(1+\frac{y}{j^{\alpha}}\right)^{2}}=\rho_{0}(y) \prod_{k=1}^{j}\left(1+\frac{y}{k^{\alpha}}\right)^{-2}, \quad j \geq 1,
$$
for fixed constants $\alpha \geq 0$ and $m \geq 0$ chosen later $(m$ large enough and $\alpha$ matching the constraints found from the estimates). The need for this family of weights will be clarified later. Let us note that locally in $y,$ this family of norms is comparable to more classical families such as
$$
\||f|\|_{\gamma, \tau, r}^{2}=\sum_{j \in \mathbb{N}}\left(\frac{\tau^{j+1}(j+1)^{r}}{(j !)^{\gamma}}\right)^{2}\left\|\partial_{x}^{j} f\right\|_{L^{2}}^{2}.
$$
For instance, for functions $f$ which are zero for $|y| \geq M,$ one has
$$
\|f\|_{\gamma, \tau, r} \leq C_{M}\|\| f\left\|_{\gamma, \tau, r}, \quad\right\||f|\left\|_{\gamma, \tau, r} \leq C_{M, \tau^{\prime}}\right\| f \|_{\gamma, \tau^{\prime}, r} \text { for any } \tau^{\prime}>\tau.
$$
\begin{Theorem}(\cite{dgv})
There exists $m$ and $\alpha$ such that: for all $0<\tau_{1}<\tau_{0}, r \in \mathbb{R},$ for all $T_{0}>0,$ for all $U^{E}$ satisfying
$$
\begin{array}{l}
\sup _{\left[0, T_{0}\right]}\left|\partial_{t} U^{E}\right|_{2, \tau_{0}, r}+\left|U^{E}\right|_{2, \tau_{0}, r}<+\infty ,\\
\sup _{\left[0, T_{0}\right]} \max _{l=0, \ldots, 3}\left\|\partial_{t}^{l}\left(\partial_{t}+U^{E} \partial_{x}\right) U^{E}\right\|_{H^{6-2 l}(\mathbb{T})}<+\infty,
\end{array}
$$
for all $U_{\text {in }}^{P}$ satisfying
$$
\begin{array}{c}
\left\|U_{\text {in }}^{P}-\left.U^{E}\right|_{t=0}\right\|_{2, \tau_{0}, r}<+\infty, \quad\left\|(1+y) \partial_{y} U_{\text {in }}^{P}\right\|_{2, \tau_{0}, r}<+\infty, \\
\left\|(1+y)^{m+6} \partial_{y} U_{\text {in }}^{P}\right\|_{H^{6}\left(\mathbb{T} \times \mathbb{R}_{+}\right)}<+\infty,
\end{array}
$$
and under usual compatibility conditions (see the last remark below), there exists $0<T \leq T_{0}$ and a unique solution $U^{P}$ of (\ref{y2.13})-(\ref{y2.14})  over $(0, T)$ with initial datum $U_{\text {in }}^{P}$ that satisfies
$$
\begin{array}{l}
\sup _{t \in[0, T]}\left\|U^{P}(t)-U^{E}(t)\right\|_{2, \tau_{1}, r}^{2}+\sup _{t \in[0, T]}\left\|(1+y) \partial_{y} U^{P}(t)\right\|_{2, \tau_{1}, r}^{2} \\
+\int_{0}^{T}\left\|(1+y) \partial_{y}^{2} U^{P}(t)\right\|_{2, \tau_{1}, r}^{2} \mathrm{~d} t<+\infty.
\end{array}
$$
\end{Theorem}

The authors in \cite{HH,LCS,WZ1} studied the local well-posedness of the Prandtl equations around others.

Hong and Hunter (\cite{HH}) in 2003 used the method of characteristics to prove the short-time existence of smooth solutions of the unsteady inviscid Prandtl equations, and presented a simple explicit solution that forms a singularity in a finite time.

Lombardo,  Cannone and   Sammartino (\cite{LCS}) in 2003 considered the mild solutions of the following Prandtl equations \eqref{h18}-\eqref{h19} on the half space and proved the short time existence and the uniqueness of the solution in the proper function space under the analyticity of the tangential variable by the  abstract Cauchy-Kowalewski theorem.

\begin{equation}\left\{\begin{array}{ll}
(\partial_{t}-\partial_{YY})u^{P}+u^{P}\partial_{x}u^{P}+v^{P}\partial_{Y}u^{P}+\partial_{x}p^{P}=0,\\
\partial_{Y}p^{P}=0,\\
\partial_{x}u^{P}+\partial_{Y}v^{P}=0,\\
u^{P}(x,Y=0,t)=v^{P}(x,Y=0,t)=0,\\
u^{P}(x,Y\rightarrow\infty,t)\rightarrow U(x,t),\\
p^{P}(x,Y\rightarrow\infty,t)\rightarrow p^{E}(x,y=0,t),\\
u^{P}(x,Y,t=0)=u_{in}^{p},
\end{array}\right.\label{h18}\end{equation}

\begin{equation}
\begin{array}{ll}
u+F(t,u)=0.
\end{array}\label{h19}
\end{equation}

\begin{Theorem}(\cite{LCS})
Suppose that there exist $ R>0$, $\rho_{0}>0$, and $\beta_{0}>0$ such
that if $0<t\leq \rho_{0}/\beta_{0}$, the following properties hold:\\
(1) $
\forall 0<\rho'<\rho\leq\rho_{0}$ and $\forall u$ such that $\{u\in X_{\rho} : \sup_{0\leq t\leq T}|u(t)|_{\rho}\leq R\}$, the map $F(t,u): [0,T]\mapsto X_{\rho'}$
 is continuous.\\
 (2) $\forall 0<\rho\leq\rho_{0}$, the function
$ F(t,0) : [0,\rho_{0}/\beta_{0}]\mapsto\{u\in X_{\rho}: \sup_{0\leq t\leq T}|u(t)|_{\rho}\leq R\} $ is continuous and
\begin{eqnarray}
|F(t,0)|_{\rho}\leq R_{0} <R.
\end{eqnarray}
(3) $\forall 0<\rho'<\rho(s)\leq \rho_{0}$ and $\forall u^{1}$ and $u^{2}\in \{u\in X_{\rho}: \sup_{0\leq t\leq T}|u(t)|_{\rho-\beta_{0}t}\leq R\}$,
\begin{eqnarray}
|F(t,u^{1})-F(t,u^{2})|_{\rho'}\leq C\int_{0}^{t}ds(\frac{|u^{1}-u^{1}|_{\rho(s)}}{\rho(s)=\rho'}+\frac{|u^{1}-u^{2}|_{\rho'}}{\sqrt{t-s}}).
\end{eqnarray}
Then  there exists $ \beta >\beta_{0}$ such that $\forall 0<\rho <\rho_{0}$, \eqref{h19} has a unique solution $u(t)\in X_{\rho}$ with $t\in [0, (\rho_{0}-\rho)/\beta]$.
Moreover, $\sup_{\rho<\rho_{0}-\beta t}|u(t)|_{\rho}\leq R$. Here $ \{X_{\rho}:0<\rho\leq \rho_{0} \} $ be a Banach scale with norms $\|\cdot \|_{\rho}$.
\end{Theorem}

\begin{Theorem}(\cite{LCS})
Suppose that $u_{1}$ and $u_{2}$ are in $K_{\beta_{0},T}^{l,\rho_{0},\mu_{0}}$ (see \eqref{-3}). Suppose $0<\rho'<\rho(s)<\rho_{0}s$ and $0<\mu'<\mu(s)<\mu_{0}$. Then the following estimate holds:
\begin{eqnarray}
|F(u^{1},t)-F(u^{2},t)|_{l,\rho',\mu'}
\leq c\int_{0}^{t}ds(\frac{|u^{1}-u^{2}|_{l,\rho(s),\mu}}{\rho(s)-\rho'}+\frac{|u^{1}-u^{2}|_{l,\rho,\mu(s)}}{\mu(s)-\mu'}+\frac{|u^{1}-u^{2}|_{l,\rho',\mu'}}{\sqrt{t-s}}).
\end{eqnarray}
\end{Theorem}

\begin{Theorem}(\cite{LCS})
Suppose $U\in K^{l,\rho_{0}}_{\beta_{0},T}$ and $u_{in}^{P}-U \in K^{l,\rho_{0},\mu_{0}}$ . Moreover, let the
compatibility conditions
\begin{eqnarray}
u_{in}^{P}(x,Y=0)=0,\\
u_{in}^{P}(x,Y\rightarrow\infty)-U \rightarrow 0,
\end{eqnarray}
hold. Then there exist $0<\rho_{1}<\rho_{0}$, $0<\mu_{1}<\mu_{0}$, $\beta_{1}>\beta_{0}>0$, and $0<T<T_{1}$ such
that \eqref{h18} admit a unique mild solution $u^{p}$ . This solution can be written as
\begin{eqnarray}
u^{P}(x,Y,t)=u(x,Y,t)+U,
\end{eqnarray}
where $u\in K^{l,\rho_{1},\mu_{1}}_{\beta_{1},T_{1}}$. Here
$$ |f |_{K^{l,\rho }_{\beta ,T} }=\sum\limits_{0\leq j\leq 1}\sum\limits_{i\leq l-j }\sup\limits_{0\leq t\leq T}|\partial_{t}^{j}\partial_{x}^{i} f(\cdot,t)|_{0,\rho-\beta t}<+\infty .$$
\end{Theorem}

Wang and Zhu  (\cite{WZ1}) in 2020 studied of the existence of a backflow point in the two-dimensional unsteady boundary layers under an adverse pressure gradient,and considered the following problem for the Prandtl boundary layer equation with nonslip boundary condition for an unsteady incompressible flow in a domain $Q_{T}=\{(t, x, y) \mid 0 \leq t<T, 0 \leq x \leq L, 0 \leq y<+\infty\}:$
\begin{equation}
\left\{\begin{array}{l}
\partial_{t} u+u \partial_{x} u+v \partial_{y} u=\partial_{y}^{2} u-\partial_{x} P, \\
\partial_{x} u+\partial_{y} v=0, \\
\left.u\right|_{y=0}=\left.v\right|_{y=0}=0, \quad \lim\limits _{y \rightarrow+\infty} u=U_{e}(t, x), \\
\left.u\right|_{t=0}=u_{0}(x, y),\left.u\right|_{x=0}=u_{1}(t, y),
\end{array}\right.\label{+2}
\end{equation}
where $(u(t, x, y), v(t, x, y))$ is the velocity field in the boundary layer, $U_{e}(t, x)$ and $P(t, x)$ are traces at the boundary $\{y=0\}$ of the tangential velocity and pressure of the Euler outer flow, respectively, interrelated through Bernoulli's law
$$
\partial_{t} U_{e}+U_{e} \partial_{x} U_{e}=-\partial_{x} P.
$$
Assume that the data satisfy the following properties
$$
u_{0}(x, y)>0, u_{1}(t, y)>0, \forall t \in[0, T), x \in[0, L], y \in[0,+\infty),
$$
and the monotonicity assumption,
$$
\partial_{y} u_{0}(x, y)>0, \partial_{y} u_{1}(t, y)>0, \forall t \in[0, T), x \in[0, L], y \in[0,+\infty).
$$
\begin{Theorem}(\cite{WZ1})
(1) Assume that the trace at the boundary of the outer Euler flow satisfies
$$
U_{e} \in C^{1}([0, T] \times[0, L]) \text { and } U_{e}(t, x)>0 \quad \forall t \in[0, T], x \in[0, L],
$$
and the uniformly adverse pressure gradient in the sense that
$$
\partial_{x} P(t, x)>0 \quad \forall t \in[0, T], x \in[0, L].
$$
Assume that $(u, v) \in C^{2}\left(Q_{T}\right)$ is the solution to the problem \eqref{+2} corresponding to the data $u_{0}, u_{1}$ satisfying the conditions above. Then, the first zero point of $\partial_{y} u(t, x, y),$ when the time evolves, should be at the boundary $\{y=0\}$ if it exists for some time $t \in(0, T)$.\\
(2) Moreover, when the initial velocity $u_{0}(x, y)$ satisfies
$$
\int_{0}^{+\infty} \int_{0}^{L} \frac{(L-x)^{\frac{3}{2}} \partial_{y} u_{0}}{\sqrt{\left(\partial_{y} u_{0}\right)^{2}+u_{0}^{2}}} d x d y \geq C_{*},
$$
for a positive constant $C_{*}$ depending only on $L, T, U_{e},$ and $\partial_{x} P,$ then there is a backflow point $\left(t^{*}, x^{*}\right) \in(0, T) \times[0, L]$ such that
$$
\left\{\begin{array}{l}
\partial_{y} u\left(t^{*}, x^{*}, 0\right)=0 ,\\
\partial_{y} u(t, x, y)>0 \quad \forall 0<t<t^{*}, x \in[0, L], y \geq 0.
\end{array}\right.
$$
Moreover, they have $\partial_{y}^{2} u\left(t^{*}, x^{*}, 0\right) \neq 0$.
\end{Theorem}

Li and Yang  (\cite{ly1}) in 2020 studied the well-posedness of the Prandtl system without monotonicity and analyticity assumption. Precisely, for any index $\sigma \in[3 / 2,2]$, they obtained the local in time well-posedness in the space of Gevrey class $G^{\sigma}$ in the tangential variable and Sobolev class in the normal variable so that the monotonicity condition on the tangential velocity is not needed to overcome the loss of tangential derivative.

The Prandtl equations introduced by Prandtl in 1904 described the behavior of the incompressible flow near a rigid wall at high Reynolds number:
\begin{equation}
\left\{\begin{array}{l}
\partial_{t} u^{P}+u^{P} \partial_{x} u^{P}+v^{P} \partial_{y} u^{P}-\partial_{y}^{2} u^{P}+\partial_{x} p=0, \quad t>0, \quad x \in \mathbb{R}, \quad y>0, \\
\partial_{x} u^{P}+\partial_{y} v^{P}=0 ,\\
\left.u^{P}\right|_{y=0}=\left.v^{P}\right|_{y=0}=0, \quad \lim\limits _{y \rightarrow+\infty} u=U(t, x), \\
\left.u^{P}\right|_{t=0}=u_{0}^{P}(x, y),
\end{array}\right.\label{+3}
\end{equation}
where $u^{P}(t, x, y)$ and $v^{P}(t, x, y)$ represent the tangential and normal velocities of the boundary layer, with $y$ being the scaled normal variable to the boundary, while $U(t, x)$ and $p(t, x)$ are the values on the boundary of the tangential velocity and pressure of the outflow satisfying the Bernoulli's law
$$
\partial_{t} U+U \partial_{x} U+\partial_{x} p=0.
$$

To have a clear presentation, the authors of \cite{ly1} constructed a solutions $u^{P}$ that is a small perturbation around a shear flow, that is, $u^{P}(t, x, y)=u^{s}(t, y)+u(t, x, y)$. For this, they supposed that the initial data $u_{0}^{P}$ can be written as
$$
u_{0}^{P}(x, y)=u_{0}^{s}(y)+u_{0}(x, y),
$$
with $u_{0}^{s}$ being independent of $x$ variable. Then one reduces the original Prandtl equation to the following two time evolutional equations, one of which is the equation for the shear flow $\left(u^{s}, 0\right)$ with $u^{s}$ solving
\begin{equation}
\left\{\begin{array}{l}
\partial_{t} u^{s}-\partial_{y}^{2} u^{s}=0 ,\\
\left.u^{s}\right|_{y=0}=0, \quad \lim\limits _{y \rightarrow+\infty} u^{s}=1 ,\\
\left.u^{s}\right|_{t=0}=u_{0}^{s},
\end{array}\right.
\end{equation}
and the another reads
\begin{equation}
\left\{\begin{array}{l}
\partial_{t} u+\left(u^{s}+u\right) \partial_{x} u+v \partial_{y}\left(u^{s}+u\right)-\partial_{y}^{2} u=0, \\
\left.u\right|_{y=0}=0, \quad \lim _{y \rightarrow+\infty} u=0, \\
\left.u\right|_{t=0}=u_{0},
\end{array}\right.
\end{equation}
where
$$
v=-\int_{0}^{y} \partial_{x} u(x, \tilde{y}) d \tilde{y} .
$$
There exists a $y_{0}>0$ such that $u_{0}^{s} \in C^{6}\left(\mathbb{R}_{+}\right)$ satisfies the following properties

(i) $\frac{d u_{0}^{s}}{d y}\left(y_{0}\right)=0$ and $\frac{d^{2} u_{0}^{s}}{d y^{2}}\left(y_{0}\right) \neq 0 .$ Moreover, there exist $0<\delta<y_{0} / 2$ and a constant $c_{0}>0$ such that
$$
\forall y \in\left[y_{0}-2 \delta, y_{0}+2 \delta\right], \quad\left|\frac{d^{2} u_{0}^{s}}{d y^{2}}(y)\right| \geq c_{0}.
$$

(ii) There exists a constant $0<c_{1}<1$ such that
$$
\forall y \in\left[0, y_{0}-\delta\right] \cup\left[y_{0}+\delta,+\infty\left[, \quad c_{1}\langle y\rangle^{-\alpha} \leq\left|\frac{d u_{0}^{s}(y)}{d y}\right| \leq c_{1}^{-1}\langle y\rangle^{-\alpha}, \right.\right.
$$
for some $\alpha>1,$ and that
$$
\forall y \geq 0, \quad\left|\frac{d^{j} u_{0}^{s}(y)}{d y^{j}}\right| \leq c_{1}^{-1}\langle y\rangle^{-\alpha-1} \text { for } 2 \leq j \leq 6.
$$

(iii) The compatibility condition holds, that is, $\left.u_{0}^{s}\right|_{y=0}=\left.\partial_{y}^{2} u_{0}^{s}\right|_{y=0}=0$ and $u_{0}^{s}(y) \rightarrow 1$ as $y \rightarrow+\infty$.

Let $\alpha$ be the number given in the assumption and let $\ell$ be a fixed number satisfying that
$$
\ell>3 / 2, \quad \alpha \leq \ell<\alpha+\frac{1}{2}.
$$
With each pair $(\rho, \sigma), \rho>0, \sigma \geq 1,$ we associate a Banach space $X_{\rho, \sigma},$ equipped with the norm $\|\cdot\|_{\rho, \sigma}$ that consists of all the smooth functions $f$ such that $\|f\|_{\rho, \sigma}<+\infty,$ where
$$
\begin{aligned}
\|f\|_{\rho, \sigma} \stackrel{\text { def }}{=} & \sup _{m \geq 6} \frac{\rho^{m-5}}{[(m-6) !]^{\sigma}}\left\|\langle y\rangle^{\ell-1} \partial_{x}^{m} f\right\|_{L^{2}}+\sup _{m \geq 6} \frac{\rho^{m-5}}{[(m-6) !]^{\sigma}}\left\|\langle y\rangle^{\ell} \partial_{x}^{m}\left(\partial_{y} f\right)\right\|_{L^{2}} \\
&+\sup _{0 \leq m \leq 5}\left(\left\|\langle y\rangle^{\ell-1} \partial_{x}^{m} f\right\|_{L^{2}}+\left\|\langle y\rangle^{\ell} \partial_{x}^{m}\left(\partial_{y} f\right)\right\|_{L^{2}}\right) \\
&+\sup _{1 \leq j \leq 4 \atop i+j \geq 6} \frac{\rho^{i+j-5}}{[(i+j-6) !]^{\sigma}}\left\|\langle y\rangle^{\ell+1} \partial_{x}^{i} \partial_{y}^{j}\left(\partial_{y} f\right)\right\|_{L^{2}}+\sup _{1 \leq j \leq 4 \atop i+j \leq 5}\left\|\langle y\rangle^{\ell+1} \partial_{x}^{i} \partial_{y}^{j}\left(\partial_{y} f\right)\right\|_{L^{2}} .
\end{aligned}
$$

\begin{Theorem}(\cite{ly1})
For $\sigma \in[3 / 2,2],$ let the initial datum $u_{0}$ belong to $X_{2 \rho_{0}, \sigma}$ for some $\rho_{0}>0$ and moreover
$$
\left\|u_{0}\right\|_{2 \rho_{0}, \sigma} \leq \eta_{0},
$$
for some $\eta_{0}>0$. Suppose that the compatibility condition holds for $u_{0}$. Then the equation \eqref{+3} admits a unique solution $u \in L^{\infty}\left([0, T] ; X_{\rho, \sigma}\right)$ for some $T>0$ and some $0<\rho<2 \rho_{0},$ provided $\eta_{0}$ is sufficiently small.
\end{Theorem}

\subsection{2D Prandtl Equations-Global Existence }

 Oleinik and  Samokhin (\cite{OS}) in 1999 considered
the following Prandtl system for a nonstationary boundary layer in an axially symmetric incompressible flow past a solid body
\begin{equation}\left\{
\begin{array}{ll}
\frac{\partial u}{\partial t}+ u\frac{\partial u}{\partial x}+v\frac{\partial u}{\partial y}
=  \nu\frac{\partial^{2} u}{\partial y^{2}} +U\frac{\partial U}{\partial x}+U \frac{dU}{d x},  \\
\frac{\partial (ru)}{\partial x}+\frac{\partial (rv)}{\partial y} =0,
\end{array}
 \label{ww21}         \right.\end{equation}
in a domain $D=\{0<t<T,0<x<X, 0<y<\infty\}$, with the boundary conditions
\begin{equation}\left\{
\begin{array}{ll}
u(0,x,y)=u_{0}(x,y), \ \ u(t,0,y)=0,\ \ u(t,x,0)=0, \ \  v(t, x,0) =v_{0}(t,x) ,  \\
u(t,x,y)\rightarrow U(t,x)  \ \  \text{as}\ \  y\rightarrow \infty.
\end{array}
 \label{ww22}         \right.\end{equation}
The function $U(t, x)$ and $r(x)$ are given functions such that $U(t,0)=0, U(t,x)>0$ for $x\geq 0$,   $r(0)=0,  r(x)>0$ for $x>0$.

Introduce  Crocco variables by
$$\tau =t, \ \  \xi=x, \ \ \eta=\frac{u(t,x,y)}{U(t,x)},$$
and a new unknown function
$$w(\tau, \xi,\eta)= \frac{u_{y}(t, x,y)}{U(t,x)}.$$
From  (\ref{ww21}), we obtain the following equation for $w(\tau, \xi,\eta)$:
\begin{eqnarray}
\nu w^{2}w_{\eta\eta}-w_{\tau}-\eta Uw_{\xi}+Aw_{\eta}+Bw=0
 \label{ww23}
\end{eqnarray}
in the domain $\Omega=\{0<\tau<T, 0<\xi<X, 0<\eta <1 \}$, with the boundary conditions
\begin{eqnarray}
w|_{\tau=0}=\frac{u_{0y}}{U}\equiv w_{0}(\xi,\eta), \ \  w|_{\eta=1}=0, \ \  (\nu w^{2}w_{\eta\eta}-v_{0}w+C)_{\eta=0}=0,
 \label{ww24}
\end{eqnarray}
where
$$A=(\eta^{2}-1)U_{x}+(\eta-1)\frac{U_{t}}{U}, \ \  B=-\eta U_{x}(\xi)+\frac{\eta r_{x}(\xi)U(\xi)}{r(\xi)}-\frac{U_{t}}{U} , \ \  C=U_{x}+\frac{U_{t}}{U}.$$

\begin{Theorem} \label{OS.2.26} (\cite{OS})
Assume that $U_{x}, U_{t}/U, Ur_{x}/r, v_{0}$ are bounded functions having bounded derivatives with respect to $t,x$ in $D$; $u_{0}(x,y)\rightarrow \infty$ as
$y\rightarrow \infty$, $u_{0}=0$ for $y=0$; $u_{0}/U, u_{0y}/U$ are continuous in $\overline{D}$; $u_{0y}>0$ for $y\geq 0,x>0$,
$$K_{1}(U(0,x)-u_{0}(x,y))\leq u_{0y}(x,y)\leq K_{2}(U(0,x)-u_{0}(x,y)) ,$$
with positive constants $K_{1}$ and $K_{2}$. Assume also that there exist bounded derivatives $u_{0y},u_{0yy},u_{0yyy},u_{0x},u_{0xy}$, and the ratios
$$\frac{u_{0yy}}{u_{0y}}, \ \ \frac{u_{0yyy}u_{0y}-u_{0yy}^{2}}{u_{0y}^{2}}$$
are bounded for $0\leq x\leq X, 0\leq y<\infty$. Let the following compatibility condition be satisfied:
$$v_{0}(0,x)u_{0y}(x,0)=-p_{x}(0,x)+\nu u_{0yy}(x,0), $$
and let
$$\left|\frac{u_{0yx}-u_{0x}u_{0yy}}{u_{0y}}+U_{x}\frac{u_{0}u_{0yy}-u_{0y}^{2}}{Uu_{0y}} \right|\leq K_{6}(U-u_{0}(x,y)).$$
Then, problem (\ref{ww21})-(\ref{ww22}) in $D$ has a unique solution $u,v$ with the following properties: $u/U,u_{y}/U$ are continuous and bounded
in $\overline{D}$;  $u_{y}/U>0$ for $y\geq 0$; $u_{y}/U\rightarrow 0$ as $y\rightarrow \infty$; $u=0$ for $y=0$;  $v$ is continuous in $D$ and
bounded for bounded $y$; the weak derivatives $u_{t},u_{x}, u_{yt}, u_{yx}, u_{yy}, u_{yyy}, u_{y}$ are bounded measurable functions in $D$;
the equations of system (\ref{ww21}) hold almost everywhere; in $D$, the functions $u_{t},u_{x}, v_{y}, u_{yy} $ are continuous
with respect to $y$. Moreover,
$$\frac{u_{ yy}}{u_{ y}}, \ \ \frac{u_{ yyy}u_{ y}-u_{ yy}^{2}}{u_{ y}^{2}}$$
are bounded and the following inequalities hold:
\begin{eqnarray*}
&& C_{1}(U(t,x)-u(t,x,y))\leq u_{ y}(t,x,y)\leq C_{2}(U(t,x)-u(t,x,y)) , \\
&& \exp(-C_{2}y)\leq 1-\frac{u(t,x,y)}{U(t,x)} \leq \exp(-C_{1}y), \\
&& \left|\frac{u_{yt}u_{y}-u_{t}u_{ yy}}{u_{ y}}+U_{t}\frac{uu_{ yy}-u_{y}^{2}}{u_{ y}U} \right|\leq C_{4}(U(t,x)-u(t,x,y)) , \\
&& \left|\frac{u_{yx}u_{y}-u_{x}u_{ yy}}{u_{ y}}+U_{x}\frac{uu_{ yy}-u_{y}^{2}}{u_{ y}U} \right|\leq C_{3}(U(t,x)-u(t,x,y)).
\end{eqnarray*}

\end{Theorem}

 Oleinik and  Samokhin (\cite{OS}) in 1999 also considered
system  (\ref{ww21}) in the domain $D=\{0<t<T,0<x<X, 0<y<\infty\}$, with the boundary and the initial conditions of the form:
\begin{equation}\left\{
\begin{array}{ll}
u(0,x,y)=u_{0}(x,y), \ \ u(t,0,y)=u_{1}(t,y),\ \ u(t,x,0)=0, \ \  v(t, x,0) =v_{0}(t,x) ,  \\
u(t,x,y)\rightarrow U(t,x)  \ \  \text{as}\ \  y\rightarrow \infty.
\end{array}
 \label{ww25}         \right.\end{equation}

\begin{Theorem} \label{OS.2.27} (\cite{OS})
Suppose that: $U_{x}, U_{t}, v_{0}, U, r_{x}/r$ are bounded and have bounded first derivatives in $D$ with respect to $t$ and $x$;
$$u_{0}(x,y)\rightarrow U(0,x), \ \ u_{1}(t,y)\rightarrow U(t,0), \ \  \text{as} \ \  y\rightarrow \infty, $$
and $U(t,x)>0, u_{0y}>0, u_{1y}>0$ for $y\geq 0$,
\begin{eqnarray*}
&& K_{1}(U(0,x)-u_{0}(x,y))\leq u_{0y}(x,y)\leq K_{2}(U(0,x)-u_{0}(x,y)) , \\
&& K_{1}^{1}(U(t,0)-u_{1}(t,y))\leq u_{1y}(x,y)\leq K_{2}^{1}(U(t,0)-u_{1}(x,y)) .
\end{eqnarray*}
Suppose also that there exist bounded weak derivatives $u_{0y}, u_{0yy}, u_{0yyy}$, $u_{0x},u_{0xy}, u_{1y},u_{1yy}, u_{1yyy}$, $u_{1t}, u_{1ty}$ and the
 functions
$$\frac{u_{i yy}}{u_{ iy}}, \ \ \frac{u_{i yyy}u_{i y}-u_{i yy}^{2}}{u_{i y}^{2}}, \ \   i=0,1,$$
are bounded; moreover,
\begin{eqnarray*}
&& \left|\frac{u_{0yx}u_{0y}-u_{0x}u_{0yy}}{u_{0y}}+U_{x}\frac{u_{0}u_{0yy}-u_{0y}^{2}}{Uu_{0y}} \right|\leq K_{5}(U-u_{0} ), \\
&& \left|\frac{u_{1yx}u_{1y}-u_{1x}u_{1yy}}{u_{1y}}+U_{t}\frac{u_{1}u_{1yy}-u_{1y}^{2}}{Uu_{1y}} \right|\leq K_{5}^{1}(U-u_{1} ) ;
\end{eqnarray*}
the function whose absolute value is estimated in the last inequality is continuous in $t$ and $y$ for small $y$. Let the following compatibility conditions
hold:
\begin{eqnarray*}
&& v_{0}(0,x)u_{0y}(x,0)=-p_{x}(0,x)+\nu u_{0yy}(x,0), \\
&& v_{0}(t,0)u_{1y}(t,0)=-p_{x}(t,0)+\nu u_{1yy}(t,0), \\
&& u_{0}=u_{1} \ \ \text{for } \ \  t=0, x=0; \ \  u_{1}=0, u_{0}=0 \ \ \text{for } \ \  y=0; \\
&& u_{1ty}(t,y)+v_{0}(t,0)u_{1yy}(t,y)-\nu u_{1yyy}(t,y)=\mathcal{O}(y)  \ \ \text{for small } \ \ y.
\end{eqnarray*}
Then problem (\ref{ww21}), (\ref{ww25}) in $D$  has a unique solution  $u,v$ with the following properties: $u,u_y$ are continuous and bounded in
 $\overline{D}$;  $u_{y} >0$ for $y\geq 0$; $u_{y} \rightarrow 0$ as $y\rightarrow \infty$; $u=0$ for $y=0$;  $v$ is continuous in $\overline{D}$ with respect
to $y$, and bounded for bounded $y$; there exist bounded weak derivatives $u_{t},u_{x}, u_{yt}, u_{yx}, u_{yy}, u_{yyy}, v_{y}$ and equations (\ref{ww21})
 hold almost everywhere in $D$. Moreover, $u_{t},u_{x}, v_{y}, u_{yy} $ are continuous in $y$; $u_{yy}/u_{y}, (u_{yyy}u_{y}-u_{yy}^{2})/u_{y}^{2}$
are bounded, and
\begin{eqnarray*}
&& C_{1}(U(t,x)-u(t,x,y))\leq u_{ y}(t,x,y)\leq C_{2}(U(t,x)-u(t,x,y)) , \\
&& \exp(-C_{2}y)\leq 1-\frac{u(t,x,y)}{U(t,x)} \leq \exp(-C_{1}y), \\
&& \left|\frac{u_{yt}u_{y}-u_{t}u_{ yy}}{u_{ y}}+U_{t}\frac{uu_{ yy}-u_{y}^{2}}{u_{ y}U} \right|\leq C_{3}(U  -u  ) , \\
&& \left|\frac{u_{yx}u_{y}-u_{x}u_{ yy}}{u_{ y}}+U_{x}\frac{uu_{ yy}-u_{y}^{2}}{u_{ y}U} \right|\leq C_{4}(U -u ).
\end{eqnarray*}
where $C_{i}$ $(i=1,2,3,4)$ are positive constants.

\end{Theorem}

\begin{Theorem} \label{OS.2.28} (\cite{OS})
Let the assumptions of Theorem \ref{OS.2.27} hold in the domain $D=\{0<t<\infty,0<x<X, 0<y<\infty\}$. Then the solution $u,v$ of problem (\ref{ww21}), (\ref{ww25})
 exists and is unique in the domain $D$ with $X>0$ depending on the data of problem (\ref{ww21}), (\ref{ww25}). This solution has
the same properties as the solution $u,v$ constructed in Theorem \ref{OS.2.27}.

\end{Theorem}

 Oleinik and  Samokhin (\cite{OS}) in 1999 also considered
  the nonstationary system of the two-dimensional boundary layer
\begin{equation}\left\{
\begin{array}{ll}
\frac{\partial u}{\partial t}+ u\frac{\partial u}{\partial x}+v\frac{\partial u}{\partial y}
=  \nu\frac{\partial^{2} u}{\partial y^{2}} -\frac{dp}{dx},  \\
\frac{\partial u}{\partial x}+\frac{\partial v}{\partial y} =0,
\end{array}
 \label{ww26}         \right.\end{equation}
in a domain $D=\{0<t<T,0<x<X, 0<y<\infty\}$, with the boundary conditions
\begin{equation}\left\{
\begin{array}{ll}
u(0,x,y)=u_{0}(x,y), \ \ u(t,0,y)=u_{1}(t,y),\ \ u(t,x,0)=0, \ \  v(t, x,0) =v_{0}(t,x) ,  \\
u(t,x,y)\rightarrow U(t,x)  \ \  \text{as}\ \  y\rightarrow \infty.
\end{array}
 \label{ww27}         \right.\end{equation}
The density of the fluid $\rho$ is assumed constant and equal to $1$. Therefore, $U_{t}+UU_{x}=-p_{x}(t,x)$. It is also assumed that
$U(t,x)>0, u_{0}>0, u_{1}>0$ for $y>0$, $u_{0y}>0, u_{1y}>0$ for $y\geq 0$.

\begin{Theorem} \label{OS.2.29} (\cite{OS})
Assume that $p(t,x), v_{0}(t,x), u_{0}(x,y), u_{1}(t,y), u_{0y}$ are sufficiency smooth and satisfy the compatibility conditions
which amount to the existence of the function to $w^{*}$ mentioned earlier. Then there is one and only one solution of problem (\ref{ww26}), (\ref{ww27})
in the domain $D$, with $X$ being arbitrary and $T$ depending on the data of problem  (\ref{ww26}), (\ref{ww27}), or $T$  being arbitrary and $X$ depending
on the data. This solution has the following properties: $u>0$ for $y>0$, $u_{y}>0$ for $y\geq 0$;  the derivatives $u_{t}, u_{x}, u_{y}, u_{yy}, v_{y}$
are continuous and bounded in $\overline{D}$; moreover, the ratios
$$\frac{u_{yy}}{u_{y}}, \ \  \frac{u_{yyy}u_{y}-u_{yy}^{2}}{u_{y}^{2}}$$
are bounded in D.

\end{Theorem}

\begin{Theorem} \label{OS.2.30} (\cite{OS})
(Global Solutions of the Prandtl System for Axially Symmetric Flows) Let $U(t,x)=ax+b(t,x)x^{3}$, where $a$ is a positive constant,
and $b(x,t)$ has bounded second order derivatives. Assume also that
$$v_{0}\leq M_{1}x, \ \ v_{0t}\geq -M_{2} x, \ \ \left|\left(U\frac{r_{x}}{r}- U_{x}\right)_{t}\right|\leq M_{3}x. $$
Assume that $u_{0}(x,y)$ satisfies the conditions: $u_{0}/U$ is continuous in $\overline{D}$, $u_{0}>0$ for $y>0,x>0$; $u_{0}(x,0)=0, u_{0}/U\rightarrow 1$
as $y\rightarrow \infty$; $u_{0y}/U>0$ for $y\geq 0$; $u_{0y}/U\rightarrow 0$ as $y\rightarrow \infty$;  moreover,
\begin{eqnarray*}
&& M_{4}(U-u_{0})\sigma_{0}\leq u_{0y}\leq  M_{5}(U-u_{0})\sigma_{0}, \\
&& -M_{6}\sigma_{0}\leq \frac{u_{0yy}}{u_{0y}}\leq   -M_{7}\sigma_{0}, \\
&& \left|\frac{u_{0yx}u_{0y}-u_{0x}u_{0 yy}}{u_{ 0y}}+\frac{U_{x}}{U}\frac{u_{0}u_{0 yy}-u_{0y}^{2}}{u_{ 0y}} \right|\leq M_{8}\sigma_{0}(U -u _{0}), \\
&& \left|\frac{u_{0yyy}u_{0y}- u_{0 yy}^{2}}{u_{ 0y}^{2}}  \right|\leq M_{9} , \\
&& \nu \frac{u_{0yyy}u_{0y}- u_{0 yy}^{2}}{u_{ 0y}}+\left[u_{0}^{2}-U^{2}\frac{U_{x}}{U}\frac{  u_{0 yy} }{u_{ 0y}}  \right]\leq  M_{10}x u_{ 0y},
\end{eqnarray*}
where $\sigma_{0}=\sqrt{-\ln \mu(1-u_{0}/U)},M_{i}\ (i=1,2,3,4,5,6,7,8,9,10),\mu$ are positive constants. Finally, let the following compatibility conditions hold:
$$v_{0}(0,x)u_{0y}(x, 0)=U_{t}(0,x)+U(0,x)U_{x}(0,x)+\nu u_{0yy}(x,0). $$

Then in the domain $D=\{0<t<\infty,0<x<X, 0<y<\infty\}$, problem (\ref{ww21})-(\ref{ww22}) has one and only one solution $u,v$ with the following properties:
$u/U, u_{y}/U$ are bounded and continuous in $\overline{D}$, $u>0$ for $y>0,x>0$; $u\rightarrow U$ as $y\rightarrow \infty$;
$u|_{y=0}=0,u|_{x=0}=0$, $u_{y}/U>0 $ for $y\geq 0$;  $u_{y}/U\rightarrow 0 $ as $y\rightarrow \infty$; $u_{y}, u_{x}, u_{yy}, u_{t}, v_{y}$ are bounded
 and continuous in $\overline{D}$ with respect to $y$; $u_{yy}/u_{y}, v$ are continuous with respect to $y$ in $\overline{D}$ and
are bounded for bounded $y$; $v|_{y=0}=v_{0}(t,x), u_{yyy}$ is bounded in $D$; $u_{yx}, u_{yt}$
are bounded for bounded $y$;  equations (\ref{ww21}) hold almost everywhere in $D$. Moreover,
\begin{eqnarray*}
&& M_{12}(U-u )\sigma \leq u_{ y}\leq  M_{11}(U-u )\sigma , \\
&& -M_{14}\sigma \leq \frac{u_{ yy}}{u_{ y}}\leq   -M_{13}\sigma ,\ \ \left|\frac{u_{ yyy}u_{ y}- u_{ yy}^{2}}{u_{ y}^{2}}  \right|\leq M_{15} ,   \\
&& \left|\frac{u_{ yx}u_{ y}-u_{ x}u_{ yy}}{u_{ y}}+\frac{U_{x}}{U}\frac{u u_{ yy}-u_{ y}^{2}}{u_{ y}} \right|\leq M_{16}\sigma (U -u  ), \\
&&-\leq M_{18}\sigma (U -u  ) \leq \frac{u_{ yt}u_{ y}-u_{t}u_{ yy}}{u_{ y}}  +\frac{U_{t}}{U}\frac{u u_{ yy}-u_{ y}^{2}}{u_{ y}}\leq M_{17}x\sigma (U -u  ) , \\
&& M_{21}U\exp(-M_{22}y^{2})\leq U-u\leq  M_{19}U\exp(-M_{20}y^{2}),
\end{eqnarray*}
where $\sigma =\sqrt{-\ln \mu(1-u /U)},M_{i}\ (i=1,2,3\cdot\cdot\cdot21),\mu$ are positive constants; $\mu<1, X>0$ depend on $U,r, v_{0}, u_{0}$.

\end{Theorem}
They have  ensured the existence of a solution for the boundary layer system for $t$ in the interval $0\leq t<\infty$ (Theorem \ref{OS.2.30}).
It is natural to examine the question of stability and asymptotic stability under certain
perturbations of the data and the parameters involved in the problem (\ref{ww21})-(\ref{ww22}) in the domain $D=\{0<t<\infty,0<x<X, 0<y<\infty\}$.
And let $u,v$ correspond to the given functions $U,r,v_{0}, u_{0}$; let $\widetilde{u}(t,x,y),\widetilde{v}(t,x,y)$ be a solution of the same problem
corresponding to $\widetilde{U},\widetilde{r},\widetilde{v}_{0}, \widetilde{u}_{0}$.
\begin{Theorem} \label{OS.2.31} (\cite{OS})
Assume that the functions
 \begin{eqnarray}
U-\widetilde{U}, \ \ U_{x}-\widetilde{U}_{x}, \ \ \frac{U_{t}}{U}-\frac{\widetilde{U}_{t}}{\widetilde{U}}, \ \ v_{0}-\widetilde{v}_{0},
\ \ \frac{r_{x}U}{r}- \frac{\widetilde{r}_{x}\widetilde{U}}{\widetilde{r}},
\label{ww28}
\end{eqnarray}
identically vanish for $t\geq t_{0}>0$, where $t_{0}=$const.$<\infty$. Then, in the domain $D$ for $0\leq y\leq y_{0}<\infty$, the following estimate holds:
$$\left|\frac{u}{U}- \frac{\widetilde{u}}{\widetilde{U}}\right|\leq K_{3} e^{-\alpha_{0}t}, $$
where $\alpha=$const.$>0$,  $K_{3}=$const.$>0$; $K_{3}$ depends on $y_{0}$.

\end{Theorem}

\begin{Theorem} \label{OS.2.32} (\cite{OS})
Assume that the functions (\ref{ww28}) uniformly converge to $0$ in the domain $D$ as $t\rightarrow \infty$. Then, for any $\delta>0$, there is a constant
$K(\delta)>0$ such that
$$\left|\frac{u}{U}- \frac{\widetilde{u}}{\widetilde{U}}\right|\leq \delta+ K(\delta) e^{-\alpha_{0}t} \ \ \text{in} \ \ D. $$

\end{Theorem}

\begin{Theorem} \label{OS.2.33} \cite{OS}
Assume that the functions (\ref{ww28}) have their absolute values $\leq \varepsilon$, and $u_{0}, \widetilde{u}_{0}$ such that $|w-\widetilde{w}|\leq \varepsilon$,
where $w$ and $\widetilde{w}$ are solutions of Corcco variables corresponding to $U,r,v_{0}, u_{0}$, and
$\widetilde{U},\widetilde{r},\widetilde{v}_{0}, \widetilde{u}_{0}$, respectively. Then
$$\left|\frac{u}{U}- \frac{\widetilde{u}}{\widetilde{U}}\right|\leq   K_{4} \varepsilon  \ \ \text{in} \ \ D, $$
for any $t$, where $K_{4}>0$ is a constant depending on $u_{0}, \widetilde{u}_{0}$ and their derivatives.

\end{Theorem}

 Oleinik and  Samokhin (\cite{OS}) in 1999   examined the boundary layer system for a nonstationary axially symmetric flow (\ref{ww21}) in the case of time-periodic outer flow
in the domain $D=\{-\infty<t<\infty,0<x<X, 0<y<\infty\}$,   with the boundary conditions
\begin{equation}\left\{
\begin{array}{ll}
  u(t,0,y)=0,\ \ u(t,x,0)=0, \ \  v(t, x,0) =v_{0}(t,x) ,  \\
u(t,x,y)\rightarrow U(t,x)  \ \  \text{as}\ \  y\rightarrow \infty,
\end{array}
 \label{ww29}         \right.\end{equation}
and the condition of periodicity
\begin{eqnarray}
u(t+T,x,y)=u(t,x,y),
\label{ww30}
\end{eqnarray}
 where $T$ is a given positive constant.

\begin{Theorem} \label{OS.2.34} (\cite{OS})
Let $U(t,x)=ax+b(t,x)x^{2}$, where $a$ is a positive constant, and $b(x,t)$ has bounded second order derivatives. Assume also that
$$v_{0}\leq K_{1}x, \ \ v_{0t}\geq -K_{2} x, \ \ \left|\left(U\frac{r_{x}}{r}- U_{x}\right)_{t}\right|\leq K_{3}x, \ \ r_{x}|_{x=0}>0; $$
$b, v_{0}, r$ are periodic in $t$ with period $T$. Then problem (\ref{ww21}), (\ref{ww29}) in $D=\{0<t<T,0<x<X, 0<y<\infty\}$ has one and only one solution
$u,v$ with the following properties: $u/U, u_{y}/U$ are bounded and continuous in $\overline{D}$, $u>0$ for $y>0,x>0$; $u\rightarrow U$ as $y\rightarrow \infty$;
$u|_{y=0}=0,u|_{x=0}=0$;  $u_{y}/U>0 $ for $y\geq 0$;  $u_{y}/U\rightarrow 0 $ as $y\rightarrow \infty$; $u_{y}, u_{x}, u_{yy}, u_{t}, v_{y}$ are continuous
 and bounded in $\overline{D}$ with respect to $y$; $u_{yy}/u_{y}, v$ are continuous  in $\overline{D}$ with respect to $y$ and
  bounded for bounded $y$; $v|_{y=0}=v_{0}(t,x), u_{yyy}$ is bounded in $D$; $u_{yx}, u_{yt}$
are bounded for bounded $y$;  equations (\ref{ww21}) hold almost everywhere in $D$, and the following inequalities are valid:
\begin{eqnarray*}
&& K_{4}(U-u )\sigma \leq u_{ y}\leq  K_{5}(U-u )\sigma , \\
&& -K_{6}\sigma \leq \frac{u_{ yy}}{u_{ y}}\leq   -K_{7}\sigma ,\ \ \left|\frac{u_{ yyy}u_{ y}- u_{ yy}^{2}}{u_{ y}^{2}}  \right|\leq K_{8} ,   \\
&& \left|\frac{u_{ yx}u_{ y}-u_{ x}u_{ yy}}{u_{ y}}+\frac{U_{x}}{U}\frac{u u_{ yy}-u_{ y}^{2}}{u_{ y}} \right|\leq K_{9}\sigma (U -u  ), \\
&&-\leq K_{10}\sigma (U -u  ) \leq \frac{u_{ yt}u_{ y}-u_{t}u_{ yy}}{u_{ y}}  +\frac{U_{t}}{U}\frac{u u_{ yy}-u_{ y}^{2}}{u_{ y}}\leq K_{11}x\sigma (U -u  ) ,
\end{eqnarray*}
where $\sigma =\sqrt{-\ln \mu(1-u /U)},K_{i}\ (i=1,2,3,4,5,6,7,8,9,10,11),\mu$ are positive constants; $0<\mu<e^{-1/2}, X $ depend on $U,r, v_{0} $.

\end{Theorem}

 Oleinik and  Samokhin (\cite{OS}) in 1999 considered
 the nonstationary system of the two-dimensional boundary layer
\begin{equation}\left\{
\begin{array}{ll}
\frac{\partial u}{\partial t}+ u\frac{\partial u}{\partial x}+v\frac{\partial u}{\partial y}
=  \nu\frac{\partial^{2} u}{\partial y^{2}} + U\frac{\partial U}{\partial x}+ \frac{\partial U}{\partial t},  \\
\frac{\partial u}{\partial x}+\frac{\partial v}{\partial y} =0
\end{array}
 \label{ww31}         \right.\end{equation}
in a domain $D=\{0<t<T,0<x<X, 0<y<\infty\}$, with the boundary conditions
\begin{equation}\left\{
\begin{array}{ll}
u(0,x,y)=0, \ \ u(t,0,y)=0,\ \ u(t,x,0)=0, \ \  v(t, x,0) =v_{0}(t,x) ,  \\
u(t,x,y)\rightarrow U(t,x)  \ \  \text{as}\ \  y\rightarrow \infty.
\end{array}
 \label{ww32}         \right.\end{equation}
 It is also assumed that $U(t,0)=0, U(0,x)=0, U>0$ for $t,x>0$, $$U(t,x)=t^{n}U_{1}(t,x), \ \ n\geq 1, $$
and $U_{1t}/U_{1}$ is a bounded function.

Introduce new independent variables in (\ref{ww31})-(\ref{ww32}),  setting
$$\tau =t^{n-\frac{1}{2}}, \ \  \xi=x, \ \ \eta=\frac{u(t,x,y)}{U(t,x)},$$
and a new unknown function
$$w(\tau, \xi,\eta)= t^{n }\frac{u_{y}(t, x,y)}{U(t,x)}.$$
From  (\ref{ww31}), we obtain the following equation for $w(\tau, \xi,\eta)$:
\begin{eqnarray}
\nu w^{2}w_{\eta\eta}-\tau^{3}(n-\frac{1}{2})w_{\tau}-\eta U\tau^{2N}w_{\xi}+n(n-1)\tau^{2}w_{\eta}+A_{1}\tau^{2N}w_{\eta}+B_{1}\tau^{2N}w=0,
 \label{ww33}
\end{eqnarray}
where $N=1+\frac{1}{2n-1}$,
$$A_{1}=(\eta^{2}-1)U_{x}+(\eta-1)\frac{U_{1t}}{U_{1}}, \ \  B_{1}=-\eta U_{x}(\xi) -\frac{U_{1t}}{U_{1}} ,  $$
in the domain $\Omega=\{0<\tau<T^{n-\frac{1}{2}}, 0<\xi<X, 0<\eta <1 \}$, with the boundary conditions
\begin{eqnarray}
w|_{\tau=0}=0, \ \  w|_{\eta=1}=0, \ \  (\nu w w_{\eta }-v_{0}w\tau^{ N}+n\tau^{2 }+C_{1}\tau^{2N})_{\eta=0}=0,
 \label{ww34}
\end{eqnarray}
where
$$  C_{1}=U_{x}+\frac{U_{1t}}{U_{1}}.$$

\begin{Theorem} \label{OS.5.35} (\cite{OS})
Let
\begin{eqnarray*}
&& U(t,x)=t^{n}U_{1}(t,x), \ \  n\geq 1,\\
&& U(t,0)=0, \ \  U_{1}(t,x)>0, \ \ x>0.
\end{eqnarray*}
Assume that $U_{1x}, u_{1t}/U_{1}, v_{0}$ have continuous first order derivatives in $t$ and $x$.

Then problem (\ref{ww31}),(\ref{ww32}) in $D=\{0<t<T,0<x<X, 0<y<\infty\}$ admits a solution $u,v$ with the following properties: $u/U, u_{y}t^{n}/U$ are
bounded and continuous in $\overline{D}$; $u(t,x,y)>0$ for $t,x>0$; $u_{y}t^{n}/U>0$ for $t>0$, $u_{y}t^{n}/U\rightarrow 0$ as $y \rightarrow 0$;
 the derivatives $u_{y}, u_{x}, u_{yy},u_{t}, v_{y}$ are  bounded and continuous in $D$ with respect to $y$;
$$|u_{y}|\leq E_{1}t^{n-\frac{1}{2}}, \ \  |u_{yy}|\leq E_{2}t^{n-1}, \ \ |u_{t}|\leq E_{3}t^{n-1}, \ \ |u_{x}|\leq E_{4}t^{n }; $$
the function $v$ is continuous in $\overline{D}$ with respect to $y$ and is bounded for bounded $y$; $t^{-n+\frac{1}{2}}u_{yx}, t^{-n+\frac{3}{2}}u_{yt}$
are bounded for bounded $y$,
$$|u_{yyy}|\leq E_{5}t^{n-\frac{3}{2}}, \ \  E_{i}=\text{const.}>0,\ (i=1,2,3,4,5);$$
equations (\ref{ww31}) hold almost everywhere in $D$. Moreover, the following estimates hold for this solution:
\begin{eqnarray*}
&& \Phi^{-1}\left(yt^{-\frac{1}{2}}(1-\alpha t^{ \frac{1}{2}})\right)U\leq u\leq  \Phi^{-1}\left(yt^{-\frac{1}{2}}(1+\beta t^{ \frac{1}{2}})\right)U,
\ \ \alpha,\beta=\text{const.}>0, \\
&& \Phi(\zeta)\equiv\int_{0}^{\zeta}(Y_{0}(s))^{-1}ds  \ \  (\Phi^{-1})\ \ \text{is the inverse of} \ \ \Phi, \\
&& U(1-e^{-\nu_{1}})\leq u\leq U(1-e^{-\nu_{2}}),  \\
&& \quad \nu_{2}=\left[\frac{M_{1}y}{2t^{1/2}(1-\beta t^{1/2})}\right]^{2}+\frac{M_{1}y\sqrt{-\ln \mu}}{t^{1/2}(1-\beta t^{1/2})} ,  \ \
 \nu_{1}=\left[\frac{M_{2}y}{2t^{1/2}(1+\alpha t^{1/2})}\right]^{2}+\frac{M_{2}y\sqrt{-\ln \mu}}{t^{1/2}(1+\alpha t^{1/2})} , \\
&& \nu M_{1}^{2}=1, \ \ \nu M_{2}^{2}=\frac{1}{2}-\delta, \ \ \delta, \varepsilon=\text{const.}>0, \\
&& 1-\frac{u}{U}=\exp\left(-\frac{1}{4\nu t}\left[y^{2}\left(1+\mathcal{O}(\sqrt{t})+ \mathcal{O}(y^{1+\varepsilon} t^{\frac{1-\varepsilon}{2}})
\right)\right]\right), \ \ \text{as} \ \  y\rightarrow \infty, t\rightarrow 0,\\
&& \left|t^{n} \frac{u_{y}}{U}-t^{n-\frac{1}{2}} Y_{0}(\frac{u}{U}) \right|\leq E_{6}t^{n}Y_{0}(\frac{u}{U}), \ \
 \left|t^{n} \frac{u_{yy}}{u_{y}}-t^{n-\frac{1}{2}} Y_{0\eta}(\frac{u}{U}) \right|\leq E_{7}t^{n}Y_{0\eta}(\frac{u}{U}),\\
&& t(u_{yyy}u_{y}-u_{yy}^{2})u_{y}^{-2}\leq -E_{8},  \ \  E_{i}=\text{const.}>0,~ (i=1,2,3,\cdot\cdot\cdot8).
\end{eqnarray*}
The solution $u,v$ of problem (\ref{ww31})-(\ref{ww32}) is unique in the class of functions.

\end{Theorem}

 Oleinik and  Samokhin (\cite{OS}) in 1999 considered
  the nonstationary system of the two-dimensional boundary layer
\begin{equation}\left\{
\begin{array}{ll}
\frac{\partial u}{\partial t}+ u\frac{\partial u}{\partial x}+v\frac{\partial u}{\partial y}
=  \nu\frac{\partial^{2} u}{\partial y^{2}} + U\frac{\partial U}{\partial x}+ \frac{\partial U}{\partial t},  \\
\frac{\partial u}{\partial x}+\frac{\partial v}{\partial y} =0,
\end{array}
 \label{ww35}         \right.\end{equation}
in a domain $D=\{0<t<T,0<x<X, 0<y<\infty\}$, with the boundary conditions
\begin{equation}\left\{
\begin{array}{ll}
u(0,x,y)=U(0,x), \ \ u(t,0,y)=0,\ \ u(t,x,0)=0, \ \  v(t, x,0) =v_{0}(t,x) ,  \\
u(t,x,y)\rightarrow U(t,x)  \ \  \text{as}\ \  y\rightarrow \infty.
\end{array}
 \label{ww36}         \right.\end{equation}
 It is also assumed that $U(t,0)=0, U(t,x) >0$ for $ x>0, t\geq 0$;  $ v_{0}(t,x), U(t,x),  U_{x}, U_{t}/U $are bounded.

Introduce new independent variables in (\ref{ww35})-(\ref{ww36}),  setting
$$\tau =\sqrt{t} , \ \  \xi=x, \ \ \eta=\frac{u(t,x,y)}{U(t,x)},$$
and a new unknown function
$$w(\tau, \xi,\eta)=  \frac{\sqrt{t}u_{y}(t, x,y)}{U(t,x)}.$$
From  (\ref{ww35}), we obtain the following equation for $w(\tau, \xi,\eta)$:
\begin{eqnarray}
\nu w^{2}w_{\eta\eta}- \frac{\tau}{2} w_{\tau}-\tau^{2}\eta U w_{\xi} +A  w_{\eta}+B w=0,
 \label{ww37}
\end{eqnarray}
where
$$A =\tau^{2}(\eta^{2}-1)U_{x}+\tau^{2}(\eta-1)\frac{U_{ t}}{U }, \ \  B =-\eta \tau^{2}U_{x}(\xi) -\tau^{2}\frac{U_{ t}}{U } +\frac{1}{2},  $$
in the domain $\Omega=\{0<\tau<\sqrt{T} , 0<\xi<X, 0<\eta <1 \}$, with the boundary conditions
\begin{eqnarray}
   w|_{\eta=1}=0, \ \  (\nu w w_{\eta }-v_{0}w\tau  +C  )_{\eta=0}=0,
 \label{ww38}
\end{eqnarray}
where
$$  C =\tau^{2}U_{x}+\tau^{2}\frac{U_{ t}}{U }.$$

\begin{Theorem} \label{OS.5.36} (\cite{OS})
Assume that the functions $U, U_{x}, U_{xx}, U_{t}/U, (U_{t}/U)_{x}, \sqrt{t}U_{xt},\sqrt{t} (U_{t}/U)_{x}, v_{0}$ are bounded, and
$$U_{x}(t,0)>0, \ \  U_{x}+\frac{U_{t}}{U}>0, \ \ -M_{1}U< \frac{U_{t}}{U}$$
for small $x$. Assume also that
$$v_{0}\leq M_{2}t^{\frac{1}{2}+\varepsilon}, \ \  |v_{0x}|\leq M_{3}t^{\frac{1}{2} }, \ \ t^{\frac{1}{2} }|v_{0t}|\leq M_{4},  $$
where $\varepsilon$ is a positive constant. Then problem (\ref{ww35})-(\ref{ww36}), in the domain $D$ with small enough $T$, has one and only one solution $u,v$
with the following properties: $u/U$ is continuous for $t=0$, $y>0$ and for $t>0, y\geq 0$; $\sqrt{t}u_{y}/U$ is bounded and continuous in $\overline{D}$;
$u\rightarrow U$ as $y\rightarrow \infty$; $u_{y}, u_{x}, u_{yy}, u_{t}, u_{y}$ bounded for $t\geq t_{0}>0$; conditions (\ref{ww36}) are satisfied
for $u, v$, and equation (\ref{ww35}) holds almost everywhere in $D$. Moreover, the following inequalities are valid:
\begin{eqnarray*}
&&  C_{8}\left(1-\frac{u}{U}\right)\sigma_{\mu}(\frac{u}{U})\leq \frac{\sqrt{t}u_{y}}{U} \leq C_{7}\left(1-\frac{u}{U}\right)\sigma_{\mu_{1}}(\frac{u}{U}) , \\
&& -C_{1} \sigma_{\mu_{1}}(\frac{u}{U})\leq \frac{\sqrt{t}u_{yy}}{u_{y}} \leq C_{12}\Phi(\frac{u}{U}),  \\
&& \left|\frac{2tu_{ yy} }{u_{ y}} \frac{uU_{t}-u_{t}U}{U^{2}}+\frac{ u_{ y } }{U}+2t\left(\frac{u_{yt}U-u_{y}U_{t}}{U^{2}}\right) \right|
\leq C_{10}\left(1-\frac{u}{U}\right)\sigma_{\mu_{1}}(\frac{u}{U}), \\
&& U(t,x)\exp \left(-\frac{C_{7}^{2}y^{2}}{4t}-\frac{C_{7}y\sqrt{-\ln \mu_{1}}}{\sqrt{t}} \right)\leq U(t,x)-u\leq
U(t,x)\exp \left(-\frac{C_{8}^{2}y^{2}}{4t}-\frac{C_{8}y\sqrt{-\ln \mu_{1}}}{\sqrt{t}} \right) ,
\end{eqnarray*}
where $\sigma_{\gamma} =\sqrt{-\ln \gamma(1-\eta)},K_{i} $, $0<\mu<e^{-1/2} <\mu_{1}<1$; $\Phi(s)=-s^{2}$ for $0\leq s\leq 1-\delta_{0}$;
$\Phi(s)=-\sigma_{\mu}(s)$ for $1-\delta_{0}<s<1, \delta_{0}>0$; $C_{7},C_{8}, C_{9}, C_{10}, C_{11},C_{12}  $ are positive constants.

\end{Theorem}

 Oleinik and Samokhin (\cite{OS}) in 1999 considered
  the boundary layer theory   for the numerical solution of these problems  by finite differences method,
the nonstationary system of the two-dimensional boundary layer
\begin{equation}\left\{
\begin{array}{ll}
\frac{\partial u}{\partial t}+ u\frac{\partial u}{\partial x}+v\frac{\partial u}{\partial y}
=  \nu\frac{\partial^{2} u}{\partial y^{2}} + U\frac{\partial U}{\partial x}+ \frac{\partial U}{\partial t},  \\
\frac{\partial u}{\partial x}+\frac{\partial v}{\partial y} =0,
\end{array}
 \label{ww39}         \right.\end{equation}
in a domain $D=\{0<t<T,0<x<X, 0<y<\infty\}$, with the boundary conditions
\begin{equation}\left\{
\begin{array}{ll}
u(0,x,y)=u_{0}(x,y), \ \ u(t,0,y)=u_{1}(t,y),\ \ u(t,x,0)=0, \ \  v(t, x,0) =v_{0}(t,x) ,  \\
u(t,x,y)\rightarrow U(t,x)  \ \  \text{as}\ \  y\rightarrow \infty.
\end{array}
 \label{ww40}         \right.\end{equation}

Introduce new independent variables in (\ref{ww39})-(\ref{ww40}),  setting
$$\tau =t, \ \  \xi=x, \ \ \eta= u(t,x,y) ,$$
and a new unknown function
$$w(\tau, \xi,\eta)=  u_{y}(t, x,y), $$
 we obtain the following equation :
\begin{eqnarray}
\nu w^{2}w_{\eta\eta}-  w_{\tau}- \eta U w_{\xi} +p_{x} w_{\eta} =0,
 \label{ww41}
\end{eqnarray}
in the domain $\Omega=\{0<\tau<T, 0<\xi<X, 0<\eta <U(\tau, \xi) \}$, with the boundary conditions
\begin{eqnarray}
   w|_{\tau=0}=u_{0y}=w_{0}(\xi,\eta), \ \  w|_{\xi=0}=u_{1y}=w_{1}(\tau,\eta), \ \ w|_{\eta=U{\tau,\xi}}=0 , \ \  (\nu w w_{\eta }-v_{0}w-p_{x}  )|_{\eta=0}=0.
 \label{ww42}
\end{eqnarray}
Denote by $\Omega'$ the union of the domain $\Omega$ and the part of its boundary
on which the initial and the boundary conditions are prescribed. A node is
said to be internal if its neighboring nodes belong to $\Omega'$. The remaining
nodes belonging to $\Omega'$ are called boundary nodes.

For any function $f(\tau, \xi,\eta)$,  denoted by $f^{m,l,k}$ its value at the node $(mh, l\sigma,k\sigma )$.
They associated a finite difference equation for $w$ which approximates equation (\ref{ww41}), namely,
\begin{eqnarray}
\begin{aligned}
&  \left(\nu(w^{m,l,k})^{2}+M\sigma\right)\frac{w^{m,l,k+1}-2w^{m,l,k}+w^{m,l,k-1} }{\sigma^{2}}-\frac{w^{m+1,l,k}-w^{m,l,k}}{h} \\
&  -k\sigma \frac{w^{m,l,k}-w^{m,l-1,k}}{\sigma}+p_{x}^{m,l,k}\frac{w^{m,l,k}-w^{m,l,k-1}}{\sigma}=0,
 \label{ww43}
\end{aligned}
\end{eqnarray}
where $M$ is a positive constant such that $M>\max |p_{x}|$.
At the boundary nodes of O' belonging to the planes $\tau=0, \xi=0, \eta=0$,
we set
\begin{equation}\left\{
\begin{array}{ll}
  w_{0}^{0,l,k}=w_{0}(l\sigma, k\sigma), \ \ w^{m,0,k}=w_{1}(mh,k\sigma),   \\
  \nu w^{m,l,0}\frac{w^{m+1,l,1}-w^{m+1,l,0}}{\sigma}- p_{x}^{m,l,0}-v_{0} ^{m,l,0}w^{m,l,0} =0.
\end{array}
 \label{ww44}         \right.\end{equation}
At every boundary node of $\Omega'$ that does not belong to any of the planes $\tau=0, \xi=0, \eta=0$ (the set of all such points is denoted
by $\Gamma_{h\sigma}$), we set
\begin{equation}
 w^{m,l,k} =0.
 \label{ww45}
\end{equation}

\begin{Theorem} \label{OS.6.37} (\cite{OS})
Let $W$ be a solution of problem (\ref{ww41})-(\ref{ww42}) with bounded second derivatives in $\Omega'$, and let $w$ be the solution of finite difference
 equations (\ref{ww43})-(\ref{ww45}). Then, for small enough $h$ and $\sigma$, the following inequality holds:
$$|W-w|\leq M_{3}(h+\sigma), $$
if either $\tau\leq T$ and $h/\sigma^{2}<1/(2\nu b_{1}^{2})$ or $\xi\leq X$ and $h/\sigma^{2}<1/(2\nu b_{2}^{2})$.

\end{Theorem}

By $f^{m,l,k}$ they denoted the value of a function $f$ at the point $(mh, l\sigma,k\sigma )$. To every internal node of $\Omega'$ with the coordinates
 $((m+1)h, l\sigma,k\sigma )$, they associated a difference equation for $w$ which approximates equation (\ref{ww41}):
\begin{eqnarray}
\begin{aligned}
&  \left(\nu(w^{m,l,k})^{2}+Mh\right)\frac{w^{m+1,l,k+1}-2w^{m+1,l,k}+w^{m+1,l,k-1} }{h^{2}}-\frac{w^{m+1,l,k}-w^{m,l,k}}{h} \\
&  -kh \frac{w^{m+1,l,k}-w^{m+1,l ,k-1}}{h}+p_{x}^{m,l,k}\frac{w^{m+1,l,k}-w^{m+1,l,k-1}}{h}=0,
 \label{ww46}
\end{aligned}
\end{eqnarray}
where $M=$const.$>\max |p_{x}|$.
To the boundary conditions (\ref{ww42}), we associate the following equations at boundary nodes:
\begin{equation}\left\{
\begin{array}{ll}
  w ^{0,l,k}=w_{0}(lh, kh), \ \ w^{m,0,k}=w_{1}(mh,kh),  \ \ \text{on}\ \  \Gamma_{h} , \\
  \nu w^{m,l,0}\frac{w^{m+1,l,1}-w^{m+1,l,0}}{h}- p_{x}^{m,l,0}-v_{0} ^{m,l,0}w^{m,l,0} =0.
\end{array}
 \label{ww47}         \right.\end{equation}
By $\Gamma_{h}$, we denote the boundary nodes of $\Omega'$ that do not belong to the planes  $\tau=0, \xi=0, \eta=0$.

\begin{Theorem} \label{OS.6.38} (\cite{OS})
 Let $W$ be the solution of problem (\ref{ww41})-(\ref{ww42}) with bounded second derivatives in $\Omega'$, and let $w$ be the solution of finite difference
 equations (\ref{ww46})-(\ref{ww47}). Then, at the nodes of $\Omega'$, for small enough $h$, the inequality
$$|W-w|\leq k_{1}h$$
holds either for $\tau\leq T$ or for $\xi\leq X$.

\end{Theorem}

 Oleinik and  Samokhin (\cite{OS}) in 1999
applied the finite difference method to construct a solution of the boundary layer system in the case of a non-stationary axially symmetric flow
(\ref{ww21})-(\ref{ww22}) in a domain $D=\{0<t<\infty,0<x<X, 0<y<\infty\}$. They  replaced equation (\ref{ww23}) and conditions (\ref{ww24})
 by the following linear system of algebraic equations with respect to $w^{m+1,l,k}$:
 \begin{eqnarray}
\begin{aligned}
&  \left(\nu(w^{m,l,k})^{2} \right)\frac{w^{m ,l,k+1}-2w^{m ,l,k}+w^{m ,l,k-1} }{\delta^{2}}-\frac{w^{m+1,l,k}-w^{m,l,k}}{h}
 -k\delta U^{m,l,k} \frac{w^{m ,l,k} -w^{m ,l-1 ,k-1}}{d}\\
 &+A^{m,l,k}\frac{w^{m,l,k}-w^{m,l,k-1}}{\delta} +B^{m,l,k} w^{m ,l,k } =0,
 \label{ww48}
\end{aligned}
\end{eqnarray}
with $m=0,1...; \ \ l=0,1,...[X/ \varepsilon]; \ \ k=1,..., J-1,$ and the following conditions
\begin{equation}\left\{
\begin{array}{ll}
w^{m+1, l,0}=\frac{1}{2}\left[w^{m+1, l,1} -\delta\nu^{-1}v_{0}^{m+1, l,0}+\sqrt{(w^{m+1, l,0}-\delta\nu^{-1}v_{0}^{m+1, l,0})^{2}
+4-\delta\nu^{-1} C^{m+1, l,0}} \right], \\
  w ^{0,l,k}=w_{0}(ld, k\delta), \ \   l=0,1,...[X/ d]; \ \ k=0,..., J;\\
  w^{m+1,l,J}=0, \ \   l=0,1,...[X/ d]; \ \ m=0,1,... .
\end{array}
 \label{ww49}         \right.\end{equation}

\begin{Theorem} \label{OS.6.39} (\cite{OS})
Let $A,B,C, v_{0}, U$ be bounded functions in $\Omega$. Assume that $U_{x}>0, C>0$ in $\Omega$ and
$$E_{1}(1-\eta)\sigma\leq w_{0}\leq E_{2}(1-\eta)\sigma, $$
where $\sigma=\sqrt{-\ln \mu(1-\eta)}, \mu<1$,  and let the  function $W(\tau,\xi,\eta)$ be a solution of problem (\ref{ww23})-(\ref{ww24}) such that
\begin{eqnarray*}
&& M_{12}(1-\eta)\sigma\leq W \leq M_{13}(1-\eta)\sigma, \\
&& |W_{\xi}|\leq M_{14}(1-\eta)\sigma, \ \  -M_{15} \sigma \leq W_{\eta} \leq -M_{16} \sigma, \\
&& - M_{17}(1-\eta)\sigma\leq W_{\tau} \leq M_{18}(1-\eta)\sigma, \\
&& - M_{19} \leq WW_{\eta\eta} \leq -M_{20} .
\end{eqnarray*}
Assume, in addition, that $W$ has bounded weak derivatives $W_{\tau\tau}, W_{\xi\xi}$ in $\Omega$, and $W_{\eta\eta } $ satisfies the inequality
$$|W_{\eta\eta }(\tau,\xi,\eta_{1})- W_{\eta\eta }(\tau,\xi,\eta_{2})|\leq \frac{M_{17}|\eta_{1}-\eta_{2}|}{(1-\eta_{2})^{2} \sigma(\eta_{2})}  $$
for any $\eta_{1}, \eta_{2}$ such that $0<\eta_{1}<\eta_{2}<1$. Then for small enough $\delta$ and any $m, l, k, mh\leq T$ (T is an arbitrary positive constant),
we have
$$|W^{m+1,l , k}-w^{m+1,l , k}|\leq M_{18}(h+d+\delta\sigma^{J-1}), $$
provided that $h\delta^{-2}\leq (2\nu a^{2}+b\delta^{2}d^{-1})^{-1}$,  where $w^{m+1,l , k}$ is a solution of (\ref{ww48})-(\ref{ww49}), and
 the constant $M_{18}$ is  independent of $h , d, \delta$, but may depend on $T$.

\end{Theorem}

 Xin and Zhang (\cite{XZ}) in 2004 obtained the existence of global weak solutions to the Prandtl equation (\ref{n0}) when the pressure is favourable  $(\partial_xp\leq 0)$,
 and the following conditions are satisfied,
\begin{equation}\left\{\begin{array}{ll}
U(x,t)>0, \ u_{0}(x,t)>0, \ u_{1}(y,t)>0, \ v_{0}(x,t)\leq 0,\\
\partial_{y}u_{0}(x,y)>0, \ \partial_{y}u_{1}(x,y)>0,\\
U(x,t)\equiv d= constant.
\end{array}\right.\label{n1}\end{equation}
The considered problem in \cite{XZ} is the following initial-boundary value problem
\begin{equation}\left\{\begin{array}{ll}
\partial_{\tau}w^{-1}+\eta U\partial_{\xi}w^{-1}+A\partial_{\eta}w^{-1}-Bw^{-1}=-v\partial_{\eta}^{2}w \\
on \ Q=\{(\xi,\eta,\tau)| 0<\tau<\infty, \ 0<\xi<L, \ 0<\eta<1\},\\
w|_{\tau=0}=\frac{\partial_{y}u_{0}}{U}\equiv w_{0}, \ \  w|_{\eta=1}=0,\\
w|_{\xi=0}=w_{1}, \ \ and \ (vw\partial_{\eta}w-v_{0}w)|_{\eta=0}=\frac{\partial_{y}U_{1}(0,y,t)}{U(0,t)},
\end{array}\right.\label{n0}\end{equation}
where $A=(1-\eta^{2})\partial_{x}U+(1-\eta)\frac{\partial_{t}U}{U}$, $B=\eta\partial_{x}U+\frac{\partial_{t}U}{U}$ and $w_{1}(\tau,\eta)=\frac{\partial_{y}U_{1}(0,y,t)}{U(0,t)}$.

To study problem \eqref{n0}, we need the following Lemmas 1.1.79-1.1.88. To this end, we consider the  problem
\begin{equation}\left\{\begin{array}{ll}
u_{t}-u^{2}u_{yy}=0, \ 0<y<d,\\
u|_{t=0}=u_{0}>0,\\
\frac{\partial u}{\partial y}|_{y=0}=v_{0}, \ u|_{y=d}=0.
\end{array}\right.\label{n2}\end{equation}
\begin{Lemma}(\cite{XZ})
Let $\phi=e^{\alpha y}\sin\alpha(d-y)$, where $\alpha=\frac{\pi}{2d}$: Suppose $v_{0}\leq 0$. Then there exist constants $\beta$ and $C_{0}$ which depend only on d $|u_{0}|_{L^{\infty}}$ and $|v_{0}|_{L^{\infty}}$ such that the solution of
\eqref{n2} satisfies
\begin{equation}
\begin{array}{ll}
|u|\leq C_{0}(d-y),
\end{array}\label{n3}
\end{equation}

\begin{equation}
\begin{array}{ll}
u\geq \theta_{0}e^{-\beta t}\phi,
\end{array}\label{n4}
\end{equation}
provided that \eqref{n3} and \eqref{n4} hold true initially, where $\theta_{0}=\min\frac{u_{0}}{\phi}$.
\end{Lemma}
\begin{Lemma}(\cite{XZ})
\begin{equation}
\begin{array}{ll}
|u_{y}^{2}|_{L^{\infty}}\leq a_{0}^{2}.
\end{array}\label{n5}
\end{equation}
\end{Lemma}

We estimate $w=u_x$ which solves the following problem:

\begin{equation}\left\{\begin{array}{ll}
w_{t}-(u^{2}+\varepsilon)w_{yy}=2uu_{yy}w+2uu_{yy}v_{0,x}(y-d)-v_{0,tx}(y-d),\\
 w|_{t=0}=u_{0,x}-v_{0,x}(y-d)\equiv w_{0},\\
\frac{\partial w}{\partial y}=0, \ w|_{y=d}=0.
\end{array}\right.\label{n'}\end{equation}

\begin{Lemma}(\cite{XZ})
Suppose that $v_{0}\leq 0$. Then there exists a constant $C_{6}$ depending only on T, d, $|u_{0,x}|_{L^{\infty}}$, $|u_{0}|_{L^{\infty}}$, $|u_{y}|_{L^{\infty}}$, $|v_{0}|_{C^{2}}$ and $\int_{0}^{d}\frac{u^{2}_{0,x}}{(y-d)^{2}}dy$, such that the solution to problem
\eqref{n'} satisfies
\begin{equation}
\begin{array}{ll}
|u_{x}(t,x,y)|\leq C_{6}.
\end{array}\label{n6}
\end{equation}

\end{Lemma}

For the initial-boundary value problem,
\begin{equation}\left\{\begin{array}{ll}
u_{t}+yu_{x}-u^{2}u_{yy}=0, \  (0,T)\times\Omega ,\\
u|_{t=0}=u_{0,y}, \  \  u|_{x=0}=u_{1,y},     \\
\frac{\partial u}{\partial y}|_{y=0}=v_{0}, \ u|_{y=1}=0,
\end{array}\right.\label{n8}\end{equation}
we have the following three lemmas
\begin{Lemma}(\cite{XZ})
For $ 0\leq t\leq T, (x,y)\in\Omega$, it holds for solution $u$ to problem \eqref{n8} that
\begin{equation}
\begin{array}{ll}
|u_{t}(t,x,y)|\leq |\overline{u}_{0}|_{L^{\infty}}+2|v_{0}|_{L^{\infty}}+|\overline{u}_{1}|_{L^{\infty}}.
\end{array}\label{n9+}
\end{equation}
\begin{equation}
\begin{array}{ll}
|u(t,x,y)|\leq |\overline{u}_{0}|_{L^{\infty}}+2|v_{0}|_{L^{\infty}}+|\overline{u}_{1}|_{L^{\infty}}.
\end{array}\label{+4}
\end{equation}
\end{Lemma}

\begin{Lemma}(\cite{XZ})
Suppose $v_{0}\leq 0$. Then there exist constants $C_{0}$ and $\beta$ depending only on
$|v_{0}|_{L^{\infty}}$, $|\overline{u}_{0}|_{L^{\infty}}$ and $|\overline{u}_{1}|_{C^{1}}$, such that the solution $u$ to problem \eqref{n8}
\begin{equation}
\begin{array}{ll}
u_{t}(t,x,y)\leq C_{0}(1-y),
\end{array}\label{n9}
\end{equation}
\begin{equation}
\begin{array}{ll}
u(t,x,y)\geq \theta_{0}\varepsilon^{-\beta t}\phi.
\end{array}\label{n10}
\end{equation}
\end{Lemma}

\begin{Lemma}(\cite{XZ})
There exists a constant $C_{3}$ depending on$|\overline{u}_{0}|_{L^{\infty}}$, $|v_{0}|_{L^{\infty}}$, and $|\overline{u}_{1}|_{C^{1}}$ such that the solution $u$ to problem \eqref{n8}
\begin{equation}
\begin{array}{ll}
 \int \int_{\Omega}|u_{x}(t,\cdot)|dx dy \leq C_{3}e^{C_{3}t}[\int\int_{\Omega}(C_{3}^{2}|\overline{u}_{0,x}|)dx dy+\frac{L}{2}|v_{0,x}|_{L^{\infty}}+\frac{1}{2}\int_{0}^{1}\frac{|\overline{u}_{1,t}|}{\overline{u}_{1}^{2}}(1-y)^{2}dy].

\end{array}\label{n12}
\end{equation}

\end{Lemma}

Now  we consider the following problem:
\begin{equation}\left\{\begin{array}{ll}
u_{t}=u^{2}u_{yy}+b(t,x,y)u, \ \ (t,x,y)\in (0,T)\times\Omega,\\
u|_{t=0}=u_{0,y}>0,\\
\frac{\partial u}{\partial y}|_{y=0}=v_{0}, \ \ u|_{y=1}=0,
\end{array}\right.\label{n13}\end{equation}

which is approximated by
\begin{equation}\left\{\begin{array}{ll}
u_{t}=(u^{2}+\varepsilon)u_{yy}+b(t,x,y)u, \ \ (t,x,y)\in (0,T)\times\Omega,\\
u|_{t=0}=u_{0,y}>0,\\
\frac{\partial u}{\partial y}|_{y=0}=v_{0}, \ \ u|_{y=1}=0.
\end{array}\right.\label{n''}\end{equation}

The following Lemmas 1.1.85-1.1.88 are concerned with problem \eqref{n''}.
\begin{Lemma}(\cite{XZ})
The solution of problem \eqref{n13} satisfies
\begin{equation}\left\{\begin{array}{ll}
|v_{0}|_{L^{\infty}}(1-y)\mp ue^{-Bt}\leq \max_{x,y}(|v_{0}|(1-y)\mp ue^{-Bt_{0}})|_{t=t_{0}},\\
|u|\leq C_{0}(1-y)e^{Bt},\\
u\geq \theta_{0} e^{-\beta t}\phi.
\end{array}\right.\label{n14}\end{equation}
\end{Lemma}

\begin{Lemma}(\cite{XZ})
There exists a constant $C_{4}$, independent of E, such that the solution $u$ to \eqref{n''}
as well as \eqref{n13} satisfies
\begin{equation}
\begin{array}{ll}
|u_{y}(t,x,y)|\leq C_{4}.
\end{array}\label{n15}
\end{equation}
\end{Lemma}

\begin{Lemma}(\cite{XZ})
There exists a constant $C_{5}$, independent of E, such that the solution $u$ to \eqref{n''}
as well as \eqref{n13} satisfies
\begin{equation}
\begin{array}{ll}
|u_{x}(t,x,y)|+|u_{t}(t,x,y)|\leq C_{5}.
\end{array}\label{n16}
\end{equation}
\end{Lemma}
\begin{Lemma}(\cite{XZ})
There exists a constant $C_{9}$, depending on $T, \Omega, B$, $|\overline{u}_{0}|_{L^{\infty}}$, $|v_{0}|_{L^{\infty}}$, and $|\overline{u}_{1}|_{C^{1}}$ such that the solution $u$ of \eqref{n''} as well as \eqref{n13} satisfies
\begin{equation}
\begin{array}{ll}
 \int \int_{\Omega}|u_{x}(t,\cdot)|dx dy \leq e^{C_{9}t}[\int\int_{\Omega}|\overline{u}_{0x}|+C_{9}]dxdy,
\end{array}\label{n17}
\end{equation}
for $t\leq T$.
\end{Lemma}

Based on Lemmas 1.1.79-1.1.88, we have the following two theorems.

\begin{Theorem}(\cite{XZ})
Assume that the data satisfy conditions \eqref{n1}. Then there exists
a weak solution $w\in BV(Q_{T})\cap L^{\infty}(Q_{T})$ to the initial value problem \eqref{n0} provided that
the pressure is favourable, i.e., $\partial_{x}p(x,t)\leq0$ holds for $t>0$, $0<x<L$.
\end{Theorem}

\begin{Theorem}(\cite{XZ})
Assume that $v_{0}\leq 0$ and the data are smooth and compatible. Then problem \eqref{n2}
 has a unique positive bounded smooth solution in $(0,T]\times\Omega$. Moreover
\begin{eqnarray}
u\in L_{p}([0,T]\times\overline{\Omega})
\end{eqnarray}
and u satisfies estimates \eqref{n3}-\eqref{n'} and \eqref{n9+}.
\end{Theorem}
From above estimates on problem  \eqref{n8}, we have the following theorem.

\begin{Theorem}(\cite{XZ})
Assume that \eqref{n1} hold and $b(t,x,y)=0$. Then problem
\eqref{n8} has a weak solution
\begin{eqnarray}
u\in L^{\infty}(0,T,BV(\Omega))\cap BV((0,t)\times\Omega),
\end{eqnarray}
which satisfies estimates \eqref{+4}-\eqref{n10}. Furthermore, its derivatives satisfy
\eqref{n12} and
\begin{equation}
\begin{array}{ll}
 \int \int_{\Omega}|u_{y}(t,\cdot)|dx dy \leq
 \int_{0}^{t}\int_{0}^{1}|\overline{u}_{1,y}|dy ds+2L|u|_{L^{\infty}}+\int \int_{\Omega}|\overline{u}_{0,y}|dx dy+\int_{0}^{t}\int \int_{\Omega}|u_{x}(s,\cdot)|dx dy ds,
\end{array}
\end{equation}
\begin{equation}
\begin{array}{ll}
 \int \int_{\Omega}|u_{t}(t,\cdot)|dx dy \leq
 \int_{0}^{1}[\frac{y(1-y)}{u(t,0,y)}+\frac{y(1-y)}{u(t,L,y)}]dy +\int_{0}^{1}|v_{0}(t,\cdot)|dx+\int \int_{\Omega}|u_{y}(t,\cdot)|dx dy.
\end{array}\label{n*}
\end{equation}
\end{Theorem}

\begin{Theorem}(\cite{XZ})
Suppose that $v_{0}\leq 0$ and the data are smooth and compatible. Then problem \eqref{n13} has a unique positive bounded smooth solution and
\begin{eqnarray}
u\in L_{p}([0,T]\times \overline{\Omega}),
\end{eqnarray}
Moreover, the solution satisfies \eqref{n14}-\eqref{n16}.
\end{Theorem}

\begin{Lemma}(\cite{XZ})
There exists a constant $C_{11}$ depending on $T, \Omega, B$, $|\overline{u}_{0}|_{L^{\infty}}$, $|v_{0}|_{L^{\infty}}$, and $|\overline{u}_{1}|_{C^{1}}$ such that the solution $u$ of problem \eqref{n''} satisfies
\begin{equation}
\begin{array}{ll}
 \int \int_{\Omega}|u_{y}(t,\cdot)|dx dy \leq e^{Bt}[\int\int_{\Omega}[|\overline{u}_{0y}|+\frac{n}{T}\Sigma_{j=1}^{i}|u_{x}(t_{j-1},\cdot)|]dxdy+C_{11}.
\end{array}\label{n17}
\end{equation}
\end{Lemma}

From Theorems 1.1.91-1.1.92 and Lemma 1.1.93, we conclude the following theorem.

\begin{Theorem}(\cite{XZ})
Under assumptions \eqref{n1}, problem \eqref{n8} has a solution
\begin{eqnarray}
u\in L^{\infty}(0,T;BV(\Omega))\cap BV((0,T)\times\Omega)
\end{eqnarray}
satisfying \eqref{n14}. Furthermore, its derivatives are estimated by \eqref{n17}, \eqref{n*} and
\begin{equation}
\begin{array}{ll}
 \int \int_{\Omega}|u_{y}(t,\cdot)|dx dy \leq
 e^{Bt}\int \int_{\Omega}[|\overline{u}_{0y}|+\int_{0}^{t}|u_{x}|ds]dx dy +C_{11}.
\end{array}\label{n*}
\end{equation}

\end{Theorem}

Xu and Zhao (\cite{XZ1}) in 2009 studied the global existence and uniqueness of weak solutions to the non-stationary boundary layer (\ref{1.1.6}) by Croco transformation.
The initial boundary value problem is
\begin{equation}
\begin{array}{ll}
w_{\tau}+\eta Uw_{\xi}+Aw_{\eta}+Bw-w^{2}w_{\eta\eta}=0 \ \ \  in \ \Omega_{T},
\end{array}\label{h*1}
\end{equation}

\begin{equation}
\begin{array}{ll}
w(0,\xi,\eta)=w_{0}(\xi,\eta)=\frac{u_{0y}}{U(0,x)}, \ w(\tau,\xi,1)=0,
\end{array}\label{h*2}
\end{equation}

\begin{equation}
\begin{array}{ll}
ww_{\eta}-v_{0}w+(U_{\xi}+\frac{U_{\tau}}{U})=0 \ \ \ on \ \eta=0,
\end{array}\label{h*3}
\end{equation}
where $\Omega_{T}=\{(\tau,\xi,\eta)|0<\tau<T,0<\xi<L,0<\eta<1\}$ and
\begin{eqnarray}
A=(1-\eta^{2})U_{\xi}+(1-\eta)\frac{U_{\tau}}{U}, \ B=\eta U_{\xi}+\frac{U_{\tau}}{U}.
\end{eqnarray}

\begin{Theorem}(\cite{XZ1})
Assume that $U=\xi^{m}U_{1}(\tau,\xi)$ $m\geq1$, $U_{1}(\tau,\xi)$, $w_{0}(\xi,\eta)$, $v_{0}(\tau,\xi)$ are appropriately smooth and
\begin{eqnarray}
U_{1}(\tau,\xi)>0, \ U_{\tau}+UU_{\xi}>0, \ v_{0}(\tau,\xi)\leq0, \  w_{0}(\tau,\xi)\geq0,
\end{eqnarray}
and there exists a constant $C_{0}$ such that
\begin{eqnarray}
C_{0}^{-1}(1-\eta)\leq w_{0}(\xi,\eta)\leq C_{0}(1-\eta),
\end{eqnarray}
where $w_{0}$ are the functions in \eqref{h*2}.
\end{Theorem}

\begin{Theorem}(\cite{XZ1})
The solution of \eqref{h*1}-\eqref{h*3} is unique.
\end{Theorem}

 Xu and  Zhang (\cite{xz1})  in 2015 proved the local-in-time existence, uniqueness and stability of solutions for the nonlinear Prandtl equations
 (\ref{1.1.6}) where $U$ is a cnstant in weighted Sobolev space with shear flow without monotonicity in the domain $\mathbb{R}_{+}^{2}$,
  and proved  the global-in-time existence of solutions to the Prandtl equations (\ref{1.1.6}) . \\
They assumed that
$u_{0}(x,y)=u_{0}^{s}(y)+\tilde{u}_{0}(x,y),$
where $u_{0}^{s}$ satisfies the following conditions: for some even integer $m\geq 6$, and $k>1$, $u_{0}^{s}\in C^{m+4}([0,+\infty)), \lim\limits_{y\rightarrow+\infty} u^{s}_{0}(y)=1$ with the compatibility conditions
$(\partial_{y}^{2p}u^{s}_{0})(0)=0$, $0\leq 2p\leq m+4$, and that there exists $0<a<+\infty, c_{1},c_{2}>0$ such
that
\begin{equation}\left\{
\begin{array}{ll}
\partial_{y}u_{0}^{s}(a) =0, \partial_{y}^{2}u_{0}^{s}(a) \neq 0, \text{and} ~ \partial_{y} u_{0}^{s}(y) \neq 0, \text{for any}\ y\in \bar{\mathbb{R}}_{+}\ \{a\}, \\
 c_{1}\langle y\rangle ^{-k}\leq |\partial_{y}u_{0}^{s}(y) |\leq c_{2}\langle y\rangle ^{-k}, \forall y\geq a+1,\\
| \partial_{y}^{ p}u^{s}_{0} (y)|\leq  c_{2}\langle y\rangle ^{-k-p+1}, \forall y\geq 0, 0\leq p\leq m+4.
\end{array}
 \label{w.28}         \right.\end{equation}
The weighted Sobolev spaces, for $\lambda\in \mathbb{R}, m\in \mathbb{N}$,
$$\|u\|_{H^{m}_{\lambda}(\mathbb{R}_{+}^{2}) }^{2}= \sum\limits_{|\alpha_{1}+\alpha_{2}|\leq m}
\int_{\mathbb{R}_{+}^{2} }\langle y\rangle ^{2\lambda+2\alpha_{2}}|\partial_{x}^{\alpha_{1}}\partial_{y}^{\alpha_{2}} u|^{2}dxdy ,$$
$$\|u\|_{L^{2}_{\lambda}(\mathbb{R}_{+}^{2}) }^{2}=
\int_{\mathbb{R}_{+}^{2} }\langle y\rangle ^{2\lambda }| u|^{2}dxdy.$$

\begin{Theorem} (\cite{xz1} )\label{t.xz1}
Let $m\geq 6$ be an even integer, $k\geq 1, 0\leq \ell\leq \frac{1}{2}, k+\ell>\frac{1}{2}$. Assume that $u^{s}_{0}$ satisfies (\ref{w.28}), the initial data $u_{0}-u_{0}^{s}\in H^{m+3}_{2k+\ell-1}( \mathbb{R}_{+}^{2})$ and $u_{0}-u_{0}^{s}$ satisfies the compatibility condition up to order $m+2$. Then there exists a $T>0$ such that if
\begin{eqnarray}
 \|(u_{0}-u_{0}^{s})\|_{H^{m+3}_{2k+\ell-1}( \mathbb{R}_{+}^{2}) } \leq \delta_{0},
\label{w.29}
\end{eqnarray}
for some $\delta_{0}>0$ small enough, the initial-boundary value problem (\ref{1.1.6})  admits a unique solution $(u,v)$ with
$$ u -u ^{s}\in L^{\infty} ([0,T];H^{m }_{ k+\ell-1}( \mathbb{R}_{+}^{2}) ),
\partial_{y} (u -u ^{s})\in L^{\infty} ([0,T];H^{m }_{ k+\ell }( \mathbb{R}_{+}^{2}) ),  $$
and
$$v \in L^{\infty} ([0,T];  L^{\infty}(\mathbb{R}_{y,+}; H^{m-1}( \mathbb{R}_{x})) ),
\partial_{y} v \in L^{\infty} ([0,T];H^{m -1}_{ k+\ell-1 }( \mathbb{R}_{+}^{2}) ) .$$
Moreover, we have the stability with respect to the initial data in the following sense: given any two initial data
$$u_{0}^{1}=u_{0}^{s}+\tilde{u}_{0}^{1}, \ \  u_{0}^{2}=u_{0}^{s}+\tilde{u}_{0}^{2}, $$
if $u_{0}^{s} $ satisfies  problem  (\ref{w.28}) and $\tilde{u}_{0}^{1},\tilde{u}_{0}^{2} $ satisfies (\ref{w.29}), then the solutions $u^{1}$ and $u^{2}$ of problem  (\ref{1.1.6})  with initial data $ u_{0}^{1}$ and $ u_{0}^{2}$ respectively satisfy,
$$ \|u^{1}-u^{2}\|_{  L^{\infty} ([0,T];H^{m-2 }_{ k+\ell -1}( \mathbb{R}_{+}^{2}) )}\leq \|u^{1}_{0}-u^{2}_{0}\|_{  H^{m  }_{ 2k+\ell -1}( \mathbb{R}_{+}^{2})}. $$

\end{Theorem}
Assume  the initial data of the shear flow  $u_{0}^{s} $ is uniformly monotonic on $\mathbb{R}_{+}$ in the following sense:
\begin{equation}\left\{
\begin{array}{ll}
 c_{1}\langle y\rangle ^{-k}\leq |\partial_{y}u_{0}^{s}(y) |\leq c_{2}\langle y\rangle ^{-k}, \forall y\geq 0,\\
| \partial_{y}^{ p}u^{s}_{0} (y)|\leq  c_{2}\langle y\rangle ^{-k-p+1}, \forall y\geq 0, 1\leq p\leq m+2,
\end{array}
 \label{w.30}         \right.\end{equation}
for certain $c_{1},c_{2}>0$, they obtained the following existence of almost global-in-time solution.
\begin{Theorem}(\cite{xz1})
Let $m\geq 6$ be an even integer, $k\geq 1, 0\leq \ell\leq \frac{1}{2}, k+\ell>\frac{3}{2}$. Assume that $u^{s}_{0}$ satisfies (\ref{w.30}), the initial data $u_{0}-u_{0}^{s}\in H^{m+3}_{2k+\ell-1}( \mathbb{R}_{+}^{2})$ and $u_{0}-u_{0}^{s}$ satisfies the compatibility condition up to order $m+2$. Then for any $T>0$,  there
exists $\delta_{0}>0$ small enough such that if
\begin{eqnarray*}
 \|(u_{0}-u_{0}^{s})\|_{H^{m+3}_{2k+\ell-1}( \mathbb{R}_{+}^{2}) } \leq \delta_{0},
\label{ }
\end{eqnarray*}
then the initial-boundary value problem (\ref{1.1.6})  admits a unique solution $(u,v)$ with
$$ u -u ^{s}\in L^{\infty} ([0,T];H^{m }_{ k+\ell-1}( \mathbb{R}_{+}^{2}) ),
\partial_{y} (u -u ^{s})\in L^{\infty} ([0,T];H^{m }_{ k+\ell }( \mathbb{R}_{+}^{2}) ),  $$
and
$$v \in L^{\infty} ([0,T];  L^{\infty}(\mathbb{R}_{y,+}; H^{m-1}( \mathbb{R}_{x})) ),
\partial_{y} v \in L^{\infty} ([0,T];H^{m -1}_{ k+\ell-1 }( \mathbb{R}_{+}^{2}) ) .$$
It also has the stability with respect to the initial data as in Theorem \ref{t.xz1}.

\end{Theorem}

 Wu (\cite{W2}) in 2016  studied the  regularity and global well-posedness of classical solutions to the following nonlinear unsteady Prandtl equations \eqref{w.47}
with Robin or Dirichlet boundary condition in half space. Under Oleinik's monotonicity assumption.   Wu (\cite{W1})   in 2015 proved the local well-posedness of the following Prandtl equations with  Robin boundary condition:
\begin{equation}\left\{
\begin{array}{ll}
u_{t}+uu_{x}+vu_{y} =u_{yy}, \\
u_{x}+v_{y}=0,\\
(u_{y}-\beta u)|_{y=0}=0,\ \ v|_{y=0}=0,\\
u(0,x,y)=u_{0}, \\
\lim\limits_{y\rightarrow+\infty}u(t,x,y)=U(t,x).
\end{array}
 \label{w.47}         \right.\end{equation}

The following Prandtl equations with Dirichlet boundary condition takes the form
\begin{equation}\left\{
\begin{array}{ll}
u_{t}+uu_{x}+vu_{y} =u_{yy}, \\
u_{x}+v_{y}=0,\\
u|_{y=0}=0,\ \ v|_{y=0}=0,\\
u(0,x,y)=u_{0}, \\
\lim\limits_{y\rightarrow+\infty}u(t,x,y)=U(t,x).
\end{array}
 \label{w.48}         \right.\end{equation}
The functional norms
$$\|f\|_{\mathcal{H}_{\ell}^{k}( \mathbb{R}^{2}_{+})}=\sum\limits_{|\alpha|+\sigma\leq k}\|(1+y)^{\ell+\sigma}\partial_{t,x}^{\alpha} \partial_{y}^{\alpha}f \|_{L^{2}_{\ell+\sigma}( \mathbb{R}^{2}_{+})}.$$

\begin{Theorem} (\cite{W2})
Considering the nonlinear unsteady Prandtl equations with Robin boundary condition (\ref{w.47}) under Oleinik's monotonicity assumption $w=u_{y}>0$.
Giving any integer $k\geq 4$, $U(t,x)\in C^{k+1}([0,+\infty)\times\mathbb{R})$ and $U(t,x)>0$, we have the following existence, uniqueness and stability results:\\
(1)~ For any fixed finite number $T\in (0,+\infty)$, there exist sufficiently small real number $0<\varepsilon_{1}= o(T^{-1})$ and suitably large real numbers $\ell_{0}>1,\delta_{\beta}>0$ such that if $\ell\geq\ell_{0},\beta\in[\delta_{\beta},+\infty)$, the compatible initial data satisfies
\begin{equation}\left\{
\begin{array}{ll}
w_{0}>0,\ \ u_{0}|_{y=0}>0, \ \  (\partial_{y}u_{0}-\beta u_{0})|_{y=0}=0,\ \ \lim\limits_{y\rightarrow+\infty} u_{0}=U|_{t=0}, \\
\|w_{0} \|_{\mathcal{H}_{\ell}^{k}( \mathbb{R}^{2}_{+})}+\frac{1}{\sqrt{\beta}}\|w_{0}|_{y=0}\|_{\mathcal{H}_{\ell}^{k}( \mathbb{R} )}\leq \varepsilon_{1},\ \
\|U(t,x)\|_{H^{k+1}([0,T]\times \mathbb{R}  )}\leq C_{0}\varepsilon_{1} ,\\
0<c_{1}(1+y)^{-\theta  }\leq w_{0}\leq c_{2}(1+y)^{-\theta  }, \theta>\frac{\ell+1}{2},
\end{array}
 \label{w.49}         \right.\end{equation}
then the Prandtl system (\ref{w.47}) admits a unique classical solution $(w,u,v)$ in  $[0,T]$ satisfying
\begin{equation}
\begin{array}{ll}
  w\in \mathcal{H}_{\ell}^{k}( \mathbb{R}^{2}_{+}), \ \ w, w_{y}\in \mathcal{H}_{\ell}^{k}([0,T]\times \mathbb{R}^{2}_{+}), \\
  u-U\in \mathcal{H}_{\ell-1}^{k}( \mathbb{R}^{2}_{+})\cap \mathcal{H}_{\ell-1}^{k}([0,T]\times \mathbb{R}^{2}_{+}), \\
  \partial_{y}^{j}u|_{y=0}\in H^{k-j}( \mathbb{R}  )\cap H^{k-j}([0,T]\times \mathbb{R}  ), \ \ 0\leq j\leq k,\\
 \partial_{t,x}^{\alpha} v+y\cdot\partial_{t,x}^{\alpha}\partial_{x}U\in L^{\infty}_{y,\ell-1}(L^{2}_{t,x}), \ \ |\alpha|\leq k-1.
\end{array}
 \label{w.50}          \end{equation}
(2)~ The classical solution to (\ref{w.47}) is stable with respect to the initial data in the following sense: for any given two initial data satisfying (\ref{w.49}), then for all $p\leq k-1$, the corresponding solutions of the Prandtl system (\ref{w.48}) satisfy
\begin{eqnarray*}
&& \|u^{1}-u^{2}\|_{\mathcal{H}_{\ell-1}^{p}( \mathbb{R}^{2}_{+})}
+\sum\limits_{j=0}^{p}\|\partial_{y}^{j}u^{1}|_{y=0}-\partial_{y}^{j}u^{2}|_{y=0} \|_{H^{p-j}( \mathbb{R}  ) }\\
&&\quad +\|w^{1}-w^{2}\|_{\mathcal{H}_{\ell }^{p}( \mathbb{R}^{2}_{+})}
 +\sum\limits_{|\alpha|\leq p-1} \|\partial_{t,x}^{\alpha}v^{1}-\partial_{t,x}^{\alpha}v^{2} \|_{ L^{\infty}_{y,\ell-1}(L^{2}_{t,x})} \\
&& \leq C(\varepsilon_{1},T)[ \|w^{1}_{0}-w^{2}_{0}\|_{\mathcal{H}_{\ell }^{p}( \mathbb{R}^{2}_{+})}
+\frac{1}{\beta-\delta_{\beta} } \|w_{0}^{1}|_{y=0}-w_{0}^{2}|_{y=0} \|_{H^{p }( \mathbb{R}  ) }^{2} ] .
\label{ }
\end{eqnarray*}
(3)~ As $\beta\rightarrow+\infty$, $\|u|_{y=0}\|_{H^{k }( \mathbb{R}  )}=\mathcal{O}(\frac{1}{\sqrt{\beta}} )$,
$\|w|_{y=0}\|_{H^{k }( \mathbb{R}  )}=\mathcal{O}( \sqrt{\beta}  )$ and $(w,u,v)$ satisfies the regularities (\ref{w.50}) uniformly.

\end{Theorem}

\begin{Theorem}( \cite{W2})
Considering the nonlinear unsteady Prandtl equations with Dirichlet boundary condition (\ref{w.48}) under Oleinik's monotonicity assumption  $w=u_{y}>0$.
Giving any integer $k\geq 6$, $U(t,x)\in C^{k+1}([0,+\infty)\times\mathbb{R})$ and $U(t,x)>0$, we have the following existence, uniqueness and stability results: \\
(1)~  For any fixed finite number $T\in (0,+\infty)$, there exist sufficiently small real number $0<\varepsilon_{1}= o(T^{-1})$ and suitably large real
 numbers $\ell_{0}>1$  such that if $\ell\geq\ell_{0}$, the compatible initial data and $U(t,x)$ satisfies
\begin{equation}\left\{
\begin{array}{ll}
w_{0}>0,\ \ u_{0}|_{y=0}>0, \ \    \lim\limits_{y\rightarrow+\infty} u_{0}=U|_{t=0}, \\
\|w_{0} \|_{\mathcal{H}_{\ell}^{k}( \mathbb{R}^{2}_{+})} \leq \varepsilon_{2},\ \
\|U(t,x)\|_{H^{k+1}([0,T]\times \mathbb{R}  )}\leq C_{0}\varepsilon_{2} ,\\
0<c_{1}(1+y)^{-\theta  }\leq w_{0}\leq c_{2}(1+y)^{-\theta  }, \theta>\frac{\ell+1}{2},
\end{array}
 \label{w.51}         \right.\end{equation}
then the Prandtl system (\ref{w.48}) admits a unique classical solution $(w,u,v)$ in $[0,T]$ satisfying
\begin{equation}
\begin{array}{ll}
  w\in \mathcal{H}_{\ell}^{k}( \mathbb{R}^{2}_{+}), \ \ w, w_{y}\in \mathcal{H}_{\ell}^{k}([0,T]\times \mathbb{R}^{2}_{+}), \\
  u-U\in \mathcal{H}_{\ell-1}^{k}( \mathbb{R}^{2}_{+})\cap \mathcal{H}_{\ell-1}^{k}([0,T]\times \mathbb{R}^{2}_{+}), \\
  \partial_{y}^{j}u|_{y=0}\in H^{k-j}( \mathbb{R}  )\cap H^{k-j}([0,T]\times \mathbb{R}  ), \ \ 0\leq j\leq k,\\
 \partial_{t,x}^{\alpha} v+y\cdot\partial_{t,x}^{\alpha}\partial_{x}U\in L^{\infty}_{y,\ell-1}(L^{2}_{t,x}), \ \ |\alpha|\leq k-1.
\end{array}
 \label{w.52}          \end{equation}
(2)~  The classical solution to (\ref{w.47}) is stable with respect to the initial data in the following sense: for any given two initial data satisfy (\ref{w.51}), then for all $p\leq k-1$, the corresponding solutions of the Prandtl system (\ref{w.48}) satisfy
\begin{eqnarray*}
&& \|u^{1}-u^{2}\|_{\mathcal{H}_{\ell-1}^{p}( \mathbb{R}^{2}_{+})}
+\sum\limits_{j=0}^{p}\|\partial_{y}^{j}u^{1}|_{y=0}-\partial_{y}^{j}u^{2}|_{y=0} \|_{H^{p-j}( \mathbb{R}  ) }\\
&&\quad +\|w^{1}-w^{2}\|_{\mathcal{H}_{\ell }^{p}( \mathbb{R}^{2}_{+})}
 +\sum\limits_{|\alpha|\leq p-1} \|\partial_{t,x}^{\alpha}v^{1}-\partial_{t,x}^{\alpha}v^{2} \|_{ L^{\infty}_{y,\ell-1}(L^{2}_{t,x})} \\
&& \leq C(\varepsilon_{2},T)  \|w^{1}_{0}-w^{2}_{0}\|_{\mathcal{H}_{\ell }^{p}( \mathbb{R}^{2}_{+})}.
\label{ }
\end{eqnarray*}

\end{Theorem}

 Ignatova and   Vicol (\cite{iv}) in 2016  proved almost global existence of solutions to the following problem (\ref{w.40}) if the Prandtl datum lies within $ \varepsilon$ of the error function erf($\frac{y}{2}$) (in
a suitable topology), then the Prandtl equations have a unique (classical in $x$ weak in $y$) solution on $[0,T_{ \varepsilon}]$, where $T_{ \varepsilon}\geq \exp(\varepsilon^{-1}/\log(\varepsilon^{-1}))$. \\

The two dimensional Prandtl boundary layer equations in the domain $\mathbb { H } = \left\{ ( x , y ) \in \mathbb { R }_{+} ^ { 2 } : y > 0 \right\}$ for the velocity field $(u^{p}, v^{p})$ takes the form
\begin{eqnarray}
\left\{\begin{array}{ll}
 \partial _  t  u ^  P - \partial _  y  ^  2  u ^  P  + u ^  P  \partial _  x  u ^  P  + v ^  P  \partial _  y u ^  P  = - \partial _  x  p ^ E  , \\  \partial _  x  u ^  P + \partial_ y  v ^  P  = 0,
 \end{array}\right.\label{w.40}\end{eqnarray}
with the boundary conditions
\begin{equation*}
\left\{\begin{array}{ll}{ \left. u ^ { P } \right| _ { y = 0 } = \left. v ^ { P } \right| _ { y = 0 } = 0 }, \\
\left. u ^ { P } \right| _ { t = 0 } = u _ { 0 } ^ { P },\\
{ \left. u ^ { P } \right| _ { y = \infty } = u ^ { E } } ,  \end{array}\right.
\end{equation*}
which are obtained by matching the Navier-Stokes no-slip boundary condition $u^{NS}$ on  $\partial \mathbb { H }$, with the Euler slip boundary condition at  $y = \infty$. The trace at $\partial \mathbb { H }$ of the Euler tangential velocity $u^{E}$ obeys Bernoulli's law.

\begin{Theorem}(\cite{iv})(Almost Global Existence). Let the Euler data be given by $u^{E}=\kappa$ and $\partial_{x}p^{E}=0$. Define
\begin{equation*}
u_{0}(x, y)=u_{0}^{P}(x, y)-\kappa \operatorname{erfc}\left(\frac{y}{2}\right),\ \ \operatorname{erf}\left(\frac{y}{2}\right)=\frac{1}{\sqrt{\pi}}\int_{0}^{y}\exp\left(-\frac{z^{2}}{4}\right)dz,
\end{equation*}
where  $\operatorname{erfc}$ is the Gauss error function. There exists a sufficiently large universal constant $C_{*}>0$ and a sufficiently small universal constant $\varepsilon_{*}>0$  such that the following holds. For any given $\varepsilon \in\left(0, \varepsilon_{*}\right]$, assume that there exists an analyticity
radius $\tau_{0}>0$ such that
\begin{equation*}
\frac{C_{*}}{\log \frac{1}{\varepsilon}} \leqq \tau_{0}^{3 / 2} \leqq \frac{1}{C_{*} \varepsilon^{3}},
\end{equation*}
and such that the function
$$g_{0}(x, y)=\partial_{y} u_{0}(x, y)+\frac{y}{2} u_{0}(x, y), $$
obeys
\begin{equation*}
\left\|g_{0}\right\|_{X_{2 \tau_{0}, 1 / 2}}:=\sum_{m \geqq 0}\left\|\exp \left(\frac{y^{2}}{8}\right) \partial_{x}^{m} g_{0}(x, y)\right\|_{L^{2}(\mathbb{H})}\left(2 \tau_{0}\right)^{m} \frac{\sqrt{m+1}}{m !} \leqq \varepsilon.
\end{equation*}
Here $$X_{2 \tau_{0}, 1 / 2}=\{\|g\|_{X_{2 \tau_{0}, 1 / 2}} <+\infty \}. $$
Then there exists a unique solution $u^{p}$ of the Prandtl boundary layer equations on $[0,T_{\varepsilon}]$, where
\begin{equation*}
T_{\varepsilon} \geqq \exp \left(\frac{\varepsilon^{-1}}{\log \left(\varepsilon^{-1}\right)}\right).
\end{equation*}
The solution $u^{p}$ is real analytic in $x$, with analyticity radius larger than  $\tau_{0} / 2$, and lies in a weighted $H^{2}$ space with respect to $y$. We emphasize that $\varepsilon \ll 1$ and $\tau_{0}=\mathcal{O}(1)$ are independent of $\kappa$.

 \end{Theorem}

Greniner,  Guo and Nguyen (\cite{ggn}) in 2016 constructed   the linearized Navier-Stokes equations about generic stationary shear flows of the boundary layer type in a regime of sufficiently large Reynolds number : $R\rightarrow \infty$.
 Greniner,  Guo and  Nguyen (\cite{ggn1}) also recalled the stability of Prandtl boundary layers of \cite{ggn}, and considered
  the linearization of the incompressible Navier-Stokes equations in the domain $\mathbb{R}\times \mathbb{R}_{+}$:
\begin{equation}\left\{
\begin{array}{ll}
v_{t}+ u_{0} \cdot \nabla v+  v\cdot \nabla u_{0}+ \nabla p= \frac{1}{R}\triangle v, \\
 \nabla\cdot v=0,\\
v|_{z=0}=0,\\
u_{0}=(U(z),0)^{T}.
\end{array}
 \label{w.14}         \right.\end{equation}
The asymptotic behavior of these branches $\alpha_{\mathrm{low}}$ and $\alpha_{\mathrm{up}}$ depends on the profile: \\
(i) for plane Poiseuille flow in a channel: $U(z)=1-z^{2}$ for $-1\leq z\leq 1$,
\begin{equation}
\alpha_{\mathrm{low}}(R)\equiv A_{1 c} R^{-1 / 7},  \quad \alpha_{\mathrm{up}}(R)\equiv A_{2 c} R^{-1 / 11},
 \label{w.15}
\end{equation}
(ii) for boundary layer profiles,
\begin{equation}
\alpha_{\mathrm{low}}(R)\equiv A_{1 c} R^{-1 / 4}, \quad \alpha_{\mathrm{up}}(R)\equiv A_{2 c} R^{-1 / 6},
 \label{w.16}
\end{equation}
(iii) for Blasius (a particular boundary layer) profile,
\begin{equation}
\alpha_{\mathrm{low}}(R)=A_{1 c} R^{-1 / 4}, \quad \alpha_{\mathrm{up}}(R)=A_{2 c} R^{-1 / 10},
 \label{w.17}
\end{equation}
where $A_{1 c},A_{2 c}$ are constants$>0$.
\begin{Theorem}(\cite{ggn})
Let $U(z)$ be an arbitrary shear profile with $U'(0)>0$ and satisfy
$$\sup\limits_{z\geq 0}|\partial_{z}^{k}(U(z)- U_{+})e^{\eta_{0}z}|\leq +\infty,\ \ \ \   k=0,\cdot\cdot\cdot,4, $$
for some constants $U_{+}$ and $\eta_{0}>0$. Let $\alpha_{\operatorname{low}}(R)$ and $\alpha_{\mathrm{up}}(R)$ be defined as in (\ref{w.16}) for
general boundary layer profiles, or defined as in (\ref{w.17}) for the Blasius profiles: those with additional assumptions:
$U''(0)=U'''(0)=0$.
Then, there is a critical Reynolds number $R_{c}$ so that for all $R\geq R_{c}$ and all $\alpha\in( \alpha_{\operatorname{low}}(R), \alpha_{\mathrm{up}}(R) )$, there exist a nontrivial triple $c(R), \widehat{v}(z,R), \widehat{p}(z,R)$, with $\operatorname{Im} c(R)>0$,
such that $v_{R}:=e^{i \alpha(y-c t)} \hat{v}(z ; R)$ and $p_{R}:=e^{i \alpha(y-c t)} \hat{p}(z ; R)$ solve the problem (\ref{w.14})
with the no-slip boundary conditions. In the case of instability, there holds the following estimate for the growth rate of the unstable solutions:
$$
\alpha \operatorname{Im} c(R) \approx R^{-1 / 2},
$$
as  $R \rightarrow \infty$.

\end{Theorem}

Xu and Zhang (\cite{xz}) in 2017 discussed the global well-posedness  of solutions for the nonlinear Prandtl boundary layer equations  on the half plane. More specifically, they used energy method to obtain the existence, uniqueness and stability of solutions in weighted Sobolev space $H^m_\mu(\mathbb{R}_+^2)\;( m\geq 6)$ when weighted functions $\mu$ is a polynomial function and proved the lifespan $T$ of solutions could be any large  when its initial datum is a perturbation around the monotonic shear profile of small size $e^{-T}$, which is the first result of  global existence of solutions of 2D Prandtl equation  in a polynomial weighted Sobolev space. In \cite{xz}, they studied the initial-boundary value problem for the Prandtl boundary layer
equation in two dimension, which reads
\begin{eqnarray}
\left\{\begin{array}{l}{\partial_{t} u+u \partial_{x} u+v \partial_{y} u+\partial_{x} p=\partial_{y}^{2} u, \quad t>0,(x, y) \in \mathbb{R}_{+}^{2}}, \\ {\partial_{x} u+\partial_{y} v=0}, \\ {\left.u\right|_{y=0}=\left.v\right|_{y=0}=0, \lim\limits _{y \rightarrow+\infty} u=U(t, x)}=1, \\ {\left.u\right|_{t=0}=u_{0}(x, y)},\end{array}\right.\label{28.1}
\end{eqnarray}
where $\mathbb{R}_{+}^{2} = \{(x, y) \in  \mathbb{R} ^{2}  ; y > 0\}$, $u(t, x, y)$  represents the tangential velocity, $v(t, x, y)$ normal
velocity.

Assume that $u^{s}_{0}$ (initial datum of shear flow) satisfies the following conditions:
\begin{eqnarray}
\left\{\begin{array}{l}{u^{s}_{0}\in C^{m+4}([0,+\infty]),\quad \lim\limits _{y \rightarrow+\infty}u^{s}_{0}(y)=1}, \\ {(\partial_{y}^{2p} u_{0}^{s})(0)=0,\quad 0\leq 2p\leq m+4}, \\ {c_{1}<y>^{-k}\leq(\partial_{y}u_{0}^{s})(y)\leq c_{2}<y>^{-k}, \quad \forall y\geq 0},\\
{|(\partial_{y}^{ p} u_{0}^{s}) (y)|\leq c_{2}<y>^{-k-p+1}, \quad \forall y\geq 0,1\leq p\leq m+4. }\end{array}\right. \label{28.2}
\end{eqnarray}
\begin{Theorem}( \cite{xz})
 Let $m\geq 6$ be an even integer, $ k > 1$ and $-\frac{1}{2}< \nu < 0$. Assume that $u^{s}_{0}$,  the initial data $\bar{ u}_{ 0} =  u_{ 0}-u^{ s}_{0}  \in  H^{m+3}_{k+\nu}(\mathbb{R}_{+}^{2})$, and $\bar{ u}_{ 0}$ satisfies the compatibility condition up to order $m +2$. Then for any $ T > 0$, there exists a $ \delta_{0} > 0$  small enough such that if
\begin{eqnarray}
\|\bar{u}_{0}\|_{ H^{m+1}_{k+\nu}(\mathbb{R}_{+}^{2})}\leq \delta_{0}, \label{28.3}
\end{eqnarray}
then the initial-boundary value problem \eqref{28.1} admits a unique solution $(u, v)$ with
\begin{eqnarray}
(u-u ^{s }) \in  L^{\infty}  ([0,T]; H^{m }_{k+\nu-\delta'}(\mathbb{R}_{+}^{2})), v  \in  L^{\infty}  ([0,T];L^{\infty}(\mathbb{R}_{y,+}; H^{m-1}(\mathbb{R}_{x}))),
\end{eqnarray}
where $\delta'>0$ satisfies $\nu+\frac{1}{2}< \delta'<\nu+1$ and $k +\nu -\delta' >\frac{1}{2}$.

Moreover, we have the stability with respect to the initial data in the following sense: given
any two initial data
\begin{eqnarray*}
u^{ 1}_{0} = u^{s}_{0} + \bar{u}^{ 1}_{0} ,\quad u^{2}_{0} = u^{s}_{0} + \bar{u}^{ 2}_{0} .
\end{eqnarray*}
If $u^{s}_{0}$ satisfies  \eqref{28.2} and $\bar{u}^{ 1}_{0},\bar{u}^{ 2}_{0} $ satisfy \eqref{28.3}, then the solutions $u^{1}$ and $u^{2}$ to problem \eqref{28.1} satisfies
\begin{eqnarray*}
\|u^{1}-u^{ 2} \|_{L^{\infty}  ([0,T]; H^{m -3}_{k+\nu-\delta'}(\mathbb{R}_{+}^{2}))}\leq C
\|u^{1}_{0}-u^{ 2}_{0} \|_{ H^{m +1}_{k+\nu }(\mathbb{R}_{+}^{2})},
\end{eqnarray*}
where the constant $C$ depends on the norm of $\partial_{y}u^{1}, \partial_{y}u^{2}$ in $L^{\infty}([0,T]; H^{m }_{k+\nu-\delta'+1}(\mathbb{R}_{+}^{2}))$.
\end{Theorem}

Paicu and Zhang   (\cite{PZ}) in 2019 proved the global well-posedness of Prandtl system with small
initial data, which is analytical in the tangential variable. The key ingredient used in the
proof is to derive sufficiently fast decay-in-time estimate of some weighted analytic energy
estimate of the tangent velocity, which is based on a Poincar$\acute{e}$ type inequality and a subtle
interplay between the tangential velocity equation and its primitive one. Their result can
be viewed as a global-in-time Cauchy-Kowalevsakya result for Prandtl system with small
analytical data.

One of the key step to rigorously justify this inviscid limit of Navier-Stokes system with Dirichelt boundary condition is to deal with the well-posedness of the following Prandtl system,
\begin{equation}
\left\{\begin{array}{l}
\partial_{t} U+U \partial_{x} U+V \partial_{y} U-\partial_{y}^{2} U+\partial_{x} p=0, \quad(t, x, y) \in \mathbb{R}_{+} \times \mathbb{R} \times \mathbb{R}_{+}, \\
\partial_{x} U+\partial_{y} V=0, \\
\left.U\right|_{y=0}=\left.V\right|_{y=0}=0 ,\quad \text { and } \quad \lim\limits _{y \rightarrow+\infty} U(t, x, y)=w(t, x), \\
\left.U\right|_{t=0}=U_{0},
\end{array}\right. \label{y2.4}
\end{equation}
where $U$ and $V$ represent the tangential and normal velocities of the boundary layer flow. $(w(t, x), p(t, x))$ are the traces of the tangential velocity and pressure of the outflow on the boundary, which satisfy Bernoulli's law:
$$
\partial_{t} w+w \partial_{x} w+\partial_{x} p=0.
$$

The authors of \cite{PZ} investigated the global existence of the solutions to the Prandtl system (\ref{y2.4}) with small data which is analytic in the tangential variable. For simplicity, they took $w(t, x)$ in (\ref{y2.4}) to be $\varepsilon f(t)$ with $f(0)=0.$ Let $R>0$ and $\chi \in C^{\infty}[0, \infty)$ with $\chi(y)=\left\{\begin{array}{ll}1 & \text { if } y \geq 2 R, \\ 0 & \text { if } y \leq R,\end{array}\right.$ they denoted $W \stackrel{\text { def }}{=} U-\varepsilon f(t) \chi(y) .$ Then $W$ solves
$$
\left\{\begin{array}{l}
\partial_{t} W+(W+\varepsilon f(t) \chi(y)) \partial_{x} W+V \partial_{y}(W+\varepsilon f(t) \chi(y))-\partial_{y}^{2} W=\varepsilon g, \\
\partial_{x} W+\partial_{y} V=0, \quad(t, x, y) \in \mathbb{R}_{+} \times \mathbb{R}_{+}^{2} ,\\
\left.W\right|_{y=0}=\left.V\right|_{y=0}=0, \quad \text { and } \quad \lim\limits _{y \rightarrow+\infty} W(t, x, y)=0, \\
\left.W\right|_{t=0}=U_{0},
\end{array}\right.
$$
where $\mathbb{R}_{+}^{2}=\mathbb{R} \times \mathbb{R}_{+}$ and $g(t, y) \stackrel{\text { def }}{=}(1-\chi(y)) f^{\prime}(t)+f(t) \chi^{\prime \prime}(y)$.

In order to get rid of the source term in the equation of  $W$, the authors of \cite{PZ} introduced $u^{s}$ via
\begin{equation}
\left\{\begin{array}{l}
\partial_{t} u^{s}-\partial_{y}^{2} u^{s}=\varepsilon g(t, y), \quad(t, y) \in \mathbb{R}_{+} \times \mathbb{R}_{+}, \\
\left.u^{s}\right|_{y=0}=0 ,\quad \text { and } \quad \lim\limits _{y \rightarrow+\infty} u^{s}=0, \\
\left.u^{s}\right|_{t=0}=0.
\end{array}\right. \label{y2.5}
\end{equation}
With $u^{s}$ being determined by the last equation,  the authors set $u \stackrel{\text { def }}{=} W-u^{s}$ and $v \stackrel{\text { def }}{=} V$. Then $(u, v)$ verifies
\begin{equation}
\left\{\begin{array}{l}
\partial_{t} u+\left(u+u^{s}+\varepsilon f(t) \chi(y)\right) \partial_{x} u+v \partial_{y}\left(u+u^{s}+\varepsilon f(t) \chi(y)\right)-\partial_{y}^{2} u=0, \\
\partial_{x} u+\partial_{y} v=0, \quad(t, x, y) \in \mathbb{R}_{+} \times \mathbb{R}_{+}^{2}, \\
\left.u\right|_{y=0}=\left.v\right|_{y=0}=0, \quad \text { and } \quad \lim\limits _{y \rightarrow \infty} u(t, x, y)=0, \\
\left.u\right|_{t=0}=u_{0} \stackrel{\text { def }}{=} U_{0} .
\end{array}\right.    \label{yy2.5}
\end{equation}

Let $s$ in $\mathbb{R}$. For $u$ in $\mathcal{S}_{h}^{\prime}\left(\mathbb{R}_{+}^{2}\right)$, which means that $u$ is in $\mathcal{S}^{\prime}\left(\mathbb{R}_{+}^{2}\right)$ and satisfies $\lim\limits _{k \rightarrow-\infty}\left\|S_{k}^{\mathrm{h}} u\right\|_{L^{\infty}}=0,$  the authors set
$$
\|u\|_{\mathcal{B}^{s, 0}} \stackrel{\text { def }}{=}\left\|\left(2^{k s}\left\|\Delta_{k}^{\mathrm{h}} u\right\|_{L_{+}^{2}}\right)_{k \in \mathbb{Z}}\right\|_{\ell^{1}(\mathbb{Z})}  .
$$

For $s \leq \frac{1}{2}$, they defined $\mathcal{B}^{s, 0}\left(\mathbb{R}_{+}^{2}\right) \stackrel{\text { def }}{=}\left\{u \in \mathcal{S}_{h}^{\prime}\left(\mathbb{R}_{+}^{2}\right) \mid\|u\|_{\mathcal{B}^{s, 0}}<\infty\right\}$.

If $\ell$ is a positive integer and if $\ell-\frac{1}{2}<s \leq \ell+\frac{1}{2}$, then they defined $\mathcal{B}^{s, 0}\left(\mathbb{R}_{+}^{2}\right)$ as the subset of distributions $u$ in $\mathcal{S}_{h}^{\prime}\left(\mathbb{R}_{+}^{2}\right)$ such that $\partial_{x}^{\ell} u$ belongs to $\mathcal{B}^{s-\ell, 0}\left(\mathbb{R}_{+}^{2}\right)$.
In order to obtain a better description of the regularizing effect of the transport-diffusion equation, they needed to use Chemin-Lerner type spaces $\widetilde{L}_{T}^{p}\left(\mathcal{B}^{s, 0}\left(\mathbb{R}_{+}^{2}\right)\right)$.

Let $p \in[1,+\infty]$ and $T_{0}, T \in[0,+\infty]$.  The authors defined $\widetilde{L}^{p}\left(T_{0}, T ; \mathcal{B}^{s, 0}\left(\mathbb{R}_{+}^{2}\right)\right)$ as the
completion of $C\left(\left[T_{0}, T\right] ; \mathcal{S}\left(\mathbb{R}_{+}^{2}\right)\right)$ by the norm
$$
\|a\|_{\tilde{L}^{p}\left(T_{0}, T ; \mathcal{B}^{s, 0}\right)} \stackrel{\text { def }}{=} \sum_{k \in \mathbb{Z}} 2^{k s}\left(\int_{T_{0}}^{T}\left\|\Delta_{k}^{\mathrm{h}} a(t)\right\|_{L_{+}^{2}}^{p} d t\right)^{\frac{1}{p}}, \ \
\Delta_{k}^{\mathrm{h}} a=\mathcal{F} ^{-1}(\varphi(2^{-k}|\xi|)\widehat{a}),
$$
$\varphi$ decays to zero sufficiently fast as $y$ approaching to $+\infty$, $\widehat{a}$  denotes the partial Fourier transform of the distribution $a$, with the usual change if $p=\infty$. In particular, when $T_{0}=0,$  the authors denoted  $\|a\|_{\widetilde{L}_{T}^{p}\left(\mathcal{B}^{s, 0}\right)} \stackrel{\text { def }}{=}$ $\|a\|_{\tilde{L}^{p}\left(0, T ; \mathcal{B}^{s, 0}\right)}$ for simplicity.

\begin{Theorem}(\cite{PZ} )
Let $\ell>3$ be an integer, $\delta>0$, and $f \in H^{1}\left(\mathbb{R}_{+}\right)$ which verifies
$$
\mathcal{C}_{f, \ell}^{2} \stackrel{\text { def }}{=} \int_{0}^{\infty}\langle t\rangle^{\frac{3}{2}}(1+\log \langle t\rangle)^{\ell}\left(f^{2}(t)+\left(f^{\prime}(t)\right)^{2}\right) d t<\infty   .
$$
Let $u_{0}=\partial_{y} \varphi_{0}$ satisfy $u_{0}(x, 0)=0$ and
$$
\mathfrak{C}\left(u_{0}\right) \stackrel{\text { def }}{=}\left\|e^{\frac{y^{2}}{8}} e^{\delta\left|D_{x}\right|}\left(\varphi_{0}, u_{0}\right)\right\|_{\mathcal{B}^{\frac{1}{2}}, 0} \leq c_{0}   ,
$$
for some $c_{0}$ sufficiently small. Then there exists $\varepsilon_{0}>0$ so that for $\varepsilon \leq \varepsilon_{0}$, the system (\ref{yy2.5}) has a unique global solution $u$ which satisfies
$$
\left\|e^{\frac{y^{2}}{8(t)}} e^{\frac{\delta}{2}\left|D_{x}\right|} u\right\|_{\tilde{L}^{\infty}\left(\mathbb{R}_{+} ; \mathcal{B}^{\frac{1}{2}, 0}\right)}+\left\|e^{\frac{y^{2}}{8(t)}} e^{\frac{\delta}{2}\left|D_{x}\right|} \partial_{y} u\right\|_{\tilde{L}^{2}\left(\mathbb{R}_{+} ; \mathcal{B}^{\frac{1}{2}, 0}\right)} \leq C\left\|e^{\delta\left|D_{x}\right|} u_{0}\right\|_{\mathcal{B}^{\frac{1}{2}}, 0} .
$$
Moreover, for $u^{s}$ being determined by (\ref{y2.5}) and for any integer $j$, there exists a positive constant $C_{j}$ so that
$$
\begin{array}{l}
\left\|e^{\frac{y^{2}}{8(t)}} e^{\frac{\delta}{2}\left|D_{x}\right|} u(t)\right\|_{\mathcal{B}^{\frac{1}{2}}, 0} \leq C_{j} \mathfrak{C}\left(u_{0}\right)\langle t\rangle^{-\frac{3}{4}} \log ^{-\frac{j}{2}}\langle t\rangle \text { and } \\
\int_{0}^{\infty}\langle t\rangle^{\frac{1}{4}}\left(\left\|e^{\frac{y^{2}}{8(t)}} \partial_{y} u^{s}(t)\right\|_{L_{\mathbf{v}}^{2}}+\left\|e^{\frac{y^{2}}{8(t)}} e^{\frac{\delta}{2}\left|D_{x}\right|} \partial_{y} u(t)\right\|_{\mathcal{B}^{\frac{1}{2}}, 0}\right) d t \leq C\left(\mathcal{C}_{f, \ell} \varepsilon+\mathfrak{C}\left(u_{0}\right)\right)        .
\end{array}
$$
Here $$e^{\frac{\delta}{2}\left|D_{x}\right|} u(t)=\mathcal{F}^{-1}_{\xi\rightarrow x}\left(e^{(\delta-\lambda\theta(t))|\xi|}   \hat{u}(t,\xi,y)\right).  $$
\end{Theorem}

Petrov and   Suslina (\cite{PS}) in 2020 obtained theorem about existence and uniqueness of the solution in the classes $\widetilde{H}^{s}(-1,1)$ with $0 \leqslant s \leqslant 1$. In particular, for $s=1$ the result is as follows: if $r^{1 / 2} f \in L^{2}$, then $r^{-1 / 2} u, r^{1 / 2} u^{\prime} \in L_{2},$ where $r(x)=1-x^{2}$.

The Prandtl equation
\begin{equation}
u(x)-p(x) \frac{1}{2 \pi} \int_{-1}^{1} \frac{u^{\prime}(t)}{t-x} d t=p(x) f(x), \quad u(-1)=u(1)=0, \label{PS.1}
\end{equation}
is one of the universal equations of mathematical physics.
It is convenient to put $V(x):=p(x)^{-1}$ and rewrite the first Prandtl equation (\ref{PS.1}) with the Dirichlet conditions as
\begin{equation}
V(x) u(x)-\frac{1}{2 \pi} \int_{-1}^{1} \frac{u^{\prime}(t)}{t-x} d t=f(x), \quad u(-1)=u(1)=0.
\end{equation}
Here and below, singular integrals are understood in the mean value sense. The authors of   \cite{PS}  assumed that $V(x)$ is a measurable function on (-1,1) satisfying the following conditions:
$$
V(x) \geqslant 0 \text { for a. e. } x \in(-1,1) ; \quad\left(1-x^{2}\right) V(x) \leqslant M<\infty.
$$
Any function $u \in H_{00}^{1 / 2}(-1,1)$ satisfies
$$
\int_{-1}^{1} \frac{|u(x)|^{2}}{1-x^{2}} d x \leqslant \frac{1}{2} \int_{-\infty}^{\infty}|\zeta \| \widehat{u}(\zeta)|^{2} d \zeta,
$$
where the norm of $ H_{00}^{1 / 2}(-1,1)$  is defined as the standard norm in $ H ^{1 / 2}(\mathbb{R})$:
$$
\|u\|_{ H_{00}^{1 / 2}}^{2} :=      \|u\|_{ H_{00}^{1 / 2}(-1,1)}^{2}=\frac{1}{2}\int_{-\infty}^{+\infty}(1+|\xi|^{2})^{\frac{1}{2}}
|\widehat{u}(\zeta)|^{s}d\zeta .$$

\begin{Theorem}(\cite{PS})
Suppose that $V(x)$ satisfies the assumptions. Then for any $f \in\left(H_{00}^{1 / 2}\right)^{*}$, there exists a unique weak solution $u \in H_{00}^{1 / 2}$ of problem (\ref{PS.1}). The solution satisfies the following estimate:
$$
\|u\|_{H_{00}^{1 / 2}} \leqslant(2 \pi+2)\|f\|_{\left(H_{00}^{1 / 2}\right)^{*}}.
$$
\end{Theorem}

  Dalibard, Paddick (\cite{DP}) in 2021 proved the global existence of solutions of the following stationary Prandtl-like equations by using the von Mises variables:
 \begin{equation}\label{w.65}
\left\{\begin{array}{ll}
\lambda_{0}\nu(y)( u\partial_{\xi}v+v\partial_{y} v)+\psi-\nu(y)^{2}\partial_{\xi}^{2}v=\psi^{0}(y),\\
(u,v)=\nabla^{\bot} \psi:=(-\partial_{y}\psi,\partial_{\xi}\psi), \\
 v|_{y=0}=v_{0}(\xi),  \\
 (u,v)|_{\xi=0}=0,\\
 \lim\limits_{\xi \rightarrow+\infty}\psi(\xi,y)=\psi^{0}(y),
  \end{array}\right.
\end{equation}
in the domain $\Omega_{Y}=\{(\xi,y)\in \mathbb{R}^{2}|\xi>0,0<y<Y\}$. The function $\psi^{0}$ is given and smooth, $\lambda_{0}$ is a strictly positive parameter, and the function $\nu$ is a smooth function with $\nu\leq 1$.

The von Mises change of variables,
$$(\xi ,y)\mapsto (\psi(\xi,y), y), $$
and define a new unknown function, $w=w(\psi, y)$ such that $w(\psi, y)=v^{2}(\psi(\xi,y), y )$. The chain rule allows us to express the derivatives:
\begin{eqnarray*}
&&\frac{\partial v}{ \partial\xi} =\frac{1}{2}\frac{\partial w}{ \partial\psi}, \ \
\frac{\partial^{2} v}{ \partial\xi^{2}} =\frac{\sqrt{w}}{2}\frac{\partial^{2} w}{ \partial\psi^{2}}, \\
&& \frac{\partial v}{ \partial y} =\frac{1}{2\sqrt{w}}\frac{\partial w}{ \partial y}+\frac{1}{2\sqrt{w}}\frac{\partial w}{ \partial\psi}\frac{\partial\psi}{ \partial y}   .
\end{eqnarray*}
Then (\ref{w.65}) turns to
 \begin{equation}\label{w.66}
\left\{\begin{array}{ll}
\lambda_{0}\nu(y) \partial_{y}w- \nu(y)^{2}\sqrt{w} \partial_{\psi}^{2}w=2(\psi^{0}(y)-\psi),\\
w|_{y=0}=w_{0}(\psi),  \\
 w|_{\psi=0}=w|_{\psi=\psi^{0}(y)}=0,
  \end{array}\right.
\end{equation}
in the domain $\Omega =\{(\psi,y)\in \mathbb{R}^{2}|0<\psi<\psi^{0}(y),0<y<Y\}$.

\begin{Theorem}  (\cite{DP}  )
Let $\psi^{0}\in \mathcal{C}^{2}([0,Y])$ such that $\inf \psi^{0}>0$, $\inf\psi_{0}'>0$. Assume that $\lambda_{0}(\psi^{0})'(y)/\nu(y)^{1/3}>\eta+9/2>\sqrt[3]{\frac{27}{4}}$ for every $y\geq 0$ and for some $\eta>0$. Let $v_{0}\in W^{2,\infty}(\mathbb{R}^{+})$, positive on $(0,+\infty)$, with $v(0)=v_{0},v_{0}'>0$, such that $\int_{0}^{+\infty}v_{0}(\xi)d \xi=\psi^{0}(0)$. We assume that $v_{0}$ satisfies the corner compatibility condition
$$\nu(0)^{2}v_{0}''(\xi)+\psi^{0}(0)=\mathcal{O}(\xi^{2})$$
as $\xi\rightarrow 0$, and the decay constraint
$$v_{0}(\xi)\sim a(0)(\psi^{0}(0)-\psi_{0}(\xi))$$
as $\xi\rightarrow+\infty$, in which $\psi_{0}(\xi)=\int_{0}^{\xi} v_{0}(s)ds$ is the initial stream function.

Then, for every $Y>0$, there exists a unique classical Lipschitz solution to the steady rotating Prandtl equations (\ref{w.65})
with $\nu>0$ in $\Omega_{Y}$ . Furthermore, $v$ satisfies the following properties:\\
(i) Behaviour close to $\xi=0$: for some $\xi_{0}, m>0$, depending on $Y$, we have $\partial_{\xi}v>m$ for $\xi\leq \xi_{0}$. \\
(ii) Behaviour close to $\xi=\infty$: for every $y, v(\xi,y)\sim a(y)(\psi^{0}(y)-\psi(\xi,y))$ as $\xi\rightarrow \infty$.

\end{Theorem}

 Wang,  Wang and   Zhang (\cite{WWZ1}) in 2021 proved the global existence and large time behavior of small solutions to 2D Prandtl system
   (\ref{1.1.6}) which $U$ is a constant for data with Gevrey 2 regularity in the $x$ variable and Sobolev regularity in the $y$ variable. They extended the global well-posedness result
in \cite{PZ} for 2D Prandtl system with analytic data to data with optimal Gevery regularity in
the sense of \cite{GD}. \\

In  \cite{WWZ1}, there exists a potential function $\varphi$  such that $u= \partial_{y} \varphi$ and $v= -\partial_{x} \varphi$. Recalled the  Gevrey form denoted by
$$G:=u+\frac{y}{2\langle t\rangle}\varphi, \ \   g:=\partial_{y}G= \partial_{y}g+\frac{y}{2\langle t\rangle} u+\frac{\varphi}{2\langle t\rangle},$$
here $\langle t\rangle=(1+t^{2})^{\frac{1}{2}}$, and denoted  by
\begin{eqnarray}
&&e^{\Phi (t,D_{x})}:=e^{\delta(t)\langle D_{x}\rangle^{\frac{1}{2} }} \ \ \text{with} \ \ \delta(t)=\delta-\lambda\theta(t)\ \  \text{and} \ \
f_{\Phi}:=e^{\Phi (t,D_{x})} f, \label{w.74}\\
&& \|f\|_{H^{s,k}_{\Psi}}:=\sum\limits_{0\leq \ell\leq k}\left(\int_{0}^{\infty}e^{2\Phi (t,x)} \|\partial_{y}^{\ell} f(\cdot, y)\|^{2}_{H^{s}_{h}} dy \right)^{ \frac{1}{2}} \ \  \text{with} \ \ \Phi(t,y)=\frac{y^{2}}{8\langle t\rangle}.
 \end{eqnarray}
\begin{Theorem}  (\cite{WWZ1})
  Let $u_{0}= \partial_{y} \varphi_{0}$ satisfy $u_{0}(x,0)=0$ and $\int_{0}^{\infty}u_{0}dy=0$. Let $G_{0}=u_{0}+\frac{y}{2\langle t\rangle}\varphi_{0}$.
For some sufficiently small but fixed $\eta\in (0,\frac{1}{6})$, we denote
\begin{eqnarray*}
\begin{array}{ll}
E(t):=&\|\langle t\rangle^{\frac{1-\eta}{4}} u_{\Phi} \|^{2}_{H^{\frac{11}{2},0}_{\Psi} }
+\|\langle t\rangle^{\frac{3-\eta}{4}} \partial_{y} u_{\Phi} \|^{2}_{H^{5,0}_{\Psi} }
+\|\langle t\rangle^{\frac{5-\eta}{4}} G_{\Phi} \|^{2}_{H^{4,0}_{\Psi} } \\
&+\|\langle t\rangle^{\frac{7-\eta}{4}} \partial_{y} G_{\Phi} \|^{2}_{H^{3,0}_{\Psi} }
+\|\langle t\rangle^{\frac{9-\eta}{4}} \partial_{y}^{2} G   \|^{2}_{H^{3,0}_{\Psi} }
+\|\langle t\rangle^{\frac{11-\eta}{4}} \partial_{y}^{3} G  \|^{2}_{H^{2,0}_{\Psi} } ,
\end{array} \end{eqnarray*}
and
\begin{eqnarray*}
\begin{array}{ll}
D(t):=&\|\langle t\rangle^{\frac{1-\eta}{4}} \partial_{y} u_{\Phi} \|^{2}_{H^{\frac{11}{2},0}_{\Psi} }
+\|\langle t\rangle^{\frac{3-\eta}{4}} \partial_{y}^{2} u_{\Phi} \|^{2}_{H^{5,0}_{\Psi} }
+\|\langle t\rangle^{\frac{5-\eta}{4}}\partial_{y} G_{\Phi} \|^{2}_{H^{4,0}_{\Psi} } \\
&+\|\langle t\rangle^{\frac{7-\eta}{4}} \partial_{y}^{2} G_{\Phi} \|^{2}_{H^{3,0}_{\Psi} }
+\|\langle t\rangle^{\frac{9-\eta}{4}} \partial_{y}^{3} G  \|^{2}_{H^{3,0}_{\Psi} }
+\|\langle t\rangle^{\frac{11-\eta}{4}} \partial_{y}^{4} G  \|^{2}_{H^{2,0}_{\Psi} } .
\end{array} \end{eqnarray*}
 Moreover, assume that
$$\| u_{\Phi} (0)\|^{2}_{H^{\frac{25}{4},0}_{\Psi} }+E(0)\leq \varepsilon ^{2}. $$
Then there exists $\varepsilon_{0},\lambda_{0}>0$ so that for $\varepsilon\leq\varepsilon_{0} $ and $\lambda\geq \lambda_{0},\theta(t)$ in (\ref{w.74})
 satisfies $\sup\limits_{t\in [0,\infty)}\theta(t)\leq \frac{\delta}{4\lambda}$, and the system (\ref{1.1.6})  has a unique global solution $u$ which satisfies
$$E(t)+c\eta\int_{0}^{t}D(t')  dt'\leq C \varepsilon ^{2}, \ \ \forall t\in \mathbb{R}. $$
\end{Theorem}

\subsection{Compressibl 2D Prandtl Equations}

  Liu and  Wang (\cite{LW2}) in 2014 studied the asymptotic behavior of circularly symmetric solutions to the following  initial boundary value problem (\ref{w.8}) of the nonisentropic compressible Navier-Stokes equations in a two-dimensional exterior domain with impermeable boundary conditions when the viscosities and the heat conduction coefficient tend to zero at the same rate. The compressible viscous flow in two space variables:
\begin{equation}\left\{
\begin{array}{ll}
\partial_{t}\rho^{\epsilon}+\nabla\cdot(\rho^{\epsilon} {\bf u}^{\epsilon} )=0,\\
\rho^{\epsilon} \{ \partial_{t} {\bf u}^{\epsilon} +({\bf u}^{\epsilon}\cdot\nabla){\bf u}^{\epsilon} \}+\nabla p(\rho^{\epsilon}, \theta^{\epsilon} )=
\epsilon \{  \mu\triangle {\bf u}^{\epsilon}+(\lambda+\mu)\nabla(\nabla\cdot {\bf u}^{\epsilon})\} ,\\
c_{V} \rho^{\epsilon} \{ \partial_{t}\theta^{\epsilon} +({\bf u}^{\epsilon}\cdot\nabla)\theta^{\epsilon} \}+  p(\rho^{\epsilon}, \theta^{\epsilon} )\nabla\cdot {\bf u}^{\epsilon}=
\epsilon \{ k\triangle \theta^{\epsilon}+\lambda(\nabla \cdot {\bf u}^{\epsilon})^{2}+ 2\mu D\cdot D\}  .
\end{array}
\label{w.8}      \right.\end{equation}
where $ \rho^{\epsilon},\theta^{\epsilon}, {\bf u}^{\epsilon}=(u ^{\epsilon}, v^{\epsilon} )^{T} ,  p(\rho^{\epsilon}, \theta^{\epsilon} )$ are the density, the absolute temperature, velocity, and the  pressure  respectively,
$$D_{ij}=\frac{1}{2}\left(\frac{\partial_{u_{i}^{\epsilon}} }{\partial_{x_{j}}}+ \frac{\partial_{u_{j}^{\epsilon}} }{\partial_{x_{i}}}\right),  D\cdot D=\sum\limits_{i,j=2}D^{2}_{ij}. $$
Take the circularly symmetric solutions of the form
$$ \rho^{\epsilon}(t,x)=\widetilde{\rho}^{\epsilon}(t,r),\ \   {\bf u}^{\epsilon}(t,x)=(  \widetilde{u}^{\epsilon},  \widetilde{v}^{\epsilon})(t,r),\ \
\theta^{\epsilon}(t,x) =  \widetilde{\theta} ^{\epsilon}(t,r), $$
and let $$U^{\epsilon}(t,x)=(\frac{1}{\widetilde{\rho}^{\epsilon}}, \widetilde{u}^{\epsilon},  \widetilde{v}^{\epsilon}, \widetilde{\theta} ^{\epsilon}    )^{T}(t,r(t,x)),$$
then $U^{\epsilon}(t,x)$ satisfies the following problem
\begin{equation}\left\{
\begin{array}{ll}
  \partial_{t}U^{\epsilon}+A(U^{\epsilon}) \partial_{x}U^{\epsilon}+Q_{1} (U^{\epsilon}) (U^{\epsilon},U^{\epsilon})-\epsilon B  (U^{\epsilon}) \partial_{x}^{2} U^{\epsilon} \\
- \epsilon Q_{2} (U^{\epsilon}) (\partial_{x}U^{\epsilon},\partial_{x}U^{\epsilon})-\epsilon Q_{3} (U^{\epsilon}) (U^{\epsilon}, U^{\epsilon})-\epsilon V(\partial_{x}U^{\epsilon}, U^{\epsilon})=0, \\
U^{\epsilon}  (0,x )=U_{0}(x)=) ( \widetilde{\tau}^{\epsilon}_{0} , \widetilde{u}^{\epsilon}_{0},  \widetilde{v}^{\epsilon}_{0}, \widetilde{\theta} ^{\epsilon} _{0}   )^{T}(x), \ \    \tau=\frac{1}{\rho}, \\
U^{\epsilon} _{1}  (t,0) =(0,v^{0}(t),\theta^{0})^{T} .
\end{array}
\label{w.9}      \right.\end{equation}

\begin{Theorem}(\cite{LW2})
Let the initial data $U_{0}(x)=(\tau_{0},u_{0},v_{0},\theta_{0})^{T}(x)$ satisfy Assumption 2, and $U_{0}-(\overline{\tau},0,0,\overline{\theta} )^{T}\in H^{s}(\mathbb{R}_{+})$ with $s>11$. Let $U^{a}(t,x)$ be an approximate solution to the problem (\ref{w.9}).
Then there exists $\delta>0$ and $C>0$, independent of $\epsilon$,  such that when
$$ \sup\limits_{t\in[0,T*],\beta\leq 2}\|\partial_{\xi}^{\beta}U^{B,0}(t,\cdot)\|_{L^{\infty}}+ \sup\limits_{t\in[0,T*] } \|\partial_{\xi} U^{B,0}(t,\cdot)\|_{L^{1}}\leq \delta, $$
holds, there exists a unique solution $U^{\epsilon}(t,x) $ to (\ref{w.9}) such that $U^{\epsilon}-U^{a}\in C([0, T*];H^{1})$ and
$$ \|U^{\epsilon}-U^{a}\|_{L^{\infty}([0,T*];L^{2})}^{2}+\epsilon^{2}\|\partial_{x}U^{\epsilon}-\partial_{x}U^{a}\|_{L^{\infty}([0,T*];L^{2})}^{2}\leq C \epsilon^{2}.$$
Moreover, we have
$$ U^{\epsilon}(t,x) -\left(U^{I,0}(t,x) +U^{B,0}(t,\frac{x}{\sqrt{\epsilon}})=\mathbb{O}(\sqrt{\epsilon} )  \right). $$

\end{Theorem}

  Wang,  Xie and  Yang (\cite{WXY}) in 2015   obtained the well-posedness of the  compressible Prandtl equations locally in time in the domain $\{(x,y)| x\in \mathbb{T}=\mathbb{R}/\mathbb{Z}, y\in[0, +\infty)  \}$
\begin{equation}\left\{
\begin{array}{ll}
u_{t}+ u u_{x}+  v u_{y}= \frac{1}{\rho(t,x)}\nu u_{yy}+p_{x}, \\
 \partial_{x}(\rho u)+ \partial_{y}(\rho v)= -\rho_{t},\\
u(0,x,y)=u_{0},\ \ u(t,x,0)=v(t,x,0)=0,\\
\lim\limits_{y\rightarrow+\infty}u(t,x,y)=U(t,x),
\end{array}
 \label{1.1.7}         \right.\end{equation}
where $U(t,x)$  satisfies Bernoulli's law (\ref{1.1.3}).

 Main assumptions (H) on the initial data are the following:\\
(H1) For a fixed integer $k_{0}>9$, the initial data $u_{0}(x,\eta)$ satisfy the compatibility
condition of the problem (\ref{1.1.7}) up to order $4k_{0}+2$;\\
(H2) monotone condition $\partial_{\eta} u_{0}(x,\eta)\geq\frac{\sigma_{0}}{(1+\eta)^{\gamma+2}}>0 $ holds for all $x\in \mathbb{T}$ and $\eta\geq 0$ with some positive constant $\sigma_{0}$ and a positive integer $\gamma\geq 2$;\\
(H3) $\| (1+\eta)^{\gamma+\alpha_{2}} D^{\alpha}(u_{0}(x,\eta)  -U(0.x))\|_{L^{2}(\mathbb{T}\times\mathbb{R}^{+})}\leq C_{0}$, where $D^{\alpha}=\partial_{x}^{\alpha_{1}}\partial_{y}^{\alpha_{2}}$
with
$\alpha=(\alpha_{1},\alpha_{2})$ and $|\alpha|= \alpha_{1}+\alpha_{2} \leq 4k_{0}+2$;\\
(H4) $\| (1+\eta)^{\gamma+2+\alpha_{2}} D^{\alpha}\partial_{\eta} u_{0} \|_{L^{\infty}(\mathbb{T}\times\mathbb{R}^{+})}\leq \frac{1}{\sigma_{0}}$  for  $|\alpha| \leq 3k_{0} $.

 \begin{Theorem}(\cite{WXY})
Suppose that the outer Euler flow is smooth for $0\leq t\leq T_{0}$, the density $\rho(t,x)$ has both positive lower and upper bounds, and the Sobolev norm
$H^{s}([0,T_{0}\times\mathbb{R}_{+}])$ of $(\rho,U,V)$ is bounded for a suitably large integer $s$. Moreover, the main assumptions (H) on the initial
data $u_{0}(x,y)$ are satisfied. Then there exists a $0<T\leq T_{0}$, such that the initial boundary-value problem (\ref{1.1.7}) has a unique classical
solution $(u,v)$ satisfying
$$ \sum\limits_{|m_{1}|+[(m_{2}+1)/2]\leq k_{0}} \|\langle\eta\rangle^{l}\partial_{t,x}^{m_{1}} \partial_{\eta}^{m_{2}} (u-U)\|_{L^{2}([0,T]
\times \mathbb{T}\times\mathbb{R}_{+})} < +\infty$$
for a fixed $l>\frac{1}{2}$ depending only on $\gamma$ given in (H) with $\langle\eta\rangle=(1+\eta)$, and
$$ \sum\limits_{|m_{1}|+[(m_{2}+1)/2]\leq k_{0}-1} \sup\limits_{\eta\in\mathbb{R}_{+}} \| \partial_{t,x}^{m_{1}} \partial_{\eta}^{m_{2}}
 (v-Vy)(\cdot,y)\|_{L^{2}([0,T]\times \mathbb{T} )} <+\infty. $$
Here
 $V(t,x)$ denotes the trace of $\partial_{y} u^{E}_{2}$ on $y=0$ for the normal velocity $u^{E}_{2}$ of outer Euler flow.
 \end{Theorem}

  Gong, Guo and   Wang (\cite{GGW5}) in 2016  studied the well-posedness for compressible steady boundary layer equations in the domain $D=\{0<t<T,0<x<X,y>0\}$ by the Crocco transformation or von Mises transformation. \\

The  unsteady boundary layer equations are  governed by
\begin{equation}\left\{
\begin{array}{ll}
\partial_{x}(\rho u)+\partial_{y}(\rho v)=-\partial_{t} \rho,\\
\partial_{t}u+ u\partial_{x} u +v\partial_{y} u  =\frac{1}{\rho}\partial_{y}^{2}u-\frac{\partial_{x} p(\rho)}{\rho},
\end{array}
 \label{w.41}         \right.\end{equation}
 associated steady boundary layer equations are  governed by
\begin{equation}\left\{
\begin{array}{ll}
\partial_{x}(\rho u)+\partial_{y}(\rho v)=0,\\
 u\partial_{x} u +v\partial_{y} u  =\frac{1}{\rho}\partial_{y}^{2}u-\frac{\partial_{x} p(\rho)}{\rho},
\end{array}
 \label{w.42}         \right.\end{equation}
with the following boundary and initial conditions
 \begin{equation}\left\{
\begin{array}{ll}
u_{x=0}=u_{1}(t,y),\\
u|_{y=0}=0,\ \ v_{y=0}=v_{0}(t,x),\\
\lim\limits_{y\rightarrow+\infty}u(t,x,y)=U(t,x),
\end{array}
 \label{w.43}         \right.\end{equation}
and
 \begin{eqnarray}
u(0,x,y)=u_{0}(x,y).
\label{w.44}
\end{eqnarray}

 \begin{Theorem}(\cite{GGW5})
 For the steady boundary value problem (\ref{w.42})-(\ref{w.43}), assume that $\rho(x)$ and $v_{0}(x)$ have second and first-order bounded derivatives respectively; $u_{1}(0)=0,u_{1}'(0)>0,u_{1}(y)>0$ for $y>0$ and $u_{1}\rightarrow U(0)>0 $ as $y\rightarrow \infty$; $u_{1},u_{1}'$ and $u_{1}''$ are bounded and H$\ddot{o}$lder continuous. Moreover, the compatibility condition holds at the origin, i.e.,
$$\frac{1}{\rho(0)}u_{1}''(y)-\frac{\partial_{x}p(0)}{\rho(0)}-v_{0}u_{1}'(y)=\mathcal{O}(y^{2})$$
in the neighborhood of $y=0$. Then, for some $X>0$ there exists a solution $(u,v)$ of the problem  (\ref{w.42})-(\ref{w.43}) in  $D=\{ 0<x<X,y>0\}$, satisfying that $u$ is bounded and continuous in $\overline{D}$, $u>0$ for $y>0$, $u_{y}>m_{1}$ in $0\leq y\leq y_{0}$ for some positive constants $m_{1}$ and $y_{0},u_{y}$ and $u_{yy}$ are bounded and continuous in $D$; $u_{x}, v$ and $v_{y}$ are bounded locally in $D$.

In particular, if $|u_{1}'(y)|\leq m_{2}\exp (-m_{3}y)$  for two positive constants $m_{2}$ and $m_{3}$, then $u_{x}$ and $v_{y}$ are bounded in $D$. Moreover, if either $\rho_{x}\leq 0, v_{0}\leq 0$ or $\rho_{x}\leq -m_{4}$ for a positive constant $m_{4}$, then the solution exists globally.
\end{Theorem}
 \begin{Theorem}(\cite{GGW5})
For the initial boundary value problem (\ref{w.41}), (\ref{w.43}) and (\ref{w.44}), let $\rho , v_{0}, u_{0},u_{1},\frac{u_{0y}}{U(0,x)}$ and $\frac{u_{1y}}{U(t,0)}$ be smooth enough and the data is sixth order compatible at the intersection of the parabolic boundary. In addition, $\partial_{y}u_{i}>0, i=0,1$,  for $y\geq 0$, and satisfy
$$K_{1}(U-u_{i})^{k}\leq \frac{\partial u_{i}}{\partial_{y}}\leq K_{2}(U-u_{i})^{k},\ \ \ \  i=0,1$$
for some $k\geq 1$ and $K_{j}>0~ (j=1,2)$. Then, there exists a unique solution $(u,v)$ to the problem (\ref{w.41}), (\ref{w.43}) and (\ref{w.44}) in $D=\{0<t<T,0<x<X,y>0\}$ , with either arbitrary large $X$ and small $T$ or arbitrary large $T$ and small $X$. In
addition, the solution satisfies that $u_{y}>0$ for $y\geq 0$, the derivatives appearing in (\ref{w.41}) are bounded and continuous in $\overline{D}$; $\frac{u_{yy}}{u_{y}}$ and $\frac{u_{yyy}u_{y}-u_{yy}^{2}}{u_{y}^{3}}$ are bounded in $D$. Furthermore,
$$K_{1}(U-u )^{k}\leq u_{y}\leq K_{2}(U-u )^{k}.$$
\end{Theorem}
Assuming that $\rho=1+\varepsilon\rho_{1}$, then  (\ref{w.41})-(\ref{w.42}) formally can be regarded as the variation of the following steady and unsteady Prandtl equations for small $\varepsilon$, respectively:
\begin{equation}\left\{
\begin{array}{ll}
  u\partial_{x} u +v\partial_{y} u  = \partial_{y}^{2}u , \\
\partial_{x} u+\partial_{y}  v=0,
\end{array}
 \label{w.45}         \right.\end{equation}
and
\begin{equation}\left\{
\begin{array}{ll}
\partial_{t} u+ u\partial_{x} u +v\partial_{y} u  = \partial_{y}^{2}u , \\
\partial_{x} u+\partial_{y}  v=0.
\end{array}
 \label{w.46}         \right.\end{equation}
 \begin{Theorem}(\cite{GGW5})
(1) Assume that $\rho_{1}(x)$ and $\rho_{1}'(x)$ are bounded, $|v_{0}^{1}-v_{0}|\leq \lambda_{1}\varepsilon ^{\frac{1}{2}}, |u_{1}^{1}(y)-u_{1}(y)|\leq \lambda_{2}\varepsilon ^{\frac{1}{2}}\min\{y, 1\}, |u_{1y}^{1}|+|u_{1y} |\leq \lambda_{3}e^{-\lambda_{4}y} $ hold for some $\lambda_{j}(1\leq j\leq 4)$. Then, the solution $u$ of (\ref{w.42})-(\ref{w.43}) converges to the solution $u^{1}$ of (\ref{w.45}) uniformly in $D$ and that
$$|u^{1}-u |\leq \lambda_{5}\varepsilon ^{\frac{1}{2}}$$
for some $\lambda_{5}>0 $.

(2) Assume that $\rho_{1}(x)$ and $\rho_{1}'(x)$ are bounded,
$|v_{0}^{1}-v_{0}|+ |U_{0}^{1}-U_{0}|\leq \lambda_{6}\varepsilon, \frac{u_{0}}{u_{0}^{1}} =\gamma(t,x,y) \frac{U_{0}}{U_{0}^{1}} $ with
$|\gamma-1|+|\gamma_{x}| +|\gamma_{y}| \leq \lambda_{7}\varepsilon, U_{0}^{1}(1-e^{-\lambda_{8}y}) $ and
$|u_{0y}^{1}|\leq \lambda_{10}e^{-\lambda_{11}y} $ hold for some $\lambda_{j}>0~(6\leq j\leq 11 )$, then the solution $u$ of (\ref{w.41}),
(\ref{w.43}) and (\ref{w.44}) converges to the solution $u^{1}$ of the corresponding problem of (\ref{w.46}) uniformly and that
$$|u-u^{1}|\leq \lambda_{12}\varepsilon$$
for a positive constant $\lambda_{12}$.

\end{Theorem}

Ding and Gong (\cite{dg}) in 2017 used the viscous splitting method to establish global existence of weak solutions to the two dimensional compressible Prandtl equations in $\Omega=\{0<t<\infty,\ 0<x<L,\ 0<y<+\infty\}$
 when the density is favourable:
\begin{eqnarray}
\left\{\begin{array}{l}{\partial_{t} u+u \partial_{x} u+v \partial_{y} u+\frac{\partial_{x} p(\rho)}{\rho}=\frac{1}{\rho} \partial_{y}^{2} u}, \\
 {\partial_{x}(\rho u)+\partial_{y}(\rho v)=-\partial_{t}\rho}, \\
 {\left.u\right|_{t=0}=u_{0}(x,y)},\\
  {\left.(u, v)\right|_{y=0}=(0,v_{0}(t,x), \lim\limits _{y \rightarrow+\infty} u=U(t, x)}=1, \\ {\left.u\right|_{x=0}=u_{1}(x, y)}.
  \end{array}\right.\label{d1}
\end{eqnarray}
Assume that \begin{eqnarray}
u_{0}(x,y)\geq 0,\ \ v_{0}(t,x)\leq 0,\ \ u_{1}(t,y)\geq 0,\ \   U(t,x)>0.\label{d2}
\end{eqnarray}
In \cite{dg}, they initially gave a positive answer to this problem with the favourable pressure $P$, which implies that
\begin{eqnarray}
\partial_{x}\rho\leq 0. \label{d3}
\end{eqnarray}
 We firstly  introduce the following Crocco transformation

\begin{eqnarray}
\tau=t,\quad \xi=x,\quad\eta=\frac{u}{U},\quad w=\frac{u_{y}}{U},
\end{eqnarray}
then the problem \eqref{d1} can be transformed into the following scalar degenerated parabolic equation in $\Omega_{1}=\{0<\tau<\infty, 0<\xi<L, 0<\eta<l\}$,
\begin{eqnarray}
\left\{\begin{array}{l}{w_{\tau} +\eta Uw  _{\xi}  +Aw_{\eta}  +Bw=\frac{1}{\rho} w^{2}w_{\eta\eta} }, \\
 {\left.w\right|_{\tau=0}=w_{0} ,\ \ \left.w\right|_{\xi=0}=w_{1} },\\
  {\left.w\right|_{\eta=1}=0,\ \  \frac{1}{\rho}w_{\eta}-v_{0}+\frac{C}{w}|_{\eta=0}=0 }.
 \end{array}\right.\label{d4}
\end{eqnarray}
Here
\begin{eqnarray*}
\left\{\begin{array}{l}
{A=(1-\eta)\frac{U_{t}}{U}+(1-\eta^{2})U_{x}, \ \  w_{i}=\frac{u_{iy}}{U},\ \ i=0,1,}\\
{B=\frac{U_{t}}{U}+\eta U_{x}-\frac{\rho_{t}}{\rho}-\frac{\rho_{t}}{\rho}\eta U,\ \ C=U_{x}+\frac{U_{t}}{U}.}
\end{array}\right.
\end{eqnarray*}
\begin{Theorem}(\cite{dg} )
Suppose the data $w _{i }~(i=0, 1)$ satisfy
\begin{eqnarray}
c_{1}(1-\eta)\leq w_{i}\leq c_{1}(1-\eta). \label{d5}
\end{eqnarray}
   Assumption condition $(\ref{d2})$ and the favourable
condition $(\ref{d3})$ hold for $(\tau, \xi)\in (0,\infty)\times (0,L)$, then problem\eqref{d4} admits a global weak solution  $w$.
\end{Theorem}

 Fan,   Ruan and  Yang (\cite{FRY}) in 2021 proved by the energy method for the local well-posedness of the  two-dimensional compressible boundary
  layer equations (\ref{1.1.7}) by energy method  in a periodic domain
  $ \mathbb{T}\times \mathbb{R}_{+} =\{(x,y)| x\in  \mathbb{R}/\mathbb{Z}, y\in[0, +\infty)  \}$ with the same method as in
    Masmoudi and Wong \cite{mw} in 2015.
The space $H^{s,\gamma}_{\sigma,2\delta}$  for $w:\Omega\rightarrow \mathbb{R}$ is defined by
$$H^{s,\gamma}_{\sigma,2\delta}:=\left\{ w| \|w\|_{H^{s,\gamma}}\leq \infty,(1+y)^{\sigma}w\geq \delta,\sum\limits_{|\alpha|\leq 2}|(1+y)^{\sigma+\alpha_{2}}D^{\alpha}w|^{2}\leq \frac{1}{\delta^{2}} , \right\}$$
where $s\geq 4,\gamma\geq 1,\sigma>\gamma +\frac{1}{2}, \delta\in (0,1), D^{\alpha}=\partial^{\alpha_{1}}_{x}\partial^{\alpha_{2}}_{y} $ with the index $\alpha=(\alpha_{1}, \alpha_{2})$, and the weighted $H^{s}$ norm $\|\cdot\|_{H^{s,\gamma}(\Omega)}$ is defined by
$$ \|w\|_{H^{s,\gamma}}^{2}:=\sum\limits_{|\alpha|\leq s}\|(1+y)^{\gamma+\alpha_{2}}D^{\alpha}w\|^{2} _{L^{2}(\Omega)}. $$
\begin{Theorem} (\cite{FRY})
Given any even integer $s\geq 4$  and real numbers $\gamma,\sigma,\delta$ satisfying $\gamma\geq 1,\sigma>\gamma+\frac{1}{2}$, and $\delta\in(0,1)$,
assume the following conditions on the initial data and the outer flow $(\rho, U)$:\\
(i) Suppose that the initial data $u_{0}-U(0,x)\in H^{s,\gamma-1}$ and $\partial_{y}u_{0}=\partial_{y}u_{0}\in H^{s,\gamma}_{\sigma,2\delta}$ satisfy the
compatibility conditions $u_{0}|_{y}=0, \lim\limits_{y\rightarrow+\infty }u_{0}|_{y=0}=U|_{t=0}$. In addition, when $s=4$,
it is further assumed that $\delta>0$ is chosen small enough such that $\|\partial_{y}u_{0}\|_{H^{s,\gamma}_{g} }\leq C\delta^{-1}$ with a generic constant $C>0$.\\
(ii) The outer flow  $(\rho, U)$ is supposed to satisfy that the density $\rho(t,x)$ has both positive lower and upper bounds, and
$$
\sup\limits_{t}\sum\limits_{l=0}^{ \frac{s}{2}+1} \|\partial_{t}^{l}U\|_{H^{s-2l+2} (\mathbb{T})}
+\sup\limits_{t}\sum\limits_{l=0}^{ \frac{s}{2}+1} \|\partial_{t}^{l}\rho\|_{H^{s-2l+2} (\mathbb{T})}
< +\infty.
$$
Then there exist a time $T:=T(s,\gamma,\sigma,\delta, \|w_{0}\|_{H^{s,\gamma}_{g}  }\rho, ,U)>0$ such that the initial-boundary value problem (\ref{1.1.7})
 has a unique classical solution $(u,v)$ satisfying
$$u-U\in L^{\infty}([0,T],H^{s,\gamma-1})\cap C([0,T]; H^{s}-w)$$
and
$$ \partial_{y} u\in L^{\infty}([0,T]; H^{s,\gamma}_{\sigma,\delta} )\cap C([0,T]; H^{s}-w), $$
where $ H^{s}-w$ is the space $H^{s}$ endowed with its weak topology.

\end{Theorem}

\section{Viscosity Limit of the Navier-Stokes Equations}

  Sammartino, and Caflisch (\cite{SC1}) in 1998 proved the existence of solutions to the Euler and Prandtl equations with analytic initial data.
 The solution will be found as a composite expansion, using the Prandtl solution near the boundary and the Euler solution far from the boundary.
 These results can  be used in (\cite {SC}) as the leading order terms in an asymptotic expansion for the solution of the Navier-Stokes equations with small viscosity.
 \begin{Theorem}(\cite {SC1})
Suppose that $u^{E}(t,x,y)$ and $u^{P}(t,x,Y)$ are solutions of the Euler and Prandtl equations  (\ref{1.1.4}) and (\ref{1.1.2}), respectively, which are analytic in the spatial
variables $x,y,Y$. Then for a short time $T$, independent of $\epsilon$, there is a solution $u^{NS}(t,x,y)$ of the Navier-Stokes equations  (\ref{1.1.1}) with
 \begin{equation*}
 u^{NS}=\left\{
\begin{array}{ll}
u^{E}+O(\epsilon),~\text{outside boundary layer}, \\
u^{P}+O(\epsilon),~\text{inside boundary layer}.
\end{array}
 \label{SC1.1}         \right.\end{equation*}

\end{Theorem}

Authors of \cite{E} in 2000 reviewed recent progresses on the analysis of Prandtl  equation and the related issue of the zero-viscosity limit for the solutions of
 the Navier-Stokes equations, and considered the following equations:

\begin{equation}\left\{\begin{array}{ll}
u_{t}^{\varepsilon}+(u^{\varepsilon}+\nabla)u^{\varepsilon}+\nabla p^{\varepsilon}=\varepsilon \triangle u^{\varepsilon},\\
\nabla\cdot u^{\varepsilon} =0 \ on \ \Omega,\\
u^{\varepsilon}\cdot n|\partial\Omega = 0, \ \ u^{\varepsilon}\cdot \tau|\partial \Omega = 0, \\
u^{\varepsilon}(x,0)=u_{0}(x).
\end{array}\right.\label{h1}\end{equation}
Here $n$ and $\tau$ are the unit normal and tangent vectors at the boundary $\partial\Omega$. They restricted
their attention to two space dimensions.

The steady Prandtl  equations to (\ref{h1}) are
\begin{equation}\left\{\begin{array}{ll}
u_{t}+uu_{x}+vu_{y}+p_{x}=u_{yy},\\
p_{y}=0,
u_{x}+v_{y}=0.
\end{array}\right.\label{h2}\end{equation}

The steady state version of \eqref{h2} has also been of considerable interest in the
engineering literature, which reads

\begin{equation}\left\{\begin{array}{ll}
uu_{x}+vu_{y}+p_{x}=u_{yy},\\
u_{x}+v_{y}=0,\\
u(0,y)=u_{0}(y).
\end{array}\right.\label{h3}\end{equation}
\begin{Theorem}(\cite{E})
Assume that $u_{0}$, $u_{0y}$, $u_{0yy}$ are bounded and H\"{o}lder continuous.
Assume also that $u_{0}$ satisfies the compatibility condition
\begin{eqnarray}
u_{0yy}-P_{x}=O{(y^{2})}
\end{eqnarray}
as $y \rightarrow 0$. Then there exists an X such that \eqref{h3} has a unique strong solution u in the domain $D_{x}=\{(x,y), 0\leq x\leq X, y\geq0 \}$. Moreover, u, $u_{x}$, $u_{y}$, $u_{yy}$ are continuous and bounded in $D_{X}$. If $P_{x}\leq 0$, then $X=+\infty$.
\end{Theorem}
\begin{Theorem}(\cite{E})
Assume that $u_{0}$ satisfies
\begin{eqnarray}
u_{0}^{2}(y)-\frac{3}{2}u_{0y}(y)\int_{0}^{y}u_{0}(z)dz \geq0
\end{eqnarray}
for $y\geq0$. Then\\
(1) There exists an $x^{\ast}$, such that \eqref{h3} has a unique strong solution on $D_{x^{\ast}}=\{(x,y), 0\leq x\leq x^{\ast}, y>0 \}$, but this strong solution cannot be extended to $x>x^{\ast}$.\\
(2) The sequences of functions $\{u_{\lambda}\}$ defined by
\begin{eqnarray}
u_{\lambda}(x,y)=\frac{1}{\lambda^{\frac{1}{2}}}u(x^{\ast}-x\lambda, \lambda^{\frac{1}{4}}y)
\end{eqnarray}
is compact in $C^{0}(D)$, $D=\{(x,y):x,y\geq0\}$.
\end{Theorem}
\begin{Theorem}(\cite{E})
Assume that $u_{0}$ satisfies
\begin{eqnarray}
u_{0y}(x,y)>0 \ for \ y\geq0.
\end{eqnarray}
Then \eqref{h2} has a unique global strong solution.
\end{Theorem}
\begin{Theorem}(\cite{E})
Assume that
\begin{eqnarray}
E(a_{0})=\int_{0}^{\infty}(\frac{1}{2}a_{0y}^{2}+\frac{1}{4}a_{0}^{3})dy < 0,
\end{eqnarray}
where $a_{0}(y)=a(y,0)$. Then global strong solutions of \eqref{h2} do not exist.
\end{Theorem}

\begin{Theorem}(\cite{E})
For any $N > 0$, there exists a profile $u_{\ast} \in C^{\infty}(\mathbb{R}_{+})$, such that
the solution of the Navier-Stokes equation $u^{\varepsilon}_{\ast}$ with initial data
\begin{eqnarray}
u^{\varepsilon}_{\ast}(x,y,0)=(u_{\ast}(\frac{y}{\sqrt{\varepsilon}}),0)
\end{eqnarray}
has the following property: for any $s>0$, and sufficiently small $\varepsilon >0$, there exists a solution
to the  Navier-Stokes equation $u^{\varepsilon}$, such that
\begin{eqnarray}
\|u^{\varepsilon}(\cdot,0)-u^{\varepsilon}_{\ast}(\cdot,0)\|_{H^{s}(\Omega)} \leq \varepsilon^{N}
\end{eqnarray}
and
\begin{eqnarray}
\|u^{\varepsilon}(\cdot,T_{\varepsilon})-u^{\varepsilon}_{\ast}(\cdot,T_{\varepsilon})\|_{H^{s}(\Omega)} \rightarrow +\infty
\end{eqnarray}
as $
\varepsilon \rightarrow0$, where
\begin{eqnarray}
T_{\varepsilon}=O(\sqrt{\varepsilon}\log\frac{1}{\varepsilon}).
\end{eqnarray}
\end{Theorem}

\begin{Theorem}(\cite{E})
Fix $T > 0$. Then $u^{
\varepsilon}(\cdot,t)\rightarrow u^{0}(\cdot,0)$ in $L^{2}(\Omega)$, uniformly for $t\in[0,T]$ if
\begin{eqnarray}
\varepsilon \int_{0}^{T}\int_{\Omega}|\nabla u^{\varepsilon}|^{2}d^{2}xdt \rightarrow 0
\end{eqnarray}
or
\begin{eqnarray}
\varepsilon \int_{0}^{T}\int_{T_{\varepsilon}}|\nabla u^{\varepsilon}|^{2}d^{2}xdt \rightarrow 0
\end{eqnarray}
\end{Theorem}
where $T_\varepsilon$ is a strip of width $O(\varepsilon)$ around $\partial \Omega$.

\begin{Theorem}(\cite{E})
Fix any $T>0$. Then $\max_{0\leq t\leq T}\|u^{\varepsilon}(\cdot,t)-u^{0}(\cdot,t)\|_{L^{2}(\Omega)}\rightarrow0$ as $\varepsilon_{\|}\rightarrow0$ if $\varepsilon_{\perp}/\varepsilon_{\|}\rightarrow0$, here $\delta=\sqrt{ \varepsilon_{\perp}}$.
\end{Theorem}

 Greniner (\cite{G})  in 2000, by giving examples of nonlinearly unstable solutions of Euler equations (\ref{1.1.4}),   proved an instability result for
 Prandtl-type boundary layers (\ref{1.1.1}).
\begin{Theorem}(\cite{G}) \label{ }(Nonlinear Instability of Euler Equations)
Let $u^{s}(y)\in C^{\infty}(\Omega)$ be a shear flow. If $u^{s}(y)$ is linearly unstable (in the sense that there exists an eigenvalue of a linearized Euler
 operator with a nonnegative real part), then it is also nonlinearly unstable in the following sense: There exist positive constants $\delta_{0}>0$
 and $C_{s}>0$~ (for $s\in\mathbb{N}$) such that for every s arbitrarily large and every $\delta>0$ arbitrarily small,  there exists a smooth solution
  $u(t,x,y)$ of problem (\ref{1.1.4}), such that
$$\|u(0,x,y)-u^{s}(y)\|_{H^{s}}\leq \delta$$
and a time
$$T_{\delta}\leq C_{s} \log \delta^{-1}+C_{s}$$
such that
$$\|u(T_{\delta},x,y)-u^{s}(y)\|_{L^{\infty}}\geq \delta_{0}$$
and
$$\|u(T_{\delta},x,y)-u^{s}(y)\|_{L^{2}}\geq \delta_{\delta_{0}}.$$
\end{Theorem}

\begin{Theorem}(\cite{G})\label{ }
 (Nonlinear Instability of Prandtl Boundary Layers) Let $\Omega=\mathbb{R}^{d-1}\times\mathbb{R}^{+}$ and $u^{s}(y)\in C^{\Omega}$ with $u^{s}(y)=0$ be a
 smooth shear layer profile, linearly unstable for Euler equations. Let $n$ be an integer, arbitrarily large. Then there exists a positive constant $\delta>0$
 such that for every s arbitrarily large and for $\nu$ small enough, there exists a sequence of solutions $u^{\nu}(t,x,y)$ of (\ref{1.1.1}) with
$$\|u^{\nu}(0,x,y)-u^{s}(\frac{y}{\sqrt{\nu}})\|_{H^{n}}\leq \nu^{n}$$
and times $T^{\nu}$ with $\lim\limits_{\nu\rightarrow 0}T^{\nu}=0$ and
$$\|\text{curl} ~u^{\nu}(T^{\nu},x,y)- \text{curl} ~u^{s}(T^{\nu},\frac{y}{\sqrt{\nu}}) \|_{L^{\infty}}\rightarrow+\infty, $$
$$\| u^{\nu}(T^{\nu},x,y)-   u^{s}(T^{\nu},\frac{y}{\sqrt{\nu}}) \|_{L^{\infty}} \geq \delta_{0} \nu^{\frac{1}{4}} $$
where $u^{s}(y)$ is the solution of
$$\partial_{t}u^{s}(y)- \partial_{yy}^{2}u^{s}(y)=0$$
with initial data $u^{s}_{0}(y)$.
\end{Theorem}

Teman (\cite{T}) in 2002 studied the local well-posed for the boundary layer of wall bounded flows in a channel at small viscosity when the boundaries are uniformly noncharacteristic.
 Denoting now $v^\epsilon$ the perturbation, i.e., $v^\epsilon=u^\epsilon -(0, 0, -U)$, the author found after
dropping the nonlinear terms:
\begin{equation}\left\{\begin{array}{ll}
\frac{\partial v^{\varepsilon}}{\partial t}-\varepsilon \Delta v^{\varepsilon}-UD_{3}v^{\varepsilon}+\nabla p^{\varepsilon}=f,\\
div v^{\varepsilon}=0,\\
v^{\varepsilon}=0 \ \ \ at\ the\ wall, i.e. ~at \ \ z=0,h.
\end{array}\right.\label{h10}\end{equation}

Setting formally $\epsilon=0$ in \eqref{h10}, the author found
\begin{equation}\left\{\begin{array}{ll}
\frac{\partial v^{0}}{\partial t}-UD_{3}v^{0}+\nabla p^{0}=f,\\
div v^{0}=0,\\
v^{0}=0 \ \ \  at \ \ z=h,
v_{3}^{0}=0 \ \ \  at \ \ z=0,
\end{array}\right.\label{h11}\end{equation}
where
\begin{equation}\begin{array}{ll}
\varphi^{\varepsilon}=curl~ \psi^{\varepsilon}=(-\frac{\partial \psi_{2}^{\varepsilon}}{\partial z},-\frac{\partial \psi_{1}^{\varepsilon}}{\partial z},-\frac{\partial \psi_{1}^{\varepsilon}}{\partial y}+\frac{\partial \psi_{2}^{\varepsilon}}{\partial x}),
\end{array}\label{h12}\end{equation}

\begin{equation}\left\{\begin{array}{ll}
\|\varphi^{\varepsilon,2}\|_{L^{\infty}(0;T;H^{k})} \leq k\varepsilon, \ \ \ \  \forall k,\\
\|\varphi^{\varepsilon,1}\|_{L^{\infty}(0;T;H^{k})} \leq k\varepsilon^{1/2-k},\\
\|\varphi^{\varepsilon,1}\|_{L^{\infty}(0;T;L^{\infty})} \leq k,\\
\|z\varphi^{\varepsilon,1}\|_{L^{\infty}(0;T;L^{\infty})} \leq k\varepsilon,\\
\|z\varphi^{\varepsilon,1}\|_{L^{\infty}(0;T;L^{2})} \leq k\varepsilon^{\frac{3}{2}},\\
\|z\nabla\varphi^{\varepsilon,1}\|_{L^{\infty}(0;T;L^{2})} \leq k\varepsilon^{\frac{3}{2}},
\end{array}\right.\label{h13}\end{equation}
where $k$ is a generic constant independent of the kinematic viscosity $\varepsilon$.

We deduce that the adjusted difference satisfies the equation
\begin{equation}\left\{\begin{array}{ll}
\frac{\partial w^{\varepsilon}}{\partial t}-\varepsilon \Delta w^{\varepsilon}-UD_{3}w^{\varepsilon}+\nabla (p^{\varepsilon}-p^{0})=-\frac{\partial \varphi^{\varepsilon}}{\partial t}+\varepsilon \Delta v^{0}-\varepsilon \Delta \varphi^{\varepsilon}-UD_{3}\varphi^{\varepsilon},\\
div w^{\varepsilon}=0,\\
w^{\varepsilon}=0 \ \ \  at \ \ z=0,h,\\
w^{\varepsilon}=0 \ \ \  at \ \ t=0.
\end{array}\right.\label{h14}\end{equation}

Recall the Navier-Stokes equations
\begin{equation}\left\{\begin{array}{ll}
\frac{\partial u^{\varepsilon}}{\partial t}-\varepsilon \Delta u^{\varepsilon}+(u\cdot\nabla )u^{\varepsilon}+\nabla p^{\varepsilon}=f,\\
div u^{\varepsilon}=0.
\end{array}\right.\label{h15}\end{equation}

The corresponding "inviscid" Euler equation $(\epsilon=0)$ is

\begin{equation}\left\{\begin{array}{ll}
\frac{\partial u^{0}}{\partial t}-(u\cdot\nabla )u^{0} +\nabla p^{0}=f,\\
div u^{0}=0,\\
u^{0}=(0,0,-U) \ \ \  at \ \ z=h,\\
u_{3}^{0}=-U \ \ \  at \ \ z=0,
\end{array}\right.\label{h16}\end{equation}
and
\begin{equation}
\begin{array}{ll}
\frac{\|v^{0}(t)\|_{L^{\infty}(z=0)}}{U}\leq \frac{e}{32} \ \ \  for \ t\leq T_{\ast}.
\end{array} \label{h17}
\end{equation}

\begin{Theorem}( \cite{T})
For $f$ an $v_{0}$ given, let $v^{\varepsilon}$ and $v^{0}$ be the solutions of \eqref{h10}
and \eqref{h11} and assume that $v^{0}$ is sufficiently regular.
Then, as $\varepsilon\rightarrow0$, $v^{\varepsilon}\rightarrow v^{0}$ is estimated by \eqref{h14} and
\begin{equation*}\left\{\begin{array}{ll}
\|w^{\varepsilon}\|_{L^{\infty}([0,T]\times \Omega)} =\|u^{\varepsilon}-u^{0}-\varphi^{\varepsilon}\|_{L^{\infty}([0,T]\times \Omega)}
\leq k\varepsilon^{\frac{3}{4}},\\
\|w^{\varepsilon}\|_{L^{\infty}(0;T;H)} \leq k\varepsilon, \\
\|w^{\varepsilon}\|_{L^{2}(0;T;(H^{1})^{3})} \leq k\varepsilon^{\frac{1}{2}}.
\end{array}\right.\end{equation*}
The difference $w^{\varepsilon}=v^{\varepsilon}-v^{0}-
\varphi^{\varepsilon}$, the corrector $\varphi^{\varepsilon}$ being given by \eqref{h12}, is
estimated by \eqref{h13}.
\end{Theorem}

\begin{Theorem}( \cite{T}) \label{T.1}
Let f and $u_{0}$ be smooth functions satisfying certain
compatibility conditions so that the inviscid problem \eqref{h16} is well-posed
(at least locally in time). Let $u^{\varepsilon}$ and $u^{0}$ be the solution of \eqref{h15} and
\eqref{h16}. Then there exists a time $T_{\ast}>0$ so that \eqref{h17} is satisfied.
Moreover, as $\varepsilon\rightarrow0$, $v^{\varepsilon}\rightarrow v^{0}$ is estimated as
\begin{equation*}\left\{\begin{array}{ll}
\|u^{\varepsilon}-u^{0}\|_{L^{\infty}(0;T_{\ast};H)} =\|v^{\varepsilon}-v^{0}\|_{L^{\infty}(0;T_{\ast};H)}\leq k\varepsilon^{\frac{1}{2}},\\
\|u^{\varepsilon}-u^{0}+\varphi^{\varepsilon}\|_{L^{2}(0;T_{\ast};V)} = \|v^{\varepsilon}-v^{0}+\varphi^{\varepsilon}\|_{L^{2}(0;T_{\ast};V) }
\leq k\varepsilon^{\frac{1}{2}}.
\end{array}\right.\end{equation*}
The adjusted difference $w^{\varepsilon}=u^{\varepsilon}-u^{0}+
\varphi^{\varepsilon}=v^{\varepsilon}-v^{0}+
\varphi^{\varepsilon}$, the corrector $\varphi^{\varepsilon}$ being given by \eqref{h12}, is
estimated by \eqref{h13}.
\end{Theorem}

\begin{Theorem}( \cite{T})
Under the assumptions of Theorem \ref{T.1}, we have also
\begin{equation*}
\begin{array}{ll}
\|u^{\varepsilon}-u^{0}+\varphi^{\varepsilon}\|
_{L^{\infty}((0,T_{\ast})\times\overline{ \Omega} )} = \|v^{\varepsilon}-v^{0}+\varphi^{\varepsilon}\|_{L^{\infty}((0,T_{\ast} )\times \overline{ \Omega})} \leq k\varepsilon^{\frac{1}{4}}.
\end{array}
\end{equation*}
\end{Theorem}

 Filho (\cite{F}) in 2007 made a survey on   results related to boundary layers and the vanishing viscosity limit for incompressible flow. The main topics covered are: a derivation of Prandtl's boundary layer equation; an outline of the rigorous
theory of Prandtl's equation, without proofs; Kato's criterion for the vanishing viscosity limit; the vanishing viscosity limit with Navier friction condition;
rigorous boundary layer theory for the Navier friction condition and boundary layers for flows in a rotating cylinder.
For the problem
\begin{equation}\left\{\begin{array}{ll}
\partial_{t}v^{1}+v\cdot\nabla v^{1}=U_{t}+UU_{x_{1}}+\partial^{2}_{x_{2}}v^{1},\\
div v=0,\\
v(x_{1},0,t)=0 \ and \  \lim\limits_{x_{2}\rightarrow\infty}v^{1}(x,t)=U(x_{1},t),\\
v(x,0)=v_{0}(x),
\end{array}\right.\label{h*}\end{equation}
the author in  \cite{F} obtained the following Theorems \ref{F.1}-\ref{F.2}.

\begin{Theorem}(\cite{F})\label{F.1}
Assume that both $U$ and $v_{0}^{1}$ are positive and that, in
addition, $\partial_{x_{2}}v_{0}^{1}\geq 0$. Then there exists a unique global strong solution of \eqref{h*}.
\end{Theorem}

 Let $u^{\nu}$ be a sequence of solutions of the incompressible Navier-Stokes equations in two space dimensions.
\begin{Theorem}(\cite{F})
Fix $T>0$. There exists a constant $c>0$ such that
$u^{\nu}\rightarrow u$ strongly in $L^{\infty}((0,T);L^{2}(\Omega))$; if and only if $\nu\int_{0}^{T}\|\nabla u^{\nu}(\cdot, t)\|_{L^{2}(\Gamma_{c\nu})}^{2}dt\rightarrow 0$
as $v\rightarrow 0$, where $\Gamma_{cv}\equiv\{1-cv<|x|<1\}$.
\end{Theorem}

Let $\omega^{\nu}=\omega^{\nu} (\cdot,t)$   to be the unique solution of the IBVP:
\begin{equation*}\left\{\begin{array}{ll}
\omega^{\nu}_{t}+u^{\nu}\cdot \nabla \omega^{\nu}=\nu\triangle \omega^{\nu} ,\\
u^{\nu}=K_{\Omega}[\omega^{\nu} ],\\
\omega^{\nu} (x,t)=(2k-\alpha)(u(x,t)\cdot \tau)~\text{for}~x\in \partial\Omega,\\
\omega(x,0)=\omega_{0}.
\end{array}\right.\label{ }\end{equation*}
\begin{Theorem}(\cite{F})
For each $p>2$, there
exists a constant $C>0$, independent of $\nu$, such that
\begin{eqnarray*}
\|\omega^{\nu}(\cdot,t)\|_{L^{p}(\Omega)}\leq C(\|\omega_{0}\|_{L^{p}(\Omega)}+\|u_{0}\|_{L^{2}(\Omega)}),
\end{eqnarray*}
with $u_{0}=K_{\Omega}[\omega_{0}]$.
\end{Theorem}

\begin{Theorem}(\cite{F})
Let $u^{v}$ be the solution of the 2D Navier-Stokes equations in the unit disk, with no slip boundary data with respect to boundary rotation with prescribed
angular velocity $\alpha \in BV(\mathbb{R})$. Assume that the initial velocity $u_{0}\in L^{2}(D)$ has
circular symmetry, $D\equiv\{|x|\leq 1\}\subset\mathbb{R}^{2}$, i.e., $u_{0}=v_{0}(|x|)x^{\perp}$. Then, $u_{0}$ is a steady solution of the 2D
Euler equation and $u^{v}$
converges strongly in $L^{\infty}([0,T],L^{2}(D))$ to $u_{0}$ as $v\rightarrow 0$. Here
$$ x^{\perp}=(-x_{2},x_{1}),\quad \nabla^{\bot}=(-\partial_{x_{2}},\partial_{x_{1}} ) .$$
\end{Theorem}

\begin{Theorem}(\cite{F})\label{F.2}
Let $\alpha\in BV(\mathbb{R})$ and assume that $\alpha$ is compactly supported in $(0,\infty)$.
Let $u_{0}\in L^{2}(D,\mathbb{R}^{2})$ be covariant under rotations and assume that $\omega_{0}=\nabla^{\bot}\cdot u_{0} \in L^{1}(D)$, $D\equiv\{|x|\leq 1\}\subset\mathbb{R}^{2}$. Then we have:
(1) There exists a constant $C>0$ such that
\begin{eqnarray*}
\int_{D}|\omega^{v}(x,t)|dx \leq C(\|\omega_{0}\|_{L^{1}}+\|\alpha\|_{BV}),
\end{eqnarray*}
for almost all time.\\
(2) For almost every $t\in \mathbb{R}$,
\begin{eqnarray*}
\alpha(t)=\int_{D}\omega^{v}(x,t)dx.
\end{eqnarray*}
(3) We also have, for any $0<a<1$,
\begin{eqnarray*}
\|\omega^{v}-\omega_{0}\|_{L^{1}({|x|\leq a})}\rightarrow 0 \ as \ v\rightarrow 0,
\end{eqnarray*}
a.e. in time. If $\omega_{0}\in C^{0}(D)$,  then the convergence is uniform in time.

\end{Theorem}

Wang, Wang and Xin (\cite{WWX}) in 2010 studied the vanishing viscosity limit for the incompressible Navier-Stokes equations with the Navier friction boundary condition. To simplify the expansion of solutions in terms of the viscosity, they  only considered the case that the slip length $\alpha$ in the Navier boundary condition is a power of the viscosity $\varepsilon$. They considered the following problem:

\begin{eqnarray}\left\{\begin{array}{ll}
u^\varepsilon_t+u^\varepsilon\nabla u^\varepsilon+\nabla p^\varepsilon=\varepsilon^{\frac{1}{2}}\partial_x^2u^\varepsilon+\varepsilon\partial_y^2u^\varepsilon,\\
\partial_xu^\varepsilon+\partial_yv^\varepsilon=0,\\
u^\varepsilon_2=0,\ \beta u^\varepsilon_1-\varepsilon^{\frac{1}{4}}\frac{\partial u_1^\varepsilon}{\partial y}=0,\ on\ y=0,\\
u^\varepsilon|_{t=0}=u_0(x,y),
\end{array}\right.\label{109.1}\end{eqnarray}
where $\beta$ is a positive constant.
\begin{Theorem} (\cite{WWX})
For the problem \eqref{109.1},  suppose that the initial data $u_0$ of
\eqref{109.1} belongs to $H^s(\Omega)$  for a fixed $s>18$, with   div $ u_0=0$ and $u_{0,2}=0$ on $y=0$. Then, in $L^\infty((0,T]\times\mathbb{R}_+^3)$, we have
\begin{eqnarray*}\left\{\begin{array}{ll}
u_1^\varepsilon(t,x,y)=u_1^{I,0}(t,x,y)+\sum_{j=1}^3\varepsilon^{j/4}(u_1^{I,j}(t,x,y)
+u_1^{B,j}(t,x,\frac{y}{\sqrt\varepsilon}))+O(\varepsilon),\\
u_2^\varepsilon(t,x,y)=\sum_{j=1}^3\varepsilon^{j/4}(u_2^{I,j}(t,x,y)
+u_2^{B,3}(t,x,\frac{y}{\sqrt\varepsilon}))+O(\varepsilon),\\
p^\varepsilon(t,x,y)=\sum_{j=1}^3\varepsilon^{j/4}(u_1^{I,j}(t,x,y)+O(\varepsilon),\\
\end{array}\right.\label{109.2}\end{eqnarray*}
for rapidly decreasing $ (u_1^{B,j},u_2^{B,3})  $ in $z\rightarrow+\infty$, where $ u_1^{I,0}$ is a solution
 for the Euler equations, $( u_1^{I,j} , u_2^{I,j})$ are solutions for the linearized Euler equations.
\end{Theorem}

The authors of \cite{WWX} further considered the following problem in $\{t,y>0,\ x\in \mathbb{R}\}$,
\begin{eqnarray}\left\{\begin{array}{ll}
u^\varepsilon_t+u^\varepsilon\nabla u^\varepsilon+\nabla p^\varepsilon=\varepsilon^{1-\frac{1}{2}\gamma}\partial_x^2u^\varepsilon+\varepsilon\partial_y^2u^\varepsilon,\\
\partial_xu^\varepsilon+\partial_yv^\varepsilon=0,\\
u^\varepsilon_2=0,\ \beta u^\varepsilon_1-\varepsilon^{\gamma}\frac{\partial u_1^\varepsilon}{\partial y}=0,\ on\ y=0,\\
u^\varepsilon|_{t=0}=u_0(x,y).
\end{array}\right.\label{109.3}\end{eqnarray}
\begin{Theorem}(\cite{WWX})
For the problem \eqref{109.3}, under certain assumptions on the regularity and compatibility conditions on the initial data $u_0$,
the solutions to \eqref{109.3} have the following expansions:
\begin{eqnarray}\left\{\begin{array}{ll}
u_1^\varepsilon(t,x,y)=\sum_{k=0}^{b-a-1}\varepsilon^{rk} u_1^{I,k}(t,x,y)+\sum_{j=b-a}^J\varepsilon^{rj}(u_1^{I,j}(t,x,y)
+u_1^{B,j}(t,x,\frac{y}{\sqrt\varepsilon}))+O(\varepsilon^{rJ}),\\
u_2^\varepsilon(t,x,y)=\sum_{k=0}^{2b-a-1}\varepsilon^{rk} u_2^{I,k}(t,x,y)+\sum_{j\geq2b-a}^J\varepsilon^{rj}(u_2^{I,j}(t,x,y)
+u_2^{B,j}(t,x,\frac{y}{\sqrt\varepsilon}))+O(\varepsilon^{rJ}),\\
p^\varepsilon(t,x,y)=\sum_{k=0}^{3b-a-1}\varepsilon^{rk} p^{I,k}(t,x,y)+\sum_{j\geq3b-a}^J\varepsilon^{rj}(p^{I,j}(t,x,y)
+p^{B,j}(t,x,\frac{y}{\sqrt\varepsilon}))+O(\varepsilon^{rJ}),
\end{array}\right.\label{109.4}\end{eqnarray}
in $L^\infty((0,T]\times \mathbb{R}^2)$ for a fixed $J\geq1$.

\end{Theorem}

Liu and  Wang (\cite{LW1}) in 2014 using curvilinear coordinate, proved that boundary layers of the Navier-Stokes equations in a bounded domain, still
satisfy the classical Prandtl equations when the viscosity goes to zero.
 \begin{Theorem}(\cite{LW1})\label{ }
 The solutions $({\bf u}^{\nu} ,p^{\nu})$ to the problem (\ref{1.1.1}) formally admit the following asymptotic expansions:
\begin{equation*}\left\{
\begin{array}{ll}
{\bf u}^{\nu} (t,x)={\bf u}^{I,0}+
u^{B,0}(t,s(x),\frac{d(x)}{\sqrt{\nu}})+\sqrt{\nu}v^{p}(t,s(x),\frac{d(x)}{\sqrt{\nu}})\overrightarrow{n}+\mathbb{O}(\sqrt{\nu}),   \\
p^{\nu}=p^{I,0}(t,x)+\mathbb{O}(\sqrt{\nu}),
\end{array}
      \right.\end{equation*}
where $v^{p}(t,s,\eta)$ is a scalar function, and $u^{B,0}(t,s,\eta), v^{p}(t,s,\eta)$ are rapidly decreasing as $\eta\rightarrow+\infty$.
 Also, $({\bf u}^{I,0},p^{I,0})(t,x)$ satisfies the boundary value problem for the incompressible Euler system  (\ref{1.1.4}); the
  term $u^{B,0}(t,s,\eta)$ satisfies
$$ u^{B,0}(t,s,\eta)\cdot \overrightarrow{n}=0 ,$$
and $(( u^{B,0}+\overline{   u ^{I,0}})\cdot \overrightarrow{\tau},v^{p})(t,s,\eta)$ satisfies the boundary value problem for the Prandtl system (\ref{1.1.6}),
where $\overline{u^{I,0}}$ satisfies the Bernoulli's law.

\end{Theorem}

Guo and Nguyen (\cite{gn}) in 2017  studied the validity of the Prandtl boundary layer theory in the inviscid limit for steady incompressible Navier-Stokes flows.
To this end, they considered the stationary incompressible Navier-Stokes equations,
\begin{eqnarray}\left\{\begin{array}{ll}
 U U_X+V U_Y+P_X =\varepsilon U_{XX}+\varepsilon U_{YY},\\
 U V_X+V V_Y+P_Y =\varepsilon V_{XX}+\varepsilon V_{YY},\\
 U_{X}+V_{Y} =0 ,\end{array}\right.\label{45.1}\end{eqnarray}
posed in a two dimensional domain $\Omega=\{(X,Y):0\leq X\leq L, 0\leq Y<+\infty\}$, with a no-slip boundary condition
\begin{equation*}\label{45.2}
\left[U(X,0), V(X,0)\right] \equiv\left[u_{e}, 0\right].
\end{equation*}
Throughout \cite{gn}, they assumed that the outside Euler flow is a shear flow
\begin{equation*}\label{45.3}
\left[U_{0}, V_{0}\right] \equiv\left[u_{e}^{0}(Y), 0\right].
\end{equation*}
They  worked with the scaled boundary layer, or Prandtl's variables:
$$x=X, \quad y=\frac{Y}{\sqrt{\varepsilon}}.$$
In these new variables, the Navier-Stokes equations \eqref{45.1} now read as
\begin{eqnarray}\left\{\begin{array}{ll}
 U^{\varepsilon} U_{x}^{\varepsilon}+V^{\varepsilon} U_{y}^{\varepsilon}+P_{x}^{\varepsilon} =U_{y y}^{\varepsilon}+\varepsilon U_{x x}^{\varepsilon}, \\
 U^{\varepsilon} V_{x}^{\varepsilon}+V^{\varepsilon} V_{y}^{\varepsilon}+\frac{P_{y}^{\varepsilon}}{\varepsilon}
 =V_{y y}^{\varepsilon}+\varepsilon V_{x x}^{\varepsilon}, \\ U_{x}^{\varepsilon}+V_{y}^{\varepsilon} =0. \end{array}\right.\label{45.4}
\end{eqnarray}

Precisely, we search for asymptotic expansions of the scaled Navier-Stokes solutions $( U^{\varepsilon}, V^{\varepsilon}, P^{\varepsilon})$
in the following form:

\begin{eqnarray}\left\{\begin{array}{ll}
 U^{\varepsilon}(x, y) =u_{e}^{0}(\sqrt{\varepsilon} y)+u_{p}^{0}(x, y)+\sqrt{\varepsilon} u_{e}^{1}(x, \sqrt{\varepsilon} y)+\sqrt{\varepsilon} u_{p}^{1}(x, y)
 +\varepsilon^{\gamma+\frac{1}{2}} u^{\varepsilon}(x, y), \\
V^{\varepsilon}(x, y) =v_{p}^{0}(x, y)+v_{e}^{1}(x, \sqrt{\varepsilon} y)+\sqrt{\varepsilon} v_{p}^{1}(x, y)+\varepsilon^{\gamma+\frac{1}{2}} v^{\varepsilon}(x, y), \\
P^{\varepsilon}(x, y) =\sqrt{\varepsilon} p_{e}^{1}(x, \sqrt{\varepsilon} y)+\sqrt{\varepsilon} p_{p}^{1}(x, y)+\varepsilon p_{p}^{2}(x, y)
+\varepsilon^{\gamma+\frac{1}{2}} p^{\varepsilon}(x, y), \end{array}\right.\label{45.5}
\end{eqnarray}
for some $\gamma>0$, $\left[u_{e}^{j}, v_{e}^{j}, p_{e}^{j}\right]$ and $\left[u_{p}^{j}, v_{p}^{j}, p_{p}^{j}\right],$ $j=0,1$,  denote the
Euler and Prandtl profiles, respectively, and $\left[u^{\varepsilon}, v^{\varepsilon}, p^{\varepsilon}\right]$ collects all the remainder solutions.
\begin{Theorem}(\cite{gn})
Let $u_b>0$ be a constant tangential velocity of the Navier-Stokes flow
on the boundary $Y=0$ and let be a given smooth and positive Euler flow so
that $\partial_{Y} u_{e}^{0}$  and its derivatives decay exponentially fast to zero at infinity. We take Euler
data $u_{b}^{1}(Y), V_{b 0}(Y), V_{b L}(Y)$, and Prandtl data $\bar{u}_0(y)$ and $\bar{u}_1(y)$ to be sufficiently
smooth and decay exponentially fast at infinity in their arguments. We assume that $\left|V_{b L}(Y)-V_{b 0}(Y)\right| \lesssim L$ for all $L$, and
$$\min _{0 \leq y \leq \infty}\left\{u_{e}^{0}(\sqrt{\varepsilon} y)+\bar{u}_{0}(y)\right\}>0.$$
Then, there exists a positive number $L_{0}$ that depends only on
 the given data so that the boundary layer expansions \eqref{45.5}, with the profiles and the remainder satisfying the boundary conditions
 for $\gamma \in\left(0, \frac{1}{4}\right)$ and $0<L \leq L_{0}$.  Precisely, $\left[U^{\varepsilon}, V^{\varepsilon}, P^{\varepsilon}\right]$
 as defined in \eqref{45.5}  is the unique solution to the Navier-Stokes equations \eqref{45.4}, so that the remainder
  solutions $\left[u^{\varepsilon}, v^{\varepsilon}\right]$ satisfy
$$
\left\|\nabla_{\varepsilon} u^{\varepsilon}\right\|_{L^{2}}+\varepsilon^{\frac{1}{2}}\left\|\nabla_{\varepsilon} v^{\varepsilon}\right\|_{L^{2}}
+\varepsilon^{\frac{\gamma}{2}}\left\|u^{\varepsilon}\right\|_{L^{\infty}}+\varepsilon^{\frac{1}{2}+\frac{\gamma}{2}}\left\|v^{\varepsilon}\right\|_{L^{\infty}}
\leq C_{0},
$$
for some constant $C_0$ that depends only on the given data. Here, $\nabla_{\varepsilon}=\left(\sqrt{\varepsilon} \partial_{x}, \partial_{y}\right)$,
 $\|\cdot\|_{L^{p}}$ denotes the usual $L^p$ norm over $[0, L] \times \mathbb{R}_{+}$ .
\end{Theorem}

Wang and Liu (\cite{WL}) in 2017 used $\lambda$-weighted energy method and the matched asymptotic expansion method to prove the well-posedness of solutions to the
boundary layer problem.  They studied the boundary layer problem and the quasineutral limit of the isentropic Euler-Poisson system in plasma for
 $t>0$, $x=(x_1,x_2,x_3)=(y,x_3)\in \mathbb{R}^3_+=\mathbb{R}^2\times \mathbb{R}_+$,
\begin{eqnarray}\left\{\begin{array}{ll}
\partial_tn^{\lambda}+div(n^{\lambda}u^{\lambda})=0,\\
\partial_tu^{\lambda}+u^{\lambda}\cdot\nabla u^{\lambda}+\nabla h(u^{\lambda})=-\nabla\phi^{\lambda},\\
-\lambda^2\Delta\phi^{\lambda}=n^{\lambda}-1,
\end{array}\right.\label{104.1}\end{eqnarray}
where $\lambda$ is the (scaled) Debye length, $(n^{\lambda},u^{\lambda})=(u_1^{\lambda},u_2^{\lambda},u_3^{\lambda})=(u_y^{\lambda},u_3^{\lambda})$,
$ \phi^{\lambda}$ denote
the electron density, the electron velocity and the electric potential, respectively. The function $h=h(n)$ is the enthalpy for the system and satisfies
$$h'(n)=\frac{p'(n)}{n},\ n>0$$
where where $p=p(n)$ is the pressure density function, and we assume
$$p(n)=a^2n^{\gamma},\ n>0,a\neq0,\ \gamma\geq1.$$
From \eqref{104.1} in the limit $\lambda\rightarrow 0$ and the fact that $n=1$, we have the following incompressible Euler system:
\begin{eqnarray}\left\{\begin{array}{ll}
div u=0,\\
\partial_tu+u\cdot\nabla u=-\phi.
\end{array}\right.\label{104.2}\end{eqnarray}
\begin{Theorem}( \cite{WL} )
Let $(u(0) ,\phi(0) )$ be a solution to \eqref{104.2} such that
$u(0)\in C^0 ([0,T],H^s (R^3_+ ))$ with s large enough. There is $\lambda_0
>0$ such that for any
$\lambda\in(0,\lambda_0]$, there exists $(n^\lambda ,u^\lambda ,\phi^\lambda )$, a unique solution to \eqref{104.1} also defined on $[0,T]$ such that
\begin{eqnarray*}
 \lim\limits_{\lambda\rightarrow0}\sup_{[0,T]}(\|n^\lambda-1\|_{L^2}+\|u^\lambda-u(0)\|_{L^2}+\|\lambda\nabla
\phi^\lambda-\lambda\nabla\phi(0)\|_{L^2})\rightarrow0.
\end{eqnarray*}
\end{Theorem}

Guo and Lyer (\cite{GL}) in 2018 investigated the asymptotic behavior of solutions to  incompressible Navier-Stokes equations on the two dimensional  as the viscosity vanishes, that is as $\varepsilon\rightarrow 0$.

They considered the steady, incompressible Navier-Stokes equations on the two dimensional domain $(X,L)\in\Omega\times(0,\infty)$, Denoting the velocity
$\textbf{U}^{N,S}=(U^{N,S},V^{N,S})$, the equations read:
\begin{eqnarray}\left\{\begin{array}{ll}
\textbf{U}^{N,S}\cdot\nabla \textbf{U}^{N,S}+\nabla P^{N,S}=\varepsilon \textbf{U}^{N,S},\\
\nabla\cdot\textbf{U}^{N,S}=0.
\end{array}\right.\label{43.1}\end{eqnarray}
The system above is taken with the no-slip boundary condition on
\begin{eqnarray*}
 (U^{N,S},V^{N,S})|_{Y=0}=0.\label{43.2}
 \end{eqnarray*}
In \cite{GL}, they addressed Prandtl's classical setup, and worked with the scaled
boundary layer variable:
$$x=X, \quad y=\frac{Y}{\sqrt{\varepsilon}}.$$
In these new variables, the Navier-Stokes equations \eqref{43.1} now read
\begin{eqnarray*}\left\{\begin{array}{ll}
 U^{\varepsilon} U_{x}^{\varepsilon}+V^{\varepsilon} U_{y}^{\varepsilon}+P_{x}^{\varepsilon} =\Delta_{\varepsilon} U^{\varepsilon}, \\ U^{\varepsilon} V_{x}^{\varepsilon}+V^{\varepsilon} V_{y}^{\varepsilon}+\frac{P_{y}^{\varepsilon}}{\varepsilon} =\Delta_{\varepsilon} V^{\varepsilon}, \\ U_{x}^{\varepsilon}+V_{y}^{\varepsilon} =0. \end{array}\right.\label{43.3}
\end{eqnarray*}

They expanded the solution in $\varepsilon$ as:

\begin{eqnarray}\left\{\begin{array}{ll}
 U^{\varepsilon} =u_{e}^{0}+u_{p}^{0}+\sum_{i=1}^n\sqrt{\varepsilon}^{i} u_{e}^{i}+\sqrt{\varepsilon} u_{p}^{i}+\varepsilon^{N_0} u^{\varepsilon}:=u_s+\varepsilon^{N_0} u^{\varepsilon}, \\
V^{\varepsilon} =v_{p}^{0}+v_{e}^{1}+\sum_{i=1}^{n-1}\sqrt{\varepsilon}^{i} (v_{p}^{i}+v_{e}^{i+1})+\varepsilon^{n/2} v_p^{n}+\varepsilon^{N_0} v^{\varepsilon}:=v_s+\varepsilon^{N_0} v^{\varepsilon}, \\
P^{\varepsilon} =P_{e}^{0}+P_{p}^{0}+\sum_{i=1}^n\sqrt{\varepsilon}^{i} P_{e}^{i}+\sqrt{\varepsilon} P_{p}^{i}+\varepsilon^{N_0} P^{\varepsilon}:=P_s+\varepsilon^{N_0} P^{\varepsilon}.  \end{array}\right.\label{43.4}
\end{eqnarray}

 At the leading order, $[u_p^0, v_p^0]$ solves the nonlinear Prandtl equations:
\begin{eqnarray}\left\{\begin{array}{ll}
\bar{u}_p^0u_{px}^0+\bar{v}_p^0u_{px}^0-u_{pyy}^0+P_{px}^0=0,\\
u_{px}^0+u_{py}^0=0,\ P_{py}^0=0,\ u_{p}^0|_{x=0}=U_{p}^0,\  u_{p}^0|_{y=0}=-u^0_e|_{Y=0}.
\end{array}\right.\label{43.5}\end{eqnarray}

It is well known that the Prandtl equations \eqref{43.5} admit the two parameter
scaling invariance:
\begin{eqnarray}
 [\bar{u}^{\lambda,\delta},\bar{v}^{\lambda,\delta}]=[\frac{\lambda^2}{\delta}\bar{u}_{p}^0(\delta x,\lambda y),\lambda\bar{v}_{p}^0(\delta x,\lambda y)],
 \label{43.6}
 \end{eqnarray}
which means that $[u_p^0, v_p^0]$ solves \eqref{43.5},  then so do $[\bar{u}^{\lambda,\delta},\bar{v}^{\lambda,\delta}]$ (with appropriately
modified initial data).

The boundary condition is taken as the following
\begin{eqnarray}\left\{\begin{array}{ll}
v_x^{(\varepsilon)}|_{x=L}=a_1^\varepsilon(y),\ v_{xx}^{(\varepsilon)}|_{x=L}=a_2^\varepsilon(y),\ v_{xxx}^{(\varepsilon)}|_{x=L}=a_3^\varepsilon(y),\\
v^{(\varepsilon)}|_{x=0}=v_0(y),\ v_y^0|+u^0=h(y)\in C^{\infty},\ h(0)=0,\\
v^{(\varepsilon)}|_{y=0}=v_y^{(\varepsilon)}|_{y=0}=u^{(\varepsilon)}|_{y=0}=0,\
v_y^{(\varepsilon)}|_{y=\infty}=0.
\end{array}\right.\label{43.7}\end{eqnarray}
Here, the $a_i^\varepsilon(y)$ are prescribed boundary data which are assumed to satisfy
\begin{eqnarray}
 \|\partial_y^ja_i^\varepsilon\{\frac{1}{\varepsilon^{1/2}}(y)(Y)^m\}\|\leq o(1)\ for \ j=0,1,2,3,4, and\ m\ lager.\label{43.8}
 \end{eqnarray}
\begin{Theorem} (\cite{GL})
 Assume boundary data are prescribed as in \eqref{43.7}-\eqref{43.8}. Assume $0<\delta<<1$ in \eqref{43.6}.
 Then let $0<\varepsilon<<L<<1$. Take $N_0=1+$ and $n=4$ in \eqref{43.4}. Then all terms in the expansion \eqref{43.4} exist and are regular,
 $\|u_s,v_s\|\leq 1.$  The remainders, $[u^{(\varepsilon)}, v^{(\varepsilon)}]$ exists uniquely in the space $X$ and satisfy

$$\|\textbf{u}^{(\varepsilon)}\|_X\leq 1.$$
The Navier-Stokes solutions satisfy
$$\|U^{N,S}-u_e^0-u_p^0\|_\infty\leq \sqrt\varepsilon\ and \ \|V^{N,S}-\sqrt\varepsilon v_e^1- \sqrt\varepsilon v_p^0\|_\infty\leq\varepsilon,$$
where \bqa X=\{(v,u^0,v^0)\in L^2(\Omega)\times L^2(\mathbb{R}_+)\times L^2(\mathbb{R}_+): \|(v,u^0,v^0)\|_X<+\infty;\non    \\
v|_{y=0}=v_y|_{y=0}=v_x|_{x=L}=v_{xx}|_{x=0}=v_{xxx}|_{x=L}=v|_{y=\infty}=0,\non   \\
v^0(0)=v^0_y(0)=\partial^k
_yv(\infty)=0\ \text{for}\ k \geq1,
u^0+v^0_y=h(y), u^0(0)=0\}.  \eqa
\end{Theorem}

We first record the properties of the leading order $(i = 0)$ layers. For the outer
Euler flow, we will take a shear flow, $[u^0_e(Y ), 0, 0]$. The derivatives of $u^0_e$
decay rapidly in $Y$ and that is bounded below, $|u_e^0|\geq 1,$

For the leading order Prandtl boundary layer, the equations are the following:
\begin{eqnarray}\left\{\begin{array}{ll}
\bar{u}_p^0u_{px}^0+\bar{v}_p^0u_{px}^0-u_{pyy}^0+P_{px}^0=0,\\
u_{px}^0+u_{py}^0=0,\ P_{py}^0=0,\ u_{p}^0|_{x=0}=U_{p}^0,\  u_{p}^0|_{y=0}=-u^0_e|_{Y=0}.
\end{array}\right.\label{43.9}\end{eqnarray}
\begin{Theorem} (\cite{GL})
Assume boundary data is prescribed satisfying $U_P^0\in C^\infty$
and exponentially decaying $|\partial_y^j(U_P^0-u_e^0(0))|$ for $j\geq0$ satisfying:
\begin{eqnarray}
\bar{U}_P^0>0\ for \ y>0,\ \partial_y\bar{U}_P^0(0)>0,\ \partial^2_y\bar{U}_P^0(0)\sim y^2\ near \ y=0.\label{43.10}
\end{eqnarray}
Then for some $L > 0$,  there exists a solution, $[\bar{u}^0_p, \bar{v}^0_p]$ to problem \eqref{43.9} satisfying, for some $y_0,m_0>0$,
$$\sup_{x\in[0,L]}\sup_{y\in(0,y_0)}|\bar{u}^0_p,\bar{v}^0_p,\partial_y\bar{u}^0_p,\partial_{yy}\bar{u}^0_p,\partial_x\bar{u}^0_p|\leq1,$$
$$\sup_{x\in[0,L]}\sup_{y\in(0,y_0)}\partial_y\bar{u}^0_p>m_0>0.$$
By evaluating the system \eqref{43.9} and $\partial_y$ of \eqref{43.9} at ${y = 0}$,  we conclude:

$$\partial^2_y\bar{u}^0_p|_{y=0}=\partial^3_y\bar{u}^0_p|_{y=0}=0.$$
\end{Theorem}
\begin{Theorem}(\cite{GL})
Assume the shear flow $u_e^0(Y)\in C^\infty$,  whose derivatives decay
rapidly. Assume \eqref{43.10}  regarding $\bar{u}^0_p|_{x=0}$,  and the conditions
\begin{eqnarray}\left\{\begin{array}{ll}
\partial^3_y\bar{v}^i_p|_{x=0}(0)=\partial_xg_1|_{x=0,y=0}=0,\\
\bar{v}^i_p|''''_{x=0}(0)=\partial_{xy}g_1|_{y=0}(x=0),\\
\bar{u}^0_{py}|_{x=0}(0)u^i_e|_{x=0}(0)-\int_0^\infty\bar{u}^0_pe^{-\int_1^y
\bar{v}^0_p}\{f^i(y)-r^i(y)\}dy=0,
 \end{array}\right.\label{43.11}\end{eqnarray}
where $r^i(y)=\bar{v}^i_p\bar{u}^0_{py}-\bar{u}^0_{p}\bar{v}^i_{py}$. We assume also standard higher order versions
of the parabolic compatibility conditions \eqref{43.11}. Let
$v^i_e|_{x=0}$,  $v^i_e|_{x=L}$, $u^i_e|_{x=0}$, be prescribed smooth and rapidly decaying Euler data. We assume on the data
standard elliptic compatibility conditions at the corners $(0,0)$ and $(L,0)$ obtained
by evaluating the equation at the corners. In addition, assume
\begin{eqnarray*}\left\{\begin{array}{ll}
v_e^1|_{x=0}\sim Y^{-m_1}\ or \ for~ some \ 0<m_1<\infty,\\
\|\partial_Y^k(v_e^i|_{x=0}-v_e^i|_{x=L})(Y)^M\|_\infty\leq L.
\end{array}\right.\label{43.12}\end{eqnarray*}
Then all profiles in $[u_s, v_s]$ exist and are smooth on $\Omega$. The following estimates hold:
\begin{eqnarray*}\left\{\begin{array}{ll}
\bar{u}_p^0>0,\ \bar{u}_{py}^0|_{y=0}>0,\ \bar{u}_{pyy}^0|_{y=0}=\bar{u}_{pyyy}^0|_{y=0}=0,\\
\|\nabla^K({u}_p^0,{v}_p^0)e^{My}\|_\infty\leq1\ for~ any \ K\geq 0,\\
\|{u}_p^1\|_{\infty}+\|\nabla^K({u}_p^ie^{My}\|_\infty+
\|\nabla^j({v}_p^ie^{My}\|_\infty\leq 1\ for ~any \ K\geq 1,\ M\geq0,\\
\|\nabla^K({u}_e^1,{v}_e^1)w_{m_1}\|_\infty\leq1\ for\ some\ fixed\ m_1>1,\\
\|\nabla^K({u}_e^i,{v}_e^i)w_{m_i}\|_\infty\leq1\ for\ some\ fixed\ m_1>1,
\end{array}\right.\label{43.13}\end{eqnarray*}
where $w_{m_i}\sim e^{m_iY}$ or $(1+Y)^{m_i}$.

\end{Theorem}

G$\acute{e}$rard-Varet and  Maekawa (\cite{GM}) in 2019 show the $H^{1}$ stability of shear flows of Prandtl type: $U^{v}=\left(U_{s}(y / \sqrt{v}), 0\right)$, in the steady two-dimensional Navier-Stokes equations, under the natural assumptions that $U_{s}(Y)>0$ for $Y>0, U_{s}(0)=0,$ and $U_{s}^{\prime}(0)>0 .$ Their result is in sharp contrast with the unsteady ones, in which at most Gevrey stability can be obtained, even under global monotonicity and concavity hypotheses.
The authors of (\cite{GM}) also considered the vanishing viscosity limit of the two-dimensional steady Navier-Stokes equations
$$
\left\{\begin{array}{ll}
v^{\nu} \cdot \nabla v^{\nu}-v \Delta v^{\nu}+\nabla q^{\nu}=g^{\nu}, & (x, y) \in \mathbb{T}_{\kappa} \times \mathbb{R}_{+} ,\\
\operatorname{div} v^{\nu}=0, & (x, y) \in \mathbb{T}_{\kappa} \times \mathbb{R}_{+}, \\
\left.v^{\nu}\right|_{y=0}=0.
\end{array}\right.
$$
Here $\mathbb{T}_{\kappa}=\mathbb{R} /(2 \pi \kappa) \mathbb{Z}, \kappa>0,$ is a torus with periodicity $2 \pi \kappa, \mathbb{R}_{+}=\{y \in \mathbb{R} \mid y>$
0\} , while $v^{\nu}=\left(v_{1}^{\nu}, v_{2}^{\nu}\right)$ and $q^{\nu}$ are, respectively, the unknown velocity field and pressure field of the fluid. The positive constant $v$ is the viscosity coefficient. The vector field $g^{\nu}$ is an external force, decaying fast enough at infinity. The usual no-slip condition is prescribed at $y=0$.

Let $U_{s}=U_{s}(Y) \in C^{2}\left(\overline{\mathbb{R}_{+}}\right)$ such that
$$
\begin{array}{c}
U_{s}(0)=0, \quad U_{s}>0 \quad \text { in } Y>0, \quad \lim _{Y \rightarrow \infty} U_{s}(Y)=U_{E}>0 ,\\
\partial_{Y} U_{s}(0)>0 ,\\
\sum_{k=1,2} \sup _{Y \geqq 0}(1+Y)^{3}\left|\partial_{Y}^{k} U_{s}(Y)\right|<\infty.
\end{array}
$$

Denoting $u^{\nu}=v^{\nu}-U^{\nu}$ as the perturbation induced by $f^{\nu}=g^{\nu}+v \partial_{y}^{2} U^{\nu},$  we get
\begin{eqnarray}
\left\{\begin{array}{ll}
U_{s}^{\nu} \partial_{x} u^{\nu}+u_{2}^{\nu} \partial_{y} U_{s}^{v} \mathbf{e}_{1}-v \Delta u^{\nu}+\nabla p^{\nu}=-u^{\nu} \cdot \nabla u^{\nu}+f^{\nu}, & (x, y) \in \mathbb{T}_{\kappa} \times \mathbb{R}_{+} ,\\
\operatorname{div} u^{\nu}=0, & (x, y) \in \mathbb{T}_{\kappa} \times \mathbb{R}_{+}, \\
\left.u^{v}\right|_{y=0}=0.
\end{array}\right.\label{+9}
\end{eqnarray}
Here $\mathbf{e}_{1}=(1,0)$.
\begin{Theorem}(\cite{GM})
There exist positive numbers $\kappa_{0}, v_{0}, \varepsilon$ such that the following statement holds for $0<\kappa \leqq \kappa_{0}$ and $0<v \leqq v_{0}:$ if $f^{v}=\mathcal{Q}_{0} f^{v}$ and $\left\|f^{\nu}\right\|_{L^{2}} \leqq$
$\varepsilon v^{\frac{3}{4}}|\log v|^{-1},$ then there exists a unique solution $\left(u^{\nu}, \nabla p^{v}\right) \in\left(X \cap W_{l o c}^{2,2}\left(\mathbb{T}_{\kappa} \times\right.\right.$
$\left.\left.\mathbb{R}_{+}\right)^{2}\right) \times L^{2}\left(\mathbb{T}_{\kappa} \times \mathbb{R}_{+}\right)^{2}$ to the problem \eqref{+9} such that
$$
\begin{array}{l}
\left\|u_{0,1}^{v}\right\|_{L^{\infty}}+v^{\frac{1}{4}}\left\|\partial_{y} u_{0,1}^{v}\right\|_{L^{2}} \\
\quad+\sum_{n \neq 0}\left\|u_{n}^{\nu}\right\|_{L^{\infty}}+v^{-\frac{1}{4}}\left\|\mathcal{Q}_{0} u^{\nu}\right\|_{L^{2}}+v^{\frac{1}{4}}\left\|\nabla \mathcal{Q}_{0} u^{\nu}\right\|_{L^{2}} \leqq \frac{C|\log v|^{\frac{1}{2}}}{v^{\frac{1}{4}}}\left\|f^{\nu}\right\|_{L^{2}},
\end{array}
$$
where $\mathcal{Q}_{0}=(I-\mathcal{P}_0)$, $P_n$ is the orthogonal projection on the $n$-th Fourier mode in variable $x$ and $I$ is the identity operator.
\end{Theorem}

Ding and Li  (\cite{DL}) in 2020 were concerned with the validity of the Prandtl boundary layer theory in the inviscid limit of the
steady incompressible Navier-Stokes equations. Under the symmetry assumption,
they established the validity of the Prandtl boundary layer expansions and the error estimates.
The convergence rate as $\varepsilon \rightarrow 0$ is also given.

In \cite{DL}, they considered the following steady incompressible Navier-Stokes equations
\begin{equation}
\left\{\begin{array}{l}
U U_{X}+V U_{Y}+P_{X}=\varepsilon U_{X X}+\varepsilon U_{Y Y}, \\
U V_{X}+V V_{Y}+P_{Y}=\varepsilon V_{X X}+\varepsilon V_{Y Y}, \\
U_{X}+V_{Y}=0,
\end{array}\right.\label{+10}
\end{equation}
in the domain
$$
\Omega:=\{(X, Y) \mid 0 \leq X \leq L, 0 \leq Y \leq 2\},
$$
with moving boundary conditions
$$
U(X, 0)=U(X, 2)=u_{b}>0, V(X, 0)=V(X, 2)=0.
$$

Under the transformation, the system \eqref{+10} can be rewritten as
\begin{equation}
\left\{\begin{array}{l}
U^{\varepsilon} U_{x}^{\varepsilon}+V^{\varepsilon} U_{y}^{\varepsilon}+P_{x}^{\varepsilon}=U_{y y}^{\varepsilon}+\varepsilon U_{x x}^{\varepsilon}, \\
U^{\varepsilon} V_{x}^{\varepsilon}+V^{\varepsilon} V_{y}^{\varepsilon}+P_{y}^{\varepsilon} / \varepsilon=V_{y y}^{\varepsilon}+\varepsilon V_{x x}^{\varepsilon}, \\
U_{x}^{\varepsilon}+V_{y}^{\varepsilon}=0,
\end{array}\right.
\end{equation}
in the domain
$$
\Omega_{\varepsilon}:=\left\{(x, y) \mid 0 \leq x \leq L, 0 \leq y \leq \frac{1}{\sqrt{\varepsilon}}\right\},
$$
with the boundary conditions
$$
\left[U^{\varepsilon}, V^{\varepsilon}\right](x, 0)=\left[u_{b}, 0\right], \quad\left[U_{y}^{\varepsilon}, V^{\varepsilon}\right]\left(x, \frac{1}{\sqrt{\varepsilon}}\right)=[0,0].
$$
In what follows, they intended to find the exact solutions $\left[U^{\varepsilon}, V^{\varepsilon}, P^{\varepsilon}\right]$ in form of
$$
\left\{\begin{array}{l}
U^{\varepsilon}(x, y)=u_{a p p}(x, y)+\varepsilon^{\gamma+\frac{1}{2}} u^{\varepsilon}(x, y), \\
V^{\varepsilon}(x, y)=v_{a p p}(x, y)+\varepsilon^{\gamma+\frac{1}{2}} v^{\varepsilon}(x, y) ,\\
P^{\varepsilon}(x, y)=p_{a p p}(x, y)+\varepsilon^{\gamma+\frac{1}{2}} p^{\varepsilon}(x, y),
\end{array}\right.
$$
where
$$
\left\{\begin{array}{l}
u_{a p p}(x, y)=u_{e}^{0}(\sqrt{\varepsilon} y)+u_{p}^{0}(x, y)+\sqrt{\varepsilon} u_{e}^{1}(x, \sqrt{\varepsilon} y)+\sqrt{\varepsilon} u_{p}^{1}(x, y) ,\\
v_{a p p}(x, y)=v_{p}^{0}(x, y)+v_{e}^{1}(x, \sqrt{\varepsilon} y)+\sqrt{\varepsilon} v_{p}^{1}(x, y), \\
p_{a p p}(x, y)=\sqrt{\varepsilon} p_{e}^{1}(x, \sqrt{\varepsilon} y)+\sqrt{\varepsilon} p_{p}^{1}(x, y)+\varepsilon p_{p}^{2}(x, y).
\end{array}\right.
$$

For convenience, denote
$z:=\sqrt{\varepsilon} y$
boundary conditions on $\{y=0\}:$
\begin{eqnarray}\left\{
\begin{array}{lll}
u_{e}^{0}(0)+u_{p}^{0}(x, 0)=u_{b}, & u_{e}^{1}(x, 0)+u_{p}^{1}(x, 0)=0, & u^{\varepsilon}(x, 0)=0 ,\\
v_{p}^{0}(x, 0)+v_{e}^{1}(x, 0)=0, & v_{p}^{1}(x, 0)=0, & v^{\varepsilon}(x, 0)=0.
\end{array}\right.\non
\end{eqnarray}
Boundary conditions on $\left\{y=\frac{1}{\sqrt{\varepsilon}}\right\}$ are as follows:
\begin{eqnarray}\left\{
\begin{array}{lll}
u_{p y}^{0}\left(x, \frac{1}{\sqrt{\varepsilon}}\right)=0, & u_{e z}^{1}(x, 1)=0, u_{p y}^{1}\left(x, \frac{1}{\sqrt{\varepsilon}}\right)=0, & u_{y}^{\varepsilon}\left(x, \frac{1}{\sqrt{\varepsilon}}\right)=0, \\
v_{p}^{0}\left(x, \frac{1}{\sqrt{\varepsilon}}\right)=0, & v_{e}^{1}(x, 1)=0, v_{p}^{1}\left(x, \frac{1}{\sqrt{\varepsilon}}\right)=0, & v^{\varepsilon}\left(x, \frac{1}{\sqrt{\varepsilon}}\right)=0.
\end{array}\right.\non
\end{eqnarray}
Boundary conditions on $\{x=0\}$ are the following:
\begin{eqnarray}\left\{
\begin{array}{lll}
u_{p}^{0}(0, y)=\bar{u}_{0}(y), \quad u_{e}^{1}(0, z)=u_{b}^{1}(z), \quad u_{p}^{1}(0, y)=\bar{u}_{1}(y), \quad u^{\varepsilon}(0, y)=0, \\
v_{e}^{1}(0, z)=V_{b 0}(z), \quad v^{\varepsilon}(0, y)=0.
\end{array}\right.\non
\end{eqnarray}
Boundary conditions on $\{x=L\}$ are:
$$
v_{e}^{1}(L, z)=V_{b L}(z), \quad\left[p^{\varepsilon}-2 \varepsilon u_{x}^{\varepsilon}, u_{y}^{\varepsilon}+\varepsilon v_{x}^{\varepsilon}\right](L, y)=0.
$$
For the existence of the Euler corrector $\left[u_{e}^{1}, v_{e}^{1}, p_{e}^{1}\right],$ it is necessary  to impose the following compatibility conditions:
$$
V_{b 0}(0)=-v_{p}^{0}(0,0), \quad V_{b L}(0)=-v_{p}^{0}(L, 0), \quad V_{b 0}(1)=V_{b L}(1)=0.
$$

\begin{Theorem}(\cite{DL})
Let $u_{b}>0$ be a constant tangential velocity of the Navier-Stokes flow on the boundary $\{Y=0\}$, and let $u_{e}^{0}(Y)$ be a smooth positive Euler flow satisfies $u_{e z}^{0}(1)=0 .$ Both the boundary condition and the compatibility condition are assumed to be true. Suppose further that the positive condition $\min _{y}\left\{u_{e}^{0}(\sqrt{\varepsilon} y)+\bar{u}_{0}(y)\right\}>0$ holds. Then there exists a constant $L_{0}>0,$ which depends only on the prescribed data, such that for $0<L \leq L_{0}$ and $\gamma \in\left(0, \frac{1}{5}\right),$ the asymptotic expansion stated before is a solution to equations on $\Omega_{\varepsilon}$ together with the corresponding boundary conditions.  The approximate solutions appearing in the expansions
are constructed in sections 2-5 of \cite{DL}, in which the remainder solutions $\left[u^{\varepsilon}, v^{\varepsilon}\right]$ satisfy the estimate
$$
\left\|\nabla_{\varepsilon} u^{\varepsilon}\right\|_{L^{2}\left(\Omega_{\varepsilon}\right)}+\left\|\nabla_{\varepsilon} v^{\varepsilon}\right\|_{L^{2}\left(\Omega_{\varepsilon}\right)}+\left\|u^{\varepsilon}\right\|_{L^{\infty}\left(\Omega_{\varepsilon}\right)}+\sqrt{\varepsilon}\left\|v^{\varepsilon}\right\|_{L^{\infty}\left(\Omega_{\varepsilon}\right)} \leq C_{0},
$$
where $\Omega_{\varepsilon}=\{(x,y)|0\leq x\leq L,\ 0\leq y\leq \frac{1}{\sqrt\varepsilon}\}.$
\end{Theorem}

Wang, Yang and  Zhang (\cite{WWZ})  in 2020 considered 2D incompressible Navier-Stokes equations in a thin domain when the depth of the domain and the viscosity coefficient converges to zero simultaneously in a related way:
\begin{equation}
\left\{\begin{array}{l}
\partial_{t} U+U \cdot \nabla U-\varepsilon^{2}\left(\partial_{x}^{2}+\eta \partial_{y}^{2}\right) U+\nabla P=0 \quad \text { in } \mathcal{S}^{\varepsilon} \times(0, \infty), \\
\operatorname{div} U=0, \\
\left.U\right|_{y=0}=\left.U\right|_{y=\varepsilon}=0,
\end{array}\right.\label{+11}
\end{equation}
where $\mathcal{S}^{\varepsilon}=\{(x, y) \in \mathbb{T} \times \mathbb{R}: 0<y<\varepsilon\}$, and $U(t, x, y), P(t, x, y)$ stand for the velocity and pressure
function respectively and $\eta$ is a positive constant independent of $\varepsilon$. The system \eqref{+11} is prescribed with the initial data of the form
$$
\left.U\right|_{t=0}=\left(u_{0}\left(x, \frac{y}{\varepsilon}\right), \varepsilon v_{0}\left(x, \frac{y}{\varepsilon}\right)\right)=U_{0}^{\varepsilon}
$$
Now rescale $(U, P)$ as follows
$$
U(t, x, y)=\left(u^{\varepsilon}\left(t, x, \frac{y}{\varepsilon}\right), \varepsilon v^{\varepsilon}\left(t, x, \frac{y}{\varepsilon}\right)\right) \quad \text { and } \quad P(t, x, y)=p^{\varepsilon}\left(t, x, \frac{y}{\varepsilon}\right).
$$

Then the system \eqref{+11} reduces to the following scaled anisotropic Navier-Stokes system:
\begin{equation}
\left\{\begin{array}{l}
\partial_{t} u^{\varepsilon}+u^{\varepsilon} \partial_{x} u^{\varepsilon}+v^{\varepsilon} \partial_{y} u^{\varepsilon}-\varepsilon^{2} \partial_{x}^{2} u^{\varepsilon}-\eta \partial_{y}^{2} u^{\varepsilon}+\partial_{x} p^{\varepsilon}=0, \quad \text { in } \mathcal{S} \times(0, \infty), \\
\varepsilon^{2}\left(\partial_{t} v^{\varepsilon}+u^{\varepsilon} \partial_{x} v^{\varepsilon}+v^{\varepsilon} \partial_{y} v^{\varepsilon}-\varepsilon^{2} \partial_{x}^{2} v^{\varepsilon}-\eta \partial_{y}^{2} v^{\varepsilon}\right)+\partial_{y} p^{\varepsilon}=0 ,\quad \text { in } \mathcal{S} \times(0, \infty), \\
\partial_{x} u^{\varepsilon}+\partial_{y} v^{\varepsilon}=0, \quad \text { in } \mathcal{S} \times(0, \infty), \\
\left.\left(u^{\varepsilon}, v^{\varepsilon}\right)\right|_{y=0,1}=0. \\
\left.\left(u^{\varepsilon}, v^{\varepsilon}\right)\right|_{t=0}=\left(u_{0}, v_{0}\right), \quad \text { in } \mathcal{S},
\end{array}\right.
\end{equation}
where $\mathcal{S}=\{(x, y) \in \mathbb{T} \times(0,1)\}$. This is a classical model in geophysical fluid, where the vertical dimension of the domain is very small compared with the horizontal dimension of the domain. For simplicity, they took $\eta=1$ in the sequel and denoted $\Delta_{\varepsilon}=\varepsilon^{2} \partial_{x}^{2}+\partial_{y}^{2}$.

Formally, taking $\varepsilon \rightarrow 0$,  we derive the hydrostatic Navier-Stokes/Prandtl system
\begin{equation}
\left\{\begin{array}{l}
\partial_{t} u^{p}+u^{p} \partial_{x} u^{p}+v^{p} \partial_{y} u^{p}-\partial_{y}^{2} u^{p}+\partial_{x} p^{p}=0 \quad \text { in } \mathcal{S} \times(0, \infty), \\
\partial_{y} p^{p}=0 \quad \text { in } \mathcal{S} \times(0, \infty), \\
\partial_{x} u^{p}+\partial_{y} v^{p}=0 \quad \text { in } \mathcal{S} \times(0, \infty), \\
\left.\left(u^{p}, v^{p}\right)\right|_{y=0,1}=0, \\
\left.u^{p}\right|_{t=0}=u_{0} \quad \text { in } \mathcal{S}.
\end{array}\right.
\end{equation}

Consider the initial data of the form
$$
u^{\varepsilon}(0, x, y)=u_{0}(x, y), \quad v^{\varepsilon}(0, x, y)=v_{0}(x, y),
$$
which satisfy the compatibility conditions
\begin{eqnarray}\left\{
\begin{array}{l}
\partial_{x} u_{0}+\partial_{y} v_{0}=0, \quad u_{0}(x, 0)=u_{0}(x, 1)=v_{0}(x, 0)=v_{0}(x, 1)=0, \\
\int_{0}^{1} \partial_{x} u_{0} d y=0,\left.\quad \partial_{y}^{2} u_{0}\right|_{y=0,1}=\int_{0}^{1}\left(-\partial_{x} u_{0}^{2}+\partial_{y}^{2} u_{0}\right) d y-\int_{\mathcal{S}} \partial_{y}^{2} u_{0},
\end{array}\right.\non\end{eqnarray}
and the convex condition
$$
\inf _{\mathcal{S}} \partial_{y}^{2} u_{0}=2 \delta_{0}>0.
$$
The authors of \cite{WWZ} further assumed that initial data fall  into the Gevrey class with the bound
$$
\left\|\partial_{y} u_{0}\right\|_{X_{\sigma, \tau_{0}}^{N_{0}}}+\left\|\partial_{y}^{3} u_{0}\right\|_{X_{\sigma, \tau_{0}}^{N_{0}-4}}=M<+\infty.
$$
Here the Gevrey class norm $\|\cdot\|_{X_{\sigma, \tau}^{r}}$ is defined by
$$
\begin{array}{r}
\|f\|_{X_{\sigma, \tau}^{r}}^{2}=\left\|e^{\tau\left\langle D_{x}\right\rangle^{\sigma}} f\right\|_{H^{r, 0}}^{2},
\text { with }\|f\|_{H^{r, s}}=\|\| f \|_{H_{x}^{r}(\mathbb{T})} \|_{H_{y}^{s}(0,1)}.
\end{array}
$$

\begin{Theorem}(\cite{WWZ})
Let the initial data $u_{0}$ satisfy the assumptions with $\sigma \in\left[\frac{8}{9}, 1\right], \tau_{0}>0$ and $N_{0} \geq 10$. Then there exist $T>0$ and a unique solution $u^{p}$ of the former equation which satisfies
$$
\begin{array}{l}
\sup _{t \in[0, T]}\left(\left\|\partial_{y} u^{p}(t)\right\|_{X_{\sigma, \tau}^{N_{0}-1}}+\left\|\partial_{y}^{3} u^{p}(t)\right\|_{X_{\sigma, \tau}^{N_{0}-5}}\right)<+\infty ,\\
\sup _{t \in[0, T] \times \mathcal{S}} \partial_{y}^{2} u^{p}>\delta_{0} .
\end{array}
$$

\end{Theorem}
\begin{Theorem}(\cite{WWZ})
Let initial data $u_{0}$ satisfies the assumptions with $\sigma \in\left[\frac{8}{9}, 1\right], \tau_{0}>0$ and $N_{0} \geq 10$. Then there exists a unique solution of the Navier-Stokes equations in $[0, T],$ which satisfies
$$
\left\|\left(u^{\varepsilon}-u^{p}, \varepsilon v^{\varepsilon}-\varepsilon v^{p}\right)\right\|_{L_{x, y}^{2} \cap L_{x, y}^{\infty}} \leq C \varepsilon^{2},
$$
where $\left(u^{p}, v^{p}\right)$ is given by the former theorem and $C>0$ is a constant independent of $\varepsilon$.
\end{Theorem}

At last, the following is the   zero viscosity limit of the 3D incompressible Navier-Stokes equations.

Fei, Tao and Zhang (\cite{FTZ1}) in 2018 investigated  the zero viscosity limit of the incompressible Navier-Stokes equations with non-slip boundary condition in $\mathbb{R}^3_+$ for the initial vorticity located away from the boundary by the energy method.
In \cite{FTZ1}, they were concerned with the zero-viscosity limit of the incompressible Navier-Stokes equations in the half-space  $\mathbb{R}^3_+$:

\begin{eqnarray}\left\{\begin{array}{ll}
u^\varepsilon_t-\varepsilon^2\Delta u^\varepsilon+u^\varepsilon\cdot\nabla_x u^\varepsilon+v^\varepsilon\cdot\nabla_y u^\varepsilon+\nabla_x p=0,\\
v^\varepsilon_t-\varepsilon^2\Delta v^\varepsilon+u^\varepsilon\cdot\nabla_x v^\varepsilon+v^\varepsilon\cdot\partial_y v^\varepsilon+\nabla_x p=0,\\
\nabla_x\cdot u^\varepsilon+\partial_yv^\varepsilon=0,\\
(u^\varepsilon,v^\varepsilon)(t,x,0)=(0,0).
\end{array}\right.\label{24.1}\end{eqnarray}
Here $x=(x_1,x_2)\in \mathbb{T}$ and $y\in \mathbb{R}$ denote the tangential and the vertical components of the space variable, respectively. $\nabla_x=(\partial_{x_1},\partial_{x_2})$ denotes the tangential gradient; $u=(u^1,u^2)(t,x,y)$ and $v=v(t,x,y)$  denote the tangential and the vertical velocities;
 $\varepsilon^2$ is the viscosity coefficient.
For simplicity, we consider the initial data of the form $$u^\varepsilon(0,x,y)=u_0(x,y),\ v^\varepsilon(0,x,y)=v_0(x,y),$$
which satisfies
\begin{eqnarray}
\nabla_x\cdot u_0^\varepsilon+\partial_yv_0^\varepsilon=0,\ u_0(x,0)=0,\ v_0(x,0)=0,\label{24.2}
\end{eqnarray}
and the initial vorticity $\omega_0=curl(u_0,v_0)$ satisfies

\begin{eqnarray}
 2d_0=dist(supp\omega_0,{y=0})>0.\label{24.3}
 \end{eqnarray}

Without loss of generality, we take $d_0=1$.

\begin{Theorem}( \cite{FTZ1})
There exist $T>0$ and $C>0$, independent of $\varepsilon$, such that for any $(u_0, v_0)\in H^{30}(\mathbb{R}^3_+)$ satisfying
\eqref{24.2}-\eqref{24.3}, there exists a unique solution $(u_\varepsilon, v_\varepsilon)$ of the Navier-Stokes equations \eqref{24.1} in $[0, T]$, which satisfies
\begin{eqnarray}\left\{\begin{array}{ll}
\sup\limits_{0\leq t\leq T}\|u^\varepsilon(t,x,y)-u^e(t,x,y)-u^p(t,x,y/\varepsilon)\|_{L^2\cap L^2(\mathbb{R}_+^3)}\leq C\varepsilon,\\
\sup\limits_{0\leq t\leq T}\|v^\varepsilon(t,x,y)-v^e(t,x,y)-v^p(t,x,y/\varepsilon)\|_{L^2\cap L^2(\mathbb{R}_+^3)}\leq C\varepsilon.
\end{array}\right.\end{eqnarray}
\end{Theorem}

\section{Ill-posedness of the 2D Prandtl Equations}

  E and  Engquist (\cite{EE}) in 1997 considered the  unsteady boundary layer equations   (\ref{1.1.6}) in the domain
 $ \mathbb{R}^{2}=\{(x,y)|x\in \mathbb{R}, y\in \mathbb{R}_{+}\}$. Let $U=p=0$, they  proved  global blowup of solutions.
  \begin{Theorem}(\cite{EE})\label{t.3}
Assume that $u_{0}$ takes the form $u_{0}(x,y)=-xb_{0}(x,y)$ where $a_{0}(\cdot)=-u_{0x}(0,\cdot)$ satisfies the conditions of nonnegative, compactly
supported initial data such that $\int_{0}^{\infty}(\frac{1}{2} a_{y}^{2}-\frac{1}{4} a^{3})dy $; then smooth solutions of (\ref{1.1.6})  do not exist globally
in time.

\end{Theorem}

Gargano,   Sammartino and   Sciacca (\cite{GSS}) in 2009 investigated the numerical evidence of the ill-posedness of the Prandtl equations in $H^1$ and led to the solution blow-up whose time can be made arbitrarily short within the class.

The authors of \cite{GSS} considered  the Cauchy problem for the Prandtl equation. This
problem is known to be well-posed for analytic data, or for data with monotonicity
properties.  G$\acute{e}$rard-Varet and   Dormy (\cite{GD}) in 2010 proved the linearly ill-posed in Sobolev type spaces.
The Prandtl equations reduce to
\begin{eqnarray}\left\{\begin{array}{ll}
u_t+u\partial_xu+v\partial_yu-\partial_{y}^2u=0,\\
\partial_xu+\partial_yv=0,\\
(u,v)|_{y=0}=0, \lim\limits_{y\rightarrow+\infty}u=0.
\end{array}\right.\label{33.1}\end{eqnarray}
Let $u_s=u_s(t,y)$ be a smooth solution of the heat equation
\begin{eqnarray}\left\{\begin{array}{ll}
\partial_t u_s-\partial_{y}^2u_s=0,\\
u_s|_{y=0}=0, \lim\limits_{y\rightarrow+\infty}u=U_0,\\
u_s|_{t=0}=U_s(y)
\end{array}\right.\label{33.2}\end{eqnarray}
with good decay as  $y\rightarrow +\infty$.  Clearly, the shear velocity profile $(u_s,v_s)=(u_s(t,y),0)$ satisfies the system \eqref{33.1}.  Consider the linearization around $(u_s,v_s)$,  that is,

\begin{eqnarray}\left\{\begin{array}{ll}
u_t+u_s\partial_xu+v\partial_yu_s-\partial_{y}^2u=0,\\
\partial_xu+\partial_yv=0,\\
(u,v)|_{y=0}=0, \lim_{y\rightarrow+\infty}u=0.
\end{array}\right.\label{33.3}\end{eqnarray}
In the sequel, \cite{GSS} studied well-posedness properties of problem \eqref{33.3},  for a certain class of velocities $u_s$. In this view, \cite{GSS} introduced the following functional spaces:

$$W_\alpha^{s,\infty}(\mathbb{R}_+):=\{f=f(y), e^{\alpha y}f\in W^{s,\infty}(\mathbb{R}_+),\ \forall \alpha,s\geq0\},$$
with $\|f\|_{W_\alpha^{s,\infty}}:=\|e^{\alpha y}f\|_{W_\alpha^{s,\infty}}$, and

$$E_{\alpha,\beta}:=\Big\{u=u(x,y)=\sum_{k\in\mathbb{Z}}\bar{u}^k(y)e^{ikx},\ \|\bar{u}^k\|_{W_\alpha^{s,\infty}}\leq C_{\alpha,\beta}e^{-\beta|k|}\Big\}, \forall \alpha,\beta>0,$$
with $\|u\|_{E_{\alpha,\beta}}:=\sup_ke^{\beta|k|}\|\bar{u}^k\|_{W_\alpha^{s,\infty}}.$

\begin{Theorem}(\cite{GD})
 Let $u_s\in C^0(\mathbb{R}_+;W_0^{4,\infty})\cap C^1(\mathbb{R}_+;W_0^{2,\infty})$.  Assume that the initial velocity has a non-degenerate critical point over $\mathbb{R}^+$. Then, there exists $\sigma>0$,  such that for all  $\delta>0$,
$$\sup_{0\leq s\leq t\leq \delta}\|e^{-\sigma(t-s)\sqrt|\partial_x|}T(t,s)\|_{L(H^m,H^{m-\mu})}=+\infty, \forall \alpha, m\geq0, \mu\in [0,\frac{1}{2}).$$

Moreover, one can find solutions $u_s$ of \eqref{33.2} and  $\sigma>0$,  such that for all  $\delta>0$.

$$\sup_{0\leq s\leq t\leq \delta}\|e^{-\sigma(t-s)\sqrt|\partial_x|}T(t,s)\|_{L(H^m,H^{m_2})}=+\infty, \forall  m_1,m_2\geq0,$$
where $H^m=H^m(\mathbb{T}_x, W_\alpha^{0,\infty}(\mathbb{R}_y^+)),\ m\geq0$.
\end{Theorem}

G$\acute{e}$rard-Varet and Dormy (\cite{GD}) in 2010  were concerned with the Cauchy problem for the Prandtl equation,  proved  that it is linearly ill-posed in Sobolev type spaces.The key of the analysis is the construction, at high tangential frequencies, of unstable quasimodes for the linearization around solutions with non-degenerate critical points.

 The authors of \cite{GD} restricted themselves to $(x, Y) \in \mathbb{T} \times \mathbb{R}_{+},$ and $\mathbf{u}^{0}=0$.  To lighten notations, they wrote $y$ instead of $Y$. The Prandtl equation comes down to
\begin{eqnarray}
\left\{\begin{aligned}
\partial_{t} u+u \partial_{x} u+v \partial_{y} u-\partial_{y}^{2} u=0, & \text { in } \mathbb{T} \times \mathbb{R}_{+} \\
\partial_{x} u+\partial_{y} v=0, & \text { in } \mathbb{T} \times \mathbb{R}_{+} \\
\left.(u, v)\right|_{y=0}=(0,0), & \lim _{y \rightarrow+\infty} u=0.
\end{aligned}\right.\label{+12}
\end{eqnarray}
Let $u_{s}=u_{s}(t, y)$ be a smooth solution of the heat equation
$$e
\partial_{t} u_{s}-\partial_{y}^{2} u_{s}=0,\left.\quad u_{s}\right|_{y=0}=0,\left.\quad u_{s}\right|_{t=0}=U_{s}
$$
with good decay as $y \rightarrow+\infty$. Clearly, the shear velocity profile $\left(u_{s}, v_{s}\right)=\left(u_{s}(t, y), 0\right)$ satisfies the system \eqref{+12}.  The authors of \cite{GD}  considered the linearization around $\left(u_{s}, v_{s}\right),$ that is,
\begin{eqnarray}
\left\{\begin{aligned}
\partial_{t} u+u_{s} \partial_{x} u+v \partial_{y} u_{s}-\partial_{y}^{2} u=0, & \text { in } \mathbb{T} \times \mathbb{R}_{+}, \\
\partial_{x} u+\partial_{y} v=0, & \text { in } \mathbb{T} \times \mathbb{R}_{+}, \\
\left.(u, v)\right|_{y=0}=(0,0), & \lim _{y \rightarrow+\infty} u=0,
\end{aligned}\label{1.2}   \right.
\end{eqnarray}
and introduced the following functional spaces:
$$
W_{\alpha}^{s, \infty}\left(\mathbb{R}_{+}\right):=\left\{f=f(y), \quad e^{\alpha y} f \in W^{s, \infty}\left(\mathbb{R}_{+}\right)\right\}, \quad \forall \alpha, s \geq 0
$$
with $\|f\|_{W_{\alpha}^{s, \infty}}:=\left\|e^{\alpha y} f\right\|_{W^{s, \infty}},$ and
$$
E_{\alpha, \beta}:=\left\{u=u(x, y)=\sum_{k \in \mathbb{Z}} \hat{u}^{k}(y) e^{i k x}, \quad\left\|\hat{u}^{k}\right\|_{W_{\alpha}^{0, \infty}} \leq C_{\alpha, \beta} e^{-\beta|k|}, \forall k\right\}, \quad \forall \alpha, \beta>0
$$
with $\|u\|_{E_{\alpha, \beta}}:=\sup _{k} e^{\beta|k|}\left\|\hat{u}^{k}\right\|_{W_{\alpha}^{0, \infty}}$.
\begin{Theorem}(\cite{GD})
\textbf{(Ill-posedness in the Sobolev Setting)}

i) Let $u_{s} \in C^{0}\left(\mathbb{R}_{+} ; W_{\alpha}^{4, \infty}\left(\mathbb{R}_{+}\right)\right) \cap C^{1}\left(\mathbb{R}_{+} ; W_{\alpha}^{2, \infty}\left(\mathbb{R}_{+}\right)\right) .$ Assume that the initial velocity has a
non-degenerate critical point over $\mathbb{R}_{+} .$ Then, there exists $\sigma>0,$ such that for all $\delta>0,$
$$
\sup _{0 \leq s \leq t \leq \delta}\left\|e^{-\sigma(t-s) \sqrt{\left|\partial_{x}\right|}} T(t, s)\right\|_{\mathcal{L}\left(H^{m}, H^{m-\mu}\right)}=+\infty, \quad \forall m \geq 0, \mu \in[0,1 / 2).
$$

ii) Moreover, one can find solutions $u_{s}$ of (\ref{1.2}) and $\sigma>0$ such that: for all $\delta>0$,
$$
\sup _{0 \leq s \leq t \leq \delta}\left\|e^{-\sigma(t-s) \sqrt{\left|\partial_{x}\right|}} T(t, s)\right\|_{\mathcal{L}\left(H^{m_{1}}, H^{m_{2}}\right)}=+\infty, \quad \forall m_{1}, m_{2} \geq 0.
$$
\end{Theorem}

G$\acute{e}$rard-Varet and  Nguyenin (\cite{GN5}) in 2012 established various ill-posedness results for the Prandtl equation. By considering perturbations of stationary shear flows, they showed that for some linearizations of the Prandtl equation and some $C^{\infty}$ initial data, local in time $C^{\infty}$ solutions do not exist. At the nonlinear level, they proved that if a flow exists in the Sobolev setting, it cannot be Lipschitz continuous. Besides ill-posedness in time, they also established some ill-posedness in space, that casts some light on the results obtained by Oleinik for monotonic data.

They were concerned in \cite{GN5} with the famous Prandtl equations:
\begin{equation}
\left\{\begin{aligned}
\partial_{t} u+u \partial_{x} u+v \partial_{y} u-\partial_{y}^{2} u+\partial_{x} P &=f, \quad y>0, \\
\partial_{x} u+\partial_{y} v &=0, \quad y>0, \\
u=v &=0, \quad y=0, \\
\lim _{y \rightarrow+\infty} u &=U(t, x).
\end{aligned}\right.
\end{equation}

They had two settings:

The initial value problem (IVP):
$(t, x, y)$ in $[0, T) \times \mathbb{T} \times \mathbb{R}_{+},\left.u\right|_{t=0}=u_{0}(x, y).$

The boundary value problem $(B V P)$:
$(t, x, y)$ in $\mathbb{T} \times[0, X) \times \mathbb{R}_{+},\left.u\right|_{x=0}=u_{1}(t, y).$

The system satisfied by the perturbation $(u, v)$ of $\left(u_{s}, 0\right)$ reads:
\begin{equation}
\left\{\begin{aligned}
\partial_{t} u+u_{s} \partial_{x} u+u_{s}^{\prime} v+u \partial_{x} u+v \partial_{y} u-\partial_{y}^{2} u &=0, \quad y>0, \\
\partial_{x} u+\partial_{y} v &=0, \quad y>0, \\
u=v &=0, \quad y=0.
\end{aligned}\right.
\end{equation}

\begin{Theorem}(\cite{GN5})
(Non-existence of Solutions for Linearized Prandtl Equations).

There exists a shear flow $u_{s}$ with $u_{s}-U \in C_{c}^{\infty}\left(\mathbb{R}_{+}\right)$ such that: for all $T>0,$ there exists an initial datum $u_{0}$ satisfying

i) $e^{y} u_{0} \in H^{\infty}\left(\mathbb{T} \times \mathbb{R}_{+}\right);$

ii) The IVP has no distributional solution u with
$$
u\in L^{\infty}((0, T); L^{2}(\mathbb{T} \times \mathbb{R}_{+})), \quad \partial_{y} u \in L^{2}((0, T)\times \mathbb{T}\times\mathbb{R}_{+}).
$$
\end{Theorem}

\begin{Theorem}(\cite{GN5})(No Lipschitz Continuity of the Flow)
There exists a shear flow $u_{s}$ with $u_{s}-U \in C_{c}^{\infty}\left(\mathbb{R}_{+}\right)$ such that: for all $m \geq 0,$ the Cauchy problem is not locally $\left(H^{m}, H^{1}\right)$ Lipschitz well-posed.
\end{Theorem}

Ding (\cite{D}) in 2012 concerned the linearized Prandtl equation around general stationary solutions with nondegenerate critical points depending on $x$ which could be considered as the time-periodic
solutions and showed some ill-posedness.

Consider the 2D Navier-Stokes equations on a half-space:
\begin{equation}
\left\{\begin{array}{l}
\partial_{t} u^{\nu}+u^{\nu} \partial_{x} u^{\nu}+v^{\nu} \partial_{y} u^{\nu}+\partial_{x} p^{\nu}-v \Delta u^{\nu}=0, \\
\partial_{t} v^{\nu}+u^{\nu} \partial_{x} v^{\nu}+v^{v} \partial_{y} v^{\nu}+\partial_{y} p^{\nu}-v \Delta v^{\nu}=0, \\
\partial_{x} u^{\nu}+\partial_{y} v^{v}=0, \\
\left.\left(u^{\nu}, v^{v}\right)\right|_{y=0}=(0,0).
\end{array}\right.
\end{equation}

When $\nu \rightarrow 0$, a natural question is that does the solution $\left(u^{v}, v^{v}\right)$ convergence to the solution of Euler equation:
\begin{equation}
\left\{\begin{array}{ll}
\partial_{t} u^{E}+u^{E} \partial_{x} u^{E}+v^{E} \partial_{y} u^{E}+\partial_{x} p^{E}=0, \\
\partial_{t} v^{E}+u^{E} \partial_{x} v^{E}+v^{E} \partial_{y} v^{E}+\partial_{y} p^{E}=0, \\
\partial_{x} u^{E}+\partial_{y} v^{E}=0, \\
v^\nu|_{y=0}=0.\end{array}\right.\non\end{equation}
As there is the no-slip condition: $\left.u^{v}\right|_{y=0}=0$ in Navier-Stokes equations, the transition from zero velocity at the boundary to the full magnitude at some distance from it take place in a very thin layer.

It is interesting to consider the linearized equation around stationary solution $(u_0,v_0)$ where $u_0$ has non-degenerate critical points, that is
\begin{eqnarray}
\left\{\begin{aligned}
\partial_{t} u+u_0 \partial_{x} u+v_0 \partial_{y} u+u \partial_{x} u_0+v \partial_{y} u_0-\partial_{y}^{2} u=0, \\
\partial_{x} u+\partial_{y} v=0,\\
u|_{y=0}=v|_{y=0}=0,\ u|_{x=0}=u_1,\ \lim\limits_{y\rightarrow+\infty}u=0.
\end{aligned}\right.\label{+13}
\end{eqnarray}

The author of \cite{D} introduced the following function spaces:
$$
W_{\alpha}^{s, \infty}:=\left\{f=f(y), e^{\alpha y} f \in W^{s, \infty}\left(\mathbb{R}_{+}\right)\right\}
$$
with the norm: $\|f\|_{W_{\alpha}^{s, \infty}}=\left\|e^{\alpha y} f\right\|_{W^{s, \infty}}$,
$$
H_{\beta}^{m}:=H^{m}\left(\mathbb{T}_{t}, W_{\beta}^{0, \infty}\left(\mathbb{R}_{y}^{+}\right)\right)
$$
and
$$
\bar{H}_{\beta}^{m}:=H_{\beta}^{m} \cap C^{1}\left(\mathbb{T}_{t}, W_{\beta}^{2, \infty}\left(\mathbb{R}_{y}^{+}\right)\right),
$$where $\mathbb{T}_{t}$ denotes the time domain.
\begin{Theorem}(\cite{D})
 Let $u_{0}-U \in C^{0}\left(\left[0, X_{0}\right) ; W_{\alpha}^{4, \infty}\left(\mathbb{R}_{+}\right)\right) \cap C^{1}\left(\left[0, X_{0}\right) ; W_{\alpha}^{2, \infty}\left(\mathbb{R}_{+}\right)\right),\left.u_{0}\right|_{x=0}$ has a non-degenerate critical point. If there
exists $X>0$, such that for every $u_{1}(t, y) \in \bar{H}_{\beta}^{m}$, with $\beta<\alpha$, the equations \eqref{+13} have a unique solution, let us denote $u(x, \cdot):=\mathfrak{X}(x, \xi) u_{1}$, where $u(x, \cdot)$ is the solution of the problem \eqref{+13} with $\left.u\right|_{x=\xi}=u_{1},$ then there exists a $ \delta>0,$ such that for every $\epsilon>0$,
$$
\sup _{0 \leqslant \xi \leqslant x \leqslant \epsilon}\left\|e^{-\delta(x-\xi) \sqrt{\left|\partial_{t}\right|}} \mathfrak{X}(x, \xi)\right\|_{\mathcal{L}\left(H_{\beta}^{m}, H_{\beta}^{m-\sigma}\right)}=\infty, \quad \forall m \geqslant 0, \sigma \in\left[0, \frac{1}{2}\right).
$$
If $\lim\limits _{y \rightarrow \infty} u_{0}(x, y)=C,$ where $C \geqslant 0$ is a constant, the result is valid for $\alpha=\beta$.

\end{Theorem}

Liu and Yang (\cite{ly}) in 2017 studied the ill-posedness of the following Prandtl equations \eqref{77.1} in Sobolev spaces around a shear flow with general decay:
\begin{eqnarray}\left\{\begin{array}{ll}
u_t+u\partial_xu+v\partial_yu-\partial_{y}^2u=0,\\
\partial_xu+\partial_yv=0,\\
(u,v)|_{y=0}=0, \lim\limits_{y\rightarrow+\infty}u=U_0.
\end{array}\right.\label{77.1}\end{eqnarray}
Note that \eqref{77.1} has a special shear flow solution $(u_s (t, y),0)$, where the function $u_s (t, y)$ is a smooth solution to the following heat equation:
\begin{eqnarray}\left\{\begin{array}{ll}
\partial_t u_s-\partial_{y}^2u_s=0,\\
u_s|_{y=0}=0, \lim\limits_{y\rightarrow+\infty}u=U_0,\\
u_s|_{t=0}=U_s(y)
\end{array}\right.\label{77.2}\end{eqnarray}
with an initial shear layer $U_s(y)$.
Denote by $T(t,s)$ the linear solution operator $$T(t,s)u_0:= u(t,\cdot).$$

\begin{Theorem} (\cite{ly}) \label{ly.1}
Let $u_s(t,y)$  be the solution of the problem \eqref{77.2} satisfying
$$u_s-U_0\in C^0(\mathbb{R}_+;W_0^{4,\infty})\cap C^1(\mathbb{R}_+;W_0^{2,\infty})$$
and assume that the initial shear layer $U_s(y)$ has a non-degenerate critical point in $\mathbb{R}_+$. Then, there exists
$\sigma > 0$ such that for all $\delta > 0$,
$$\sup_{0\leq s\leq t\leq \delta}\|e^{-\sigma(t-s)\sqrt{|\partial_x|}}T(t,s)\|_{\mathcal{L}(H^m_\alpha,H^{m-\mu}_0)}=+\infty, \forall \alpha, m\geq0, \mu\in [0,\frac{1}{2}).$$
\end{Theorem}
\begin{Theorem}
Under the assumptions of Theorem \ref{ly.1}, it holds that for any $\delta>0$ and $\alpha,m\geq0$
$$\sup_{0\leq s\leq t\leq\delta}\|T(t,s)\|_{\mathcal{L}(H^m_\alpha,H^0_0)}=+\infty.$$
\end{Theorem}

 Kukavica,   Vicol and  Wang (\cite{KVW}) in 2017  proved a numerical conjecture by rigorously establishing the finite time blowup of the boundary layer thickness.They considered the 2D Prandtl boundary layer equations for the unknown velocity field $(u, v) = (u(t, x, y), v(t, x, y))$:
\begin{eqnarray}\left\{\begin{array}{ll}
u_t+u\partial_xu+v\partial_yu-\partial_{y}^2u+\partial_xp^{E}=0,\\
\partial_xu+\partial_yv=0,\\
(u,v)|_{y=0}=0, \lim\limits_{y\rightarrow+\infty}u=U^E,\\
u|_{t=0}=u_0.
\end{array}\right.\label{57.1}\end{eqnarray}
They (\cite{KVW})   wanted to prove the formation of finite time singularities in the
Prandtl boundary layer equations \eqref{57.1} when the underlying Euler flow is not trivial, i.e.,
when $U^E\neq0$. For this purpose, consider the Euler trace
\begin{eqnarray}\left\{\begin{array}{ll}
U^E=ksinx,\\
-\partial_xP^E=\frac{k^2}{2}sin(2x).
\end{array}\right.\label{57.2}\end{eqnarray}
\begin{Theorem}(\cite{KVW})
Consider the Cauchy problem for the
Prandtl equations \eqref{57.1}, with boundary conditions at $y = \infty$ matching \eqref{57.2},
with $k\neq0$. There exists a large class of initial conditions $(u_0, v_0)$ which are real-analytic
in $x$ and $y$, such that the unique real-analytic solution $(u, v)$ to the Prandtl equations \eqref{57.1} blows up in a finite time.
\end{Theorem}

Guo and Nguyen (\cite{GN}) in 2018 concerned nonlinear ill-posedness of the Prandtl equation and an invalidity of asymptotic boundary layer expansions of incompressible fluid flows
near a solid boundary. Their analysis was built upon recent remarkable linear ill-posedness results established by G$\acute{e}$rard-Varet and Dormy (\cite{GD}) and an analysis by
Guo and Tice (\cite{guotice}). They showed that the asymptotic boundary layer expansion is not
valid for nonmonotonic shear layer flows in Sobolev spaces. They also introduced
a notion of weak well-posedness and proved that the nonlinear Prandtl equation is
not well-posed in this sense near nonstationary and nonmonotonic shear flows.
On the other hand, they were able to verify that Oleinik's monotonic solutions are
well-posed.

The boundary layer or Prandtl equations for $(u, v)$ then reads:
\begin{eqnarray}
\left\{\begin{array}{ll}
\partial_{t} u+u \partial_{x} u+v \partial_{Y} u-\partial_{Y}^{2} u+\partial_{x} P=0, & Y>0 ,\\
\partial_{x} u+\partial_{Y} v=0, & Y>0, \\
\left.u\right|_{t=0}=u_{0}(x, y), \\
\left.u\right|_{Y=0}=\left.v\right|_{Y=0}=0, \\
\lim\limits _{Y \rightarrow+\infty} u=U(t, x)
\end{array}\right.
\end{eqnarray}
where $U=u^{0}(t, x, 0)$ and $P=P(t, x)$ are the normal velocity and pressure describing the Euler flow just outside the boundary layer and satisfy the Bernoulli's  equation
$$
\partial_{t} U+U \partial_{x} U+\partial_{x} P=0.
$$

Introduce the standard Sobolev spaces $L^{2}$ and $H^{m}, m \geq 0,$ with the usual norms:
$$
\|u\|_{L_{x, Y}^{2}}:=\left(\int_{\mathbb{T} \times \mathbb{R}_{+}}|u|^{2} d x d Y\right)^{1 / 2}, \quad\|u\|_{H_{x, Y}^{m}}:=\sum_{k=0}^{m} \sum_{i+j=k}\left\|\partial_{x}^{i} \partial_{Y}^{j} u\right\|_{L^{2}}.
$$
Regarding the validity of the asymptotic boundary layer expansion, we can write
$$
\begin{aligned}
\left(\begin{array}{c}
u^{v} \\
v^{\nu}
\end{array}\right)(t, x, y)=&\left(\begin{array}{c}
u^{0}-\left.u^{0}\right|_{y=0} \\
0
\end{array}\right)(y)+\left(\begin{array}{c}
u_{s} \\
0
\end{array}\right)\left(t, \frac{y}{\sqrt{v}}\right) \\
&+(\sqrt{v})^{\gamma}\left(\begin{array}{c}
\tilde{u}^{\nu} \\
\sqrt{v} \tilde{v}^{v}
\end{array}\right)\left(t, x, \frac{y}{\sqrt{v}}\right),
\end{aligned}
$$
and
$$
p^{\nu}(t, x, y)=(\sqrt{v})^{\gamma} \tilde{p}^{\nu}\left(t, x, \frac{y}{\sqrt{v}}\right),
$$
for shear flows $u_{s}$ and for some $\gamma>0$, where $\left(u^{0}(y), 0\right)^{T}$ is the Euler flow.

\begin{Theorem}(\cite{GN})(Invalidity of Asymptotic Expansions)
The expansion is not valid in the sense of validity of asymptotic expansions for any $\gamma>0$.
\end{Theorem}

Sun (\cite{S2}) in 2019 obtained a blow-up criterion for classical solutions to the two-dimensional Prandtl equations under a monotonicity assumption in weighted Sobolev spaces. The proof is based on nonlinear energy estimates inspired by Masmoudi and Wong (\cite{mw}) and maximum principles for parabolic equations.

Consider the following Prandtl system:
\begin{equation}
\left\{\begin{array}{ll}
\partial_{t} u+u \partial_{x} u+v \partial_{y} u-\partial_{y}^{2} u=\partial_{t} U+U \partial_{x} U, & \text { in }[0, T] \times \mathbb{T} \times \mathbb{R}_{+}, \\
\partial_{x} u+\partial_{y} v=0, & \text { in }[0, T] \times \mathbb{T} \times \mathbb{R}_{+}, \\
\left.u\right|_{t=0}=u_{0}, \quad & \text { on } \mathbb{T} \times \mathbb{R}_{+}, \\
\left.u\right|_{y=0}=\left.v\right|_{y=0}=0, & \text { on }[0, T] \times \mathbb{T}, \\
\lim\limits _{y \rightarrow \infty} u(t, x, y)=U(t, x), & \text { for all }(t, x) \in[0, T] \times \mathbb{T}.
\end{array}\right.
\end{equation}

Masmoudi and Wong (\cite{mw}) introduced the following space:
$$
H_{\sigma, \delta}^{s, \gamma}:=\left\{\omega: \mathbb{T} \times\left.\mathbb{R}_{+} \rightarrow \mathbb{R}\left|\|\omega\|_{H^{s}, \nu}<+\infty,(1+y)^{\sigma} \omega \geq \delta, \sum_{|\alpha| \leq 2}\right|(1+y)^{\sigma+\alpha_{2}} D^{\alpha} \omega\right|^{2} \leq \frac{1}{\delta^{2}}\right\},
$$
where $s \geq 4, \gamma \geq 1, \sigma>\gamma+\frac{1}{2}, \delta \in(0,1), D^{\alpha}=\partial_{x}^{\alpha_{1}} \partial_{y}^{\alpha_{2}},$ and the weighted $H^{s}$ norm $\|\cdot\|_{H^{s}, \nu}$ is defined by
$$
\|\omega\|_{H^{s y}}^{2}:=\sum_{|\alpha| \leq s}\left\|(1+y)^{\gamma+\alpha_{2}} D^{\alpha} \omega\right\|_{L^{2}\left(\mathbb{T} \times \mathbb{R}_{+}\right)}^{2}.
$$
Moreover, another weighted norm $\|\cdot\|_{H_{g}^{s y}}$ introduced in \cite{mw}, which is almost equivalent to $\|\cdot\|_{H^{s}, \gamma},$ reads as
$$
\|\omega\|_{H_{g}^{s, y}}^{2}:=\left\|(1+y)^{\gamma} g_{s}\right\|_{L^{2}\left(\mathbb{T} \times \mathbb{R}_{+}\right)}^{2}+\sum_{|\alpha| \leq s, \alpha_{1} \leq s-1}\left\|(1+y)^{\gamma+\alpha_{2}} D^{\alpha} \omega\right\|_{L^{2}\left(\mathbb{T} \times \mathbb{R}_{+}\right)}^{2},
$$
where
$$
g_{s}:=\partial_{x}^{s} \omega-\frac{\partial_{y} \omega}{\omega} \partial_{x}^{s}(u-U), \quad u(t, x, y):=\int_{0}^{y} \omega(t, x, \tilde{y}) d \tilde{y}.
$$

\begin{Theorem}(\cite{S2}) \label{S2.1}
Let $s \geq 6$ be an even integer, $\gamma \geq 1, \sigma>\gamma+\frac{1}{2}$ and $0<\delta<\frac{1}{2}$. Suppose that the outer flow $U$ satisfies
$$
\sup _{t}\|U\|_{s+9, \infty}:=\sup _{t} \sum_{l=0}^{\left[\frac{s+9}{2}\right]}\left\|\partial_{t}^{l} U\right\|_{W^{s-2 l+, \infty}(\mathbb{T})}<+\infty .
$$
Assume that $u_{0}-\left.U\right|_{t=0} \in H^{s, \gamma-1}$ and the initial vorticity $\omega_{0}:=\partial_{y} u_{0} \in H_{\sigma, 2 \delta}^{s, \gamma} .$ Then, there exist a time $T:=T\left(s, \gamma, \sigma, \delta,\left\|\omega_{0}\right\|_{H_{g}^{s,}}, U\right)>0$ and a
unique classical solution $(u, v)$ to the Prandtl system such that
$$
u-U \in L^{\infty}\left([0, T] ; H^{s, \gamma-1}\right) \cap C\left([0, T] ; H^{s}-\omega\right),
$$
and the vorticity
$$
\omega:=\partial_{y} u \in L^{\infty}\left([0, T] ; H_{\sigma, \delta}^{s, \gamma}\right) \cap C\left([0, T] ; H^{s}-\omega\right),
$$
where $H^{s}-\omega$ is the space $H^{s}$ endowed with its weak topology.
Let $T^{*}$ be a maximal existence time of the classical solution such that $u-U \in H^{s, \gamma-1}$ and $\omega \in H_{\sigma, \delta}^{s, \gamma}$.
\end{Theorem}

\begin{Theorem}\cite{S2}
Under the conditions of Theorem \ref{S2.1}, and $\left[\sigma-\frac{\gamma+1}{2}\right] \geq 5 .$ If $T^{*}<\infty,$ then
$$
\lim _{t \rightarrow T^{*}}\left(\left\|\frac{1}{(1+y)^{\sigma} \partial_{y} u}\right\|_{L^{\infty}\left(\mathbb{T} \times \mathbb{R}_{+}\right)}+\|\nabla u\|_{L^{\infty}\left(\mathbb{T} \times \mathbb{R}_{+}\right)}+\left\|\frac{v_{x}}{1+y}\right\|_{L^{\infty}\left(\mathbb{T} \times \mathbb{R}_{+}\right)}\right)=\infty.
$$
Let us stress that $H_{\sigma, \delta}^{s, \gamma}$ is not a linear space, but a subset of linear function space $H^{s, \gamma} .$ That is, the function $\omega \in H_{\sigma, \delta}^{s, \gamma}$ satisfies not only the regularity condition
$$
\|\omega\|_{H^{s, \gamma}}<\infty,
$$
but also the monotonicity and decay conditions
$$
(1+y)^{\sigma} \omega \geq \delta, \quad \sum_{|\alpha| \leq 2}\left|(1+y)^{\sigma+\alpha_{2}} D^{\alpha} \omega\right|^{2} \leq \frac{1}{\delta^{2}}.
$$
Different from general blow up criterion, this blow-up criterion in Theorem 1.3.12 pays attention to both the regularity condition and the monotonicity and decay conditions. Combining this with the minimum principle, we arrive at a blow-up criterion in terms of the regularity condition of vorticity on the boundary.
\end{Theorem}

\begin{Theorem}(\cite{S2})
Under the conditions of Theorem \ref{S2.1}, and $\left[\sigma-\frac{\gamma+1}{2}\right] \geq 5 .$ If $T^{*}<\infty,$ then
$$
\lim _{t \rightarrow T^{*}}\left(\left\|\frac{1}{\left.\omega\right|_{y=0}}\right\|_{L^{\infty}(\mathbb{T})}+\|\nabla u\|_{L^{\infty}\left(\mathbb{T} \times \mathbb{R}^{+}\right)}+\left\|\frac{v_{x}}{1+y}\right\|_{L^{\infty}\left(\mathbb{T} \times \mathbb{R}^{+}\right)}\right)=\infty.
$$
\end{Theorem}

\section{3D Prandtl Equations}
In this section, we survey results on the 3D Prandtl equations, taking the form as
\begin{equation}\left\{
\begin{array}{ll}
\partial _{t}u+(u\partial _{x} +v\partial _{y} +w\partial_{z})u+\partial _{x}P=\partial _{z}^{2}u,\\
\partial _{t}v+(u\partial _{x} +v\partial _{y} +w\partial_{z})v+\partial _{y}P=\partial _{z}^{2}v,\\
\partial _{x}u+\partial _{y}v+\partial_{z}w=0,\\
(u,v)|_{t=0}=(u_{0}(x,y,z),v_{0}(x,y,z)),\\
(u,v,w)|_{z=0}=0,\\
  \lim\limits_{z\rightarrow+\infty}(u,v)=(U(t,x,y), V(t,x,y)),
\end{array}
 \label{w.100}         \right.\end{equation}
 where $(U(t, x,y), V(t,x,y)) $ satisfies
\begin{eqnarray}
\left\{\begin{array}{l}
{\partial_{t} U+U \partial_{x} U+V \partial_{y} U +\partial_{x} p=0    }, \\
{\partial_{t} V+U \partial_{x} V+V \partial_{y} V +\partial_{y} p=0    }.\end{array}\right.
\label{w.101}
\end{eqnarray}

\subsection{3D Prandtl Equations-Local Existence  }

  Oleinik and Samokhin (\cite{OS}) in 1999 considered
the system of 3D nonstationary boundary layer
\begin{equation}\left\{
\begin{array}{ll}
\frac{\partial u}{\partial t}+u\frac{\partial u}{\partial x}+v\frac{\partial u}{\partial y}+w\frac{\partial u}{\partial z}
 -\nu\frac{\partial^{2} u}{\partial z^{2}} =\frac{1}{\rho}\frac{\partial p}{\partial x},  \\
\frac{\partial v}{\partial t}+u\frac{\partial v}{\partial x}+v\frac{\partial v}{\partial y}+w\frac{\partial v}{\partial z}
-\nu\frac{\partial^{2} v}{\partial z^{2}} =\frac{1}{\rho}\frac{\partial p}{\partial y},\\
\frac{\partial u}{\partial x}+\frac{\partial v}{\partial y}+\frac{\partial w}{\partial z}=0,
\end{array}
 \label{ww10}         \right.\end{equation}
in a domain $G=\{0<t<T, 0<x,X,0<y<Y,0<z<Z\}$ with the conditions
\begin{equation}\left\{
\begin{array}{ll}
u(t,x,y,0)=0,\ \ v(t,x,y,0)=0, u( x,0) =0, \ \  w(t,x,y,0)=w_{0}(t,x,y)\leq 0,  \\
u=u_{0}, \ \ v=v_{0}, \ \  \text{on}\ \  \Gamma,
\end{array}
 \label{ww11}         \right.\end{equation}
where $\Gamma$ is the part of the boundary of $G$ belonging to the planes $t=0,x=0,y=0,z=Z$.

\begin{Theorem} \label{OS.2.13} (\cite{OS})
Assume that problem (\ref{ww10})-(\ref{ww11}) has a solution $(u,v,w)\in C^{2}(\overline{G})$. Then $u\geq 0, w\geq 0$ everywhere in $\overline{G}$;
and $u=0,w=0$ only for $z=0$.

\end{Theorem}

\begin{Theorem} \label{OS.2.14} (\cite{OS})
Let  $(u,v,w)\in C^{2}(\overline{G})$ be a solution of problem (\ref{ww10})-(\ref{ww11}). Then
\begin{eqnarray*}
&& 0\leq u\leq \max u_{0}+t\max (-\rho^{-1}p_{x}), \\
&& 0\leq v\leq \max v_{0}+t\max (-\rho^{-1}p_{y}).
\end{eqnarray*}

\end{Theorem}

\begin{Theorem} \label{OS.2.15} \cite{OS}
Let $(u,v,w)\in C^{2}(\overline{G})$ be a solution of problem (\ref{ww10})-(\ref{ww11}). Then the following inequalities hold in the domain $G$:
$$\min \left\{ \inf\limits_{z>0} \frac{v_{0}}{u_{0}},\min \frac{p_{y}}{p_{x}} \right\}\leq \frac{v}{u}
\leq \max \left\{ \sup\limits_{z>0} \frac{v_{0}}{u_{0}},\max \frac{p_{y}}{p_{x}} \right\} .$$
\end{Theorem}

 Under the assumption of Oleinik monotonicity, Liu, Wang and Yang (\cite{lwy}) in 2017 investigated the local well-posedness of solutions to the 3D Prandtl equations by Crocco transformation.
In \cite{lwy}, they studied the well-posedness for the following
problem of the Prandtl equations in the domain $Q_T=\{0<t\leq T,(x,y)\in D, z>0\}$,
\begin{eqnarray}\left\{\begin{array}{ll}
\partial_tu+(u\partial_xu+v\partial_yu+w\partial_zu)-\partial_x^2u=-\partial_xP,\\
\partial_tv+(u\partial_xv+v\partial_yv+w\partial_zv)-\partial_z^2v=-\partial_yP,\\
\partial_xu+\partial_xv+\partial_xw=0,\\
(u,w)|_{z=0}=0, \lim_{z\rightarrow\infty}=(U(t,x,y),KU(t,x,y)),\\
(u,v)|_{\partial Q_T^-}=(u_1(t,x,y,z),k(x,y)u_1(t,x,y,z)),\\
(u,v)|_{t=0}=(u_0(t,x,y,z),k(x,y)u_0(t,x,y,z))
\end{array}\right.\label{72.1}\end{eqnarray}
where $\partial Q_T^-=(0,T]\times\gamma_-\times \mathbb{R}_+$ with $\gamma_-=\{(x,y)\in\partial D|(1,k(x,y))\cdot \vec{n}(x,y)<0\}$.
The main results on the well-posedness of the initial boundary value problem \eqref{72.1} are given in the following theorem.
\begin{Theorem}( \cite{lwy})
Under the above assumption $(H)$ with $k\in C^{10}(D)$ and
$(U,p)\in C^{10}((0,T]\times D)$ for a fixed $T>0$, assume that
$$u_0\in C^{10}(\mathbb{R}^+_z\times D), u_1\in C^{15}(\partial Q_T^-)$$
have the following properties:
(1) $\partial_zu_0>0, u_1>0$ for all $z\geq 0$, and there is a constant $C_0 > 0$ such that
$$C_0^{-1}(U(0,x,y)-u_0(x,y,z))\leq\partial_zu_0(x,y,z)\leq C_0(U(0,x,y)-u_0(x,y,z)),$$

and

$$C_0^{-1}(U(t,x,y)-u_1(x,y,z))\leq\partial_zu_1(x,y,z)\leq C_0(U(t,x,y)-u_1(x,y,z)),\ on\  \partial Q_T^-.$$
(2)  the compatibility conditions of the problem \eqref{72.1} hold up to order 6 at $\{t=0\}\cap\partial Q_T^-$, and up to order 3 (4 resp.) at $\{t=0\}\cap\{z=0\}$($\{t=0\}\cap\{z=\infty\}$) and $\{z=0\}\cap\partial Q_T^-$($\{z=\infty\}\cap\partial Q_T^-$ resp).
Then, there exist $0 < T_0 \leq T$ and a unique classical solution $(u, v, w)$ to the problem
 \eqref{72.1} in the domain $Q_{T_0}$ . Moreover, the solution is linearly stable with respect to any three-dimensional smooth perturbation.
\end{Theorem}

Li and Yang (\cite{ly3}) in 2019 studied the 3D Prandtl equations without
any monotonicity condition on the velocity field, and proved that when one tangential
component of the velocity field has a single curve of non-degenerate critical points with
respect to the normal variable, the system is locally well-posed in the Gevrey function
space with Gevrey index.

Denote by $(u, v)$ the tangential component and by $w$ the vertical component of the velocity field, then the 3D Prandtl system in $\Omega$ reads as
\begin{equation}
\left\{\begin{array}{ll}
\partial_{t} u+\left(u \partial_{x}+v \partial_{y}+w \partial_{z}\right) u-\partial_{z}^{2} u+\partial_{x} p=0, \quad t>0, \quad(x, y, z) \in \Omega, \\
\partial_{t} v+\left(u \partial_{x}+v \partial_{y}+w \partial_{z}\right) v-\partial_{z}^{2} v+\partial_{y} p=0, \quad t>0,(x, y, z) \in \Omega, \\
\partial_{x} u+\partial_{y} v+\partial_{z} w=0, \quad t>0,(x, y, z) \in \Omega, \\
\left.u\right|_{z=0}=\left.v\right|_{z=0}=\left.w\right|_{z=0}=0, \quad \lim _{z \rightarrow+\infty}(u, v)=(U(t, x, y), V(t, x, y)) ,\\
\left.u\right|_{t=0}=u_{0},\left.\quad v\right|_{t=0}=v_{0}, \quad(x, y, z) \in \Omega,
\end{array}\right.\label{+14}
\end{equation}
where $(U(t, x, y), V(t, x, y))$ and $p(t, x, y)$ are the boundary traces of the tangential velocity field and pressure of the outer flow, satisfying Bernoulli's law
\begin{equation}
\left\{\begin{array}{ll}
\partial_{t} U+U \partial_{x} U+V \partial_{y} U+\partial_{x} p=0, \\
\partial_{t} V+U \partial_{x} V+V \partial_{y} V+\partial_{y} p=0.
\end{array}\right.\label{+15}
\end{equation}

Let $U, V$ be the data given in the equation. With each pair $(\rho, \sigma), \rho>0$ and $\sigma \geq 1$, a Banach space $X_{\rho, \sigma}$ consists of all smooth vector functions $(u, v)$ such that $\|(u, v)\|_{\rho, \sigma}<+\infty$, where the Gevrey norm $\|\cdot\|_{\rho, \sigma}$ is defined below. For each multi-index $\alpha=\left(\alpha_{1}, \alpha_{2}\right) \in \mathbb{Z}_{+}^{2},$ denote $\partial^{\alpha}=\partial_{x}^{\alpha_{1}} \partial_{y}^{\alpha_{2}}$ and
$$
\psi=\partial_{z}(u-U)=\partial_{z} u, \quad \eta=\partial_{z}(v-V)=\partial_{z} v.
$$
Then the Gevrey norm is defined by
$$
\begin{aligned}
\|(u, v)\|_{\rho, \sigma}=& \sup _{|\alpha| \geq 7} \frac{\rho^{|\alpha|-6}}{[(|\alpha|-7) !]^{\sigma}}\left(\left\|\langle z\rangle^{\ell-1} \partial^{\alpha}(u-U)\right\|_{L^{2}}+\left\|\langle z\rangle^{\kappa} \partial^{\alpha}(v-V)\right\|_{L^{2}}\right) \\
&+\sup _{|\alpha| \leq 6}\left(\left\|\langle z)^{\ell-1} \partial^{\alpha}(u-U)\right\|_{L^{2}}+\left\|\langle z\rangle^{\kappa} \partial^{\alpha}(v-V)\right\|_{L^{2}}\right) \\
&+\sup _{|\alpha| \geq 7} \frac{\rho^{|\alpha|-6}}{[(|\alpha|-7) !]^{\alpha}}\left\|\langle z\rangle^{\ell} \partial^{\alpha} \psi\right\|_{L^{2}}+\sup _{|\alpha| \leq 6}\left\|\langle z\rangle^{\ell} \partial^{\alpha} \psi\right\|_{L^{2}} \\
&+\sup _{|\alpha| \geq 7} \frac{\rho^{|\alpha|-5}}{[(|\alpha|-6) !]^{\sigma}}|\alpha|\left\|\langle z\rangle^{\kappa+2} \partial^{\alpha} \eta\right\|_{L^{2}}+\sup _{|\alpha| \leq 6}\left\|\langle z\rangle^{\kappa+2} \partial^{\alpha} \eta\right\|_{L^{2}} \\
&+\sup _{1 \leq j \leq 4 \atop|\alpha|+j \geq 7} \frac{\rho^{|\alpha|+j-6}}{[(|\alpha|+j-7) !]^{\sigma}}\left\|\langle z\rangle^{\ell+1} \partial^{\alpha} \partial_{z}^{j} \psi\right\|_{L^{2}}+\sup _{1<j \leq 4 \atop|\alpha|+j \leq 6}\left\|\langle z\rangle^{\ell+1} \partial^{\alpha} \partial_{z}^{j} \psi\right\|_{L^{2}} \\
+& \sup _{1 \leq j \leq 4 \atop|\alpha|+j \geq 7} \frac{\rho^{|\alpha|+j-5}}{[(|\alpha|+j-6) !]^{\sigma}}|\alpha|\left\|\langle z\rangle^{\kappa+2} \partial^{\alpha} \partial_{z}^{j} \eta\right\|_{L^{2}}+\sup _{1 \leq \leq \leq \atop|\alpha|+j \leq 6}\left\|\langle z\rangle^{\kappa+2} \partial^{\alpha} \partial_{z}^{j} \eta\right\|_{L^{2}}, \\
\end{aligned}
$$
where we use $L^{2}$ instead of $L^{2}(\Omega)$ without confusion. Moreover, define another Gevrey space $Y_{\rho, \sigma}$ consisting of smooth functions $F(x, y)$ such that $\|F\|_{\rho, \sigma}<+\infty$, where
$$
\|F\|_{\rho, \sigma}=\sup _{|\alpha| \geq 0} \frac{\rho^{|\alpha|}}{(|\alpha| !)^{\sigma}}\left\|\partial^{\alpha} F\right\|_{L^{2}\left(\mathbb{T}^{2}\right)}.
$$

The initial data $\left(u_{0}, v_{0}\right)$ satisfy the following compatibility conditions
\begin{equation}
\left\{\begin{array}{l}
\left.\left(u_{0}, v_{0}\right)\right|_{z=0}=(0,0), \lim \limits_{z \rightarrow+\infty}\left(u_{0}, v_{0}\right) \mid=(U, V) \text { and }\left.\left(\partial_{z} \psi_{0}, \partial_{z} \eta_{0}\right)\right|_{z=0}=\left(\partial_{x} p, \partial_{y} p\right) ,\\
\left.\partial_{z}^{3} \psi_{0}\right|_{z=0}=\left.\psi_{0}\left(\partial_{x} \psi_{0}-\partial_{y} \eta_{0}\right)\right|_{z=0}+\left.2 \eta_{0} \partial_{y} \psi_{0}\right|_{z=0}+\partial_{t} \partial_{x} p, \\
\left.\partial_{z}^{3} \eta_{0}\right|_{z=0}=\left.\eta_{0}\left(\partial_{y} \eta_{0}-\partial_{x} \psi_{0}\right)\right|_{z=0}+\left.2 \psi_{0} \partial_{x} \eta_{0}\right|_{z=0}+\partial_{t} \partial_{y} p,
\end{array}\right.\label{+16}
\end{equation}
where $\psi_{0}=\partial_{z} u_{0}$ and $\eta_{0}=\partial_{z} v_{0}$.

\begin{Theorem}(\cite{ly3})
Let $1<\sigma \leq 2,$ under the compatibility conditions,  suppose $U, V, p \in$ $Y_{2 \rho_{0}, \sigma}$ and $\left(u_{0}, v_{0}\right) \in X_{2 \rho_{0}, \sigma}$ for some $\rho_{0}>0 .$ Moreover, suppose $u_{0}$ satisfies the assumptions \eqref{+16}. Then the Prandtl system admits a unique solution $(u, v) \in L^{\infty}\left([0, T] ; X_{\rho, \sigma}\right)$ for some $T>0$ and some $0<\rho<2 \rho_{0}$.
\end{Theorem}

Li,  Masmoudi and   Yang  (\cite{LMD}) in 2020 established the well-posedness in Gevrey function space with
optimal class of regularity 2 for the 3D Prandtl system without any structural assumption.

Denote $\mathbb{R}_{+}^{3}=\left\{(x, y, z) \in \mathbb{R}^{3} ; z>0\right\}$ and let $(u, v)$ be the tangential component and $w$ be the vertical component of the velocity field. Then the 3D Prandtl system in $\mathbb{R}_{+}^{3}$ reads as
\begin{equation}
\left\{\begin{array}{l}
\left(\partial_{t}+u \partial_{x}+v \partial_{y}+w \partial_{z}-\partial_{z}^{2}\right) u+\partial_{x} p=0, \quad t>0,(x, y, z) \in \mathbb{R}_{+}^{3}, \\
\left(\partial_{t}+u \partial_{x}+v \partial_{y}+w \partial_{z}-\partial_{z}^{2}\right) v+\partial_{y} p=0, \quad t>0,(x, y, z) \in \mathbb{R}_{+}^{3}, \\
\partial_{x} u+\partial_{y} v+\partial_{z} w=0, \quad t>0,(x, y, z) \in \mathbb{R}_{+}^{3} ,\\
\left.u\right|_{z=0}=\left.v\right|_{z=0}=\left.w\right|_{z=0}=0, \quad \lim \limits_{z \rightarrow+\infty}(u, v)=(U(t, x, y), V(t, x, y)), \\
\left.u\right|_{t=0}=u_{0},\left.\quad v\right|_{t=0}=v_{0}, \quad(x, y, z) \in \mathbb{R}_{+}^{3},
\end{array}\right.
\end{equation}
where $(U(t, x, y), V(t, x, y))$ and $p(t, x, y)$ are the boundary traces of the tangential velocity field and pressure of the outer flow, satisfying
\begin{equation}
\left\{\begin{array}{l}
\partial_{t} U+U \partial_{x} U+V \partial_{y} U+\partial_{x} p=0 ,\\
\partial_{t} V+U \partial_{x} V+V \partial_{y} V+\partial_{y} p=0,
\end{array}\right.\label{+17}
\end{equation}
Here, $p, U, V$ are given functions determined by the Euler flow.
Then for the zero outer flow, the Prandtl system can be written as
$$
\left\{\begin{array}{l}
\left(\partial_{t}+u \partial_{x}+v \partial_{y}+w \partial_{z}-\partial_{z}^{2}\right) u=0, \quad t>0,(x, y, z) \in \mathbb{R}_{+}^{3}, \\
\left(\partial_{t}+u \partial_{x}+v \partial_{y}+w \partial_{z}-\partial_{z}^{2}\right) v=0, \quad t>0,(x, y, z) \in \mathbb{R}_{+}^{3}, \\
\left.(u, v)\right|_{z=0}=(0,0), \quad \lim\limits _{z \rightarrow+\infty}(u, v)=(0,0), \\
\left.(u, v)\right|_{t=0}=\left(u_{0}, v_{0}\right), \quad(x, y, z) \in \mathbb{R}_{+}^{3},
\end{array}\right.
$$
with
$$
w(t, x, y, z)=-\int_{0}^{z} \partial_{x} u(t, x, y, \tilde{z}) d \tilde{z}-\int_{0}^{z} \partial_{y} v(t, x, y, \tilde{z}) d \tilde{z}.
$$

We will use without confusion $\|\cdot\|_{L^{2}}$ and $(\cdot, \cdot)_{L^{2}}$ to denote the norm and inner product of $L^{2}=L^{2}\left(\mathbb{R}_{+}^{3}\right),$ and use the notations $\|\cdot\|_{L^{2}\left(\mathbb{R}_{x, y}^{2}\right)}$ and $(\cdot, \cdot)_{L^{2}\left(\mathbb{R}_{x, y}^{2}\right)}$ when the variables are specified. Similarly for $L^{\infty}$. Moreover,  we also use $L_{x, y}^{\infty}\left(L_{z}^{2}\right)=L^{\infty}\left(\mathbb{R}^{2} ; L^{2}\left(\mathbb{R}_{+}\right)\right)$ to stands for the
classical Sobolev space, so does the Sobolev space $L_{x, y}^{2}\left(L_{z}^{\infty}\right) .$ In the following discussion by $\partial^{\alpha}$, we always mean $\partial^{\alpha}=\partial_{x}^{\alpha_{1}} \partial_{y}^{\alpha_{2}}$ with each multi-index $\alpha=$ $\left(\alpha_{1}, \alpha_{2}\right) \in \mathbb{Z}_{+}^{2}$.

Let $\ell>1 / 2$ be a given number. With each pair $(\rho, \sigma), \rho>0$ and $\sigma \geq 1$, a Banach space $X_{\rho, \sigma}$ consists of all smooth vector-valued functions $(u, v)$ such that the Gevrey norm $\|(u, v)\|_{\rho, \sigma}<+\infty,$ where $\|\cdot\|_{\rho, \sigma}$ is defined below. Recalling $\partial^{\alpha}=\partial_{x}^{\alpha_{1}} \partial_{y}^{\alpha_{2}}$, we define
$$
\begin{aligned}
\|(u, v)\|_{\rho, \sigma}=& \sup _{0 \leq j \leq 5 \atop|\alpha|+j \geq 7} \frac{\rho^{|\alpha|+j-7}}{[(|\alpha|+j-7) !]^{\sigma}}\left(\left\|\langle z\rangle^{\ell+j} \partial^{\alpha} \partial_{z}^{j} u\right\|_{L^{2}}+\left\|\langle z\rangle^{\ell+j} \partial^{\alpha} \partial_{z}^{j} v\right\|_{L^{2}}\right) \\
&+\sup _{0 \leq j \leq 5 \atop|\alpha|+j \leq 6}\left(\left\|\langle z\rangle^{\ell+j} \partial^{\alpha} \partial_{z}^{j} u\right\|_{L^{2}}+\left\|\langle z\rangle^{\ell+j} \partial^{\alpha} \partial_{z}^{j} v\right\|_{L^{2}}\right),
\end{aligned}
$$
where $\langle z\rangle=\left(1+|z|^{2}\right)^{1 / 2}$. We call $\sigma$ the Gevrey index.

The initial data $\left(u_{0}, v_{0}\right)$ satisfy the following compatibility conditions
$$
\left\{\begin{array}{l}
\left.\left(u_{0}, v_{0}\right)\right|_{z=0}=(0,0),\ \   \lim \limits_{z \rightarrow+\infty}\left(u_{0}, v_{0}\right)\left|=(0,0),\left(\partial_{z}^{2} u_{0}, \partial_{z}^{2} v_{0}\right)\right|_{z=0}=(0,0), \\
\left.\partial_{z}^{4} u_{0}\right|_{z=0}=\left.\left(\partial_{z} u_{0}\right)\left(\partial_{x} \partial_{z} u_{0}-\partial_{y} \partial_{z} v_{0}\right)\right|_{z=0}+\left.2\left(\partial_{z} v_{0}\right) \partial_{y} \partial_{z} u_{0}\right|_{z=0}, \\
\left.\partial_{z}^{4} v_{0}\right|_{z=0}=\left.\left(\partial_{z} v_{0}\right)\left(\partial_{y} \partial_{z} v_{0}-\partial_{x} \partial_{z} u_{0}\right)\right|_{z=0}+\left.2\left(\partial_{z} u_{0}\right) \partial_{x} \partial_{z} u_{0}\right|_{z=0}.
\end{array}\right.
$$

\begin{Theorem}(\cite{LMD})
Let $1<\sigma \leq 2$. Suppose the initial datum $\left(u_{0}, v_{0}\right)$ belongs to $X_{2 \rho_{0}, \sigma}$ for some $\rho_{0}>0,$ and satisfies the compatibility condition.
Then the  system \eqref{+17} admits a unique solution $(u, v) \in L^{\infty}\left([0, T] ; X_{\rho, \sigma}\right)$ for some $T>0$ and some $0<\rho<2 \rho_{0}$.
\end{Theorem}

\subsection{3D Prandtl Equations-Global Existence}
In this subsection, we survey the global existence of solutions to the 3D Prandtl equations.
 Liu, Wang and Yang (\cite{lwy1}) in 2016 proved the global existence of weak solutions to the  Prandtl equations (\ref{w.100})  under special structure  takes the form:
    $$(u(t,x,y,z), k(x,y)u(t,x,y,z), w(t,x,y,z)),$$
in the domain $Q=\{t>0, (x,y)\in D , z>0\},~D \subseteq \mathbb{R}^{2}$,  by the Crocco transformation,
$$\xi=x,\ \ \eta=y,\ \ \zeta=\frac{u(t,x,y,z)}{U(t,x,y)}, W(t,\xi, \eta,\zeta)=\frac{\partial_{z}u(t,x,y,z)}{U(t,x,y)},$$
the problem   (\ref{w.100})  becomes
\begin{eqnarray}
\left\{\begin{array}{l}
{L(W)\triangleq \partial_{t} W+\zeta U( \partial_{\xi}+k  \partial_{\eta})W+ A\partial_{\zeta}W+BW-W^{2}\partial_{\zeta}^{2}W=0    }, \\
{ W|_{\zeta=1}=0,\ \ W\partial_{\zeta}W|_{\zeta=0}=\frac{p_{x}}{U} },\\
{W|_{\Gamma_{-}}=W_{1}(t,\xi,\eta,\zeta)\triangleq  \frac{\partial_{z}u_{1}}{U}    },\\
{W|_{t=0}=W_{0}( \xi,\eta,\zeta)\triangleq  \frac{\partial_{z}u_{0}}{U}    },
\end{array}\right.\label{w.56}
\end{eqnarray}
where
\begin{eqnarray*}
\left\{\begin{array}{l}
{\Omega_{T}=\{(t,\xi,\eta,\zeta)| 0<t<T, (\xi,\eta)\in \mathbb{R}^{2}, 0\leq \zeta<1\}   }, \\
{  \Gamma_{-}= \{(t,\xi,\eta,\zeta)| 0<t<T, (\xi,\eta)\in \gamma_{-}, 0\leq \zeta<1\}     },
\end{array}\right.
\end{eqnarray*}
and
$$A=-\zeta(1-\zeta)\frac{U_{t}}{U}-(1-\zeta^{2})\frac{p_{x}}{U} ,\ \ B=\frac{U_{t}}{U}+\zeta(U_{x}+kU_{y})-\partial_{y}k\cdot\zeta U. $$
Assume that
\begin{eqnarray*}
\partial_{z}u>0, \ \ \partial_{z}u_{0}>0, \ \ \text{for}~z>0,
\end{eqnarray*}
and the favourable pressure condition holds:
\begin{eqnarray}
p_{\xi}(t,\xi,\eta)\leq 0,\ \  \forall t>0, (\xi,\eta)\in \mathbb{R}^{2}.
\label{w.57}
\end{eqnarray}
\begin{Theorem}(\cite{lwy1})
For the problem (\ref{w.56}), and any $T>0$, assume that $ k\in C^{2}(D), U\in C^{2}((0,T)\times D), p_{x}\in C^{1}((0,T)\times D)$ satisfies (\ref{w.57}), and the initial boundary data $W_{0}\in C^{1}(\Omega), W_{1}\in C^{3}(\Gamma_{-})$ satisfy
\begin{eqnarray*}
C_{0}^{-1}(1-\zeta)\leq W_{0}, W_{1}\leq C_{0}(1-\zeta),
\end{eqnarray*}
for a positive constant $C_{0}$. Then, there exists a weak solution $W(t,\xi,\eta,\zeta)\in L^{\infty}(0, T; BV(\Omega))$ to the problem (\ref{w.56}) in the sense of  weak solution.
\end{Theorem}

 Zhang,  Zhang (\cite{zz}) in 2016 investigated the long time existence and uniqueness of small solutions for
the Prandtl system with small initial data by shear flow with  $\mathbb{R}_{+} \times \mathbb{R}_{+}^{d},$ $d=2,3$,
\begin{equation}\label{w.58}
\left\{\begin{array}{l}{\partial_{t} u+u \cdot \nabla_{\mathrm{h}} u+v \partial_{y} u-\partial_{y y} u+\nabla_{h} p=0, \quad(t, x, y) \in \mathbb{R}_{+} \times \mathbb{R}^{d-1} \times \mathbb{R}_{+}}, \\ {\operatorname{div}_{\mathrm{h}} u+\partial_{y} v=0} ,\\ {\left.u\right|_{y=0}=0,\left.\quad v\right|_{y=0}=0, \quad \text { and } \quad \lim\limits _{y \rightarrow \infty} u(t, x, y)=U(t, x)}, \\ {\left.u\right|_{t=0}=u_{0},}\end{array}\right.
\end{equation}
where $u=u, \nabla_{\mathrm{h}}=\operatorname{div}_{\mathrm{h}}=\partial_{x}$ for $d=2$, and  $u=\left(u^{1}, u^{2}\right), \nabla_{\mathrm{h}}=\left(\partial_{x_{1}}, \partial_{x_{2}}\right), \operatorname{div}_{\mathrm{h}}=\partial_{x_{1}}+\partial_{x_{2}}$ for $d=3$. For simplicity, let $U=\varepsilon \mathbf{e}$ for some unit vector $\mathbf{e} \in \mathbb{R}^{d-1}$.

Let $u^{s}(t, y)$ be determined by
\begin{equation}\label{w.59}
\left\{\begin{array}{ll}{\partial_{t} u^{s}-\partial_{y y} u^{s}=0, (t, y) \in \mathbb{R}_{+} \times \mathbb{R}_{+}} ,\\
{ u^{s}|_{y=0}=0, \quad \lim \limits_{y \rightarrow \infty} u^{s}(t, y)=\varepsilon \mathbf{e}} ,\\ {u^{s}|_{t=0}=\varepsilon \chi(y) \mathbf{e}},\end{array}\right.
\end{equation}
where $\chi(y) \in C^{\infty}(\mathbb{R})$, and $\chi(y)=0$ for $y\leq1$ and $\chi(y)=1$ for $y\leq2$, set $u=u^{s}+w$. Then $w$ satisfies
\begin{equation}\label{w.60}
\left\{\begin{array}{l}{\partial_{t} w+\left(w+u^{s}\right) \cdot \nabla_{\mathrm{h}} w-\int_{0}^{y} \operatorname{div}_{\mathrm{h}} w d y^{\prime} \partial_{y} w-\int_{0}^{y} \operatorname{div}_{\mathrm{h}} w d y^{\prime} \partial_{y} u^{s}-\partial_{y y} w=0}, \\ {\left.w\right|_{y=0}=0, \quad   \quad \lim\limits_{y \rightarrow \infty} w=0}, \\ {\left.w\right|_{t=0}=u_{0}-\varepsilon \chi \mathbf{e} \stackrel{\text { def }}{=} w_{0}}. \end{array}\right.
\end{equation}
 Zhang and  Zhang (\cite{zz}) introduced a key quantity  $\theta(t)$ to describe the evolution of the analytic band of $w$:
\begin{equation}\label{w.61}
\left\{\begin{array}{l}{\dot{\theta}(t)=\langle t\rangle^{\frac{1}{4}}\left(\left\|e^{\Psi} \partial_{y} w_{\Phi}(t)\right\|_{\mathcal{B}^{\frac{d-1}{2}, 0}}+\left\|e^{\Psi} \partial_{y} u^{s}(t)\right\|_{L_{v}^{2}}\right)}, \\ {\left.\theta\right|_{t=0}=0}.\end{array}\right.
\end{equation}
Here $\langle t\rangle= 1+t$, the phase function $\Phi$ is defined by
\begin{equation}\label{w.62}
\Phi(t, \xi) \stackrel{\text { def }}{=}(\delta-\lambda \theta(t))|\xi|,
\end{equation}
and the weighted function $\Psi(t, y)$ is determined by
\begin{equation}\label{w.63}
\Psi(t, y) \stackrel{\text { def }}{=} \frac{1+y^{2}}{8\langle t\rangle}.
\end{equation}
\begin{Theorem}(\cite{zz})
Let $\delta>0$ and $\varepsilon$ be a sufficiently small positive constant. Assume that  $w_{0}$ satisfies
$$
\left\|e^{\frac{1+y^{2}}{8\langle t\rangle }} e^{\delta|D_{x}|} w_{0}\right\|_{\mathcal{B}^{\frac{d-1}{2}}, 0} \leq \varepsilon,
$$
then there exits a positive time $T_{\varepsilon}$ which is of size greater than $\varepsilon^{-\frac{4}{3}}$ so that the problem (\ref{w.60}) has a
unique solution w which satisfies
\begin{eqnarray*}
 \|e^{\Psi }w_{\Phi }\|_{ \widetilde{L}_{t}^{\infty}\left(\mathcal{B}^{\frac{d-1}{2}, 0}\right)}
+c\sqrt{\lambda} \|e^{\Psi }w_{\Phi }\|_{ \widetilde{L}_{t, \dot{\theta}(t)}^{2}\left(\mathcal{B}^{\frac{d }{2}, 0}\right)}
+\|e^{\Psi }\partial_{y}w_{\Phi }\|_{ \widetilde{L}_{t}^{2}\left(\mathcal{B}^{\frac{d-1}{2}, 0}\right)}
 \leq \|e^{\frac{1+y^{2}}{8}} e^{\delta |D_{x}|}w_{0}\|_{\mathcal{B}^{\frac{d-1}{2}, 0}}   \ \   \text{for any}\ \ t\in[0,T_{\varepsilon}].
\end{eqnarray*}
In particular, we have
\begin{eqnarray*}
e^{\Psi(t, y)} e^{\Phi(t, D)} w \in \widetilde{L}_{T}^{\infty}\left(\mathcal{B}^{\frac{d-1}{2}, 0}\right),\left.\quad e^{\Psi(t, y)} e^{\Phi(t, D)}\right|_{\partial y} w \in \widetilde{L}_{T}^{2}\left(\mathcal{B}^{\frac{d-1}{2}, 0}\right),
\end{eqnarray*}
for any $T\leq T_{\varepsilon}$, and where the functions $\Psi(t, y), \ \Phi(t, \xi)$ are determined by (\ref{w.61})-(\ref{w.63})  respectively.
Here define $  \widetilde{L}_{T}^{p}\left(\mathcal{B}^{s, 0}\right)$ as the completion of $C([0,T];\mathcal{S}(\mathbb{R}^{d}_{+}))$ by the norm
\begin{eqnarray*}
&&\|a\|_{ \widetilde{L}_{T}^{p}\left(\mathcal{B}^{s, 0}\right) }:=\sum\limits_{k\in\mathbb{Z}}2^{ks}\left(
\int_{0}^{T}\|\triangle ^{h}_{k}a(t)\|_{L^{2}_{+}}dt  \right)^{\frac{1}{p}},\\
&&\|a\|_{ \widetilde{L}_{t,f}^{p}\left(\mathcal{B}^{s, 0}\right) }:=\sum\limits_{k\in\mathbb{Z}}2^{ks}\left(
\int_{0}^{t} f(t') \|\triangle ^{h}_{k}a(t)\|_{L^{2}_{+}}dt  \right)^{\frac{1}{p}},
\end{eqnarray*}
with the usual change if $p=\infty$.
Here
$$e^{\delta|D_{x}|} w_{0}=\mathcal{F}^{-1}_{\xi\rightarrow x}e^{ (\delta-\lambda \theta(t))|\xi|} \hat{w}_{0}(t,\xi,y) .$$

\end{Theorem}

Luo and Xin (\cite{LX}) in 2018  obtained H$\ddot{o}$lder continuous weak solutions to the 3D Prandtl equation  with tangential velocities transverse to the outflow $U$.
Consider the following 3D Prandtl system
\begin{eqnarray}\left\{\begin{array}{ll}
u_t+u\partial_xu+v\partial_yu-\partial_{y}^2u+\nabla_xp=0,\\
\nabla_xu+\partial_yv=0,\\
(u,v)|_{y=0}=0,\quad  \lim\limits_{y\rightarrow+\infty}u=U,\\
u|_{t=0}=u_0.
\end{array}\right.\label{82.1}\end{eqnarray}
Here $x=(x_1,x_2)\in \mathbb{T}$ and $y\in \mathbb{R}$ denote the tangential and the vertical components of the space variable, respectively. $\nabla_x=(\partial_{x_1},\partial_{x_2})$ denotes the tangential gradient; $u=(u^1,u^2)(t,x,y)$ and $v=v(t,x,y)$  denote the tangential and the vertical velocities;
$U(t,x)$ and $p(t,x)$  denote the tangential velocity and the pressure on the boundary $\{y=+\infty\}$ of the outer Euler flow, respectively, which satisfy
$$\partial_tU+U\cdot\nabla_x U+\nabla_x p=0.$$
\begin{Theorem}(\cite{LX}) \label{LX1}
Suppose that $(u_C, v_C )$ is a classical solution to the system \eqref{82.1} and
$(u,v)$ is a smooth perturbation of $(u_C, v_C)$ such that the difference $(u-u_C, v-v_C)\in C_c^\infty(\mathbb{R}_+\times\mathbb{T}^2\times\mathbb{R}_+)$ has compact support and satisfies
$$\nabla_x {u}+\partial_y {v}=0,\ in \ \mathbb{R}_+\times\mathbb{T}^2\times\mathbb{R}_+,$$
$$\int_{\mathbb{T}^2}(u-u_C)(t,x,y)dx=0,\ for\ (t,y)\in \mathbb{R}_+\times\mathbb{R}_+.$$
Let $\rho>0$  be a given positive number such that
$$\overline{N(supp_{t,y}(u-u_C, {v}-v_C);\rho,\rho^{1/2})}\subset\mathbb{R}_+\times\mathbb{R}_+,$$
where $supp_{t,y}$ denotes the projection of the support to $\mathbb{R}_+\times\mathbb{R}_+$.  Then there exists a sequence of H$\ddot{o}$lder continuous weak solutions $\{(u_k,v_k)\}_{=1}^\infty$ to the system \eqref{82.1}
and a sequence of positive numbers $\{C_k\}$  satisfying the estimates
$$[u_k-{u},v_k-{v}]_{\frac{1}{21}-\varepsilon,\frac{1}{10}-\varepsilon}\leq C_k,$$
and
$$supp_{t,y}(u_k-{u},v_k-{v})\subset\overline{N(supp_{t,y}(u-u_C, {v}-v_C);\rho,\rho^{1/2})}.$$
Furthermore, $(u_k,v_k)\rightharpoonup(u,v)$ in the weak- $\ast$  topology on $L^\infty(\mathbb{R}_+\times\mathbb{T}^2\times\mathbb{R}_+).$
\end{Theorem}
Here denote
\begin{eqnarray*}
&&N(E; \rho,\rho')=\{(t,y):|t-t_{0}|<\rho, |y-y_{0}|<\rho',\text{for ~ some }~(t_{0},y_{0})\in E\},\\
&& [f]_{\alpha,\beta(\Omega)}:=\sup\limits_{(t,x,y),(t',x',y')\in \Omega}\frac{|f(t,x,y)-f(t',x'y')|}
{|t-t'|^{\alpha}+ |x-x'|^{\alpha}+|y-y'|^{\beta} }.
\end{eqnarray*}
The constructions in the proof of   Theorem  $\ref{LX1}$  can also be adapted to some other
models with vertical viscosities. In particular, consider the system
\begin{eqnarray}\left\{\begin{array}{ll}
u_t+u\cdot\nabla u-\partial_{y}^2u+\nabla_xp=0,\\
\nabla_xu+\partial_yv=0,\\
u|_{t=0}=u_0.
\end{array}\right.\label{82.2}\end{eqnarray}
\begin{Theorem}(\cite{LX})
 There exists a non-trivial H$\ddot{o}$lder continuous weak solutions $u$ to the
system \eqref{82.2} which is supported in a compact time interval, with
$$[u]_{\frac{1}{21}-\varepsilon,\frac{1}{10}-\varepsilon}<+\infty.$$
\end{Theorem}

Lin and Zhang (\cite{LZ1}) in 2020  proved the almost global existence of classical solutions to the 3D Prandtl system with the initial data which lie within $\varepsilon$ of a stable shear flow. Using anisotropic Littlewood-Paley energy estimates in tangentially analytic norms and introducing new linearly-good unknowns, they proved that the 3D Prandtl system has a unique solution with the lifespan of which is greater than $\exp \left(\varepsilon^{-1} / \log \left(\varepsilon^{-1}\right)\right) .$
They considered the following Prandtl boundary layer equations in $\mathbb{R}_{+} \times \mathbb{R}_{+}^{3}$ :
\begin{equation}
\left\{\begin{array}{l}
\partial_{t} u^{p}+\left(u^{p} \partial_{x}+v^{p} \partial_{y}+w^{p} \partial_{z}\right) u^{p}+\partial_{x} p^{E}=\partial_{z}^{2} u^{p}, \\
\partial_{t} v^{p}+\left(u^{p} \partial_{x}+v^{p} \partial_{y}+w^{p} \partial_{z}\right) v^{p}+\partial_{y} p^{E}=\partial_{z}^{2} v^{p}, \\
\partial_{x} u^{p}+\partial_{y} v^{p}+\partial_{z} w^{p}=0, \\
\left.\left(u^{p}, v^{p}, w^{p}\right)\right|_{z=0}=(0,0,0), \\
\lim\limits _{z \rightarrow+\infty}\left(u^{p}, v^{p}\right)=\left(U^{E}(t, x, y), V^{E}(t, x, y)\right), \\
\left.\left(u^{p}, v^{p}\right)\right|_{t=0}=\left(u_{0}(x, y, z), v_{0}(x, y, z)\right).
\end{array}\right.\label{LZ.1}
\end{equation}
Here and in what follows, $(t, x, y, z) \in \mathbb{R}_{+} \times \mathbb{R} \times \mathbb{R} \times \mathbb{R}_{+},\left(u^{p}, v^{p}\right)$ and $w^{p}$ denote the tangential and normal velocity of the boundary layer flow, the initial data $\left(u_{0}, v_{0}\right):=$ $\left(u_{0}(x, y, z), v_{0}(x, y, z)\right)$ and the far-field $\left(U^{E}(t, x, y), V^{E}(t, x, y)\right)$ are given. Furthermore, $\left(U^{E}(t, x, y), V^{E}(t, x, y)\right)$ and the given scalar pressure $p^{E}(t, x, y)$ are the tangential velocity field and pressure on the boundary $\{z=0\}$ of the Euler flow, satisfying
\begin{equation}
\left\{\begin{array}{l}
\partial_{t} U^{E}+U^{E} \partial_{x} U^{E}+V^{E} \partial_{y} U^{E}+\partial_{x} p^{E}=0, \\
\partial_{t} V^{E}+U^{E} \partial_{x} V^{E}+V^{E} \partial_{y} V^{E}+\partial_{y} p^{E}=0. \quad t>0,(x, y) \in \mathbb{R}^{2}.
\end{array}\right.\label{LZ.2}
\end{equation}

Set
\begin{equation}
\lim _{z \rightarrow+\infty}\left(u^{p}, v^{p}\right)=\left(\kappa_{1}, \kappa_{2}\right)  .\label{LZ.3}
\end{equation}
Write $u^{p}(t, x, y, z)$ and $v^{p}(t, x, y, z)$ as perturbations $u(t, x, y, z)$ and $v(t, x, y, z)$ of the lifts $\kappa_{1} \varphi(t, z)$ and $\kappa_{2} \varphi(t, z)$ respectively via
$$
u^{p}(t, x, y, z)=\kappa_{1} \varphi(t, z)+u(t, x, y, z), \quad v^{p}(t, x, y, z)=\kappa_{2} \varphi(t, z)+v(t, x, y, z)
$$
where
$$
\varphi(t, z)=\frac{1}{\sqrt{\pi}} \int_{0}^{z / \sqrt{\langle t\rangle}} \exp \left(-\frac{\widetilde{z}^{2}}{4}\right) d \widetilde{z}
$$
and denote by $\langle t\rangle=1+t$. Then introduce new linearly-good unknowns
\begin{equation*}
\left\{\begin{array}{l}
g_{1}(t, x, y, z)=\partial_{z} u(t, x, y, z)+\frac{z}{2\langle t\rangle} u(t, x, y, z), \\
g_{2}(t, x, y, z)=\partial_{z} v(t, x, y, z)+\frac{z}{2\langle t\rangle} v(t, x, y, z).
\end{array}\right.
\end{equation*}
Through some simple calculations,  the equations for the good unknowns $g=$ $\left(g_{1}, g_{2}\right)^{T}$ in $\left\{t>0,(x, y) \in \mathbb{R}^{2}, z \in \mathbb{R}_{+}\right\}$ are
\begin{equation}
\left\{\begin{array}{l}
\partial_{t} g-\partial_{z}^{2} g+\frac{1}{\langle t\rangle} g+\kappa_{1} \partial_{z} \varphi \nabla_{h}^{\perp} v+\kappa_{1} \varphi \partial_{x} g+g_{1} \nabla_{h}^{\perp} v-\frac{z}{2\langle t\rangle} u \nabla_{h}^{\perp} v+u \partial_{x} g, \\
-\kappa_{2} \partial_{z} \varphi \nabla_{h}^{\perp} u+\kappa_{2} \varphi \partial_{y} g-g_{2} \nabla_{h}^{\perp} u+\frac{z}{2\langle t\rangle} v \nabla_{h}^{\perp} u+v \partial_{y} g+w \partial_{z} g-\frac{1}{2\langle t\rangle} w \mathbf{u}=0, \\
u=U\left(g_{1}\right), \quad v=U\left(g_{2}\right), \quad \text { and } \quad w=W\left(g_{1}, g_{2}\right), \\
\left.\partial_{z} g\right|_{z=0}=\lim \limits_{z \rightarrow+\infty} g=0, \\
\left.g\right|_{t=0}=\left(g_{10}, g_{20}\right),
\end{array}\right.
\end{equation}
where
\begin{equation*}
\left\{\begin{array}{l}
\mathbf{u}=(u, v)^{T}, \quad \nabla_{h}=\left(\begin{array}{c}
\partial_{x} \\
\partial_{y}
\end{array}\right), \quad \nabla_{h}^{\perp}=\left(\begin{array}{c}
-\partial_{y}, \\
\partial_{x}
\end{array}\right), \\
U\left(g_{i}\right) \doteq \exp \left(\frac{-z^{2}}{4\langle t\rangle}\right) \int_{0}^{z} g_{i}(t, x, y, \widetilde{z}) \exp \left(\frac{\widetilde{z}^{2}}{4\langle t\rangle}\right) d \widetilde{z}, \quad i=1,2, \\
W\left(g_{1}, g_{2}\right) \doteq-\int_{0}^{z}\left[U\left(\partial_{x} g_{1}\right)+U (\partial_{y} g_{2} )\right] d \widetilde{z}.
\end{array}\right.
\end{equation*}

For any local bounded function $\Phi$ on $\mathbb{R}_{+} \times \mathbb{R}^{2}$,   define
$$
\begin{array}{l}
u_{\Phi}(t, x, y, z)=\mathscr{F}_{\xi \rightarrow(x, y)}^{-1}\left(e^{\Phi(t, \xi)} \widehat{u}(t, \xi, z)\right), \\
\widetilde{u_{\Phi}}(t, x, y, z)=\mathscr{F}_{\xi \rightarrow(x, y)}^{-1}\left(\left|e^{\Phi(t, \xi)} \widehat{u}(t, \xi, z)\right|\right),
\end{array}
$$
where and in all that follows, $\mathscr{F} u$ and $\widehat{u}$ always denote the partial Fourier transform of the distribution $u$ with respect to $(x, y)$ variables, $\widehat{u}(\xi, z)=\mathscr{F}_{(x, y) \rightarrow \xi}(u)(\xi, z)$, and $\xi=\left(\xi_{1}, \xi_{2}\right) .$ Define the phase function $\Phi$ as follows
$$
\Phi(t, \xi) \doteq\left(\tau_{0}-\lambda \theta(t)\right)|\xi|,
$$
where $\lambda$ is a large enough positive constant which will be chosen precisely later,  and $\theta(t)$ is a key quantity to describe the analytic band of $\left(g_{1}, g_{2}\right)$ :
$$
\left\{\begin{array}{l}
\frac{d}{d t} \theta(t)=\left\|\left(g_{1 \Phi}, g_{2 \Phi}\right)\right\|_{B_{1, \alpha}}+\left(\left|\kappa_{1}\right|+\left|\kappa_{2}\right|\right)\langle t\rangle^{\frac{1}{4}}\left\|e^{\psi} \partial_{z} \varphi\right\|_{L_{z}^{2}}, \\
\left.\theta\right|_{t=0}=0.
\end{array}\right.
$$
Let $s$ be in $\mathbb{R}$. For $g$ in $\mathscr{S}_{h}^{\prime}\left(\mathbb{R}_{+}^{3}\right)$, which means that $g$ is in $\mathscr{S}^{\prime}\left(\mathbb{R}_{+}^{3}\right)$ and satisfies $\lim\limits _{k \rightarrow-\infty}\left\|S_{k}^{h} g\right\|_{L^{\infty}}=0$, we set
$$
\left\|g_{\Phi}\right\|_{X_{s, \alpha}}=\left\|\left(2^{k s}\left\|e^{\psi} \Delta_{k}^{h} g_{\Phi}\right\|_{L_{+}^{2}}\right)_{k}\right\|_{\ell^{1}(\mathbb{Z})} .
$$
For compactness of notation, for a function $g_{\Phi}$ such that $g_{\Phi}, \zeta g_{\Phi}, \partial_{z} g_{\Phi} \in X_{s, \alpha}$,   use the time weighted norm
$$
\left\|g_{\Phi}\right\|_{B_{s, \alpha}}=\langle t\rangle^{\frac{1}{4}}\left\|g_{\Phi}\right\|_{X_{s, \alpha}}+\langle t\rangle^{\frac{1}{4}}\left\|\zeta g_{\Phi}\right\|_{X_{s, \alpha}}+\langle t\rangle^{\frac{3}{4}}\left\|\partial_{z} g_{\Phi}\right\|_{X_{s, \alpha}},
$$
where $\zeta=\zeta(t, z)=\frac{z}{\langle t\rangle^{1 / 2}}$, the heat self-similar variable.
In order to obtain a better description of the regularizing effect of the transport diffusion equation,   use Chemin-Lerner type spaces $\widetilde{L}_{T}^{p}\left(X_{s, \alpha}\left(\mathbb{R}_{+}^{3}\right)\right)$, with the norm and
$$
\|g\|_{L_{T}^{p}\left(X_{s, \alpha}\right)} \doteq \sum_{k \in \mathbb{Z}} 2^{k s}\left(\int_{0}^{T}\left\|e^{\psi} \Delta_{k}^{h} g_{\Phi}(t)\right\|_{L_{+}^{2}}^{p} d t\right)^{\frac{1}{p}},\quad
e^{\psi}=\exp \left(\frac{\alpha z^{2}}{4\langle t\rangle}\right) ,
$$
with the usual change if $p=\infty$.
Let $f(t) \in L_{l o c}^{1}\left(\mathbb{R}_{+}\right)$ be a nonnegative function. Define
$$
\begin{aligned}
\|g\|_{\tilde{L}_{t, f}^{p}\left(X_{s, \alpha}\right)} \doteq & \sum_{k \in \mathbb{Z}} 2^{k s}\left(\int_{0}^{t} f\left(t^{\prime}\right)\left\|e^{\psi} \Delta_{k}^{h} g_{\Phi}\left(t^{\prime}\right)\right\|_{L_{+}^{2}}^{p} d t^{\prime}\right)^{\frac{1}{p}}, \\
\|g\|_{L_{t, f}^{p}\left(B_{s, \alpha}\right)} & \doteq \sum_{k \in \mathbb{Z}} 2^{k s}\left[\int _ { 0 } ^ { t } f ( t ^ { \prime } ) \left(\left|t^{\prime}\right\rangle^{\frac{p}{4}}\left(\left\|e^{\psi} \Delta_{k}^{h} g_{\Phi}\right\|_{L_{+}^{2}}^{p}+\left\|\zeta e^{\psi} \Delta_{k}^{h} g_{\Phi}\right\|_{L_{+}^{2}}^{p}\right)\right.\right.\\
&\left.\left.+\left\langle t^{\prime}\right\rangle^{\frac{3 p}{4}}\left\|e^{\psi} \Delta_{k}^{h} \partial_{z} g_{\Phi}\right\|_{L_{+}^{2}}^{p}\right) d t^{\prime}\right]^{\frac{1}{p}},
\end{aligned}
$$
where $p<\infty$. When $p=\infty$, define
$$
\|g\| \widetilde{L}_{T, f(t)}^{\infty}\left(X_{s, \alpha}\right) \doteq \sum_{k \in \mathbb{Z}} 2^{k s} \sup _{t \in[0, T]} f(t)\left\|e^{\psi} \Delta_{k}^{h} g_{\Phi}(t)\right\|_{L_{+}^{2}}.
$$
\begin{Theorem}(\cite{LZ1})
There exist $C_{*}>0$ and $\varepsilon_{*}>0$ such that for any $\varepsilon \in\left(0, \varepsilon_{*}\right]$, when $e^{\tau_{0}\left|D_{h}\right|} g_{i 0} \in$ $X_{1, \frac{1}{2}}, i=1,2$, and
$$
\left\|e^{\tau_{0}\left|D_{h}\right|} g_{10}\right\|_{X_{1, \frac{1}{2}}}+\left\|e^{\tau_{0}\left|D_{h}\right|} g_{20}\right\|_{X_{1, \frac{1}{2}}} \leq \varepsilon,
$$
with
$$
\tau_{0} \geq \frac{C_{*}}{\ln \frac{1}{\varepsilon}}+C_{*}\left(\left|\kappa_{1}\right|+\left|\kappa_{2}\right|\right) \exp \left(\frac{1}{\varepsilon \ln \frac{1}{\varepsilon}}\right),
$$
the system has a unique solution $\left(g_{1}, g_{2}\right)$ on $\left[0, T_{\varepsilon}\right]$, where
$$
T_{\varepsilon} \geq \exp \left(\frac{\varepsilon^{-1}}{\ln \left(\varepsilon^{-1}\right)}\right).
$$
Here $|D_{h}|$ denotes the Fourier multiplier with the symbol $|\xi|$.
Furthermore, the phase function $\Phi(t, \xi)$  of the solution $\left(g_{1}, g_{2}\right)$ satisfies
$$
\Phi(t, \xi) \geq \frac{3 \tau_{0}}{4}|\xi|,
$$
for all $t \in\left[0, T_{\varepsilon}\right]$, and the solution $\left(g_{1}, g_{2}\right)$ obeys the bounds
$$
\begin{array}{l}
\left\|\left(g_{1 \Phi}, g_{2 \Phi}\right)\right\|_{\tilde{L}_{T_{\varepsilon}, f_{3}(t)}^{\infty}\left(X_{1, \alpha}\right)} \leq C \varepsilon, \\
\left\|\left(g_{1 \Phi}, g_{2 \Phi}\right)\right\|_{\widetilde{L}_{T_{\varepsilon}, f_{2}(t)}^{2}\left(B_{1, \alpha}\right)} \leq \frac{\varepsilon C \sqrt{C_{0}}}{\sqrt{\delta}},
\end{array}
$$
where $\delta=\varepsilon \ln \frac{1}{\varepsilon}, \alpha=\frac{1-\delta}{2}, f_{2}(t)=\langle t\rangle^{1-2 \delta}, f_{3}(t)=\langle t\rangle^{\frac{5}{4}-\delta}$, and $C_{0}$ is a positive constant.
\end{Theorem}

\subsection{Ill-posedness of the 3D Prandtl Equations}

 Liu,  Wang and   Yang (\cite{lwy2}) in 2016 gave an instability (ill-posedness) criterion for the 3D Prandtl equations, which shows that the monotonicity condition is not sufficient for the well-posedness of the 3D Prandtl equations, in contrast to the classical well-posedness theory of the 2D Prandtl equations under the Oleinik monotonicity assumption. \\

They considered the  3D Prandtl equations (\ref{w.100}) in $\Omega \triangleq\left\{(t, x, y, z): t>0,(x, y) \in \mathbb{T}^{2}, z \in \mathbb{R}_{+}\right\}$.  Let $u^{s}(t, z)$ and $v^{s}(t, z)$  be smooth solutions of the heat equations:
 \begin{equation}\label{w.64}
\left\{\begin{array}{ll}
{\partial_{t} u^{s}-\partial_{z}^{2} u^{s}=0,} \ \ {\partial_{t} v^{s}-\partial_{z}^{2} v^{s}=0}, \\
 {\left.\left(u^{s}, v^{s}\right)\right|_{z=0}=0,} \ \ {\lim\limits_{z \rightarrow+\infty}\left(u^{s}, v^{s}\right)=\left(U_{0}, V_{0}\right)}, \\
 {\left.\left(u^{s}, v^{s}\right)\right|_{t=0}=\left(U_{s}, V_{s}\right)(z)}.  \end{array}\right.
\end{equation}

\begin{Theorem}(\cite{lwy2} ) \label{lwy2.t}
Let $\left(u^{s}, v^{s}\right)(t, z)$ be the solution of the problem (\ref{w.64}) satisfying
$$
\begin{aligned}\left(u^{s}-U_{0}, v^{s}-V_{0}\right) \in & C^{0}\left(\mathbb{R}_{+} ; W_{\alpha}^{4, \infty}\left(\mathbb{R}_{+}\right) \cap H_{\alpha}^{4}\left(\mathbb{R}_{+}\right)\right)  \cap C^{1}\left(\mathbb{R}_{+} ; W_{\alpha}^{2, \infty}\left(\mathbb{R}_{+}\right) \cap H_{\alpha}^{2}\left(\mathbb{R}_{+}\right)\right). \end{aligned}
$$
Assume that the initial data of the problem (\ref{w.64}) satisfies that there exists a $ z_{0}>0$, such that
$$
  V_{s}^{\prime}\left(z_{0}\right) U_{s}^{\prime \prime}\left(z_{0}\right) \neq U_{s}^{\prime}\left(z_{0}\right) V_{s}^{\prime \prime}\left(z_{0}\right).
$$
Then we have the following two instability statements.

(i) There exists a $\sigma>0$ such that for all $\delta>0, m>0$ and $\mu \in\left[0, \frac{1}{4}\right)$,
$$
\sup _{0 \leqq s \leqq t \leqq \delta}\left\|e^{-\sigma(t-s) \sqrt{| \partial \mathscr{T}| }} T(t, s)\right\|_{\mathscr{L}\left(\mathscr{H}_{\alpha}, \mathscr{H}_{\alpha}^{m-\mu}\right)}=+\infty,
$$
where the operator $\partial_{\mathscr{T}}$ represents the tangential derivative $\partial_{x}$ or $\partial_{y}$,

(ii) There exists an initial shear layer $\left(U_{s}, V_{s}\right)$  to the problem (\ref{w.64}) and $\sigma>0$, such that for all $\delta>0$,
\begin{equation*}
\sup _{0 \leqq s \leqq t \leqq \delta}\left\|e^{-\sigma(t-s) \sqrt{|\partial \mathscr{T}|}} T(t, s)\right\|_{\mathscr{L}\left(\mathscr{H}_{\alpha}^{m_{1}}, \mathscr{H}_{\alpha}^{m_{2}}\right)} =+\infty, \quad \forall m_{1}, m_{2}>0.
\end{equation*}
\end{Theorem}
Here denote by
 \begin{equation*}
\left\{\begin{array}{ll}
   T(t,s)(u_{0},v_{0}):=(u,v)(t,\cdot), \quad \mathscr{H}_{\alpha}^{m}:=H^{m}\left(\mathbb{T}^{2}_{x,y};L^{2}_{\alpha}(\mathbb{R}_{+})\right),\\
  \left\| T(t, s)\right\|_{\mathscr{L}\left(\mathscr{H}_{\alpha}^{m_{1}}, \mathscr{H}_{\alpha}^{m_{2}}\right)} :=
\sup\limits_{(u_{0},v_{0})\in E_{\alpha,\beta}}
\frac{\| T(t,s)(u_{0},v_{0})\|_{\mathscr{H}_{\alpha}^{m_{2}} } }   {\| (u_{0},v_{0})\|_{\mathscr{H}_{\alpha}^{m_{1}} } }\in \mathbb{R}_{+}\cup\{\infty\}, \\
  L^{2}_{\alpha}(\mathbb{R}_{+}):=\{f=f(z), z\in\mathbb{R}_{+}; \|f\|_{L^{2}_{\alpha}}=\|e^{\alpha z}f\|_{L^{2}}<\infty \}, \\
 H^{m}_{\alpha}(\mathbb{R}_{+}):=\{f=f(z), z\in\mathbb{R}_{+}; \|f\|_{H^{m}_{\alpha}}=\|e^{\alpha z}f\|_{H^{m}}<\infty \}.
\end{array}\right.
\end{equation*}

\begin{Theorem} (\cite{lwy2}) \label{lwy2.d}
The problem (\ref{w.100}) with the initial data $\left.(u, v)\right|_{t=0}=\left(u_{0}, v_{0}\right)(x, y, z)$ is locally well-posed, if there exist positive continuous functions $T(\cdot, \cdot), C(\cdot, \cdot)$, some $\alpha>0$ and integer $m \geq 1$ such that for any initial data $\left(u_{0}^{1}, v_{0}^{1}\right)$ and $\left(u_{0}^{2}, v_{0}^{2}\right)$ with
$$
\left(u_{0}^{1}-U_{0}, v_{0}^{1}-V_{0}\right) \in \mathscr{H}_{\alpha}^{m}, \quad\left(u_{0}^{2}-U_{0}, v_{0}^{2}-V_{0}\right) \in \mathscr{H}_{\alpha}^{m},
$$
there are unique distributional solutions $\left(u ^{1}, v ^{1}\right)$ and $\left(u ^{2}, v ^{2}\right)$ satisfying that for
$(i=1, 2)$, $\left.\left(u^{i}, v^{i}\right)\right|_{t=0}=\left(u_{0}^{i}, v_{0}^{i}\right)$ and
\begin{eqnarray*}
\left(u^{i}-U_{0}, v^{i}-V_{0}\right) \in L^{\infty}\left(0, T ; L^{2}\left(\mathbb{T}^{2} \times \mathbb{R}_{+}\right)\right)
 \cap L^{2}\left(0, T ; H^{1}\left(\mathbb{T}^{2} \times \mathbb{R}_{+}\right)\right),
\end{eqnarray*}
and the following estimate holds
$$
\begin{array}{l}{\left\|\left(u^{1}, v^{1}\right)-\left(u^{2}, v^{2}\right)\right\|_{L^{\infty}\left(0, T ; L^{2}\left(\mathbb{T}^{2} \times \mathbb{R}_{+}\right)\right)}+\left\|\left(u^{1}, v^{1}\right)-\left(u^{2}, v^{2}\right)\right\|_{L^{2}\left(0, T ; H^{1}\left(\mathbb{T}^{2} \times \mathbb{R}_{+}\right)\right)}} \\ {\quad \leqq C\left\|\left(u_{0}^{1}-u_{0}^{2}, v_{0}^{1}-V_{0}^{2}\right)\right\|_{ \mathscr{H}_{\alpha}^{m}}}.\end{array}
$$
where both of $C$ and $T$ depend on $\left(u_{0}^{1}-U_{0}, v_{0}^{1}-V_{0}\right)_{\mathscr{H}_{\alpha}^{m}}$ and $\left(u_{0}^{2}-U_{0}, v_{0}^{2}-V_{0}\right) \in \mathscr{H}_{\alpha}^{m}$.
\end{Theorem}

\begin{Theorem}(\cite{lwy2})
Under the same assumption as given in Theorem \ref{lwy2.t}, the problem (\ref{w.100}) of the 3D nonlinear Prandtl equations is not locally well-posed in the sense of Definition \ref{lwy2.d}.
\end{Theorem}

\section{MHD Boundary Layer Equations}

Liu, Xie and Yang (\cite{lxy2}) in 2019 studied the high Reynolds numbers limit for the MHD
system with Prandtl boundary layer expansion when no-slip boundary condition is imposed on velocity field
and perfect conducting boundary condition on magnetic field. Under the assumption that the viscosity and
resistivity coefficients are of the same order and the initial tangential magnetic field on the boundary is not
degenerate, they justified the validity of the Prandtl boundary layer expansion and gave a $L^{\infty}$ estimate on the
error by multi-scale analysis.
Consider the following 2D incompressible viscous MHD equations in the domain $\left\{(t, x, y) \mid t>0, x \in \mathbb{T}, y \in \mathbb{R}_{+}\right\}$,

\begin{equation}
\left\{\begin{array}{l}
\partial_{t} \mathbf{u}^{\epsilon}+\left(\mathbf{u}^{\epsilon} \cdot \nabla\right) \mathbf{u}^{\epsilon}+\nabla p^{\epsilon}-\left(\mathbf{H}^{\epsilon} \cdot \nabla\right) \mathbf{H}^{\epsilon}=\mu \epsilon \Delta \mathbf{u}^{\epsilon}, \\
\partial_{t} \mathbf{H}^{\epsilon}+\left(\mathbf{u}^{\epsilon} \cdot \nabla\right) \mathbf{H}^{\epsilon}-\left(\mathbf{H}^{\epsilon} \cdot \nabla\right) \mathbf{u}^{\epsilon}=\kappa \epsilon \Delta \mathbf{H}^{\epsilon} ,\\
\nabla \cdot \mathbf{u}^{\epsilon}=0, \quad \nabla \cdot \mathbf{H}^{\epsilon}=0.
\end{array}\right.  \label{y2.6}
\end{equation}

Here $\mathbf{u}^{\epsilon}=\left(u^{\epsilon}, v^{\epsilon}\right)$ and $\mathbf{H}^{\epsilon}=\left(h^{\epsilon}, g^{\epsilon}\right)$ stand for the velocity field and magnetic field respectively, $p^{\epsilon}$ denotes the total pressure, the tangential variable is periodic: $x \in \mathbb{T},$ and the normal variable $y \in \mathbb{R}_{+} .$ Also assume $\mu, \kappa$ are positive constants, and the viscosity and resistivity coefficients are of the same order in a small parameter $\epsilon$. The initial data of the system (\ref{y2.6}) is given by
$$
\left.\left(\mathbf{u}^{\epsilon}, \mathbf{H}^{\epsilon}\right)\right|_{t=0}=\left(\mathbf{u}_{0}, \mathbf{H}_{0}\right)(x, y)=\left(u_{0}, v_{0}, h_{0}, g_{0}\right)(x, y),
$$
independent of $\epsilon$. The no-slip boundary condition is imposed on velocity field, and the perfectly conducting boundary condition on magnetic field:
$$
\left.\mathbf{u}^{\epsilon}\right|_{y=0}=\mathbf{0},\left.\quad\left(\partial_{y} h^{\epsilon}, g^{\epsilon}\right)\right|_{y=0}=\mathbf{0}.
$$
\begin{Theorem}(\cite{lxy2})
Let $m \geqslant 36$ be an integer. Let the initial data $\left(\mathbf{u}_{\mathbf{0}}, \mathbf{H}_{\mathbf{0}}\right)(x, y) \in H^{m}\left(\mathbb{T} \times \mathbb{R}_{+}\right)$ satisfy:

i) $\mathbf{u}_{\mathbf{0}}=\left(u_{0}, v_{0}\right)$ is a divergence free vector field vanishing on the boundary, and $\mathbf{H}_{\mathbf{0}}=\left(h_{0}, g_{0}\right)$ is a divergence free vector field tangent to the boundary;

ii) there exists a small $\epsilon_{0}>0$ such that the following ``strong" compatibility conditions hold for any $\epsilon \in\left[0, \epsilon_{0}\right],$
$$
\left.\partial_{t}^{i}\left(u^{\epsilon}, v^{\epsilon}, g^{\epsilon}\right)(0)\right|_{y=0}=\mathbf{0},\left.\partial_{t}^{i-1} \partial_{y} h^{\epsilon}(0)\right|_{y=0}=0, \quad 1 \leqslant i \leqslant\left[\frac{m}{2}\right]-3,
$$
where $[k], k \in \mathbb{R}$,  stands for the largest integer less than or equal to $k$, and $\partial_{t}^{i}\left(u^{\epsilon}, v^{\epsilon}, h^{\epsilon}, g^{\epsilon}\right)(0)$ is the $i$-th time derivative at $\{t=0\}$ of any solution of the equation, as calculated from the equation to yield an expression in terms of derivatives of $\left(\mathbf{u}_{\mathbf{0}}, \mathbf{H}_{\mathbf{0}}\right)$;

iii) the initial tangential magnetic filed is non-degenerate:
$h_{0}(x, 0) \geqslant \delta_{0}>0 $ for some constant $\delta_{0}$.
Then there exists $T_{*}>0$ independent of $\epsilon$ such that, the problem admits a solution $\left(\mathbf{u}^{\epsilon}, \mathbf{H}^{\epsilon}\right)$ in the time interval $\left[0, T_{*}\right],$ and there exists a smooth solution $\left(\mathbf{u}^{\mathbf{0}}, \mathbf{H}^{\mathbf{0}}\right)(t, x, y) \in C\left(\left[0, T_{*}\right], H^{m}\right)$ to the problem with the initial data $\left(\mathbf{u}_{\mathbf{0}}, \mathbf{H}_{\mathbf{0}}\right)$, and a boundary layer profile
$$
\left(u_{b}^{0}, v_{b}^{0}, h_{b}^{0}, g_{b}^{0}\right) \in C\left(\left[0, T_{*}\right] \times \mathbb{T} \times \mathbb{R}_{+}\right),
$$
such that for any arbitrarily small $\sigma>0$,
$$
\sup _{0 \leqslant t \leqslant T_{*}}\left\|\left(\mathbf{u}^{\epsilon}, \mathbf{H}^{\epsilon}\right)(t, x, y)-\left(\mathbf{u}^{\mathbf{0}}, \mathbf{H}^{\mathbf{0}}\right)(t, x, y)-\left(u_{b}^{0}, \sqrt{\epsilon} v_{b}^{0}, h_{b}^{0}, \sqrt{\epsilon} g_{b}^{0}\right)\left(t, x, \frac{y}{\sqrt{\epsilon}}\right)\right\|_{L_{x y}^{\infty}} \leqslant C \epsilon^{3 / 8-\sigma},
$$
where the constant $C>0$ is independent of $\epsilon$.
\end{Theorem}

Wang and Xin (\cite{WX}) in 2017 identified a non-trivial class of initial data for which the uniform stability of the Prandtl's type
boundary layers  can be  established, and  proved rigorously that the solutions to the viscous and diffusive incompressible MHD systems converges strongly to the superposition of the solution to the ideal MHD systems with a Prandtl's type boundary layer corrector.

To this end, they considered  zero viscosity-diffusion vanishing inviscid limit and zero magnetic diffusion vanishing limit for the 3D/2D incompressible viscous and diffusive magnetohydrodynamic (MHD) systems with Dirichlet boundary (no-slip characteristic) boundary conditions,
\begin{eqnarray}\left\{\begin{array}{ll}
u^\varepsilon_t-\varepsilon_1\Delta u^\varepsilon+u^\varepsilon\cdot\nabla u^\varepsilon+\nabla p=b^\varepsilon\cdot\nabla b^\varepsilon,\\
b^\varepsilon_t-\varepsilon_2\Delta b^\varepsilon+u^\varepsilon\cdot\nabla b^\varepsilon=b^\varepsilon\cdot\nabla u^\varepsilon,\\
div u^\varepsilon=div b^\varepsilon=0,\\
u^\varepsilon=b^\varepsilon=0,\ in \ \partial\Omega\times(0,T),\\
u^\varepsilon(t=0)=u_0^\varepsilon,\ b^\varepsilon(t=0)=b_0^\varepsilon.
\end{array}\right.\label{113.1}\end{eqnarray}
Letting $\varepsilon \rightarrow0$ in \eqref{113.1}, we obtain  formally the following inviscid
MHD system
\begin{eqnarray}\left\{\begin{array}{ll}
\partial_tu^{0}+u^{0}\cdot\nabla u^{0}+\nabla p^{0}=b^{0}\cdot \nabla b^{0}, \ \Omega\times (0,T),\\
\partial_tb^{0}+u^{0}\cdot\nabla b^{0}0=b^{0}\cdot\nabla u^{0},\\
div u^{0}=div b^{0}=0,\\
u^{0}=0, b^{0}=0,   \partial\Omega\times(0,T),\\
(u^{0},b^{0})|_{t=0}=(u_0^{0},b_0^{0}).
\end{array}\right.\label{113.2}\end{eqnarray}
\begin{Theorem} ( \cite{WX})
Let $(u^0, p^0, b^0)$ be the solution to the incompressible ideal MHD system \eqref{113.2}. Assume that $(u_\varepsilon^0, b_\varepsilon^0)$ strongly converges in $L^2(\Omega)$ to $(u^0_0, b^0_0)$, where $(u^0_0, b^0_0)\in H^s(\Omega)$ satisfies $divu^0_0=divb^0_0=0$,
and $u^0_0\cdot n=b^0_0\cdot n=0$. Assume that $u^0_0(x,y,z)=b^0_0(x,y,z)$ or $u^0_0(x,y,z)=-b^0_0(x,y,z)$. Furthermore, assume that $\varepsilon,\varepsilon_1,\varepsilon_2$ satisfy the following convergence:
\begin{eqnarray}
 \frac{\varepsilon_1+\varepsilon_2}{\sqrt\varepsilon}\rightarrow0,\ \frac{(\varepsilon_1-\varepsilon_2)^2}{\sqrt\varepsilon\varepsilon(\varepsilon_1+\varepsilon_2)}\rightarrow0,\
\frac{(\varepsilon_1-\varepsilon_2)^2}{\varepsilon(\varepsilon_1+\varepsilon_2)}\leq C\min\{\varepsilon_1,\varepsilon_2\},\label{113.3}
\end{eqnarray}
for some constant $C>0$, independent of $\varepsilon,\varepsilon_1,\varepsilon_2$ as $\varepsilon\rightarrow 0,\ \varepsilon_1\rightarrow 0,\ \varepsilon_1\rightarrow 0$. Then there exists a global Leray-Hopf weak solutions $(u^\varepsilon, p^\varepsilon, b^\varepsilon )$ of the system \eqref{113.1} such that
\begin{eqnarray}
 (u^\varepsilon-u^0,  b^\varepsilon-b^0)\rightarrow(0,0)\ in \ L^\infty(0,T;L^2(\Omega)),\label{113.4}
\end{eqnarray}
for any $T : 0 < T < \infty$, as viscosity coefficient $\varepsilon_1\rightarrow0$ and diffusion coefficient $\varepsilon_2\rightarrow0$.
Moreover, if
\begin{eqnarray}
 \|(u_0^\varepsilon-u_0^0,  b_0^\varepsilon-b_0^0)\rightarrow(0,0)\|_{L^2(\Omega)}\leq C\varepsilon^k,\ k\geq 1,\label{113.4}
 \end{eqnarray}
then there exists $C(T)$, independent of $\varepsilon$ such that, for $0<t\leq T<+\infty$,
\begin{eqnarray}
 \|(u^\varepsilon-u^0,  b^\varepsilon-b^0)\rightarrow(0,0)\|_{L^2(\Omega)}\leq C(T)(\varepsilon^{k-1}+\varepsilon_1^2+\varepsilon_2^2+\frac{\varepsilon_1+\varepsilon_2}{\sqrt\varepsilon}
+\frac{(\varepsilon_1-\varepsilon_2)^2}{\sqrt\varepsilon\varepsilon(\varepsilon_1+\varepsilon_2)}),\ k\geq 1.\label{113.5}
\end{eqnarray}
Furthermore, we have the following stronger $L^\infty$ convergence results for the viscous and diffusive incompressible 2D MHD systems:
Assume that $u_1^0(x,z=0)=b_1^0(x,z=0)=constant$ and $\varepsilon_1=\varepsilon_2$ or that $\varepsilon,\varepsilon_1,\varepsilon_2$ satisfy suitable relations (stated below). If, for some suitably large $k>2$,
\begin{eqnarray}
 \|(u_0^\varepsilon-u_0^0-u_B^{\varepsilon}(t=0),  b_0^\varepsilon-b_0^0-b_B^{\varepsilon}(t=0))\|_{H^s(\Omega)}\leq C\varepsilon^k,\ s\geq 3,\label{113.6}
 \end{eqnarray}
then there exists a $C(T)>0$, independent of $\varepsilon$, such that, for $0\leq  t \leq T < \infty$,
\begin{eqnarray}
 \|(u^\varepsilon-u^0-u_B^{\varepsilon},  b^\varepsilon-b^0-b_B^{\varepsilon})\|_{L^\infty(\Omega)}\leq C(T)\varepsilon^{1/2},\ if\ \varepsilon_=\varepsilon_1
 =\varepsilon_2 .\label{113.7}
 \end{eqnarray}
\end{Theorem}

Wang, Wang, Liu, Wang (\cite{wwlw}) in 2017 studied the boundary layer problem and zero viscosity-diffusion limit of
the following  initial boundary value problem for the incompressible viscous and diffusive magnetohydrodynamic
(MHD) system with (no-slip characteristic) Dirichlet boundary conditions and  proved that the corresponding Prandtl's type boundary layer are stable
 with respect to small viscosity-diffusion coefficients:
\begin{eqnarray}\left\{\begin{array}{ll}
\partial_tu^{\epsilon}+u^{\epsilon}\cdot\nabla u^{\epsilon}-\epsilon\Delta u^{\epsilon}=b^{\epsilon}\cdot \nabla b^{\epsilon},in\ \Omega\times (0,T),\\
\partial_tb^{\epsilon}+u^{\epsilon}\cdot\nabla b^{\epsilon}-\epsilon\Delta b^{\epsilon}=b^{\epsilon}\cdot\nabla u^{\epsilon},\\
div u^{\epsilon}=div b^{\epsilon}=0,\\
u^{\epsilon}=0, b^{\epsilon}=0, in~ \partial\Omega\times(0,T),\\
(u^{\epsilon},b^{\epsilon})|_{t=0}=(u_0^{\epsilon},b_0^{\epsilon})
\end{array}\right.\label{108.1}\end{eqnarray}
and
\begin{eqnarray}\left\{\begin{array}{ll}
\partial_tu^{\eta,\nu}+u^{\eta,\nu}\cdot\nabla u^{\eta,\nu}-\nu_1\partial_z^2u^{\eta,\nu}-\eta_1\Delta_{x,y} u^{\eta,\nu}+\nabla p^{\eta,\nu}=b^{\eta,\nu}\cdot \nabla b^{\eta,\nu},in\ \Omega\times (0,T),\\
\partial_tb^{\eta,\nu}+u^{\eta,\nu}\cdot\nabla b^{\eta,\nu}-\nu_2\partial_z^2b^{\eta,\nu}-\eta_2\Delta_{x,y} b^{\eta,\nu}=b^{\eta,\nu}\cdot\nabla u^{\eta,\nu},\\
div u^{\eta,\nu}=div b^{\eta,\nu}=0,\\
u^{\eta,\nu}=0, b^{\eta,\nu}=0, in~ \partial\Omega\times(0,T),\\
(u^{\eta,\nu},b^{\eta,\nu})|_{t=0}=(u_0^{\eta,\nu},b_0^{\eta,\nu})
\end{array}\right.\label{108.2}\end{eqnarray}

First, letting $\epsilon\rightarrow 0$ in \eqref{108.1}, we formally obtain  the following 3D incompressible ideal MHD system
\begin{eqnarray}\left\{\begin{array}{ll}
\partial_tu^{0}+u^{0}\cdot\nabla u^{0}+\nabla p^{0}=b^{0}\cdot \nabla b^{0},in\ \Omega\times (0,T),\\
\partial_tb^{0}+u^{0}\cdot\nabla b^{0}0=b^{0}\cdot\nabla u^{0},\\
div u^{0}=div b^{0}=0,\\
u^{0}=0, b^{0}=0, in ~ \partial\Omega\times(0,T),\\
(u^{0},b^{0})|_{t=0}=(u_0^{0},b_0^{0}).
\end{array}\right.\label{108.3}\end{eqnarray}
Secondly, setting $\nu=(\eta,\nu)\rightarrow 0$ (which is called the anisotropic MHD limit in the following) in \eqref{108.1} one formally gets the following three-dimensional incompressible anisotropic MHD system
\begin{eqnarray}\left\{\begin{array}{ll}
\partial_tu^{\eta,0}+u^{\eta,0}\cdot\nabla u^{\eta,0}-\eta_1\Delta_{x,y} u^{\eta,0}+\nabla p^{\eta,0}=b^{\eta,0}\cdot \nabla b^{\eta,0},in\ \Omega\times (0,T),\\
\partial_tb^{\eta,0}+u^{\eta,0}\cdot\nabla b^{\eta,0}-\eta_2\Delta_{x,y} b^{\eta,0}=b^{\eta,0}\cdot\nabla u^{\eta,0},\\
div u^{\eta,0}=div b^{\eta,0}=0,\\
u_3^{\eta,0}=0, b_3^{\eta,0}=0, in~ \partial\Omega\times(0,T),\\
(u^{\eta,0},b^{\eta,0})|_{t=0}=(u_0^{\eta,0},b_0^{\eta,0}).
\end{array}\right.\label{108.4}\end{eqnarray}

For the MHD system \eqref{108.1}, the authors of \cite{wwlw} had the following stability result of Prandtl's type
boundary layer for a class of special initial data, which yields to the appearance of the Prandtl's type boundary layer.
\begin{Theorem}( \cite{wwlw})
Let $(u^0,p^0,b^0)(x,y,z,t)$  be
the solution of the incompressible ideal MHD equations \eqref{108.3} with the initial data
$(u_0^0,b_0^0)(x,y,z)$ satisfying $u_0^0,b_0^0\in H^s, s>\frac{3}{2}+1$, $div u^{0}=div b^{0}=0$ and
$u_0^0\cdot n|_{\partial\Omega}=b_0^0\cdot n|_{\partial\Omega}=0$. Assume that $(u_0^\epsilon, b_0^\epsilon)$  strongly converges in $(u_0^0, b_0^0)$.
 Also, assume that $u_0^0(x,y,z)= b_0^0(x,y,z)$ or $u_0^0(x,y,z)= -b_0^0(x,y,z)$.  Then there exists a global Leray-Hopf
weak solutions $(u^\epsilon , p^\epsilon, b^\epsilon )$ of the problem \eqref{108.1} such that
$$(u^\epsilon-u^0,b^\epsilon-b^0)\rightarrow (0,0),~ in ~\ L^\infty(0,T;L^2(\Omega))$$
for any $T : 0 < T < \infty$, as viscosity and diffusion coefficient $ \epsilon\rightarrow 0$.

Here, in fact,  $(u^0,p^0,b^0)(x,y,z,t)=(u^0_0,0,b^0_0)(x,y,z)$ is a special steady state solution of the ideal MHD system \eqref{108.3}.   Moreover, if
$$\|(u^\epsilon_0-u^0_0,b^\epsilon_0-b^0_0)\|_{L^2(\Omega)}\leq C\epsilon^\kappa,\kappa>1,$$
then there exists a $C(T) > 0$, independent of $\epsilon$, such that
$$\|(u^\epsilon-u^0,b^\epsilon-b^0)\|_{L^2(\Omega)}\leq C(T)\epsilon^{\min\{\kappa-1,\frac{1}{2}\}}.$$
\end{Theorem}

For anisotropic viscous and diffusive MHD system \eqref{108.2} with the different horizontal and vertical viscosities and magnetic diffusions, the authors of \cite{wwlw} obtained the following result.
\begin{Theorem}( \cite{wwlw})
 (The Ideal MHD Limit) Assume that $(u_0^{\eta,\nu},b_0^{\eta,\nu})$  strongly converges in $L^2(\Omega)$ to $(u_0^{0,0},b_0^{0,0})$ where $u_0^{ 0,0 },b_0^{ 0,0 }\in H^s, s>\frac{3}{2}+1$.
Assume that
$$\eta_1\rightarrow 0,\eta_2\rightarrow 0, \nu_1\rightarrow 0,\nu_2\rightarrow 0,$$
$$\frac{\nu_1}{\sqrt{\min\{\eta_1,\eta_2\}\min\{\nu_1,\nu_2\}}} \rightarrow 0 \ and \ \frac{\nu_2}{\sqrt{\min\{\eta_1,\eta_2\}\min\{\nu_1,\nu_2\}}}\rightarrow 0.$$
Assume that $(u^{0,0},p^{0,0},b^{0,0})$  is the smooth solution to the system \eqref{108.3} defined on $[0,T_\ast)$ with $0<T_\ast\leq \infty$. Then there exists  a global Leray-Hopf weak solution $(u^{\eta,\nu},p^{\eta,\nu},b^{\eta,\nu})$ of the system \eqref{108.2} such that $$(u^{\eta,\nu},b^{\eta,\nu})-(u^{0,0},b^{0,0})\rightarrow (0,0),\ in\ \ L^\infty(0,T;L^2),$$
for any $T:0<T<T^\ast$.

Moreover, if
$$\|(u^{\eta,\nu}_0-u^{0,0}_0,b^{\eta,\nu}_0-b^{0,0}_0)\|_{L^2(\Omega)}\leq C\sqrt{\alpha},$$
where $\alpha=\max\left\{\nu_1+\nu_2+\eta_1+\eta_2,\sqrt{\min\{\eta_1,\eta_2\}\min\{\nu_1,\nu_2\}},
\frac{\nu_1}{\sqrt{\min\{\eta_1,\eta_2\}\min\{\nu_1,\nu_2\}}},
\frac{\nu_2}{\sqrt{\min\{\eta_1,\eta_2\}\min\{\nu_1,\nu_2\}}}\right\}$.
Then there exists a $C=C(T)>0$, independent of $\eta,\nu$, such that
$$\|(u^{\eta,\nu}-u^{0,0},b^{\eta,\nu}-b^{0,0})\|_{L^2(\Omega)}\leq C\sqrt{\alpha}.$$
\end{Theorem}

\begin{Theorem}( \cite{wwlw}) (The Anisotropic MHD Limit) Let $\eta_1>0,\eta_2>0$. Assume that $(u_0^{\eta,\nu},b_0^{\eta,\nu})$  strongly converges in $L^2(\Omega)$ to $(u_0^{\eta,0},b_0^{\eta,0})$ where $u_0^{ \eta,0 },b_0^{ \eta,0 }\in H^s, s>\frac{3}{2}+2$.
Assume that
$$ \nu_1\rightarrow 0,\nu_2\rightarrow 0,$$
$$\frac{\nu_1}{\sqrt{\min\{\nu_1,\nu_2\}}} \rightarrow 0 \ and \ \frac{\nu_2}{\sqrt{\min\{\nu_1,\nu_2\}}}\rightarrow 0.$$
Assume that
$(u^{\eta,0},p^{\eta,0},b^{\eta,0})$  is the smooth solution to the system \eqref{108.4} defined on $[0,T_\ast)$ with $0<T_\ast\leq \infty$. Then there exists a global Leray-Hopf weak solution  $(u^{\eta,\nu},p^{\eta,\nu},b^{\eta,\nu})$ of the system \eqref{108.2} such that $$(u^{\eta,\nu},b^{\eta,\nu})-(u^{\eta,0},b^{\eta,0})\rightarrow (0,0),\ in\ \ L^\infty(0,T;L^2),$$
for any $T:0<T<T^\ast$ as $\nu=(\nu_1,\nu_2)\rightarrow(0,0).$

Moreover, if
$$\|(u^{\eta,\nu}_0-u^{\eta,0}_0,b^{\eta,\nu}_0-b^{\eta,0}_0)\|_{L^2(\Omega)}\leq C\sqrt{\beta},$$
where $\beta=\max\left\{\nu_1+\nu_2,\sqrt{\min\{\nu_1,\nu_2\}},
\frac{\nu_1}{\sqrt{\min\{\nu_1,\nu_2\}}},
\frac{\nu_2}{\sqrt{\min\{\nu_1,\nu_2\}}}\right\}$.

Then there exists a $C=C(T)>0$, independent of $\eta,\nu$, such that
$$\|(u^{\eta,\nu}-u^{\eta,0},b^{\eta,\nu}-b^{\eta,0})\|_{L^2(\Omega)}\leq C\sqrt{\beta}.$$
\end{Theorem}

G$\acute{e}$rand-Varet and Prestipino (\cite{gvp}) in 2017 provided a systematic derivation of boundary layer models in magnetohydrodynamics (MHD), through an asymptotic analysis of
the incompressible MHD system and recovered classical linear models, related to the famous Hartmann and Shercliff layers, as well as nonlinear ones,
that called magnetic Prandtl models. An appropriate starting point to describe such dynamics is the classical incompressible MHD system:
\begin{eqnarray}
\left\{\begin{array}{l}{\partial_{t} \mathbf{u} + \mathbf{u}  \cdot \nabla \mathbf{u} +\nabla p - \frac{1}{Re} \triangle \mathbf{u}=S\mathbf{b}  \cdot \nabla \mathbf{b}}, \\ {\partial_{t} \mathbf{b} -\text{curl}~(\mathbf{u} \times\mathbf{b} )+ \frac{1}{Rm}\text{curl}~\text{curl}\mathbf{b}=0}, \\
 {\nabla \cdot \mathbf{u} =0, \quad \nabla \cdot \mathbf{b} =0,\ \ t>0, \ \ x\in\Omega ,}\end{array}\right. \label{d6.1}
\end{eqnarray}
with
$$S=\frac{B_{0}^{2}}{\mu\rho U^{2}}\frac{Ha^{2}}{ReRm},\ \ Ha=B_{0}L(\frac{\sigma}{\eta})^{\frac{1}{2}}.$$
They considered the solutions of the system (\ref{d6.1}) behaving like:
\begin{eqnarray*}
\left\{\begin{array}{l}
\mathbf{u}\approx\left(u_{x}'(t,x,y,\lambda^{-1}z),u_{y}'(t,x,y,\lambda^{-1}z),\lambda u_{z}'(t,x,y,\lambda^{-1}z)\right),\\
 \mathbf{b}\approx \mathbf(e)_{z} +\delta \left(b_{x}'(t,x,y,\lambda^{-1}z),b_{y}'(t,x,y,\lambda^{-1}z),\lambda b_{z}'(t,x,y,\lambda^{-1}z)\right).
\end{array}\right.
\end{eqnarray*}
The 2D mixed Prandtl/Hartmann regime is as follows
\begin{eqnarray}
\left\{\begin{array}{l}{\partial_{t} u+u \partial_{x} u+v \partial_{z} u-\frac{Ha^{2}}{Re}\partial_{z}^{2} u+\frac{Ha^{2}}{Re}u=-\partial_{x}p}, \\
{\partial_{x} u+\partial_{z} v=0}, \\
{\left.u\right|_{z=0}=\left.v\right|_{z=0}=0, \lim\limits _{z \rightarrow+\infty} u=U(t, x)} ,  \end{array}\right.\label{6.2}
\end{eqnarray}
and the 2D Mixed Prandtl/Shercliff regime is in the following
\begin{eqnarray}
\left\{\begin{array}{l}{\partial_{t} u+u \partial_{x} u+v \partial_{z} u-\frac{Ha}{Re}\partial_{z}^{2} u=\frac{Ha}{Re}\partial_{x}b-\partial_{x}p}, \\
{\partial_{x} u+\partial_{z}^{2} b=0}, \\
{\partial_{x} u+\partial_{z} v=0}, \\
{\left.u\right|_{z=0}=\left.v\right|_{z=0}=\left.b\right|_{z=0}=0, \lim\limits _{z \rightarrow+\infty} u=U(t, x) , \lim\limits _{b \rightarrow+\infty} u=B(t, x)} .\end{array}\right.\label{6.3}
\end{eqnarray}
The linearized system of the problem $(\ref{6.3})$  reads as
\begin{eqnarray}
\left\{\begin{array}{l}{\partial_{t} u+U \partial_{x} u+vU'-\frac{Ha}{Re}\partial_{z}^{2} u=\frac{Ha}{Re}\partial_{x}b} , \\
{\partial_{x} u+\partial_{z}^{2} b=0}, \\
{\partial_{x} u+\partial_{z} v=0}, \\
{\left.u\right|_{z=0}=\left.v\right|_{z=0}=\left.b\right|_{z=0}=0, \lim\limits _{z \rightarrow+\infty} u=  \lim\limits _{b \rightarrow+\infty} u=(0,0)} .  \end{array}\right.\label{6.4}
\end{eqnarray}
\begin{Theorem}(\cite{gvp})
   Assume that $U \in  W^{ 2,\infty} (\mathbb{R}_{ +} ), zU',zU'' \in L^{\infty}(\mathbb{R}_{ +})$. Let $u _{0}\in L^{2}(\Omega)$ such that $w_{0}=\partial_{z}u_{0}\in L^{2}(\Omega)$, $u_{0}|_{z=0}=0$.  Let $b _{0}\in L^{2}_{loc}(\Omega)$    such that $\partial_{z}b_{0}\in L^{2}(\Omega)$, $b_{0}|_{z=0}=0$ and with zero average in x. Then there exists a unique solution $(u,v,b)$ of the problem $(\ref{6.4})$ satisfying $(u,b)|_{ t=0} = (u_{ 0} ,b_{ 0} )$.
\end{Theorem}
The 2D nonlinear MHD layer $(Re\sim Rm\sim Ha)$ is the following:
\begin{eqnarray}
\left\{\begin{array}{l}{\partial_{t} u+u \partial_{x} u+v \partial_{z} u- \partial_{z}^{2} u= S\mathbf{b}  \cdot \nabla \mathbf{b}-\partial_{x}p}, \\
{\partial_{t} \mathbf{b} -\nabla^{T}(\mathbf{u} \times\mathbf{b} )- \frac{Re}{Rm}\partial_{z}^{2}\mathbf{b}  =0}, \\
{\partial_{x} u+\partial_{z} v= div \mathbf{b}  =0}, \\
{\left.u\right|_{z=0}=\left.v\right|_{z=0}=\left.b\right|_{z=0}=0, \lim\limits _{z \rightarrow+\infty} u=U(t, x) , \lim\limits _{b \rightarrow+\infty} b=B(t, x)} ,  \end{array}\right.\label{6.5}
\end{eqnarray}
whose linearized system reads as:
\begin{eqnarray}
\left\{\begin{array}{l}{\partial_{t} u+U \partial_{x} u+U'v- \partial_{z}^{2} u= S\partial_{x}b   }, \\
{\partial_{t} \mathbf{b} -\nabla^{T}(v-Uc )- \frac{Re}{Rm}\partial_{z}^{2}\mathbf{b}  =0}, \\
{\partial_{x} u+\partial_{z} v= div \mathbf{b}  =0}, \\
{\left.u\right|_{z=0}=\left.v\right|_{z=0}=\left.b\right|_{z=0}=0, \lim\limits _{z \rightarrow+\infty} u= \lim\limits _{b \rightarrow+\infty} b=0}. \end{array}\right.\label{6.6}
\end{eqnarray}
\begin{Theorem} (\cite{gvp})
  Assume that $U \in  W^{ 3,\infty} (\mathbb{R}_{ +} ), zU',zU'',zU''' \in L^{\infty}(\mathbb{R}_{ +})$. Let $u _{0}\in L^{2}(\Omega)$, $\phi_{0}\in  L^{2}_{loc}(\Omega)$, such that  $b_{0}=\partial_{y}\phi_{0}\in L^{2}(\Omega)$, $\phi_{0}|_{z=0}=0$ and with zero average in x. Then there exists a unique solution   of the problem   $(\ref{6.6})$ satisfying $(u,b)|_{ t=0} = (0 ,-\nabla^{T}\phi_{ 0} )$.
\end{Theorem}

\subsection{MHD Boundary Layer Equations-Local Existence   }

  Oleinik and Samokhin (\cite{OS}) in 1999 considered
the boundary layer theory which, in fact, can also be applied to some problems connected with the boundary layer in magnetohydrodynamics (MHD).
They considered the stationary boundary layer system for a plane-parallel conducting flow in the presence of magnetic and electric fields
in the domain $D =\{0<x<X, 0<y< \infty \}$  by the von Mises variables,  that has the form:
\begin{equation}\left\{
\begin{array}{ll}
  u \frac{\partial u }{\partial x}+v \frac{\partial u }{\partial y}
=  \nu \frac{\partial^{2} u }{\partial y^{2}} - \frac{\partial p }{\partial x}-jB, \\
\frac{\partial u }{\partial x}+\frac{\partial v }{\partial y} =0
\end{array}
 \label{ww74}         \right.\end{equation}
where $j=\sigma(E+B)$ is the $z$-component of the electric current vector; ${\bf E}=(0,0,E)$ is the electric field; ${\bf B}=(0,B,0)$ is the magnetic field;
$p,B,E$ are given functions that depend only on $x$. System (\ref{ww74}) is supplemented with the conditions
 \begin{equation}
\begin{array}{ll}
u (0,y)=u_{0}(y), \ \ u (x,0)=0, \ \ v (x,0)=v_{0}(x), \\
u (x,y)\rightarrow U (x),  \ \ \text{as}\ \   y\rightarrow \infty.
\end{array}
 \label{ww75}         \end{equation}
The velocity $U(x)$ of the outer flow is related to the pressure $p(x)$ and the electromagnetic field components by
$$U\frac{dU}{dx}=-\frac{dp}{dx}-\sigma EB-\sigma B^{2}U. $$
The above equation is a generalization of the Bernoulli's law. Taking into account this relation, we can rewrite (\ref{ww74}) in the form
\begin{equation}\left\{
\begin{array}{ll}
  u \frac{\partial u }{\partial x}+v \frac{\partial u }{\partial y}
=  \nu \frac{\partial^{2} u }{\partial y^{2}} +U \frac{d U }{d x}+\sigma B^{2}(U-u), \\
\frac{\partial u }{\partial x}+\frac{\partial v }{\partial y} =0.
\end{array}
 \label{ww76}         \right.\end{equation}

\begin{Theorem} \label{OS.9.50} (\cite{OS})
For some $X>0$, problem (\ref{ww75})-(\ref{ww76}) in $D$ admits a solution $u,v$ with the following properties: $u$ continuous and bounded
in $\overline{D}$;  $u>0$ for $y>0$; $\partial u/\partial y>m>0$ for $0\leq y\leq y_{0}$; $m,y_{0}=$const.;
$\partial u/\partial y, \partial^{2} u/\partial y^{2}$ are bounded and continuous in $D$; $\partial u/\partial x, v, \partial v/\partial y$ are
continuous and bounded in any finite subdomain of $D$. If
$$\frac{d U }{d x}\geq 0 \ \ \text{and} \ \ v_{0}(x)\leq 0  \ \ \text{or} \ \ U\left(\frac{d U }{d x}+\sigma B^{2}\right)\geq \alpha=\text{const.}>0, $$
then the solution of problem (\ref{ww75})-(\ref{ww76})  exists for any $X>0$.

\end{Theorem}

\begin{Theorem} \label{OS.9.51} (\cite{OS})
Let $u(x,y), v(x,y)$ be a solution of problem (\ref{ww75})-(\ref{ww76})  with the following properties: $u,v$ are continuous in $\overline{D}$; $0<u<C_{0}$ for
$y>0$;
\begin{eqnarray}
&& k_{1}y\leq u\leq k_{2}y \ \ \text{for}\ \  0<y<y_{0}, \label{ww77}   \\
&& \partial^{2} u/\partial y^{2}\leq k_{3} \ \ \text{in} \ \  D,
 \label{ww78}
\end{eqnarray}
where $C_{0}, k_{i}, (i=1,2,3), y_{0}$ positive constants. Such a solution of problem (\ref{ww75})-(\ref{ww76})  is unique.
If $U(dU/dx+\sigma B^{2})\geq 0$ and $v_{0}(x)\leq 0$,  then the conditions (\ref{ww77})-(\ref{ww78}) can be dropped. If, however, $v_{0}(x)$ is allowed to take
positive values, then the condition (\ref{ww77})  should be retained.

\end{Theorem}

\begin{Theorem} \label{OS.9.52}( \cite{OS})
Assume  the following conditions holds
 $$v_{0}(x)\leq 0, \ \ 0\leq\frac{dU}{dx}\leq \frac{M_{0}}{(x+1)^{1+\gamma_{0}}}, \ \ \gamma_{0}>0, \ \ \lim\limits_{x\rightarrow\infty}U(x)=U_{0}.$$
  Let $u(x,y), v(x,y)$ be the solution of problem (\ref{ww75})-(\ref{ww76})  in  $D =\{0<x<\infty, 0<y< \infty \}$.  Assume also that
$$\sigma B^{2}(x)=\frac{U_{0}b(x)}{2(x+1)^{1+\delta_{1}}}, \ \ v_{0}(x)= \frac{\sqrt{\nu U_{0}}q(x)}{\sqrt{2}(x+1)^{\frac{1}{2}+\delta_{2}}}, $$
where $-1\leq \delta_{1}\leq \infty, -1/2\leq \delta_{2}\leq \infty$; $b(x),q(x)$ are bounded functions. If $\delta_{1}>0$,
$b(x)\geq 0$.
 In the case $\delta_{1}\leq 0$, if it is assumed that $b(x)>0$ and, for some $N=$const.$>0$, we have
$$|b(x)-N|\leq \frac{M_{1}}{(x+1)^{1+\delta_{3}} }, \ \ \delta_{3}>0. $$
If $\delta_{2}>0$, we assume that $q(x)\leq 0$; if $\delta_{2}\leq 0$, we assume that $q(x)<0$ and, for some $\alpha=$const.$\leq 0$, the following inequality holds
$$|q(x)-\alpha|\leq \frac{M_{2}}{(x+1)^{1+\delta_{4}} }, \ \ \delta_{4}>0. $$
For definiteness, it is also assumed that $B(x)\geq 0$ and that for $\delta_{1}<0,\delta_{2}<0$,  we have
$$\frac{B(x)}{v_{0}(x)}=\mathcal{O}((x+1)^{1+\delta_{5}} ) \ \ \text{as} \ \  x\rightarrow \infty . $$

Under these assumptions, it is claimed that
$$|u(x,y)-\phi(x,y)|\rightarrow  0 \ \ \text{as} \ \  x\rightarrow \infty$$
uniformly with respect to $y$, where the function $\phi(x,y)$ is defined as follows: \\
(1)  $\phi(x,y)=u_{0,0}(x,y)$ if $\delta_{1}>0, \delta_{2}>0$;\\
(2)  $\phi(x,y)=u_{0,N}(x,y)$ if $\delta_{1}=0, \delta_{2}>0$;\\
(3)  $\phi(x,y)=u_{\alpha,0}(x,y)$ if $\delta_{1}>0, \delta_{2}=0$;\\
(4)  $\phi(x,y)=u_{\alpha,N}(x,y)$ if $\delta_{1}=\delta_{2}=0$;\\
(5) if $-1\leq\delta_{1}<0$ and either $\delta_{2}\geq 0$ or $\delta_{2}<0$ and $\delta_{5}> 0$, then
$$\phi(x,y)=h_{1}(x,y)=U_{0}\left(1-\exp\left\{-\frac{y\sqrt{U_{0}N}}{\sqrt{2\nu(x+1)^{1+\delta_{2}}}} \right\}\right);  $$
(6) if $-\frac{1}{2}\leq\delta_{2}<0$ and either $\delta_{1}\geq 0$ or $\delta_{1}<0$ and $\delta_{5}< 0$, then
$$\phi(x,y)=h_{2}(x,y)=U_{0}\left(1-\exp\left\{ \frac{\alpha y\sqrt{U_{0} }}{\sqrt{2\nu(x+1)^{1+2\delta_{2}}}} \right\}\right);  $$
(7) if $-1\leq\delta_{1}<0, -\frac{1}{2}\leq\delta_{2}<0, \delta_{5}= 0$, then
$$\phi(x,y)=h_{3}(x,y)=U_{0}\left(1-\exp\left\{  C(x)y \right\}\right),  $$
where
$$C(x)=\frac{ \sqrt{U_{0} }(\alpha -\sqrt{\alpha^{2}+4N})}{2\sqrt{2\nu(x+1)^{1+2\delta_{2}}}}.  $$

\end{Theorem}

\begin{Theorem} \label{OS.9.53} (\cite{OS})
Suppose the assumptions of Theorem \ref{OS.9.52} hold for $U(x), v_{0}(x), b(x)$, and   $(u(x,y), v(x,y))$ is a solution of problem (\ref{ww75})-(\ref{ww76}). If
$v_{0}(x)$ and $B(x)$ ensure the convergence $|u(x,y)-\phi(x,y )|\rightarrow 0$, uniformly in $y$, as $x\rightarrow\infty $
(where $\phi(x,y)$ is the same as in Theorem \ref{OS.9.52}), then
$$\left|\tau-\nu \frac{\partial\phi}{\partial y}\right|\rightarrow 0\ \  \text{as} \ \ x\rightarrow\infty . $$
\end{Theorem}

  Oleinik and Samokhin (\cite{OS}) in 1999 also considered
 the MHD boundary in pseudo-plastic fluids as follows:
\begin{equation}\left\{
\begin{array}{ll}
  u \frac{\partial u }{\partial x}+v \frac{\partial u }{\partial y}
=  \nu \frac{\partial}{\partial y }\left(|\frac{\partial u }{\partial y}|^{n-1}\frac{\partial u }{\partial y} \right)+U(x)\frac{d U }{d x}
+\sigma B^{2}(U-u), \ \  n>0, \\
\frac{\partial u }{\partial x}+\frac{\partial v }{\partial y} =0,
\end{array}
 \label{ww79}         \right.\end{equation}
in the domain $D =\{0<x<X, 0<y< \infty \}$, with the boundary conditions
 \begin{equation}
\begin{array}{ll}
u (0,y)=u_{0}(y), \ \ u (x,0)=0, \ \ v (x,0)=v_{0}(x), \\
u (x,y)\rightarrow U (x),  \ \ \text{as}\ \   y\rightarrow+\infty,
\end{array}
 \label{ww80}         \end{equation}
where $U(x), B(x), u_{0}(x), v_{0}(x)$ are given functions.
Assume that
$$U(x) =xV(x), \ \  v_{0}(x)=x^{\frac{n-1}{n+1}}v_{1}(x), $$
where$V(x)>0, v_{1}(x)$ are bounded functions. In  (\ref{ww79}), introduce  new independent variables
$$\xi=x, \ \  \eta=\frac{u}{U}, \ \  w(\xi, \eta)=\frac{|u_{y}|^{n-1}u_{y}}{x^{\frac{n-1}{n+1}}U}, $$
then  problem (\ref{ww79})-(\ref{ww80}) reduces to the equation
 \begin{equation}
\begin{array}{ll}
\nu n V^{\frac{1-n}{n}}|w|^{\frac{1+n}{n} } w_{\eta\eta} -\eta\xi Vw_{\xi}+Aw_{\eta}+Bw=0,
\end{array}
 \label{ww81}         \end{equation}
in the domain $\Omega=\{  0<\xi<X, 0<\eta <1 \}$, with the boundary conditions
\begin{eqnarray}
   w|_{\eta=1}=0, \ \  (\nu w |w|^{\frac{1-n}{n} }_{\eta }-v_{1}w|w|^{\frac{1-n}{n} }+C )|_{\eta=0}=0,
 \label{ww82}
\end{eqnarray}
 where
 $$A =(\eta^{2}-1)(V+\xi V_{x})+(\eta-1)\sigma B^{2} , \ \  B =-\eta \left(\frac{2nV}{1+n}  +\xi V_{x}\right)-\sigma B^{2},
 \ \  C=V^{\frac{ n-1}{n}}(V+\xi V_{x}+\sigma B^{2}). $$

\begin{Theorem} \label{OS.9.54} (\cite{OS})
Assume that
$$U(x) =x(a+xa_{1}(x)), \ \  v_{0}(x)=x^{\frac{n-1}{n+1}}(b+xb_{1}(x)), $$
$a=$const.$>0$, $b,d_{0}=$const.; $U>0$ for $x>0$; $a_{1}, a_{1x},a_{1xx}, b_{1}, b_{1x}$ are bounded for $0\leq x\leq X$. Then problem (\ref{ww79})-(\ref{ww80})
in the domain $D$, with $X$ depending on $U,v_{0},B,n,\nu$, admits a solution $(u(x,y), v(x,y))$ with the following properties: $u/U$ is bounded and continuous
 in $\overline{D}$;  $v$ is continuous in $D$; $v$ is continuous with respect to $y$ in $\overline{D}$ and bounded for bounded $y$;
 $u(x,y)>0$ for $y>0$ and $x>0$; the derivatives $u_{x}, u_{y}, u_{yy}, u_{y}$ are continuous and bounded in $D$ and continuous with respect to $y$ in
 $\overline{D}$; $u_{y}>0$ for $y\geq 0$ and $x>0$; $u$ and $v$ satisfy equations (\ref{ww79}) at every point of $D$; $u(x,y)\rightarrow U(x)$ as
  $y\rightarrow\infty$, uniformly on the segment [0,X];  the boundary conditions (\ref{ww80}) are satisfied. Moreover,
$u_{y}^{n}/\left(x^{\frac{n-1}{n+1}}U\right)$ is bounded and continuous in  $\overline{D}$, $u_{y}\rightarrow 0$ as  $y\rightarrow\infty$;
$u_{yy}^{n}/\left(x^{\frac{n-1}{n+1}}u_{y}^{2-n}\right)$ is continuous with respect to $y$ in $\overline{D}$;  $u_{xy}, u_{yyy}$ are continuous in $D$;
the following inequalities hold:
\begin{eqnarray*}
&& x^{\frac{n-1}{n+1}}U(x)Y(\frac{u}{U})e^{-K_{1}x} \leq u_{y}^{n}\leq x^{\frac{n-1}{n+1}}U(x)Y(\frac{u}{U})e^{ K_{3}x}  , \\
&&  -K_{6}\leq \frac{n(n-2)(u_{yy})^{2}+ nu_{y}u_{yyy}}{x^{\frac{n-1}{n}} U^{\frac{n-1}{n}} u_{y}^{3-n}}\leq -K_{7} , \\
&& \left(M_{2}^{\frac{1}{n}} \frac{1-n}{n+1} e^{\frac{K_{2}}{n}}x^{\frac{n-1}{n(n+1)}}U^{\frac{ 1-n}{n }}y +1\right)^{\frac{n+1}{n-1}}\leq 1-\frac{u}{U}
\leq \left(M_{1}^{\frac{1}{n}} \frac{1-n}{n+1} e^{-\frac{K_{1}}{n}}x^{\frac{n-1}{n(n+1)}}U^{\frac{ 1-n}{n }}y +1\right)^{\frac{n+1}{n-1}}.
\end{eqnarray*} The solution of
problem (\ref{ww79})-(\ref{ww80}) with the above properties is unique.

\end{Theorem}

Now consider a more general problem, with $U(x)=\mathcal{O}(x^{m})$ as $x\rightarrow 0$, $m>0$, that is, system (\ref{ww79}) with the boundary conditions
 \begin{equation}\left\{
\begin{array}{ll}
u (0,y)=0, \ \ u (x,0)=0, \ \ v (x,0)=v_{0}(x), \\
u (x,y)\rightarrow U (x),  \ \ \text{as}\ \   y\rightarrow+\infty.
\end{array}\right.
 \label{ww83}         \end{equation}
Assume that
$$U(x) =x^{m}V(x), \ \  v_{0}(x)=x^{\frac{2mn-n-m}{n+1}}v_{1}(x), \ \  d^{2}(x)=\sigma B(x)^{2}, \ \ d(x)=x^{\frac{m-1}{2}}W(x), $$
where $V(x)>0$  for $0\leq x\leq X, m>0$; $W(x)>0$ for $0\leq x\leq X $;   moreover, the functions $V(x), W(x), v_{1}(x)$ are smooth, $V(0)=a>0, v_{1}(0)=b,
W(0)=d_{0}$. In  (\ref{ww79}),  introduce  new independent variables
$$\xi=x, \ \  \eta=\frac{u}{U}, \ \  w(\xi, \eta)=\frac{|u_{y}|^{n-1}u_{y}}{x^{\frac{2mn-n-m}{n+1}}U}. $$

\begin{Theorem} \label{OS.9.55} (\cite{OS})
 Under the assumptions
\begin{eqnarray*}
&& x(U_{x}+\sigma B(x)^{2})+\frac{2mn-n-m}{n+1}U>0, \ \ \frac{(3m-1)na}{n+1}+\sigma B^{2}(0)>0, \\
&& V(x) =x(a+xa_{1}(x)), W(x)=d_{0}+xd_{1}(x), \ \  v_{1}(x)= b+xb_{1}(x ),
\end{eqnarray*}
where $a >0$ and the functions $a_{1}, a_{1x},a_{1xx},  b_{1}, b_{1x}, b_{1}, b_{1x}$ are bounded for $0\leq x\leq X$. Then the problem  (\ref{ww79}),  (\ref{ww83})
in the domain $D$,  with $X$ depending on $U,v_{0},d,n,\nu$ has a solution $u(x,y), v(x,y)$ with the following properties: $u/U$ is bounded and
continuous in $\overline{D}$;  $v$  is continuous in  $D$; $v$  is continuous with respect to
$y$ in $\overline{D}$ and bounded for bounded $y$; $u(x,y)>0$ for $y>0$ and $x>0$; the derivatives $u_{x}, u_{y}, u_{yy}, u_{y}$ are continuous and bounded in $D$ and continuous with respect to $y$ in
 $\overline{D}$; $u_{y}>0$ for $y\geq 0$ and $x>0$; equations  (\ref{ww79}) hold for $u,v$ at each point of $D$; $u(x,y)\rightarrow U(x)$ as
  $y\rightarrow\infty$, uniformly on the segment $[0,X]$, the boundary conditions (\ref{ww83}) hold for $u,v$. Moreover,
$u_{y}^{n}/\left(x^{\frac{2mn-n-m}{n+1}}U(x)\right)$ is bounded and continuous in $\overline{D}$,  $u_{y}\rightarrow 0$ as  $y\rightarrow\infty$;
$u_{yy}^{n}/\left(x^{\frac{n-1}{n+1}}u_{y}^{2-n}\right)$ is continuous with respect to $y$ in $\overline{D}$;  $u_{xy}, u_{yyy}$
are continuous in $D$; the following inequalities hold:
\begin{eqnarray*}
&& x^{\frac{2n-m-n}{n+1}}U(x)Y(\frac{u}{U})e^{-K_{8}x} \leq u_{y}^{n}\leq x^{\frac{2mn-m-n}{n+1}}U(x)Y(\frac{u}{U})e^{ K_{9}x}  , \\
&&  -K_{13}\leq \frac{n(n-2)(u_{yy})^{2}+ nu_{y}u_{yyy}}{x^{\frac{2mn-m-n}{n}} U^{\frac{n-1}{n}} u_{y}^{3-n}}\leq -K_{14} , \\
&& \left(M_{8}^{\frac{1}{n}} \frac{1-n}{n+1} e^{\frac{K_{9}}{n}}x^{\frac{2mn-m-n}{n(n+1)}}U^{\frac{ 1-n}{n }}y +1\right)^{\frac{n+1}{n-1}}\leq 1-\frac{u}{U}
\leq \left(M_{7}^{\frac{1}{n}} \frac{1-n}{n+1} e^{-\frac{K_{8}}{n}}x^{\frac{2mn-m-n}{n(n+1)}}U^{\frac{ 1-n}{n }}y +1\right)^{\frac{n+1}{n-1}}.
\end{eqnarray*}
 The solution of
problem (\ref{ww79}),  (\ref{ww83}) with the above properties is unique.

\end{Theorem}

  Oleinik and  Samokhin (\cite{OS}) in 1999 continued to consider
 the system of equations of the MHD boundary layer for a dilatable fluid:
\begin{equation}\left\{
\begin{array}{ll}
  u \frac{\partial u }{\partial x}+v \frac{\partial u }{\partial y}
= \frac{k}{\rho} \frac{\partial}{\partial y }\left(|\frac{\partial u }{\partial y}|^{n-1}\frac{\partial u }{\partial y} \right)+U(x)\frac{d U }{d x}
+\frac{\sigma B^{2}(x)}{\rho}(U-u), \ \  n>1, \\
\frac{\partial u }{\partial x}+\frac{\partial v }{\partial y} =0
\end{array}
 \label{ww86}         \right.\end{equation}
in the domain $D =\{0<x<X, 0<y< \infty \}$, with the boundary conditions
 \begin{equation}\left\{
\begin{array}{ll}
u (0,y)=u_{0}(y), \ \ u (x,0)=0, \ \ v (x,0)=0, \\
u (x,y)\rightarrow U (x),  \ \ \text{as}\ \   y\rightarrow+\infty.
\end{array}
 \label{ww87}   \right.      \end{equation}
Assume that
$$U( x) =C_{0}x^{m} , \ \  B( x)=bx^{\frac{m-1}{2}} , $$
where  $C_{0}=$const.$>0$, $m=$const., $b=$const..   In  (\ref{ww84}),  introduce  new independent variables
$$  \eta=yx^{\frac{m-n}{n+1}}\left(\frac{r\rho C_{0}^{2-n}}{n(n+1)k}\right)^{\frac{1}{n+1}}, \ \
\psi(x,y)=C_{0}x^{\frac{r}{n+1}} \left(\frac{n(n+1)k}{r\rho C_{0}^{2-n}}\right)^{\frac{1}{n+1}} f(\eta),  $$
and take the flow function of the form
$$u=\partial\psi/ \partial y, \ \  v=-\partial\psi/ \partial x,  $$
then problem (\ref{ww86})-(\ref{ww87}) reduces to the following problem
\begin{eqnarray}
|f''|^{n-1}f'''+f f''+\beta(1-f'^{2})+N(1-f')=0,
 \label{ww88}
\end{eqnarray}
with the boundary conditions
\begin{eqnarray}
f(0)=0, \ \  f'(0)=0, \ \  f'(\eta)\rightarrow 1\ \ \text{as} \ \  \eta\rightarrow \infty,
 \label{ww89}
\end{eqnarray}
where $\beta=m(n+1)/r, N=(n+1)\sigma b^{2}/(C_{0}\rho r)$.

\begin{Theorem} \label{OS.9.57} (\cite{OS})
Let $f(\eta)$ be a solution of problem (\ref{ww88})-(\ref{ww89}) such that $f'(\eta)\rightarrow C$ as $\eta\rightarrow\infty$.
 Assume that $\eta^{-1}(f'')^{n}\rightarrow 0$ as $\eta \rightarrow \infty$. Then $C=1$.

\end{Theorem}

  Oleinik and  Samokhin (\cite{OS}) in 1999 considered
 the dilatable fluid in a transversed magnetic field system of equations
\begin{equation}\left\{
\begin{array}{ll}
\frac{\partial u }{\partial t}+  u \frac{\partial u }{\partial x}+v \frac{\partial u }{\partial y}
=  \nu \frac{\partial}{\partial y }\left(|\frac{\partial u }{\partial y}|^{n-1}\frac{\partial u }{\partial y} \right)+U(x)\frac{d U }{d x}
+d^{2} (U-u) +\frac{\partial U }{\partial t}, \ \  n>1, \\
\frac{\partial u }{\partial x}+\frac{\partial v }{\partial y} =0
\end{array}
 \label{ww90}         \right.\end{equation}
 in the domain $D =\{  0<x<X, 0<y< \infty \}$, with the boundary conditions
 \begin{equation}\left\{
\begin{array}{ll}
  u ( 0,y)=0, \ \ u (x,0)=0, \ \ v (x,0)=v_{0}( x), \\
u ( x,y)\rightarrow U ( x),  \ \ \text{as}\ \   y\rightarrow \infty.
\end{array}
 \label{ww91}  \right.       \end{equation}

\begin{Theorem} \label{OS.9.58} (\cite{OS})
Suppose that $U'(x), U''(x), v_{0 }'(x), d^{2}_{x}$ are continuous functions and $U(0)=0$, $U'(0)>0$,
$$U'(x)+d^{2}(x)/2\geq a_{0}>0, \ \  -M_{0}x^{2n-1}\leq v_{0}(x)\leq 0, $$
for $0\leq x\leq X$, where $M_{0}=$const.$>0$. Then there exist $u(x,y), v(x,y)$ with the following properties: $u,v$ have locally square summable derivatives
involved in (\ref{ww90}); equations (\ref{ww90}) hold almost everywhere; $u$ is continuous in $\overline{D}$;  $v$ is square summable on any compact set inside $D$;
$u(x,0)=0, u(0,y)=0, u(x,y)\rightarrow U(x)$ as $y\rightarrow\infty$; $u(x,y)>0$ for $x>0$ and $y>0$;
$0\leq u_{y}\leq M_{1}, M_{1}=$const.$>0$; $u_{y}>0$ for $y=0$ and $x>0$; $u_{x}, v$ are continuous in $y$ for $y=0$; $v(x,0)=v_{0}(x)$ for $0<x<X$;
$u(x,y)\equiv U(x)$ for $y\geq C_{0}$, where $C_{0}>0$ is a constant depending on $X$ and the other data of the problem.
If $d^{2}(x)/U^{2n+1}(x)\geq a_{1}>0$, then $C_{0}$ can be chosen independent of $X$. The solution of problem (\ref{ww90})-(\ref{ww91}) with these properties
is unique.

\end{Theorem}

  Oleinik and Samokhin (\cite{OS}) in 1999 considered
  the stationary MHD boundary layer in a rapidly oscillating magnetic field, the system reads as
\begin{equation}\left\{
\begin{array}{ll}
   u \frac{\partial u }{\partial x}+v \frac{\partial u }{\partial y}
=  \nu  \frac{\partial^{2} u }{\partial y^{2}} +d^{2}(x)(U(x)-u) +U(x)\frac{d U }{d x}, \\
\frac{\partial u }{\partial x}+\frac{\partial v }{\partial y} =0,
\end{array}
 \label{ww94}         \right.\end{equation}
 in the domain $D =\{  0<x<X, 0<y< \infty \}$, with the boundary conditions
 \begin{equation}\left\{
\begin{array}{ll}
  u ( 0,y)=u_{0}(y), \ \ u (x,0)=0, \ \ v (x,0)=v_{0}( x), \\
u ( x,y)\rightarrow U ( x),  \ \ \text{as}\ \   y\rightarrow \infty.
\end{array}
 \label{ww95}    \right.     \end{equation}

\begin{Theorem} \label{OS.10.60} (\cite{OS})
Suppose that $U'(x)\geq 0$ for $0\leq x\leq X$; $U(x), u_{0}(y), v_{0}(x), d^{2}(x,\varepsilon^{-1}x)$ satisfy the compatibility conditions for any
$\varepsilon\in (0,\varepsilon_{0}]$. Then for $d(x)=d(x,\varepsilon^{-1}x)$ and $\varepsilon \rightarrow 0$, the solution of problem   (\ref{ww94})-(\ref{ww95})
 in the sense that
\begin{eqnarray*}
&& u_{\varepsilon}\rightarrow \overline{u} \ \ \text{uniformly for} \ \ 0\leq x\leq X, \ \  0\leq y\leq y_{1}<\infty, \\
&& v_{\varepsilon}\rightarrow \overline{v} \ \ \text{uniformly for} \ \ x_{0}\leq x\leq X, \ \  0\leq y\leq y_{1} , \\
&& \text{grad}~u_{\varepsilon}\rightarrow \text{grad}~\overline{u} \ \ \text{weakly in} \ \  L_{2}(\{0\leq x\leq X, 0\leq y\leq y_{1}\}) ,
\end{eqnarray*}
where $y_{1}>0$ and $x_{0}\in (0,X)$ are arbitrary.

\end{Theorem}

Gao, Guo, Huang  (\cite{GGH}) in 2017 studied  the local existence and uniqueness of solutions to the 2D MHD boundary layer where the initial tangential magnetic field has a lower bound $\delta>0$ plays an important role.   They studied the following initial boundary value problem for the nonlinear MHD boundary layer equations
\begin{eqnarray}
\left\{\begin{array}{l}
{\partial_{t} u_{1}+u_{1} \partial_{x} u_{1}+u_{2} \partial_{y} u_{1} -h_{1} \partial_{x} h_{1}-h_{2} \partial_{y} h_{1}=\mu \partial_{y}((\theta+1)\partial_{y}u_{1}) -P_{x},} \\
{\partial_{t} \theta+u_{1} \partial_{x} \theta+u_{2} \partial_{y} \theta_{1}-\kappa\partial_y^2\theta=\mu(\theta+1)(\partial_{y}u_{1})^2 +\nu(\partial_yh_1)^2,}\\
{\partial_{x} u_{1}+\partial_{y} u_{2}=0,\ \ \partial_{x} h_{1}+\partial_{y} h_{2}=0, }\\
{u_{1}|_{t=0}=u_{10}(x,y),\ \ h_{1}|_{t=0}=h_{10}(x,y), }\\
{(u_{1},u_{2},\partial_{y}h_{1}, h_{2})|_{y=0}=\mathbf{0},\ \ \lim\limits_{y\rightarrow +\infty}}(u_{1},\theta, h_{1})=(U,\Theta,H)(t,x).
\end{array}\right. \label{38.1}
\end{eqnarray}
\begin{Theorem}(\cite{GGH})
 Let $ m\geq 5$ be an integer, and $l\geq 0$ a real number. Assume that the outer flow $(U,\theta,H,P _{x})(t,x) $ satisfies that for some $T>0$,
\begin{eqnarray}
M_{0}:=\sum\limits_{i=0}^{2m+2}\left(\sup\limits_{0\leq t\leq T}\|\partial_{t}^{i}(U,\theta,H,P)(t,\cdot)\|_{H^{2m+2-i}(\mathbb{T}_{x})}+
\|\partial_{t}^{i}(U,H,P) \|_{L^{2}(0,T;H^{2m+2-i}(\mathbb{T}_{x}))}  \right)<+\infty.\label{38.2}
\end{eqnarray}
Also, we suppose the initial data $(u_{10},\theta_{10},h_{10})(x,y)$ satisfies
\begin{eqnarray}
 \left( u_{10}(x,y)-U(0,x),\theta_{10}(x,y)-\Theta(0,x),  h_{10}(x,y)-H(0,x) \right)\in H^{3m+2}_{l}(\Omega),
\end{eqnarray}
and the compatibility conditions up to m-th order. Moreover, there exists a sufficiently small constant $\delta_{0}>0$ such that
\begin{eqnarray}
 \left|<y>^{l+1}\partial_{y}^{i} (u_{10} ,  h_{10})(x,y) \right|\leq (2\delta_{0})^{-1},\ \  h_{10}(x,y)\geq 2\delta_{0}, \ \ for \ \ i=1,2,\ \ (x,y)\in \Omega.
\end{eqnarray}
Then, there exist a postive time $0<T_{*}\leq T$ and a unique solution $(u_{1},u_{2}, h_{1}, h_{2})$ to the initial boundary
value problem \eqref{38.1}, such that
\begin{eqnarray}
  (u_{1 }-U ,\theta-\Theta,  h_{1 } -H ) \in \bigcap\limits_{i=0}^{m} W^{i,\infty}(0,T_{*}; H^{m-i}_{l}(\Omega)),
\end{eqnarray}
and
\begin{eqnarray}
  (u_{2 }+U_{x}y ,  h_{2 } +H_{x}y ) \in \bigcap\limits_{i=0}^{m-1} W^{i,\infty}(0,T_{*}; H^{m-1-i}_{l}(\Omega)),\\
    (\partial_{y}u_{2 }+U_{x} ,\partial_{y}  h_{2 }+H_{x}) \in \bigcap\limits_{i=0}^{m-1} W^{i,\infty}(0,T_{*}; H^{m-i}_{l}(\Omega)).
\end{eqnarray}
Moreover, if $l>\frac{1}{2}$,
\begin{eqnarray}
  (u_{2 }+U_{x}y ,  h_{2 } +H_{x}y ) \in \bigcap\limits_{i=0}^{m-1} W^{i,\infty}(0,T_{*}; L^{\infty}(\mathbb{R}_{y,+};H^{m-1-i}(\mathbb{T}_{x}))),\\
    (\partial_{y}u_{2 }+U_{x} ,\partial_{y}  h_{2 }+H_{x}) \in \bigcap\limits_{i=0}^{m-1} W^{i,\infty}(0,T_{*}; H^{m-i}_{l}(\Omega)),
\end{eqnarray}
where $<y>= 1+y $, and
$$H^{m}_{l}(\Omega):=\left\{f(x,y):\Omega\rightarrow \mathbb{R},   \|f\|_{H^{m}_{l}(\Omega) } ^{2} =
\sum\limits_{m_{1}+m_{2}\leq m} \|<y>^{l+m_{2}}\partial_{x}^{m_{1}} \partial_{y}^{m_{2}}f \|_{L^{2}(\Omega)}^{2}<+\infty \right\}. $$
\end{Theorem}

Xie and Yang (\cite{xy}) in 2019 considered the lifespan of solution to the MHD boundary layer system as an analytic
perturbation of general shear flow. 
Since there is no restriction on the strength of the shear flow and the lifespan estimate is
larger than the one obtained for the classical Prandtl system in this setting, it reveals the stabilizing effect of
the magnetic field on the electrically conducting fluid near the boundary.

Consider the high Reynolds number limit to the MHD system near a no-slip boundary, when both of the Reynolds number and the magnetic Reynolds number have the same order in two space dimensions. Precisely, they considered the MHD system in the domain $\left\{(x, Y) \mid x \in \mathbb{R}, Y \in \mathbb{R}_{+}\right\}$ with $Y=0$ being the boundary,
\begin{equation}
\left\{\begin{array}{l}
\partial_{t} u^{\varepsilon}+\left(u^{\varepsilon} \partial_{x}+v^{\varepsilon} \partial_{Y}\right) u^{\varepsilon}+\partial_{x} p^{\varepsilon}-\left(h^{\varepsilon} \partial_{x}+g^{\varepsilon} \partial_{Y}\right) h^{\varepsilon}=\varepsilon\left(\partial_{x}^{2} u^{\varepsilon}+\partial_{Y}^{2} u^{\varepsilon}\right) ,\\
\partial_{t} v^{\varepsilon}+\left(u^{\varepsilon} \partial_{x}+v^{\varepsilon} \partial_{Y}\right) v^{\varepsilon}+\partial_{Y} p^{\varepsilon}-\left(h^{\varepsilon} \partial_{x}+g^{\varepsilon} \partial_{Y}\right) g^{\varepsilon}=\varepsilon\left(\partial_{x}^{2} v^{\varepsilon}+\partial_{Y}^{2} v^{\varepsilon}\right) ,\\
\partial_{t} h^{\varepsilon}+\left(u^{\varepsilon} \partial_{x}+v^{\varepsilon} \partial_{Y}\right) h^{\varepsilon}-\left(h^{\varepsilon} \partial_{x}+g^{\varepsilon} \partial_{Y}\right) u^{\varepsilon}=\kappa \varepsilon\left(\partial_{x}^{2} h^{\varepsilon}+\partial_{Y}^{2} h^{\varepsilon}\right), \\
\partial_{t} g^{\varepsilon}+\left(u^{\varepsilon} \partial_{x}+v^{\varepsilon} \partial_{Y}\right) g^{\varepsilon}-\left(h^{\varepsilon} \partial_{x}+g^{\varepsilon} \partial_{Y}\right) v^{\varepsilon}=\kappa \varepsilon\left(\partial_{x}^{2} g^{\varepsilon}+\partial_{Y}^{2} g^{\varepsilon}\right), \\
\partial_{x} u^{\varepsilon}+\partial_{Y} v^{\varepsilon}=0, \quad \partial_{x} h^{\varepsilon}+\partial_{Y} g^{\varepsilon}=0,
\end{array}\right. \label{y2.8}
\end{equation}
where both the viscosity and resistivity coefficients are denoted by a small positive parameter $\varepsilon,\left(u^{\varepsilon}, v^{\varepsilon}\right)$ and $\left(h^{\varepsilon}, g^{\varepsilon}\right)$ represent the velocity and the magnetic field respectively. The no-slip boundary condition is imposed on the velocity field
\begin{equation}
\left.\left(u^{\varepsilon}, v^{\varepsilon}\right)\right|_{Y=0}=\mathbf{0}. \label{yy2.9}
\end{equation}
Note that the system (\ref{y2.8})-(\ref{yy2.9}) governs the fluid behavior in the leading order
of approximation near the boundary  derived in
\begin{equation}
\left\{\begin{array}{l}
\partial_{t} u_{1}+  u_{1} \partial_{x}u_{1} +u_{2}  \partial_{y}u_{1}=b_{1} \partial_{x}b_{1} +b_{2} \partial_{y}b_{1} +\partial_{y}^{2}u_{1} ,\\
\partial_{t} b_{1}+  \partial_{y}(u_{2}b_{1}-u_{1}b_{2})=k\partial_{y}^{2}b_{1} ,\\
 \partial_{x}u_{1} +  \partial_{y}u_{2}=0 ,\quad   \partial_{x}b_{1} +  \partial_{y}b_{2}=0 , \\
u_{1}(t,x,y)|_{t=0}=u_{0}(x,y) , \quad b_{1}(t,x,y)|_{t=0}=b_{0}(x,y) ,\\
\lim\limits_{y\rightarrow\infty} u_{1}=\bar{u}, \quad  \lim\limits_{y\rightarrow\infty} b_{1}=\bar{b}.
\end{array}\right. \label{yy2.8}
\end{equation}
Assume the shear flow $u^{s} (t,y)$ has the following properties:
$$(  H ) \quad   \quad
\| \partial_{y}^{i}u_{s}(t,\cdot) \|_{L_{y}^{\infty}}\leq \frac{C}{\langle t\rangle^{\frac{i}{2}}}, \quad
\int_{0}^{\infty}|\partial_{t}u_{s}(t,y) |dy\leq C,  \quad
\|\theta_{\alpha} \partial_{y}^{2}u_{s}(t,\cdot) \|_{L_{y}^{2}}\leq \frac{C}{\langle t\rangle^{\frac{3}{4}}},
 \quad \langle t\rangle=1+t  .$$
To define the function space of the solution considered in \cite{xy}, the following Gaussian weighted function $\theta_{\alpha}$ will be used:
$$
\theta_{\alpha}(t, y)=\exp \left(\frac{\alpha z(t, y)^{2}}{4}\right), \quad \text { with } z(t, y)=\frac{y}{\sqrt{\langle t\rangle}} \text { and } \alpha \in[1 / 4,1 / 2] .
$$
With this and
$$
M_{m}=\frac{\sqrt{m+1}}{m !},
$$
define the Sobolev weighted semi-norms by
\begin{eqnarray*}
&& X_{m}=X_{m}(f, \tau)=\left\|\theta_{\alpha} \partial_{x}^{m} f\right\|_{L^{2}} \tau^{m} M_{m}, \quad D_{m}=D_{m}(f, \tau)=\left\|\theta_{\alpha} \partial_{y} \partial_{x}^{m} f\right\|_{L^{2}} \tau^{m} M_{m},\\
 &&Y_{m}=Y_{m}(f, \tau)=\left\|\theta_{\alpha} \partial_{x}^{m} f\right\|_{L^{2}} \tau^{m-1} m M_{m} .
 \end{eqnarray*}

Then the following space of analytic functions in the tangential variable $x$ and weighted Sobolev in the normal variable $y$ is defined by
$$
X_{\tau, \alpha}=\left\{f(t, x, y) \in L^{2}\left(\mathbb{H} ; \theta_{\alpha} d x d y\right):\|f\|_{X_{\tau, \alpha}}<\infty\right\},
$$
with $\tau>0$ and the norm
$$
\|f\|_{X_{\tau, \alpha}}=\sum_{m \geq 0} X_{m}(f, \tau) .
$$
In addition, the following two semi-norms will also be used:
$$
\|f\|_{D_{\tau, \alpha}}=\sum_{m \geq 0} D_{m}(f, \tau)=\left\|\partial_{y} f\right\|_{X_{\tau, \alpha}}, \quad\|f\|_{Y_{\tau, \alpha}}=\sum_{m \geq 1} Y_{m}(f, \tau).
$$

\begin{Theorem}(\cite{xy})
For any $\lambda \in[3 / 2,2),$ there exists a small positive constant $\varepsilon_{*}$ depending on $2-\lambda$. Under the assumption $(H)$ on the backgroud shear flow $\left(u_{s}(t, y), 0, \bar{b}, 0\right)$ with $\bar{b} \neq 0$, assume the initial data $u_{0}$ and $b_{0}$ satisfy
$$
\left\|u_{0}-u_{s}(0, y)\right\|_{X_{2 \tau_{0}, 1 / 2}} \leq \varepsilon, \quad\left\|b_{0}-\bar{b}\right\|_{X_{2 \tau_{0}, 1 / 2}} \leq \varepsilon,
$$
for some given $\varepsilon \in\left(0, \varepsilon_{*}\right] .$ Then there exists a unique solution $\left(u_{1}, u_{2}, b_{1}, b_{2}\right)$ to the $M H D$ boundary layer equations $(\ref{yy2.8})$ such that
$$
\left(u_{1}-u_{s}(t, y), b_{1}-\bar{b}\right) \in X_{\tau, \alpha}, \quad \alpha \in[1 / 4,1 / 2],
$$
with analyticity radius $\tau$ larger than $\tau_{0} / 4$ in the time interval $\left[0, T_{\varepsilon}\right] .$ And the lifespan $T_{\varepsilon}$ has the following low bound estimate
$$
T_{\varepsilon} \geq C \varepsilon^{-\lambda},
$$
where the constant $C>0$ is independent of $\varepsilon$.
\end{Theorem}

Huang, Liu and Yang (\cite{hly}) in 2019 were concerned with the motion of electrically conducting fluid governed by the 2D non-isentropic viscous compressible MHD system on the half plane with no-slip condition on
the velocity field, perfectly conducting wall condition on the magnetic field and Dirichlet boundary condition on the temperature on the boundary, and  obtained the local-in-time well-posedness of the following  boundary layer system (\ref{y2.9}) in
weighted Sobolev spaces:
\begin{equation}
\left\{\begin{array}{l}
\partial_{t} \rho^{e}+\nabla \cdot\left(\rho^{e} \mathbf{u}^{e}\right)=0, \\
\rho^{e}\left(\partial_{t} \mathbf{u}^{e}+\left(\mathbf{u}^{e} \cdot \nabla\right) \mathbf{u}^{e}\right)+\nabla p\left(\rho^{e}, \theta^{e}\right)-\left(\nabla \times \mathbf{H}^{e}\right) \times \mathbf{H}^{e}=0 ,\\
c_{V} \rho^{e}\left(\partial_{t} \theta^{e}+\left(\mathbf{u}^{e} \cdot \nabla\right) \theta^{e}\right)+p\left(\rho^{e}, \theta^{e}\right) \nabla \cdot \mathbf{u}^{e}=0, \\
\partial_{t} \mathbf{H}^{e}-\nabla \times\left(u^{e} \times \mathbf{H}^{e}\right)=0, \quad \nabla \cdot \mathbf{H}^{e}=0,
\end{array}\right. \label{y2.9}
\end{equation}
with the following  natural and sufficient boundary conditions:
$$
\left.\left(u_{2}^{e}, h_{2}^{e}\right)\right|_{y=0}=\mathbf{0}.
$$
\begin{Theorem}(\cite{hly})
For the initial-boundary value problem  with a smooth outflow $(U, \Theta, H, P)(t, x),$ assume that the initial data $\left(u_{1,0}, \theta_{0}, h_{1,0}\right)(x, y)$ and the boundary value $\theta^{*}(t, x)$ are smooth, compatible and satisfy
$$
\theta^{*}(t, x), \theta_{0}(x, y), h_{1,0}(x, y) \geq 2 \delta, \quad \frac{1}{2}\left(h_{1,0}(x, y)\right)^{2} \leq P(0, x)-2 \delta,
$$
for $t \in[0, T],(x, y) \in \mathbb{T} \times \mathbb{R}_{+}$ with some constant $\delta>0 .$ Then there exists a $T_{*} \in(0, T]$ and a unique classical solution $\left(u_{1}, u_{2}, \theta, h_{1}, h_{2}\right)$ to the problem (\ref{y2.9}) in the domain \\ $D_{T_{*}}:=\left\{(t, x, y) \mid t \in\left[0, T_{*}\right], x \in \mathbb{T}, y \in \mathbb{R}_{+}\right\}$ with the following properties:\\
(1). $\theta(t, x, y), h_{1}(t, x, y) \geq \delta$ and $\left(h_{1}(x, y)\right)^{2} / 2 \leq P(t, x)-\delta$ in $D_{T_{*}}$;\\
(2). $\left(u_{1}, \theta, h_{1}\right), \partial_{y}\left(u_{1}, \theta, h_{1}\right)$ and $\partial_{y}^{2}\left(u_{1}, \theta, h_{1}\right)$ are continuous and bounded in $D_{T_{*}}, \partial_{t}\left(u_{1}, \theta, h_{1}\right)$
and $\partial_{x}\left(u_{1}, \theta, h_{1}\right)$ are continuous and bounded in any compact set of $D_{T_{*}}$;\\
(3). $\left(u_{2}, h_{2}\right)$ and $\partial_{y}\left(u_{2}, h_{2}\right)$ are continuous and bounded in any compact set of $D_{T_{*}}$.
\end{Theorem}

Liu, Xie and Yang (\cite{lxy1}) in 2019 studied the well-posedness theory for the MHD boundary layer. The boundary
layer equations are governed by the Prandtl-type equations that are derived from
the incompressible MHD system with non-slip boundary condition on the velocity and perfectly conducting condition on the magnetic field. Under the assumption that the initial tangential magnetic field is not zero, they established the local-in time existence, uniqueness of solutions for the nonlinear MHD boundary layer
equations.

One important problem about magnetohydrodynamics (MHD) is to understand the high Reynolds number limits in a domain with boundary. To this end, they considered the following initial boundary value problem for the 2D viscous MHD equations in a periodic domain $\{(t, x, y): t \in$ $\left.[0, T], x \in \mathbb{T}, y \in \mathbb{R}_{+}\right\}:$
\begin{equation}
\left\{\begin{array}{l}
\partial_{t} \mathbf{u}^{\epsilon}+\left(\mathbf{u}^{\epsilon} \cdot \nabla\right) \mathbf{u}^{\epsilon}-\left(\mathbf{H}^{\epsilon} \cdot \nabla\right) \mathbf{H}^{\epsilon}+\nabla p^{\epsilon}=\mu \epsilon \Delta \mathbf{u}^{\epsilon}, \\
\partial_{t} \mathbf{H}^{\epsilon}-\nabla \times\left(\mathbf{u}^{\epsilon} \times \mathbf{H}^{\epsilon}\right)=\kappa \in \Delta \mathbf{H}^{\epsilon}, \\
\nabla \cdot \mathbf{u}^{\epsilon}=0, \quad \nabla \cdot \mathbf{H}^{\epsilon}=0.
\end{array}\right. \label{L1}
\end{equation}
Here,   assume the viscosity and resistivity coefficients have the same order of a small parameter $\epsilon $.  $\mathbf{u}^{\epsilon}=\left(u_{1}^{\epsilon}, u_{2}^{\epsilon}\right)$ denotes the velocity vector, $\mathbf{H}^{\epsilon}=\left(h_{1}^{\epsilon}, h_{2}^{\epsilon}\right)$
denotes the magnetic field, and $p^{\epsilon}=\tilde{p}^{\epsilon}+\left|\mathbf{H}^{\epsilon}\right|^{2} / 2$ denotes the total pressure with $\tilde{p}^{\epsilon}$ the pressure of the fluid. On the boundary, the non-slip boundary condition is imposed on velocity field
$$
\left.\mathbf{u}^{\epsilon}\right|_{y=0}=\mathbf{0},
$$
and the perfectly conducting boundary condition on magnetic field
$$
\left.h_{2}^{\epsilon}\right|_{y=0}=\left.\partial_{y} h_{1}^{\epsilon}\right|_{y=0}=0.
$$

Denote
$$
\Omega:=\{(x, y): x \in \mathbb{T}, y \in \mathbb{R}+\}.
$$
For any $l \in \mathbb{R},$ denote by $L_{l}^{2}(\Omega)$ the weighted Lebesgue space with respect to the spatial variables:
$$
\begin{aligned}
L_{l}^{2}(\Omega):=\left\{  f(x, y): \Omega \rightarrow \mathbb{R},
 \|f\|_{L_{l}^{2}(\Omega)}:=\left(\int_{\Omega}\langle y\rangle^{2 l}|f(x, y)|^{2} dx d y\right)^{\frac{1}{2}}<+\infty\right\},
\end{aligned}
$$
with $\langle y\rangle=1+y .$ Then, for any given $m \in \mathbb{N},$ denote by $H_{l}^{m}(\Omega)$, the weighted Sobolev spaces,
$$
H_{l}^{m}(\Omega):=\left\{f(x, y): \Omega \rightarrow \mathbb{R},\|f\|_{H_{l}^{m}(\Omega)}<+\infty\right\},
$$
with the norm
$$
\|f\|_{H_{l}^{m}(\Omega)}:=\left(\sum_{m_{1}+m_{2} \leq m}\left\|\langle y\rangle^{l+m_{2}} \partial_{x}^{m_{1}} \partial_{y}^{m_{2}} f\right\|_{L^{2}(\Omega)}^{2}\right)^{\frac{1}{2}}.
$$
\begin{Theorem}(\cite{lxy1})
Let $m \geq 5$ be an integer and $l \geq 0$ a real number. Assume that the outer flow $\left(U, H, P_{x}\right)(t, x)$ satisfies that for some $T>0,$
$$
M_{0}:=\sum_{i=0}^{2 m+2}\left(\sup _{0 \leq t \leq T}\left\|\partial_{t}^{i}(U, H, P)(t, \cdot)\right\|_{H^{2 m+2-i}\left(\mathbb{T}_{x}\right)}\right)<+\infty.
$$
Also,   assume the initial data $\left(u_{10}, h_{10}\right)(x, y)$ satisfies
$$
\left(u_{10}(x, y)-U(0, x), h_{10}(x, y)-H(0, x)\right) \in H_{l}^{3 m+2}(\Omega),
$$
and the compatibility conditions up to $m$-th order. Moreover, there exists a sufficiently small constant $\delta_{0}>0$ such that
$\left|\langle y\rangle^{l+1} \partial_{y}^{i}\left(u_{10}, h_{10}\right)(x, y)\right| \leq\left(2 \delta_{0}\right)^{-1}  $ for $i=1,2,(x, y) \in \Omega$,  $h_{10}(x, y) \geq 2 \delta_{0}$.
Then, there exist a positive time $0<T_{*} \leq T$ and a unique solution $\left(u_{1}, u_{2}, h_{1}, h_{2}\right)$ to the initial boundary value problem (\ref{L1}) such that
$$
\left(u_{1}-U, h_{1}-H\right) \in \bigcap_{i=0}^{m} W^{i, \infty}\left(0, T_{*} ; H_{l}^{m-i}(\Omega)\right),
$$
and
$$
\begin{aligned}
\left(u_{2}+U_{x} y, h_{2}+H_{x} y\right) & \in \bigcap_{i=0}^{m-1} W^{i, \infty}\left(0, T_{*} ; H_{-1}^{m-1-i}(\Omega)\right) ,\\
\left(\partial_{y} u_{2}+U_{x}, \partial_{y} h_{2}+H_{x}\right) & \in \bigcap_{i=0}^{m-1} W^{i, \infty}\left(0, T_{*} ; H_{l}^{m-1-i}(\Omega)\right) .
\end{aligned}
$$
Moreover, if $l>\frac{1}{2},$
$$
\left(u_{2}+U_{x} y, h_{2}+H_{x} y\right) \in \bigcap_{i=0}^{m-1} W^{i, \infty}\left(0, T_{*} ; L^{\infty}\left(\mathbb{R}_{y,+} ; H^{m-1-i}\left(\mathbb{T}_{x}\right)\right)\right).
$$
\end{Theorem}

Chen, Ren, Wang and Zhang (\cite{CRWZ1}) in 2020 studied the long time well-posedness of 2D MHD boundary layer equation, and proved that if the initial data satisfies
$$
\left\|\left(u_{0}, h_{0}-1\right)\right\|_{H_{\mu}^{3,0} \cap H_{\mu}^{1,2}} \leq \varepsilon,
$$
then the life span of the solution is at least of order $\varepsilon^{2-\eta}$ for $\eta>0$.

Let $h(t, x, y)=1+\widetilde{h}(t, x, y) .$ Then $(u, \widetilde{h})$ satisfies the following system
\begin{equation}
\left\{\begin{array}{l}
\partial_{t} u+u \partial_{x} u+v \partial_{y} u-h \partial_{x} \tilde{h}-g \partial_{y} \tilde{h}-\partial_{y}^{2} u=0, \\
\partial_{t} \widetilde{h}+u \partial_{x} \widetilde{h}+v \partial_{y} \widetilde{h}-h \partial_{x} u-g \partial_{y} u-\partial_{y}^{2} \widetilde{h}=0 ,\\
\partial_{x} u+\partial_{y} v=0, \quad \partial_{x} \widetilde{h}+\partial_{y} g=0 ,\\
\left.\left(u, v, \partial_{y} \widetilde{h}, g\right)\right|_{y=0}=0 \quad \text { and } \quad \lim\limits _{y \rightarrow+\infty}(u, \widetilde{h})=(0,0), \\
\left.(u, \widetilde{h})\right|_{t=0}=\left(u_{0}, \widetilde{h}_{0}\right) .
\end{array}\right.
\end{equation}
Now  introduce the following weighted Sobolev space. For $k, \ell \in \mathbb{N},$ the space $H_{\omega}^{k, \ell}\left(\mathbb{R}_{+}^{2}\right)$ consists of all functions $f \in L_{\omega}^{2}$ satisfying
$$
\|f\|_{H_{\omega}^{k, \ell}}^{2} \stackrel{\text { def }}{=} \sum_{\alpha=0}^{k} \sum_{\beta=0}^{\ell}\left\|\partial_{x}^{\alpha} \partial_{y}^{\beta} f\right\|_{L_{\omega}^{2}}^{2}<+\infty,
$$
where $\|f\|_{L_{\omega}^{p}}=\|\omega(y) f(x, y)\|_{L^{p}}$ with $\omega(y)$ a positive weight function.

\begin{Theorem}(\cite{CRWZ1})
Let $\mu=\exp \left(\frac{1+y^{2}}{8(t)}\right)$ with $\langle t\rangle=1+t .$ For any $\eta \in(0,1),$ there exists an $\varepsilon>0$ such that if the initial data $\left(u_{0}, \widetilde{h}_{0}\right)$ satisfies
$$
\left\|\left(u_{0}, \widetilde{h}_{0}\right)\right\|_{H_{\mu}^{3,0}}+\left\|\left(u_{0}, \widetilde{h}_{0}\right)\right\|_{H_{\mu}^{1,2}} \leq \varepsilon,
$$
then there exists a time $T_{\varepsilon} \geq \varepsilon^{-(2-\eta)}$ such that the system has a solution $(u, \widetilde{h})$ on $\left[0, T_{\varepsilon}\right]$, which satisfies
$$
(u, \widetilde{h}) \in L^{\infty}\left(\left[0, T_{\varepsilon}\right] ; H_{\mu}^{3,0}\left(\mathbb{R}_{+}^{2}\right) \cap H_{\mu}^{1,2}\left(\mathbb{R}_{+}^{2}\right)\right) \cap L^{2}\left(\left[0, T_{\varepsilon}\right] ; H_{\mu}^{3,1}\left(\mathbb{R}_{+}^{2}\right) \cap H_{\mu}^{1,3}\left(\mathbb{R}_{+}^{2}\right)\right).
$$
\end{Theorem}

Liu, Wang, Xie and Yang (\cite{LWXY}) in 2020 were concerned with the magnetic effect on
the Sobolev solvability of boundary layer equations for the 2D
incompressible MHD system without resistivity. The MHD
boundary layer is described by the Prandtl type equations
derived from the incompressible viscous MHD system without
resistivity under the no-slip boundary condition on the
velocity. Assuming that the initial tangential magnetic field
does not degenerate, a local-in-time well-posedness in Sobolev
spaces is proved without the monotonicity condition on
the velocity field. Moreover, the authors of  \cite{LWXY} showed that if the tangential
magnetic field of shear layer is degenerate at one point, then
the linearized MHD boundary layer system around the shear
layer profile is ill-posed in the Gevrey function space provided.

Consider the initial-boundary value problem for the following  2D  MHD  boundary layer equations in the domain $\left\{(t, x, y): t \in[0, T], x \in \mathbb{T}, y \in \mathbb{R}_{+}\right\}:$
\begin{equation}
\left\{\begin{array}{l}
\partial_{t} u+u \partial_{x} u+v \partial_{y} u+\partial_{x} p-\partial_{y}^{2} u-b_{1} \partial_{x} b_{1}-b_{2} \partial_{y} b_{1}=0, \\
\partial_{t} b_{1}+u \partial_{x} b_{1}+v \partial_{y} b_{1}-b_{1} \partial_{x} u-b_{2} \partial_{y} u=0 ,\\
\partial_{x} u+\partial_{y} v=0, \quad \partial_{x} b_{1}+\partial_{y} b_{2}=0, \\
\left.\left(u, v, b_{2}\right)\right|_{y=0}=0, \quad \lim\limits _{y \rightarrow+\infty}\left(u, b_{1}\right)(t, x, y)=(U, B)(t, x) ,\\
\left.\left(u, b_{1}\right)\right|_{t=0}=\left(u_{0}, b_{0}\right)(x, y),
\end{array}\right. \label{y3.10}
\end{equation}
where $\mathbb{T}$ stands for a torus or a periodic domain, $\mathbb{R}_{+}=[0,+\infty),(u, v)$ and $\left(b_{1}, b_{2}\right)$ are the velocity and magnetic boundary layer functions respectively, and the known functions $U, B$ and $p$ satisfy the Bernoulli law:
$$
\left\{\begin{array}{l}
\partial_{t} U+U \partial_{x} U+\partial_{x} p=B \partial_{x} B, \\
\partial_{t} B+U \partial_{x} B=B \partial_{x} U.
\end{array}\right.
$$

\begin{Theorem}(\cite{LWXY})
Suppose that the outflow $(U, p, B)(t, x)$ in the equation is smooth, and the initial data and boundary conditions are smooth, compatible and satisfy
$$
b_{0}(x, y) \geq \delta_{0},
$$
for some positive constant $\delta_{0}$. Then there exists a time $T$, such that the initial-boundary value problem (\ref{y3.10}) admits a unique smooth solution $\left(u, v, b_{1}, b_{2}\right)(t, x, y)$ satisfying
$$
b_{1}(t, x, y) \geq \delta_{0} / 2,
$$
for all $t \in[0, T],(x, y) \in \mathbb{T} \times \mathbb{R}_{+}$.

\end{Theorem}

Li and Yang (\cite{LY}) in 2020 established the well-posedness of the MHD boundary layer system in Gevrey function space
without any structural assumption. Compared to the classical Prandtl equation, the loss of tangential
derivative comes from both the velocity and magnetic fields that are coupled with each other. By observing
a new type of cancellation mechanism in the system for overcoming the loss derivative degeneracy, they showed
that the MHD boundary layer system is well-posed with Gevrey index up to 3/2 in both two and three
dimensional spaces.
Consider the system
\begin{equation}
\left\{\begin{array}{l}
\left(\partial_{t}+\vec{u} \cdot \nabla-\nu \partial_{z}^{2}\right) u_{h}-(\vec{f} \cdot \nabla) f_{h}=0, \\
\partial_{t} \vec{f}-\nabla \times(\vec{u} \times \vec{f})-\mu \partial_{z}^{2} \vec{f}=0 ,\\
\operatorname{div} \vec{u}=\operatorname{div} \vec{f}=0, \\
\left.\vec{u}\right|_{z=0}=\left.\left(\partial_{z} f_{h}, f_{z}\right)\right|_{z=0}=\mathbf{0},\left.\quad\left(u_{h}, f_{h}\right)\right|_{z \rightarrow+\infty}=\mathbf{0}, \\
\left.u_{h}\right|_{t=0}=u_{h, 0},\left.\quad f_{h}\right|_{t=0}=f_{h, 0}.
\end{array}\right. \label{L2}
\end{equation}
Let $\ell \geq 1$ be a given number. With a given integer $N \geq 0$ and a pair $(\rho, \sigma), \rho>0$ and $\sigma \geq 1,$ a Banach space $X_{\rho, \sigma, N}$ consists of all smooth vector-valued functions $\boldsymbol{A}=\boldsymbol{A}\left(x_{h}, z\right)$ with $\left(x_{h}, z\right) \in \mathbb{R}_{+}^{d}$ such that the Gevrey norm $\|\boldsymbol{A}\|_{\rho, \sigma, N}<+\infty,$ where $\|\cdot\|_{\rho, \sigma, N}$ is defined below. Denote $\partial_{x_{h}}^{\alpha}=\partial_{x_{1}}^{\alpha_{1}} \cdots \partial_{x_{d-1}}^{\alpha_{d-1}}$
and define
$$
\|\boldsymbol{A}\|_{\rho, \sigma, N}=\sup _{0 \leq j \leq N \atop|\alpha|+j \geq 7} \frac{\rho^{|\alpha|+j-7}}{[(|\alpha|+j-7) !]^{\sigma}}\left\|\langle z\rangle^{\ell+j} \partial_{x_{h}}^{\alpha} \partial_{z}^{j} \boldsymbol{A}\right\|_{L^{2}\left(\mathbb{R}_{+}^{d}\right)}+\sup _{0 \leq j \leq N \atop|\alpha|+j \leq 6}\left\|\langle z\rangle^{\ell+j} \partial_{x_{h}}^{\alpha} \partial_{z}^{j} \boldsymbol{A}\right\|_{L^{2}\left(\mathbb{R}_{+}^{d}\right)},
$$
where $\langle z\rangle=\left(1+|z|^{2}\right)^{1 / 2}$ and
$$
\|\boldsymbol{A}\|_{L^{2}\left(\mathbb{R}_{+}^{d}\right)} \stackrel{\text { def }}{=}\left(\sum_{1 \leq j \leq k}\left\|A_{j}\right\|_{L^{2}\left(\mathbb{R}_{+}^{d}\right)}^{2}\right)^{1 / 2},
$$
for $\boldsymbol{A}=\left(A_{1}, \cdots, A_{k}\right)$. Here, $\sigma$ is the Gevrey index.

\begin{Theorem}(\cite{LY})
Let the dimension $d=2$ or 3. Suppose the initial data $\left(u_{h, 0}, f_{h, 0}\right)$ in the system (\ref{L2}) belong to $X_{2 \rho_{0}, \sigma, 8}$ for some $1<\sigma \leq 3 / 2$ and some $0<\rho_{0} \leq 1,$ compatible with the boundary condition. Then the system (\ref{L2}) admits a unique solution $\left(u_{h}, f_{h}\right) \in L^{\infty}\left([0, T] ; X_{\rho, \sigma, 4}\right)$ for some $T>0$ and some $0<\rho<2 \rho_{0}$.
\end{Theorem}

\subsection{MHD Boundary Layer Equations-Global Existence  }
In this subsection, we survey the results on the 2D MHD boundary layer equations.

Xie and Yang (\cite{xy1}) in 2018 proved global existence of solutions with analytic regularity to
the 2D MHD boundary layer equations in the mixed Prandtl and Hartmann regime.
The following mixed Prandtl and Hartmann boundary layer equations from the classical incompressible MHD system,
\begin{eqnarray}
\left\{\begin{array}{l}{\partial_{t} u_{1}+u_{1} \partial_{x} u_{1}+u_{2} \partial_{y} u_{1}=\partial_{y} b_{1}+\partial_{y}^{2} u_{1}}, \\ {\partial_{y} u_{1}+\partial_{y}^{2} b_{1}=0} ,\\ {\partial_{x} u_{1}+\partial_{y} u_{2}=0, \quad x \in \mathbb{R}, y \in \mathbb{R}_{+}}, \\
u_{ 1} (t = 0,x,y) = u _{10} (x,y), \\
u _{1 }|_{ y=0} = 0,\quad  u_{ 2} | _{y=0} = 0,\\
\lim\limits_{y\rightarrow +\infty}u_{ 1} = \bar{u}, \quad \lim\limits_{y\rightarrow +\infty}b_{ 1} = \bar{b}.
\end{array}. \right.\label{118.1}
\end{eqnarray}
Here, $(u_1,u_2)$ denotes the velocity field of the boundary layer and $b_1$ is the corresponding tangential magnetic component.
Integrating the equation of $\eqref{118.1}_2$ over $y$ yields
\begin{eqnarray}
-u _{1 }(t,x,y) + \bar{u} = \partial_{y} b_{ 1} .  \label{118.2}
\end{eqnarray}
Thus, the equations \eqref{118.1} can be rewritten as
\begin{eqnarray}
\left\{\begin{array}{l}{\partial_{t} u_{1}+u_{1} \partial_{x} u_{1}+u_{2} \partial_{y} u_{1}=-u _{1 }(t,x,y) + \bar{u} +\partial_{y}^{2} u_{1}}, \\ {\partial_{x} u_{1}+\partial_{y}  u_{1}=0} . \end{array}\right.\label{118.3}
\end{eqnarray}
Recall that the classical Hartmann boundary layer is given by
\begin{eqnarray}
u_{ 1} = (1- e^{ -y} )\bar{u}, \quad u_{ 2} = 0,\label{118.4}
\end{eqnarray}
which is a steady solution to  the system \eqref{118.3}. Without loss of generality, set $\bar{u} = 1$ and denote the
perturbation by $(u,v)$:
\begin{eqnarray}
u_{ 1} = (1-e^{ -y} ) + u, \quad u_{ 2} = v. \label{118.5}
\end{eqnarray}
Obviously, $(u,v)$ satisfies
\begin{eqnarray}
\left\{\begin{array}{l}{\partial_{t} u+\left(1-e^{-y}+u\right) \partial_{x} u+v \partial_{y}\left(-e^{-y}+u\right)=-u+\partial_{y}^{2} u} ,\\ {\partial_{x} u+\partial_{y} v=0},\\
u _{0} (x,y) = u_{ 10} (x,y)- (1 - e^{-y} ),\\
u _{1 }|_{ y=0} = 0,\quad  u_{ 2} | _{y=0} = 0.
\end{array}\right. \label{118.6}
\end{eqnarray}

For some $r > 1$, denote an analytic weight $M_m$ by
$$M_{m}=\frac{(m + 1) r}{m!}.$$
Define
\begin{eqnarray}
\left\{\begin{array}{l}{ X_{m}  =\left\|e^{\alpha y} \partial_{x}^{m} g\right\|_{L^{2}\left(\mathbb{R}_{+}^{2}\right)} \tau^{m} M_{m},  \ \ Z_{m}=\left\|e^{\alpha y} \partial_{y} \partial_{x}^{m} g\right\|_{L^{2}\left(\mathbb{R}_{+}^{2}\right)} \tau^{m} M_{m},} \\{ Y_{m} =\left\|e^{\alpha y} \partial_{x}^{m} g\right\|_{L^{2}\left(\mathbb{R}_{+}^{2}\right)} \tau^{m-1 / 2} m^{1 / 2} M_{m}, \quad D_{m}=\left\|e^{\alpha y} \partial_{x}^{m} g\right\|_{L_{y}^{\infty} L_{x}^{2}} \tau^{m} M_{m},} \label{118.7}\end{array}\right.
\end{eqnarray}
and
\begin{eqnarray}
\|g\|_{X_{\tau, \alpha}^{\prime}}^{2}=\sum_{m \geq 0} X_{m}^{2}, \quad\|g\|_{Z_{\tau, \alpha}^{2}}^{2}=\sum_{m \geq 0} Z_{m}^{2},\\
 \|g\|_{Y_{\tau, \alpha}^{\prime}}^{2}=\sum_{m \geq 0} Y_{m}^{2}, \quad\|g\|_{D_{\tau, \alpha}^{2}}^{2}=\sum_{m \geq 0} D_{m}^{2}. \label{118.9}
 \end{eqnarray}
Here, $\tau$ denotes the analytic radius.
\begin{Theorem}(\cite{xy1})
 Let the initial data $u _{10} (x,y)$ be a small perturbation of the Hartmann profile
$(\bar{u}(1-e ^{-y }),0)$ satisfying the compatibility conditions and
\begin{eqnarray}
\|\partial_{ y} u _{10} + u_{10} - \bar{ u}\|_{X_{\tau_{0}, \alpha}^{r}}  \leq \delta_{0},\label{118.10}
\end{eqnarray}
with $r > 1, 0 < \alpha<\frac{\sqrt{2}}{2}$  for some small constant $\delta_{0}>0$. Then there exists a unique global-in-
time solution $(u_1,u_2 )$ to the problem \eqref{118.1} satisfying
\begin{eqnarray}
 \|g\|_{X_{\tau(t), \alpha}^{r}} \leq e^{-2\left(1-2 \alpha^{2}\right) t} \delta_{0}, \quad \text{with }\quad \tau(t)>\tau_{0} / 2,
\end{eqnarray}
 for all time $t\geq0$, where $g =\partial_{y} u _{1} + u _{1} - \bar{ u}$.
\end{Theorem}

Chen, Ren, Wang and Zhang (\cite{CRWZ}) in 2020 proved the global well-posedness of the 2D magnetic Prandtl model in the mixed Prandtl/Hartmann
regime when the initial data is a small perturbation of the Hartmann layer in Sobolev spaces.
They considered a 2D magnetic Prandtl model in $\mathbb{R}_{+}^{2}$ :
\begin{equation}
\left\{\begin{array}{l}
\partial_{t} u+u \partial_{x} u+v \partial_{y} u-\partial_{y}^{2} u=\partial_{y} b, \\
\partial_{y} u+\partial_{y}^{2} b=0, \\
\partial_{x} u+\partial_{y} v=0, \\
\left.u\right|_{y=0}=\left.v\right|_{y=0}=\left.b\right|_{y=0}=0, \\
\lim\limits _{y \rightarrow+\infty} u(t, x, y)=u^{\infty}, \quad \lim\limits _{y \rightarrow+\infty} b(t, x, y)=b^{\infty},
\end{array}\right.
\end{equation}
where $(u, v)$ denotes the velocity field and $b$ denotes the tangential magnetic field.

Thanks to the boundary conditions and $\partial_{y} u+\partial_{y}^{2} b=0$, we get
$$
\partial_{y} b=-\left(u-u^{\infty}\right)
$$
which, substituted into the first equation, gives  a Prandtl type equation with an extra damping term $\left(u-u^{\infty}\right)$ :
\begin{equation}
\left\{\begin{array}{l}
\partial_{t} u+u \partial_{x} u+v \partial_{y} u-\partial_{y}^{2} u+\left(u-u^{\infty}\right)=0, \\
\partial_{x} u+\partial_{y} v=0 ,\\
\left.u\right|_{y=0}=\left.v\right|_{y=0}=0 \quad \text { and } \quad \lim \limits_{y \rightarrow+\infty} u(t, x, y)=u^{\infty}.
\end{array}\right.
\end{equation}
Now introduce the following weighted Sobolev space. For $k, \ell \in  \mathbb{N},$ the space $H_{\omega}^{k, \ell}$ consists of all functions $f \in L_{\omega}^{2}$ satisfying
$$
\|f\|_{H_{\omega}^{k, \ell}}^{2} := \sum_{\alpha=0}^{k} \sum_{\beta=0}^{\ell}\left\|\partial_{x}^{\alpha} \partial_{y}^{\beta} f\right\|_{L_{\omega}^{2}}^{2}<+\infty,
$$
where $\|f\|_{L_{\omega}^{p}}=\|\omega(y) f(x, y)\|_{L^{p}}$ with $\omega(y)$ being a positive weight function.

The authors of  \cite{CRWZ} took the
global stability of the Hartmann layer $u^{s}=u^{\infty}(1-e^{-y}) $, without loss of generality, let $u^{\infty}=1$ and introduce the perturbation $\widetilde{u}=u-u^{s}$, which satisfies
\begin{equation}
\left\{\begin{array}{l}
\partial_{t} \widetilde{u}+u \partial_{x} \widetilde{u}+v \partial_{y} u-\partial_{y}^{2}\widetilde{u}+\widetilde{u}=0, \\
\partial_{x} \widetilde{u}+\partial_{y} v=0 ,\\
\left.\widetilde{u}\right|_{y=0}=\left.v\right|_{y=0}=0 \quad \text { and } \quad \lim \limits_{y \rightarrow+\infty} \widetilde{u}(t, x, y)=0,\\
\widetilde{u}|_{t=0}=\widetilde{u}_{0}(x,y).
\end{array}\right. \label{L3}
\end{equation}
\begin{Theorem}(\cite{CRWZ})
Assume that $\widetilde{u}_{0} \in H_{\mu}^{3,1} \cap H_{\mu}^{1,2}\left(\mu(y)=e^{\frac{y}{2}}\right)$ with $\left.\tilde{u}_{0}\right|_{y=0}=0$ satisfies
$$\left\|\tilde{u}_{0}\right\|_{H_{\mu}^{3,1}}+\left\|\tilde{u}_{0}\right\|_{H_{\mu}^{1,2}} \leqslant \varepsilon, \quad
\left|\partial_{x} \tilde{u}_{0}(x, y)\right|+\left|\partial_{x} \partial_{y} \tilde{u}_{0}(x, y)\right|+\left|\partial_{y}^{2} \tilde{u}_{0}(x, y)\right| \leqslant \varepsilon e^{-y} \quad \text { for }(x, y) \in \mathbb{R}_{+}^{2}.
$$

Then there exists $\varepsilon_{0}>0$, such that for any $\varepsilon \in\left(0, \varepsilon_{0}\right),$ the magnetic Prandtl model (\ref{L3}) has a unique global in time solution. Moreover, there holds that for any $(t, x, y) \in[0, \infty) \times \mathbb{R}_{+}^{2},$
$$
\begin{array}{l}
c e^{-y} \leqslant \partial_{y} u(t, x, y) \leqslant C e^{-y} ,\\
\left|\partial_{x} \widetilde{u}(t, x, y)\right|+\left|\partial_{x} \partial_{y} \tilde{u}(t, x, y)\right|+\left|\partial_{y}^{2} \tilde{u}(t, x, y)\right| \leqslant C \varepsilon e^{-y} ,\\
e^{\frac{t}{4}}\left(\|\widetilde{u}(t)\|_{H_{\mu}^{3,1}}^{2}+\|\tilde{u}(t)\|_{H_{\mu}^{1,2}}^{2}\right)+\int_{0}^{t} e^{\frac{s}{4}}\|\tilde{u}(s)\|_{H_{\mu}^{1,3}}^{2} d s \leqslant C \varepsilon^{2}.
\end{array}
$$
\end{Theorem}

  Gong and   Wang (\cite{GGW}) in 2021  obtained the global existence of a weak solution to the following  mixed Prandtl-Hartmann boundary layer problem (\ref{w.67}) under the Crocco transformation, and   the existence of a back-flow point of the Prandtl-Hartmann boundary layer when the initial velocity satisfies certain growth condition  in the domain $D_{T}=\{0\leq t <T, 0\leq x\leq X, y\geq 0\}$,
 \begin{equation}\label{w.67}
\left\{\begin{array}{ll}
 \partial_{t} u +u \partial_{x} u+v \partial_{y} u= \partial_{y}^{2} u+\sigma B_{2} \partial_{y} b_{1}-p_{x}, \\
 \sigma B_{2} \partial_{y} u+\partial_{y}^{2} b_{1}=0, \\
 \partial_{x} u+\partial_{y} v=0,\\
u|_{t=0}=u_{0}(x,y),  \ \  u|_{x=0}=u_{1}(t,y), \\
(u,v)|_{y=0}= (0,0), \ \
  \lim\limits_{y \rightarrow+\infty}\left(u, b_{1}\right)=\left(U(t,,x), B_{1}(t,x)\right) , \end{array}\right.
\end{equation}
where $(u,v)$ is the velocity field, and $b_{1} $ is the tangential magnetic component. Here $\sigma>0$ is the magnetic conductivity parameter, the given functions $U,p,\overrightarrow{B}(t,x)=(B_{1},B_{2})$ are the tangential velocity, hydrodynamic pressure and magnetic field of the outer flow, respectively, which obey the Bernoulli's  law   (\ref{1.1.3}).  The pressure is favourable for the classical Prandtl equation, a.e.,
 \begin{equation}
  \partial_{x} p\leq 0.
\label{w.71}
 \end{equation}
Assume  that
 \begin{equation}
 U(t,x)\geq 0, \ \ u_{0}(x,y)>0, \ \ u_{1}(x,y)>0, \ \ \forall 0\leq x\leq X, \ \  y\geq 0, \ \ t\geq 0,
\label{w.72}
 \end{equation}
and the initial tangential velocity $u_{0}(x,y)$ and the upstream velocity $u_{1}(t,y)$ satisfy the following monotonicity condition
 \begin{equation}
\partial_{y}u_{0}(x,y)>0, \ \ \partial_{y}u_{1}(x,y)>0, \ \ \forall 0\leq x\leq X, \ \  y\geq 0, \ \ t\geq 0,
\label{w.73}
 \end{equation}
and
  \begin{equation}
\sigma B_{2}^{2} U- \partial_{x} p\geq 0.
\label{w.68}
 \end{equation}
The Crocco transformation is
$$\tau=t,\ \ \xi=x,\ \ \eta=\frac{u(t,x,y)}{U(t,x)} , \ \ w(\tau,\xi,\eta)=\frac{u_{y}(t,x,y)}{U(t,x)} ,$$
and denoting by $\alpha=\sigma B_{2}^{2}(t,x)$ and $\Omega=\{0\leq \xi\leq X, 0\leq\eta<1\}$, the original problem (\ref{w.67}) can be
transformed into the following initial-boundary value problem in $Q_{T}=[0,T)\times \Omega$,
 \begin{equation}\label{w.69}
\left\{\begin{array}{ll}
 w_{\tau}+\eta Uw_{\xi}+Aw_{\eta}+B w=w^{2}w_{\eta\eta}, \\
w|_{\tau=0}=w_{0} ,  \ \  w|_{\xi=0}=w_{1} , \\
w|_{\eta=1}= 0, \ \  (w_{\eta}+\frac{C}{w})|_{\eta=0}= 0, \end{array}\right.
\end{equation}
where
\begin{eqnarray*}
&&A=(1-\eta)  \left(\frac{U_{t}}{U}+\alpha  \right) +(1-\eta^{2})U_{x}, \ \ w_{0}=\frac{u_{0y}}{U}, \ \  w_{1}=\frac{u_{1y}}{U}, \\
&&B=\frac{U_{t}}{U}+ \eta U_{x}+\alpha, \ \  C=\alpha+   U_{x}+ \frac{U_{t}}{U}\geq 0.
\end{eqnarray*}

\begin{Definition}(\cite{GGW}) \label{GGW.d}
 A function $w\in L^{\infty} (0,T; BV(\Omega))\cap L^{\infty}(Q_{T})$ is a weak solution to the problem (\ref{w.69}), if the following conditions hold: \\
(1) There exist some positive constants $c_{1}$ and  $c_{2}$ such that
 \begin{equation}
c_{1}(1-\eta)\leq w(\tau,\xi,\eta)\leq c_{2}(1-\eta), \ \ \forall (\tau,\xi,\eta)\in Q_{T};
\label{w.70}
 \end{equation}
(2) $w(\tau,\xi,\eta)$ satisfies the problem (\ref{w.69}) in the weak sense:
\begin{eqnarray*}
&& \int_{Q_{T}}\left(\frac{1}{w}(\phi_{\tau}+\eta(U\phi)_{\xi}+B\phi)- w\phi_{\eta\eta}\right) d\xi d\eta d\tau
+ \int_{0}^{T}\int_{0}^{X}(A\frac{\phi}{w}) (\tau,\xi,0)  d\xi   d\tau \\
&& =-\int_{0}^{T}\int_{0}^{1} \eta\frac{U(\tau, 0)}{w_{1}} \phi(\tau, 0,\eta) d\eta d\tau
-\int_{\Omega} \frac{1}{w_{0}}\phi(0, \xi,\eta) d\xi d\eta + \int_{0}^{T}\int_{0}^{X}( \frac{C\phi}{w}) (\tau,\xi,0)  d\xi   d\tau
\end{eqnarray*}
for any test function $\phi\in C^{\infty}( \overline{Q}_{T})$ satisfying
$$\phi|_{\tau=T}=\phi|_{\xi=T}=\phi|_{\eta=1}=\partial_{\eta}\phi|_{\eta=0}=0.  $$
\end{Definition}
\begin{Theorem}(\cite{GGW})  (Well-posedness)
 For any $T>0$, assume that
$$ U\in C^{2}([0,T]\times[0,X]), \alpha\in C^{1}([0,T]\times[0,X]), $$
and the initial and boundary data $w_{0}\in C^{1}(\Omega), w_{1}\in C^{3}([0,T]\times[0,1])$ satisfy the compatibility condition. If the condition (\ref{w.68}) holds for all $\tau\geq 0$ and $0\leq\xi\leq L$, then the problem (\ref{w.69}) admits a global weak solution $w\in L^{\infty}((0,T); BV(\Omega))$ in the sense of Definition \ref{GGW.d}.
\end{Theorem}
The condition (\ref{w.68}) is fulfilled if the pressure $p(t,x)$ is favourable in the sense of (\ref{w.71}). When
the pressure $p(t,x)$ is not favourable, the condition  (\ref{w.68}) may still hold as long as $UB_{2}^{2}(t,x)$ is properly
large. Thus this result shows that the magnetic field has certain stabilizing effect on hydrodynamic flow.
When investigating the occurrence of a back-flow point in the Prandtl-Hartmann boundary layer, we
have the following result.
\begin{Theorem}(\cite{GGW})
(1) For any fixed $T>0$ and $X>0$, assume that the tangential velocity of the outer flow $U(t,x)$, and the initial and boundary data $u_{0}(x,y), u_{1}(t,y)$ satisfy the conditions (\ref{w.72}) and (\ref{w.73}), and
$$U\in C^{1}([0,T]\times [0,X]), \ \ B_{2}\in C^{1} ([0,T]\times[0,X] ). $$
Then, for the solution $(u, v,b_{1})\in C^{2}(D_{T})$ to the problem (\ref{w.67}), a first zero point of $\partial_{y} u(t,x,y)$, when the time evolves, must appear at the boundary $y=0$ if it exists for some time $0<t<T$. \\
(2) Furthermore, assume that
$$C(t,x)=\sigma B_{2}^{2}(t,x) +U_{x}+\frac{U_{t}}{U}\leq 0 , \ \ \forall(t,x)\in [0,T]\times [0,X] ,$$
and for some constant $C_{*}$ depending on $X,T,U$ and $B_{2}$, if the initial velocity $u_{0}(x,y)$ satisfies
$$\int_{0}^{X}\int_{0}^{+\infty} \frac{u_{0y}(X-x)^{\frac{3}{2}} }{\sqrt{u_{oy^{2}+u^{2}_{0}}} }dxdy \geq C_{*},  $$
then there exists a back-flow point $(t^{*}, x^{*})\in (0,T)\times[0,X]$, such that
 \begin{equation*}
\left\{\begin{array}{ll}
\partial_{y} u(t^{*}, x^{*})=0, \\
\partial_{y} u(t , x ,y)>0 ,  \ \ \forall 0<t<t^{*}, \ \ x\in [0,X],\ \  y\geq 0. \end{array}\right.
\end{equation*}
\end{Theorem}

 Ding,  Ji and  Lin (\cite{DJL}) in 2021 established  the boundary layer theory for 2D steady viscous incompressible MHD equations under the assumption of a moving boundary at $\{y=0\}$, and convergence rates in Sobolev sense.\\

The following steady viscous incompressible magnetohydrodynamics system in $\Omega=[0,L]\times\mathbb{R}_{+}$ with moving boundary conditions on velocity field and the perfect conducting boundary conditions on magnetic field at $\{y=0\}$:
 \begin{equation}\label{w.76}
\left\{\begin{array}{ll}
(u\partial_{x}+v\partial_{y})u- (h\partial_{x}+g\partial_{y})h+p_{x}=\nu\varepsilon (\partial_{xx}+ \partial_{yy})u, \\
(u\partial_{x}+v\partial_{y})v- (h\partial_{x}+g\partial_{y})g+p_{y}=\nu\varepsilon (\partial_{xx}+ \partial_{yy})v, \\
(u\partial_{x}+v\partial_{y})h- (h\partial_{x}+g\partial_{y})u =k\varepsilon (\partial_{xx}+ \partial_{yy})h, \\
(u\partial_{x}+v\partial_{y})g- (h\partial_{x}+g\partial_{y})v =k\varepsilon (\partial_{xx}+ \partial_{yy})g, \\
\partial_{x} u+\partial_{y} v=0, \ \  \partial_{x} h+\partial_{y} g=0, \\
(u,v,\partial_{y}h, g)|_{y=0}=(u_{b}, 0,0, 0),
\end{array}\right.
\end{equation}
where $(u,v)$ and $(h,g)$ are velocity and magnetic field respectively, and the given constant $u_{b}$ stands for the moving speed of the plate. Assume that the viscosity and resistivity coefficients have the same order of a small parameter $\varepsilon>0$. Suppose  that the outer ideal MHD flows are prescribed by
$$(u^{0}_{e}, v^{0}_{e}, h^{0}_{e},g^{0}_{e},p^{0}_{e})(x,y)$$
as $\varepsilon\rightarrow 0$, $( u^{0}_{e} , h^{0}_{e})(x,0)= ( \overline{u}^{0}_{e} , \overline{h}^{0}_{e}) $.

Following the idea of Prandt with the scaling boundary layer variable $(x', y')$:
$$x'=x,\ \ y'=\frac{y}{\sqrt{\varepsilon}}, $$
and letting
 $$(u^{\varepsilon},v^{\varepsilon},g^{\varepsilon},h^{\varepsilon},p^{\varepsilon} )(x',y')
=(u,\frac{v}{\sqrt{\varepsilon}},h,\frac{g}{\sqrt{\varepsilon}}, p )(x,y), $$
we can reduce the problem (\ref{w.76}) to the form
 \begin{equation}
\left\{\begin{array}{ll}
(u^{\varepsilon}\partial_{x}+v^{\varepsilon}\partial_{y})u^{\varepsilon}-(h^{\varepsilon}\partial_{x}+g^{\varepsilon}\partial_{y})h^{\varepsilon}
+p^{\varepsilon}_{x}=\nu\varepsilon \partial_{xx}u^{\varepsilon}+ \nu \partial_{yy} u^{\varepsilon}, \\
(u^{\varepsilon}\partial_{x}+v^{\varepsilon}\partial_{y})v^{\varepsilon}- (h^{\varepsilon}\partial_{x}+g^{\varepsilon}\partial_{y})g^{\varepsilon}
+\frac{p^{\varepsilon}_{y}}{ \varepsilon }=\nu\varepsilon \partial_{xx} v^{\varepsilon}+ \nu\partial_{yy} v^{\varepsilon}, \\
(u^{\varepsilon}\partial_{x}+v^{\varepsilon}\partial_{y})h^{\varepsilon}- (h^{\varepsilon}\partial_{x}+g^{\varepsilon}\partial_{y})u^{\varepsilon}
 =k\varepsilon \partial_{xx} h^{\varepsilon}+ k\partial_{yy} h^{\varepsilon}, \\
(u^{\varepsilon}\partial_{x}+v^{\varepsilon}\partial_{y})g^{\varepsilon}- (h^{\varepsilon}\partial_{x}+g^{\varepsilon}\partial_{y})v^{\varepsilon} =
k\varepsilon  \partial_{xx}g^{\varepsilon}+k \partial_{yy} g^{\varepsilon}, \\
\partial_{x} u^{\varepsilon}+\partial_{y} v^{\varepsilon}=0, \ \  \partial_{x} h^{\varepsilon}+\partial_{y} g^{\varepsilon}=0, \\
(u^{\varepsilon},v^{\varepsilon},\partial_{y}h^{\varepsilon}, g^{\varepsilon})|_{y=0}=(u_{b}, 0,0, 0).
\label{w.77}
\end{array}\right.
\end{equation}
Define the following norm of $\mathcal{X}$
\begin{eqnarray*}
\begin{aligned}
\|u^{\varepsilon}, v^{\varepsilon},h^{\varepsilon},g^{\varepsilon}\|_{\mathcal{X} }
:=&\|\{u^{\varepsilon}_{y}, v^{\varepsilon}_{y},\sqrt{\varepsilon}(h^{\varepsilon}_{x},g^{\varepsilon}_{x})\}\cdot y\|_{L^{2}}
+\| u^{\varepsilon}_{y}, v^{\varepsilon}_{y},\sqrt{\varepsilon}(h^{\varepsilon}_{x},g^{\varepsilon}_{x}) \cdot y\|_{L^{2}} \\
+&\|\{u^{\varepsilon}_{yy}, v^{\varepsilon}_{yy},\sqrt{\varepsilon}(h^{\varepsilon}_{xy},g^{\varepsilon}_{xy}), \varepsilon (h^{\varepsilon}_{xx},g^{\varepsilon}_{xx}) \}\cdot y\|_{L^{2}} \\
+&\|u^{\varepsilon}, v^{\varepsilon},h^{\varepsilon},g^{\varepsilon}\|_{B}
+\varepsilon^{\frac{\gamma}{2}} \| u^{\varepsilon}_{y}, v^{\varepsilon}_{y},\sqrt{\varepsilon}(h^{\varepsilon} ,g^{\varepsilon} )
\|_{L^{\infty}},
\end{aligned}
\end{eqnarray*}
where the boundary term is defined by
$$ \|u^{\varepsilon}, v^{\varepsilon},h^{\varepsilon},g^{\varepsilon}\|_{B}
:=\|\{u^{\varepsilon}_{y},  \sqrt{\varepsilon}(u^{\varepsilon}_{x},h^{\varepsilon}_{x})\}\cdot y\|_{L^{2}(x=L)}
+\| \sqrt{\varepsilon}(u^{\varepsilon}_{x},h^{\varepsilon}_{x}) \|_{L^{2}(x=L)} .$$
\begin{Theorem}(\cite{DJL}) \label{DJL.t}
Let $u_{b}>0$ be a constant tangential velocity of the viscous MHD flows on the boundary $\{y=0\}$, and the given smooth non-shear ideal MHD flows $(u^{0}_{e}, v^{0}_{e}, h^{0}_{e},g^{0}_{e},p^{0}_{e})(x,y)$ satisfy the following hypotheses:
 \begin{equation*}
 \begin{array}{ll}
0<c_{0}\leq h^{0}_{e}\ll u^{0}_{e} \leq C_{0}\leq \infty, \\
\|\frac{v^{0}_{e}}{y}\|_{L^{\infty}}\ll 1, \\
\|y^{k}\nabla^{m}(v^{0}_{e} ,g^{0}_{e} )\|_{L^{\infty}} \leq\infty \ \ \text{for sufficiently large} \ \ k, m\geq 0,\\
\|y^{k}\nabla^{m}(u^{0}_{e} ,h^{0}_{e} )\|_{L^{\infty}} \leq\infty  \ \ \text{for sufficiently large} \ \ k\geq 0, m\geq 1,\\
\|\langle y\rangle \partial_{y}(u^{0}_{e} ,h^{0}_{e} )\|_{L^{\infty}} \leq\delta_{0}  \ \ \text{ suitable small} \ \ \delta_{0} \geq 0.
\end{array}
\end{equation*}
In addition, let $m\geq 5$ be an integer, the outer ideal MHD flows are assumed to enjoy
$$M_{0}:=\sum\limits^{m+2}_{i=0} \|(\overline{u}^{0}_{e} ,\overline{h}^{0}_{e}, \overline{p}^{0}_{e} )(x) \|_{H^{m+2-i}(0,L)}<+\infty.$$
Moreover, for some positive constants $\vartheta_{0},\eta_{0}$ and small $\sigma_{0}$, suppose that
 \begin{equation*}
\left\{\begin{array}{ll}
 \overline{u}^{0}_{e}+u^{0}_{p}(0,y)>\overline{h}^{0}_{e}>\overline{h}^{0}_{e}+h^{0}_{p}(0,y)\geq \vartheta_{0}, \\
|u^{0}_{e}(0,y)+u^{0}_{p}(0,y)|\gg |h^{0}_{e}(0,y)+h^{0}_{p}(0,y)|\geq \eta_{0}, \\
|\langle y\rangle^{l+1}\partial_{y}(u^{0}_{p},h^{0}_{p})(0,y) |\leq \frac{1}{2}\sigma_{0}, \\
|\langle y\rangle^{l+1}\partial_{y}^{2}(u^{0}_{p},h^{0}_{p})(0,y) |\leq \frac{1}{2}\vartheta_{0}^{-1}.
 \end{array}\right.
\end{equation*}
Then there exist the remainder solutions $(u^{\varepsilon}, v^{\varepsilon},h^{\varepsilon},g^{\varepsilon})$  in the space $\mathcal{X}$ satisfying
$$\|(u^{\varepsilon}, v^{\varepsilon},h^{\varepsilon},g^{\varepsilon}) \|_{\mathcal{X}} \lesssim 1  $$
in $[0,L]\times[0,+\infty)$, where the positive number $L$ is sufficiently small relative to universal constant.

\end{Theorem}

\begin{Theorem}(\cite{DJL})
Under the hypotheses stated in Theorem \ref{DJL.t} with the profile $(u^{0}_{p},h^{0}_{p} )$, it holds that
$$\|(u-u^{0}_{e}-u^{0}_{p}, h-h^{0}_{e}-h^{0}_{p})\|_{L^{\infty}}+\|(v-v^{0}_{e}v, g-g^{0}_{e} )\|_{L^{\infty}}  \lesssim \sqrt{\varepsilon} $$
where $(u^{0}_{p},h^{0}_{p} )$ is the solution of Prandtl equations  (\ref{w.77}).
\end{Theorem}

Liu and Zhang (\cite{LZ2}) in 2021  proved the global existence and the large time decay estimate of solutions to the 2D MHD boundary layer equations (\ref{w.400}) with small initial data  in the upper space $\mathbb{R}^{2}_{+} =\{(x,y ):x\in\mathbb{R},y\in \mathbb{R}_{+} \} $:
  \begin{equation}\label{w.400}
\left\{\begin{array}{ll}
 \partial_{t} u_{1}-\partial_{y}^{2} u_{1}+u_{1}\partial_{x}u_{1}+u_{2}\partial_{y}u_{1}+\partial_{x}p
 =b_{1}\partial_{x}b_{1}+b_{2}\partial_{y}b_{1},\\
\partial_{t} b_{1}-k\partial_{y}^{2} b_{1}+u_{1}\partial_{x}b_{1}+u_{2}\partial_{y}b_{1}
 =b_{1}\partial_{x}u_{1}+b_{2}\partial_{y}u_{1},\\
u_{1}|_{y=0}=u_{2}|_{y=0}=0,\quad \partial_{y}b_{1}|_{y=0}=b_{2}|_{y=0}=0,\\
\lim\limits_{y\rightarrow+\infty}u_{1}=U_{1},\quad \lim\limits_{y\rightarrow+\infty}b_{1}=B_{1},\\
u_{1}|_{t=0}=u_{1,0},\quad b_{1}|_{t=0}=b_{1,0},
 \end{array}\right.
\end{equation}
where $(u_{1},u_{2})$ and $(b_{1},b_{2})$ represent the velocity of fluid and the magnetic field respectively, $k>0$ is a constant which represents the difference between the Reynolds number and the magnetic
Reynolds number, $(U_{1},B_{1},p)(t,x)$ are the traces of the tangential fields and pressure of the
outflow on the boundary, which satisfy Bernoulli's law:
\begin{equation}\label{w.401}
\left\{\begin{array}{ll}
 \partial_{t} U_{1} +U_{1}\partial_{x}U_{1}+ \partial_{x}p =B_{1}\partial_{x}B_{1} ,\\
\partial_{t} B_{1} +U_{1}\partial_{x}B_{1}  =B_{1}\partial_{x}U_{1}.
 \end{array}\right.
\end{equation}
Take a cut-off function $\mathcal{X}\in C^{\infty}[0,\infty]$ with
\begin{equation*}
\mathcal{X}(y)=\left\{\begin{array}{ll}
 y\quad \text{if } y\geq 2,\\
0\quad \text{if } y\leq 2,
 \end{array}\right.
\end{equation*}
 and make the following function transformations:
 \begin{equation}\label{w.402}
\left\{\begin{array}{ll}
u=u_{1}-\mathcal{X}'(y)U\quad \text{and}\quad v=u_{2}+\mathcal{X} (y)\partial_{x}U ,\\
b=b_{1}-\mathcal{X}'(y)B-\overline{B}_{x}\quad \text{and}\quad h=b_{2}+\mathcal{X} (y)\partial_{x}B ,
 \end{array}\right.
\end{equation}
where $U=U_{1}$ and $B=B_{1}-\overline{B}_{x}$.
Then in view of (\ref{w.400}) and (\ref{w.401}), $(u,v,b,h)$ solves
\begin{eqnarray}\label{w.403}
\left\{\begin{array}{ll}
 \partial_{t} u -\partial_{y}^{2} u -\overline{B}_{k}\partial_{x}b+u \partial_{x}u -b \partial_{x}b+v\partial_{y}u-h \partial_{y}b
 +\mathcal{X}'(U  \partial_{x}u- B\partial_{x}b)+\mathcal{X}'(\partial_{x}U u- \partial_{x} B b) \\
\quad \quad+\mathcal{X} (- \partial_{x}U  \partial_{y}u+ \partial_{x} B\partial_{y}b)+\mathcal{X}''( U v- B h)=m_{U},\\
\partial_{t} b -k\partial_{y}^{2}b -\overline{B}_{k}\partial_{x}u+u \partial_{x}b -b \partial_{x}u+v\partial_{y}b-h \partial_{y}u
 +\mathcal{X}'(U  \partial_{x}b- B\partial_{x}u)+\mathcal{X}'(\partial_{x}B u- \partial_{x}U b)\\
\quad \quad +\mathcal{X} (- \partial_{x}U  \partial_{y}b+ \partial_{x} B\partial_{y}u)+\mathcal{X}''( B v- U h)=m_{B},\\
\partial_{x}u+\partial_{y}v =0,\quad \partial_{x}b+\partial_{y}h =0,\\
u |_{y=0}=v|_{y=0}=0,\quad \partial_{y}b |_{y=0}=h|_{y=0}=0,\\
\lim\limits_{y\rightarrow+\infty}u =0,\quad \lim\limits_{y\rightarrow+\infty}b =0,\\
\lim\limits_{y\rightarrow+\infty}v =0,\quad \lim\limits_{y\rightarrow+\infty}h=0,\\
u |_{t=0}=u_{ 0}=u_{1,0}-\mathcal{X}'U_{0},\quad b |_{t=0}=b_{ 0}=b_{1,0}-\mathcal{X}'B_{0}-\overline{B}_{k},
 \end{array}\right.
\end{eqnarray}
where $(U_{0}(x), B_{0}(x))=(U(0,x), B(0,x))$, and the terms $(m_{U},m_{B})$ are given by
 \begin{equation}\label{w.404}
\left\{\begin{array}{ll}
m_{U}=(1-\mathcal{X}')(\partial_{t}U-\overline{B}_{k}\partial_{x}B)+\mathcal{X}'''U+(1-\mathcal{X}'^{2}+\mathcal{X}\mathcal{X}'')
(U\partial_{x}U-B\partial_{x}B),\\
m_{B}=(1-\mathcal{X}')(\partial_{t}B-\overline{B}_{k}\partial_{x}U)+\mathcal{X}'''B+(1-\mathcal{X}'^{2}-\mathcal{X}\mathcal{X}'')
(U\partial_{x}B-B\partial_{x}U).
 \end{array}\right.
\end{equation}
Note that there exist two potential functions $(\phi,\varphi)$ so that $(u,b)=\partial_{y}(\phi,\varphi)$ and $(v,h)=\partial_{x}(\phi,\varphi)$, then  integrating $(u, b)$ equations of (\ref{w.403}) with respect to $y$ variable over $[y,\infty]$, we  get
 \begin{equation}\label{ww.403}
\left\{\begin{array}{ll}
 \partial_{t} \varphi -\partial_{y}^{2}\varphi -\overline{B}_{k}\partial_{x}\psi+u \partial_{x}\varphi -b \partial_{x}\psi
 +2\int_{y}^{\infty} ( \partial_{x} \varphi \partial_{y}u-\partial_{x}\psi\partial_{y}b)dy'
 +\mathcal{X}'(U  \partial_{x}\varphi- B\partial_{x}\psi)\\
\quad \quad +2\int_{y}^{\infty}\mathcal{X}''(U  \partial_{x}\varphi- B\partial_{x}\psi)dy'
 +\mathcal{X}'(\partial_{x}U \varphi- \partial_{x} B \psi) \\
\quad \quad+\mathcal{X} (- \partial_{x}U  \partial_{y}\varphi+ \partial_{x} B\partial_{y}\psi) =m_{U},\\
\partial_{t}\psi-k\partial_{y}^{2}\psi -\overline{B}_{k}\partial_{x}\varphi+u \partial_{x}\psi -b \partial_{x}\varphi
+\mathcal{X}'(U  \partial_{x}\varphi- B\partial_{x}\varphi)\\
\quad \quad +\mathcal{X} (- \partial_{x}U   b+ \partial_{x} B u) =m_{B},\\
\varphi |_{y=0}=\psi|_{y=0}=0,\quad \partial_{y}b |_{y=0}=h|_{y=0}=0,\\
\lim\limits_{y\rightarrow+\infty}u =0,\quad \lim\limits_{y\rightarrow+\infty}\varphi =\lim\limits_{y\rightarrow+\infty}\psi =0,\\
\varphi |_{t=0}=\varphi_{ 0}=-\int_{y}^{\infty}u_{0}dy',\quad \psi |_{t=0}=\psi_{ 0}=-\int_{y}^{\infty}b_{0}dy'.
 \end{array}\right.
\end{equation}
Define
$$\Delta^{h}_{k}a=\mathcal{F}^{-1}(\varphi2^{-k}|\xi| \hat{a}),\quad S^{h}_{k}a=\mathcal{F}^{-1}(\mathcal{X}2^{-k}|\xi| \hat{a}), $$
where $\mathcal{F}^{-1}$ and $\hat{a}$  denote the partial Fourier transform of the distribution $a$ with respect to $x$  variable, and
\begin{eqnarray}
&& \|u \|_{\mathcal{B}^{s,0} }:=\|(2^{ks} \|\Delta^{h}_{k}a(t)\|_{L^{2}_{+}} )_{k\in \mathbb{Z}}\|_{L^{1}(\mathbb{Z})},\\
&& \|\|_{\tilde{L}^{p}(T_{0}, T; \mathcal{B}^{s,0})}:=\sum\limits_{k\in \mathbb{Z}}2^{ks}\left(\int_{T_{0}}^{T} \|\Delta^{h}_{k}a(t)\|_{L^{2}_{+}}^{p}dt\right)^{\frac{1}{p}}.
\end{eqnarray}
Here the phase function $\Phi$ is defined by
$$\Phi(t,\xi):=(\delta-\lambda\theta(t))|\xi|, $$
and the weighted function $\Psi(t,y)$ is determined by
$$ \Psi(t,y):=\frac{y^{2}}{8\langle t\rangle}, \quad  \langle t\rangle =1+t, $$
which satisfies
$$\partial_{t} \Psi+2(\partial_{y}\Psi)^{2}=0. $$
\begin{Theorem}(\cite{LZ2})
 Let $k\in[0,12]$, $\overline{B}_{k}=\left\{\begin{array}{ll}
 1, ~if ~k=1, \\
0, ~otherwise ,
 \end{array}\right. $   and $\varepsilon,\delta>0$. Assume that the far field states $(U,B)$ satisfy
$$\|\langle t\rangle^{\frac{9}{4}}e^{\delta|D_{x}|}(U,B)\|_{\widetilde{L}^{\infty}(\mathbb{R}^{+};\mathcal{B}_{h}^{\frac{3}{2}} )}
+\|\langle t\rangle^{\frac{7}{4}}e^{\delta|D_{x}|}(\partial_{t}U,\partial_{t}B,U,B)\|_{\widetilde{L}^{\infty}(\mathbb{R}^{+};
\mathcal{B}_{h}^{\frac{1}{2}} )} \leq \varepsilon.$$
Let the initial data $(u_{0},b_{0},\varphi_{0},\psi_{0})$   satisfy the compatibility
condition:  $u_{0}|_{y=0}=\partial_{y}b_{0}|_{y=0}=0$, $\int_{0}^{\infty}u_{0}dy= \int_{0}^{\infty}b_{0}dy=0$ and
$$\|e^{\frac{y^{2}}{8}} e^{\delta|D_{x}|}(u_{0},b_{0},\varphi_{0},\psi_{0})\|_{ \mathcal{B} ^{\frac{1}{2},0} }< \infty ,\quad
\|e^{\frac{y^{2}}{8}} e^{\delta|D_{x}|}(G_{0},H_{0} )\|_{ \mathcal{B} ^{\frac{1}{2},0} }\leq \sqrt{\varepsilon} ,$$
where $G_{0}=u_{0}+\frac{y}{2\langle t\rangle}\varphi_{0}$ and $H_{0}=b_{0}+\frac{y}{2k\langle t\rangle}\psi_{0}$. Then there exist positive constants $\lambda\varsigma $ and $\varepsilon_{0}(\lambda,\varsigma,k\delta)$ so that for $\varepsilon\leq \varepsilon_{0}$ and $l_{k}=\frac{k(2-k)}{4}\in[0,\frac{1}{4}]$, the system (\ref{ww.403}) has a unique global
solution $(u,b)$ which satisfies $\sup\limits_{t\in[0,\infty]}\theta(t)\leq\frac{\delta}{2\lambda}$, and
$$\| e^{\Psi}(u,b)\Phi\|_{\widetilde{L}^{\infty}_{t}( \mathcal{B} ^{\frac{1}{2},0} )}
+\sqrt{l_{k}}\| e^{\Psi} \partial_{y}(u,b)\Phi\|_{\widetilde{L}^{2}_{t}( \mathcal{B} ^{\frac{1}{2},0} )}
 \leq \| e^{\frac{y^{2}}{8}}e^{\delta|D_{x}|}(u_{0},b_{0}) \|_{ \mathcal{B} ^{\frac{1}{2},0} }
+C\sqrt{\varepsilon}.$$
Furthermore, for any $t>0$ and $\gamma\in[0,1]$, there holds that
\begin{eqnarray*}
&&\|\langle \tau\rangle^{\frac{1}{2}+l_{k}-\varsigma\varepsilon}e^{\Psi}(u,b)\Psi\|_{\widetilde{L}^{\infty}_{t}( \mathcal{B}^{\frac{1}{2},0} )}
+\|\langle \tau\rangle^{\frac{1}{2}+l_{k}-\varsigma\varepsilon}e^{\Psi}\partial_{y}(u,b)\Psi\|_{\widetilde{L}^{2}_{[\frac{t}{2},t]}( \mathcal{B}^{\frac{1}{2},0} )} \\
&&\quad \leq C(\|e^{\frac{y^{2}}{8}} e^{\delta|D_{x}|}(u_{0},b_{0},\varphi_{0},\psi_{0})\|_{ \mathcal{B} ^{\frac{1}{2},0} }+ \sqrt{\varepsilon}),\\
&&\|\langle \tau\rangle^{\frac{1}{2}+l_{k}-\varsigma\varepsilon}e^{\Psi}(G,H)\Psi\|_{\widetilde{L}^{\infty}_{t}( \mathcal{B}^{\frac{1}{2},0} )}
+\|\langle \tau\rangle^{\frac{1}{2}+l_{k}-\varsigma\varepsilon}e^{\Psi}\partial_{y}(G,H)\Psi\|_{\widetilde{L}^{2}_{[\frac{t}{2},t]}( \mathcal{B}^{\frac{1}{2},0} )}\\
&&\quad\leq C\sqrt{\varepsilon}(1+\|e^{\frac{y^{2}}{8}} e^{\delta|D_{x}|}(u_{0},b_{0},\varphi_{0},\psi_{0})\|_{ \mathcal{B} ^{\frac{1}{2},0} }),\\
&&\|\langle \tau\rangle^{\frac{1}{2}+l_{k}-\varsigma\varepsilon}e^{\gamma\Psi}(u,b)\Psi\|_{\widetilde{L}^{\infty}_{t}( \mathcal{B}^{\frac{1}{2},0} )}
+\|\langle \tau\rangle^{\frac{1}{2}+l_{k}-\varsigma\varepsilon}e^{\gamma\Psi}\partial_{y}(u,b)\Psi\|_{\widetilde{L}^{2}_{[\frac{t}{2},t]}( \mathcal{B}^{\frac{1}{2},0} )}\\
&&\quad \leq C\sqrt{\varepsilon}(1+\|e^{\frac{y^{2}}{8}} e^{\delta|D_{x}|}(u_{0},b_{0},\varphi_{0},\psi_{0})\|_{ \mathcal{B} ^{\frac{1}{2},0} }).
\end{eqnarray*}
\end{Theorem}

\subsection{Compressible MHD Boundary Layer Equations}

Jiang and Zhang (\cite{JZ}) in 2017 studied the non-resistive limit of the global solutions with large data and the global well-posedness of the following  compressible non-resistive MHD equations:
\begin{eqnarray}\left\{\begin{array}{ll}
\rho_t+(\rho u)_x=0,\\
(\rho u)_t+(\rho u^2+P(\rho)+\frac{1}{2}b^2)_x=\lambda u_{xx},\\
b_t+(ub)_x=\nu b_{xx},
\end{array}\right.\label{51.1}\end{eqnarray}
with the following initial and boundary conditions:
\begin{eqnarray}\left\{\begin{array}{ll}
(\rho,u,b)(x,0) = (\rho_0,u_0 ,b_0 )(x),x\in[0,1],\\
u(0,t) = u(1,t) = 0, b(0,t) = b_1 (t), b(1,t) = b_2 (t), t\geq0,
\end{array}\right.\label{51.2}\end{eqnarray}
where $P(\rho)=A\rho^{\gamma}$, $A>0$, $\gamma>1$.

The one-dimensional compressible, viscous, non-resistive MHD equations are as follows:
\begin{eqnarray}\left\{\begin{array}{ll}
\bar{\rho}_t+(\bar{\rho} \bar{u})_x=0,\\
(\bar{\rho}\bar{ u})_t+(\bar{\rho} \bar{u}^2+P(\bar{\rho})+\frac{1}{2}\bar{b}^2)_x=\lambda \bar{u}_{xx},\\
\bar{b}_t+(\bar{u}\bar{b})_x=\nu \bar{b}_{xx},
\end{array}\right.\label{51.3}\end{eqnarray}
which subjects to the following initial and boundary conditions:
\begin{eqnarray}\left\{\begin{array}{ll}
(\bar{\rho},\bar{u},\bar{b})(x,0) = ({\rho}_0,{u}_0 ,{b}_0 )(x),x\in[0,1],\\
\bar{u}(0,t) = \bar{u}(1,t) = 0, t\geq0.
\end{array}\right.\label{51.4}\end{eqnarray}
\begin{Theorem}(\cite{JZ})
Assume that $P(\bar{\rho})=A\bar{\rho}^{\gamma}$, with $A>0$, $\gamma>1$, and that the initial data $(\rho_0, u_0, b_0 )$
satisfy
\begin{eqnarray}
 \inf_{x\in[0,1]}\rho_0(x)\geq 0,\ (\rho_0,b_0)\in H^1,\ u_0\in H_0^1,\ u_0(0)=u_0(1)=0. \label{51.5}
\end{eqnarray}
Then for any $0 < T<\infty$ , there exist a positive constant $C$ and a unique global strong
solution $(\bar{\rho},\bar{u},\bar{b})$ to the initial-boundary value problem \eqref{51.3}-\eqref{51.4} on $(0,1)\times[0,T]$ ,
such that
\begin{eqnarray}
0<C^{-1}\leq \bar{\rho}\leq C,\ \forall (x,t)\in(0,1)\times[0,T], \label{51.6}
\end{eqnarray}
and
\begin{eqnarray}
 (\bar{\rho},\bar{b})\in H^1,\ \bar{u}\in H^1_0, \ (\bar{\rho}_t,\bar{b}_t)\in L^2,\ (\bar{u}_{t}, \bar{u}_{xx})\in L^2(0,T;L^2(\Omega)).\label{51.7}
 \end{eqnarray}
Assume that $P(\bar{\rho})=A\bar{\rho}^{\gamma}$, with $A>0$, $\gamma>1$. Moreover, in addition to \eqref{51.5}, suppose
that for any given $T\in(0,\infty)$,
\begin{eqnarray}
 (b_1,b_2)\in C^1([0,T]),\ b_0(0)=b_1(0),\ b_0(1)=b_2(0).\label{51.8}
 \end{eqnarray}
Then for each fixed $\nu>0$ , there exist a positive constant $C$ and a unique global strong
solution   to the initial-boundary value problem \eqref{51.1}-\eqref{51.2} on $(0,1)\times[0,T]$,
such that
\begin{eqnarray*}
0<C^{-1}\leq {\rho}\leq C,\ \forall (x,t)\in(0,1)\times[0,T],
\end{eqnarray*}
and
\begin{eqnarray*}
&& \sup\limits_{0\leq t\leq T}(\|u_x\|_{L^2}^2+\|b\|_{L^\infty}^2+\nu^{1/2}\|\rho_x\|_{L^2}^2+\nu^{1/2}\|b_x\|_{L^2}^2) \\
&&\quad +\int_0^T\Big(\|u_t\|_{L^2}^2+\nu^{1/2}\|u_{xx}+\nu^{3/2}\|b_{xx}\|_{L^2}^2\Big)ds\leq C.
\end{eqnarray*}
Moreover, as $\nu\rightarrow 0$ ,
\begin{eqnarray*}\left\{\begin{array}{ll}
(\rho,u,b)\rightarrow(\bar{\rho},\bar{u},\bar{b}),\ \ strongly \ in\ L^\infty(0,T;L^2),\\
\nu b_x\rightarrow0,\ u_x\rightarrow\bar{u}_x\ \ strongly \ in\ L^\infty(0,T;L^2)
\end{array}\right. \end{eqnarray*}
and
\begin{eqnarray*}
 \sup\limits_{0\leq t\leq T}(\|u-\bar{u}\|_{L^2}^2+\|b-\bar{b}\|_{L^2}^2+\|\rho-\bar{\rho}\|_{L^2}^2)
+\int_0^T\|u_x-\bar{u}_x\|_{L^2}^2ds\leq C\nu^{1/2}.
\end{eqnarray*}
Here, $C>0$ is a positive constant independent of $\nu$.
\end{Theorem}
\begin{Theorem}(\cite{JZ})
In addition to \eqref{51.5} and \eqref{51.8}, assume further that $(\rho_0, u_0, b_0 )\in H^2$.  Then any function $\delta(\nu)\in(0,1/2)$, satisfying
\begin{eqnarray*}
 \delta(\nu)\rightarrow 0\ and\ \frac{\delta(\nu)}{\nu^{1/2}}\rightarrow \infty,\ as\ \nu\rightarrow0,
 \end{eqnarray*}
 such that
\begin{eqnarray*}
 \lim\limits_{\nu\rightarrow 0}\Big(\|\rho-\bar{\rho}\|^2_{L^\infty(0,T;C(\bar{\Omega}_{\delta(\nu)}))}+
 \|b-\bar{b}\|^2_{L^\infty(0,T;C(\bar{\Omega}_{\delta(\nu)}))}\Big)=0,
 \end{eqnarray*}
 and
\begin{eqnarray*}
 \lim\limits_{\nu\rightarrow 0}\Big(\|\rho-\bar{\rho}\|^2_{L^\infty(0,T;C(\bar{\Omega}))}+
 \|b-\bar{b}\|^2_{L^\infty(0,T;C(\bar{\Omega}))}\Big)=0,
 \end{eqnarray*}
provided $b_i(t)\neq \bar{b}_i(t)( i = 1,2 )$, where $\bar{b}_1(t)$ and $\bar{b}_2(t)$ denote the boundary values of $\bar{b}(x,t)$ on the boundaries $x = 0,1$ , respectively. Here $ C(\bar{\Omega}_{\delta(\nu)}):= \{ x\in \Omega
|\delta<x<1-\delta \}$.
\end{Theorem}

Under the assumption that initial tangential magnetic field is not zero and density is a small perturbation of the outer constant flow in supernorm,
Gao, Huang and Yao (\cite{ghy}) in 2018 established the local-in-time existence and uniqueness of inhomogeneous incompressible MHD boundary layer equation
 in weighted   Sobolev spaces by the energy method:
 \begin{eqnarray}\left\{\begin{array}{ll}
\partial_t\rho+u_1\partial_x\rho+u_2\partial_y\rho=0,\\
\rho \partial_tu_1+\rho u_1\partial_xu_1+\rho u_2\partial_yu_1-\kappa\partial_{y}^2u_1=h_1\partial_xh_1+h_2\partial_yh_1=0,\\
\partial_th_1-\partial_y(u_2h_1-u_1h_2)-\kappa\partial_{y}^2h_1=0,\\
\partial_xu_1+\partial_yu_2=0,\partial_xh_1+\partial_yh_2=0,
\end{array}\right.\label{29.1}\end{eqnarray}
where the density $\rho:=\rho(t,x,y)$,  velocity field $(u_1,u_2)$, the magnetic field $(h_1,h_2)$ are unknown functions. The boundary conditions for equation \eqref{29.1} are
given by
\begin{eqnarray}\left\{\begin{array}{ll}
(u_1,u_2,\partial_yh_1,h_2)|_{y=0}=0, \\ \lim\limits_{y\rightarrow+\infty}\rho=\lim_{y\rightarrow+\infty}u_1=\lim_{y\rightarrow+\infty}h_1=1.
\end{array}\right.\label{29.2}\end{eqnarray}
We supplement the MHD boundary layer equations \eqref{29.1} with the initial data
\begin{eqnarray}
(\rho,u_1,h_1)(0,x,y)=(\rho_0,u_{10},h_{10})(x,y),\label{29.3}
\end{eqnarray}
satisfying
\begin{eqnarray}
 0<m_0\leq\rho_0\leq M_0<+\infty,\label{29.4}
 \end{eqnarray}
and
\begin{eqnarray}
 \|(\rho_0,u_{10},h_{10})\|_{B_l^m}\leq C_0<+\infty,\label{29.5}
 \end{eqnarray}
where $m_0 ,M_0 ,C_0 > 0$ are positive constants.
Define the functional space $B_l^m$  for a pair of function $(\rho,u_1,h_1)(t,x,y)$ as  follows:
\begin{eqnarray*}
B_l^m=\left\{(\rho-1,u_1-1,h_1-1)\in L^\infty([0,T];L^2_l(\Omega)):esssup_{0\leq t\leq T}\|(\rho,u_1,h_1)\|_{B^m_l}<+\infty\right\}, \\ B_{BL,app}^{m,l}=\left\{(\rho-1,u_1-1,h_1-1)\in H_l^{4k}|\partial_t(\rho,u_1,h_1),k=1,2,...,m    \right\}
\end{eqnarray*}
  and
 \begin{eqnarray}\left\{\begin{array}{ll}
  B_{BL}^{m,l}=\text{the closure of}\ B_{BL,app}^{m,l} \text{in the norm } \|\cdot\|_{B_l^m},\\
  L^{2}_{l}(\Omega)=\left\{ f(x,y):\Omega\rightarrow\mathbb{R}, \|f\|_{L^{2}_{l}(\Omega)}^{2}:=\int_{\Omega}\langle y\rangle  ^{2l}|f|^{2}dxdy<+\infty, \langle y\rangle  =1+y \right\}, \\
  L^{\infty}_{l}(\Omega)=\left\{ f(x,y):\Omega\rightarrow\mathbb{R}, \|f\|_{L^{\infty}_{l}(\Omega)} := \sup\limits_{(x,y)\in\Omega }|\langle y\rangle  ^{ l} f| dxdy<+\infty, \langle y\rangle  =1+y \right\}.
\end{array}\right.\label{29.111} \end{eqnarray}
\begin{Theorem}(\cite{ghy})
 Let $m \geq 5$ be an integer and $l\geq 2$ be a real number. Assume the initial
data $(\rho_0 ,u_{10} ,h_{10} ) \in B^{m,l}_{BL}$
given in \eqref{29.111} and satisfying \eqref{29.4}-\eqref{29.5}. Moreover, there exists a small constant $\delta_0 > 0$ such that
$$h_{10}\geq2\delta_0, \ for \ all \ (x,y)\in\Omega$$
and
$$\|\rho_0-1\|_{L_0^\infty(\Omega)}\leq\frac{2l-1}{16}\delta^2_0,\ \ \|\partial_yu_{10}\|_{L^\infty_1(\Omega)}\leq(2\delta_0)^{-1}.$$
Then, there exist a positive time $0<T_0$ and a unique solution $(\rho,u_1,u_2,h_1,h_2)$
 to the initial boundary value problem \eqref{29.1}-\eqref{29.2}, such that
\begin{eqnarray*}
&&\sup_{0\leq t\leq T_0}{\|(\rho-1,u_1-1,h_1-1)(t)\|_{H^m}^2+\|\partial_y(\rho,u_1,h_1)(t)\|_{H_1^{m-1}}^2
+\|\partial_y\rho(t)\|_{H_1^{1,\infty}}^2}  \\
&& +\int_0^{T_0}\|\partial_y(\sqrt\mu u_1,\sqrt\kappa h_1)(t)\|^2_{H_1^m}dt+\int_0^{T_0}\|\partial_y(\sqrt\mu u_1,\sqrt\kappa h_1)(t)\|^2_{H_1^{m-1}}dt\leq C_0<+\infty,
\end{eqnarray*}
here $C_{0}$ depends only on $l,\delta_{0}$.
\end{Theorem}

Ye and Zhang (\cite{YZ}) in  2020 considered that the 3D MHD flow with spatial variables $(x, y, z)$ being  moving only in the longitudinal direction $x$ and uniformly in the transverse direction $(y, z) .$ Then, based on the specific choice of dependent variables $\rho=\rho(x, t), \boldsymbol{u}=(u, \boldsymbol{w})(x, t), \boldsymbol{B}=(b, \boldsymbol{b})(x, t)$ and
$\theta=\theta(x, t),$ the compressible viscous heat-conducting MHD equations can be written in the following form:
\begin{equation}
\left\{\begin{array}{l}
\rho_{t}+(\rho u)_{x}=0 ,\\
(\rho u)_{t}+\left(\rho u^{2}+R \rho \theta+\frac{|b|^{2}}{2}\right)_{x}=\left(\lambda u_{x}\right)_{x}, \\
(\rho w)_{t}+(\rho u \boldsymbol{w}-\boldsymbol{b})_{x}=\left(\mu \boldsymbol{w}_{x}\right)_{x}, \\
\boldsymbol{b}_{t}+(u \boldsymbol{b}-\boldsymbol{w})_{x}=\left(\nu \boldsymbol{b}_{x}\right)_{x} ,\\
c_{V}(\rho \theta)_{t}+c_{V}(\rho u \theta)_{x}+R \rho \theta u_{x}=\left(\kappa(\theta) \theta_{x}\right)_{x}+\lambda u_{x}^{2}+\mu\left|\boldsymbol{w}_{x}\right|^{2}+v\left|\boldsymbol{b}_{x}\right|^{2}, \label{YZ1}
\end{array}\right.
\end{equation}
where the unknown functions $\rho, u \in \mathbb{R}, \boldsymbol{w}=\left(w^{1}, w^{2}\right) \in \mathbb{R}^{2}, \boldsymbol{b}=\left(b^{1}, b^{2}\right) \in \mathbb{R}^{2}$ and $\theta$
are the density of the fluid, the longitudinal velocity, the transverse velocity, the transverse magnetic field, and the temperature, respectively. The constants $\lambda>0, \mu \geq 0$ are the bulk and shear viscosity coefficients respectively, $v>0$ is the resistivity coefficient, $\kappa=\kappa(\theta)>0$ is the heat-conductivity coefficient, and $R, c_{V}>0$ are the perfect gas constant and the specific heat at constant volume, respectively.

\begin{Theorem}(\cite{YZ})
For any given $0<T<\infty$, assume that
$$
\left\{\begin{array}{l}
\rho_{0}>0, \quad \theta_{0}>0, \quad\left(\rho_{0}, u_{0}, \boldsymbol{w}_{0}, \boldsymbol{b}_{0}, \theta_{0}\right) \in H^{2}, \quad\left(\boldsymbol{w}_{1}, \boldsymbol{w}_{2}\right)(t) \in C^{2}([0, T]), \\
\left.\left(u_{0}, \boldsymbol{b}_{0}\right)\right|_{x=0,1}=0,\left.\quad \theta_{0 x}\right|_{x=0,1}=0,\left.\quad \boldsymbol{w}_{0}\right|_{x=0}=\boldsymbol{w}_{1}(0),\left.\quad \boldsymbol{w}_{0}\right|_{x=1}=\boldsymbol{w}_{2}(0).
\end{array}\right.
$$
Assume also that there exist two positive numbers $\underline{\kappa}$ and $\bar{\kappa}$ such that for some $q>0$,
$$
\kappa(\cdot) \in C^{2}(0, \infty) \text { and } \underline{\kappa} \xi^{q} \leq \kappa(\xi) \leq \bar{\kappa} \xi^{q}, \quad \forall \xi \geq 0.
$$
Then there exists a global unique strong solution $(\rho, u, \boldsymbol{w}, \boldsymbol{b}, \theta)~(\operatorname{resp} .(\bar{\rho}, \bar{u}, \overline{\boldsymbol{w}}, \overline{\boldsymbol{b}}, \bar{\theta}))$ to the
problem (\ref{YZ1}) with fixed $\mu>0$ on $(0,1) \times[0, T]$, such that
$$
\begin{array}{l}
\sup _{0 \leq t \leq T}\left(\|(\rho-\bar{\rho}, u-\bar{u}, \boldsymbol{b}-\overline{\boldsymbol{b}}, \theta-\bar{\theta})(t)\|_{H^{1}}^{2}+\|(\boldsymbol{w}-\overline{\boldsymbol{w}})(t)\|_{L^{2}}^{2}\right)  \\
+\int_{0}^{T}\left(\left\|(u-\bar{u}, \boldsymbol{b}-\overline{\boldsymbol{b}}, \theta-\bar{\theta})_{t}\right\|_{L^{2}}^{2}+\|(u-\bar{u}, \theta-\bar{\theta})\|_{H^{2}}^{2}\right) d t \leq C \mu^{1 / 2},
\end{array}
$$
where $C>0$ is a positive constant independent of $\mu .$
\end{Theorem}

\subsection{Ill-posedness of the Boundary Layer of MHD Equations}

Liu, Xie and Yang (\cite{lxy}) in  2018 investigated the sufficient degeneracy in the tangential magnetic field at a non-degenerate critical point of the tangential
 velocity field of shear flow, which indeed yields instability as for the classical Prandtl equations without magnetic field.
 They also partially showed the necessity of the non-degeneracy in the tangential magnetic field for the stability of the boundary layer of MHD in
 2D at least in Sobolev spaces:
\begin{eqnarray}\left\{\begin{array}{ll}
u_t+u_s\partial_xu+v\partial_yu_s-b_s\partial_xb-g\partial_yb_s=\kappa\partial_{y}^2u,\\
b_t-b_s\partial_xu+v\partial_yb_s+u_s\partial_xb-g\partial_yu_s=\nu\partial_{y}^2b,\\
\partial_xu+\partial_yv=0,\partial_xb+\partial_yg=0,\\
(u,v,b,g)|_{y=0}=0, \quad \lim\limits_{y\rightarrow+\infty}u=\lim\limits_{y\rightarrow+\infty}b=0,\\
u(0,x,y)=u_0(x,y),b(0,x,y)=b_0(x,y).
\end{array}\right.\label{73.1}\end{eqnarray}

\begin{Theorem}( \cite{lxy}) (Ill-posedness in the Sobolev Type Spaces)\\
(1) ~Let
$$u_s-\bar{u}\in \bigcap_0^1 C^0(\mathbb{R}_+;W_\alpha^{4-2i,\infty}), b_s-\bar{b}\in \bigcap_0^2 C^0(\mathbb{R}_+;W_\alpha^{6-2i,\infty}).$$
Suppose the initial velocity $u^s_0$ has a non-degenerate critical point $a>0$, i.e.,
$$u'_{s0}(a)=0,u''_{s0}(a)\neq0.$$
Moreover
$$\partial_y^jb_{s0}(a)=0, \ \ j=0,...6.$$
Then, there exists $\sigma > 0$, such that for any $\gamma > 0$,
$$\sup_{0\leq s\leq t\leq \gamma}\|e^{-\sigma(t-s)\sqrt{|\partial_x|}}T(t,s)\|_{L(H^m,H^{m-\mu})}=+\infty, \forall m\geq0, \mu\in [0,\frac{1}{2}).$$

(2)~Moreover, there exist some background shear flow $(u_s(t,y),b_s(t,y))$ and $\sigma > 0$ such that for any $\gamma > 0$, it holds that
\begin{eqnarray*}
\sup\limits_{0\leq s\leq t\leq \gamma}  \|e^{-\sigma(t-s)\sqrt{|\partial_x|}}T(t,s)\|_{L(H^m_1,H^ m_2)}=+\infty , \forall m_{1}, m_{2}\geq 0.
\end{eqnarray*}
\end{Theorem}
\begin{Theorem}(\cite{lxy})(Well-posedness in the Sobolev Type Spaces)
Let $u_s-\bar{u}\in C^1(\mathbb{R}_+;W_\alpha^{1,\infty})$, $b_s=u_s~(\text{or} -u_s)$ and $\kappa=\nu$. Then the linearized operator $T(t,s)$  satisfies
$$\|T(t,s)\|_{L(H^{m+1},H^m)}=\sup_{(u_0,b_0)\in K_{\alpha,\beta}}\frac{\|T(t,s)(u_0,b_0)\|_{H^m}}{\|(u_0,b_0)\|_{H^{m+1}}}\leq C<\infty,\forall m\geq0.$$
Here, the constant $C>0$ depends only on time $t$ and the background shear flows $u_s$ and $b_s$.
\end{Theorem}


\section{Other Boundary Layer Equations}
 Korpela (\cite{K}) in 1985  used  the Dorodnitsyn-Howarth transformation   across the boundary layer,  which reduces to the Blasius equation:
\begin{equation}
 f'''+ff''=0, \ \ \ \  f(0)=f'(0)=f'(\infty)-1=0,
 \label{ }
\end{equation}
the solution of which for small values of $\eta$ is
$$f=\frac{1}{2}\alpha\eta^{2}-\frac{1}{5!}\alpha^{2}\eta^{5}+\frac{11}{8!}\alpha^{3}\eta^{8}-... .$$
From the coefficients, it is obvious that the large Prandtl number  expansion breaks down in the range occupied by gases.
Thus the result has little utility  in aeronautical applications.

 Joseph and  Lefloch  (\cite{JL}) in 1999 considered  the initial-boundary-value problem of conservation laws
\begin{equation}
\partial_{t} u + \partial_{x}f(u)=0 , \ \ \ \  x>0,\ \ \ \ t>0 .
 \label{w.6}
\end{equation}
They derived conditions satisfied by the boundary layer, which take the form of a family of boundary entropy inequalities and a boundary-layer equation.
\begin{Theorem}( \cite{JL})
 The following statements hold for all convex entropy pairs $U, F$ associated with the system  (\ref{w.6}), all functions $\eta\in BV(\mathbb{R}_{+})$, and any bounded interval $(T_{1},T_{2})$.\\
1) When $\eta(t)\geq 0$, the distribution
$$y\mapsto\int_{T_{1}}^{T_{2}}\langle \mu_{y,t},F\rangle \eta(t)dt-\frac{\partial}{\partial y}\int_{T_{1}}^{T_{2}}\langle \mu_{y,t},U\rangle \eta(t)dt$$
is in fact a function of locally bounded variation and thus is defined pointwise as a right-continuous function. There exists a Young measure $\mu_{0,t}$,  such that the following limit exists and is given by $\mu_{0,t}$:
$$\lim\limits_{y\rightarrow 0^{+}}\int_{T_{1}}^{T_{2}}\langle \mu_{y,t},U\rangle \eta(t)dt=\int_{T_{1}}^{T_{2}}\langle \mu_{0,t},U\rangle \eta(t)dt.$$
When $\eta(t)\geq 0$, the function
$$x\mapsto\int_{T_{1}}^{T_{2}}\langle \nu_{x,t},F\rangle \eta(t)dt $$
has locally bounded variation. There exists a Young measure  $\nu_{0,t}$, the  ``trace"  of  $\nu_{x,t}$ at $x=0$, such that the following limit exists and is given by $\nu_{0,t}$:
$$\lim\limits_{x\rightarrow 0^{+}}\int_{T_{1}}^{T_{2}}\langle \nu_{x,t},U\rangle \eta(t)dt=\int_{T_{1}}^{T_{2}}\langle \nu_{0,t},F\rangle \eta(t)dt.$$
When $(U,F)=(u_{j},f_{j})$, $1\leq J\leq N$, all of the results above still hold when the function $\eta$ has no specific sign.\\
2) For all $0<y_{1}<y_{2}$ and for $t\in \mathbb{R}_{+}$, the following inequalities hold in the
sense of distributions:
\begin{eqnarray*}
\begin{aligned}
F(u_{B})+\nabla U(u_{B})\left(\langle \nu_{0,t},f\rangle-f(u_{B})\right)&\geq  \langle \mu_{y_{1},t},F\rangle-\partial_{y} \mu_{y_{1},t},U\rangle \\
&\geq   \langle \mu_{y_{2},t},F\rangle-\partial_{y} \mu_{y_{2},t},U\rangle\\
&\geq \langle \nu_{0,t},F\rangle  .
\end{aligned}
\end{eqnarray*}
3) Moreover,
\begin{eqnarray*}
\mu_{0,t}=\delta_{uB(t)} \ \ \ \  \text{a.e. } t\in \mathbb{R}_{+}
 \end{eqnarray*}
and, when $\eta\geq 0$,
$$\lim\limits_{y\rightarrow \infty}\left(\int_{T_{1}}^{T_{2}}\langle \mu_{y,t},F\rangle \eta(t)dt -\frac{d}{\partial y}\int_{T_{1}}^{T_{2}}\langle \mu_{y,t},U\rangle \eta(t)dt\right)\geq \int_{T_{1}}^{T_{2}}\langle \nu_{0,t},F\rangle \eta(t)dt.$$
Here $\langle \cdot,\cdot \rangle$ denotes inner product of $ BV(\mathbb{R}_{+})$.
\end{Theorem}

 Oleinik and Samokhin (\cite{OS}) in 1999 considered the  boundary layer in Non-Newtonian flows.
They considered the stationary plane-parallel symmetric flow of a pseudo-plastic fluid in the boundary layer is described by the system
 \begin{equation}\left\{
\begin{array}{ll}
  u \frac{\partial u }{\partial x}+v \frac{\partial u }{\partial y}
=  \nu \frac{\partial}{\partial y }\left(|\frac{\partial u }{\partial y}|^{n-1}\frac{\partial u }{\partial y} \right)+U(x)\frac{d U }{d x}, \ \  0<n<1, \\
\frac{\partial u }{\partial x}+\frac{\partial v }{\partial y} =0,
\end{array}
 \label{ww61}         \right.\end{equation}
in the domain $D =\{0<x<X, 0<y< \infty \}$, with the boundary conditions
 \begin{equation}
\begin{array}{ll}
u (0,y)=0, \ \ u (x,0)=0, \ \ v (x,0)=0, \\
u (x,y)\rightarrow U (x),  \ \ \text{as}\ \   y\rightarrow+\infty.
\end{array}
 \label{ww62}         \end{equation}
Assume  that $U(0)=0, U(x)>0$ for $x>0$; $\nu>0$, and
$$U(x) =xV(x), \ \  v_{0}(x)=x^{\frac{n-1}{n+1}}v_{1}(x), $$
where $V(x)>0, v_{1}(x)$ are bounded functions. In  (\ref{ww61}), introduce  new independent variables
$$\xi=x, \ \  \eta=\frac{u}{U}, \ \  w(\xi, \eta)=\frac{|u_{y}|^{n-1}u_{y}}{x^{\frac{n-1}{n+1}}U}, $$
then problem (\ref{ww61})-(\ref{ww62}) reduces to the equation
 \begin{equation}
\begin{array}{ll}
\nu n V^{\frac{1-n}{n}}|w|^{\frac{1-n}{n} } w_{\eta\eta} -\eta\xi Vw_{\xi}+Aw_{\eta}+Bw=0,
\end{array}
 \label{ww63}         \end{equation}
in the domain $\Omega=\{  0<\xi<X, 0<\eta <1 \}$, with the boundary conditions
\begin{eqnarray}
   w|_{\eta=1}=0, \ \  (\nu w |w|^{\frac{1+n}{n} }_{\eta }-v_{1}w|w|^{\frac{1-n}{n} }+C )|_{\eta=0}=0,
 \label{ww64}
\end{eqnarray}
 where
 $$A =(\eta^{2}-1)(V+\xi V_{x}) , \ \  B =-\eta \left(\frac{2n}{1+n} V+\xi V_{x}\right),   \ \  C=V^{\frac{ n-1}{n}}(V+\xi V_{x}). $$

 \begin{Theorem} \label{OS.8.43} (\cite{OS})
Assume that
$$U(x) =x(a+xa_{1}(x)), \ \  v_{0}(x)=x^{\frac{n-1}{n+1}}(b+xb_{1}(x)), $$
where $a=$const.$>0$, $b=$const.; $U>0$ for $x>0$; $a_{1}, a_{1x},a_{1xx}, b_{1}, b_{1x}$ are bounded. Then problem (\ref{ww61})-(\ref{ww62}) in the domain $D$,
with $X$ depending on $U, v_{0},n$, admits one and only one solution $u,v $ with the following properties: $u_{y}>0$ for $y\geq 0$ and $x>0$;
$$\frac{u}{U}, \ \ \frac{u_{y}^{n}}{x^{\frac{n-1}{n+1}}U}$$
are bounded and continuous in $\overline{D}$; $u>0$ for $y>0$ and $x>0$;
$$u\rightarrow U, \ \ \frac{u_{y}^{n}}{x^{\frac{n-1}{n+1}}U} \rightarrow 0\ \ \text{as} \ \ y\rightarrow\infty; $$
$u_{x}, u_{y}, u_{yy}$ are bounded and continuous in $D$;  $v$ is continuous in  $\overline{D}$ with respect to $y$ for $x>0$;
$v$ is continuous in $x$  and $y$  inside $D$;
$$\frac{u_{yy} }{x^{\frac{n-1}{n+1}}u_{y}^{2-n}} $$
is continuous in $D$ with respect to $y$; the following estimates hold:
\begin{eqnarray*}
&& x^{\frac{n-1}{n+1}}UY(\frac{u}{U})e^{-C_{3}x} \leq u_{y}\leq x^{\frac{n-1}{n+1}}UY(\frac{u}{U})e^{ C_{4}x}  , \\
&&  x^{\frac{n-1}{n+1}} Y_{n}(\frac{u}{U})e^{ C_{5}x} \leq \frac{u_{yy}}{u_{y}^{2-n}}\leq x^{\frac{n-1}{n+1}}Y_{n}(\frac{u}{U})e^{ -C_{7}x} , \\
&& \left(M_{2}^{\frac{1}{n}} \frac{n-1}{n+1} e^{\frac{C_{4}}{n}}x^{\frac{n-1}{n(n+1)}}U^{\frac{ 1-n}{n }}y +1\right)^{\frac{n+1}{n-1}}\leq 1-\frac{u}{U}
\leq \left(M_{1}^{\frac{1}{n}} \frac{n-1}{n+1} e^{-\frac{C_{3}}{n}}x^{\frac{n-1}{n(n+1)}}U^{\frac{ 1-n}{n }}y +1\right)^{\frac{n+1}{n-1}}.
\end{eqnarray*}

\end{Theorem}

  Oleinik and Samokhin (\cite{OS})  in 1999 considered
  system (\ref{ww61}) in the domain $D =\{0<x<X, 0<y< \infty \}$, whose the boundary conditions  are more complicated than those in (\ref{ww62}),
 \begin{equation}
\begin{array}{ll}
u (0,y)=u_{0}(y), \ \ u (x,0)=0, \ \ v (x,0)=v_{0}(x), \\
u (x,y)\rightarrow U (x),  \ \ \text{as}\ \   y\rightarrow+\infty.
\end{array}
 \label{ww65}         \end{equation}

\begin{Theorem} \label{OS.8.44} (\cite{OS})
Let $U(x), U_{x}(x), v_{0}(x)$ be continuous functions on the interval $0\leq x\leq X$, and let $ U_{x}(x)\geq$const.,  and assume also that $u_{0}(y)$ is
continuous for $0\leq y< \infty$; $u_{0}(0)=0$; $u_{0}(y)>0$ for $y>0$; $u_{0}(y)\rightarrow U(0)$ as $y\rightarrow \infty$; the derivative $u_{0y}(y)$
is continuous; $ u_{0y}(y)>0$; there exists a weak derivative $u_{0yy}(y)$ and
\begin{eqnarray*}
&& K_{1}^{\frac{1}{n}}U^{\frac{1}{n}}(0)\left(1-\frac{u_{0}(y)}{U(0)}\right)^{\frac{2}{n+1}} \leq u_{0y}(y)
\leq K_{2}^{\frac{1}{n}}U^{\frac{1}{n}}(0)\left(1-\frac{u_{0}(y)}{U(0)}\right)^{\frac{2}{n+1}} , \\
&& \int_{0}^{1} u_{0y}^{2n-3}(y)u_{0yy}^{2}(y)\left(1-\frac{u_{0}(y)}{U(0)}\right)^{\frac{2(1-n)}{n+1}}dy < \infty.
\end{eqnarray*}
Then problem (\ref{ww61}), (\ref{ww65}) has a weak solution $u(x,y), v(x,y)$ in the following sense: $u(x,y)$ is bounded, measurable, and continuous
in $\overline{D}$ with respect to $y$; $u>0$ for $y>0$; $u=0$ for $y=0$; $u_{y}(x, y)$ is bounded in $D$,
\begin{eqnarray*}
&& K_{11}^{\frac{1}{n}}U^{\frac{1}{n}}(x)\left(1-\frac{u(x,y)}{U(x)}\right)^{\frac{2}{n+1}} \leq u_{ y}(x,y)
\leq K_{12}^{\frac{1}{n}}U^{\frac{1}{n}}(x)\left(1-\frac{u(x,y)}{U(x)}\right)^{\frac{2}{n+1}} , \\
&& U(x)-U(x)\left(1+K_{12}^{\frac{1}{n}}\frac{ 1-n }{n+1}U^{\frac{1-n}{n }}(x)y \right)^{\frac{n+1}{n-1}}\leq u(x,y)
\leq  U(x)-U(x)\left(1+K_{12}^{\frac{1}{n}}\frac{ 1-n }{n+1}U^{\frac{1-n}{n }}(x)y \right)^{\frac{n+1}{n-1}} ;
\end{eqnarray*}
then the weak derivatives $u_{yy}, u_{yyy}$ exist and satisfy the inequalities
\begin{eqnarray*}
&& \int_{D}  u_{ y}^{2n-3}(y)u_{ yy}^{2}(y)\left(1-\frac{u }{U }\right)^{\frac{2(1-n)}{n+1}}dxdy < \infty, \\
&& \int_{D}  u^{2}\left(1-\frac{u }{U }\right)^{\frac{6-2n}{n+1}} \left( \frac{nu_{y}u_{yyy}+n(n-1)u_{yy} }{u_{y}^{4-n} }\right)^{2}u_{y}  dxdy < \infty ;
\end{eqnarray*}
$$\int_{D\cap\{(x,y):0\leq x\leq X-\varepsilon\}} v(x,y)(1+y)^{\frac{n+3}{n-1}}dxdy < \infty , \ \ \varepsilon >0, $$
here $v(x,y)$ is measurable in $D$ and equations  (\ref{ww61}) and the boundary conditions (\ref{ww65}) for $x=0,y=0$ hold in terms of the following integral identities
\begin{eqnarray*}
&& \int_{D} \left[\nu U^{\frac{1-n}{n}}\frac{u_{y}\phi_{yy}-\phi_{y}u_{yy}}{u_{y}^{2-n}} -u\phi_{x}U^{\frac{1-n}{n}}
+\left(\frac{u^{2}}{U^{2}}-1\right)\frac{\phi_{y}}{u_{y}}U^{\frac{1}{n}}U_{x} +^{\frac{n-1}{n}}u\phi U^{\frac{1-2n}{n}}U_{x} \right] dxdy  \\
&&\quad -\int_{0}^{X}U^{\frac{1-n}{n}}(x)\phi(x,0)v_{0}(x)dx -\int_{0}^{\infty}u_{0}(y)\phi(0,y)dy=0 , \\
&& \int_{D}   \left(\nu nu_{y}^{n-1}u_{yy}\phi+ \frac{1}{2}u^{2}\phi_{x}-vu_{y}\phi+UU_{x}\phi \right) dxdy
+ \int_{0}^{\infty}\frac{1}{2}u^{2}_{0}(y)\phi(0,y)U^{\frac{1-n}{n}}(0)dy=0 ,
\end{eqnarray*}
for any $\phi(x,y)$ such that the derivatives $\phi_{x}, \phi_{y}, \phi_{yy}$ are continuous and bounded in  $\overline{D}$, and $\phi(X,y)=0, \phi_{y}(x,0)=0$,
$$|\phi(x,y)|\leq K_{13}(1+K_{14}y)^{\frac{n+1}{n-1}}  \ \  \text{as} \ \  y\rightarrow \infty . $$
The solution $(u,v)$ of problem  (\ref{ww61}), (\ref{ww65}) with all these properties is unique.

\end{Theorem}

  Oleinik and   Samokhin (\cite{OS}) in 1999 considered
  the nonstationary motion of a pseudo-plastic fluid in a symmetric boundary layer is described by a system of the form $(0<n<1)$:
 \begin{equation}\left\{
\begin{array}{ll}
\frac{\partial u }{\partial t}+  u \frac{\partial u }{\partial x}+v \frac{\partial u }{\partial y}
=  \nu \frac{\partial}{\partial y }\left(|\frac{\partial u }{\partial y}|^{n-1}\frac{\partial u }{\partial y} \right)+\frac{\partial U }{\partial t}+
U(x)\frac{d U }{d x},  \\
\frac{\partial u }{\partial x}+\frac{\partial v }{\partial y} =0,
\end{array}
 \label{ww66}         \right.\end{equation}
in the domain $D =\{0<t<\infty, 0<x<X, 0<y< \infty \}$, with the boundary conditions
 \begin{equation}
\begin{array}{ll}
u (0,x,y)=u_{0}(x,y), \ \ u (t,0,y)=0, \ \ u (t,x,0)=0, \ \ v (t,x,0)=v_{0}(t,x), \\
u (t,x,y)\rightarrow U (x,y),  \ \ \text{as}\ \   y\rightarrow \infty.
\end{array}
 \label{ww67}         \end{equation}
Assume that
$$U(x) =xV(t,x), \ \ V(t,x)>0, \ \   v_{0}(t,x)=x^{\frac{n-1}{n+1}}v_{1}(x), $$
where $V  , v_{1}(x), V_{x}, V_{t}$ are bounded functions. In  (\ref{ww66}),  introduce new independent variables
$$\tau=t, \ \ \xi=x, \ \  \eta=\frac{u}{U}, \ \  w(\tau, \xi, \eta)=\frac{|u_{y}|^{n-1}u_{y}}{x^{\frac{n-1}{n+1}}U}, $$
then  problem (\ref{ww66})-(\ref{ww67}) reduces to the equation
 \begin{equation}
\begin{array}{ll}
\nu n V^{\frac{1-n}{n}}|w|^{\frac{1+n}{n} } w_{\eta\eta}-w_{\tau} -\eta\xi Vw_{\xi}+Aw_{\eta}+Bw=0,
\end{array}
 \label{ww68}         \end{equation}
in the domain $\Omega=\{0<\tau<\infty,  0<\xi<X, 0<\eta <1 \}$, with the boundary conditions
\begin{eqnarray}
 w|_{\tau=0}=w_{0}=\frac{|u_{0y}|^{n-1}u_{0y}}{x^{\frac{n-1}{n+1}}U}, \ \   w|_{\eta=1}=0, \ \
   (\nu w |w|^{\frac{1-n}{n} }_{\eta }-v_{1}w|w|^{\frac{1-n}{n} }+C )|_{\eta=0}=0,
 \label{ww69}
\end{eqnarray}
 where
 $$A =(\eta^{2}-1)(V+\xi V_{x})+(\eta-1)\frac{V_{t}}{V} , \ \  B =-\eta \left(\frac{2n}{1+n} V+\xi V_{x}\right)-\frac{V_{t}}{V},   \ \
  C=V^{\frac{ n-1}{n}}(V+\xi V_{x}+\frac{V_{t}}{V}). $$

\begin{Theorem} \label{OS.8.45} (\cite{OS})
Assume that
$$U(t,x) =x(a+xa_{1}(t,x)), \ \  v_{0}(t,x)=x^{\frac{n-1}{n+1}}(b+xb_{1}(t,x)), $$
where $a=$const.$>0$, $b=$const., $a_{1}(t,x)$ has bounded second order derivatives, $b_{1}(t,x)$ has bounded first order derivatives for
 $0\leq t<\infty, 0\leq x\leq X$. Let the function $u_{0}(x,y)$ be such that
$$\frac{|u_{0y}|^{n-1}u_{0y}}{x^{\frac{n-1}{n+1}}U}$$
satisfies the conditions
\begin{eqnarray*}
&& Ye^{-M_{4}\xi}\leq w_{0}\leq  Ye^{ M_{5}\xi}, \ \  |w_{0\xi}|\leq M_{11}(1-\eta)^{\frac{2n}{1+n}}, \ \
 Y_{\eta}e^ {M_{12}\xi}\leq w_{0\eta}\leq  Y_{\eta}e^{- M_{15}\xi}, \\
&& -M_{14}(1-\eta)^{\frac{2n}{1+n}}\leq \nu nV^{\frac{1-n}{n}}(0,x)w_{0}^{\frac{1+n}{n}}w_{0\eta\eta}+A(0,\xi,\eta)w_{0\eta}+B(0,\xi,\eta)w_{0}
\leq M_{15}x(1-\eta)^{\frac{2n}{1+n}}, \\
&&  (\nu w_{0\eta}  w_{0} ^{\frac{1-n}{n} } -v_{1}w_{0}^{\frac{1-n}{n} }+C )|_{\tau=0, \eta=0}=0.
\end{eqnarray*}
Then for some $X$ depending on $U,u_{0},v_{0}$, problem (\ref{ww66})-(\ref{ww67}) in $D$ has one and only one solution $u,v$ with the following properties:
$u>0$ for $y>0$ and $x>0$; $u\rightarrow U$ as $y\rightarrow \infty$; $u_{y}>0$ for $y\geq 0$;
$$\frac{u}{U}, \ \ \frac{u_{y}^{n}}{x^{\frac{n-1}{n+1}}U}   $$
are bounded and continuous in $\overline{D}$;
$$ \frac{u_{y}^{n}}{x^{\frac{n-1}{n+1}}U} \rightarrow 0\ \ \text{as} \ \ y\rightarrow\infty; $$
$u_{y}, u_{x}, u_{t}, u_{yy},v_{y}$ are bounded and continuous in $D$  with respect to $y$; $v$ is continuous in  $\overline{D}$ with respect to $y$ and bounded
for bounded $y$; equations (\ref{ww66})  hold almost everywhere in $D$; the following inequalities are satisfied:
\begin{eqnarray*}
&& x^{\frac{n-1}{n+1}}UY(\frac{u}{U})e^{-C_{3}x} \leq u_{y}^{n}\leq x^{\frac{n-1}{n+1}}UY(\frac{u}{U})e^{ C_{4}x}  , \\
&& U -U \left(1+M_{1}^{\frac{1}{n}}\frac{ 1-n }{n+1}\exp\left(-\frac{C_{3}x}{n}\right)x^{\frac{n-1}{n(n+1)}} U^{\frac{1-n}{n }}(x)y \right)^{\frac{n+1}{n-1}}
\leq u(x,y) \\
&& \quad
\leq U -U \left(1+M_{2}^{\frac{1}{n}}\frac{ 1-n }{n+1}\exp\left(\frac{C_{4}x}{n}\right)x^{\frac{n-1}{n(n+1)}} U^{\frac{1-n}{n }}(x)y \right)^{\frac{n+1}{n-1}},
\end{eqnarray*}
where $M_{1}$ and $M_{2}$ are positive constants.

\end{Theorem}

Assume that
$$U(x) =x\widetilde{V}(t,x), \ \ \widetilde{V}(t,x)>0,   $$
where $\widetilde{V} =V$ is  bounded functions, and $\widetilde{V}(t,x)$ is the corresponding solution $\widetilde{u}(t,x)$.
 \begin{Theorem} \label{OS.8.46}( \cite{OS})
Assume that the functions
\begin{eqnarray}
V-\widetilde{V}, \ \ V_{x}-\widetilde{V}_{x}, \ \ \frac{V_{t}}{V}-\frac{\widetilde{V}_{t}}{\widetilde{V}}, \ \  v_{1}-\widetilde{v}_{1}
 \label{ww70}
\end{eqnarray}
vanish for $t\geq t_{0}=$const.$>0$. Then
$$\left|\frac{u}{V}-\frac{\widetilde{u}}{\widetilde{V}} \right|\leq x K_{3}e^{-\alpha t}\ \ \text{in} \ \ D$$
for $0\leq t\leq y_{0}<\infty$, where $\alpha, K$ are positive constants.

\end{Theorem}

\begin{Theorem} \label{OS.8.47} (\cite{OS})
Assume that the functions  (\ref{ww70}) have their absolute values less than $\varepsilon$, and $V,\widetilde{V},u_{0}, \widetilde{u}_{0}$ are such that
$|w-\widetilde{w}|\leq \varepsilon$, where $u$ and $\widetilde{u}$ are the solutions of problem  (\ref{ww66})-(\ref{ww67}) with $V,v_{1}, u_{0}$ and
$\widetilde{V},\widetilde{v}_{1}, \widetilde{u}_{0}$  respectively. Then
$$\left|\frac{u}{V}-\frac{\widetilde{u}}{\widetilde{V}} \right|\leq x K_{6}\varepsilon \ \ \text{in} \ \ D$$
for $y\leq y_{0}<\infty$ and all $t\geq 0$, where $K_{6}$ is a positive constant.

\end{Theorem}

  Oleinik and   Samokhin (\cite{OS}) in 1999 considered
 the continuation problem for the following  stationary boundary layer equations in a dilatable medium by the von Mises variables:
\begin{equation}\left\{
\begin{array}{ll}
  u \frac{\partial u }{\partial x}+v \frac{\partial u }{\partial y}
=  \nu \frac{\partial}{\partial y }\left(|\frac{\partial u }{\partial y}|^{n-1}\frac{\partial u }{\partial y} \right)+U(x)\frac{d U }{d x}, \ \  1<n<\infty, \\
\frac{\partial u }{\partial x}+\frac{\partial v }{\partial y} =0,
\end{array}
 \label{ww71}         \right.\end{equation}
in the domain $D =\{0<x<X, 0<y< \infty \}$, with the boundary conditions
 \begin{equation}\left\{
\begin{array}{ll}
u (0,y)=u_{0}(y), \ \ u (x,0)=0, \ \ v (x,0)=v_{0}(x), \\
u (x,y)\rightarrow U (x),  \ \ \text{as}\ \   y\rightarrow+\infty.
\end{array}
 \label{ww72}    \right.     \end{equation}

\begin{Theorem} \label{OS.8.48} (\cite{OS})
 Assume that $u_{0}(y)>0$ for $y>0$, $u_{0}=0$, $u'_{0}(0)>0$, $u_{0}(y)\rightarrow  U(0)\neq 0 $ as $y\rightarrow \infty$,
$U(x)>0, U_{x}(x)$ and $v_{0}(x)$ are continuously differentiable on $[0,X], u_{0}(y), u'_{0}(y), u''_{0}(y)$ are bounded for $0\leq y\leq \infty$ and
satisfy the H$\ddot{o}$lder condition; the compatibility condition
$$\nu n\left|\frac{du_{0}}{dy}\right|^{n-1}\frac{d^{2}u_{0}}{dy^{2}}-v_{0}(0)\frac{du_{0}}{dy}+U(0)U_{x}(0)=\mathcal{O}(y^{2})$$
holds as $y\rightarrow 0$. Then, for some $X>0$, there exists a weak solution $(u(x,y), v(x,y))$ of problem (\ref{ww71})-(\ref{ww72}).
 If $U_{x}\geq 0$ and $v_{0}\leq 0$, or $U_{x}(x)>0$, then the solution exists for any $X>0$. If we assume, in addition, that
$U_{x}>0, v_{0}(x)\leq 0, u_{0}(y)\equiv U(0)$ for $y\geq y_{1}>0$, then there is a $y_{2}\geq y_{1}$ such that $u(x,y)\equiv U(x)$ for
 $y \geq y_{2}$ and this solution is unique.

\end{Theorem}

  Oleinik and  Samokhin (\cite{OS}) in 1999 considered
 system (\ref{ww71}) in the domain $D =\{0<x<X, 0<y< \infty \}$ by the von Mises variables, with the boundary conditions
 \begin{equation}\left\{
\begin{array}{ll}
u (0,y)=0, \ \ u (x,0)=0, \ \ v (x,0)=v_{0}(x), \\
u (x,y)\rightarrow U (x),  \ \ \text{as}\ \   y\rightarrow+\infty.
\end{array}
 \label{ww73}      \right.   \end{equation}

\begin{Theorem} \label{OS.8.49} (\cite{OS})
 Suppose that $U(x), U_{x}(x), U_{xx}(x), v_{0}(x), v_{0x}(x)$ are continuous for $0\leq x\leq X$, and let $U_{x}>0$. Then problem (\ref{ww71}), (\ref{ww73})
has a weak solution $u(x,y), v(x,y)$ with the following properties: $0\leq u_{y}\leq M_{1}, M_{1}=$const.$>0$; $u_{y}>0$ for $y=0$ and $x>0$; the function
$u_{y}^{n-1}(u_{yy})^{2}$ is summable in any finite portion of the domain $\Omega$; for any $x_{1}>0$ and $0< x_{1}\leq x\leq X$,
the functions $(u_{x}(x,y), v(x,y))$ are continuous at $y=0$ and $v(x,0)=v_{0}(x)$. The solution $(u (x,y), v(x,y))$ with these properties is unique and $u(x,y)\equiv U(x)$ for
$$y\geq \int_{0}^{Cx}\frac{ds}{u(x,s)}, $$
where $C$  is a positive constant depending on $U(x)$ and $v_{0}(x)$.

\end{Theorem}

Hou, Liu, Wang and Wang (\cite{hlww}) in 2003 studied  the stability of boundary layer solutions for a viscous hyperbolic system transformed via a Cole-Hopf transformation
  from a singular chemotactic system modeling the initiation of tumor angiogenesis. They  proved the stability of boundary layer solutions and identify the
   precise structure of boundary layer solutions.

 The prototypical chemotaxis model, known as the Keller-Segel (KS) mode as follows
\begin{eqnarray}
\left\{\begin{array}{l}{u_{t}=\left[D u_{x}-\chi u(\ln c)_{x}\right]_{x}}, \\ {c_{t}=\varepsilon c_{x x}+\mu u c},\end{array}\right.\label{47.1}
\end{eqnarray}
where $u(x,t)$ and $c(x,t)$  denote the cell density and chemical (signal) concentration
at position $x$ and time $t$, respectively. The function $\chi$ is called the chemotactic sensitivity function accounting for the signal response mechanism. $D >0$ and $\varepsilon>0$ are cell and chemical diffusion coefficients, respectively.

The common approach currently used to overcome this singularity is the Cole-Hopf type transformation
\begin{eqnarray}
 v=-\frac{\sqrt{\chi \mu}}{\mu}(\ln c)_{x}=-\frac{\sqrt{\chi \mu}}{\mu} \frac{c_{x}}{c},
\end{eqnarray}
which transforms the model \eqref{47.1} into a nonsingular system of conservation laws
\begin{eqnarray}
\left\{\begin{array}{l}{u_{t}-(u v)_{x}=u_{x x}}, \\ {v_{t}-\left(u-\frac{\varepsilon}{\chi} v^{2}\right)_{x}=\frac{\varepsilon}{D} v_{x x}},
\quad\quad (x,t)\in (0,1)\times (0,+\infty),
\\ {(u, v)(x, 0)=\left(u_{0}, v_{0}\right)(x)}.
\end{array}\right.\label{47.2}
\end{eqnarray}
 For illustration, let's first consider the following initial-boundary value problem of system \eqref{47.2} in an interval $(0,1)$:
\begin{eqnarray}
\left\{\begin{array}{l}{u_{t}-(u v)_{x}=u_{x x}}, \\ {v_{t}-\left(u-\varepsilon |v|^{2}\right)_{x}=\varepsilon v_{x x}},
\quad\quad (x,t)\in (0,1)\times (0,+\infty),
\\ {(u, v)(x, 0)=\left(u_{0}, v_{0}\right)(x)},\\
{u_x|_{x=0,1}=v|_{x=0,1}=0}.\end{array}\right.\label{47.2}
\end{eqnarray}
 Precisely, we postulate that the initial data $(u_0,v_0 )\in H^3(0,1)\times H^3(0,1)$ satisfying
\begin{eqnarray}\left\{\begin{array}{ll}
(u_0,v_0)|_{x=0,1}=(\bar{u},\bar{v}),\\
u_{0x}|_{x=0,1}=0,\\
|(u_0v_0)_x+u_{0xx}|_{x=0,1}=0.
\end{array}\right.\label{47.3}\end{eqnarray}
\begin{Theorem} (\cite{hlww})
Assume that $(u_0,v_0 )\in H^3(0,1)\times H^3(0,1)$ with $u_0\geq 0$
 satisfies the compatibility conditions \eqref{47.3}. Denote by $(u^\varepsilon,v^\varepsilon)$  the unique global solution of problem \eqref{47.2} with $\varepsilon\geq 0$. Then as $\varepsilon\rightarrow0$,  the following asymptotic expansions hold in space $L^\infty([0,1]\times[0,T])$ for any fixed $0<T<\infty$,
 $$u^\varepsilon(x,t)=u^0(x,t)+O(\varepsilon^{1/2}),$$
 $$v^\varepsilon(x,t)=v^0(x,t)+v^{B,0}(\frac{x}{\sqrt\varepsilon},t)+v^{b,0}(\frac{x-1}{\sqrt\varepsilon},t)+O(\varepsilon^{1/2})$$
 where $(u^0,v^0)=(u^{I,0},v^{I,0})$  denotes the outer layer profile and the inner layer profile $(v^{B,0},v^{b,0})$ is given by
 \begin{eqnarray*}\left\{\begin{array}{ll}
 v^{B,0}=\int_0^t\int_{-\infty}^0\frac{1}{\sqrt{\pi(t-s)}}e^{-(\frac{(z-y)^2}{4(t-s)})
 +\bar{u}(t-s)}[\bar{u}(\bar{v}-v^0(0,s))-v^0(0,s)]dyds, \\
 v^{b,0}=\int_0^t\int_{-\infty}^0\frac{1}{\sqrt{\pi(t-s)}}e^{-(\frac{(z-y)^2}{4(t-s)})
 +\bar{u}(t-s)}[\bar{u}(\bar{v}-v^0(1,s))-v^0(1,s)]dyds.
 \end{array}\right.
\end{eqnarray*}
 \end{Theorem}

 The counterpart of the original system \eqref{47.1} in $[0,1]$ corresponding to the initial-boundary value problem of the transformed system \eqref{47.2} reads as follows:
 \begin{eqnarray}
\left\{\begin{array}{l}{u_{t}=(u _{x}-u(lnc)_{x})_x}, \\ {c_{t}=\left(uc-\varepsilon c_{xx}\right)},
\quad\quad (x,t)\in (0,1)\times (0,+\infty),
\\ {(u, v)(x, 0)=\left(u_{0}, v_{0}\right)(x)},\\
{u_x|_{x=0,1}=\bar{u}, \frac{c_x}{c}|_{x=0,1}=-\bar{c}},\ if\ \varepsilon>0,\\
u|_{x=0,1}=\bar{u}.\end{array}\right.\label{47.4}
\end{eqnarray}

 \begin{Theorem}(\cite{hlww})
  Suppose that the initial data $(u_0,ln c_0 )\in H^3\times H^4$ with $u_0\geq 0$, $c_0\geq 0$  and the compatibility conditions \eqref{47.3} with $v_0=-(ln c_0)_x$ and $\bar{v}=\bar{c}$. Let $(u^\varepsilon,v^\varepsilon)$  be the unique global solution of problem  \eqref{47.4} with $\varepsilon\geq 0$. Then for any fixed $0<T<\infty$,  we have in space
 $L^\infty([0,1]\times(0,T))$ that
 $$u^\varepsilon(x,t)=u^0(x,t)+O(\varepsilon^{1/2}),$$
 $$c^\varepsilon(x,t)=c^0(x,t)+O(\varepsilon^{1/2}),$$
 and
 $$c^\varepsilon(x,t)=c^0(x,t)-c^0(x,t)\Big(c^{B,0}(\frac{x}{\sqrt\varepsilon},t)
 +v^{b,0}(\frac{x-1}{\sqrt\varepsilon},t)\Big)+O(\varepsilon^{1/2}). $$

 \end{Theorem}

  Renardy and   Wang (\cite{RW}) in 2014  considered the flow of the Maxwell fluid in the  boundary layer and proved the well-posedness of these equations.
They  started with the upper convected Maxwell model in dimensionless form:
\begin{equation}\left\{
\begin{array}{ll}
  \partial_{t}{\bf v}+( {\bf v}\cdot\nabla ){\bf v}=\frac{1}{R} \nabla\cdot {\bf T}-\nabla p, \\
 {\bf T}_{t}+( {\bf v}\cdot \nabla) {\bf T}- (\nabla  {\bf v}) {\bf T}- {\bf T}(\nabla  {\bf v})^{T}+\frac{1}{W} {\bf T}=\frac{1}{W}(\nabla  {\bf v}+(\nabla {\bf v})^{T}).
\end{array}
\label{w.10}      \right.\end{equation}
Here $v$ denotes the velocity, $T$ the extra stress and $p$ the pressure, where $R$ is the Reynolds number and $W$ is the Weissenberg number measuring the
 elasticity.
Set ${\bf S}={\bf T}+\frac{1}{W^{2}}$, they assumed a flat boundary given by $y=0$,   set
\begin{equation*}
t'=t,\ \ x'=x,\ \ y'=Wy, \ \ u'=u,\ \  v'=Wv,\\
S'_{11}=S_{11}, \ \ S'_{12}=WS_{12}, \ \ S'_{22}=W^{2}S_{22}, \ \  p'=p.
\end{equation*}
Then the momentum equations reduce to the problem
\begin{equation}\left\{
\begin{array}{ll}
 u_{t} +uu_{x}+vu_{y}=E(S_{11}x+ S_{11}y) - p_{x}, \\
p_{y}=0,
\end{array}
\label{w.11}      \right.\end{equation}
here $E=\frac{R}{W}$. Let $$x=x(\xi_{1},\xi_{2},t ),\ \  y=y(\xi_{1},\xi_{2},t ).$$

 \begin{Theorem}(\cite{RW})
 Consider the boundary layer system (\ref{w.11}) on the domain $-\infty<\xi_{1}<\infty$, $0\leq \xi_{2}<\infty$, for $t\geq 0$.
Let $m\geq 1$ be an integer. We assume that  their derivatives up to order $m$ are bounde for $(\xi_{1},\xi_{2},t )\in (-\infty,\infty)\times[0,\infty)\times[0,T]$. Moreover they are all periodic in $\xi_{1}$ with period $1$. Let $\Omega=[0,1]\times[0,\infty),Q=\Omega\times[0,T]$, and denote by $H^{m}_{p}(Q),m\in\mathbb{N}$, the space of all periodic (in $\xi_{1}$) functions
which have $H^{m}$ regularity. If $\Phi\in H^{m}_{p}(Q)$, then there exists some $T'\in (0,T]$ such that a unique solution
of the problem  (\ref{w.11}) exists and satisfies $ \in H^{m}_{p}(Q)$.

\end{Theorem}

 Renardy (\cite{r}) in 2015 proved the local well-posedness of the Prandtl boundary layer equations for the upper convected Maxwell fluid (\ref{w.10}).   He   considered boundary layers in the case where the Reynolds number tends to infinity, but the Weissenberg number remains finite.
As before, Renardy and  Wang (\cite{RW}) in 2014 considered boundary layers for the upper convected Maxwell model in a limit where both the Reynolds and Weissenberg number become infinite.

Let $X(\xi,t)$ denote the outer flow solution, i.e.,  the solution of the equation
$$\frac {\partial ^ { 2 } X } { \partial t ^ { 2 } } = -P(X,t), $$
with initial conditions
 $$X ( \xi , 0 ) = \xi , \quad \frac { \partial X } { \partial t } ( \xi , 0 ) = \lim _{ \eta \rightarrow \infty } u _ { 0 } ( \xi , \eta ).$$
Let $\alpha>0$ be a constant. We shall look for solutions of the form
\begin{equation*}
\begin{aligned}
& x(\xi,\eta,t)=X(\xi,t)+e^{-\alpha \eta}(Z(\xi,t)+w(\xi,\eta,t)), \\
& C (\xi,\eta,t)= \left(  \begin{array} { c } { e^{-2\alpha \eta} (G_{11}(\xi,t)+H_{11} (\xi,\eta,t)\ \ \ \
e^{- \alpha \eta} (G_{12}(\xi,t)+H_{12} (\xi,\eta,t)   ) }\\
( e^{- \alpha \eta} (G_{12}(\xi,t)+H_{12} (\xi,\eta,t)\ \ \ \  (G_{22}(\xi,t)+H_{22} (\xi,\eta,t)   )  )
   \end{array} \right),
\end{aligned}
\end{equation*}
where it is assumed that $w$ and $H$ tend to zero as $\eta\rightarrow\infty$. With this ansatz, we find the
equations
\begin{equation*}\left\{
\begin{aligned} \frac { \partial ^ { 2 } Z } { \partial t ^ { 2 } } & = \alpha ^ { 2 } G _ { 22 } Z - \alpha G _ { 12 } \frac { \partial X } { \partial \xi } - \frac { \partial P } { \partial X } ( X , t ) Z ,\\
\frac { \partial G } { \partial t } & = - \frac { 1 } { W }  { G } + \frac { 1 } { W ^ { 2 } } \left( \begin{array} { c } { \alpha ^ { 2 } Z ^ { 2 } \ \ \  \alpha Z \frac { \partial X } { \partial \xi }   } \\
{ \alpha Z \frac { \partial X } { \partial \xi } \ \ \left( \frac { \partial X } { \partial \xi } \right) ^ { 2 } }   \end{array} \right).\end{aligned}\right.
\end{equation*}

After a series of transformations, then the momentum equations   becomes
\begin{equation}
\begin{aligned} \frac
{ \partial ^ { 2 } w } { \partial t ^ { 2 } }
= & e ^ { - \alpha \eta } \frac { \partial } { \partial \xi } \left( \left( G _ { 11 } + H _ { 11 } \right) \frac { \partial } { \partial \xi } \left( X + e ^ { - \alpha \eta } ( Z + w ) \right) \right) \\
 & + e ^ { - \alpha \eta } \frac { \partial } { \partial \xi } \left( \left( G _ { 12 } + H _ { 12 } \right)
 \left( \frac { \partial } { \partial \eta } - \alpha \right) ( Z + w ) \right) + \left( \frac { \partial } { \partial \eta } - \alpha \right) H _ { 12 } \frac { \partial X } { \partial \xi } \\
  & + \left( \frac { \partial } { \partial \eta } - \alpha \right) \left(
  \left( G _ { 12 } + H _ { 12 } \right) e ^ { - \alpha \eta } \frac { \partial } { \partial \xi } ( Z + w ) \right) \\
 & - \alpha \left( \left( \frac { \partial } { \partial \eta } - \alpha \right) H _ { 22 } \right) Z
  + e ^ { \alpha \eta } \frac { \partial } { \partial \eta } \left( \left( G _ { 22 } + H _ { 22 } \right) e ^ { - \alpha \eta } \left( \frac { \partial w } { \partial \eta } - \alpha w \right) \right) \\
 & + e ^ { \alpha \eta } \left( - P \left( X + e ^ { - \alpha \eta } ( Z + w ) , t \right)
 + P ( X , t ) \right) + \frac { \partial P } { \partial X } ( X , t ) Z,
 \label{w.31}
\end{aligned}
\end{equation}

\begin{equation}
\begin{aligned}
\frac { \partial H _ { 11 } } { \partial t } = & - \frac { 1 } { W } H _ { 11 } - \frac { 2 } { W ^ { 2 } } \alpha Z \left( \frac { \partial w } { \partial \eta } - \alpha w \right) + \frac { 1 } { W ^ { 2 } } \left( \frac { \partial w } { \partial \eta } - \alpha w \right) ^ { 2 } ,\\
\frac { \partial H _ { 12 } } { \partial t } = & - \frac { 1 } { W } H _ { 12 } - \frac { 1 } { W ^ { 2 } } \frac { \partial X } { \partial \xi } \left( \frac { \partial w } { \partial \eta } - \alpha w \right) \\ & - \frac { 1 } { W ^ { 2 } } e ^ { - \alpha \eta } \left( \frac { \partial Z } { \partial \xi } + \frac { \partial w } { \partial \xi } \right) \left( - \alpha Z + \frac { \partial w } { \partial \eta } - \alpha w \right), \\
 \frac { \partial H _ { 22 } } { \partial t } = & - \frac { 1 } { W } H _ { 22 } + \frac { 2 } { W ^ { 2 } } e ^ { - \alpha \eta } \frac { \partial X } { \partial \xi } \left( \frac { \partial Z } { \partial \xi } + \frac { \partial w } { \partial \xi } \right) + \frac { 1 } { W ^ { 2 } } e ^ { - 2 \alpha \eta } \left( \frac { \partial Z } { \partial \xi } + \frac { \partial w } { \partial \xi } \right) ^ { 2 }.
 \label{w.32}
\end{aligned}
\end{equation}

\begin{Theorem}(\cite{r})
Consider   problem (\ref{w.31})-(\ref{w.32}) with initial conditions $\omega(\xi,\eta,0)=0$, \\
 $\omega_{t}(\xi,\eta,0)= w _ { 1 } ( \xi , \eta )$,
$H ( \xi , \eta , 0 ) =  H  _ { 0 } ( \xi , \eta )$. Assume that the functions $P, X, G$ and $G$ are smooth
and periodic in $X$ or, respectively, $\xi$ with period $L$. Assume, moreover, that  $\omega_{1}$, $H_{0}\in \overline{\omega}_{n-1}$,
where $n\geq 5$, and that appropriate compatibility conditions between the initial and boundary data for w hold. Finally, assume that $G(\xi,0)+H(\xi,\eta,0)$
 is strictly and uniformly positive definite. Then, for some $T>0$, there exists a unique solution with the regularity
 \begin{eqnarray*}
\left\{\begin{array}{l}
{ w \in L ^ { \infty } \left( [ 0 , T ] , V _ { n } \right) \cap W ^ { 1 , \infty }
\left( [ 0 , T ] , W _ { n - 1 } \right) \cap \bigcap\limits_ { k = 2 } ^ { n } W ^ { k , \infty } \left( [ 0 , T ] , H _ { n - k } \right) }, \\ { { H } \in L ^ { \infty } \left( [ 0 , T ] , W _ { n - 1 } \right) \cap \bigcap \limits_ { k = 1 } ^ { n - 1 } W ^ { k , \infty } \left( [ 0 , T ] , H _ { n - 1 - k } \right) }.
 \end{array}\right.
\end{eqnarray*}
\end{Theorem}

Li and Wang (\cite{LW3}) in 2018 investigated the global solutions to  initial-boundary value problem for the following nonlinear evolution system (\ref{61.1}) on the domain $[0,1]\times[0,\infty)$ with damping and diffusion and gave the convergence rates and boundary layer thickness:
\begin{eqnarray}\left\{\begin{array}{ll}
\phi^\alpha_t=-(k-1)\alpha\phi^\alpha-k\alpha\theta^\alpha_x+\alpha\phi^\alpha_{xx}, 0\leq x\leq 1,\ t>0,\\
\theta^\alpha_t=-(1-\beta)\theta^\alpha+\mu\alpha\phi^\alpha_x+2\phi^\alpha\theta^\alpha
+\beta\theta^\alpha_{xx},\\
(\phi^\alpha,\theta^\alpha)(x,0)=(\phi^\alpha_0,\theta^\alpha_0)(x),\ 0\leq x\leq 1,\\
(\phi^\alpha_x,\theta^\alpha)(0,t)=(\phi^\alpha_x,\theta^\alpha)(1,t)=(0,0),\ t\geq0,
\end{array}\right.\label{61.1}\end{eqnarray}
where $k,\alpha,\beta$ and $\mu$ are positive constants with $k>1$ and $0<\beta<1$. The limit problem of the vanishing parameter $\alpha\rightarrow 0^+$ in \eqref{61.1} is the following
\begin{eqnarray}\left\{\begin{array}{ll}
\phi^0_t=0,\\
\theta^0_t=-(1-\beta)\theta^0+2\phi^0\theta^0_{xx},\ 0<x<1,\ t>0,\\
(\phi^0,\theta^0)(x,0)=(\phi_1(x),\phi_2(x)) ,\ 0\leq x\leq 1,\\
\theta^0(0,t)=\theta^0(1,t)=0,\ t\geq0.
\end{array}\right.\label{61.2}\end{eqnarray}
We first  show that the initial-boundary
value problem \eqref{61.1} admits a unique global smooth solution $(\phi^\alpha,\theta^\alpha)$, which is stated in the next theorem.
\begin{Theorem}( \cite{LW3} )
Suppose that the initial data satisfy the conditions:
$(\phi^0,\theta^0)\in H^2$, $(\phi^0,\theta^0)(0,t)=(\phi^0,\theta^0)(1,t)=(0,0)$ and
$\|(\phi^0,\theta^0)\|_{H^1}$  is sufficiently small, then there exists a unique
solution $(\phi^\alpha,\theta^\alpha)$ to the initial boundary value problem \eqref{61.1} satisfying
$$\theta^\alpha\in L^\infty(0,T;H^2)\cap L^2(0,T;H^2),\ \phi^\alpha\in L^\infty(0,T;H^2)\cap L^2(0,T;H^1),$$
where the norms are all uniform in $\alpha$.
\end{Theorem}
The second result in the following shows that the initial-boundary
value problem \eqref{61.2} admits a unique global smooth solution $(\phi^0,\theta^0)$.
\begin{Theorem}( \cite{LW3} )
Suppose that the initial data satisfy the conditions:
$(\phi_1(x),\phi_2(x))\in H^2$, $\phi_1(0,t)=0, ~\phi_2(1,t)\leq\varepsilon$ and
$\varepsilon$   is sufficiently small, then there exists a unique
solution $(\phi^0,\theta^0)$ to the initial boundary value problem \eqref{61.2} satisfying
$$\phi^0\in L^\infty(0,T;H^2),\ \theta^0\in L^\infty(0,T;H^2)\cap L^2(0,T;H^2),\ \theta^0_t\in L^\infty(0,T;H^2)\cap L^2(0,T;H^1).$$
\end{Theorem}
Furthermore, Li and Wang (\cite{LW3}) gave the convergence rates and boundary layer thickness.

\begin{Theorem} ( \cite{LW3} )
Under the assumptions of the above two theorems, suppose that
$$\|\phi_0(x)-\phi_1(x)\|_{H^1}^2\leq C\alpha,\quad  \|\theta_0(x)-\phi_2(x)\|_{H^1}^2\leq C\alpha,$$
then any function $\delta(\alpha)$ satisfying the condition  $\delta(\alpha)\rightarrow0$
and $\frac{\alpha}{\delta(\alpha)}\rightarrow0$ as $\alpha\rightarrow0^+$ is a BL-thickness such that
$$\|\theta^\alpha-\theta^0\|_{L^\infty(0,T;L^\infty[0,1]})\leq C\alpha,\quad
|\phi^\alpha-\phi^0\|_{L^\infty(0,T;L^\infty[\delta,1-\delta]})\leq C\frac{\alpha}{\delta(\alpha)},$$
$$\lim\inf\limits_{\alpha\rightarrow0}\|\phi^\alpha-\phi^0\|_{L^\infty(0,T;L^\infty[0,1]} )>0.$$
Consequently,
$$\lim\limits_{\alpha\rightarrow0}\|\theta^\alpha-\theta^0\|_{L^\infty(0,T;L^\infty[0,1]} )=0,\quad \lim\limits_{\alpha\rightarrow0}\|\phi^\alpha-\phi^0\|_{L^\infty(0,T;L^\infty[\delta,1-\delta]} )=0. $$
\end{Theorem}

 Oleinik and  Samokhin (\cite{OS}) in 1999  considered
  the nonstationary symmetric flow of an electrically conducting pseudo-plastic fluid within the boundary layer,
\begin{equation}\left\{
\begin{array}{ll}
\frac{\partial u }{\partial t}+  u \frac{\partial u }{\partial x}+v \frac{\partial u }{\partial y}
=  \nu \frac{\partial}{\partial y }\left(|\frac{\partial u }{\partial y}|^{n-1}\frac{\partial u }{\partial y} \right)+U(x)\frac{d U }{d x}
+d^{2} (U-u) +\frac{\partial U }{\partial t}, \ \  0<n<1, \\
\frac{\partial u }{\partial x}+\frac{\partial v }{\partial y} =0.
\end{array}
 \label{ww84}         \right.\end{equation}
Here $d^{2} =\sigma B^{2}(x) $ in the domain $D =\{0<t<\infty, 0<x<X, 0<y< \infty \}$, with the boundary conditions
 \begin{equation}
\begin{array}{ll}
u (0,x,y)=u_{0}(x,y), \ \ u (t,0,y)=0, \ \ u (x,0)=0, \ \ v (x,0)=v_{0}(t,x), \\
u (t,x,y)\rightarrow U (t,x),  \ \ \text{as}\ \   y\rightarrow \infty.
\end{array}
 \label{ww85}         \end{equation}
Assume that
$$U(t,x) =x^{m}V(t,x), \ \  v_{0}(t,x)=x^{\frac{2mn-n-m}{n+1}}v_{1}(t,x),\ \    d(t,x)=x^{\frac{m-1}{2}}W(t,x), $$
where the functions $V,W,v_{1}$ are bounded in $D$, $V(t,0)=a=$const.$>0$, $v_{1}(t,0)=b=$const., $W(t,0)=d_{0}=$const.. The assumptions about the
smoothness of $V,W,v_{1}$, and the properties of their derivatives will be stated later. In  (\ref{ww84}), introduce  new independent variables
$$\tau=t, \ \  \xi=x, \ \  \eta=\frac{u}{U}, \ \  w(\tau,\xi, \eta)=\frac{|u_{y}|^{n-1}u_{y}}{x^{\frac{2mn-n-m}{n+1}}U}. $$

\begin{Theorem} \label{OS.9.56} (\cite{OS})
Suppose that $ U(t,x) =x^{m}(a+xa_{1}(t,x))$, $d(t,x)=x^{\frac{m-1}{2}}(d_{0}+xd_{1}(t,x))$, $ v_{0}(t,x)=x^{\frac{2mn-n-m}{n+1}}( b+xb_{1}(t,x ))$,
$m>0$, $a_{1}(t,x),d_{1}(t,x)$ have bounded second order derivatives; $b_{1}(t,x)$ has bounded first order derivatives; $d_{0}^{2}+(3m-1)na/(n+1)>0$;
the functions $u_{0}(x,y)$ and $v_{0}(t,x)$ satisfy
$$w_{0}(\xi,\eta)=|u_{0y}|^{n-1}u_{0y}\left(\xi^{\frac{2mn-n-m}{n+1}} U(0,\xi)\right)^{-1} .$$
Then problem  (\ref{ww84})-(\ref{ww85}) in $D$, with $X$ depending on $U,v_{0}, d,\nu, m,n$, admits a solution $u(t,x,y)$, $v(t,x,y)$ with the following
 properties: $u/U$ is continuous in $\overline{D}$; $v$ is continuous in $\overline{D}$ with respect to $y$ and bounded for bounded $y$,
the derivatives $u_{y},u_{yy}, u_{x},u_{t}, v_{y}$ are continuous in $D$ with respect to $y$; equations (\ref{ww84}) hold for $u,v$ almost everywhere in $D$;
$u,v$ satisfy the boundary conditions (\ref{ww85}). Moreover, $u>0$ for $y>0$ and $x>0$; $u_{y}>0$ for $y\leq 0$; $u_{y}\rightarrow 0$ as
$y\rightarrow \infty$; $u_{yyy}$ is bounded in $D$; the following inequalities hold:
\begin{eqnarray*}
&& x^{\frac{2n-m-n}{n+1}}U(t,x)Y(\frac{u}{U})e^{-K_{15}x} \leq u_{y}^{n}(t,x,y)\leq x^{\frac{2mn-m-n}{n+1}}U(t,x)Y(\frac{u}{U})e^{ K_{16}x}  , \\
&&  -K_{22}\leq \frac{n(n-2)(u_{yy})^{2}+ nu_{y}u_{yyy}}{x^{\frac{2mn-m-n}{n}} U^{\frac{1-n}{n}} u_{y}^{3-n}}\leq -K_{23} , \\
&& \left(M_{8}^{\frac{1}{n}} \frac{1-n}{n+1} e^{\frac{K_{14}}{n}}x^{\frac{2mn-m-n}{n(n+1)}}U^{\frac{ 1-n}{n }}(t,x)y +1\right)^{\frac{n+1}{n-1}}\leq 1-\frac{u}{U}
\leq \left(M_{7}^{\frac{1}{n}} \frac{1-n}{n+1} e^{-\frac{K_{15}}{n}}x^{\frac{2mn-m-n}{n(n+1)}}U^{\frac{ 1-n}{n }}(t,x)y +1\right)^{\frac{n+1}{n-1}}.
\end{eqnarray*}
The solution of problem (\ref{ww84})-(\ref{ww85}) with the above properties is unique.

\end{Theorem}

\chapter{Global Well-posedness of Solutions to the 2D  Prandtl-Hartmann Equations in An Analytic Framework}

In this chapter, we shall consider the 2D Prandtl-Hartmann  equations on the half plane and prove the  global existence and uniqueness of solutions to the 2D  Prandtl-Hartmann  equations by using the classical energy methods in an analytic framework. We prove that the lifespan of the  solutions to the 2D  Prandtl-Hartmann  equations can be extended up to $T_\varepsilon$ (see Theorem 2.2.1) when  the strength of the perturbation is of the order of $\varepsilon$. The difficulty of solving the Prandtl-Hartmann equations in the analytic framework is the loss of $x$-derivative in the term $v\partial_yu$. To overcome this difficulty, we show the Gaussian weighted Poincar$\acute{e}$ inequality (see Lemma 2.2.3). Compared to the existence and uniqueness of solutions to the classical Prandtl equations where the monotonicity condition of the tangential velocity, playing a key role, is required there, while in this chapter it is not needed for the 2D Prandtl-Hartmann equations in an analytic framework. Besides, compared with the existence and uniqueness of solutions to the 2D MHD boundary layer where the initial tangential magnetic field has a lower bound, also playing an important role there, which is also not needed in this chapter for 2D Prandtl-Hartmann  equations  in an analytic framework.
The content of this chapter is chosen from \cite{dongqin}.
\section{Introduction}
\setcounter{equation}{0}
The following  Prandtl-Hartmann  equations were derived in \cite{gvp} by the classical 2D incompressible MHD system
 under a transverse magnetic field on the boundary for the parameters of the system in the mixed Prandtl-Hartmann regime. In this chapter, we shall
investigate the global existence and uniqueness of solutions to the following initial boundary value problem for the 2D  Prandtl-Hartmann  system in the upper half plane $\mathbb{R}_+^2=\{(x,y):x\in\mathbb{R},\ y\in\mathbb{R}_+\}$,  which reads as
\begin{eqnarray}\left\{\begin{array}{ll}
\partial_t u_1+u_1\partial_xu_1+u_2\partial_yu_1-\partial^2_yu_1+u_1=0,\\
\partial_xu_1+\partial_yu_2=0,
\end{array}\right.\label{2.1.1}\end{eqnarray}
where the velocity field $(u_1(t,x,y),u_2(t,x,y))$ of the boundary layer is a unknown vector function.

The system \eqref{2.1.1} subjects to the  initial data and boundary conditions
\begin{eqnarray}\left\{\begin{array}{ll}
u_1(t,x,y)|_{t=0}=u_{10}(x,y),\\
u_1|_{y=0}=u_2|_{y=0}=0.
\end{array}\right.\label{2.1.2}\end{eqnarray}

The far field is represented by a positive constant $\bar{u}$,
\bqa \lim\limits_{y\rightarrow+\infty}u_1(t,x,y)=\bar{u}.\label{2.1.3}\eqa
\section{Global Wellposedness of Solutions}
\setcounter{equation}{0}

As in \cite{iv,xy}, we use the following Gaussian weighted function $\theta_\alpha$ to introduce the function space of the solutions,
$$\theta_\alpha=\exp{\big(\frac{\alpha z^2}{4}\big)}\ \ \text{with}\ \ z=\frac{ y}{\sqrt{\langle t\rangle}},\ \langle t\rangle=(1+t),\ \alpha\in[1/4,1/2],$$
which, combined with

$$M_m=\frac{\sqrt{(m+1)}}{m!},$$
defines the Sobolev weighted semi-norms by
\begin{eqnarray}\left\{\begin{array}{ll}
X_m=X_m(f,\tau)=\|\theta_\alpha\partial_x^mf\|_{L^2(\mathbb{R}_+^2)}\tau^mM_m,\\
D_m=D_m(f,\tau)=\|\theta_\alpha\partial_y\partial_x^mf\|_{L^2(\mathbb{R}_+^2)}\tau^mM_m,\\
Y_m=Y_m(f,\tau)=\|\theta_\alpha\partial_x^mf\|_{L^2(\mathbb{R}_+^2)}\tau^{m-1}mM_m.
\end{array}\right.\label{2.2.1}\end{eqnarray}

To define the following analytic function space  in the tangential variable $x$ and weighted Sobolev space
in the normal variable $y$  by
$$X_{\tau,\alpha}=\Big\{f(t,x,y)\in L^2(\mathbb{R}_+^2;\theta_\alpha dxdy):\|f\|_{X_{\tau,\alpha}}<\infty\Big\},$$
with $\tau>0$ and the norms
\begin{eqnarray}\left\{\begin{array}{ll}
\|f\|_{X_{\tau,\alpha}}=\sum\limits_{m\geq 0}X_m(f,\tau),\\
\|f\|_{D_{\tau,\alpha}}=\sum\limits_{m\geq 0}D_m(f,\tau),\\
\|f\|_{Y_{\tau,\alpha}}=\sum\limits_{m\geq 1}Y_m(f,\tau).
\end{array}\right.\non\end{eqnarray}
\begin{Theorem}
Let $\frac{1}{4}\leq \alpha\leq \frac{1}{2}$ and assume that the initial datum $u_1(0)$ satisfies
\bqa \|u_1(0)-\bar{u}\|_{X_{\tau_0,\alpha}}\leq \varepsilon,\label{2.2.2}\eqa
for some small constant $0<\varepsilon\ll 1$ and assume there exists an analyticity radius $\tau_0$ such that
\bqa\max\{\frac{7\varepsilon}{8C_1},128\sqrt2C\varepsilon,\frac{256\sqrt2CC_2^{\frac{3}{4}}}{\alpha}\}\leq\tau_0^{\frac{3}{2}}.\label{2.2.3}\eqa
Then there exists a unique solution $(u_1,u_2)$
to the damping Prandtl equations \eqref{2.1.1}-\eqref{2.1.3} such that
\bqa \|u\|_{X_{\tau,\alpha}}\leq e^{-t}\varepsilon,\label{2.2.4}\eqa
with analyticity radius $\tau$ larger than $\frac{\tau_0}{4}$ in the time interval $[0,T_\varepsilon]$,
where $u=u_1-\bar{u}$, $T_\varepsilon=C_2\varepsilon^{-\frac{4}{3}}-1$ and $\bar{u}$ is given by \eqref{2.1.3}, $C_1=\frac{3C(\alpha+2)}{2\alpha}$ and $C_2=\Big(\frac{7}{8C_1}\tau_0^{\frac{3}{2}}\Big)^\frac{4}{3}.$
\end{Theorem}
The above theorem is the main part of this chapter, and its proof follows from the discussion in Sections 2.2.1 and 2.2.2.

Before  proving the theorem, we first list some important inequalities as follows. The following main estimates on the functions in the norms defined in
the previous section will be used frequently, whose proofs can be found in \cite{h,xy,iv}.
\begin{Lemma}(Agmon inequality \cite{iv}) Let $f\in H^1(\mathbb{R}_+^2)$, then
\bqa\|f\|_{L_x^\infty L^2_y}\leq C\|f\|^{\frac{1}{2}}_{L^2(\mathbb{R}_+^2)}\|f\|^{\frac{1}{2}}_{H^1(\mathbb{R}_+^2)}.\non\eqa
\end{Lemma}
\begin{Lemma}(\cite{h,iv})
Let $f$ be a function such that $f|_{y=0}(or\ \partial_yf|_{y=0})$ and $f|_{y=\infty}=0$. Then,
for  $m\geq0$ and $t\geq0$, there holds that
\begin{eqnarray}\left\{\begin{array}{ll}
 \frac{\alpha}{\langle t\rangle}\|\theta_\alpha\partial_x^mf\|^2_{L^2_y}\leq
\|\theta_\alpha\partial_y\partial_x^mf\|^2_{L^2_y},\\
\sum\limits_{m\geq0}\frac{\|\theta_\alpha\partial_y\partial_x^mf\|^2_{L^2(\mathbb{R}_+^2)}}{\|\theta_\alpha\partial_x^mf\|_{L^2(\mathbb{R}_+^2)}}\tau^mM_m\geq
\frac{\alpha}{2\langle t\rangle^{\frac{1}{2}}}\|f\|_{D_{\tau,\alpha}}+\frac{\alpha(2-\sqrt\alpha)}{\langle t\rangle}\|f\|_{X_{\tau,\alpha}},
\end{array}\right.\label{2.2.5}\end{eqnarray}
with $\alpha\in[1/4,1/2]$, $\langle t\rangle=1+t$.
\end{Lemma}
\textbf{Proof}. For the convenience of the readers, we will prove the second inequality as follows. Actually, by the first inequality, we have
\bqa \frac{\|\theta_\alpha\partial_y\partial_x^mf\|^2_{L^2(\mathbb{R}_+^2)}}{\|\theta_\alpha\partial_x^mf\|_{L^2(\mathbb{R}_+^2)}}&\geq&
\frac{\sqrt\alpha}{4}\frac{\|\theta_\alpha\partial_y\partial_x^mf\|^2_{L^2(\mathbb{R}_+^2)}}{\|\theta_\alpha\partial_x^mf\|_{L^2(\mathbb{R}_+^2)}}+
\frac{2-\frac{\sqrt\alpha}{2}}{2}\frac{\sqrt\alpha}{\sqrt{\langle t\rangle}}\|\theta_\alpha\partial_y\partial_x^mf\|_{L^2(\mathbb{R}_+^2)}\non\\
&\geq&
\frac{\sqrt\alpha}{4}\frac{\|\theta_\alpha\partial_y\partial_x^mf\|^2_{L^2(\mathbb{R}_+^2)}}{\|\theta_\alpha\partial_x^mf\|_{L^2(\mathbb{R}_+^2)}}+
\frac{\alpha}{4\sqrt{\langle t\rangle}}\|\theta_\alpha\partial_y\partial_x^mf\|_{L^2(\mathbb{R}_+^2)}
+
\frac{\alpha(2-\sqrt\alpha)}{2\langle t\rangle}\|\theta_\alpha\partial_x^mf\|_{L^2(\mathbb{R}_+^2)}\non\\
&\geq&\frac{\alpha}{2\sqrt{\langle t\rangle}}\|\theta_\alpha\partial_y\partial_x^mf\|_{L^2(\mathbb{R}_+^2)}
+
\frac{\alpha(2-\sqrt\alpha)}{2\langle t\rangle}\|\theta_\alpha\partial_x^mf\|_{L^2(\mathbb{R}_+^2)}.\non\eqa
\hfill $\Box$

\begin{Lemma}(\cite{xy}) Let $f$ be a function such that $f|_{y=0}$. Then, it holds that
\begin{eqnarray}\left\{\begin{array}{ll}
\|f\|_{L^2_xL^\infty_y}\leq C\|\theta_\alpha f\|^{\frac{1}{2}}_{L^2(\mathbb{R}_+^2)}\|\theta_\alpha \partial_yf\|^{\frac{1}{2}}_{L^2(\mathbb{R}_+^2)},\\
\|f\|_{L_{xy}^\infty}\leq C\|\theta_\alpha f\|^{\frac{1}{4}}_{L^2(\mathbb{R}_+^2)}\|\theta_\alpha\partial_x f\|^{\frac{1}{4}}_{L^2(\mathbb{R}_+^2)}\|\theta_\alpha \partial_yf\|^{\frac{1}{4}}_{L^2(\mathbb{R}_+^2)}\|\theta_\alpha \partial_x\partial_yf\|^{\frac{1}{4}}_{L^2(\mathbb{R}_+^2)}.
\end{array}\right.\label{2.2.6}\end{eqnarray}
\end{Lemma}

\subsection{Uniform Estimates}

In this subsection, we will prove the existence of  solutions for the 2D Prandtl-Hartmann  regime in an analytic framework defined in Section 2.2  by the standard energy method.

Now, we are ready to establish  uniform estimates for the solutions to system \eqref{2.1.1}-\eqref{2.1.3}. For our purpose, we first rewrite
the solutions $(u_1,u_2)$ to the problem \eqref{2.1.1}-\eqref{2.1.3} as a perturbation $(u,v)$ around the point $(\bar{u},0)$ by standing for
\begin{eqnarray}\left\{\begin{array}{ll}
u_1=u+\bar{u},\\
u_2=v.
\end{array}\right.\label{2.3.1}\end{eqnarray}
Then system \eqref{2.1.1} becomes
\begin{eqnarray}\left\{\begin{array}{ll}
\partial_t u+(\bar{u}+u)\partial_xu+v\partial_y(u+\bar{u})-\partial^2_yu+u+\bar{u}=0,\\
\partial_xu+\partial_yv=0,
\end{array}\right.\label{2.3.2}\end{eqnarray}
which subjects to the initial datum and boundary data
\begin{eqnarray}\left\{\begin{array}{ll}
u(t,x,y)|_{t=0}=u_0(x,y)-\bar{u},\\
u|_{y=0}=v|_{y=0}=0,
\end{array}\right.\label{2.3.3}\end{eqnarray}
and the corresponding far field condition
\bqa \lim\limits_{y\rightarrow+\infty}u=0.\label{2.3.4}\eqa
Next, we will establish the estimates of solutions to the system \eqref{2.3.2}-\eqref{2.3.4}. For $m\geq 0$, applying the operator $\partial^m_x$ on $\eqref{2.3.2}_1$, multiplying
the resulting equation by $\theta^2_\alpha\partial_x^mu$ and integrating it by parts over $\mathbb{R}_+^2$,  we derive that
\bqa\int_{\mathbb{R}_+^2}\partial^m_x\big(\partial_t u+(u+\bar{u})\partial_xu+v\partial_yu-\partial^2_yu+u+\bar{u}\big)
\theta^2_\alpha\partial_x^mudxdy=0.\label{2.3.5}\eqa
We now deal with each term in \eqref{2.3.5} as follows. For the first term, the second term and the fourth term, we have, respectively
\bqa&&\int_{\mathbb{R}_+^2}\partial^m_x\partial_t u\theta^2_\alpha\partial_x^mudxdy=
\frac{1}{2}\int_{\mathbb{R}_+^2}\partial_t(\partial^m_xu)^2\theta^2_\alpha dxdy\non\\
&=&\frac{1}{2}\frac{d}{dt}\int_{\mathbb{R}_+^2}(\partial^m_xu)^2\theta^2_\alpha dxdy-\int_{\mathbb{R}_+^2}(\partial^m_xu)^2\theta_\alpha\frac{d}{dt}\theta_\alpha dxdy\non\\
&=&\frac{1}{2}\frac{d}{dt}\|\theta_\alpha\partial^m_xu\|^2_{L^2(\mathbb{R}_+^2)}+\frac{\alpha}{4\langle t\rangle}\|\theta_\alpha z\partial^m_xu\|^2_{L^2(\mathbb{R}_+^2)},\label{2.3.6}\eqa
\bqa \int_{\mathbb{R}_+^2}\partial^m_x(u+\bar{u})\theta^2_\alpha\partial_x^mudxdy=\|\theta_\alpha\partial^m_xu\|^2_{L^2(\mathbb{R}_+^2)},\label{2.3.7}\eqa

and
\bqa &&-\int_{\mathbb{R}_+^2}\partial^2_y\partial^m_x u\theta^2_\alpha\partial_x^mudxdy\non\\
&=&
\|\theta_\alpha\partial^m_x\partial_yu\|^2_{L^2(\mathbb{R}_+^2)}+\int_{\mathbb{R}_+^2}\partial_y\partial^m_x u\partial_y(\theta^2_\alpha)\partial_x^mudxdy\non\\
&=&\|\theta_\alpha\partial^m_x\partial_yu\|^2_{L^2(\mathbb{R}_+^2)}-\frac{1}{2}\int_{\mathbb{R}_+^2}\partial_y(\partial^m_x u)^2\partial_y(\theta^2_\alpha)dxdy\non\\
&=&\|\theta_\alpha\partial^m_x\partial_yu\|^2_{L^2(\mathbb{R}_+^2)}-\frac{\alpha}{2{\langle t\rangle}}\|\theta_\alpha\partial^m_xu\|^2_{L^2(\mathbb{R}_+^2)}
-\frac{\alpha^2}{2\langle t\rangle}\|\theta_\alpha z\partial^m_xu\|^2_{L^2(\mathbb{R}_+^2)},\label{2.3.8}\eqa
where we have used the fact $\partial^2_y(\theta^2_\alpha)=\frac{\alpha}{{\langle t\rangle}}\theta^2_\alpha+\frac{\alpha^2}{\langle t\rangle}z^2\theta^2_\alpha$ and boundary conditions $\partial_x^mu|_{y=0}=0$.
Next, we build the estimate of the nonlinear terms in \eqref{2.3.5} as follows
\bqa&&R_1=\int_{\mathbb{R}_+^2}\partial^m_x ((u+\bar{u})\partial_xu)\theta^2_\alpha\partial_x^mudxdy=\sum\limits_{i=0}^m\binom mi\int_{\mathbb{R}_+^2}\partial^{m-i}_x u\partial^{i+1}_xu\theta^2_\alpha\partial_x^mudxdy\non\\
&&\qquad\leq \sum\limits_{i=0}^{[m/2]}\binom mi\|\partial^{m-i}_x u\|_{L^2_xL_y^\infty}\|\theta_\alpha\partial^{i+1}_xu\|_{L^\infty_xL_y^2}\|\theta_\alpha\partial_x^mu\|_{L^2(\mathbb{R}_+^2)}\non\\
&&\qquad\qquad+\sum\limits_{i=[m/2]+1}^{m}\binom mi\|\partial^{m-i}_x u\|_{L_{xy}^\infty}\|\theta_\alpha\partial^{i+1}_xu\|_{L^2(\mathbb{R}_+^2)}\|\theta_\alpha\partial_x^mu\|_{L^2(\mathbb{R}_+^2)}.\label{2.3.9}\eqa
For $0\leq i\leq [m/2]$, using  Lemma 2.2.4 and Lemma 2.2.2 (the Agmon inequality) in $x$, we deduce
\bqa \|\partial^{m-i}_x u\|_{L^2_xL_y^\infty}\leq C\|\theta_\alpha \partial^{m-i}_x u\|^{\frac{1}{2}}_{L^2(\mathbb{R}_+^2)}\|\theta_\alpha \partial_y\partial^{m-i}_x u\|^{\frac{1}{2}}_{L^2(\mathbb{R}_+^2)},\label{2.3.10}\eqa

and

\bqa \|\theta_\alpha\partial^{i+1}_xu\|_{L^\infty_xL_y^2}\leq C\|\theta_\alpha\partial^{i+1}_xu\|^{\frac{1}{2}}_{L^2(\mathbb{R}_+^2)}\|\theta_\alpha\partial^{i+2}_xu\|^{\frac{1}{2}}_{L^2(\mathbb{R}_+^2)}.\label{2.3.11}\eqa
For $ [m/2]\leq i\leq m$, by using Lemma 2.2.4, we arrive at
\bqa \|\partial^{m-i}_x u\|_{L_{xy}^\infty}
&\leq& C\|\theta_\alpha \partial^{m-i}_x u\|^{\frac{1}{4}}_{L^2(\mathbb{R}_+^2)}\|\theta_\alpha \partial^{m-i+1}_x u\|^{\frac{1}{4}}_{L^2(\mathbb{R}_+^2)}\non\\
&&\times\|\theta_\alpha \partial_y\partial^{m-i}_x u\|^{\frac{1}{4}}_{L^2(\mathbb{R}_+^2)}\|\theta_\alpha \partial^{m-i+1}_x\partial_y u\|^{\frac{1}{4}}_{L^2(\mathbb{R}_+^2)}.\label{2.3.12}\eqa
Accordingly,
\bqa \frac{|R_1|\tau^mM_m}{\|\theta_\alpha\partial_x^mu\|_{L^2(\mathbb{R}_+^2)}}
&\leq& \frac{C}{\tau(t)^{\frac{1}{2}}}\Big(\sum\limits_{i=0}^{[m/2]}\binom miX^{\frac{1}{2}}_{m-i}D^{\frac{1}{2}}_{m-i}Y^{\frac{1}{2}}_{i+1}Y^{\frac{1}{2}}_{i+2}\non\\
&&+\sum\limits_{i=[m/2]+1}^{m}\binom miX^{\frac{1}{4}}_{m-i}X^{\frac{1}{4}}_{m-i+1}D^{\frac{1}{4}}_{m-i}D^{\frac{1}{4}}_{m-i+1}Y_{i+1}\Big).\label{2.3.13}\eqa
Similarly, we also have
\bqa&&R_2=\int_{\mathbb{R}_+^2}\partial^m_x (v\partial_yu)\theta^2_\alpha\partial_x^mudxdy=\sum\limits_{i=0}^m\binom mi\int_{\mathbb{R}_+^2}\partial^{m-i}_x v\partial^{i}_x\partial_yu\theta^2_\alpha\partial_x^mudxdy\non\\
&&\qquad\leq \sum\limits_{i=0}^{[m/2]}\binom mi\|\partial^{m-i}_x v\|_{L^2_xL_y^\infty}\|\theta_\alpha\partial^{i}_x\partial_yu\|_{L^\infty_xL_y^2}\|\theta_\alpha\partial_x^mu\|_{L^2(\mathbb{R}_+^2)}\non\\
&&\qquad\qquad+\sum\limits_{i=[m/2]+1}^{m}\binom mi\|\partial^{m-i}_x v\|_{L_{xy}^\infty}\|\theta_\alpha\partial^{i}_x\partial_yu\|_{L^2(\mathbb{R}_+^2)}\|\theta_\alpha\partial_x^mu\|_{L^2(\mathbb{R}_+^2)}.\label{2.3.14}\eqa
For $0\leq i\leq [m/2]$, using  Lemmas 2.2.3-2.2.4 and Lemma 2.2.2 (the Agmon inequality) in $x$, we attain
\bqa \|\partial^{m-i}_x v\|_{L^2_xL_y^\infty}&\leq& C\|\theta_\alpha \partial^{m-i}_x v\|^{\frac{1}{2}}_{L^2(\mathbb{R}_+^2)}\|\theta_\alpha \partial_y\partial^{m-i}_x v\|^{\frac{1}{2}}_{L^2(\mathbb{R}_+^2)}\non\\
&\leq&C \langle t\rangle^{\frac{1}{4}}\|\theta_\alpha \partial^{m-i+1}_x u\|^{\frac{1}{2}}_{L^2(\mathbb{R}_+^2)}\|\theta_\alpha\partial^{m-i+1}_x u\|^{\frac{1}{2}}_{L^2(\mathbb{R}_+^2)}\non\\
&=& C\langle t\rangle^{\frac{1}{4}}\|\theta_\alpha\partial^{m-i+1}_x u\|_{L^2(\mathbb{R}_+^2)},\label{2.3.15}\eqa
and

\bqa \|\theta_\alpha\partial^{i}_x\partial_yu\|_{L^\infty_xL_y^2}\leq C \|\theta_\alpha\partial^{i}_x\partial_yu\|^{\frac{1}{2}}_{L^2(\mathbb{R}_+^2)}\|\theta_\alpha\partial^{i+1}_x\partial_yu\|^{\frac{1}{2}}_{L^2(\mathbb{R}_+^2)}.\label{2.3.16}\eqa
For $ [m/2]\leq i\leq m$, using Lemmas 2.2.3-2.2.4, we find
\bqa \|\partial^{m-i}_x v\|_{L_{xy}^\infty}&\leq& C\|\theta_\alpha \partial^{m-i}_x v\|^{\frac{1}{4}}_{L^2(\mathbb{R}_+^2)}\|\theta_\alpha \partial^{m-i+1}_x v\|^{\frac{1}{4}}_{L^2(\mathbb{R}_+^2)}\|\theta_\alpha \partial_y\partial^{m-i}_x v\|^{\frac{1}{4}}_{L^2(\mathbb{R}_+^2)}\|\theta_\alpha \partial^{m-i+1}_x\partial_y v\|^{\frac{1}{4}}_{L^2(\mathbb{R}_+^2)}\non\\
&\leq&C\langle t\rangle^{\frac{1}{4}}\|\theta_\alpha \partial^{m-i+1}_x u\|^{\frac{1}{2}}_{L^2(\mathbb{R}_+^2)}\|\theta_\alpha\partial^{m-i+2}_x u\|^{\frac{1}{2}}_{L^2(\mathbb{R}_+^2)}.\label{2.3.17}\eqa
Hence,
\bqa \frac{|R_2|\tau^mM_m}{\|\theta_\alpha\partial_x^mu\|_{L^2(\mathbb{R}_+^2)}}
&\leq& \frac{C}{\tau(t)^{\frac{1}{2}}}\Big(\sum\limits_{i=0}^{[m/2]}\binom mi\langle t\rangle^{\frac{1}{4}}Y_{m-i+1}D^{\frac{1}{2}}_{i}D^{\frac{1}{2}}_{i+1}\non\\
&&+\sum\limits_{i=[m/2]+1}^{m}\binom mi\langle t\rangle^{\frac{1}{4}}Y^{\frac{1}{2}}_{m-i+1}Y^{\frac{1}{2}}_{m-i+2}D_{i}\Big).\label{2.3.18}\eqa
Inserting \eqref{2.3.6}-\eqref{2.3.8}, \eqref{2.3.13}, \eqref{2.3.18} into \eqref{2.3.5} and summing over $m\geq0$ and choosing $\alpha\leq \frac{1}{2}$, we obtain
\bqa &&\frac{d}{dt}\sum\limits_{m\geq0}\|\theta_\alpha\partial_x^mu\|_{L^2(\mathbb{R}_+^2)}\tau^mM_m+\sum\limits_{m\geq0}\frac{\alpha(1-2\alpha)}{4\langle t\rangle}\tau^mM_m\frac{\|\theta_\alpha z\partial^m_xu\|^2_{L^2(\mathbb{R}_+^2)}}{\|\theta_\alpha\partial_x^mu\|_{L^2(\mathbb{R}_+^2)}}\non\\
&&\quad+(1-\frac{\alpha}{2\langle t\rangle})\sum\limits_{m\geq0}\|\theta_\alpha\partial_x^mu\|_{L^2(\mathbb{R}_+^2)}\tau^mM_m
+\sum\limits_{m\geq0}\tau^mM_m\frac{\|\theta_\alpha\partial^m_x\partial_yu\|^2_{L^2(\mathbb{R}_+^2)}}{\|\theta_\alpha\partial_x^mu\|_{L^2(\mathbb{R}_+^2)}}\non\\
&&\leq\frac{C}{\tau(t)^{\frac{1}{2}}}\Big((\sum_{m\geq 0} \|\theta_\alpha\partial_x^mu\|_{L^2(\mathbb{R}_+^2)}\tau^{m}M_m)^\frac{1}{2}(\sum_{m\geq 0} \|\theta_\alpha\partial_y\partial_x^mu\|_{L^2(\mathbb{R}_+^2)}\tau^{m}M_m)^\frac{1}{2}\non\\
&&\quad\times\sum_{m\geq 0}\|\theta_\alpha\partial_x^mu\|_{L^2(\mathbb{R}_+^2)}\tau^{m-1}mM_m\non\\
&&\quad+\langle t\rangle^{\frac{1}{4}}\sum_{m\geq 0}\|\theta_\alpha\partial_y\partial_x^mu\|_{L^2(\mathbb{R}_+^2)}\tau^{m}M_m\sum_{m\geq 0}\|\theta_\alpha\partial_x^mu\|_{L^2(\mathbb{R}_+^2)}\tau^{m-1}mM_m\Big)\non\\
&&\quad+\dot{\tau}(t)\sum\limits_{m\geq0}\|\theta_\alpha\partial_x^mu\|_{L^2(\mathbb{R}_+^2)}\tau^{m-1}mM_m\non\\
&&\leq\frac{C}{\tau(t)^{\frac{1}{2}}}\Big(\big(\sum_{m\geq 0} \|\theta_\alpha\partial_x^mu\|_{L^2(\mathbb{R}_+^2)}\tau^{m}M_m+\sum_{m\geq 0} \|\theta_\alpha\partial_y\partial_x^mu\|_{L^2(\mathbb{R}_+^2)}\tau^{m}M_m\big)\non\\
&&\quad\times\sum_{m\geq 0}\|\theta_\alpha\partial_x^mu\|_{L^2(\mathbb{R}_+^2)}\tau^{m-1}mM_m\non\\
&&\quad+\langle t\rangle^{\frac{1}{4}}\sum_{m\geq 0}\|\theta_\alpha\partial_y\partial_x^mu\|_{L^2(\mathbb{R}_+^2)}\tau^{m}M_m\sum_{m\geq 0}\|\theta_\alpha\partial_x^mu\|_{L^2(\mathbb{R}_+^2)}\tau^{m-1}mM_m\Big)\non\\
&&\quad+\dot{\tau}(t)\sum\limits_{m\geq0}\|\theta_\alpha\partial_x^mu\|_{L^2(\mathbb{R}_+^2)}\tau^{m-1}mM_m,\label{2.3.19}\eqa
where we have used the fact that for any sequences $\{a_j\}_{j\geq 0}$ and $\{b_j\}_{j\geq 0}$
\bqa \sum\limits_{m\geq 0}\sum\limits_{j\geq 0}^{m}a_jb_{m-j}\leq \sum\limits_{j\geq 0}a_j\sum\limits_{j\geq 0}b_j,\non\eqa
and the discrete H$\ddot{o}$lder's inequality for any sequences $\{a_{ij}\}_{i,j\geq 0}$ and $\{b_{ij}\}_{i,j\geq 0}$,
\bqa \sum\limits_{j=1}^{m}a_{1j}a_{2j}\cdot\cdot\cdot a_{mj}\leq \Big(\sum\limits_{j=1}^m(a_{1j})^{p_1}\Big)^{\frac{1}{p_1}}\cdot\cdot\cdot
\big(\sum\limits_{j=1}^m(a_{mj})^{p_m}\big)^{\frac{1}{p_m}},\non\eqa
where $\{p_i\}_{i\leq m}$ are positive integers and satisfy $\sum\limits_{i=1}^{m}\frac{1}{p_i}=1$.

It follows from \eqref{2.3.19} that
\bqa &&\frac{d}{dt}\|u\|_{X_{\tau,\alpha}}+(1-\frac{\alpha}{2\langle t\rangle})\|u\|_{X_{\tau,\alpha}}
+\sum\limits_{m\geq0}\tau^mM_m\frac{\|\theta_\alpha\partial^m_x\partial_yu\|^2_{L^2(\mathbb{R}_+^2)}}{\|\theta_\alpha\partial_x^mu\|_{L^2(\mathbb{R}_+^2)}}\non\\
&&\qquad\leq \Big(\dot{\tau}(t)+\frac{C}{\tau(t)^{\frac{1}{2}}}\big(\|u\|_{X_{\tau,\alpha}}+\langle t\rangle^{\frac{1}{4}}\|u\|_{D_{\tau,\alpha}}\big)\Big)\|u\|_{Y_{\tau,\alpha}}.\label{2.3.20}\eqa

Then, using Lemma 2.2.3, we discover
\bqa \sum\limits_{m\geq0}\frac{\|\theta_\alpha\partial_y\partial_x^mu\|^2_{L^2(\mathbb{R}_+^2)}}{\|\theta_\alpha\partial_x^mu\|_{L^2(\mathbb{R}_+^2)}}\tau^mM_m\geq
\frac{\alpha}{2\langle t\rangle^{\frac{1}{2}}}\|u\|_{D_{\tau,\alpha}}+\frac{\alpha(2-\sqrt\alpha)}{2\langle t\rangle}\|u\|_{X_{\tau, \alpha}},\label{2.3.21}\eqa
which, combined with \eqref{2.3.20}, gives
\bqa &&\frac{d}{dt}\|u\|_{X_{\tau,\alpha}}+\big(1+\frac{\alpha(1-\sqrt\alpha)}{2\langle t\rangle}\big)\|u\|_{X_{\tau,\alpha}}+
\frac{\alpha}{2\langle t\rangle^{\frac{1}{2}}}\|u\|_{D_{\tau,\alpha}}\non\\
&&\qquad\leq\Big(\dot{\tau}(t)+\frac{C}{\tau(t)^{\frac{1}{2}}}\big(\|u\|_{X_{\tau,\alpha}}+\langle t\rangle^{\frac{1}{4}}\|u\|_{D_{\tau,\alpha}}\big)\Big)\|u\|_{Y_{\tau,\alpha}}.\label{2.3.22}\eqa
We will choose suitable function $\tau(t)$ such that the following ordinary differential equation holds
\bqa \frac{d}{dt}({\tau}(t))^{\frac{3}{2}}+\frac{3C}{2}\big(\|u\|_{X_{\tau,\alpha}}+\langle t\rangle^{\frac{1}{4}}\|u\|_{D_{\tau,\alpha}}\big)=0.\label{2.3.23}\eqa
Inserting \eqref{2.3.23} into  \eqref{2.3.22} and applying $\alpha\in[1/4, 1/2]$, we arrive at

\bqa &&\frac{d}{dt}\|u\|_{X_{\tau,\alpha}}+\|u\|_{X_{\tau,\alpha}}+
\frac{\alpha}{2\langle t\rangle^{\frac{1}{2}}}\|u\|_{D_{\tau,\alpha}}\non\\
&&\quad\leq \frac{2}{3\tau(t)^{\frac{1}{2}}}\Big(\frac{d}{dt}({\tau}(t))^{\frac{3}{2}}+\frac{3C}{2}\big(\|u\|_{X_{\tau,\alpha}}+\langle t\rangle^{\frac{1}{4}}\|u\|_{D_{\tau,\alpha}}\big)\Big)\|u\|_{Y_{\tau,\alpha}}\non\\
&&\quad=0,\non\eqa
which gives
\bqa &&\frac{d}{dt}\big(e^t\|u\|_{X_{\tau,\alpha}}\big)+\frac{\alpha e^t}{2\langle t\rangle^{\frac{1}{2}}}\|u\|_{D_{\tau,\alpha}}\leq 0.\non\eqa

Integrating above inequality in time $t$ yields
\bqa   \|u\|_{X_{\tau,\alpha}}+\int_0^t\frac{\alpha e^{-(t-s)}}{2\langle s\rangle^{\frac{1}{2}}}\|u\|_{D_{\tau,\alpha}}ds\leq e^{-t}\|u(0)\|_{X_{\tau_0,\alpha}},\non\eqa
implying
\bqa \|u\|_{X_{\tau,\alpha}}+\frac{\alpha e^{-t}}{2}\int_0^t\frac{1}{\langle s\rangle^{\frac{1}{2}}}\|u\|_{D_{\tau,\alpha}}ds\leq e^{-t}\|u(0)\|_{X_{\tau_0,\alpha}}.\label{2.3.24}\eqa
From \eqref{2.3.23},  we easily derive for all $t\geq 0$ and $\alpha\in[1/4,1/2]$,
\bqa {\tau}(t)^{\frac{3}{2}}-{\tau}(0)^{\frac{3}{2}}&=&-\frac{3C}{2}\int_0^t\big(\|u\|_{X_{\tau,\alpha}}+\langle s\rangle^{\frac{1}{4}}\|u\|_{D_{\tau,\alpha}}\big)ds\non\\
&\geq&-\frac{3C\|u(0)\|_{X_{\tau_0,\alpha}}}{2}\Big( 1+\frac{2}{\alpha}\langle t\rangle^{\frac{3}{4}}\Big)\non\\
&\geq&-\frac{3C(\alpha+2)\|u(0)\|_{X_{\tau_0,\alpha}}}{2\alpha}\langle t\rangle^{\frac{3}{4}}\non\\
&\geq&-C_1\varepsilon\langle t\rangle^{\frac{3}{4}},\label{2.3.25}\eqa
where the positive constant $C_1=\frac{3C(\alpha+2)}{2\alpha}$.

Therefore, we choose the initial perturbation datum is so small enough  that
\bqa\tau(t)\geq \frac{\tau_0}{4},\label{2.3.26}\eqa
in the time interval $[0,T_\varepsilon]$, where $T_\varepsilon$ satisfies
\bqa T_\varepsilon=C_2\varepsilon^{-\frac{4}{3}}-1\non,\eqa
with $C_2=\Big(\frac{7}{8C_1}\tau_0^{\frac{3}{2}}\Big)^\frac{4}{3}.$
The proof is thus complete.\hfill$\Box$

We have completed the proof of the existence of solutions to the 2D Prandtl-Hartmann equations \eqref{2.3.2}-\eqref{2.3.4} with damping
term and derived the lifespan $T_\varepsilon$ of solutions.
\subsection{Uniqueness of Solutions}

In this subsection, we shall prove the uniqueness of solutions to equations
\eqref{2.3.2}-\eqref{2.3.4}. Let  $u_1(t,x,y)$ and $u_2(t,x,y)$ are two solutions to the equations \eqref{2.3.2}-\eqref{2.3.4} with same initial data $u_1(0,x,y)=u_2(0,x,y)\in X_{\tau_0,\alpha}$. Represent the tangential radii
of analytic regularity of $u_1$ and $u_2$ by $\tau_1(t)$ and $\tau_2(t)$ respectively, which satisfy the bounds in Theorem 2.2.1.
Define $\tau(t)$ to solve
\bqa \dot{\tau}(t)+\frac{C}{\tau(t)^{\frac{1}{2}}}\big(\|u_1\|_{X_{\tau_1,\alpha}}+\langle t\rangle^{\frac{1}{4}}\|u_1\|_{D_{\tau_1,\alpha}}\big)=0,\ \tau(0)=\frac{\tau_0}{8}.\label{2.4.1}\eqa
From the above equation, we can deduce
\bqa \dot{\tau}(t)+\frac{C}{\tau(t)^{\frac{1}{2}}}\big(\|u_1\|_{X_{\tau,\alpha}}+\langle t\rangle^{\frac{1}{4}}\|u_1\|_{D_{\tau,\alpha}}\big)\leq \dot{\tau}(t)+\frac{C}{\tau(t)^{\frac{1}{2}}}\big(\|u_1\|_{X_{\tau_1,\alpha}}+\langle t\rangle^{\frac{1}{4}}\|u_1\|_{D_{\tau_1,\alpha}}\big)= 0,\label{2.4.2}\eqa
where we have used the fact that $\tau(t)\leq\tau_1(t)$ and the norms $X_{\tau,\alpha}$ and $D_{\tau,\alpha}$ are increasing with respect to $\tau$.

In view of the estimate \eqref{2.3.24} for $u_1$ and the lower bounds
\eqref{2.3.26} for $\tau_1(t)$ and $\tau_2(t)$, we conclude
\bqa \frac{\tau_0}{32}\leq \tau(t)\leq\frac{\tau_0}{8}\leq \frac{\min\{\tau_1(t),\tau_2(t)\}}{2},\label{2.4.3}\eqa
for all $t\in[0,T_\varepsilon]$.

On the one hand, from \cite{iv}, the following inequality holds
\bqa \|u\|_{Y_{\tau,\alpha}}\leq \tau^{-1} \|u\|_{X_{2\tau,\alpha}}.\label{2.4.4}\eqa

On the other hand, using \eqref{2.3.24} and \eqref{2.4.3}-\eqref{2.4.4}, we also have
\bqa \|u_2\|_{Y_{\tau,\alpha}}\leq \frac{1}{\tau(t)}\|u_2\|_{X_{2\tau,\alpha}}\leq
\frac{1}{\tau(t)}\|u_2\|_{X_{\tau_2,\alpha}}\leq \frac{\varepsilon }{\tau(t)}.\label{2.4.5}\eqa

Setting $u=u_1-u_2$, then the system \eqref{2.3.2} reads as
\begin{eqnarray}\left\{\begin{array}{ll}
\partial_t u+(\bar{u}+u_1)\partial_xu+u\partial_xu_2+v\partial_y(u_1+\bar{u})+v_2\partial_yu-\partial^2_yu+u=0,\\
\partial_xu+\partial_yv=0.
\end{array}\right.\label{2.4.6}\end{eqnarray}
For $m\geq 0$, applying the operator $\partial^m_x$ on $\eqref{2.4.6}_1$, multiplying
the resulting equation by $\theta^2_\alpha\partial_x^mu$ and integrating it by parts over $\mathbb{R}_+^2$,  we derive that
\bqa\int_{\mathbb{R}_+^2}\partial^m_x\big(\partial_t u+u_1\partial_xu+u\partial_xu_2+v\partial_yu_1+v_2\partial_yu-\partial^2_yu+u\big)
\theta^2_\alpha\partial_x^mudxdy=0.\label{2.4.7}\eqa
Now, we deal with this integral equality. Similar to the estimates given in Section 2.2.1, we can deduce
\bqa\int_{\mathbb{R}_+^2}\partial_t \partial^m_xu\theta^2_\alpha\partial_x^mudxdy=
\frac{1}{2}\frac{d}{dt}\|\theta_\alpha\partial^m_xu\|^2_{L^2(\mathbb{R}_+^2)}+\frac{\alpha}{4\langle t\rangle}\|\theta_\alpha z\partial^m_xu\|^2_{L^2(\mathbb{R}_+^2)},\label{2.4.8}\eqa
\bqa \int_{\mathbb{R}_+^2}\partial^m_xu\theta^2_\alpha\partial_x^mudxdy=\|\theta_\alpha\partial^m_xu\|^2_{L^2(\mathbb{R}_+^2)},\label{2.4.9}\eqa

and
\begin{align} &-\int_{\mathbb{R}_+^2}\partial^2_y\partial^m_x u\theta^2_\alpha\partial_x^mu dxdy\nonumber\\
=&\|\theta_\alpha\partial^m_x\partial_yu\|^2_{L^2(\mathbb{R}_+^2)}-\frac{\alpha}{2{\langle t\rangle}}\|\theta_\alpha\partial^m_xu\|^2_{L^2(\mathbb{R}_+^2)}
-\frac{\alpha^2}{2\langle t\rangle}\|\theta_\alpha z\partial^m_xu\|^2_{L^2(\mathbb{R}_+^2)}.\label{2.4.10}\end{align}

Next, we  estimate the nonlinear terms in \eqref{2.4.7} in which the first term can be estimated as follows,
\bqa&&\tilde{R}_1=\int_{\mathbb{R}_+^2}\partial^m_x (u_1\partial_xu)\theta^2_\alpha\partial_x^mudxdy=\sum\limits_{i=0}^m\binom mi\int_{\mathbb{R}_+^2}\partial^{m-i}_x u_1\partial^{i+1}_xu\theta^2_\alpha\partial_x^mudxdy\non\\
&&\qquad\leq \sum\limits_{i=0}^{[m/2]}\binom mi\|\partial^{m-i}_x u_1\|_{L^2_xL_y^\infty}\|\theta_\alpha\partial^{i+1}_xu\|_{L^\infty_xL_y^2}\|\theta_\alpha\partial_x^mu\|_{L^2(\mathbb{R}_+^2)}\non\\
&&\qquad+\sum\limits_{i=[m/2]+1}^{m}\binom mi\|\partial^{m-i}_x u_1\|_{L_{xy}^\infty}\|\theta_\alpha\partial^{i+1}_xu\|_{L^2(\mathbb{R}_+^2)}\|\theta_\alpha\partial_x^mu\|_{L^2(\mathbb{R}_+^2)}.\label{2.4.11}\eqa

For $0\leq i\leq [m/2]$, using  Lemma 2.2.4 and Lemma 2.2.2 (the Agmon inequality) in $x$, we deduce
\bqa \|\partial^{m-i}_x u_1\|_{L^2_xL_y^\infty}\leq C\|\theta_\alpha \partial^{m-i}_x u_1\|^{\frac{1}{2}}_{L^2(\mathbb{R}_+^2)}\|\theta_\alpha \partial_y\partial^{m-i}_x u_1\|^{\frac{1}{2}}_{L^2(\mathbb{R}_+^2)},\label{2.4.12}\eqa
and
\bqa \|\theta_\alpha\partial^{i+1}_xu\|_{L^\infty_xL_y^2}\leq C\|\theta_\alpha\partial^{i+1}_xu\|^{\frac{1}{2}}_{L^2(\mathbb{R}_+^2)}\|\theta_\alpha\partial^{i+2}_xu\|^{\frac{1}{2}}_{L^2(\mathbb{R}_+^2)}.\label{2.4.13}\eqa

For $ [m/2]\leq i\leq m$, using Lemma 2.2.4, we arrive at
\bqa \|\partial^{m-i}_x u_1\|_{L_{xy}^\infty}&\leq& C\|\theta_\alpha \partial^{m-i}_x u_1\|^{\frac{1}{4}}_{L^2(\mathbb{R}_+^2)}\|\theta_\alpha \partial^{m-i+1}_x u_1\|^{\frac{1}{4}}_{L^2(\mathbb{R}_+^2)}\non\\
&&\quad\times\|\theta_\alpha \partial_y\partial^{m-i}_x u_1\|^{\frac{1}{4}}_{L^2(\mathbb{R}_+^2)}\|\theta_\alpha \partial^{m-i+1}_x\partial_y u_1\|^{\frac{1}{4}}_{L^2(\mathbb{R}_+^2)}.\label{2.4.14}\eqa

Therefore,
\bqa \frac{|\tilde{R}_1|\tau^mM_m}{\|\theta_\alpha\partial_x^mu\|_{L^2(\mathbb{R}_+^2)}}
&\leq& \frac{C}{\tau(t)^{\frac{1}{2}}}\Big(\sum\limits_{i=0}^{[m/2]}\binom mi\tilde{X}^{\frac{1}{2}}_{m-i}\tilde{D}^{\frac{1}{2}}_{m-i}Y^{\frac{1}{2}}_{i+1}Y^{\frac{1}{2}}_{i+2}\non\\
&+&\sum\limits_{i=[m/2]+1}^{m}\binom mi\tilde{X}^{\frac{1}{4}}_{m-i}\tilde{X}^{\frac{1}{4}}_{m-i+1}\tilde{D}^{\frac{1}{4}}_{m-i}\tilde{D}^{\frac{1}{4}}_{m-i+1}Y_{i+1}\Big).\label{2.4.15}\eqa
From now on, we use $\tilde{X}_i$, $\tilde{Y}_i$, $\tilde{D}_i$ to stand for the semi-norms of function $u_1$ and $\bar{X}_i$, $\bar{Y}_i$, $\bar{D}_i$ to stand for the semi-norms of function $u_2$.
Similarly,
\bqa \frac{|\tilde{R}_2|\tau^mM_m}{\|\theta_\alpha\partial_x^mu\|_{L^2(\mathbb{R}_+^2)}}
&\leq& \frac{C}{\tau(t)^{\frac{1}{2}}}\Big(\sum\limits_{i=0}^{[m/2]}\binom miX^{\frac{1}{2}}_{m-i}D^{\frac{1}{2}}_{m-i}\bar{Y}^{\frac{1}{2}}_{i+1}\bar{Y}^{\frac{1}{2}}_{i+2}\non\\
&&+\sum\limits_{i=[m/2]+1}^{m}\binom miX^{\frac{1}{4}}_{m-i}X^{\frac{1}{4}}_{m-i+1}D^{\frac{1}{4}}_{m-i}D^{\frac{1}{4}}_{m-i+1}\bar{Y}_{i+1}\Big),\label{2.4.16}\eqa
\bqa \frac{|\tilde{R}_3|\tau^mM_m}{\|\theta_\alpha\partial_x^mu\|_{L^2(\mathbb{R}_+^2)}}
&\leq& \frac{C}{\tau(t)^{\frac{1}{2}}}\Big(\sum\limits_{i=0}^{[m/2]}\binom mi\langle t\rangle^{\frac{1}{4}}Y_{m-i+1}\tilde{D}^{\frac{1}{2}}_{i}\tilde{D}^{\frac{1}{2}}_{i+1}\non\\
&&+\sum\limits_{i=[m/2]+1}^{m}\binom mi\langle t\rangle^{\frac{1}{4}}Y^{\frac{1}{2}}_{m-i+1}Y^{\frac{1}{2}}_{m-i+2}\tilde{D}_{i}\Big),\label{2.4.17}\eqa
and
\bqa \frac{|\tilde{R}_4|\tau^mM_m}{\|\theta_\alpha\partial_x^mu\|_{L^2(\mathbb{R}_+^2)}}
&\leq& \frac{C}{\tau(t)^{\frac{1}{2}}}\Big(\sum\limits_{i=0}^{[m/2]}\binom mi\langle t\rangle^{\frac{1}{4}}\bar{Y}_{m-i+1}D^{\frac{1}{2}}_{i}D^{\frac{1}{2}}_{i+1}\non\\
&&+\sum\limits_{i=[m/2]+1}^{m}\binom mi\langle t\rangle^{\frac{1}{4}}\bar{Y}^{\frac{1}{2}}_{m-i+1}Y^{\frac{1}{2}}_{m-i+2}D_{i}\Big).\label{2.4.18}\eqa
Similarly to \eqref{2.3.19}, inserting \eqref{2.4.8}-\eqref{2.4.18} into \eqref{2.4.7} and summing over $m\geq0$, we attain
\bqa &&\frac{d}{dt}\sum\limits_{m\geq0}\|\theta_\alpha\partial_x^mu\|_{L^2(\mathbb{R}_+^2)}\tau^mM_m+\sum\limits_{m\geq0}\frac{\alpha(1-2\alpha)}{4\langle t\rangle}\tau^mM_m\frac{\|\theta_\alpha z\partial^m_xu\|^2_{L^2(\mathbb{R}_+^2)}}{\|\theta_\alpha\partial_x^mu\|_{L^2(\mathbb{R}_+^2)}}\non\\
&&+(1-\frac{\alpha}{2\langle t\rangle})\sum\limits_{m\geq0}\|\theta_\alpha\partial_x^mu\|_{L^2(\mathbb{R}_+^2)}\tau^mM_m
+\sum\limits_{m\geq0}\tau^mM_m\frac{\|\theta_\alpha\partial^m_x\partial_yu\|^2_{L^2(\mathbb{R}_+^2)}}{\|\theta_\alpha\partial_x^mu\|_{L^2(\mathbb{R}_+^2)}}\non\\
&\leq&\frac{C}{\tau(t)^{\frac{1}{2}}}\Big(\big(\sum_{m\geq 0} \|\theta_\alpha\partial_x^mu_1\|_{L^2(\mathbb{R}_+^2)}\tau^{m}M_m+\sum_{m\geq 0} \|\theta_\alpha\partial_y\partial_x^mu_1\|_{L^2(\mathbb{R}_+^2)}\tau^{m}M_m\big)\non\\
&&\qquad\times\sum_{m\geq 0}\|\theta_\alpha\partial_x^mu\|_{L^2(\mathbb{R}_+^2)}\tau^{m-1}mM_m\non\\
&&+\langle t\rangle^{\frac{1}{4}}\sum_{m\geq 0}\|\theta_\alpha\partial_y\partial_x^mu_1\|_{L^2(\mathbb{R}_+^2)}\tau^{m}M_m\sum_{m\geq 0}\|\theta_\alpha\partial_x^mu\|_{L^2(\mathbb{R}_+^2)}\tau^{m-1}mM_m\Big)\non\\
&&+\frac{C}{\tau(t)^{\frac{1}{2}}}\Big(\big(\sum_{m\geq 0} \|\theta_\alpha\partial_x^mu\|_{L^2(\mathbb{R}_+^2)}\tau^{m}M_m+\sum_{m\geq 0} \|\theta_\alpha\partial_y\partial_x^mu\|_{L^2(\mathbb{R}_+^2)}\tau^{m}M_m\big)\non\\
&&\qquad\times\sum_{m\geq 0}\|\theta_\alpha\partial_x^mu_2\|_{L^2(\mathbb{R}_+^2)}\tau^{m-1}mM_m\non\\
&&+\langle t\rangle^{\frac{1}{4}}\sum_{m\geq 0}\|\theta_\alpha\partial_y\partial_x^mu\|_{L^2(\mathbb{R}_+^2)}\tau^{m}M_m\sum_{m\geq 0}\|\theta_\alpha\partial_x^mu_2\|_{L^2(\mathbb{R}_+^2)}\tau^{m-1}mM_m\Big)\non\\
&&+\dot{\tau}(t)\sum\limits_{m\geq0}\|\theta_\alpha\partial_x^mu\|_{L^2(\mathbb{R}_+^2)}\tau^{m-1}mM_m.\label{2.4.19}\eqa

From \eqref{2.4.19} and using Lemma 2.2.3, for $\alpha\in[1/4,1/2]$, we can deduce that
\bqa &&\frac{d}{dt}\|u\|_{X_{\tau,\alpha}}+\|u\|_{X_{\tau,\alpha}}+
\frac{\alpha}{2\langle t\rangle^{\frac{1}{2}}}\|u\|_{D_{\tau,\alpha}}\non\\
&&\qquad\leq\Big(\dot{\tau}(t)+\frac{C}{\tau(t)^{\frac{1}{2}}}\big(\|u_1\|_{X_{\tau,\alpha}}+\langle t\rangle^{\frac{1}{4}}\|u_1\|_{D_{\tau,\alpha}}\big)\Big)\|u\|_{Y_{\tau,\alpha}}+\frac{C}{\tau(t)^{\frac{1}{2}}}\big(\|u\|_{X_{\tau,\alpha}}+\langle t\rangle^{\frac{1}{4}}\|u\|_{D_{\tau,\alpha}}\big)\|u_2\|_{Y_{\tau,\alpha}}\non\\
&&\qquad\leq \frac{C\varepsilon}{\tau^{\frac{3}{2}}(t)}\big(\|u\|_{X_{\tau,\alpha}}+\langle t\rangle^{\frac{1}{4}}\|u\|_{D_{\tau,\alpha}}\big).\label{2.4.20}\eqa
According to the inequality \eqref{2.2.3}, we can derive  the following inequalities,
\bqa\frac{C\varepsilon}{\tau^{\frac{3}{2}}(t)}\leq \frac{128\sqrt2C\varepsilon}{\tau_0^{\frac{3}{2}}}\leq 1,\ \frac{C\varepsilon\langle t\rangle^{\frac{3}{4}}}{\tau^{\frac{3}{2}}(t)}\leq\frac{128\sqrt2CC_2^{\frac{3}{4}}}{\tau_0^{\frac{3}{2}}}\leq \frac{\alpha}{2}.\label{2.4.21}\eqa
Inserting  \eqref{2.4.21}  into \eqref{2.4.20}   that
\bqa &&\frac{d}{dt}\|u\|_{X_{\tau,\alpha}}+\delta\|u\|_{X_{\tau,\alpha}}\leq 0,\label{2.4.22}\eqa
where $\delta\in(0,1)$ is a small constant. This implies the uniqueness of solutions since $u(0)=0$.
The proof of uniqueness is thus complete.\hfill$\Box$

Combining \eqref{2.3.24} with \eqref{2.4.22}, we complete the proof of Theorem 2.2.1.\hfill$\Box$

\section{Bibliographic Comments}
Let us  briefly review some known results of the problem \eqref{2.1.1}-\eqref{2.1.3}. Especially, when the damping term $u_1$ does not exist, the system \eqref{2.1.1}-\eqref{2.1.3} reduces to the classical Prandtl equations which was firstly introduced formally by Prandtl (\cite{prandtl}) in 1904.
This system  is the foundation of the boundary layer equations, which describes that the away from the boundary part can be considered
as general ideal fluid, but the near a rigid wall part is deeply affected by the viscous force.
Formally, the asymptotic limit of the Navier-Stokes equations can be denoted by the Prandtl
equations within the boundary layer and by the Euler equations away from boundary. About sixty years later, Oleinik (\cite{O}) established  the first systematic
work in strictly mathematical analysis, in which she pointed out that the local-in-time well-posedness of the Prandtl system  can be proved in 2D by using the Crocco transformation  under the monotonicity condition on the tangential velocity field in the normal variable to the boundary. This result together with an expanded introduction to the boundary layer theory was showed in Oleinik-Samokhin's classical book \cite{OS}. In addition to Oleinik$'$s monotonicity assumption
on the tangential velocity field, Xin and Zhang (\cite{XZ}) obtained the existence of global weak solutions to the Prandtl equation \eqref{3.0.1a} when the pressure is favourable $(\partial_xp\leq 0)$ where the 2D Prandtl system \eqref{3.0.1a} takes the form
\begin{eqnarray}\left\{\begin{array}{ll}
\partial_tu+u\partial_xu+v\partial_yu+\partial_xp=\partial^2_{y}u, (x,y,t)\in \mathbb{R}_+^2\times\mathbb{R}_+,\\
\partial_xu+\partial_yv=0,\ (x,y,t)\in \mathbb{R}_+^2\times\mathbb{R}_+,\\
(u,v)|_{y=0}=0, \lim\limits_{y\rightarrow+\infty}u=U, \ \text{on}\ \mathbb{R}^2\times\mathbb{R}_+, \\
u|_{t=0}=u_0(x,y), \text{on}\ \mathbb{R}_+^2,
\end{array}\right.
\label{3.0.1a}\end{eqnarray}

There have been some results in the  analytic framework  with analytic radii $\tau(t)$ for 2D   boundary layer equations. First, authors of \cite{KMVW} investigated the local well-posedness of solutions to the 2D Prandtl and hydrostatic Euler equations by energy methods, which was the original result  for boundary layer problem in the analytic framework  with analytic radii $\tau(t)$. Kukavica and Vicol (\cite{KV}) considered the local well-posedness of solutions to  the  2D Prandtl boundary layer equations with  general initial data by  using analytic energy estimates in the tangential variables. Later on, Ignatova and Vicol (\cite{iv}) investigated the almost global existence for the two-dimensional Prandtl
equations when the small initial datum lies in an exponential weighted space by applying analytic energy estimates in the tangential variables and got the unique solution can be extended at least up to $T_\varepsilon\geq \exp(\varepsilon^{-1}/\log(\varepsilon^{-1}))$.

There have also been a lot of results on the 2D  Prandtl equation in different frameworks, for example, weighted Sobolev space framework  and Gevery framework, etc. Under the Oleinik$'$s monotonicity assumption, some authors in \cite{AWXY,mw,xz} proved the well-posedness of solution for 2D Prandtl equations by using energy method and constructing a new unknown function to eliminate the troublesome term from the convection term.  Qin and Wang (\cite{qinwang1}) investigated the local existence of solutions to 2D magnetic Prandtl model in the Prandtl-Hartmann regime by the energy method, see also Chapter 5 of this book. So far, besides the well-posedness of solutions to Prandtl equation  in weighted Sobolev space, there are so many results on the existence theory of solutions, for example, authors in \cite{gvm,dgv,cwz,lwx,YL,PZ}
investigated the well-posedness of the Prandtl equations  for the initial data in Gevrey framework. There have also been some  results for the Prandtl equations by using the methods of different from the energy method. For interested readers, we refer to \cite{LCS,PZ,zz,CLS,CLS1,CLS2,GD} and the references therein for the recent progress.

Xu and Zhao (\cite{XZ1}) obtained the global existence and uniqueness of weak solutions to the 2D Prandtl equations \eqref{3.0.1} (see, Chapter 3) with non-stationary boundary layer by the Crocco transformation, and assumed the outer flow has a structure $U(t,x,y)=x^mU_1(t,x,y)$ where the $x\in[0,L]$, $m\geq1$ and $U_1(t,x,y)>0$ and obtained the existence and uniqueness of weak solutions to the Prandtl equations  for any time $T$ and any interval $L$.

Motivated by \cite{xy,iv}, we investigate the global well-posedness of solutions to the 2D Prandtl-Hartmann equations \eqref{2.1.1}-\eqref{2.1.3} in an analytic framework by the standard energy method, the proof of which can be divided into two sections. First, we prove the existence of solutions for equations \eqref{2.1.1}-\eqref{2.1.3} by using the classical energy method. The difficulty of solving this problem \eqref{2.1.1}-\eqref{2.1.3} in the analytic framework is the loss of $x$-derivative in the term $v\partial_yu$. To overcome this difficulty, we show the Gaussian weighted Poincar$\acute{e}$ inequality (see Lemma 2.2.3). In other words, the Poincar$\acute{e}$ inequality does not holds in unbounded domains. Then, similar to the estimates of the existence of solutions, we verify the uniqueness of solutions for equations \eqref{2.1.1}-\eqref{2.1.3} as usual. Compared to the existence and uniqueness of solutions to the classical Prandtl equations
where the monotonicity condition of the tangential velocity, which is not needed for the 2D Prandtl-Hartmann  equations \eqref{2.1.1}-\eqref{2.1.3}  in analytic framework, plays a key role. Besides,
compared with the existence and uniqueness of solutions to the 2D MHD boundary
layer where the initial tangential magnetic field has a lower bound, which is also not needed for the 2D Prandtl-Hartmann  equations \eqref{2.1.1}-\eqref{2.1.3}  in analytic framework, this condition plays an important
role.

In addition to the above some results for the 2D Prandtl problem \eqref{3.0.1a}, there are some results on the vanishing viscosity limit for the 2D Prandtl problems \eqref{3.0.1a}, we refer to \cite{CS,GM,gmm,gn,GL,m,SC1} and the references therein for the recent progress.
So far, besides the existence of solutions to the Prandtl equations \eqref{3.0.1a} in weighted Sobolev spaces, there are so many results on the existence theory of solutions, which  have been investigated  for the  Prandtl boundary layer system \eqref{3.0.1a} when the initial data belong to the following  different suitable function spaces: \\
1) The well-posedness of solutions were investigated for the  Prandtl equations \eqref{3.0.1a} when  the  initial data with respect to the tangential variable are analytic, see \cite{LCS,KV,zz,iv,CLS,KMVW,CLS1,CLS2}.\\
2) The authors in \cite{gvm,dgv,cwz,lwx,YL,PZ}
investigated the well-posedness of the Prandtl equations \eqref{3.0.1a} for the initial data with Gevrey class in the horizontal variable $x$ and  Sobolev regularity in normal variable $y$. Therefore, it is not necessary to use the monotonicity condition of tangential velocity to overcome the loss of tangential derivative.

It is noteworthy that the above results are mainly devoted to the 2D Prandtl equations and there  also have been some results (\cite{lin,FTZ1,lwy,lwy1,lwy2,LX}) for the 3D Prandtl equations as follows,
\begin{equation}
\left\{\begin{array}{l}{\partial_{t} u+\left(u \partial_{x}+v \partial_{y}+w \partial_{z}\right) u-\partial_{z}^{2} u+\partial_xp=0},\ \text{in}\ Q,\\ {\partial_{t} v+\left(u \partial_{x}+v \partial_{y}+w \partial_{z}\right) v-\partial_{z}^{2} v+\partial_yp=0},\ \text{in}\ Q, \\ {\partial_{x} u+\partial_{y} v+\partial_{z} w=0},\ \text{in}\ Q, \\ {\left.(u, v, w)\right|_{z=0}=0, \quad \lim\limits_{z \rightarrow+\infty}(u, v)=\left(U(t,x,y), V(t,x,y)\right)},\\
(u,v,w)(t)=(u_0,v_0,w_0),\ \text{on}\  D\times \mathbb{R}_+, \end{array}\right.\label{3.159}
\end{equation}
where $Q=\{t>0,(x,y)\in D\subset\mathbb{R}^2,z>0\}$. Functions  $(U(t,x,y),V(t,x,y))$ and $p(t,x,y)$  are the
values on the boundary of the Euler's tangential velocity and Euler's pressure of the outflow,
which satisfy the Bernoulli's law,
\begin{equation}
\left\{\begin{array}{l}
\partial_tU(t,x,y)+U(t,x,y)\partial_xU(t,x,y)+V(t,x,y)\partial_yU(t,x,y)+\partial_xp(t,x,y)=0,\\
\partial_tV(t,x,y)+U(t,x,y)\partial_xV(t,x,y)+V(t,x,y)\partial_yV(t,x,y)+\partial_yp(t,x,y)=0.
\end{array}\right.\non
\end{equation}

Under a special structure on  the velocity of the outer Euler  flow,  Liu, Wang and Yang (\cite{lwy}) studied  the  local well-posedness of classical solutions of the 3D Prandtl equations \eqref{3.159} by using the Crocco transformation. The local classical solutions structure are same as the structure of the velocity of the outer Euler flow. Moreover, they proved that this structured flow is linearly stable for any smooth 3D perturbation. A continuous result was obtained with the same structure conditions (\cite{lwy}) for the 3D Prandtl equations \eqref{3.159} in \cite{lwy1}, they proved  the global existence of weak solutions  by the Crocco transformation to the 3D Prandtl equations \eqref{3.159}.
Li and Yang (\cite{ly3})  investigated the locally well-posedness in the Gevrey function
space with Gevrey index in $[1,2]$ to 3D Prandtl equations \eqref{3.159} without any monotonicity condition on the velocity field and improved  some new observation of cancellation mechanism in the three space dimensional system.

For general MHD equations, there have  been some results in the analytic framework. Xie and Yang (\cite{xy1}) studied the global existence of solutions  to the 2D MHD boundary layer equations in the mixed Prandtl and Hartmann regime when initial datum is a small perturbation of the Hartmann profile, and proved the solution in an analytic norm is an exponential decay in time. Recently, Liu and Zhang (\cite{LZ2}) established the global existence and the asymptotic decay estimate of solutions to the 2D MHD boundary layer equations with small initial data. Inspired by \cite{iv}, Xie and Yang (\cite{xy}) investigated the  lifespan of  solution to the 2D MHD boundary layer system  by using the cancellation mechanism and obtained that the lifespan of solution has a lower bound. Motivated by \cite{KV}, Lin and Zhang (\cite{lz1}) studied the local existence of  solutions for 2D incompressible magneto-micropolar boundary layer equations when the initial data are analytic in the $x$ variable by  using the energy methods. There have also been some  results for the MHD equations by using the methods of the different from energy method. For interested readers, we refer to \cite{lxy,LXY,lxy1,LWXY,hly,ghy} and the references therein for the recent progress. Liu, Xie and Yang (\cite{lxy})
studied the well-posedness of solutions to the 2D MHD in an analytic framework and inspired by \cite{GD} on the classical Prandtl equations, they  proved that if the  tangential  magnetic field is degenerate sufficiently, then the non-degenerate critical point in the tangential velocity does not prevent the formation of singularity. Liu, Xie and Yang (\cite{lxy1}) investigated  the local existence and uniqueness of solutions in weighed Sobolev space for the two-dimensional nonlinear MHD boundary layer equations  by using energy method under the assumption that the initial tangential magnetic field has a lower bound $\delta_0>0$.
As a continuation of \cite{lxy1}, the same authors in \cite{LXY} proved the validity of the Prandtl boundary layer expansion and gave a $L^\infty$ estimate on the
error by multi-scale analysis under the assumptions that both the viscosity and
resistivity coefficients with same order and the initial tangential magnetic field on the boundary is not degenerate.
Liu, Wang, Xie and Yang (\cite{LWXY}) proved the local well-posedness of the 2D MHD  boundary layer equations in Sobolev spaces, and used the initial tangential magnetic field
has a lower bound $\delta_0>0$ instead of  the monotonicity condition of the velocity field. Finally, they got the linear instability of the 2D MHD
boundary layer when the tangential magnetic field is degenerate at one point. So far, besides the well-posedness of solutions in the frameworks of Sobolev space and analytic class, there have been some results on the vanishing limits for the incompressible MHD systems, we refer to \cite{WX,wwlw,ww} and the references therein for the recent progress.

For the following 2D compressible  MHD boundary layer, Huang, Liu and Yang (\cite{hly}) attained the local well-posedness of solutions to the 2D MHD system in weighted Sobolev spaces by applying  the classical iteration scheme under the tangential magnetic field is the non-degeneracy instead of monotonicity of velocity filed and assuming the viscosity, heat conductivity and magnetic diffusivity coefficients tend to zero in the same rate,
\begin{equation}
\left\{\begin{array}{l}{\partial_{t} \rho+\nabla \cdot(\rho \mathbf{u})=0}, \\
{\rho\left(\partial_{t} \mathbf{u}+(\mathbf{u} \cdot \nabla) \mathbf{u}\right)+\nabla p(\rho, \theta)-(\nabla \times \mathbf{H}) \times \mathbf{H}}  \\
{=\mu \Delta \mathbf{u}+(\lambda+\mu) \nabla(\nabla \cdot \mathbf{u})} , \\ {c_{V} \rho\left(\partial_{t} \theta+(\mathbf{u} \cdot \nabla) \theta\right)+p(\rho, \theta) \nabla \cdot \mathbf{u}} \\ {=\kappa \Delta \theta+\lambda(\nabla \cdot \mathbf{u})^{2}+2 \mu|\mathfrak{D}(\mathbf{u})|^{2}+\nu|\nabla \times \mathbf{H}|^{2}}, \\ {\partial_{t} \mathbf{H}-\nabla \times(u \times \mathbf{H})=\nu \Delta \mathbf{H}, \quad \nabla \cdot \mathbf{H}=0}, \non\end{array}\right.
\end{equation}
where unknown functions $\rho$, $\textbf{u}$ and $\textbf{H}$ stand for the density, velocity, absolute temperature and magnetic field respectively,
and $p(\rho, \theta)$ represents  the pressure of fluid. The constant $c_{V}$ is the specific heat capacity, $\lambda$ and $\mu$  are the viscosity coefficients satisfying
$\mu>0$, $\lambda+\mu>0$, $\kappa$  is the heat conductivity coefficient and $\nu$ is the magnetic diffusivity coefficient. The deformation tensor $\mathfrak{D}(\mathbf{u})$ is given by
$$\mathfrak{D}(\mathbf{u})=\frac{\nabla \textbf{u}+(\nabla \textbf{u})'}{2}.$$

For the following 2D incompressible MHD boundary layer equations, Gao, Huang and Yao (\cite{ghy}) investigated the local well-posedness of solutions in weighted conormal Sobolev spaces to the 2D MHD boundary layer equations  with any large initial data  by the energy method under the assumption that initial tangential magnetic field has a lower bound $\delta>0$,

\begin{equation}
\left\{\begin{array}{l}{\partial_{t} \rho+\nabla \cdot(\rho \mathbf{u})=0}, \\
{\rho\left(\partial_{t} \mathbf{u}+(\mathbf{u} \cdot \nabla) \mathbf{u}\right)+\nabla p-\mu \Delta \mathbf{u}=(\textbf{H}\cdot\nabla)  \mathbf{H}} , \\
{\partial_{t} \mathbf{H}-\nabla \times(u \times \mathbf{H})=\kappa \Delta \mathbf{H}},\\
 {\nabla \cdot \mathbf{u}=0,\ \ \nabla \cdot \mathbf{H}=0}. \non\end{array}\right.
\end{equation}

\chapter{Local Existence  of Solutions to the 2D Prandtl Equations in A Weighted  Sobolev Space}
In this chapter, we shall investigate  two-dimensional nonlinear Prandtl equations on the half plane and prove the local existence  of solutions  by energy methods  in an exponential weighted Sobolev space. We use  the skill of cancellation mechanism and construct a new unknown function to overcome some difficulties respectively. The content of this chapter is chosen from \cite{qindong}.

\section{Introduction}
\setcounter{equation}{0}
In this chapter, we investigate the local existence  of the initial boundary value problem for the 2D Prandtl system which
reads as follows
\begin{eqnarray}\left\{\begin{array}{ll}
\partial_tu+u\partial_xu+v\partial_yu+\partial_xp=\partial^2_{y}u, (x,y,t)\in \mathbb{R}_+^2\times\mathbb{R}_+,\\
\partial_xu+\partial_yv=0,\ (x,y,t)\in \mathbb{R}_+^2\times\mathbb{R}_+,\\
(u,v)|_{y=0}=0, \lim\limits_{y\rightarrow+\infty}u=U, \ \text{on}\ \mathbb{R}^2\times\mathbb{R}_+, \\
u|_{t=0}=u_0(x,y), \text{on}\ \mathbb{R}_+^2,
\end{array}\right.\label{3.0.1}\end{eqnarray}
where $\mathbb{R}_+^2=\{(x,y)\in\mathbb{R}^2; y>0\}$ and $\mathbb{R}_+=[0, +\infty)$, $u(x,y,t)$ and $v(x,y,t)$ stand for the tangential velocity and normal
velocity, respectively. Functions $U(x,t)$ and $p(x,t)$    are the
values on the boundary of the Euler's tangential velocity and Euler's pressure of the outflow respectively,
which satisfy the Bernoulli's law,
$$\partial_tU(x,t)+U(x,t)\partial_xU(x,t)+\partial_xp(x,t)=0.$$

\section{Local Wellposedness of Solutions}
\setcounter{equation}{0}

As a preparation, we define some weighted Sobolev spaces and Sobolev norms as follows,
$$\|f\|^2_{L^2_\mu}=\|e^\mu f\|^2_{L^2},\ \ \|f\|^2_{H^{m,m-1}_\mu}=\sum\limits_{0\leq\alpha+\beta\leq m,\ \alpha\leq m-1}\|\partial_x^\alpha\partial_y^\beta f\|^2_{L^2_\mu},$$
$$\|f\|^2_{H^m_\mu}=\|f\|^2_{H^{m,m-1}_\mu}+\|\partial_x^m f\|^2_{L^2_\mu},\ \ D^m_\mu(\mathbb{R}^2_+)=\Big\{u:\partial_y u\in H^m_\mu(\mathbb{R}^2_+)\Big\},$$
where $\mu=\frac{y}{4}.$
If $f|_{y=0}=0$, then we get the following inequality,
\bqa \|f\|_{L^\infty(\mathbb{R}_+^2)}\leq C\|\partial_y f\|_{L^2_\mu}^{\frac{1}{2}}\|\partial_x\partial_y f\|_{L^2_\mu}^{\frac{1}{2}}.\label{3.0.2}\eqa

Suppose that $u_0^s$ satisfies the following conditions:
\begin{eqnarray}\left\{\begin{array}{ll}
\partial_y^{2j}u_0^s(0)=\tilde{u_0}|_{y=0}=0,\ 0\leq j\leq 2,\\
C^{-1}e^{-\frac{y}{4}}\leq \partial_y u^s_0(y), \ \text{for}\ y\in[0,+\infty),\\
|\partial^p_y u^s_0(y)|\leq Ce^{-\frac{y}{4}}, \ \text{for}\ y\in[0,+\infty),\ 0\leq p\leq 5.\\
\end{array}\right.\label{3.0.3}\end{eqnarray}

We can now state our main result of this chapter as follows.
\begin{Theorem}
Assume that $u^s_0(y)$ satisfies \eqref{3.0.3}, the initial datum $\tilde{u}_0(y)=u_0(y)-u_0^s(y)\in D^4_\mu(\mathbb{R}^2_+)$
satisfies the compatibility condition  up to order 6. Then there exists $0\leq T\leq T_{\max}=\sup\{t\in[0,T], 1-C\|\tilde{u}_0\|_{D^4_\mu(\mathbb{R}^2_+)}t-\frac{t^2}{2}>0\}$,
if $$\|\tilde{u}_0(y)\|_{D^4_\mu(\mathbb{R}_+^2)}\leq \eta,$$
for some $\eta\in [0,
\sqrt2(e^{-\frac{1}{2}}-\frac{1}{2})^{\frac{1}{2}}]$, then the system \eqref{3.0.1} admits a  solution $(\tilde{u},\tilde{v})$ with
$$\tilde{u}(t,x,y)=u(t,x,y)-u^s(t,y)\in L^\infty([0,T];D^4_\mu(\mathbb{R}_+^2)),\ \ \tilde{v}(t,x,y)\in L^\infty([0,T];L^\infty(\mathbb{R}_{+,y};H^3(\mathbb{R}_x))).$$
\end{Theorem}

\subsection{The Decay of Shear Flow $u^s(t,y)$ in $y$}
\setcounter{equation}{0}
For the simplicity, we consider the  case of a uniform  outflow $U(x,t)=1$, which implies pressure $\partial_xp(x,t)=0$.
Then the Prandtl equation \eqref{3.0.1} can be rewritten as
\begin{eqnarray}\left\{\begin{array}{ll}
u_t+u\partial_xu+v\partial_yu=\partial_{yy}u, t>0, x\in \mathbb{R}, y>0,\\
\partial_xu+\partial_yv=0,\\
(u,v)|_{y=0}=0, \lim\limits_{y\rightarrow+\infty}u=1,\\
u|_{t=0}=u_0(x,y).
\end{array}\right.\label{3.1.1}\end{eqnarray}

We consider the datum around a shear flow, i.e.,
$$u_0(x,y)=u_0^s(y)+\tilde{u}_0(x,y).$$
First of all, a shear flow $(u^s(t,y), 0)$ is a trivial solution to the equations \eqref{3.1.1} with $u^s(t,y)$ being the solution of heat equation
\begin{eqnarray}\left\{\begin{array}{ll}
\partial_t u^s-\partial_{y}^2u^s=0,\\
u^s|_{y=0}=0,\ \text{and}\ \lim\limits_{y\rightarrow+\infty}u^s=1,\\
u^s|_{t=0}=u^s_0(y).
\end{array}\right.\label{3.1.2}\end{eqnarray}
We rewrite the solution to \eqref{3.1.1} as a perturbation $u$ of the $u^s(t,y)$ by denoting $u(t,x,y)=u^s(t,y)+\tilde{u}(t,x,y)$, then $\tilde{u}(t,x,y)$ satisfies the following equations
\begin{eqnarray}\left\{\begin{array}{ll}
\tilde{u}_t+(u^s+\tilde{u})\partial_x\tilde{u}+\tilde{v}(\partial_yu^s+\partial_y\tilde{u})-\partial_{y}^2\tilde{u}=0,\\
\partial_x\tilde{u}+\partial_y\tilde{v}=0,\\
(\tilde{u},\tilde{v})|_{y=0}=0\ \text{and}\ \lim\limits_{y\rightarrow+\infty}\tilde{u}=0,\\
\tilde{u}|_{t=0}=\tilde{u}_0.\end{array}\right.\label{3.1.3}\end{eqnarray}
The properties of a monotonic shear flow $u^s(t,y)$ will  be first given and we also investigate the decay of shear flow in $y$ as follows.
\begin{Lemma} Assume that initial datum $u_0^s(y)$ satisfies \eqref{3.0.3},
then for any $T_1 > 0$, there exists a constant $C>0$ such that the solution $u^s(t,y)$ of the initial boundary value problem \eqref{3.1.2} satisfies
\begin{eqnarray}\left\{\begin{array}{ll}
C^{-1}e^{-\frac{y}{4}}\leq \partial_yu^s(t,x,y),\  \text{for}\ (t,x,y)\in [0,T_1]\times \mathbb{R}_+^2,\\
|\partial_y^p u^s(t,x,y)|\leq Ce^{-\frac{y}{4}},\  \text{for}\ (t,x,y)\in [0,T_1]\times \mathbb{R}_+^2,
\end{array}\right.\label{3.106}\end{eqnarray}
for some $0\leq p \leq 5$, where a constant $C>0$ depends on $T_1$.
\end{Lemma}

\textbf{Proof}. To deduce the pointwise estimates of the solutions, we can use the representation of the solution of \eqref{3.1.2},

\bqa u^s(t,y)&=&\frac{1}{2\sqrt{\pi t}}\int_0^\infty\Big( e^{-\frac{(y-\bar{y})^2}{4t}}-e^{-\frac{(y+\bar{y})^2}{4t}}\Big) u^s_0(\bar{y})d\bar{y}\non\\
&=&\frac{1}{\sqrt\pi}\Big(\int_{-\frac{y}{2\sqrt t}}^{+\infty}e^{-\xi^2}u_0^s(2\sqrt t\xi+y)d\xi-\int_{\frac{y}{2\sqrt t}}^{+\infty}e^{-\xi^2}u_0^s(2\sqrt t\xi-y)d\xi\Big).\label{}\eqa

By using $\partial^{2j}_y u^s_0(0)=0$ for $0\leq j\leq 2$, it follows
\bqa \partial^p_yu^s(t,y)&=&\frac{1}{\sqrt\pi}\Big(\int_{-\frac{y}{2\sqrt t}}^{+\infty}e^{-\xi^2}(\partial^p_yu_0^s)(2\sqrt t\xi+y)d\xi+(-1)^{p+1}\int_{\frac{y}{2\sqrt t}}^{+\infty}e^{-\xi^2}(\partial^p_yu_0^s)(2\sqrt t\xi-y)d\xi\Big)\non\\
&=&\frac{1}{2\sqrt{\pi t}}\int_0^\infty\Big( e^{-\frac{(y-\bar{y})^2}{4t}}+(-1)^{p+1}e^{-\frac{(y+\bar{y})^2}{4t}}\Big)\partial_y^p u^s_0(\bar{y})d\bar{y}:=A+B.\label{}\eqa

We now estimate $A$ and $B$ by using  the monotonicity of $\partial^p_y u^s_0(y)$, respectively. We can derive that
\bqa A&=&\frac{1}{2\sqrt{\pi t}}\int_0^\infty e^{-\frac{(y-\bar{y})^2}{4t}}\partial^p_y u^s_0(\bar{y})d\bar{y}\non\\
&\leq&C e^{-\frac{y}{4}}\frac{1}{2\sqrt{\pi t}}\int_{0}^{+\infty} e^{-\frac{(y-\bar{y})^2}{4t}}e^{\frac{(y-\bar{y})^2}{8t}}e^{\frac{t}{8}}d\bar{y}\non\\
&\leq&C e^{-\frac{y}{4}}\frac{1}{2\sqrt{\pi t}}\int_{0}^{+\infty} e^{-\frac{(y-\bar{y})^2}{8t}}d\bar{y}\non\\
&\leq&Ce^{-\frac{y}{4}}.\label{3.aaa}\eqa

Next, we establish the estimate of $B$. The result is obvious for $|y|\leq T_1$. Thus, we assume $y\geq T_1\geq t$.
Thanks to $\partial^p_y u^s_0(y)\leq Ce^{-\frac{y}{4}}$, it follows that
\bqa B&=&(-1)^{p+1}\frac{1}{2\sqrt{\pi t}}\int_0^{+\infty}e^{-\frac{(y+\bar{y})^2}{4t}} \partial^p_{\bar{y}}u^s_0(\bar{y})d\bar{y}\non\\
&\leq&C \frac{1}{2\sqrt{\pi t}} \int_0^{+\infty } e^{-\frac{y^2}{4t}}e^{-\frac{2y\bar{y}+\bar{y}^2}{4t}}e^{-\bar{y}}d\bar{y}\non\\
&\leq&Ce^{-\frac{y}{4}}\frac{1}{2\sqrt{\pi t}}\int_{0}^{+\infty} e^{-\frac{(\bar{y})^2}{4t}}e^{-\bar{y}}d\bar{y}\non\\
&\leq&Ce^{-\frac{y}{4}+t}\frac{1}{\sqrt{\pi }}\int_{\sqrt t}^{+\infty}e^{-x^2}dx\non\\
&\leq&Ce^{-\frac{y}{4}}.\label{3.1.8}\eqa

For $p=1$, we consider the estimate of another side, since the monotonic assumptions of $\partial_y u^s_0(y)$, we can conclude that
\bqa \partial_yu^s(t,y)
&=&\frac{1}{2\sqrt{\pi t}}\int_0^\infty\Big( e^{-\frac{(y-\bar{y})^2}{4t}}+e^{-\frac{(y+\bar{y})^2}{4t}}\Big)\partial_{\bar{y}}^p u^s_0(\bar{y})d\bar{\bar{y}}\non\\
&\geq&\frac{1}{2\sqrt{\pi t}}\int_0^\infty e^{-\frac{(y-\bar{y})^2}{4t}}\partial_{\bar{y}}^p u^s_0(\bar{y})d\bar{\bar{y}}\non\\
&=&\frac{1}{\sqrt\pi}\int_{\frac{y}{2\sqrt t}}^{+\infty}e^{-\xi^2}(\partial^p_\xi u_0^s)(2\sqrt t\xi+y)d\xi\non\\
&\geq& \frac{C^{-1}}{\sqrt\pi}\int_{0}^{1}e^{-\xi^2}e^{-\frac{1}{4}(2\sqrt t\xi+y)}d\xi\non\\
&\geq&C^{-1} e^{-\frac{1}{4}(2\sqrt t+y)}\non\\
 &\geq& C^{-1}(T_1) e^{-\frac{y}{4}},\ \text{for}\ y\in[0,+\infty).\label{3.200}\eqa
The proof is thus complete.\hfill $\Box$

\begin{Remark}
From this lemma, we know that the decay of the gradient of shear flow is the same as  that of the gradient of initial data of shear flow.
The gradient of shear flow is the polynomial decay when the gradient of initial data of shear flow is the polynomial decay.
\end{Remark}

In this position, we show the precise version of the compatibility condition for the nonlinear system \eqref{3.1.3}.
Since authors in \cite{xz} studied the well-posedness of Prandtl equations when the weighted function is a polynomial function, we investigate the local existence of Prandtl equations when the weighted function is an exponential function in this chapter. Thus, the compatibility conditions and  Propositions 3.2.1-3.2.3 below are the same as those in \cite{xz}.
\begin{Proposition} (\cite{xz})
Assume that $\tilde{u}(t,x,y)$ is a smooth solutions of the system \eqref{3.1.3}, then the initial datum $\tilde{u}_0(x,y)$
have to satisfy the following compatibility conditions up to order 6:
\begin{eqnarray}\left\{\begin{array}{ll}
 \tilde{u}_0(x,0)=0,\ \ (\partial_y^2\tilde{u}_0)(x,0)=0, \ \forall x\in\mathbb{R},\\
 (\partial_y^4\tilde{u}_0)(x,0)=\big(\partial_yu_0^s(0)+\partial_y\tilde{u}_0(x,0)\big)\partial_x\partial_y\tilde{u}_0(x,0), \ \forall x\in\mathbb{R},\\
 (\partial_y^6\tilde{u}_0)(x,0)=\Big(\partial^3_yu_0^s(0)+\partial^3_y\tilde{u}_0(x,0)\Big)\partial_x\partial_y\tilde{u}_0(x,0)\\
 +\sum\limits_{1\leq j\leq 3}C_4^j\Big(\partial_y^j(u^s_0+\tilde{u}_0)\partial_y^{4-j}\partial_x\tilde{u}_0-(\partial_y^{j-1}\partial_x\tilde{u}_0)\partial_y^{4-j}(\partial_yu^s_0+\partial_y\tilde{u}_0)\Big)(x,0).
 \end{array}\right.\label{3.1.9}\end{eqnarray}
\end{Proposition}

\subsection{The Approximate Solutions}

In this subsection, we investigate the existence of a parabolic regular equation  by using classical energy method
to prove the existence of solutions of the Prandtl equation.

 Specifically, for a small parameter $ \varepsilon\in(0,1]$, we study the existence of solution $\tilde{u}_\varepsilon$ to the regularized Prandtl equation,
\begin{eqnarray}\left\{\begin{array}{ll}
\partial_t \tilde{u}_\epsilon+(u^s+\tilde{u}_\epsilon)\partial_x\tilde{u}_\epsilon+\tilde{v}_\epsilon(u_y^s+\partial_y\tilde{u}_\epsilon)
-\partial_{y}^2\tilde{u}_\epsilon-\epsilon\partial_x^2\tilde{u}_\epsilon=0,\ \text{in}\ \mathbb{R}_+^2\times \mathbb{R}_+, \\
\partial_x\tilde{u}_\varepsilon+\partial_y\tilde{v}_\varepsilon=0, \ \text{in}\ \mathbb{R}_+^2\times \mathbb{R}_+,\\
(\tilde{u}_\epsilon,\tilde{v}_\varepsilon)|_{y=0}=0,\ \text{and}\ \lim\limits_{y\rightarrow+\infty}\tilde{u}_\epsilon=0, \text{on}\ \mathbb{R}\times\mathbb{R}_+,\\
\tilde{u}_\epsilon|_{t=0}=\tilde{u}_{0}+\varepsilon \mu_\varepsilon, \text{on}\ \mathbb{R}^2_+,
\end{array}\right.\label{3.2.1}\end{eqnarray}
where $\varepsilon\mu_\varepsilon$ is corrector and $\tilde{u}_0+\varepsilon\mu_\varepsilon$ satisfies the compatibility condition for the regularized system \eqref{3.2.1}, $\tilde{u}_0$ is the initial datum of the original problem \eqref{3.1.3}.

We show the boundary data of the solutions for the regularized system \eqref{3.2.1} which give also the accurate edition of the compatibility conditions
for the system \eqref{3.2.1}.

\begin{Proposition}
Assume that $\tilde{u}_0(x,y)$  satisfies the  compatibility conditions \eqref{3.1.9} for the system \eqref{3.2.1}, and $\mu_\varepsilon\in H^5_\mu(\mathbb{R}^2_+)$
 such that $\tilde{u}_0+\varepsilon\mu_\varepsilon$ satisfies the compatibility conditions  up to order 6 for the regularized system \eqref{3.2.1}.
 If $\tilde{u}_\varepsilon(t,x,y)\in L^\infty([0,T];D_\mu^4(\mathbb{R}^2_+))\cap Lip([0,T];D_\mu^2(\mathbb{R}^2_+))$
 is a solution of the system \eqref{3.2.1}, then we have
\begin{eqnarray}\left\{\begin{array}{ll}
 \tilde{u}_\varepsilon(t,x,0)=0,\ \ (\partial_y^2\tilde{u}_\varepsilon)(t,x,0)=0, \ \forall (t,x)\in[0,T]\times\mathbb{R},\\
 (\partial_y^4\tilde{u}_\varepsilon)(t,x,0)=\big(\partial_yu_\varepsilon^s(t,0)+\partial_y\tilde{u}_\varepsilon(t,x,0)\big)
 \partial_x\partial_y\tilde{u}_\varepsilon(t,x,0), \ \forall (t,x)\in[0,T]\times\mathbb{R},\\
 (\partial_y^6\tilde{u}_\varepsilon)(t,x,0)=\Big(\partial^3_yu_\varepsilon^s(t,0)+\partial^3_y\tilde{u}_\varepsilon(t,x,0)\Big)
 \partial_x\partial_y\tilde{u}_\varepsilon(t,x,0)\\
 +\sum\limits_{1\leq j\leq 3}C_4^j\Big(\partial_y^j(u^s_\varepsilon(t,0)+\tilde{u}_\varepsilon)\partial_y^{4-j}\partial_x\tilde{u}_\varepsilon-
 (\partial_y^{j-1}\partial_x\tilde{u}_\varepsilon(t,x,0))\partial_y^{4-j}(\partial_yu^s_\varepsilon(t,0)+\partial_y\tilde{u}_\varepsilon(t,x,0))\Big)\\
 -2\varepsilon\partial_x\partial_y\tilde{u}_\varepsilon(t,x,0)\partial_x^2\partial_y\tilde{u}_\varepsilon(t,x,0).
 \end{array}\right.\label{}\end{eqnarray}
\end{Proposition}
\textbf{Proof}. The proof is similar to that of Proposition 3.1 in \cite{xz}, so we omit it here.\hfill $\Box$

\begin{Proposition}
Assume that $\tilde{u}_0(x,y)$  satisfies the  compatibility conditions \eqref{3.1.9} for the system \eqref{3.1.3}, and $\partial_y\tilde{u}_0(x,y)\in H^4_\mu(\mathbb{R}^2_+)$ , then there exists  a constant $\varepsilon_0> 0$ such that for any $\varepsilon\in(0, \varepsilon_0]$, there exists a $\mu_\varepsilon\in H^5_\mu(\mathbb{R}^2_+)$
 such that $\tilde{u}_0+\varepsilon\mu_\varepsilon$ satisfies the compatibility conditions  up to order 6 for the regularized system \eqref{3.2.1}.
 Moreover,
 $$\|\partial_y\tilde{u}_{0,\varepsilon}\|^2_{H^4_\mu(\mathbb{R}^2_+)}\leq2\|\partial_y\tilde{u}_0\|^2_{H^4_\mu(\mathbb{R}^2_+)},$$
 and
 $$\lim\limits_{\varepsilon\rightarrow 0^+}\|\partial_y\tilde{u}_{0,\varepsilon}-\partial_y\tilde{u}_0\|^2_{H^4_\mu(\mathbb{R}^2_+)}=0.$$
\end{Proposition}
\textbf{Proof}. The proof is similar to that of Corollary 3.3 in \cite{xz}, thus we omit it here. \hfill $\Box$

We shall now prove the existence of approximate solutions of the system \eqref{3.2.1} by using vorticity equation $w_\varepsilon=\partial_y\tilde{u}_\varepsilon $, which reads
\begin{eqnarray}\left\{\begin{array}{ll}
\partial_t \tilde{w}_\epsilon+(u^s+\tilde{u}_\epsilon)\partial_x\tilde{w}_\epsilon+\tilde{v}_\epsilon(u_{yy}^s+\partial_y\tilde{w}_\epsilon)
-\partial_{y}^2\tilde{w}_\epsilon-\epsilon\partial_x^2\tilde{w}_\epsilon=0\ \text{in}\ \mathbb{R}_+^2\times \mathbb{R}_+,\\
\partial_y\tilde{w}_\varepsilon=0 \ \text{on}\ \mathbb{R}\times\mathbb{R}_+,\\
\tilde{w}_\epsilon|_{t=0}=\tilde{w}_{0}+\varepsilon \partial_y\mu_\varepsilon\ \text{on}\ \mathbb{R}_+^2.
\end{array}\right.\label{3.2.2}\end{eqnarray}
Using $\eqref{3.2.2}_1$ and boundary conditions $\eqref{3.2.2}_2$, we can get following the higher  derivative of boundary data,
\bqa \partial^3_y\tilde{w}_\varepsilon|_{y=0}=((u^s_y+\tilde{w}_\varepsilon)\partial_x\tilde{w}_\varepsilon)|_{y=0},\non\eqa
and
\bqa \partial^5_y\tilde{w}_\varepsilon|_{y=0}&&=(\partial^3_yu^s+\partial^2_y\tilde{w}_\varepsilon
+\varepsilon\partial^2_x\tilde{w}_\varepsilon)\partial_x\tilde{w}_\varepsilon|_{y=0}\non\\
&&-(u^s_y+\tilde{w}_\varepsilon)(\partial_x\partial^2_y\tilde{w}_\varepsilon+\varepsilon\partial^3_x\tilde{w}_\varepsilon)|_{y=0}
-\partial_x\partial_y\tilde{w}_\varepsilon(u^s_y+\tilde{w}_\varepsilon)|_{y=0}\non\\
&&+\sum_{1\leq j\leq 4}C_4^j\Big(\partial^j_y(u^s+\tilde{u}_\varepsilon)\partial^{4-j}_y\partial_x \tilde{u}_\varepsilon-
\partial^{j-1}_y\partial_x \tilde{u}_\varepsilon\partial^{4-j}_y(u^s+\tilde{w}_\varepsilon)\Big)|_{y=0}\non\\
&&-\varepsilon\partial^2_x\Big((u^s_y+\tilde{w}_\varepsilon)\partial_x\tilde{w}_\varepsilon\Big)|_{y=0},\non\eqa
where $\tilde{u}_\varepsilon=-\int_y^{+\infty}\tilde{w}_\varepsilon(t,x,\tilde{y})d\tilde{y}$ and $\tilde{v}_\varepsilon=-\int_0^{y}\partial_x\tilde{u}_\varepsilon(t,x,\tilde{y})d\tilde{y}$.

We can state  the existence of approximate solutions as follows.
\begin{Theorem} Let $\partial_y\tilde{u}_{0}\in H^4_\mu(\mathbb{R}_+^2)$  satisfies the compatibility conditions up to order 6 for system \eqref{3.1.3} and assume $u_0^s(y)$ satisfies \eqref{3.0.3},  then for any $\varepsilon\in (0, 1]$ and $\eta\in [0,
\sqrt2(e^{-\frac{1}{2}}-\frac{1}{2})^{\frac{1}{2}}]$, there exists a $T_\varepsilon>0$
such that if
$$\|\tilde{w}_{0}\|_{H^4_\mu(\mathbb{R}_+^2)}\leq \eta,$$
then the system  \eqref{3.2.2} admits a  solution $$\tilde{w}_\varepsilon\in L^\infty([0,T_\varepsilon];H^4_\mu(\mathbb{R}_+^2)),$$
which satisfies
\bqa\|\tilde{w}_\varepsilon\|^2_{L^\infty([0,T_\varepsilon];H^4_\mu(\mathbb{R}_+^2))}\leq2
\|\tilde{w}_{0,\varepsilon}\|^2_{L^\infty([0,T_\varepsilon];H^4_\mu(\mathbb{R}_+^2))}
\leq
4\|\tilde{w}_{0}\|^2_{L^\infty([0,T_\varepsilon];H^4_\mu(\mathbb{R}_+^2))}.\label{3.2.11}\eqa
\end{Theorem}

First, we will use the following  lemma,  since it is helpful to deal with boundary value.
\begin{Lemma}(\cite{adam})
Let $1< p< +\infty$. If $U\in W^{m,p}(\mathbb{R}^{n+1})$, then its trace $u$ belongs
to the space $B=B^{m-\frac{1}{p};p,p}(\mathbb{R}^{n})$ and
$$\|u\|_{B}\leq K\|U\|_{W^{m,p}(\mathbb{R}^{n+1})},$$
with the constant $K>0$ independent of $U$.
\end{Lemma}

\begin{Corollary}
Let $1< p< +\infty$. If $U\in W^{m,p}(\mathbb{R}^{n+1})$, then its trace $u$ belongs
to the space $W^{m-1,p}(\mathbb{R}^{n})$ and
$$\|u\|_{W^{m-1,p}(\mathbb{R}^{n})}\leq K\|U\|_{W^{m,p}(\mathbb{R}^{n+1})},$$
with the constant $K>0$ independent of $U$.
\end{Corollary}
\textbf{Proof}. Since $1< p< +\infty$, it follows from the fact $W^{m-1,p}(\mathbb{R}^{n})=F^{m-1}_{p,2}(\mathbb{R}^{n})$
and the embedding theorem in \cite{Hans}, $B^{m-1,p,1}(\mathbb{R}^{n})\hookrightarrow F^{m-1}_{p,2}(\mathbb{R}^{n})$
that $B^{m-\frac{1}{p};p,p}(\mathbb{R}^{n})\hookrightarrow B^{m-1,p,1}(\mathbb{R}^{n})\hookrightarrow W^{m-1,p}(\mathbb{R}^{n})$,
which gives
$$\|u\|_{W^{m-1,p}(\mathbb{R}^{n})}\leq C\|U\|_{B^{m-\frac{1}{p};p,p}(\mathbb{R}^{n})},$$
which, together with Lemma 3.2.4, gives the desired inequality. \hfill $\Box$

 Next, we will  use the following lemmas to prove  Theorem 3.2.1.
\begin{Lemma}Let $m=1,2,3,4$ and under the assumptions of Theorem 3.2.1, for a smooth solution $\tilde{w}_\varepsilon$ to equations \eqref{3.2.2}, then
\bqa&& \frac{d}{dt}\|\tilde{w}_\varepsilon\|^2_{H^{m,m-1}_\mu(\mathbb{R}^2_+)}+\varepsilon\|\partial_x\tilde{w}_{\varepsilon }\|^2_{H^{m,m-1}_\mu(\mathbb{R}^2_+)}+\|\partial_y\tilde{w}_{\varepsilon}\|^2_{H^{m,m-1}_\mu(\mathbb{R}^2_+)}\non\\
&&\qquad\leq C(\|\tilde{w}_\varepsilon\|^2_{H^4_\mu(\mathbb{R}^2_+)}+\|\tilde{w}_\varepsilon\|^4_{H^4_\mu(\mathbb{R}^2_+)}),\label{3.2.44}\eqa
where positive constant $C$ is independent of $\varepsilon$.
\end{Lemma}
\textbf{Proof}. For $m=\alpha=\beta=0$,  multiplying \eqref{3.2.2} by $e^{2\mu}\tilde{w}_\varepsilon$ and integrating it by parts over $\mathbb{R}_+^2$, we can conclude that
\bqa&& \frac{d}{dt}\|\tilde{w}_\varepsilon\|^2_{L^2_\mu}+\varepsilon\|\partial_x\tilde{w}_{\varepsilon }\|^2_{L^2_\mu}+\|\partial_y\tilde{w}_{\varepsilon}\|^2_{L^2_\mu}\non\\
&&\quad\leq\frac{1}{8}\|\partial_y\tilde{w}_\varepsilon\|^2_{L^2_\mu}+C\|\tilde{w}_\varepsilon\|^2_{L^2_\mu}+
\|\partial_x\tilde{u}_\varepsilon\|_{L^2_\mu}\|\tilde{w}_\varepsilon\|^2_{L^2_\mu}+
\|\partial_x\tilde{u}_\varepsilon\|_{L^2_\mu}\|\tilde{w}_\varepsilon\|_{L^2_\mu}\non\\
&&\qquad+
\|\partial_x\tilde{u}_\varepsilon\|_{L^2_\mu}\|\partial_y\tilde{w}_\varepsilon\|_{L^2_\mu}\|\tilde{w}_\varepsilon\|_{L^2_\mu}\non\\
&&\quad\leq\frac{1}{2}\|\partial_y\tilde{w}_\varepsilon\|^2_{L^2_\mu}+C(\|\tilde{w}_\varepsilon\|^2_{L^2_\mu}+\|\tilde{w}_\varepsilon\|^4_{L^2_\mu}
+\|\partial_x\tilde{w}_\varepsilon\|^2_{L^2_\mu}\non\\
&&\qquad+\|\partial^2_x\tilde{w}_\varepsilon\|^2_{L^2_\mu}
+\|\partial^2_x\tilde{w}_\varepsilon\|^4_{L^2_\mu}+\|\partial^2_x\tilde{w}_\varepsilon\|^4_{L^2_\mu}),\non\eqa
which implies
\bqa&& \frac{d}{dt}\|\tilde{w}_\varepsilon\|^2_{L^2_\mu}+\varepsilon\|\partial_x\tilde{w}_{\varepsilon }\|^2_{L^2_\mu}+\|\partial_y\tilde{w}_{\varepsilon}\|^2_{L^2_\mu}
\leq C(\|\tilde{w}_\varepsilon\|^2_{H^2_\mu}+\|\tilde{w}_\varepsilon\|^4_{H^2_\mu}).\label{3.2.45}\eqa
For $\alpha+\beta=m$ and $\alpha\leq m-1$ respectively, applying the operator $\partial^m$ on \eqref{3.2.2} , multiplying
the resulting equation by $e^{2\mu}\partial^m\tilde{w}$ and integrating it by parts over $\mathbb{R}_+^2$,  we conclude that
\bqa&& \frac{d}{dt}\|\partial^m\tilde{w}_\varepsilon\|^2_{L^2_\mu}+\varepsilon\|\partial_x\partial^m\tilde{w}_{\varepsilon }\|^2_{L^2_\mu}+\|\partial_y\partial^m\tilde{w}_{\varepsilon}\|^2_{L^2_\mu}=\int_{\mathbb{R}}\partial_y\partial^m\tilde{w}_{\varepsilon}
\partial^m\tilde{w}_{\varepsilon}|_{y=0}dx\non\\
&&-\int_{\mathbb{R}^2_+}e^{2\mu}\partial_y\partial^m\tilde{w}_{\varepsilon}
\partial^m\tilde{w}_{\varepsilon}dxdy-
\int_{\mathbb{R}^2_+}e^{2\mu}\partial^m\Big((u^s+\tilde{u}_\epsilon)\partial_x\tilde{w}_\epsilon+\tilde{v}_\epsilon(u_{yy}^s+\partial_y\tilde{w}_\epsilon)\Big)
\partial^m\tilde{w}_{\varepsilon}dxdy\non\\
&&=j_1+j_2+j_3+j_4.\label{3.2.46}\eqa
Although $\alpha$ and $\beta$ have many various cases, there is a same estimate on $j_2$ by the Cauchy-Schwarz inequality, that is,
$$|j_2|\leq\frac{1}{8}\|\partial_y\partial^m\tilde{w}_{\varepsilon}\|^2_{L^2_\mu}+C\|\partial^m\tilde{w}_{\varepsilon}\|^2_{L^2_\mu}.$$
Firstly, we use the boundary conditions of $\partial_y\tilde{w}_\varepsilon$ and $\tilde{v}_\varepsilon$ to deal with $j_1$. After a complicated
analysis and calculation, we can get following cases which do not vanish on the boundary of $y=0$.

Case 1: $\alpha=1,\ \beta=2$. We use the boundary value of $\partial_y^3\tilde{w}$ on $y=0$ and Lemma 3.2.4 to conclude that
\bqa \Big|\int_{\mathbb{R}}\partial_x\partial^3_y\tilde{w}_{\varepsilon}
\partial_x\partial^2_y\tilde{w}_{\varepsilon}|_{y=0}dx\Big|&\leq& \|\partial_x((u^s_y+\tilde{w}_\varepsilon)\partial_x\tilde{w}_\varepsilon)|_{y=0}\|_{L^2_\mu}
\|\partial_x\partial^2_y\tilde{w}_{\varepsilon}|_{y=0}\|_{L^2_\mu}\non\\
&\leq&\|\partial_x\tilde{w}_\varepsilon|_{y=0}\|^2_{L^2_\mu}
\|\partial_x\partial^2_y\tilde{w}_{\varepsilon}|_{y=0}\|_{L^2_\mu}\non\\
&\leq&C\|\partial_x\partial_y\tilde{w}_\varepsilon\|^2_{L^2_\mu}
\|\partial_x\partial^3_y\tilde{w}_{\varepsilon}\|_{L^2_\mu}\non\\
&\leq&C\|\tilde{w}_\varepsilon\|^4_{H^4_\mu}+C\|\tilde{w}_{\varepsilon}\|^2_{H^4_\mu}.\label{3.2.47}\eqa

Case 2: $\alpha=0,\ \beta=4$. We use the boundary value of $\partial_y^5\tilde{w}$ on $y=0$ and Lemma 3.2.4 to deduce that
\bqa \Big|\int_{\mathbb{R}}\partial^5_y\tilde{w}_{\varepsilon}
\partial^4_y\tilde{w}_{\varepsilon}|_{y=0}dx\Big|&\leq& \|\partial^5_y\tilde{w}_{\varepsilon}|_{y=0}\|_{L^2_\mu}
\|\partial^5_y\tilde{w}_{\varepsilon}\|_{L^2_\mu}\non\\
&\leq&C\|\tilde{w}_{\varepsilon}\|^2_{H^4_\mu}+C\|\tilde{w}_{\varepsilon}\|^4_{H^4_\mu}+\frac{1}{8}\|\partial_y^5\tilde{w}_{\varepsilon}\|^2_{L^2_\mu}.
\label{3.2.48}\eqa

Case 3: $\alpha=1,\ \beta=3$. We use the boundary value of $\partial_y^3\tilde{w}$ on $y=0$, Lemma 3.2.4 and integrate it by part with respect to $x$ variable to obtain that
\bqa \Big|\int_{\mathbb{R}}\partial_x\partial^4_y\tilde{w}_{\varepsilon}
\partial_x\partial^3_y\tilde{w}_{\varepsilon}|_{y=0}dx\Big|&=&\Big|\int_{\mathbb{R}}\partial^4_y\tilde{w}_{\varepsilon}
\partial^2_x\partial^3_y\tilde{w}_{\varepsilon}|_{y=0}dx\Big|\non\\
&\leq& \|\partial^2_x((u^s_y+\tilde{w}_\varepsilon)\partial_x\tilde{w}_\varepsilon)|_{y=0}\|_{L^2_\mu}
\|\partial^5_y\tilde{w}_{\varepsilon}\|_{L^2_\mu}\non\\
&\leq&C\|\tilde{w}_{\varepsilon}\|^4_{H^4_\mu}+\frac{1}{8}\|\partial_y^5\tilde{w}_{\varepsilon}\|^2_{L^2_\mu}.\label{3.2.49}\eqa

Case 4: $\alpha=2,\ \beta=2$. We use the boundary value of $\partial_y^3\tilde{w}$ on $y=0$, Lemma 3.2.4 and integrate it by part with respect to $x$ variable to derive that
\bqa \Big|\int_{\mathbb{R}}\partial^2_x\partial^3_y\tilde{w}_{\varepsilon}
\partial^2_x\partial^2_y\tilde{w}_{\varepsilon}|_{y=0}dx\Big|
&\leq& \|\partial^2_x((u^s_y+\tilde{w}_\varepsilon)\partial_x\tilde{w}_\varepsilon)|_{y=0}\|_{L^2_\mu}
\|\partial^2_x\partial^3_y\tilde{w}_{\varepsilon}\|_{L^2_\mu}\non\\
&\leq&C\|\tilde{w}_{\varepsilon}\|^4_{H^4_\mu}+\frac{1}{8}\|\partial^2_x\partial_y^3\tilde{w}_{\varepsilon}\|^2_{L^2_\mu}.\label{3.50}\eqa
We can use the estimates of commutator operator to deal with
remaining terms $j_3$ and $j_4$. For the third term $j_3$, we can write it as following pattern and use integration by part with respect to $x$ variable and use Lemma 3.2.2 to conclude that
\bqa&&-\int_{\mathbb{R}^2_+}e^{2\mu}\partial^m\Big((u^s+\tilde{u}_\epsilon)\partial_x\tilde{w}_\epsilon\Big)\partial^m\tilde{w}_\epsilon dxdy
=-\int_{\mathbb{R}^2_+}e^{2\mu}(u^s+\tilde{u}_\epsilon)\partial_x\partial^m\tilde{w}_\epsilon\partial^m\tilde{w}_\epsilon dxdy \non\\
&&\qquad+\int_{\mathbb{R}^2_+}e^{2\mu}[(u^s+\tilde{u}_\epsilon),\partial^m]\partial_x\tilde{w}_\epsilon\partial^m\tilde{w}_\epsilon dxdy\non\eqa
\bqa&&\leq  \|\partial_x\tilde{u}_{\varepsilon}\|_{L^\infty}\|\tilde{w}_{\varepsilon}\|^2_{H^4_\mu}+C\|\tilde{w}_{\varepsilon}\|^2_{H^m_\mu}+C
\|\tilde{w}_{\varepsilon}\|^3_{H^4_\mu}\non\\
&&\leq C\|\tilde{w}_{\varepsilon}\|^2_{H^4_\mu}+C\|\tilde{w}_{\varepsilon}\|^3_{H^4_\mu},\label{3.2.51}
\eqa
where we have used the fact of the estimate of commutator operator, which has been obtained by using the Sobolev inequality,
\bqa \int_{\mathbb{R}^2_+}e^{2\mu}[(u^s+\tilde{u}_\epsilon),\partial^m]\partial_x\tilde{w}_\epsilon\partial^m\tilde{w}_\epsilon dxdy
&=&\int_{\mathbb{R}^2_+}e^{2\mu}\sum_{n\leq m,\ 1\leq n}C_m^n\partial^n(u^s+\tilde{u}_\epsilon)\partial^{m-n}\partial_x\tilde{w}_\epsilon\partial^m\tilde{w}_\epsilon dxdy\non\\
&\leq& C(\|\tilde{w}_{\varepsilon}\|_{H^4_\mu}+C\|\tilde{w}_{\varepsilon}\|^2_{H^4_\mu})\|\tilde{w}_{\varepsilon}\|_{H^4_\mu}\non\\
&\leq&C\|\tilde{w}_{\varepsilon}\|^2_{H^4_\mu}+C\|\tilde{w}_{\varepsilon}\|^3_{H^4_\mu},\non \eqa
where positive constant $C$ is independent of $\varepsilon$.

For the last  term $j_4$, the estimate is similar to that of \eqref{3.2.51},
\bqa &&\Big|\int_{\mathbb{R}^2_+}e^{2\mu}\partial^m\Big(\tilde{v}_\epsilon(u_{yy}^s+\partial_y\tilde{w}_\epsilon)\Big)
\partial^m\tilde{w}_{\varepsilon}dxdy\Big|\non\\
&&=\Big|\int_{\mathbb{R}^2_+}e^{2\mu}\Big(\tilde{v}_\epsilon\partial_y\partial^m\tilde{w}_\epsilon-[\tilde{v}_\epsilon,\partial^m]\partial_y\tilde{w}_\epsilon
+\partial^m(\tilde{v}_\epsilon u_{yy}^s)\Big)\partial^m\tilde{w}_{\varepsilon}dxdy\Big|\non\\
&&\leq \|\tilde{v}_\epsilon\|_{L^\infty}\|\partial_y\partial^m\tilde{w}_\epsilon\|_{L^2_\mu}\|\partial^m\tilde{w}_\epsilon\|_{L^2_\mu}
+C\|\partial^m\tilde{w}_\epsilon\|^3_{L^2_\mu}+C\|\partial^m\tilde{w}_\epsilon\|^2_{L^2_\mu}\non\\
&&\leq \frac{1}{8}\|\partial_y\tilde{w}_\epsilon\|^2_{H^4_\mu}+C\|\tilde{w}_\epsilon\|^4_{H^4_\mu}
+C\|\tilde{w}_\epsilon\|^3_{H^4_\mu}+C\|\tilde{w}_\epsilon\|^2_{H^4_\mu}\non\\
&&\leq \frac{1}{8}\|\partial_y\tilde{w}_\epsilon\|^2_{H^4_\mu}+C\|\tilde{w}_\epsilon\|^4_{H^4_\mu}+C\|\tilde{w}_\epsilon\|^2_{H^4_\mu}.\label{3.2.52}\eqa

Thus combining \eqref{3.2.45}-\eqref{3.2.52}, we can derive the desired result \eqref{3.2.44}. The proof is thus complete.

In order to ensure the completion of energy estimate, we need to establish the estimate of term $\partial_x^m\tilde{w}_\varepsilon$,
which will be done in the following lemma.
\begin{Lemma}Let $m=1,2,3,4$ and under the assumptions of Theorem 3.2.1, for a smooth solution $\tilde{w}_\varepsilon$ to equations \eqref{3.2.2}, then
\bqa&& \frac{d}{dt}\|\partial^m_x\tilde{w}_\varepsilon\|^2_{L^2_\mu(\mathbb{R}^2_+)}+\varepsilon\|\partial^{m+1}_x\tilde{w}_{\varepsilon }\|^2_{L^2_\mu(\mathbb{R}^2_+)}+\|\partial_y\partial^m_x\tilde{w}_{\varepsilon}\|^2_{L^2_\mu(\mathbb{R}^2_+)}\non\\
&&\qquad\leq C(\|\tilde{w}_\varepsilon\|^2_{H^4_\mu(\mathbb{R}^2_+)}+\|\tilde{w}_\varepsilon\|^3_{H^4_\mu(\mathbb{R}^2_+)})+
\frac{C}{\varepsilon}(\|\tilde{w}_\varepsilon\|^2_{H^4_\mu(\mathbb{R}^2_+)}+\|\tilde{w}_\varepsilon\|^4_{H^4_\mu(\mathbb{R}^2_+)}),\label{3.2.53}\eqa
where positive constant $C$ is independent of $\varepsilon$.
\end{Lemma}
\textbf{Proof}. Applying the operator $\partial^m_x$ on \eqref{3.2.2}, multiplying
the resulting equation by $e^{2\mu}\partial_x^m\tilde{w}$ and integrating it by parts over $\mathbb{R}_+^2$,  we derive that
\bqa&& \frac{d}{dt}\|\partial_x^m\tilde{w}_\varepsilon\|^2_{L^2_\mu}+\varepsilon\|\partial^{m+1}_x\tilde{w}_{\varepsilon }\|^2_{L^2_\mu}+\|\partial_y\partial_x^m\tilde{w}_{\varepsilon}\|^2_{L^2_\mu}\non\\
&&\leq C(\|\tilde{w}_\varepsilon\|^2_{H^4_\mu}+\|\tilde{w}_\varepsilon\|^3_{H^4_\mu})+\Big|\int_{\mathbb{R}^2_+}e^{2\mu}\partial_x^m\Big(\tilde{v}_\epsilon(u_{yy}^s+\partial_y\tilde{w}_\epsilon)\Big)
\partial_x^m\tilde{w}_{\varepsilon}dxdy\Big|,\label{3.2.54}\eqa
where we have used the fact
$$(\partial_y\partial_x^m\tilde{w}_\varepsilon)(t,x,0)=0,\ (t,x)\in[0,T_\varepsilon]\times\mathbb{R}.$$
Next, we will deal with the last term on the right-hand side of \eqref{3.2.54}, which is also the main difficulty to  investigate Prandtl equations.  Using
the Leibniz formula, Sobolev inequality and Lemma 3.2.2, we can get
\bqa &&\Big|\int_{\mathbb{R}^2_+}e^{2\mu}\partial_x^m\Big(\tilde{v}_\epsilon(u_{yy}^s+\partial_y\tilde{w}_\epsilon)\Big)
\partial_x^m\tilde{w}_{\varepsilon}dxdy\Big|\leq\Big|\int_{\mathbb{R}^2_+}e^{2\mu}\tilde{v}_\epsilon\partial_y\partial_x^m\tilde{w}_\epsilon
\partial_x^m\tilde{w}_{\varepsilon}dxdy\Big|\non\\
&&\qquad+\Big|\int_{\mathbb{R}^2_+}e^{2\mu}\sum_{1\leq j\leq m-1}C_m^j\partial_x^j\tilde{v}_\epsilon\partial_y\partial_x^{m-j}\tilde{w}_\epsilon
\partial_x^m\tilde{w}_{\varepsilon}dxdy\Big|+\Big|\int_{\mathbb{R}^2_+}e^{2\mu}\partial_x^m\tilde{v}_\epsilon(u_{yy}^s+\partial_y\tilde{w}_\epsilon)
\partial_x^m\tilde{w}_{\varepsilon}dxdy\Big|\non\\
&&\leq C\|\partial_x\tilde{u}\|_{L^\infty}\|\tilde{w}_\varepsilon\|^2_{H^4_\mu}+
C\|\tilde{v}\|_{L^\infty}\|\tilde{w}_\varepsilon\|^2_{H^4_\mu}+C\|\tilde{w}_\varepsilon\|^3_{H^4_\mu}
+\|\partial_x^m\tilde{v}_\epsilon u_{yy}^s\|_{L^2_\mu}\|\tilde{w}_\varepsilon\|_{H^4_\mu}\non\\
&&\qquad+C\|\partial_x^m\tilde{v}_\epsilon\|_{L^2(\mathbb{R}_x;L^\infty(\mathbb{R}^+))}\|\partial_y\tilde{w}_\epsilon\|_{L^\infty(\mathbb{R}_x;L^2_\mu)}
\|\tilde{w}_\varepsilon\|_{H^4_\mu}\non\\
&&\leq C\|\tilde{w}_\varepsilon\|^3_{H^4_\mu}+\|\partial_x^{m+1}\tilde{w}_\varepsilon\|_{L^2_\mu}\|\tilde{w}_\varepsilon\|_{H^4_\mu}+
\|\partial_x^{m+1}\tilde{w}_\varepsilon\|_{L^2_\mu}\|\tilde{w}_\varepsilon\|^2_{H^4_\mu}\non\\
&&\leq C\|\tilde{w}_\varepsilon\|^3_{H^4_\mu}+ \frac{C}{\varepsilon}(\|\tilde{w}_\varepsilon\|^2_{H^4_\mu}+\|\tilde{w}_\varepsilon\|^4_{H^4_\mu})
+\frac{\varepsilon}{4}\|\partial_x^{m+1}\tilde{w}_\varepsilon\|^2_{L^2_\mu}.\label{3.2.55}\eqa
Inserting \eqref{3.2.55} into \eqref{3.2.54}, we can derive the desired result \eqref{3.2.53}. The proof is thus complete.\hfill $\Box$

\textbf{Proof of Theorem 3.2.1}. Combining \eqref{3.2.44} and \eqref{3.2.53},  we have for any $\varepsilon\in(0,1]$
\bqa&& \frac{d}{dt}\|\tilde{w}_\varepsilon\|^2_{H^4_\mu(\mathbb{R}^2_+)}\leq
\frac{C}{\varepsilon}(\|\tilde{w}_\varepsilon\|^2_{H^4_\mu(\mathbb{R}^2_+)}+\|\tilde{w}_\varepsilon\|^4_{H^4_\mu(\mathbb{R}^2_+)}),\label{3.2.56}\eqa
where a positive constant $C$ is independent of $\varepsilon$.

Let $y=\|\tilde{w}_\varepsilon\|^2_{H^4_\mu(\mathbb{R}^2_+)}$.  We can derive from \eqref{3.2.56},
\bqa&& \frac{dy}{dt}\leq
\frac{C}{\varepsilon}(y+y^2),\non\eqa
i.e.,
\bqa \Big(-\frac{1}{y}\Big)_t\leq \frac{C}{\varepsilon}\Big(-\big(-\frac{1}{y}\big)+1\Big).\non\eqa

Integrating above inequality over $[0,t]$, we arrive at
\bqa -\frac{1}{y}\leq -\frac{1}{y(0)}e^{-\frac{C}{\varepsilon}t}+(1-e^{-\frac{C}{\varepsilon}t}),\non\eqa
which implies
\bqa y&\leq&\frac{y_0}{e^{-\frac{C}{\varepsilon}t}-(1-e^{-\frac{C}{\varepsilon}t})y_0}\non\\
&\leq&\frac{y_0}{e^{-\frac{C}{\varepsilon}t}-\frac{C}{\varepsilon}ty_0},\non\eqa
where we used the fact $e^{-\frac{C}{\varepsilon}t}-1\geq-\frac{C}{\varepsilon}t$ as $0\leq t \leq \frac{\varepsilon}{2C}.$

Indeed, let $f(x)=e^x-1-x$ for $x<0$, thus $f'(x)=e^x-1<0$ for $[-\frac{1}{2},0]$. We have $f(0)\leq f(x)\leq f(-\frac{1}{2})$,
i.e., $0\leq e^x-1-x\leq e^{-\frac{1}{2}}-\frac{1}{2}$.
Taking $x=-\frac{C}{\varepsilon}t$, namely, $-\frac{1}{2}\leq -\frac{C}{\varepsilon}t\leq 0$, i.e., $0\leq t \leq \frac{\varepsilon}{2C}$,
 we get $f(-\frac{C}{\varepsilon}t)=e^{-\frac{C}{\varepsilon}t}-1-(-\frac{C}{\varepsilon}t)\geq 0$ for $0\leq t \leq \frac{\varepsilon}{2C}$.
Therefore, we can obtain
\bqa \|\tilde{w}_\varepsilon\|^2_{H^4_\mu(\mathbb{R}^2_+)}\leq \frac{\|\tilde{w}_\varepsilon(0)\|^2_{H^4_\mu(\mathbb{R}^2_+)}}{e^{-\frac{C}{\varepsilon}t}
-\frac{C}{\varepsilon}t\|\tilde{w}_\varepsilon(0)\|^2_{H^4_\mu(\mathbb{R}^2_+)}},\ 0<t\leq \min\{T_{\varepsilon},\frac{\varepsilon}{2C}\},\label{3.105}\eqa
where we have chosen $T_{\varepsilon}>0$ so small that
\bqa \Big(e^{-\frac{C}{\varepsilon}T_{\varepsilon}}
-\frac{C}{\varepsilon}T_{\varepsilon}\eta^2\Big)^{-1}\leq2,\ \text{when}\ 0\leq T_{\varepsilon}\leq \frac{\varepsilon}{2C}. \label{3.3.53}\eqa
In fact, since $(e^{-\frac{C}{\varepsilon}T_{\varepsilon}}
-\frac{C}{\varepsilon}T_{\varepsilon}\eta^2)^{-1}$ for $0\leq T_{\varepsilon}\leq \frac{\varepsilon}{2C}$ is a monotonically increasing function in $T_{\varepsilon}$ and $\eta\in [0,
\sqrt2(e^{-\frac{1}{2}}-\frac{1}{2})^{\frac{1}{2}}]$, \eqref{3.3.53} follows immediately.

Thus we can conclude from \eqref{3.105} and \eqref{3.3.53} that for any $\|\tilde{w}_\varepsilon(0)\|_{H^4_\mu(\mathbb{R}^2_+)}\leq \eta$ and $0<\varepsilon\leq \varepsilon_0,$
$$\|\tilde{w}_\varepsilon(t)\|^2_{H^4_\mu(\mathbb{R}^2_+)}\leq2\|\tilde{w}_\varepsilon(0)\|^2_{H^4_\mu(\mathbb{R}^2_+)}\leq 4\|\tilde{w}(0)\|^2_{H^4_\mu(\mathbb{R}^2_+)}, \ 0<t\leq T_{\varepsilon}.$$
The proof is thus complete.\hfill$\Box$

We know that the solutions of the Prandtl equation  depend on $\varepsilon$ from \eqref{3.2.53}, we need to improve the results of Lemma 3.2.6 and establish
an uniform estimate with respect to $\varepsilon$. We now follow the approach in \cite{xz} to obtain the new estimate \eqref{3.156} below.
\subsection{Formal Transformations and Uniform Estimates}

As for the estimate \eqref{3.2.53}  depends on  $\varepsilon$, we have to choose a new way to obtain an estimate which is independent of $\varepsilon$  instead of \eqref{3.2.53}. To simplify the notations, from now on, we abandon the notation tilde and sub-index $\varepsilon$, namely,
$$u=\tilde{u}_\varepsilon,\ v=\tilde{v}_\varepsilon,\ w=\tilde{w}_\varepsilon.$$

Let $w\in L^\infty([0,T];H^4_\mu(\mathbb{R}^2_+))$ be a classical solution of \eqref{3.2.2} which satisfies the following a priori estimate
\bqa \|w\|_{L^\infty([0,T];H^4_\mu(\mathbb{R}^2_+))}\leq 2\eta.\label{3.3.1}\eqa
We know that
\bqa \|e^{\frac{y}{4}}w\|_{L^\infty([0,T]\times\mathbb{R}^2_+)}&=&\Big|\int_0^y(e^{\frac{y}{4}}w)_ydy\Big|
\leq\Big(\int_0^ye^{-\frac{y}{2}}dy\Big)^{\frac{1}{2}}\Big(\int_0^ye^{\frac{y}{2}}((e^{\frac{y}{4}}w)_y)^2dy\Big)^{\frac{1}{2}}\non\\
&\leq& C\Big(\int_0^ye^{\frac{y}{2}}((e^{\frac{y}{4}}w)_y)^2dy\Big)_{L^\infty(\mathbb{R}^+_x)}^{\frac{1}{2}}\non\\
&\leq& C\Big(\int_0^\infty\int_0^ye^{\frac{y}{2}}((e^{\frac{y}{4}}w)_y)^2dxdy\Big)^{\frac{1}{2}}
\Big(\int_0^\infty\int_0^ye^{\frac{y}{2}}((e^{\frac{y}{4}}w)_{xy})^2dxdy\Big)^{\frac{1}{2}}\non\\
&\leq&C\|w\|_{H^4_{\mu}},\label{3.3.2}\eqa
which implies
$$|\partial_yu(t,x,y)|=|w(t,x,y)|\leq C\eta e^{-\frac{y}{4}},\ \forall(t,x,y)\in[0,T]\times\mathbb{R}^2_+.$$
Now choose $\eta>0$ is so small that
$$ C\eta\leq\frac{C^{-1}}{4}.$$
Then we have
\bqa \frac{C^{-1}}{4}e^{-\frac{y}{4}}\leq|u^s_y+u_y|\leq4Ce^{-\frac{y}{4}}, (t,x,y)\in[0,T]\times\mathbb{R}^2_+.\label{3.3.3}\eqa
Under the conditions \eqref{3.3.2}-\eqref{3.3.3}, system \eqref{3.2.1}  reads as follows,  for $ 0\leq m\leq 4,$
\begin{eqnarray}\left\{\begin{array}{ll}
\partial_t g_m+(u^s+u)\partial_xg_m-\partial_{y}^2g_m-\varepsilon\partial_x^2g_m-2\varepsilon(\partial_x\partial^{-1}_yg_m)\partial_y\eta_1\overset{\bigtriangleup}{=}M_m,\\
\partial_yg_m=0,\\
g_m|_{t=0}=g_{m,0},
\end{array}\right.\label{3.3.4}\end{eqnarray}
with $M_m=\sum^6_{i=1}M^m_i,$
\begin{eqnarray}\left\{\begin{array}{ll}
 M^m_1=-(u^s+u)(g_m\eta_1+(\partial^{-1}_yg_m)\partial_y\eta_1),\\
M^m_2=2(\partial_yg_m)\eta_2+2g_m(\partial_y\eta_2-2\eta_2^2)-8(\partial^{-1}_yg_m)\eta_2\partial_y\eta_2,\\
M^m_3=\varepsilon\big(2(\partial_xg_m)\eta_1-2g_m\eta_1^2-4(\partial_y^{-1}g_m)\eta_1\partial_y\eta_1\big),\\
M^m_4=\partial_y\Big(\partial_y^{-1}g_m\frac{(u^s+u)w_x+v(w_y+u^s_{yy}))}{u^s_y+u_y}\Big),\\
M^m_5=-\partial_y\Big(\frac{\sum_j^4C_4^j\partial_x^ju\cdot\partial_x^{5-j}u}{u^s_y+u_y}\Big),\\
M^m_6=-\partial_y\Big(\frac{\sum_j^4C_4^j\partial_x^jw\cdot\partial_x^{4-j}v}{u^s_y+u_y}\Big),
\end{array}\right.\non\end{eqnarray}
where $$g_m=\Big(\frac{\partial_x^mu}{u^s_y+u_y}\Big)_y, \ \eta_1=\frac{u_{xy}}{u^s_y+\tilde{u}_y}, \ \eta_2=\frac{u^s_{yy}+u_{yy}}{u^s_y+\tilde{u}_y},\ \forall(t,x,y)\in[0,T]\times\mathbb{R}^2_+.$$
\begin{Lemma}
If $\tilde{w}_0\in H^4_\mu(\mathbb{R}^2_+)$ and  \eqref{3.3.1}-\eqref{3.3.2} hold, then we have $g_m(0)\in L^2_{\frac{\mu}{2}}(\mathbb{R}^2_+)$, and
$$\|g_m(0)\|_{L^2_{\frac{\mu}{2}}(\mathbb{R}^2_+)}\leq C\|\tilde{w}_0\|_{H^4_\mu(\mathbb{R}^2_+)}.$$
\end{Lemma}
\textbf{Proof}. Actually,
$$g_m(0)=\Big(\frac{\partial_x^m\tilde{u}_0}{u_{0,y}^s+\tilde{u}_{0,y}}\Big)_y=\frac{\partial_y\partial_x^m\tilde{u}_0}{u_{0,y}^s
+\tilde{u}_{0,y}}-\frac{\partial_x^n\tilde{u}_0}{u_{0,y}^s+\tilde{u}_{0,y}}\eta_2(0).$$
From \eqref{3.3.3}, we can derive
$$e^{\frac{y}{4}}g_m(0)\leq e^{\frac{y}{2}}\partial_x^m\tilde{w}_0+(1+e^{\frac{y}{4}})\partial_x^m\tilde{u}_0,$$
implying
$$\|g_m(0)\|_{L^2_{\frac{\mu}{2}}(\mathbb{R}^2_+)}\leq C\|\partial_x^m\tilde{w}_0\|_{L^2_\mu(\mathbb{R}^2_+)}\leq C\|\tilde{w}_0\|_{H^4_\mu(\mathbb{R}^2_+)}.$$
This completes the proof of the lemma.\hfill$\Box$
\begin{Lemma}
Let $w\in L^\infty([0,T];H^4_\mu(\mathbb{R}^2_+))$ and  \eqref{3.3.1}-\eqref{3.3.2} hold. Assume $u_0^s(y)$ satisfies \eqref{3.0.3} and $g_m$ satisfies the equation \eqref{3.3.4} for $1\leq m\leq 4$, then we have the following
estimates, for $t\in[0,T]$,
\bqa \frac{d}{dt}\sum_{m=1}^4\|g_m\|^2_{L^2_{\frac{\mu}{2}}(\mathbb{R}^2_+)}+\sum_{m=1}^4\|\partial_yg_m\|^2_{L^2_{\frac{\mu}{2}}(\mathbb{R}^2_+)}
+\varepsilon\sum_{m=1}^4\|\partial_xg_m\|^2_{L^2_{\frac{\mu}{2}}(\mathbb{R}^2_+)}\non\\
\qquad \qquad \leq C\Big(\sum_{m=1}^4\|g_m\|^2_{L^2_{\frac{\mu}{2}}(\mathbb{R}^2_+)}+\|w\|^2_{H^4_\mu}\Big)\label{3.3.5},\eqa
where a constant $C>0$ is independent of $\varepsilon$.
\end{Lemma}
\textbf{Proof}. Multiplying \eqref{3.3.4} by $e^{\frac{y}{4}}g_m$, integrating it by part over $\mathbb{R}^2_+$ and applying the boundary values and Lemma 3.2.2,
 we have
\bqa &&\frac{1}{2}\frac{d}{dt}\|g_m\|^2_{L^2_{\frac{\mu}{2}}(\mathbb{R}^2_+)}+\|\partial_yg_m\|^2_{L^2_{\frac{\mu}{2}}(\mathbb{R}^2_+)}+
\varepsilon \|\partial_xg_m\|^2_{L^2_{\frac{\mu}{2}}(\mathbb{R}^2_+)}\non\\
&&=-\int_{\mathbb{R}^2_+}(u^s+u)\partial_xg_me^{\frac{y}{4}}g_mdxdy-\frac{1}{4}\int_{\mathbb{R}^2_+}\partial_yg_me^{\frac{y}{4}}g_mdxdy\non\\
&&\qquad-2\varepsilon
\int_{\mathbb{R}^2_+}(\partial_x\partial^{-1}_yg_n)\partial_y\eta_1e^{\frac{y}{4}}g_mdxdy+\sum_{j=1}^6\Big|\int_{\mathbb{R}^2_+}M^m_j
e^{\frac{y}{4}}g_mdxdy\Big|\non\\
&&\leq C\|u_x\|_{L^\infty(\mathbb{R}^2_+)}\|g_m\|^2_{L^2_{\frac{\mu}{2}}(\mathbb{R}^2_+)}+\frac{1}{4}\|\partial_yg_m\|^2_{L^2_{\frac{\mu}{2}}(\mathbb{R}^2_+)}
+C\|g_m\|^2_{L^2_{\frac{\mu}{2}}(\mathbb{R}^2_+)}\non\\
&&\qquad+2\varepsilon\int_{\mathbb{R}^2_+}(\partial^{-1}_yg_n)\partial_y\eta_1e^{\frac{y}{4}}\partial_xg_mdxdy
+2\varepsilon\int_{\mathbb{R}^2_+}(\partial^{-1}_yg_n)\partial_y\partial_x\eta_1e^{\frac{y}{4}}g_mdxdy\non\\
&&\qquad+\sum_{j=1}^6\Big|\int_{\mathbb{R}^2_+}M^m_j
e^{\frac{y}{4}}g_mdxdy\Big|\non\\
&&\leq C\|w\|_{H^2_\mu(\mathbb{R}^2_+)}\|g_m\|^2_{L^2_{\frac{\mu}{2}}(\mathbb{R}^2_+)}+\frac{1}{4}\|\partial_yg_m\|^2_{L^2_{\frac{\mu}{2}}(\mathbb{R}^2_+)}
+C\|g_m\|^2_{L^2_{\frac{\mu}{2}}(\mathbb{R}^2_+)}\non\\
&&\qquad+\varepsilon\|\partial_y^{-1}g_m\partial_y\eta_1\|
^2_{L^2_{\frac{\mu}{2}(\mathbb{R}^2_+)}}+\frac{\varepsilon}{8}\|\partial_xg_m\|^2_{L^2_{\frac{\mu}{2}}(\mathbb{R}^2_+)}+
\varepsilon\|\partial_y^{-1}g_m\partial_y\partial_x\eta_1\|
^2_{L^2_{\frac{\mu}{2}(\mathbb{R}^2_+)}}\non\\
&&\qquad+\varepsilon\|g_m\|^2_{L^2_{\frac{\mu}{2}}(\mathbb{R}^2_+)}+\sum_{j=1}^6\Big|\int_{\mathbb{R}^2_+}M^m_j
e^{\frac{y}{4}}g_mdxdy\Big|,\non
\eqa
 for $0<\varepsilon\leq 1$, which implies
\bqa &&\frac{d}{dt}\|g_m\|^2_{L^2_{\frac{\mu}{2}}(\mathbb{R}^2_+)}+\|\partial_yg_m\|^2_{L^2_{\frac{\mu}{2}}(\mathbb{R}^2_+)}+
\varepsilon \|\partial_xg_m\|^2_{L^2_{\frac{\mu}{2}}(\mathbb{R}^2_+)}\non\\
&&\leq C\|g_m\|^2_{L^2_{\frac{\mu}{2}}(\mathbb{R}^2_+)}+\|\partial_y^{-1}g_m\partial_y\eta_1\|
^2_{L^2_{\frac{\mu}{2}}(\mathbb{R}^2_+)}+\|\partial_y^{-1}g_m\partial_y\partial_x\eta_1\|
^2_{L^2_{\frac{\mu}{2}}(\mathbb{R}^2_+)}\non\\
&&\qquad+2\sum_{j=1}^6\Big|\int_{\mathbb{R}^2_+}M^m_j
e^{\frac{y}{4}}g_mdxdy\Big|.\non
\label{3.3.6*}\eqa
We use the similar method as in \cite{xz} to deal with right-hand side of terms in above inequality.
Firstly, we derive the estimate of second and third terms from \eqref{3.3.1} and \eqref{3.3.3},
\bqa\big|\eta_1\big|&=&\Big|\frac{u_{xy}}{u^s_y+\tilde{u}_y}\Big|\leq C e^{-\frac{y}{4}},\
\big|\partial_x\eta_1\big|=\Big|\frac{u_{xxy}}{u^s_y+\tilde{u}_y}-\frac{u_{xy}\tilde{u}_{xy}}{(u^s_y+\tilde{u}_y)^2}\Big|\leq C e^{-\frac{y}{4}},\non\\
\big|\partial_y\eta_1\big|&=&\Big|\frac{u_{xyy}}{u^s_y+\tilde{u}_y}-\frac{u_{xy}\tilde{u}_{yy}}{(u^s_y+\tilde{u}_y)^2}\Big|\leq C e^{-\frac{y}{4}},\
\big|\partial_x\partial_y\eta_1\big|\leq C e^{-\frac{y}{4}}.\non
\eqa
Then we have $$\|\partial_y^{-1}g_m\partial_y\eta_1\|
^2_{L^2_{\frac{\mu}{2}}(\mathbb{R}^2_+)}\leq C \int_{\mathbb{R}^2_+}e^{-\frac{y}{2}}\Big(\int_0^yg_m(t,x,\tilde{y})d\tilde{y}\Big)^2dxdy\leq C\|g_m\|^2_{L^2_{\frac{\mu}{2}}},$$
and
$$\|\partial_y^{-1}g_m\partial_y\partial_x\eta_1\|
^2_{L^2_{\frac{\mu}{2}}(\mathbb{R}^2_+)}\leq C \int_{\mathbb{R}^2_+}e^{-\frac{y}{2}}\Big(\int_0^yg_m(t,x,\tilde{y})d\tilde{y}\Big)^2dxdy\leq C\|g_m\|^2_{L^2_{\frac{\mu}{2}}}.$$

Then we estimate the terms of $M^m_j, j=1,2,\cdot\cdot\cdot,6$. For the term $M^m_1$, we have
\bqa\Big|\int_{\mathbb{R}^2_+}M^m_1e^{\frac{y}{4}}g_mdxdy\Big| &=&\Big|\int_{\mathbb{R}^2_+}(u^s+u)(g_m\eta_1+(\partial^{-1}_yg_m)\partial_y\eta_1)e^{\frac{y}{4}}g_mdxdy\Big|\non\\
&\leq&C\|g_m\|^2_{L^2_{\frac{\mu}{2}}}+\|\partial_y^{-1}g_m\partial_y\eta_1\|
^2_{L^2_{\frac{\mu}{2}}(\mathbb{R}^2_+)}\non\\
&\leq&C\|g_m\|^2_{L^2_{\frac{\mu}{2}}}.\label{3.149}\eqa
We
use the similar method to that of $\eta_1$ above to get decay rate of $\eta_2$,
\bqa\big|\eta_2\big|\leq C(1+ e^{-\frac{y}{4}}),\
\big|\partial_x\eta_2\big|\leq C e^{-\frac{y}{4}},\non\\
\big|\partial_y\eta_2\big|\leq C e^{-\frac{y}{4}},\
\big|\partial_x\partial_y\eta_2\big|\leq C e^{-\frac{y}{4}}.\non
\eqa
Thus we get
\bqa\Big|\int_{\mathbb{R}^2_+}M^m_2e^{\frac{y}{4}}g_mdxdy\Big|&\leq& C\Big|\int_{\mathbb{R}^2_+}\Big((\partial_yg_m)\eta_2+2g_m(\partial_y\eta_2-2\eta_2^2)-8(\partial^{-1}_yg_m)\eta_2\partial_y\eta_2\Big)
e^{\frac{y}{4}}g_mdxdy\Big|\non\\
&\leq& C\|g_m\|^2_{L^2_{\frac{\mu}{2}}}+\frac{1}{4}\|\partial_yg_m\|^2_{L^2_{\frac{\mu}{2}}}+\|\partial_y^{-1}g_m\partial_y\eta_2\|
^2_{L^2_{\frac{\mu}{2}}(\mathbb{R}^2_+)}\non\\
&\leq& C\|g_m\|^2_{L^2_{\frac{\mu}{2}}}+\frac{1}{4}\|\partial_yg_m\|^2_{L^2_{\frac{\mu}{2}}}.\label{3.150}\eqa
Similarly as above for $M_2^m$,  we get for any $\varepsilon\in(0,1]$,
\bqa\Big|\int_{\mathbb{R}^2_+}M^m_3e^{\frac{y}{4}}g_mdxdy\Big|\leq C\|g_m\|^2_{L^2_{\frac{\mu}{2}}}+\frac{\varepsilon}{4}\|\partial_xg_m\|^2_{L^2_{\frac{\mu}{2}}}.\label{3.151}\eqa
For $M_4^m$,  using a direct computation and  \eqref{3.3.1}-\eqref{3.3.3}, we have
\bqa\Big|\int_{\mathbb{R}^2_+}M^m_4e^{\frac{y}{4}}g_mdxdy\Big|\leq C\|g_m\|^2_{L^2_{\frac{\mu}{2}}}.\label{3.152}\eqa

Noting that  $M_5^m$  can be written as
\bqa M_5^m=g_4\partial_xu+\sum_{1\leq m\leq 3}C_4^m\partial_xug_{5-m}
+\partial_y^{-1}g_4\partial_xw+\sum_{1\leq m\leq 3}C_4^m\partial^m_xw\partial_y^{-1}g_{5-m},\non\eqa
we have
\bqa\Big|\int_{\mathbb{R}^2_+}M^m_5e^{\frac{y}{4}}g_mdxdy\Big|&\leq&C\|g_m\|_{L^2_{\frac{\mu}{2}}}
(\|g_4\partial_xu\|_{L^2_{\frac{\mu}{2}}}+\|\partial_y^{-1}g_4\partial_xw\|_{L^2_{\frac{\mu}{2}}})\non\\
&&+C\|g_m\|_{L^2_{\frac{\mu}{2}}}\sum_{1\leq m\leq 3}(\|\partial_xug_{5-m}\|_{L^2_{\frac{\mu}{2}}}+\|\partial^m_xw\partial_y^{-1}g_{5-m}\|_{L^2_{\frac{\mu}{2}}})\non\\
&\leq&C\|g_m\|_{L^2_{\frac{\mu}{2}}}
(\|g_4\|_{L^2_{\frac{\mu}{2}}}\|\partial_xu\|_{L^\infty}+\|g_4\|_{L^2_{\frac{\mu}{2}}})\non\\
&&+C\|g_m\|_{L^2_{\frac{\mu}{2}}}\sum_{1\leq m\leq 3}(\|g_{5-m}\|_{L^2_{\frac{\mu}{2}}}\|\partial_xu\|_{L^\infty}+\|g_{5-m}\|_{L^2_{\frac{\mu}{2}}})\non\\
&\leq&C\|g_m\|_{L^2_{\frac{\mu}{2}}}
(\|g_4\|_{L^2_{\frac{\mu}{2}}}\|\partial_xw\|_{L^2_\mu(\mathbb{R}_+^2)}+\|g_4\|_{L^2_{\frac{\mu}{2}}})\non\\
&&+C\|g_m\|_{L^2_{\frac{\mu}{2}}}\sum_{1\leq m\leq 3}(\|g_{5-m}\|_{L^2_{\frac{\mu}{2}}}\|\partial_xw\|_{L^2_\mu(\mathbb{R}_+^2)}+\|g_{5-m}\|_{L^2_{\frac{\mu}{2}}})\non\\
&\leq&C\sum_{m=1}^4\|g_m\|^2_{L^2_{\frac{\mu}{2}}},\label{3.153}\eqa
where we have used the fact
\bqa \|e^{\frac{y}{4}}\partial^m_xw\|_{L^\infty([0,T]\times\mathbb{R}^2_+)}&=&\int_0^y(e^{\frac{y}{4}}\partial^m_xw)_ydy
\leq\Big(\int_0^ye^{-\frac{y}{2}}dy\Big)^{\frac{1}{2}}\Big(\int_0^ye^{\frac{y}{2}}((e^{\frac{y}{4}}\partial^m_xw)_y)^2dy\Big)^{\frac{1}{2}}\non\\
&\leq& C\Big(\int_0^ye^{\frac{y}{2}}((e^{\frac{y}{4}}w)_y)^2dy\Big)_{L^\infty(\mathbb{R}^+_x)}^{\frac{1}{2}}\non\\
&\leq& C\Big(\int_0^\infty\int_0^ye^{\frac{y}{2}}((e^{\frac{y}{4}}\partial^m_xw)_y)^2dxdy\Big)^{\frac{1}{2}}
\Big(\int_0^\infty\int_0^ye^{\frac{y}{2}}((e^{\frac{y}{4}}\partial^m_xw)_{xy})^2dxdy\Big)^{\frac{1}{2}}\non\\
&\leq&C\|w\|_{H^4_{\mu}}.\non\eqa

Therefore it follows $$|\partial^m_xw(x,y,t)|\leq Ce^{-\frac{y}{4}}, \forall(x,y,t)\in\mathbb{R}^2_+\times[0,T].$$
Last, we deal with the term $M^m_6$ by using integration by part and boundary condition,
\bqa&-&\int_{\mathbb{R}^2_+}M^m_6e^{\frac{y}{4}}g_mdxdy\non\\
&=&\frac{1}{4}\int_{\mathbb{R}^2_+}\Big(\frac{\sum_j^4C_4^j\partial_x^jw\cdot\partial_x^{4-j}v}{u^s_y+u_y}\Big)e^{\frac{y}{4}}g_mdxdy
+\int_{\mathbb{R}^2_+}\Big(\frac{\sum_j^4C_4^j\partial_x^jw\cdot\partial_x^{4-j}v}{u^s_y+u_y}\Big)e^{\frac{y}{4}}\partial_yg_mdxdy\non\\
&\leq&C\int_{\mathbb{R}^2_+}\sum_{j=1}^4\partial_x^jw\partial_x^{4-i}vg_mdxdy+C\int_{\mathbb{R}^2_+}\sum_{j=1}^4\partial_x^jw\partial_x^{4-i}v\partial_yg_mdxdy\non\\
&\leq&C(\|\partial_x^4w\|_{L^2_\mu}\|v\|_{L^\infty}+\sum_{j=1}^3\|\partial_x^jw\|_{L^\infty(R_x;L^2(R_y))}
\|\partial_x^{4-j}v\|_{L^\infty(R_y;L^2(R_x))})(\|g_m\|_{L^2_\mu}+\|\partial_yg_m\|_{L^2_\mu})\non\\
&\leq&C\|w\|_{H^4_\mu}(\|g_m\|_{L^2_\mu}+\|\partial_yg_m\|_{L^2_\mu})\non\\
&\leq&C\|w\|^2_{H^4_\mu}+C\|g_m\|^2_{L^2_\mu}+\frac{1}{4}\|\partial_yg_m\|^2_{L^2_\mu}.\label{3.154}\eqa

Thus, inserting \eqref{3.149}-\eqref{3.154} into \eqref{3.3.6*} completes the proof of \eqref{3.3.5}  .\hfill$\Box$

In this position, we give the energy estimate for the sequence of approximate solutions. Let's go back to  the notations with
tilde and the sub-index $\varepsilon$ and $g_m^\varepsilon$ is the function defined by $\tilde{u}_\varepsilon$.

We now collect some established  estimates \eqref{3.2.44} and \eqref{3.3.5} to be used,
\bqa&& \frac{d}{dt}\sum_{m=1}^4\|\tilde{w}_\varepsilon\|^2_{H^{m,m-1}_\mu(\mathbb{R}^2_+)}+\varepsilon\sum_{m=1}^4\|\partial_x\tilde{w}_{\varepsilon }\|^2_{H^{m,m-1}_\mu(\mathbb{R}^2_+)}+\sum_{m=1}^4\|\partial_y\tilde{w}_{\varepsilon}\|^2_{H^{m,m-1}_\mu(\mathbb{R}^2_+)}\non\\
&&\qquad\leq C(\|\tilde{w}_\varepsilon\|^2_{H^4_\mu(\mathbb{R}^2_+)}+\|\tilde{w}_\varepsilon\|^4_{H^4_\mu(\mathbb{R}^2_+)}),\label{3.3.6}\eqa
\bqa \frac{d}{dt}\sum_{m=1}^4\|g^\varepsilon_m\|^2_{L^2_{\frac{\mu}{2}}(\mathbb{R}^2_+)}+\sum_{m=1}^4\|\partial_yg^\varepsilon_m\|^2_{L^2_{\frac{\mu}{2}}(\mathbb{R}^2_+)}
+\varepsilon\sum_{m=1}^4\|\partial_xg^\varepsilon_m\|^2_{L^2_{\frac{\mu}{2}}(\mathbb{R}^2_+)}\non\\
\qquad \qquad \leq C\Big(\sum_{m=1}^4\|g^\varepsilon_m\|^2_{L^2_{\frac{\mu}{2}}(\mathbb{R}^2_+)}+\|\tilde{w}_\varepsilon\|^2_{H^4_\mu}\Big)\label{3.3.7}.\eqa
\begin{Lemma}Let $m=1,2,3,4$. Then we have
\bqa h_m^\varepsilon(g,w)(0)=\sum_{m=1}^4\|g^\varepsilon_m\|^2_{L^2_{\frac{\mu}{2}}(\mathbb{R}^2_+)}+\|\tilde{w}_\varepsilon(0)\|^2_{H^{m,m-1}_\mu(\mathbb{R}^2_+)}
\leq C\|\tilde{u}_0\|^2_{D^4_\mu(\mathbb{R}^2_+)},\label{3.155}\eqa
where a constant $C>0$ is independent of $\varepsilon$.
\end{Lemma}
\textbf{Proof}. For any $1\leq m\leq 4$, applying  Lemma 3.2.7 and $\tilde{u}_\varepsilon(0)=\tilde{u}(0)$, we conclude that
$$h_m^\varepsilon(g,w)(0)\leq C\|\tilde{w}_0\|^2_{H^4_\mu(\mathbb{R}^2_+)}\leq C\|\tilde{u}_0\|^2_{D^4_\mu(\mathbb{R}^2_+)},$$
which is just \eqref{3.155}.\hfill$\Box$
\begin{Lemma}
Let $m=1,2,3,4$, the following estimate holds
\bqa\|\partial_x^m\tilde{w}_\varepsilon\|^2_{L^2_\mu(\mathbb{R}^2_+)}\leq C\|g^\varepsilon_m\|^2_{L^2_{\frac{\mu}{2}}(\mathbb{R}^2_+)},\label{3.156}\eqa
where a constant $C>0$ is independent of $\varepsilon$.
\end{Lemma}
\textbf{Proof}.
By the definition, we know that
$$\partial_x^m\tilde{u}_\varepsilon(t,x,y)=(u^s_y+\tilde{w}_\varepsilon)\int_0^yg^\varepsilon_m(t,x,z)dz,\ y\geq 0.$$
Thus,
$$\partial_x^m\tilde{w}_\varepsilon(t,x,y)=(u^s_{yy}+(\tilde{w}_\varepsilon)_y)\int_0^yg^\varepsilon_m(t,x,z)dz
+(u^s_y+\tilde{w}_\varepsilon)g^\varepsilon_m(t,x,y),\ y\geq 0,$$
which gives
\bqa\|\partial_x^m\tilde{w}_\varepsilon\|^2_{L^2_\mu(\mathbb{R}^2_+)}\leq C\|g^\varepsilon_m\|^2_{L^2(\mathbb{R}^2_+)}+ C\|g^\varepsilon_m\|^2_{L^2_{\frac{\mu}{2}}(\mathbb{R}^2_+)}\leq C\|g^\varepsilon_m\|^2_{L^2_{\frac{\mu}{2}}(\mathbb{R}^2_+)}.\nonumber\eqa\hfill$\Box$

\begin{Theorem} Assume  $\tilde{u}_0\in D^4_\mu(\mathbb{R}^2_+)$
 satisfies the compatibility condition \eqref{3.1.9}.
Suppose that $\tilde{w}_\varepsilon\in L^\infty([0,T];H^4_\mu(\mathbb{R}^2_+))$ is a solution to \eqref{3.2.2} such that
\bqa\|\tilde{w}_\varepsilon\|_{L^\infty([0,T];H^4_\mu(\mathbb{R}^2_+))}\leq 2\eta,\label{157}\eqa
where $0\leq T\leq T_{\max}=\sup\{t\in[0,T], 1-C\|\tilde{u}_0\|_{D^4_\mu(\mathbb{R}^2_+)}t-\frac{t^2}{2}>0\}$.  Then we have

\bqa \|\tilde{w}_\varepsilon\|^2_{H^{4}_\mu(\mathbb{R}^2_+)}\leq
(\|\tilde{u}_0\|^2_{D^4_\mu(\mathbb{R}^2_+)}+1)\Big\{1-C\|\tilde{u}_0\|^2_{D^4_\mu(\mathbb{R}^2_+)}t-\frac{t^2}{2}\Big\}^{-1},\label{3.3.8}\eqa
where a constant $C>0$ is independent of $\varepsilon\in (0,1].$

\end{Theorem}

\textbf{Proof}. Applying \eqref{3.3.6}-\eqref{3.3.7}, for $m=1,2,3,4$, we get
\bqa \sum_{i=1}^4\big(\|g^\varepsilon_m\|^2_{L^2_{\frac{\mu}{2}}(\mathbb{R}^2_+)}+\|\tilde{w}_\varepsilon\|^2_{H^{m,m-1}_\mu(\mathbb{R}^2_+)}\big)\leq
(\|\tilde{u}_0\|^2_{D^4_\mu(\mathbb{R}^2_+)}+t)\non\\
\qquad+C\int_0^t\Big(\sum_{i=1}^4(\|g^\varepsilon_m\|^2_{L^2_{\frac{\mu}{2}}(\mathbb{R}^2_+)}+\|\tilde{w}_\varepsilon\|^2_{H^{m,m-1}_\mu(\mathbb{R}^2_+)})\Big)^2ds.\non\eqa
Applying nonlinear Gronwall inequality (Theorem 2, P362, \cite{mpf}) to above inequality, we have
\bqa \sum_{i=1}^4\big(\|g^\varepsilon_m\|^2_{L^2_{\frac{\mu}{2}}(\mathbb{R}^2_+)}+\|\tilde{w}_\varepsilon\|^2_{H^{m,m-1}_\mu(\mathbb{R}^2_+)}\big)\leq
(\|\tilde{u}_0\|^2_{D^4_\mu(\mathbb{R}^2_+)}+t)\Big\{1-C\|\tilde{u}_0\|^2_{D^4_\mu(\mathbb{R}^2_+)}t-\frac{t^2}{2}\Big\}^{-1}.\label{3.3.9}\eqa

Combining \eqref{3.3.9} and Lemma 3.2.10, we have, for $t\in[0,T]$,
\bqa \|\tilde{w}_\varepsilon\|^2_{H^{4}_\mu(\mathbb{R}^2_+)}\leq
(\|\tilde{u}_0\|^2_{D^4_\mu(\mathbb{R}^2_+)}+t)\Big\{1-C\|\tilde{u}_0\|^2_{D^4_\mu(\mathbb{R}^2_+)}t-\frac{t^2}{2}\Big\}^{-1},\label{3.158}\eqa
where a constant $C>0$  is independent of $\varepsilon\in(0,1]$.
Thus \eqref{3.3.8} follows immediately.\hfill$\Box$

\subsection{ Convergence and Consistency }

Using evolution equation \eqref{3.2.2} and  uniform  $H^4_\mu$ bound in \eqref{3.2.11}, we conclude that $\partial_t\tilde{w}^\varepsilon$
is uniformly (in $\varepsilon$) bounded in $L^\infty([0,T];H^2_\mu)$. By the Lions-Aubin Lemma and the compact embedding of
$H^4_\mu$ in $H^{4-\delta}_{\mu, loc}$, for $0<\delta<1$. Then taking a subsequence, as $\varepsilon_k\rightarrow 0^+$,
$$\tilde{w}^{\varepsilon_k}\stackrel{\ast}{\rightharpoonup}\tilde{w}\quad \text{in}\quad L^\infty([0,T];H^4_\mu)
\quad \text{and}\ \tilde{w}^{\varepsilon_k}\rightarrow \tilde{w}\quad \text{in}\quad C([0,T];H^{4-\delta}_{\mu, loc}),$$
where $\tilde{w}=\partial_y \tilde{u}\in L^\infty([0,T];H^4_\mu)$. Applying the local uniform convergence of
$\partial_x^k u^{\varepsilon_k}$, we have following the pointwise convergence of $v^{\varepsilon_k}$: as $\varepsilon_k\rightarrow 0^+,$
\bqa v^{\varepsilon_k}=-\int_0^y \partial_x \tilde{u}^{\varepsilon_k}dy\rightarrow-\int_0^y \partial_x \tilde{u}dy=:v.\label{3.2.57}\eqa

In fact, we also have $$\lim\limits_{y\rightarrow 0}\tilde{u}(t,x,y)=\lim\limits_{y\rightarrow 0}\int_0^y\tilde{w}(t,x,\tilde{y})d\tilde{y}$$
and
$$\tilde{v}=-\int_0^y\tilde{u}_xd\tilde{y}\in L^\infty([0,T_1];L^\infty(\mathbb{R}_{+,y};H^3(\mathbb{R}_x))).$$

We have proven that $\tilde{w}$ is a classical solution to the following vorticity Prandtl equation
\begin{eqnarray}\left\{\begin{array}{ll}
\partial_t \tilde{w}+(u^s+\tilde{u})\partial_x\tilde{w}+\tilde{v}(u_{yy}^s+\partial_y\tilde{w})
=\partial_{y}^2\tilde{w}\ \text{in}\ \mathbb{R}_+^2\times\mathbb{R}_+,\\
\partial_y\tilde{w}=0\ \text{on}\ \mathbb{R}\times\mathbb{R}_+,\\
\tilde{w}|_{t=0}=\tilde{w}_{0}\ \text{on}\ \mathbb{R}_+^2,
\end{array}\right.\non\end{eqnarray}
and $(\tilde{u},\tilde{v})$ is a classical solution to problem \eqref{3.1.3}. Lastly, $(u,v)=(u^s+\tilde{u}, \tilde{v})$ is a classical
solution to \eqref{3.1.1}. Therefore we have proved  Theorem 3.2.1.\hfill $\Box$

\section{Bibliographic Comments}
Since the energy method is suitable for the Navier-Stokes equations, there have been some results on the well-posedness of the 2D Prandtl equations  \eqref{3.0.1} by using energy method under a monotonicity condition on the tangential velocity field. Recently, under a monotonicity condition on the tangential velocity field, Alexandre et al (\cite{AWXY}) used  a direct energy method  instead of using the Crocco transformation to investigate the local well-posedness theory of the solution to the 2D Prandtl equations \eqref{3.0.1} in a weighted Sobolev space. Later on, under the Oleinik$'$s monotonicity assumption, Masmoudi and Wong (\cite{mw}) proved the local existence and uniqueness for the 2D Prandtl system \eqref{3.0.1} in a weighted Sobolev space by the energy method. They used a new method of nonlinear energy estimate  instead of using the Crocco transformation  for the Prandtl system. This new  estimate is based on a cancellation property in the convection terms to overcome the loss of $x$-derivative in the
tangential direction, which is valid under the monotonicity assumption.  The difference from \cite{mw} is that the initial datum is a small perturbation of a monotonic shear flow in \cite{AWXY}. Based on the results of   the local well-posedness for Prandtl equations \eqref{3.0.1} in \cite{AWXY,mw}, Xu and Zhang (\cite{xz}) discussed the global well-posedness  of solutions for the nonlinear Prandtl boundary layer equations \eqref{3.0.1} on the half plane. More specifically, they used the energy method to obtain the existence, uniqueness and stability of solutions in the weighted Sobolev space $H^m_\mu(\mathbb{R}_+^2)\;( m\geq 6)$ when weighted functions $\mu$ is a polynomial function and proved the lifespan $T$ of solutions could be any large  when its initial datum is a perturbation around the monotonic shear profile of small size $e^{-T}$, which is the first result of  global existence of solutions of 2D Prandtl equations \eqref{3.0.1} in a polynomial weighted Sobolev space.

Using the paralinearization technique, Chen, Wang and Zhang (\cite{cwz1}) proved  the local well-posedness of the Prandtl equations \eqref{3.0.1} for monotonic data in the anisotropic Sobolev space with an exponential weighted and low regularity and provided a new  way for the zero-viscosity limit problem of the Navier-Stokes equations with the non-slip boundary condition, who obtained the local  well-posedness of Prandtl equations \eqref{3.0.1} when weighted function $\mu$ is an exponential function in $H^{3,1}_\mu(\mathbb{R}_+^2)\cap H^{1,2}_\mu(\mathbb{R}_+^2)$. To be more precise,  our motivation in this chapter is   to use the energy method to prove the existence  of solutions of Prandtl
equations \eqref{3.0.1} with the non-slip boundary condition when the weighted function $\mu$ is  an exponential function  in $H^4_\mu(\mathbb{R}_+^2)$, which is still an open problem in Sobolev spaces before the present paper \cite{qindong} solved it.

To deal with our problem \eqref{3.0.1}, we have encountered some difficulties. The first difficulty arises from  the loss of horizontal direction $x$-derivative by the
term $v\partial_yu$. The loss of $x$-derivative and the absence of any horizontal diffusion, it is necessary to use  a type of cancellation mechanism in the analysis. To overcome this difficulty, inspired by recent results in \cite{AWXY,mw,xz}, we construct a new unknown function $g_n=\big(\frac{\partial_x^n u}{u_y^s+u_y}\big)_y$ to avoid the loss of $x$ derivative and establish the estimates of $g_n$ in $L^2$-norm, and also give the relationship between the quantities $g_n$ and $\partial^n_x\tilde{w}$ in $L^2$-norm
 by Lemma 3.2.7.  Using the boundary condition of vertical velocity and the divergence free, we obtain the representation formula of vertical velocity $$v(x,y,t)=-\int_0^y\partial_xu(x,\tilde{y},t)d\tilde{y}.$$
The second difficulty is that vertical velocity $v$ creates a loss of $x$-variable, so that the standard energy estimates  fail. To overcome this one, we first use the inequality \eqref{3.0.2} to take the $L^\infty$ norm of normal velocity $v(t,x,y)$, which plays a crucial role in the estimates of this chapter since the estimate of normal velocity $v(t,x,y)$  can not be controlled directly by the norm of the vorticity $w=\partial_y u$, then we use incompressible equation $\eqref{3.0.1}_2$.

In order to prove our desired result (Theorem 3.2.1) on the local-in-time existence of solutions to the 2D Prandtl equations \eqref{3.0.1}, inspired by the ideas of  \cite{mw} and \cite{xz}, we  investigate the parabolic regularized Prandtl equations \eqref{3.2.1}  by constructing an approximate scheme, which keeps the nonlinear
term of the original Prandtl equation \eqref{3.0.1} and nonlinear cancellation properties. Therefore, the  local-in-time existence of solutions
of the original equation \eqref{3.0.1} follows by uniform energy estimates of the approximate solutions. The first duty of this chapter is  to derive uniform energy estimates of the approximates solutions. The second one is to construct a new unknown function $g_n$ and to establish the relationship between the norms of $g_n$ and $\partial_x^m\tilde{w}$ from  Lemma 3.2.7, then we can obtain a priori estimate of solutions independent of $\varepsilon$.

\chapter{Local Well-posedness of Solutions to the 2D Mixed Prandtl Equations in  A Sobolev Space Without Monotonicity and Lower Bound}
In this chapter, we shall investigate the 2D  Prandtl-Shercliff regime equations on the half plane and prove the local existence and uniqueness of solutions  for any initial datum by using the classical energy methods in a Sobolev space. Compared to the existence and uniqueness of solutions to the classical Prandtl
equations where the monotonicity condition of the tangential velocity plays
a key role, this monotonicity condition is not needed for the 2D mixed Prandtl equations. Besides, compared with the existence and uniqueness of solutions to the 2D MHD boundary layer where the initial tangential magnetic field has a lower bound plays
an important role, this lower bound condition is also not needed for the 2D mixed Prandtl equations. In other words, we need neither the monotonicity condition of the tangential velocity nor the initial tangential magnetic field has a lower bound and for any initial datum in this chapter.  As far as we have learned, this is the first result of the 2D mixed Prandtl-Shercliff regime equations in a Sobolev space. The content of this chapter is chosen from \cite{qindong1}.

\section{Introduction}
\setcounter{equation}{0}
The following  Prandtl-Shercliff regime equations were derived in \cite{gvp} by the classical 2D incompressible MHD system
when the physical parameters such as Reynolds number, magnetic Reynolds number and Hartmann number satisfy some constraints in the high Reynolds numbers limit. In this chapter, we shall investigate the local existence and uniqueness of solutions to the following initial boundary value problem for the 2D  Prandtl-Shercliff regime  system in the upper half plane $\mathbb{R}_+^2=\{(x,y):x\in\mathbb{R},\ y\in\mathbb{R}_+\}$,  which reads as
\begin{eqnarray}\left\{\begin{array}{ll}
\partial_t u+u\partial_xu+v\partial_yu-\partial^2_yu=\partial_x b,\\
\partial_xu+\partial_y^2b=0,\\
\partial_xu+\partial_yv=0,
\end{array}\right.\label{4.1.1}\end{eqnarray}
where $(u,v)$ denotes the velocity field of the boundary layer and $b$ stands for the corresponding tangential magnetic component, respectively.

The initial datum of \eqref{4.1.1} is given by
\bqa u|_{t=0}=u_{0}.\label{4.0.2}\eqa

The no-slip boundary conditions are imposed on the velocity  field and magnetic filed
\bqa (u,v,b)|_{y=0}=\textbf{0}.\label{4.0.3}\eqa

Far fields boundary conditions take the form
\bqa \lim\limits_{y\rightarrow+\infty}(u,b)=(\bar{u},\bar{b}),\label{4.0.4}\eqa
where $(\bar{u},\bar{b})$ is a pair of positive constants.

\section{Local Wellposedness of Solutions}
\setcounter{equation}{0}
First, we define some  Sobolev spaces  which will be used throughout the chapter as follows,
$$\|f\|^2_{H^{m,m-1}(\mathbb{R}_+^2)}=\sum\limits_{0\leq\alpha+\beta\leq m,\ \alpha\leq m-1}\|\partial_x^\alpha\partial_y^\beta f\|^2_{L^2(\mathbb{R}_+^2)},$$
$$\|f\|^2_{H^m(\mathbb{R}_+^2)}=\|f\|^2_{H^{m,m-1}(\mathbb{R}_+^2)}+\|\partial_x^m f\|^2_{L^2(\mathbb{R}_+^2)}.$$
Next, we will state the main result as follows.
\begin{Theorem}
Let $|\gamma|\leq m=4$ and $\partial_yu_0\in H^4(\mathbb{R}_+^2)$. Then there exists a constant\\ $0<T=\min\big\{\|w(0)\|^2_{L^2(\mathbb{R}_+^2)},\frac{3}{32C\|w(0)\|^4_{H^4(\mathbb{R}_+^2)}}\big\}$ such that if  $$\|\partial_yu_0\|_{H^4(\mathbb{R}_+^2)}\leq C,$$
then for $0\leq t\leq T$, the systems \eqref{4.1.1} admits a unique solution $(u,v,b)$ in $[0,T]$ such that
$$ u\in L^\infty([0,T]; H^4(\mathbb{R}_+^2)),\ (v,b)\in L^\infty([0,T]; L^\infty(\mathbb{R}_{+}\times H^3(\mathbb{R}))),$$
$$ u\in L^2([0,T]; H^5(\mathbb{R}_+^2)),$$
where  $C$ is a positive constant.
\end{Theorem}

Before  proving  Theorem 4.2.1, we first list an important inequality as follows.
\begin{Lemma}(Agmon inequality \cite{iv}) Let $f\in H^1(\mathbb{R}_+^2)$, then
\bqa\|f\|_{L_x^\infty L^2_y}\leq C\|f\|^{\frac{1}{2}}_{L^2(\mathbb{R}_+^2)}\|f\|^{\frac{1}{2}}_{H^1(\mathbb{R}_+^2)}.\non\eqa
\end{Lemma}

\begin{Lemma}(\cite{mpf})
Let $\eta(t)$ be nonnegative, absolutely continues function on $[0,T]$, which
satisfies for a.e. $t \in(0,T)$, the integral inequality
$$\eta(t)\leq a(t)+b(t)\int_0^tk(s)\eta^m(s)ds,$$
where $m \geq 2$ is a positive integer and $a(t)$, $b(t)$, $k(t)$ are positive, continuous
functions on $[0,T]$, and $\frac{a(t)}{b(t)}$ is a nondecreasing
function. Then for $0 \leq t \leq \alpha_m$,
$$\eta(t)\leq a(t)\Big\{1-
(m-1)\int_0^tk(s)b(s)a^{m-1}(s)ds\Big\}^{\frac{1}{1-m}},$$
where $\alpha_m$ is defined to be the supreme of those $t \in [0,T]$ satisfying
$$\alpha_m=\sup\Big\{t\in[0,T]:1-
(m-1)\int_0^tk(s)b(s)a^{m-1}(s)ds>0 \Big\}.$$
\end{Lemma}

Next, we will prove Theorem 4.2.1 by  the following a series of lemmas. First, we show the boundary value of the $(u,b)$ of higher order partial derivatives in $y$ on $y=0$ for  systems \eqref{4.1.1}.
\begin{Lemma}

Assume that $(u,b)$ is a smooth solution of the system \eqref{4.1.1}, then the boundary values of $(u,b)$ satisfy
\begin{eqnarray}\left\{\begin{array}{ll}
u(t,x,0)=0,\ b(t,x,0)=0,\\
\partial_y^2u(t,x,0)=0,\ \ \partial_y^2b(t,x,0)=0,\\
\partial^4_{y}u(t,x,0)
=\partial_x\partial_yu\partial_yu(t,x,0),\\
\partial^4_{y}b(t,x,0)=0,\\
\partial^6_{y}u(t,x,0)
=(4\partial_x\partial^3_yu+\partial^2_x\partial_yb)\partial_yu(t,x,0)+(\partial_x\partial_yb-\partial^3_{y}u)\partial_x\partial_yu(t,x,0),\\ \partial^6_{y}b(t,x,0)=-\Big(\partial^2_x\partial_yu\partial_yu+(\partial_x\partial_yu)^2\Big)(t,x,0).
\end{array}\right.\label{4.1.2}\end{eqnarray}
\end{Lemma}
\textbf{Proof}. Looking back the boundary conditions in \eqref{4.0.3},
\begin{eqnarray}
u(t,x,0)=0, \ v(t,x,0)=0,\ b(t,x,0)=0, \  (t,x)\in [0,T]\times\mathbb{R},
\label{4.1.3}\end{eqnarray}

By use of the equations $(4.1.1)_{3}$, we have,
\begin{eqnarray*}u(t,x,y)=-\int_{-\infty}^x\partial_yv(t,\tilde{x},y)d\tilde{x}\ \text{and}\ \partial^n_xv(t,x,y)=-\int_0^y\partial^{n+1}_xu(t,x,\tilde{y})d\tilde{y},\end{eqnarray*}
thus, by (4.1.3), we get $\partial_yv(t,x,0)=0$, $\partial_xu(t,x,0)=0$ and $\partial^n_xv(t,x,0)=0$.

We know that
\begin{eqnarray*}u_x(t,x,y)=\int_{-\infty}^x\partial^2_xu(t,\tilde{x},y)d\tilde{x}\end{eqnarray*}
which, combined with $\partial_xu(t,x,0)=0$, we have $\partial^2_xu(t,x,0)=0$. Similarly, we can derive $\partial^3_xu(t,x,0)=0$
and $\partial^4_xu(t,x,0)=0$.

By $(4.1.1)_2$, we give
\begin{eqnarray*}b(t,x,y)&=&\int_0^y\int_0^{\tilde{y}}\partial^2_yb(t,x,\bar{y})d\bar{y}d\tilde{y}-\int_0^y\partial_yb(t,x,0)d\tilde{y}\\
&=&-\int_0^y\int_0^{\tilde{y}}\partial_xu(t,x,\bar{y})d\bar{y}d\tilde{y}-\int_0^y\partial_yb(t,x,0)d\tilde{y},\end{eqnarray*}
which gives $\partial^n_xb(t,x,0)=0$. Therefore, we can deduce that
\begin{eqnarray}
\partial^n_x(u,v,b)(t,x,0)=0, \ \ (t,x)\in [0,T]\times\mathbb{R}, 0\leq n\leq 4.\label{4.1.4}\end{eqnarray}

Applying \eqref{4.1.1} and the boundary conditions \eqref{4.1.4}, we have
\bqa \partial^2_{y}u|_{y=0}=0,\ \ \partial^2_{y}b|_{y=0}=0.\label{4.1.5} \eqa
Besides, we also have
\bqa\partial^n_x \partial^2_{y}u|_{y=0}=0,\ \ \partial^n_x\partial^2_{y}b|_{y=0}=0,  \ \ (t,x)\in [0,T]\times\mathbb{R},\ \ 0\leq n\leq 4.\label{4.1.6} \eqa
Differentiating the equation of $\eqref{4.1.1}_1$ with respect to $y$,
\bqa  \partial_t\partial_yu-\partial^3_{y}u+\partial_y
\big(u\partial_xu+v\partial_yu\big)=
\partial_x\partial_yb,\label{4.1.7}\eqa
and using the boundary conditions \eqref{4.1.3}-\eqref{4.1.5}, then we have
\bqa \partial_t\partial_yu|_{y=0}=\partial^3_{y}u|_{y=0}+\partial_x\partial_yb|_{y=0}.\label{4.1.8}\eqa
Differentiating \eqref{4.1.1} with respect to $y$ twice, we arrive at
\bqa\partial_t\partial^2_yu-\partial^4_{y}u+\partial^2_y
\big(u\partial_xu+v\partial_yu\big)
=\partial_x\partial^2_yb.\label{4.1.9}\eqa
Now appealing to the Leibniz formula
\bqa\partial^2_y
\big(u\partial_xu+v\partial_yu\big)
=\partial^2_yu\partial_xu+\partial^2_yv\partial_yu+u\partial_x\partial^2_yu+v\partial^3_yu+2\partial_yu\partial_x\partial_yu+2\partial_yv\partial^2_yu,\label{4.1.10}\eqa
we conclude,
\bqa \partial^4_{y}u(t,x,0)=\partial_x\partial_yu\partial_yu(t,x,0).\label{4.1.11}\eqa

Differentiating \eqref{4.1.11} with respect to $t$ and using \eqref{4.1.8}, we get
\bqa \partial_t\partial^4_{y}u(t,x,0)
&=&\partial_t\partial_{y}u\partial_x\partial_yu(t,x,0)+\partial_yu\partial_t\partial^2_{xy}u(t,x,0)\non\\
&=&(\partial^3_{y}u+\partial_x\partial_yb)\partial_x\partial_yu(t,x,0)+(\partial_x\partial^3_{y}u+\partial^2_x\partial_yb)\partial_yu(t,x,0).
\label{4.1.12}\eqa

Similarly, differentiating \eqref{4.1.1} four times with respect to $y$, we get
\bqa\partial_t\partial^4_yu-\partial^6_{y}u+\partial^4_y
\big(u\partial_xu+v\partial_yu\big)
=\partial_x\partial^4_yb,\label{4.1.13}\eqa
and invoking the Leibniz formula, we deduce
\bqa\partial^4_y
\big(u\partial_xu+v\partial_yu\big)
&=&\partial^4_yu\partial_xu+\partial^4_yv\partial_yu+4\partial^3_yu\partial_x\partial_yu+4\partial^3_yv\partial^2_yu\non\\
&&+6\partial^2_yu\partial_x\partial^2_yu+6\partial^2_yv\partial^3_yu+4\partial_yu\partial_x\partial^3_yu+4\partial_yv\partial^4_yu\non\\
&&+u\partial_x\partial^4_yu+v\partial^5_yu,\label{}\eqa
which, combined with \eqref{4.1.12}, gives
\bqa \partial^6_{y}u|_{y=0}
=(4\partial_x\partial^3_yu+\partial^2_x\partial_yb)\partial_yu(t,x,0)+(\partial_x\partial_yb-\partial^3_{y}u)\partial_x\partial_yu(t,x,0).\label{4.1.15}\eqa
Analogously, we can arrive at
\bqa \partial^6_{y}b(t,x,0)=-\Big(\partial^2_x\partial_yu\partial_yu+(\partial_x\partial_yu)^2\Big)(t,x,0).\label{4.1.16}\eqa
The proof is therefore complete.\hfill$\Box$

\subsection{ Existence of Solutions}

In this subsection, we shall investigate the existence of solutions of the 2D mixed Prandtl equation  by using classical energy method. We are going to investigate the existence of solutions to the equations \eqref{4.1.1} by using the following system of vorticity $w=\partial_yu$,  which reads as
\begin{eqnarray}\left\{\begin{array}{ll}
\partial_t w+u\partial_xw+v\partial_yw-\partial^2_yw=\partial_x\partial_y b,\ \ (x,y,t)\in\mathbb{R}_+^2\times T, \\
\partial_xu+\partial_y^2b=0,\ \ (x,y,t)\in\mathbb{R}_+^2\times T, \\
\partial_xu+\partial_yv=0,\ \ (x,y,t)\in\mathbb{R}_+^2\times T, \\
w|_{t=0}=w(0)=\partial_yu_0,(x,y)\in\mathbb{R}_+^2,\\
\partial_yw|_{y=0}=v|_{y=0}=b|_{y=0}=0, \ \ (x,t)\in\mathbb{R}\times T.
\end{array}\right.\label{4.2.1}\end{eqnarray}
First, we will use the following  lemma,  which is helpful to deal with the boundary value.
\begin{Lemma}(\cite{adam})
Let $1< p< +\infty$. If $U\in W^{m,p}(\mathbb{R}^{n+1})$, then its trace $u$ belongs
to the space $B=B^{m-\frac{1}{p};p,p}(\mathbb{R}^{n})$ and
$$\|u\|_{B}\leq K\|U\|_{W^{m,p}(\mathbb{R}^{n+1})},$$
with the constant $K>0$ independent of $U$.
\end{Lemma}
\textbf{Proof}. See Lemma 3.2.4.\hfill $\Box$
\begin{Corollary}
Let $1< p< +\infty$. If $U\in W^{m,p}(\mathbb{R}^{n+1})$, then its trace $u$ belongs
to the space $W^{m-1,p}(\mathbb{R}^{n})$ and
$$\|u\|_{W^{m-1,p}(\mathbb{R}^{n})}\leq K\|U\|_{W^{m,p}(\mathbb{R}^{n+1})},$$
with the constant $K>0$ independent of $U$.
\end{Corollary}
\textbf{Proof}. See Corollary 3.2.1.\hfill $\Box$

We will prove the existence of solutions to problem \eqref{4.2.1} by the standard energy method and complete the proof by the following two lemmas, where the first one is devoted to the estimate with loss of $x$-derivative of solutions $w$ to problem \eqref{4.2.1}.

\begin{Lemma}Let $|\gamma|\leq m=4$ and under the assumptions of Theorem 4.2.1, for a smooth solution $w$ to system \eqref{4.2.1}, then
\bqa&& \frac{d}{dt}\|w\|^2_{H^{m,m-1}(\mathbb{R}^2_+)}+\|\partial_xu \|^2_{H^{m,m-1}(\mathbb{R}^2_+)}+\|\partial_yw\|^2_{H^{m,m-1}(\mathbb{R}^2_+)}\non\\
&&\qquad\leq C(1+\|w\|^6_{H^4(\mathbb{R}^2_+)}),\label{4.2.4}\eqa
where  $C$ is a positive constant.
\end{Lemma}
\textbf{Proof}.
For $|\gamma|=\alpha=\beta=0$,  multiplying \eqref{4.2.1} by $w$ and integrating it by parts over $\mathbb{R}_+^2$, we conclude that
\bqa &&\frac{d}{dt}\|w\|^2_{L^2(\mathbb{R}_+^2)}+\|\partial_yw\|^2_{L^2(\mathbb{R}_+^2)}=\int_{\mathbb{R}_+^2}\partial_x\partial_y b wdydx=\int_{\mathbb{R}_+^2}\partial_x\partial^2_y b udydx=-\|\partial_xu\|^2_{L^2(\mathbb{R}_+^2)},\non\eqa
which implies
\bqa&& \frac{d}{dt}\|w\|^2_{L^2(\mathbb{R}_+^2)}+\|\partial_xu\|^2_{L^2(\mathbb{R}_+^2)}+\|\partial_yw\|^2_{L^2(\mathbb{R}_+^2)}=0.\label{4.2.3}\eqa
For $1\leq\alpha+\beta=m$ and $\alpha\leq m-1$ respectively, applying the operator $\partial^\gamma=\partial_x^\alpha\partial_y^\beta$ on $\eqref{4.2.1}_1$ , multiplying
the resulting equation by $\partial^\gamma w$ and integrating it by parts over $\mathbb{R}_+^2$,  we arrive at
\bqa&& \frac{d}{dt}\sum_{|\gamma|\leq m,\ \alpha\leq m-1}\|\partial^\gamma w\|^2_{L^2(\mathbb{R}_+^2)}+\sum_{|\gamma|\leq m,\ \alpha\leq m-1}\|\partial_y\partial^\gamma w\|^2_{L^2(\mathbb{R}_+^2)}\non\\
&&=\sum_{|\gamma|\leq m,\ \alpha\leq m-1}\Big(\int_{\mathbb{R}}\partial_y\partial^\gamma w
\partial^\gamma w|_{y=0}dx+\int_{\mathbb{R}_+^2}\Big([\partial^\gamma, u]\partial_xw+[\partial^\gamma, v]\partial_yw\Big)\partial^\gamma w dxdy\non\\
&&\qquad+\int_{\mathbb{R}_+^2}\partial_x\partial_y \partial^\gamma b \partial^\gamma \partial_yudxdy\Big)\non\\
&&=\sum_{i=1}^4J_i.\label{4.2.5}\eqa
Actually,  for the commutator, we can write it as follows
\bqa [\partial^\gamma, u]\partial_xw=\sum_{\gamma_1\leq\gamma,\ 1\leq |\gamma_1|}C_\gamma^{\gamma_1} \partial^{\gamma_1} u \partial^{\gamma-\gamma_1}\partial_xw.\non\eqa
Then, for $|\gamma|\leq4$, applying Lemma 4.2.2 (the Agmon inequality), we deduce
\bqa \|[\partial^\gamma, u]\partial_xw\|_{L^2(\mathbb{R}_+^2)}&\leq&C\sum_{\gamma_1\leq\gamma,\ 1\leq |\gamma_1|}\|\partial^{\gamma_1}u\|_{L^\infty(\mathbb{R}\times L^2(\mathbb{R}_+))}\|\partial^{\gamma-\gamma_1+1}w\|_{L^2(\mathbb{R}\times L^\infty(\mathbb{R}_+))}\non\\
&\leq& C\sum_{\gamma_1\leq\gamma,\ 1\leq |\gamma_1|}\|\partial^{\gamma_1}u\|^\frac{1}{2}_{L^2(\mathbb{R}_+^2)}\|\partial_x\partial^{\gamma_1}u\|^\frac{1}{2}_{L^2(\mathbb{R}_+^2)}\|\partial^\gamma w\|^\frac{1}{2}_{L^2(\mathbb{R}_+^2)}\|\partial_y\partial^\gamma w\|^\frac{1}{2}_{L^2(\mathbb{R}_+^2)}\non\\
&\leq& C\sum_{\gamma_1\leq\gamma,\ 1\leq |\gamma_1|}\|\partial^{\gamma_1}u\|_{L^2(\mathbb{R}_+^2)}\|\partial_x\partial^{\gamma_1}u\|_{L^2(\mathbb{R}_+^2)}\non\\
&&\quad+\|\partial^\gamma w\|_{L^2(\mathbb{R}_+^2)}\|\partial_y\partial^\gamma w\|_{L^2(\mathbb{R}_+^2)}.\label{4.2.6}\eqa
Similarly, we also have
\bqa &&\|[\partial^\gamma, v]\partial_yw\|_{L^2(\mathbb{R}_+^2)}\|\partial^\gamma w\|_{L^2(\mathbb{R}_+^2)}\non\\
&&\quad\leq\sum_{\gamma_1\leq\gamma,\ 1\leq |\gamma_1|}\|\partial^{\gamma_1}v\|_{L^2(\mathbb{R}\times L^\infty(\mathbb{R}_+))}\|\partial^{\gamma-\gamma_1}\partial_yw\|_{L^\infty(\mathbb{R}\times L^2(\mathbb{R}_+))}
\|\partial^\gamma w\|_{L^2(\mathbb{R}_+^2)}\non\\
&&\quad\leq C\sum_{\gamma_1\leq\gamma,\ 1\leq |\gamma_1|}\|\partial_x\partial^{\gamma_1}u\|_{L^2(\mathbb{R}_+^2)}\|\partial^{\gamma}w\|^\frac{3}{2}_{L^2(\mathbb{R}_+^2)}
\|\partial_x\partial_y\partial^{\gamma-\gamma_1} w\|^\frac{1}{2}_{L^2(\mathbb{R}_+^2)}\non\\
&&\quad\leq \frac{1}{4}\sum_{\gamma_1\leq\gamma,\ 1\leq |\gamma_1|}\Big(\|\partial_x\partial^{\gamma_1}u\|^2_{L^2(\mathbb{R}_+^2)}+\frac{1}{4}\|\partial_y\partial^\gamma w\|^2_{L^2(\mathbb{R}_+^2)}\Big)+C\|\partial^\gamma w\|^6_{L^2(\mathbb{R}_+^2)}.\label{4.2.7}\eqa
Firstly, although $\alpha$ and $\beta$ have many various cases, there are  same estimates on $j_2$ and $j_3$ by using \eqref{4.2.6}-\eqref{4.2.7} and the Cauchy-Schwarz inequality, that is,
\bqa|j_2|+|j_3|\leq \sum_{|\gamma|\leq m,\ \alpha\leq m-1}\Big(C\|\partial^{\gamma}w\|^6_{L^2(\mathbb{R}_+^2)}+\frac{1}{4}\|\partial_x\partial^{\gamma}u\|^2_{L^2(\mathbb{R}_+^2)}+\frac{1}{4}\|\partial_y\partial^\gamma w\|^2_{L^2(\mathbb{R}_+^2)}\Big).\label{4.10}\eqa
Next, we use the boundary conditions of $\partial_yw$ and $v$ to deal with $j_1$. After a complicated
analysis and calculation, we can get following cases that do not vanish on the boundary  $y=0$.

Case 1: $\alpha=1,\ \beta=2$. We use the boundary value of $\partial_y^3{w}$ on $y=0$, the Young inequality and Lemma 4.2.5 to conclude that
\bqa \Big|\int_{\mathbb{R}}\partial_x\partial^3_yw
\partial_x\partial^2_yw|_{y=0}dx\Big|&\leq& \|\partial_x(\partial_x\partial_yu\partial_yu))|_{y=0}\|_{L^2(\mathbb{R})}
\|\partial_x\partial^2_yw|_{y=0}\|_{L^2(\mathbb{R})}\non\\
&\leq&\|\partial_xw|_{y=0}\|^2_{L^2(\mathbb{R})}
\|\partial_x\partial^2_yw|_{y=0}\|_{L^2(\mathbb{R})}\non\\
&\leq&C\|\partial_x\partial_yw\|^2_{L^2(\mathbb{R}_+^2)}
\|\partial_x\partial^3_yw\|_{L^2(\mathbb{R}_+^2)}\non\\
&\leq&C\|w\|^6_{H^4(\mathbb{R}_+^2)}+C.\non\eqa

Case 2: $\alpha=0,\ \beta=4$. We use the boundary value of $\partial_y^5{w}$ on $y=0$, the Young inequality and Lemma 4.2.5 to deduce that
\bqa \Big|\int_{\mathbb{R}}\partial^5_yw
\partial^4_yw|_{y=0}dx\Big|&\leq& \|\Big((\partial_x\partial^3_yu+\partial^2_x\partial_yb)\partial_yu+(\partial_x\partial_yb+\partial^3_yu)\partial_x\partial_yu\Big)|_{y=0}\|_{L^2(\mathbb{R})}
\|\partial^5_yw\|_{L^2(\mathbb{R}_+^2)}\non\\
&\leq&C\|w\|^4_{H^4(\mathbb{R}_+^2)}+\frac{1}{12}\|\partial_y^5w\|^2_{L^2(\mathbb{R}_+^2)}.
\non\eqa

Case 3: $\alpha=1,\ \beta=3$. We use the boundary value of $\partial_y^3w$ on $y=0$, the Young inequality, Lemma 4.2.5 and integrate it by part with respect to $x$ variable to obtain that
\bqa \Big|\int_{\mathbb{R}}\partial_x\partial^4_yw
\partial_x\partial^3_yw|_{y=0}dx\Big|&=&\Big|\int_{\mathbb{R}}\partial^4_yw
\partial^2_x\partial^3_yw|_{y=0}dx\Big|\non\\
&\leq& \|\partial^2_x(\partial_x\partial_yu\partial_yu)|_{y=0}\|_{L^2(\mathbb{R})}
\|\partial^5_yw\|_{L^2(\mathbb{R}_+^2)}\non\\
&\leq&C\|w\|^4_{H^4(\mathbb{R}_+^2)}+\frac{1}{12}\|\partial_y^5w\|^2_{L^2(\mathbb{R}_+^2)}.\non\eqa

Case 4: $\alpha=2,\ \beta=2$. We use the boundary value of $\partial_y^3w$ on $y=0$, the Young inequality, Lemma 4.2.5 and integrate it by part with respect to $x$ variable to derive that
\bqa \Big|\int_{\mathbb{R}}\partial^2_x\partial^3_yw
\partial^2_x\partial^2_yw|_{y=0}dx\Big|
&\leq& \|\partial^2_x(\partial_x\partial_yu\partial_yu)|_{y=0}\|_{L^2(\mathbb{R})}
\|\partial^2_x\partial^3_yw\|_{L^2(\mathbb{R}_+^2)}\non\\
&\leq&C\|w\|^4_{H^4(\mathbb{R}_+^2)}+\frac{1}{12}\|\partial^2_x\partial_y^3w\|^2_{L^2(\mathbb{R}_+^2)}.\non\eqa
Summing up, we can derive that
\bqa|j_1|\leq C+C\|w\|^6_{H^4(\mathbb{R}_+^2)}+\frac{1}{4}\|\partial_yw\|^2_{H^4(\mathbb{R}_+^2)}.\label{4.2.8}\eqa
Finally, we deal with the term $j_4$ by integrating it by parts over $\mathbb{R}_+^2$,
\bqa \int_{\mathbb{R}_+^2}\partial_x\partial_y \partial^\gamma b \partial^\gamma \partial_yudxdy&=&
\int_{\mathbb{R}}\partial_x\partial_y \partial^\gamma b \partial^\gamma u|_{y=0}dx-\int_{\mathbb{R}_+^2}\partial_x\partial^2_y \partial^\gamma b \partial^\gamma udxdy\non\\
&=&\int_{\mathbb{R}}\partial_x\partial_y \partial^\gamma b \partial^\gamma u|_{y=0}dx-\|\partial_x\partial^\gamma u\|^2_{L^2(\mathbb{R}_+^2)}.\label{4.2.9}\eqa
Similarly to the estimate on $j_1$, we need to discuss the term $\int_{\mathbb{R}}\partial_x\partial_y \partial^\gamma b \partial^\gamma u|_{y=0}dx$ for the different types of  $\alpha$ and $\beta$, we can get the following case  not vanishing on the boundary of $y=0$ by a complicated analysis and calculation,
\bqa \Big|\int_{\mathbb{R}}\partial_x\partial^5_y  b \partial^4_yu|_{y=0}dx\Big|&=&\Big|\int_{\mathbb{R}}\partial_x\partial^3_yu\partial_x\partial^4_yu|_{y=0}dx\Big|\non\\
&\leq&C\|\partial_x\partial^3_yu|_{y=0}\|_{L^2(\mathbb{R})}\|\partial_x\partial^4_yu|_{y=0}\|_{L^2(\mathbb{R})}\non\\
&\leq& C\|w\|^2_{H^4(\mathbb{R}_+^2)}+\frac{1}{4}\|\partial_yw\|^2_{H^4(\mathbb{R}_+^2)},\label{4.2.10}\eqa
where we  have used  Lemma 4.2.5.

Therefore, we attain the following inequality,
\bqa |j_4|\leq C+ C\|w\|^6_{H^4(\mathbb{R}_+^2)}+\frac{1}{4}\|\partial_yw\|^2_{H^4(\mathbb{R}_+^2)}-\sum_{|\gamma|\leq m,\ \alpha\leq m-1}\|\partial_x\partial^\gamma u\|^2_{L^2(\mathbb{R}_+^2)}.\label{4.2.11}\eqa
Inserting above inequalities \eqref{4.10}-\eqref{4.2.8} and \eqref{4.2.11} into \eqref{4.2.5}, we conclude
\bqa \sum_{|\gamma|\leq m,\ \alpha\leq m-1}\Big(\frac{d}{dt}\|\partial^\gamma w\|^2_{L^2(\mathbb{R}_+^2)}+\|\partial_x\partial^\gamma u\|^2_{L^2(\mathbb{R}_+^2)}+\|\partial_y\partial^\gamma w\|^2_{L^2(\mathbb{R}_+^2)}\Big)\leq C+C\|w\|^6_{H^4(\mathbb{R}_+^2)}.\non\eqa
The proof is thus complete. \hfill$\Box$

We have completed the estimate with loss of $x$-derivative of solutions $w$ to problem \eqref{4.2.1} in previous section. In order to ensure the completion of energy estimate, we need to establish the estimate of term $\partial_x^\alpha w$,
which will be done in the following lemma.
\begin{Lemma}Let $\alpha\leq m=4$ and under the assumptions of Theorem 4.2.1, for a smooth solution $w$ to problem \eqref{4.2.1}, then
\bqa&& \frac{d}{dt}\sum_{\alpha\leq m}\|\partial^\alpha_xw\|^2_{L^2(\mathbb{R}^2_+)}+\sum_{\alpha\leq m}\|\partial_x\partial_x^\alpha u\|^2_{L^2(\mathbb{R}_+^2)}+\sum_{\alpha\leq m}\|\partial_y\partial^\alpha_xw\|^2_{L^2(\mathbb{R}^2_+)}\leq C\| w\|^6_{H^4(\mathbb{R}^2_+)},\label{4.2.12b}\eqa
where  $C$ is a positive constant.
\end{Lemma}
\textbf{Proof}. Applying the operator $\partial^\alpha_x$ on $\eqref{4.2.1}_1$, multiplying
the resulting equation by $\partial_x^\alpha w$ and integrating it by parts over $\mathbb{R}_+^2$,  we arrive at
\bqa&& \frac{d}{dt}\|\partial_x^\alpha w\|^2_{L^2(\mathbb{R}^2_+)}+\|\partial_x\partial_x^\alpha u\|^2_{L^2(\mathbb{R}^2_+)}+\|\partial_y\partial_x^\alpha w\|^2_{L^2(\mathbb{R}^2_+)}\non\\
&&=\int_{\mathbb{R}_+^2}\Big([\partial_x^\alpha, u]\partial_xw+[\partial_x^\alpha, v]\partial_yw\Big)\partial_x^\alpha w dxdy,\label{4.2.12}\eqa
where we have used the fact
$$(\partial_x^m\partial_yw)(t,x,0)=0,\ (t,x)\in[0,T]\times\mathbb{R}.$$
In fact,  for the commutator, we can write it as follows
\bqa [\partial_x^\alpha, u]\partial_xw=\sum_{\alpha_1\leq\alpha,\ 1\leq |\alpha_1|}C_\alpha^{\alpha_1} \partial_x^{\alpha_1} u \partial_x^{\alpha-\alpha_1}\partial_xw.\non\eqa
Then, for $|\alpha|\leq4$, applying Lemma 4.2.2 (the Agmon inequality) and the Young inequality, we deduce
\bqa \|[\partial_x^\alpha, u]\partial_xw\|_{L^2(\mathbb{R}_+^2)}&\leq& \sum_{\alpha_1\leq\alpha,\ 1\leq |\alpha_1|}\|\partial_x^{\alpha_1}u\|_{L^\infty(\mathbb{R}\times L^2(\mathbb{R}_+))}\|\partial_x^{\alpha-\alpha_1+1}w\|_{L^2(\mathbb{R}\times L^\infty(\mathbb{R}_+))}\non\\ &\leq& C\sum_{\alpha_1\leq\alpha,\ 1\leq |\alpha_1|}\|\partial_x^\alpha u\|^{\frac{1}{2}}_{L^2(\mathbb{R}_+^2)}\|\partial_x\partial_x^\alpha u\|^{\frac{1}{2}}_{L^2(\mathbb{R}_+^2)}\|\partial_x^{\alpha-\alpha_1+1}w\|^{\frac{1}{2}}_{L^2(\mathbb{R}_+^2)}
\|\partial_x^{\alpha-\alpha_1+1}\partial_yw\|^{\frac{1}{2}}_{L^2(\mathbb{R}_+^2)}\non\\
&\leq& C\sum_{\alpha_1\leq\alpha,\ 1\leq |\alpha_1|}\|\partial_x^\alpha u\|_{L^2(\mathbb{R}_+^2)}\|\partial_x\partial_x^\alpha u\|_{L^2(\mathbb{R}_+^2)}+C\|\partial_y\partial_x^\alpha w\|_{L^2(\mathbb{R}_+^2)}\|\partial_x^\alpha w\|_{L^2(\mathbb{R}_+^2)}.\non\eqa
Hence, \bqa \|[\partial_x^\alpha, u]\partial_xw\|_{L^2(\mathbb{R}_+^2)}\|\partial_x^\alpha w\|_{L^2(\mathbb{R}_+^2)}\leq
\frac{1}{4}\|\partial_x\partial_x^\alpha u\|^2_{L^2(\mathbb{R}_+^2)}+\frac{1}{4}\|\partial_y\partial_x^\alpha w\|^2_{L^2(\mathbb{R}_+^2)}+C\| w\|^6_{H^4(\mathbb{R}_+^2)}.\label{4.2.13}\eqa

Similarly, we also have
\bqa &&\|[\partial_x^\alpha, v]\partial_yw\|_{L^2(\mathbb{R}_+^2)}\|\partial_x^\alpha w\|_{L^2(\mathbb{R}_+^2)}\non\\
&&\quad\leq
\sum_{\alpha_1\leq\alpha,\ 1\leq |\alpha_1|}\|\partial_x^{\alpha_1}v\|_{L^2(\mathbb{R}\times L^\infty(\mathbb{R}_+))}
\|\partial_x^{\alpha-\alpha_1}\partial_yw\|_{L^\infty(\mathbb{R}\times L^2(\mathbb{R}_+))}\|\partial_x^\alpha w\|_{L^2(\mathbb{R}_+^2)}\non\\
&&\quad\leq C\sum_{\alpha_1\leq\alpha,\ 1\leq |\alpha_1|}\|\partial_x\partial_x^\alpha u\|_{L^2(\mathbb{R}_+^2)}\|\partial_x^{\alpha-\alpha_1}\partial_yw\|^{\frac{1}{2}}_{L^2(\mathbb{R}_+^2)}
\|\partial_x^{\alpha-\alpha_1+1}\partial_yw\|^{\frac{1}{2}}_{L^2(\mathbb{R}_+^2)}\|\partial_x^\alpha w\|_{L^2(\mathbb{R}_+^2)}\non\\
&&\quad\leq\frac{1}{4}\|\partial_x\partial_x^\alpha u\|^2_{L^2(\mathbb{R}_+^2)}+\frac{1}{4}\|\partial_y\partial_x^\alpha w\|^2_{L^2(\mathbb{R}_+^2)}+C\| w\|^6_{H^4(\mathbb{R}_+^2)}.\label{4.2.14}\eqa
Inserting \eqref{4.2.13}-\eqref{4.2.14} into \eqref{4.2.12} and summing over $\alpha$, we can conclude that
\bqa\frac{d}{dt}\sum_{\alpha\leq m}\|\partial_x^\alpha w\|^2_{L^2(\mathbb{R}_+^2)}+\sum_{\alpha\leq m}\|\partial_x\partial_x^\alpha u\|^2_{L^2(\mathbb{R}_+^2)}+\sum_{\alpha\leq m}\|\partial_y\partial_x^\alpha w\|^2_{L^2(\mathbb{R}_+^2)}\leq C\| w\|^6_{H^4(\mathbb{R}_+^2)}.\non\eqa
The proof is thus complete. \hfill$\Box$

Up to now, we have derived the main estimates of the solution $u$. Now combining Lemma 4.2.6 with Lemma 4.2.7 with $m=4$, we have the following estimate
\bqa \frac{d}{dt}\|w\|^2_{H^4(\mathbb{R}_+^2)}+\|\partial_xu \|^2_{H^4(\mathbb{R}_+^2)}+\|\partial_yw\|^2_{H^4(\mathbb{R}_+^2)}\leq C(1+\|w\|^6_{H^4(\mathbb{R}_+^2)}).\label{4.2.16}\eqa
Applying Lemma 4.2.3 to above inequality \eqref{4.2.16}, we have
\bqa \|w\|^2_{H^4(\mathbb{R}_+^2)}\leq\{\|w(0)\|^2_{H^4(\mathbb{R}_+^2)}+t\}\{1-2C(\|w(0)\|^2_{H^4(\mathbb{R}^2_+)}+t)^2t\}^{-\frac{1}{2}},\label{4.2.17}\eqa
if $1-2C(\|w(0)\|^2_{H^4(\mathbb{R}^2_+)}+t)^2t>0$.

Taking $T=\min\big\{\|w(0)\|^2_{H^4(\mathbb{R}_+^2)},\frac{3}{32C\|w(0)\|^4_{H^4(\mathbb{R}^2_+)}}\big\}$ such that
$$\|w(0)\|^2_{H^4(\mathbb{R}_+^2)}+t\leq 2\|w(0)\|^2_{H^4(\mathbb{R}_+^2)}, \ \{1-2C(\|w(0)\|^2_{H^4(\mathbb{R}^2_+)}+t)^2t\}^{-\frac{1}{2}}\leq 2,$$
which, combined with \eqref{4.2.17},  gives
\bqa \|w\|^2_{H^4(\mathbb{R}_+^2)}\leq 4\|w(0)\|^2_{H^4(\mathbb{R}_+^2)}.\label{4.a2.17}\eqa
Hence
\bqa \|w\|^2_{H^4(\mathbb{R}_+^2)}+\int_0^t(\|\partial_xu \|^2_{H^4(\mathbb{R}_+^2)}+\|\partial_yw\|^2_{H^4(\mathbb{R}_+^2)})ds\leq C_1,\ 0\leq t\leq T,\label{4.2.18}\eqa
where constant $C_1>0$ depends on $T$ and $\|w(0)\|^2_{H^4(\mathbb{R}_+^2)}$.

According to equations $\eqref{4.2.1}_{2,3}$ and \eqref{4.2.17}, we can conclude $\|(v,b)\|_{L^\infty([0,T];L^\infty({\mathbb{R}_+}\times H^3(\mathbb{R})))}\leq C$.
We therefore have proved the  existence of solutions in Theorem 4.2.1.\hfill$\Box$

\subsection{The Uniqueness of Local Solutions}
In this subsection, we will prove the uniqueness of solutions to the 2D mixed Prandtl equations.
Let $(w_1,b_1)$ and $(w_2,b_2)$ be two solutions to problem \eqref{4.2.1} with same initial data $w_1(0,x,y)=w_2(0,x,y)\in H^4(\mathbb{R}_+^2)$. Set $w=w_1-w_2$ and $b=b_1-b_2$. Then
\begin{eqnarray}\left\{\begin{array}{ll}
\partial_t w+u\partial_xw_1+v\partial_yw_1+u_2\partial_xw+v_2\partial_yw-\partial^2_yw=\partial_x\partial_y b,\ \ (x,y,t)\in\mathbb{R}_+^2\times T,\\
\partial_xu+\partial_y^2b=0,\ \ (x,y,t)\in\mathbb{R}_+^2\times T,\\
\partial_xu+\partial_yv=0,\ \ (x,y,t)\in\mathbb{R}_+^2\times T, \\
\partial_yu(0,x,y)=0,\ \ (x,y)\in\mathbb{R}_+^2, \\
\partial_yw|_{y=0}=v|_{y=0}=b|_{y=0}=0,\ \ (x,y,t)\in\mathbb{R}_+^2\times T.
\end{array}\right.\label{4.3.1}\end{eqnarray}

Similarly to Lemma 4.2.4, we can derive the following  initial data and boundary values,
\begin{eqnarray}\left\{\begin{array}{ll}
\partial_yu(0,x,y)=0,\ (x,y)\in\mathbb{R}_+^2,\\
u(t,x,0)=0,\ b(t,x,0)=0,\\
\partial_y^2u(t,x,0)=0,\ \ \partial_y^2b(t,x,0)=0,\\
\partial^4_{y}u(t,x,0)=\partial_x\partial_yu_1\partial_yu(t,x,0)+\partial_x\partial_yu\partial_yu_2(t,x,0),\\
\partial^4_{y}b(t,x,0)=0,\\
\partial^6_{y}u|_{y=0}
=(4\partial_x\partial^3_yu_1+\partial^2_x\partial_yb_1)\partial_yu(t,x,0)+(\partial_x\partial_yb_1-\partial^3_{y}u_1)\partial_x\partial_yu(t,x,0),\\
+(4\partial_x\partial^3_yu+\partial^2_x\partial_yb)\partial_yu_2(t,x,0)+(\partial_x\partial_yb-\partial^3_{y}u)\partial_x\partial_yu_2(t,x,0),\\
\partial^6_{y}b(t,x,0)=-\Big(\partial^2_x\partial_yu_1\partial_yu+\partial_x\partial_yu_1\partial_x\partial_yu+\partial^2_x\partial_yu\partial_yu_2
+\partial_x\partial_yu\partial_x\partial_yu_2\Big)(t,x,0).
\end{array}\right.\label{4.3.2}\end{eqnarray}
\textbf{Proof}. Looking back the boundary condition in $\eqref{4.3.1}_5$,
\begin{eqnarray}
u(t,x,0)=0, \ v(t,x,0)=0,\ b(t,x,0)=0, \  (t,x)\in [0,T]\times\mathbb{R},
\label{4.3.3}\end{eqnarray}
thus the following results are obvious
\begin{eqnarray}
\partial^n_x(u,b)(t,x,0)=0, \ \partial^n_xv(t,x,0)=0, \ \ (t,x)\in [0,T]\times\mathbb{R}, 0\leq n\leq 4.\label{4.3.4}\end{eqnarray}
Appealing to $\eqref{4.3.1}_1$ and the boundary condition \eqref{4.3.4}, we have
\bqa \partial^2_{y}u|_{y=0}=0,\ \ \partial^2_{y}b|_{y=0}=0.\label{4.3.5} \eqa
Besides, we also have
\bqa\partial^n_x \partial^2_{y}u|_{y=0}=0,\ \ \partial^n_x\partial^2_{y}b|_{y=0}=0,  \ \ (t,x)\in [0,T]\times\mathbb{R},\quad 0\leq n\leq 4.\label{4.3.6} \eqa
Differentiating the equation of $\eqref{4.3.1}_1$ with $y$, we get
\bqa  \partial_t\partial_yu-\partial^3_{y}u+\partial_y
\big(u\partial_xu_1+v\partial_yu_1+u_2\partial_xu+v_2\partial_yu\big)=
\partial_x\partial_yb,\label{4.3.7}\eqa
and using the boundary values \eqref{4.3.4}, then we have
\bqa \partial_t\partial_yu|_{y=0}=\partial^3_{y}u|_{y=0}+\partial_x\partial_yb|_{y=0}.\label{4.3.8}\eqa
Differentiating $\eqref{4.3.1}_1$ with respect to $y$ twice, it follows
\bqa\partial_t\partial^2_yu-\partial^4_{y}u+\partial^2_y
\big(u\partial_xu_1+v\partial_yu_1+u_2\partial_xu+v_2\partial_yu\big)
=\partial_x\partial^2_yb.\label{4.3.9}\eqa
Applying the Leibniz formula, we arrive at
\bqa&&\partial^2_y
\big(u\partial_xu_1+v\partial_yu_1+u_2\partial_xu+v_2\partial_yu\big)\non\\
&&=\partial^2_yu\partial_xu_1+\partial^2_yv\partial_yu_1+u\partial_x\partial^2_yu_1+v\partial^3_yu_1+2\partial_yu\partial_x\partial_yu_1+2\partial_yv\partial^2_yu_1\non\\
&&+\partial^2_yu_2\partial_xu+\partial^2_yv_2\partial_yu+u_2\partial_x\partial^2_yu+v_2\partial^3_yu+2\partial_yu_2\partial_x\partial_yu
+2\partial_yv_2\partial^2_yu.\label{4.3.10}\eqa
Therefore,
\bqa \partial^4_{y}u(t,x,0)=\partial_x\partial_yu_1\partial_yu(t,x,0)+\partial_x\partial_yu\partial_yu_2(t,x,0).\label{4.3.11}\eqa
Differentiating \eqref{4.3.11} with respect to $t$ and using \eqref{4.3.8}, we get
\bqa \partial_t\partial^4_{y}u(t,x,0)
&=&\partial_t\partial_{y}u\partial_x\partial_yu_1(t,x,0)+\partial_yu\partial_t\partial_x\partial_yu_1(t,x,0)\non\\
&&+\partial_t\partial_{y}u_2\partial_x\partial_yu(t,x,0)+\partial_yu_2\partial_t\partial_x\partial_yu(t,x,0)\non\\
&=&(\partial^3_{y}u+\partial_x\partial_yb)\partial_x\partial_yu_1(t,x,0)+(\partial_x\partial^3_{y}u_1+\partial^2_x\partial_yb_1)\partial_yu(t,x,0)\non\\
&&+(\partial^3_{y}u_2+\partial_x\partial_yb_2)\partial_x\partial_yu(t,x,0)+(\partial_x\partial^3_{y}u+\partial^2_x\partial_yb)\partial_yu_2(t,x,0).
\label{4.3.12}\eqa
Similarly, differentiating $\eqref{4.3.1}_1$ with respect to $y$ four times, we have
\bqa\partial_t\partial^4_yu-\partial^6_{y}u+\partial^4_y
\big(u\partial_xu_1+v\partial_yu_1+u_2\partial_xu+v_2\partial_yu\big)
=\partial_x\partial^4_yb,\label{4.3.13}\eqa
and applying the Leibniz formula, we conclude
\bqa&&\partial^4_y
\big(u\partial_xu_1+v\partial_yu_1+u_2\partial_xu+v_2\partial_yu\big)\non\\
&&=\partial^4_yu\partial_xu_1+\partial^4_yv\partial_yu_1+4\partial^3_yu\partial_x\partial_yu_1+4\partial^3_yv\partial^2_yu_1
+6\partial^2_yu\partial_x\partial^2_yu_1+6\partial^2_yv\partial^3_yu_1\non\\
&&+4\partial_yu\partial_x\partial^3_yu_1+4\partial_yv\partial^4_yu_1+u\partial_x\partial^4_yu_1+v\partial^5_yu_1+\partial^4_yu_2\partial_xu+\partial^4_yv_2\partial_yu
+4\partial^3_yu_2\partial_x\partial_yu\non\\
&&+4\partial^3_yv_2\partial^2_yu+6\partial^2_yu_2\partial_x\partial^2_yu+6\partial^2_yv_2\partial^3_yu+4\partial_yu_2\partial_x\partial^3_yu+4\partial_yv_2\partial^4_yu
+u_2\partial_x\partial^4_yu+v_2\partial^5_yu,\non\\
\label{4.3.14}\eqa
 which, combined with \eqref{4.3.12}, yields to
\bqa \partial^6_{y}u|_{y=0}
&=&(4\partial_x\partial^3_yu_1+\partial^2_x\partial_yb_1)\partial_yu(t,x,0)+(\partial_x\partial_yb_1-\partial^3_{y}u_1)\partial_x\partial_yu(t,x,0)\non\\
&&+(4\partial_x\partial^3_yu+\partial^2_x\partial_yb)\partial_yu_2(t,x,0)+(\partial_x\partial_yb-\partial^3_{y}u)\partial_x\partial_yu_2(t,x,0).\label{4.3.15}\eqa
Analogously, we can conclude
\bqa \partial^6_{y}b(t,x,0)=-\Big(\partial^2_x\partial_yu_1\partial_yu+\partial_x\partial_yu_1\partial_x\partial_yu+\partial^2_x\partial_yu\partial_yu_2
+\partial_x\partial_yu\partial_x\partial_yu_2\Big)(t,x,0).\label{4.3.16}\eqa
The proof is thus complete.\hfill$\Box$
\begin{Lemma}Let $|\gamma|\leq m-1=3$ and assume $u_1$ and $u_2$ are two solutions obtained in  Theorem 4.2.1 with respect to the initial data $u_1(0,x,y)$ and $u_2(0,x,y)$, then
\bqa&& \frac{d}{dt}\|w\|^2_{H^{m-1,m-2}(\mathbb{R}^2_+)}+\|\partial_xu \|^2_{H^{m-1,m-2}(\mathbb{R}^2_+)}+\|\partial_yw\|^2_{H^{m-1,m-2}(\mathbb{R}^2_+)}\non\\
&&\qquad\leq C(1+\|u_2\|^2_{H^{m-1,m-2}(\mathbb{R}^2_+)}+\|\partial_y w_1\|^2_{H^4(\mathbb{R}^2_+)}+\|\partial_y w_2\|^2_{H^4(\mathbb{R}^2_+)})\|w\|^2_{H^3(\mathbb{R}^2_+)},\label{4.3.17}\eqa
where  $C$ is a positive constant.
\end{Lemma}
\textbf{Proof}.  The proof of the lemma is similar to that of Lemma 4.2.6
for $0\leq\alpha+\beta=m-1$ and $\alpha\leq m-2$ respectively. Applying the operator $\partial^\gamma=\partial^\alpha_x\partial^\beta_y$ on $\eqref{4.3.1}_1$ , multiplying
the resulting equation by $\partial^\gamma w$ and integrating it by parts over $\mathbb{R}_+^2$,  we arrive at
\bqa&& \frac{d}{dt}\|\partial^\gamma w\|^2_{L^2(\mathbb{R}^2_+)}+\|\partial_y\partial^\gamma w\|^2_{L^2(\mathbb{R}^2_+)}\non\\
&&=\int_{\mathbb{R}}\partial_y\partial^\gamma w
\partial^\gamma w|_{y=0}dx+\int_{\mathbb{R}_+^2}\Big([\partial^\gamma, u_2]\partial_xw+[\partial^\gamma, v_2]\partial_yw\Big)\partial^\gamma w dxdy\non\\
&&\qquad+\int_{\mathbb{R}_+^2}\Big(\partial^\gamma(u\partial_xw_1)+\partial^\gamma(v\partial_yw_1)\Big)\partial^\gamma w dxdy+\int_{\mathbb{R}_+^2}\partial_x\partial_y \partial^\gamma b \partial^\gamma \partial_yudxdy\non\\
&&=\sum_{i=1}^6I_i.\label{4.3.18}\eqa
Actually, we can write  the commutator as follows
\bqa [\partial^\gamma, u_2]\partial_xw=\sum_{\gamma_1\leq\gamma,\ 1\leq |\gamma_1|}C_\gamma^{\gamma_1} \partial^{\gamma_1} u_2 \partial^{\gamma-\gamma_1}\partial_xw.\non\eqa
Then, for $|\gamma|\leq3$, applying the Sobolev inequality, we deduce
\bqa \|[\partial^\gamma, u_2]\partial_xw\|_{L^2(\mathbb{R}^2_+)}&\leq&C\sum_{\gamma_1\leq\gamma,\ 1\leq |\gamma_1|}\|\partial^{\gamma_1}u_2\|_{L^\infty(\mathbb{R}\times L^2(\mathbb{R}_+))}\|\partial^{\gamma-\gamma_1+1}w\|_{L^2(\mathbb{R}\times L^\infty(\mathbb{R}_+))}\non\\
&\leq& C\sum_{\gamma_1\leq\gamma,\ 1\leq |\gamma_1|}\|\partial_x\partial^{\gamma_1}u_2\|_{L^2(\mathbb{R}^2_+)}\|\partial_y\partial^\gamma w\|_{L^2(\mathbb{R}^2_+)}.\label{4.3.19}\eqa
Similarly, we also have
\bqa \|[\partial^\gamma, v_2]\partial_yw\|_{L^2(\mathbb{R}^2_+)}\leq C\sum_{\gamma_1\leq\gamma,\ 1\leq |\gamma_1|}\|\partial_x\partial^{\gamma_1}u_2\|_{L^2(\mathbb{R}^2_+)}\|\partial_y\partial^\gamma w\|_{L^2(\mathbb{R}^2_+)}.\label{4.3.20}\eqa
Firstly, although $\alpha$ and $\beta$ have many various cases, there are  same estimates on $I_2$ and $I_3$ by using \eqref{4.3.19}-\eqref{4.3.20} and the Cauchy-Schwarz inequality, that is,
\bqa |I_2|+|I_3|\leq\frac{1}{4}\|\partial_y\partial^\gamma w\|^2_{L^2(\mathbb{R}^2_+)}+C\|\partial_x\partial^{\gamma}u_2\|^2_{L^2(\mathbb{R}^2_+)}\|w\|^2_{H^3(\mathbb{R}^2_+)}\label{4.c1}.\eqa
For the items $I_4$ and $I_5$,  we can deduce the following estimates by using the Sobolev inequality,
\bqa |I_4|&\leq&C\sum_{\gamma_1\leq\gamma}\|\partial^{\gamma_1}u\|_{L^\infty(\mathbb{R}\times L^2(\mathbb{R}_+))}
\|\partial^{\gamma-\gamma_1}\partial_xw_1\|_{L^2(\mathbb{R}\times L^\infty(\mathbb{R}_+))}\non\\
&\leq& C\sum_{\gamma_1\leq\gamma}\|\partial_x\partial^{\gamma_1}u\|_{L^2(\mathbb{R}^2_+)}\|\partial_x\partial_y\partial^{\gamma-\gamma_1} w_1\|_{L^2(\mathbb{R}^2_+)}.\non\eqa
Analogously,
\bqa|I_5|\leq C\sum_{\gamma_1\leq\gamma}\|\partial_x\partial^{\gamma_1}u\|_{L^2(\mathbb{R}^2_+)}\|\partial_x\partial_y\partial^{\gamma-\gamma_1} w_1\|_{L^2(\mathbb{R}^2_+)}.\non\eqa
Hence,
\bqa |I_4|+|I_5|
\leq\frac{1}{4} \sum_{\gamma_1\leq\gamma}\|\partial_x\partial^{\gamma_1}u\|^2_{L^2(\mathbb{R}^2_+)}+C(\|\partial_y w_1\|^2_{H^4(\mathbb{R}^2_+)}+\|\partial_y w_2\|^2_{H^4(\mathbb{R}^2_+)})\|\partial^\gamma w\|^2_{L^2(\mathbb{R}^2_+)}.\label{4.c2}\eqa
Next, we use the boundary conditions of $w$ and $v$ to deal with $I_1$. After a complicated
analysis and calculation, we can get following case do not vanish on the boundary  $y=0$.

Case: $\alpha=1,\ \beta=2$. We use the boundary value of $\partial_y^3{w}$ on $y=0$ and Lemma 4.2.5 to conclude that
\bqa \Big|\int_{\mathbb{R}}\partial_x\partial^3_yw
\partial_x\partial^2_yw|_{y=0}dx\Big|
\leq C(1+\|w_1\|^2_{H^4(\mathbb{R}^2_+)}+\|w_2\|^2_{H^4(\mathbb{R}^2_+)})\|w\|^2_{H^3(\mathbb{R}^2_+)}.\label{4.c3}\eqa

Finally, we deal with the term $I_6$ by integrating it by parts over $\mathbb{R}_+^2$. Similarly to the estimate on $I_1$, we need to discuss the term $\int_{\mathbb{R}}\partial_x\partial_y \partial^\gamma b \partial^\gamma u|_{y=0}dx$ for the different types of  $\alpha$ and $\beta$. However,  the boundary values of $y=0$  vanish by a complicated analysis and calculation,
\bqa \int_{\mathbb{R}_+^2}\partial_x\partial_y \partial^\gamma b \partial^\gamma \partial_yudxdy&=&
\int_{\mathbb{R}}\partial_x\partial_y \partial^\gamma b \partial^\gamma u|_{y=0}dx-\int_{\mathbb{R}_+^2}\partial_x\partial^2_y \partial^\gamma b \partial^\gamma udxdy\non\\
&=&-\|\partial_x\partial^\gamma u\|^2_{L^2(\mathbb{R}^2_+)}.\label{4.3.21}\eqa

Inserting \eqref{4.c1}-\eqref{4.3.21} into \eqref{4.3.18} and summing $\gamma$, we can derive the desired result \eqref{4.3.17}.
Thus the proof is complete.\hfill $\Box$

\begin{Lemma}Let $|\alpha|\leq m-1=3$ and assume $u_1$, $u_2$ be two solutions obtained in  Theorem 4.2.1 with respect to the initial data $u_1(0,x,y)$, $u_2(0,x,y)$, then we have
\bqa&& \frac{d}{dt}\sum_{\alpha\leq m}\|\partial^\alpha_xw\|^2_{L^2(\mathbb{R}^2_+)}+\sum_{\alpha\leq m}\|\partial_x\partial_x^\alpha u\|^2_{L^2(\mathbb{R}^2_+)}+\sum_{\alpha\leq m}\|\partial_y\partial^\alpha_xw\|^2_{L^2(\mathbb{R}^2_+)}\non\\
&&\leq C(\|\partial_y w_1\|^2_{H^4(\mathbb{R}^2_+)}+\|\partial_y w_2\|^2_{H^4(\mathbb{R}^2_+)}+\|\partial_x\partial_x^\alpha u_2\|^2_{L^2(\mathbb{R}^2_+)})\| w\|^2_{H^3(\mathbb{R}^2_+)},\label{4.3.22}\eqa
where  $C$ is a positive constant.
\end{Lemma}
\textbf{Proof}. Acting the operator $\partial^\alpha_x$ on \eqref{4.3.1}, multiplying
the resulting equation by $\partial_x^\alpha w$ and integrating it by parts over $\mathbb{R}_+^2$,  we arrive at
\bqa&& \frac{d}{dt}\|\partial_x^\alpha w\|^2_{L^2(\mathbb{R}^2_+)}+\|\partial_x\partial_x^\alpha u\|^2_{L^2(\mathbb{R}^2_+)}+\|\partial_y\partial_x^\alpha w\|^2_{L^2(\mathbb{R}^2_+)}\non\\
&&=\int_{\mathbb{R}_+^2}\Big([\partial_x^\alpha, u_2]\partial_xw+[\partial_x^\alpha, v_2]\partial_yw+\partial_x^\alpha(u\partial_xw_1)+\partial_x^\alpha( v\partial_yw_1)\Big)\partial_x^\alpha w dxdy,\label{4.3.23}\eqa
where we have used the fact
$$(\partial_x^\alpha\partial_yw)(t,x,0)=0,\ (t,x)\in[0,T]\times\mathbb{R}.$$
In fact,  for the commutator, we can write  as follows
\bqa [\partial_x^\alpha, u_2]\partial_xw=\sum_{\alpha_1\leq\alpha,\ 1\leq |\alpha_1|}C_\alpha^{\alpha_1} \partial_x^{\alpha_1} u_2 \partial_x^{\alpha-\alpha_1}\partial_xw.\non\eqa
Then, for $|\alpha|\leq3$, applying the Sobolev inequality, we deduce
\bqa \|[\partial_x^\alpha, u_2]\partial_xw\|_{L^2(\mathbb{R}^2_+)}&\leq& \sum_{\alpha_1\leq\alpha,\ 1\leq |\alpha_1|}\|\partial_x^{\alpha_1}u\|_{L^\infty(\mathbb{R}\times L^2(\mathbb{R}_+))}\|\partial_x^{\alpha-\alpha_1+1}w\|_{L^2(\mathbb{R}\times L^\infty(\mathbb{R}_+))}\non\\
&\leq&C\sum_{\alpha_1\leq\alpha,\ 1\leq |\alpha_1|} \|\partial_x\partial_x^\alpha u_2\|_{L^2(\mathbb{R}^2_+)}\|\partial_y\partial_x^\alpha w\|_{L^2(\mathbb{R}^2_+)}.\non\eqa
Similarly, we also have
\bqa \|[\partial_x^\alpha, v_2]\partial_yw\|_{L^2(\mathbb{R}^2_+)}\leq C\|\partial_x\partial_x^\alpha u_2\|_{L^2(\mathbb{R}^2_+)}\|\partial_y\partial_x^\alpha w\|_{L^2(\mathbb{R}^2_+)}.\non\eqa
Hence, \bqa&&(\|[\partial_x^\alpha, u_2]\partial_xw\|_{L^2(\mathbb{R}^2_+)}+\|[\partial_x^\alpha, v_2]\partial_yw\|_{L^2(\mathbb{R}^2_+)})\|\partial_x^\alpha w\|^2_{L^2(\mathbb{R}^2_+)}\non\\
&&\leq \frac{1}{4}
\|\partial_y\partial_x^\alpha w\|_{L^2(\mathbb{R}^2_+)}+C\|\partial_x\partial_x^\alpha u_2\|^2_{L^2(\mathbb{R}^2_+)}\|\partial_x^\alpha w\|^2_{L^2(\mathbb{R}^2_+)},\label{4.3.24}\eqa
and
\bqa &&(\|\partial_x^\alpha(u\partial_xw_1)\|_{L^2(\mathbb{R}^2_+)}+ \|\partial_x^\alpha( v\partial_yw_1)\|_{L^2(\mathbb{R}^2_+)})\|\partial_x^\alpha w\|^2_{L^2(\mathbb{R}^2_+)}\non\\
&&\leq\frac{1}{4} \|\partial_x\partial^{\alpha}u\|_{L^2(\mathbb{R}^2_+)}+C(\|\partial_y w_1\|^2_{H^4(\mathbb{R}^2_+)}+\|\partial_y w_2\|^2_{H^4(\mathbb{R}^2_+)})\|\partial_x^\alpha w\|^2_{L^2(\mathbb{R}^2_+)}.\label{4.3.25}\eqa
Inserting \eqref{4.3.24}-\eqref{4.3.25} into \eqref{4.3.23} and summing $\alpha$, we can conclude that
\bqa&&\frac{d}{dt}\sum_{\alpha\leq m}\|\partial_x^\alpha w\|^2_{L^2(\mathbb{R}^2_+)}+\sum_{\alpha\leq m}\|\partial_x\partial_x^\alpha u\|^2_{L^2(\mathbb{R}^2_+)}+\sum_{\alpha\leq m}\|\partial_y\partial_x^\alpha w\|^2_{L^2(\mathbb{R}^2_+)}\non\\
&&\qquad\leq C(\|\partial_y w_1\|^2_{H^4(\mathbb{R}^2_+)}+\|\partial_y w_2\|^2_{H^4(\mathbb{R}^2_+)}+\|\partial_x\partial_x^\alpha u_2\|^2_{L^2(\mathbb{R}^2_+)})\| w\|^2_{H^3(\mathbb{R}^2_+)}.\non\eqa
The proof is thus complete.\hfill $\Box$

Finally, combining Lemma 4.2.8 with Lemma 4.2.9, we can derive the following inequality
\bqa&& \frac{d}{dt}\|w\|^2_{H^3(\mathbb{R}^2_+)}+\|\partial_xu \|^2_{H^3(\mathbb{R}^2_+)}+\|\partial_yw\|^2_{H^3(\mathbb{R}^2_+)}\non\\
&&\qquad\leq C(1+\| u_2\|^2_{H^4(\mathbb{R}^2_+)}+\|\partial_y w_1\|^2_{H^4(\mathbb{R}^2_+)}+\|\partial_y w_2\|^2_{H^4(\mathbb{R}^2_+)})\|w\|^2_{H^3(\mathbb{R}^2_+)}.\label{4.3.26}\eqa
Applying the Gronwall inequality to \eqref{4.3.26} and using \eqref{4.2.18}, we can obtain
\bqa \|w\|^2_{H^3(\mathbb{R}^2_+)}&\leq& \|w_0\|^2_{H^3(\mathbb{R}^2_+)}\exp\Big(\int_0^tC\big(1+\|u_2\|^2_{H^4(\mathbb{R}^2_+)}+\|\partial_y w_1\|^2_{H^4(\mathbb{R}^2_+)}+\|\partial_y w_2\|^2_{H^4(\mathbb{R}^2_+)}ds\big)\Big)\non\\
&=&0,\ 0\leq t\leq T,\label{4.3.27}\eqa
which yields $w_1(t,x,y)=w_2(t,x,y)$. The proof of uniqueness of solutions is thus complete.\hfill$\Box$

Combining \eqref{4.2.18} with \eqref{4.3.27}, we can complete the proof of Theorem 4.2.1.\hfill$\Box$

\section{Bibliographic Comments}
In this section, let us  briefly review some known
results to the problem \eqref{4.1.1}. For mixed Prandtl  boundary equations, Xie and Yang (\cite{xy1}) investigated the global existence and uniqueness of solutions with analytic regularity of the 2D mixed Prandtl  boundary layer  equations by using the standard energy methods in analytic space. As far as we have learned, there have been no result for the 2D mixed Prandtl-Shercliff regime equations, this is the first result of the $2D$ mixed Prandtl-Shercliff regime equations in a Sobolev space. Compared to the existence and uniqueness of solutions to the classical Prandtl equations where the monotonicity condition of the tangential velocity plays a key role, this monotonicity condition is not needed for the 2D mixed Prandtl equations. Compared with the existence and uniqueness of solutions to the 2D MHD boundary layer where the initial tangential magnetic field has a lower bound $\delta>0$ plays an important role, this lower bound condition is also not needed for the 2D mixed Prandtl equations. In other words, we need neither the monotonicity condition of the tangential velocity nor the initial tangential magnetic field has a lower bound $\delta>0$ and for any initial datum in this chapter. Besides, compared with the existence and uniqueness of solutions to the classical 2D Prandtl and MHD equations where the initial data are small enough and the vorticity satisfies decay condition, these  conditions are also not needed for the 2D mixed Prandtl equations. Similarly to the classical Prandtl equations, the difficulty of solving the problem \eqref{4.1.1} in the Sobolev framework is the loss of $x$-derivative in the term $v\partial_yu$. We can insert the equations $\eqref{4.1.1}_2$ into $\eqref{4.1.1}_1$ in calculating to eliminate this difficulty.  In this chapter, we investigate the existence and uniqueness of solutions to the 2D mixed Prandtl-Shercliff regime system for any initial datum by using the direct energy methods.

\chapter{Local well-posedness    of solutions to 2D magnetic Prandtl model in the Prandtl-Hartmann regime}
In this chapter, we shall consider the 2D  magnetic Prandtl equation in the Prandtl-Hartmann regime in a periodic domain and prove the local existence and uniqueness of solutions by the energy method in a polynomial weighted  Sobolev space. On the one hand, we have noted that the $x$-derivative of the pressure $P$ plays a key role in all known results on the existence and uniqueness of solutions to the Prandtl-Hartmann regime equations, in which the case of favorable $P$ ( $\partial_x P<0$) or the case of $\partial_x P=0$ (led by constant outer flow $U=\text{constant} $) was only considered. While in this chapter, we have no any restriction on the sign of $\partial_x P$, which has generalized all previous results and gives definitely rise to a difficulty in mathematical treatments. To overcome this difficulty, we shall use the skill of cancellation mechanism which is valid under the monotonicity assumption.
One the other hand, we shall consider  the general outer flow $U\neq\text{constant}$, leading to the boundary data at $y=0$ being much more complicated. To deal with these boundary data, some more delicate estimates and mathematical induction method will be used. Moreover, our result has given us a physical understanding that the outer flow has a stabilizing effect on the Prandtl-Hartmann regime boundary layer in mathematics. The content of this chapter is selected from \cite{qinwang2}.
\section{Introduction}
The system of magnetohydrodynamics is a fundamental system to describe the fluid under the influence of electromagnetic field. The following mixed Prandtl and Hartmann boundary layer equations arise from the incompressible MHD system when the physical parameters such as Reynolds number, magnetic Reynolds number and the Hartmann number satisfy some constraints in the high Reynolds numbers limit, which were proposed in \cite{GSS,gvp}.

In this chapter, we shall consider the  Prandtl-Hartmann regime equations  in
a periodic domain $\mathbb{T}\times \mathbb{R}^{+}:=\{(t,x,y)\big|t>0, x\in \mathbb{R}/\mathbb{Z}, y\in\mathbb{R_{+}}\}$:
\begin{equation}\left\{
\begin{array}{ll}
\partial _{t}u+u\partial _{x}u+v\partial _{y}u= \partial _{y}b+\partial _{y}^{2}u-\partial _{x}P,\\
\partial _{y}u+\partial _{y}^{2}b=0, \\
\partial _{x}u+\partial _{y}v=0,\\
\lim\limits_{y\rightarrow+\infty}u(t,x,y)=U(t,x),\quad  \lim\limits_{y\rightarrow+\infty}b(t,x,y)=B(t,x).
\end{array}
 \label{5.1.1}         \right.\end{equation}
Here
$(u,v)=(u(t,x,y),v(t,x,y))$ denotes the velocity field, $b(t,x,y)$ is the corresponding tangential magnetic component.  The given scalar pressure $P:=P(x,t)$ and the outer flow $U$ satisfy
the well-known Bernoulli's law:
\begin{eqnarray*}
\partial_{t}U+U\partial_{x}U+\partial_{x}P=0 .
\end{eqnarray*}
The initial data and no-slip boundary condition are imposed by
\begin{equation}\left\{
\begin{array}{ll}
u(0,x,y)=u_{0}(x,y),\\
u(t,x,0)=0,\ \ v(t,x,0)=0,\ \ b(t,x,0)=0.
\end{array}
 \label{5.1.2}         \right.\end{equation}

Noticing $(\ref{5.1.1})_{4}$, and integrating $(\ref{5.1.1})_{2}$ in $y$ over $[y, +\infty)$, we obtain
\begin{eqnarray}
\partial_{y}b=U-u ,\quad  b=B+\int_{y}^{+\infty}(u-U)dy. \label{5.1.3}
\end{eqnarray}
Inserting $(\ref{5.1.3})$ into $(\ref{5.1.1})$, then we arrive at the following equations
\begin{equation}\left\{
\begin{array}{ll}
\partial _{t}u+u\partial _{x}u+v\partial _{y}u=(U-u) +\partial _{y}^{2}u-\partial _{x}P,\\
\partial _{x}u+\partial _{y}v=0.
\end{array}
 \label{5.1.4}         \right.\end{equation}
Let the vorticity $w=\partial_{y}u$, then equations $(\ref{5.1.4})$ reduce to the following vorticity system
\begin{equation}\left\{
\begin{array}{ll}
\partial _{t}w+u\partial _{x}w+v\partial _{y}w=-w+\partial _{y}^{2}w,\\
w(0,x,y)=\partial _{y} u_{0},\\
\partial _{y}w|_{y=0}=\partial _{x}P-U.
\end{array}
 \label{5.1.5}         \right.\end{equation}

In this chapter, we shall consider system $(\ref{5.1.1})$ under the Oleinik's monotonicity assumption $w=\partial_{y}u>0$.

The chapter is structured as follows. Section 5.1 is the introduction, which contains many inequalities used in this chapter. In Section 5.2, we shall give the uniform estimates with the weighted norm, our proof is based on using estimates of regularized parabolic equation and maximal principle. The difficulty of this part is the estimates on the boundary at $y=0$, and we shall give the regular pattern eventually by mathematical induction method. In Section 5.3, we shall prove local existence and uniqueness of solutions to the Prandtl system.

We now introduce the weighted Sobolev space and define the space $H^{s,\gamma}_{\sigma,\delta}$ by
$$H^{s,\gamma}_{\sigma,\delta}:=\left\{w:\mathbb{T}\times\mathbb{R_{+}}\rightarrow\mathbb{R}:
\|w\|_{H^{s,\gamma}}<+\infty, (1+y)^{\sigma}|w|\geq\delta,
\sum \limits_{|\alpha|\leq 2}  |(1+y)^{\sigma+\alpha_{2}}D^{\alpha}w | ^{2}\leq\frac{1}{\delta^{2}} \right\},$$
where $D^{\alpha}:=\partial_{x}^{\alpha_{1}}\partial_{y}^{\alpha_{2}}$, $\alpha_{1}+\alpha_{2}=s$,  $s\geq 4$, $\gamma\geq 1$, $\sigma>\gamma+\frac{1}{2}$ and $\delta\in (0,1)$. We define the norm as
 $$\|w\|_{H^{s,\gamma}}^{2}:=\sum \limits_{|\alpha|\leq s} \|(1+y)^{\gamma+\alpha_{2}}D^{\alpha}w\|_{L^{2}}^{2}$$
and
$$\|w\|_{H^{s,\gamma}_{g}}^{2}:=\|(1+y)^{\gamma}g_{s}\|^{2}_{L^{2}}+
\sum\limits_{\substack{ |\alpha|\leq s\\ \alpha_{1}\leq s-1}}  \|(1+y)^{\gamma+\alpha_{2}}D^{\alpha}w\|_{L^{2}}^{2},$$
here
$$g_{s}:=\partial_{x}^{s}w-\frac{\partial_{y}w}{w}\partial_{x}^{s}(u-U).$$
In this chapter, for convenience, we simply write
$$\iint\cdot:=\int_{\mathbb{T}}\int_{\mathbb{R_{+}}}\cdot dxdy.$$

Last, we state our main result as follows.
\begin{Theorem}\label{t1.1}
Given any   integer  $s \geq 4$, and real numbers $\gamma, \sigma, \delta$ satisfying $\gamma\geq 1$, $\sigma>\gamma+\frac{1}{2}$ and $\delta\in (0,1)$.   Assume the following conditions hold,

(i) suppose that the initial data $u_0-U(0,x) \in H^{s, \gamma-1}$ and $\partial_y u_0 \in H^{s,\gamma}_{\sigma, 2\delta}$ satisfy the compatibility conditions $u_0|_{y=0}$ and $\lim \limits_{y \rightarrow +\infty} u=U|_{t=0}$. In addition, when $s = 4$,   assume  that $\delta \geq 0$ is chosen small enough such that $\|\omega_0\|_{H^{s,\gamma} } \leq C \delta^{-1}$ with a generic constant $C>0$.

(ii) the outer flow $U$  satisfies
 \begin{eqnarray}
\sup \limits_{t} \sum\limits_{l=0}^{\frac{s}{2}+1}\|\partial_t^l U\|_{H^{s-2l+2}(\mathbb{T})} < + \infty.
\label{5.1.8}
\end{eqnarray}

Then there exists a time $T := T(s, \gamma, \sigma, \delta,  \|w_0\|_{H^{s,\gamma} },U)$ such that the initial boundary value problem
(\ref{5.1.1})-(\ref{5.1.2}) has a unique local classical solution $(u, v, b)$ satisfying
$$u-U \in L^{\infty}([0,T];H^{s,\gamma-1}) \cap C([0,T]; H^{s}-w)$$
and
$$ \partial_y u \in L^{\infty}([0,T];H^{s,\gamma}_{\sigma, \delta}) \cap C([0,T]; H^{s}-w),$$
where $H^s-w$ is the space $H^s$ endowed with its weak topology.
\end{Theorem}

\subsection{Preliminaries}
In this subsection, we will give some lemmas to our main result.
\begin{Lemma}(\cite{mw})\label{p5.4.1}
Let $s\geq 4$ be an integer, $\gamma\geq 1$, $\sigma>\gamma+\frac{1}{2}$ and $\delta\in (0,1)$. Then for any $w\in H^{s,\gamma}_{\sigma,\delta}(\mathbb{T}\times\mathbb{R}_{+})$,
we have the following inequality
 \begin{eqnarray}
C_{\delta }\|w\|_{H^{s,\gamma}_{g}}\leq  \|w\|_{H^{s,\gamma} }+ \|u-U\|_{H^{s,\gamma-1} } \leq C_{s, \gamma,\sigma,\delta  }\left(\|w\|_{H^{s,\gamma}_{g}}+\|\partial_{x}^{s}U\|_{{L^2}(\mathbb{T})}\right)
  ,
\label{5.4.1}
\end{eqnarray}
where $C_{s, \gamma,\sigma,\delta  }>0$ is a constant and only depends on $s, \gamma,\sigma,\delta  $.
\end{Lemma}

\begin{Lemma}(Hardy type inequality, Lemma B.1 in \cite{mw})\label{y5.4.1}
Let $f:\mathbb{T} \times\mathbb{R}_{+}\rightarrow \mathbb{R}$, \\
(i)  if $\lambda>-\frac{1}{2}$  and $ \lim \limits_{y\rightarrow +\infty }f(x,y ) =0$,  then
 \begin{eqnarray}
 \|(1+y)^{\lambda}f\|_{L^{2}(\mathbb{T} \times\mathbb{R}_{+})} \leq 
\frac{2}{2\lambda+1} \|(1+y)^{\lambda+1}\partial_{y}f\|_{L^{2}(\mathbb{T} \times\mathbb{R}_{+})};
\label{5.4.10}
\end{eqnarray}
(ii)    if $\lambda<-\frac{1}{2}$  and $ f(x,y )|_{y=0}=0$,  then
 \begin{eqnarray}
 \|(1+y)^{\lambda}f\|_{L^{2}(\mathbb{T} \times\mathbb{R}_{+})} \leq 
 -\frac{2}{2\lambda+1} \|(1+y)^{\lambda+1}\partial_{y}f\|_{L^{2}(\mathbb{T} \times\mathbb{R}_{+})}.
\label{5.4.11}
\end{eqnarray}
\end{Lemma}

\begin{Lemma}(\cite{mw})\label{y5.4.2}
Let $s\geq 4$ be an integer, $\gamma\geq 1$, $\sigma>\gamma+\frac{1}{2}$ and $\delta\in (0,1)$. For any $w\in H^{s,\gamma}_{\sigma,\delta}(\mathbb{T}\times\mathbb{R}_{+})$, then  for $k=0,1,2,\cdots, s$,
we have
 \begin{eqnarray}
 &&\|(1+y)^{\gamma} g_{k}\|_{L^{2}}\leq\|(1+y)^{\gamma}\partial _{x}^{k}w\|_{L^{2}}+\delta^{-2}\|(1+y)^{\gamma-1}\partial _{x}^{k}(u-U)\|_{L^{2}}  ,
 \label{5.4.5}
\end{eqnarray}
and
 \begin{eqnarray}
 &&\|(1+y)^{\gamma}\partial _{x}^{k}w\|_{L^{2}}+\|(1+y)^{\gamma-1}\partial _{x}^{k}(u-U)\|_{L^{2}}\leq C_{\gamma,\sigma,\delta} \left(\|\partial_{x}^{k}U\|_{{L^2}(\mathbb{T})}+\|(1+y)^{\gamma }g_{k}\|_{L^{2}}\right),\label{5.4.4}
\end{eqnarray}
where $C_{s, \gamma,\sigma,\delta  }>0$ is a constant and only depends on $s, \gamma,\sigma,\delta  $.
\end{Lemma}

\begin{Lemma}(Sobolev type inequality, Lemma B.2 in  \cite{mw})\label{y5.4.3}
Let $f:\mathbb{T} \times\mathbb{R}_{+}\rightarrow \mathbb{R}$.  Then there exists a constant $C>0$ such that
 \begin{eqnarray}
 \|f(x,y )\|_{L^{\infty}(\mathbb{T} \times \mathbb{R}_{+})}    \leq C
\left(\|f(x,y)\|_{L^{2}(\mathbb{T}\times \mathbb{R}_{+} )}   +\|\partial_{x}f(x,y)\|_{L^{2}(\mathbb{T}\times \mathbb{R}_{+} )}  +\|\partial_{y }^{2}f(x,y)\|_{L^{2}(\mathbb{T}\times \mathbb{R}_{+} )} \right).
\label{5.4.12}
\end{eqnarray}

\end{Lemma}

\begin{Lemma}(\cite{mw})\label{y5.4.4}
Let $s\geq 4$ is an integer, $\gamma\geq 1$, $\sigma>\gamma+\frac{1}{2}$ and $\delta\in (0,1)$. Then for any $w\in H^{s,\gamma}_{\sigma,\delta}(\mathbb{T}\times\mathbb{R}_{+})$,
we have the following inequalities: \\
$(i)$ for $k=0, 1, 2,  \cdots, s-1 $,
\begin{eqnarray}
 \| \frac{\partial_{x}^{k} v +y\partial_x^{k+1}U}{1+y} \|_{L^{2}} \leq C_{s, \gamma,\sigma,\delta  } (\|w\|_{H^{s,\gamma}_{g}}+\|\partial_{x}^{s}U\|_{{L^2}(\mathbb{T})});
\label{5.4.13}
\end{eqnarray}
(ii) for $k=0, 1, 2,  \cdots, s $,
\begin{eqnarray}
 \|  (1+y)^{\gamma-1}\partial_x^{k}(u-U) \|_{L^{2}} \leq C_{s, \gamma,\sigma,\delta  } (\|w\|_{H^{s,\gamma}_{g}}+\|\partial_{x}^{s}U\|_{{L^2}(\mathbb{T})});
\label{5.4.14}
\end{eqnarray}
(iii) for $k=0, 1, 2,  \cdots, s-2 $,
\begin{eqnarray}
 \|  \frac{\partial_{x}^{k} v }{1+y}  \|_{L^{\infty}} \leq C_{s, \gamma,\sigma,\delta  } (\|w\|_{H^{s,\gamma}_{g}}+\|\partial_{x}^{s}U\|_{{L^2}(\mathbb{T})});
\label{5.4.15}
\end{eqnarray}
(iv) for $k=0, 1, 2,  \cdots, s-1 $,
\begin{eqnarray}
 \|  \partial_{x}^{k}u  \|_{L^{\infty}} \leq C_{s, \gamma,\sigma,\delta  }  (\|w\|_{H^{s,\gamma}_{g}}+\|\partial_{x}^{s}U\|_{{L^2}(\mathbb{T})});
\label{5.4.16}
\end{eqnarray}
(v) for $k=0, 1, 2,  \cdots, s-2$,
\begin{eqnarray}
 \|  (1+y)^{\gamma+k_{\alpha_{2}} }D^{k  }w  \|_{L^{\infty}} \leq C_{s, \gamma,\sigma,\delta  } \|w\|_{H^{s,\gamma}_{g}};
\label{5.4.17}
\end{eqnarray}
(vi) for all $|\alpha| \leq s$,
\begin{eqnarray}
\|(1+y)^{\gamma+\alpha_2}D^{\alpha} w\|_{L^2} \leq
\left\{
    \begin{array}{ll}
        C_{s, \gamma,\sigma,\delta  }  (\|w\|_{H^{s,\gamma}_{g}}+\|\partial_{x}^{s}U\|_{{L^2}(\mathbb{T})}) & if~ \alpha=(s,0), \\
       \|w\|_{H^{s,\gamma}_{g}} & if ~ \alpha \neq (s,0);
    \end{array}
\right.\label{5.4.52}
\end{eqnarray}
(vii) for all $k=0, 1, 2,  \cdots, s$,
\begin{eqnarray}
\|(1+y)^{\gamma}g_k\|_{L^2(\mathbb{T})} \leq
\left\{
    \begin{array}{ll}
        C_{s, \gamma,\sigma,\delta  }  (\|w\|_{H^{s,\gamma}_{g}}+\|\partial_{x}^{s}U\|_{{L^2}(\mathbb{T})}) & if~ k=0, 1, 2,  \cdots, s-1, \\
       \|w\|_{H^{s,\gamma}_{g}} & if ~ k=s,
    \end{array}
\right.\label{5.4.51}
\end{eqnarray}
where $C_{s, \gamma,\sigma,\delta  }>0$ is a constant and only depends on $s, \gamma,\sigma,\delta  $.
\end{Lemma}

\begin{Lemma}(Maximum principle for parabolic equations, Lemma E.1 in  \cite{mw})\label{y5.1.6}
Let $\epsilon \geq 0$. If $H \in C([0,T]; C^{2}(\mathbb{T}\times \mathbb{R}^{+}) \cap C^{1}([0,T]; C^{0}(\mathbb{T}\times \mathbb{R}^{+}))$ is a bounded function that satisfies the differential inequality
\begin{eqnarray*}
\left\{\partial_{t}+b_{1}\partial_{x}+b_2\partial_{y}-\epsilon^{2}\partial_{x}^{2}-\partial_{y}^{2} \right\}H \leq fH \ \ \ in \ [0,T]  \times \mathbb{T} \times \mathbb{R}^{+},
\end{eqnarray*}
where the coefficients $b_{1}$, $b_{2}$, and $f$ are continuous and satisfy
\begin{eqnarray*}
\left\|\frac{b_{2}}{1+y}\right\|_{L^{\infty}([0,T] \times \mathbb{T}\times \mathbb{R}^{+})} < +\infty  \ \ and \ \ \|f\|_{L^{\infty}([0,T] \times \mathbb{T} \times \mathbb{R}^{+})} \leq \lambda,
\end{eqnarray*}
then for any $t \in [0,T]$,
\begin{eqnarray*}
\sup \limits_{\mathbb{T} \times \mathbb{R}^{+}} H(t) \leq \max\left\{e^{\lambda t}\|H(0)\|_{L^{\infty}(\mathbb{T} \times \mathbb{R}^{+})}, \max \limits_{\tau \in [0,T]}\left\{e^{\lambda(t-\tau)}\|H(\tau)\big{|}_{y=0}\|_{L^{\infty}(\mathbb{T})}\right\}\right\}.
\end{eqnarray*}
\end{Lemma}

\begin{Lemma}(Minimum principle for parabolic equations, Lemma E.2 in  \cite{mw})\label{y5.1.7}
Let $\epsilon \geq 0$. If $H \in C([0,T]; C^{2}(\mathbb{T}\times \mathbb{R}^{+}) \cap C^{1}([0,T]; C^{0}(\mathbb{T}\times \mathbb{R}^{+}))$ is a bounded function with
 \begin{eqnarray*}
\kappa(t) :=\min \left\{ \min \limits_{\mathbb{T} \times \mathbb{R}^{+}}H(0), \min \limits_{[0,T] \times \mathbb{T}} H \big{|}_{y=0} \right\} \geq 0
\end{eqnarray*}
and satisfies
\begin{eqnarray*}
\left\{\partial_{t}+b_{1}\partial_{x}+b_2\partial_{y}-\epsilon^{2}\partial_{x}^{2}-\partial_{y}^{2} \right\}H = fH
\end{eqnarray*}
where the coefficients $b_{1}$, $b_{2}$, and $f$ are continuous and satisfy
\begin{eqnarray*}
\left\|\frac{b_{2}}{1+y}\right\|_{L^{\infty}([0,T] \times \mathbb{T}\times \mathbb{R}^{+})} <  +\infty  \ \ and \ \ \|f\|_{L^{\infty}([0,T] \times \mathbb{T} \times \mathbb{R}^{+})} \leq \lambda,
\end{eqnarray*}
then for any $t \in [0,T]$,
\begin{eqnarray*}
\sup \limits_{\mathbb{T} \times \mathbb{R}^{+}} H(t) \geq (1-\lambda t e^{\lambda t}) \kappa(t).
\end{eqnarray*}

\end{Lemma}
\section{Uniform estimates on the regularized system }
In this section, we will estimate the norm $\|w\|_{H^{s,\gamma}_{g}}^{2}$  by the energy method. We consider the regularized equations to the problem $(\ref{5.1.4})$ for any $\epsilon$,
\begin{equation}\left\{
\begin{array}{ll}
\partial _{t}u^\epsilon+u^\epsilon\partial _{x}u^\epsilon+v^\epsilon\partial _{y}u^\epsilon=(U-u^\epsilon)+\epsilon^2\partial _{x}^{2}u^\epsilon +\partial _{y}^{2}u^\epsilon-\partial _{x}P^\epsilon,\\
\partial _{x}u^\epsilon+\partial _{y}v^\epsilon=0,\\
u^{\epsilon}|_{t=0}=u_0,\\
u^\epsilon|_{y=0}=v^\epsilon |_{y=0}=0,\\
\lim\limits_{y\rightarrow+\infty}u^\epsilon(t,x,y)=U(t,x),
\label{5.2.0001}
\end{array}
         \right.\end{equation}
where $P^\epsilon$ and $U$ satisfy a regularized Bernoulli's law:
\begin{eqnarray}
\partial_{t}U+U\partial_{x}U-\epsilon^2\partial^2_xU+\partial_{x}P^\epsilon=0 .
\end{eqnarray}
Then, the regularized vorticity $w^\epsilon=\partial_y u^\epsilon$ satisfies the following regularized vorticity system for any $\epsilon>0$,
\begin{equation}\left\{
\begin{array}{ll}
\partial _{t}w^\epsilon+u^\epsilon\partial _{x}w^\epsilon+v^\epsilon\partial _{y}w^\epsilon=-w^\epsilon+\epsilon^2\partial _{x}^{2}w^\epsilon +\partial _{y}^{2}w^\epsilon,\\
{w^\epsilon}|_{t=0}=w_0=\partial _{y}u_0,\\
\partial_y{w^\epsilon}|_{y=0}=\partial_x P^\epsilon -U,
\label{5.2.0002}
\end{array}
         \right.\end{equation}
where the velocity field $(u^\epsilon,v^\epsilon)$ is given
\begin{eqnarray}
u^\epsilon(t,x,y)=U-\int^{+\infty}_y w(t,x,\widetilde{y})d \widetilde{y}, \quad v^\epsilon(t,x,y)=-\int^{y}_0 \partial_x u(t,x,\widetilde{y})d \widetilde{y}.
\end{eqnarray}
From now on, we drop the superscript $\epsilon$ for simplicity of notations.
\subsection{Estimates on $D^{\alpha}w$ }\label{su5.1}
In this subsection, we will estimate the norm on $D^{\alpha}w$ with weight $(1+y)^{ \gamma+ \alpha_{2}}$, $\alpha=(\alpha_{1},\alpha_{2})$, $|\alpha|\leq s$, and $\alpha_{1}$ should be restricted as $\alpha_{1}\leq s-1$. Otherwise, it is impossible to directly estimate the norm on $D^{\alpha}w$ with $\partial _{x}^{s}w$.
\begin{Lemma}(Reduction of boundary data)\label{y5.2.1}
If $w$ solves $(\ref{5.2.0001})$ and $(\ref{5.2.0002})$, then on the boundary at $y=0$, we have,
\begin{equation}\left\{
\begin{array}{ll}
\partial _{y}w|_{y=0}=\partial _{x}P-U(t,x),\\
\partial _{y}^{3}w|_{y=0}= (\partial _{t}-\epsilon^2\partial_x^2)(\partial _{x}P-U(t,x))+(\partial _{x}P-U(t,x))+w\partial _{x}w|_{y=0}.
\end{array}
         \right.\end{equation}
For any $5 \leq 2k+1 \leq s$, there are some constants $C_k$, $C_{\wedge_{\alpha},k,l,\rho^1,\rho^2,...,\rho^j}$,   not depending on $\epsilon$ or $(u,v,w)$, such that
\begin{eqnarray}
\begin{aligned}
\partial_y^{2k+1}w|_{y=0}&=C_k\sum\limits_{s=0}^{k}(\partial_t-\epsilon^2\partial^2_x)^s(\partial_x P -U)\\
&\quad+\sum\limits_{l=0}^{k-1}\epsilon^{2l}\sum\limits_{j=2}^{\max\{2,k-l\}}\sum\limits_{\rho \in A^j_{k,l}}C_{\wedge_{\alpha},k,l,\rho^1,\rho^2,...,\rho^j}\prod^{j}_{i=1}D^{\rho^i}w \big{|}_{y=0}
\label{5.2.005}
\end{aligned}
\end{eqnarray}
where $A^j_{k,l}:=\{\rho:=(\rho^1,\rho^2,...,\rho^j)\in \mathbb{N}^{2j};3\sum\limits_{i=1}^{j}\sum\limits_{i=1}^{j}\rho^i_2=2k+4l+1,\sum\limits_{i=1}^{j}\rho^i_1 \leq k+2l-1,\sum\limits_{i=1}^{j}\rho_2^j\leq 2k-2l-2, ~and~ |\rho^i| \leq 2k-l-1 ~for~ all~ i=1,2,...,j\}$.
\end{Lemma}
\textbf{Proof}.
According to equation $(\ref{5.2.0001})$ and the boundary condition $(\ref{5.1.2})_{2}$, we get
$$\partial _{y}w|_{y=0}=\partial _{x}P-U(t,x)=K,$$
which implies
\begin{eqnarray}
\begin{aligned}
\partial _{y}^{3}w \big{|}_{y=0}&= \left\{\partial _{t}\partial _{y}w+\partial _{y}(u\partial _{x}w)+\partial _{y}(v\partial _{y}w)+\partial _{y}w-\epsilon^2\partial_y\partial_x^2w\right\} \big{|}_{y=0}\\
&=(\partial _{t}-\epsilon^2\partial_x^2)K+K+w\partial _{x}w \big{|}_{y=0},
\end{aligned}
\end{eqnarray}
and
\begin{eqnarray}
\begin{aligned}
\partial _{y}^{2n+1}w \big{|}_{y=0}&= \left\{\left(\partial _{t}-\epsilon^2\partial_x^2w+1\right)\partial _{y}^{2n-1}w+\partial _{y}^{2n-1}(u\partial _{x}w+v\partial _{y}w)\right\} \big{|}_{y=0}.
 \label{5.2.15}
\end{aligned}
\end{eqnarray}
Hence, for $n=2$,
\begin{eqnarray}
\partial _{y}^{5}w \big{|}_{y=0}
&=& (\partial _{t}-\epsilon^2\partial_x^2)\partial_y^{3}w|_{y=0}+\partial _{y}^{3}w|_{y=0}+(2\partial_y w\partial_x\partial_yw+3w\partial_x\partial_y^2w-2\partial_xw\partial_y^2w) \big{|}_{y=0}\nonumber\\
&=&(\partial _{t}-\epsilon^2\partial_x^2)^2K+(\partial _{t}-\epsilon^2\partial_x^2)K+(\partial _{t}-\epsilon^2\partial_x^2)(w\partial_x w)|_{y=0}+(\partial _{t}-\epsilon^2\partial_x^2)K\nonumber\\
&\quad&+K+w\partial _{x}w|_{y=0}+(2\partial_y w\partial_x\partial_yw+3w\partial_x\partial_y^2w-2\partial_xw\partial_y^2w) \big{|}_{y=0},
 \label{5.2.16}
\end{eqnarray}
here
\begin{equation}\left\{
\begin{array}{ll}
\partial _{t} w|_{y=0}=( \partial _{y}^{2}w-u\partial _{x}w-v\partial _{y}w-w-\epsilon^2\partial_x^2 w)|_{y=0}=(\partial _{y}^{2}w-w+\epsilon^2\partial_x^2 w)\big{|}_{y=0}, \\
\partial _{x}\partial _{t} w|_{y=0} =(\partial _{x}\partial _{y}^{2}w-\partial _{x}w+\epsilon^2\partial_x^3 w)\big{|}_{y=0},
\label{5.2.17}
\end{array}
         \right.\end{equation}
\begin{eqnarray}
(\partial _{t}-\epsilon^2\partial_x^2)(w\partial_x w)\big{|}_{y=0}=
(w\partial_x\partial_y^2w+\partial_xw\partial_y^2w-2w\partial_xw-2\epsilon^2\partial_x\partial_x^2w)\big{|}_{y=0}.
 \label{5.2.18}
\end{eqnarray}
It follows from $(\ref{5.2.16})$, $(\ref{5.2.17})$ and $(\ref{5.2.18})$ that
\begin{eqnarray}
\begin{aligned}
\partial _{y}^{5}w\big{|}_{y=0}
&=C_n\sum\limits_{i=0}^{2}(\partial _{t}-\epsilon^2\partial_x^2)^i K+(2\partial_y w\partial_x\partial_yw+4w\partial_x\partial_y^2w-\partial_xw\partial_y^2w)\big{|}_{y=0} \\
&\quad-w\partial _{x}w \big{|}_{y=0}-2\epsilon^2\partial_xw\partial_x^2w\big{|}_{y=0}.
 \label{5.2.19}
\end{aligned}
\end{eqnarray}
For $n=3$,
\begin{eqnarray}
\partial _{y}^{7}w\big{|}_{y=0}
&=&(\partial _{t}-\epsilon^2\partial_x^2)\partial_y^{5}w\big{|}_{y=0}+\partial_y^{5}w\big{|}_{y=0}+\partial_y^5(u\partial _{x}w+v\partial _{y}w)\big{|}_{y=0}\nonumber \\
&=&C_n\sum\limits_{i=1}^{3}(\partial _{t}-\epsilon^2\partial_x^2)^iK+(\partial _{t}-\epsilon^2\partial_x^2)(2\partial_y w\partial_x\partial_yw+4w\partial_x\partial_y^2w-\partial_xw\partial_y^2w)\big{|}_{y=0}\nonumber \\
&\quad&-(\partial _{t}-\epsilon^2\partial_x^2)(w\partial_x w)\big{|}_{y=0}-(\partial _{t}-\epsilon^2\partial_x^2)2\epsilon^2\partial_xw\partial_x^2w\big{|}_{y=0}\nonumber\\
&\quad&+\partial_y^{5}w\big{|}_{y=0}+\partial_y^5(u\partial _{x}w+v\partial _{y}w)\big{|}_{y=0}.
 \label{5.2.20}
\end{eqnarray}
Since the first term and last four terms on the right-hand side  of $(\ref{5.2.20})$ are in the desired form,  we only need to calculate the second term on the right-hand side of $(\ref{5.2.20})$,
\begin{eqnarray}
&&(\partial _{t}-\epsilon^2\partial_x^2)(2\partial_y w\partial_x\partial_yw+4w\partial_x\partial_y^2w-\partial_xw\partial_y^2w)\big{|}_{y=0}\nonumber\\
&&= (-\partial _{x}\partial _{t}w\partial _{y}^{2}w-\partial _{x}w \partial _{y}^{2}\partial _{t}w+2\partial _{x}\partial _{y}\partial _{t}w\partial _{y}w+2\partial _{x}\partial _{y}w\partial _{y}\partial _{t}w
+4w\partial _{x}\partial _{y}^{2}\partial _{t}w +4\partial _{t}w \partial _{x}\partial _{y}^{2}w )\big{|}_{y=0}\nonumber\\
&&\quad-\epsilon^2\partial_x^2(2\partial_y w\partial_x\partial_yw+4w\partial_x\partial_y^2w-\partial_xw\partial_y^2w)\big{|}_{y=0},\label{5.2.21}
\end{eqnarray}
 here
 \begin{eqnarray}
\begin{aligned}
&-(\partial _{x}\partial _{t}w\partial _{y}^{2}w)\big{|}_{y=0}=-\partial _{y}^{2}w(\partial _{x}\partial _{y}^{2}w-\partial _{x}w+\epsilon^2\partial_x^3 w)\big{|}_{y=0},
 \label{5.2.22}
\end{aligned}
\end{eqnarray}
 \begin{eqnarray}
-\partial _{y}^{2}\partial _{t} w\big{|}_{y=0}
&=& -\partial _{y}^{2}(\partial _{y}^{2}w-u\partial _{x}w-v\partial _{y}w-w+\epsilon^2\partial_x^2w)\big{|}_{y=0}\nonumber\\
&=&-(\partial _{y}^{4}w-\partial _{y}^{2}w+2w\partial _{x}\partial _{y}w+\epsilon^2\partial _{y}^{2}\partial_x^2w)\big{|}_{y=0},
\end{eqnarray}
which implies
 \begin{eqnarray}
\begin{aligned}
\partial _{x}w\partial _{y}^{2}\partial _{t} w\big{|}_{y=0}&=\sum\limits_{j=0}^{4} \wedge_{\alpha}\partial _{x}  \partial _{y}^{j} w  \partial _{y}^{4-j} w
+\sum\limits_{j=0}^{2 } \wedge_{\alpha}\partial _{x}  \partial _{y}^{j} w  \partial _{y}^{2-j} w+\sum\limits_{\rho_{1}+\rho_{2}+\rho_{3}=3} \wedge_{\alpha} D^{\rho_{1}}w D^{\rho_{2}}w D^{\rho_{3}}w \big{|}_{y=0}\\
&\quad -\epsilon^2\partial _{x}w\partial _{y}^{2}\partial_x^2w\big{|}_{y=0}.
\end{aligned}
\end{eqnarray}
Since $ v=\partial _{x}\partial _{y}^{-1}w$,
\begin{eqnarray}
\partial _{x}\partial _{y}\partial _{t} w\big{|}_{y=0}
&=&\partial _{x}\partial _{y}(\partial _{y}^{2}w-u\partial _{x}w-v\partial _{y}w-w+\epsilon^2\partial_x^2w)\big{|}_{y=0}\nonumber\\
&=&\{\partial _{x}\partial _{y}^{3}w-\partial _{x}\partial _{y}w+\epsilon^2\partial_y\partial_x^3w-\partial _{x}(w\partial _{x}w+u\partial _{x}\partial _{y}w+\partial _{y}v\partial _{y}w+v\partial _{y}^2w)\}\big{|}_{y=0}\nonumber\\
&=&(\partial _{x}\partial _{y}^{3}w-\partial _{x}\partial _{y}w+\epsilon^2\partial_y\partial_x^3w-\partial _{x}w\partial _{x}w-w\partial _{x}^2w-w\partial _{y}\partial _{x}w)\big{|}_{y=0}
\end{eqnarray}
which implies
 \begin{eqnarray}
\begin{aligned}
2\partial _{x}\partial _{y}\partial _{t}w\partial _{y} w|_{y=0}&\in \left(\sum\limits_{j=0}^{4} \wedge_{\alpha}\partial _{x}  \partial _{y}^{j} w  \partial _{y}^{4-j} w
+\sum\limits_{j=0}^{2 } \wedge_{\alpha}\partial _{x}  \partial _{y}^{j} w  \partial _{y}^{2-j} w+\sum\limits_{\rho_{1}+\rho_{2}+\rho_{3}=3} \wedge_{\alpha} D^{\rho_{1}}w D^{\rho_{2}}w D^{\rho_{3}}w \right)|_{y=0}\\
&\quad+2\epsilon^2\partial_y\partial_x^3w\partial_yw|_{y=0}.
\end{aligned}
\end{eqnarray}
Similarly,
 \begin{eqnarray}
\begin{aligned}
&\left(2\partial _{x}\partial _{y}w\partial _{y}\partial _{t}w
+4w\partial _{x}\partial _{y}^{2}\partial _{t}w +4\partial _{t}w \partial _{x}\partial _{y}^{2}w\right) |_{y=0}\\
&\in \left(\sum\limits_{j=0}^{4} \wedge_{\alpha}\partial _{x}  \partial _{y}^{j} w  \partial _{y}^{4-j} w
+\sum\limits_{j=0}^{2 } \wedge_{\alpha}\partial _{x}  \partial _{y}^{j} w  \partial _{y}^{2-j} w+\sum\limits_{\rho_{1}+\rho_{2}+\rho_{3}=3} \wedge_{\alpha} D^{\rho_{1}}w D^{\rho_{2}}w D^{\rho_{3}}w \right)|_{y=0}\\
&\quad+2\epsilon^2\partial _{x}\partial _{y}w\partial _{y}\partial _{x}^2w|_{y=0}+4\epsilon^2w\partial _{x}^3\partial _{y}^{2}w|_{y=0}+4\epsilon^2\partial _{x}^2w \partial _{x}\partial _{y}^{2}w|_{y=0}. \label{5.2.23}
\end{aligned}
\end{eqnarray}

Thus it follows from $(\ref{5.2.20})$, $(\ref{5.2.21})$, $(\ref{5.2.22})$ -$(\ref{5.2.23})$ that
\begin{eqnarray}
\begin{aligned}
\partial _{y}^{7}w|_{y=0}
&=\sum\limits_{i=0}^{3}(\partial _{t}-\epsilon^2\partial_x^2)^iK+(1+\epsilon^2\partial^2_x+\epsilon^4\partial^4_x)(w\partial _{x}w)|_{y=0}+\sum\limits_{j=0}^{0 } \partial _{x}  \partial _{y}^{0} w  \partial _{y}^{0-j} w|_{y=0}\\
&\quad+\sum\limits_{j=0}^{2 } \wedge_{\alpha}\partial _{x}  \partial _{y}^{j} w  \partial _{y}^{2-j} w|_{y=0} +\sum\limits_{j=0}^{4 } \wedge_{\alpha}\partial _{x}  \partial _{y}^{j} w  \partial _{y}^{4-j} w|_{y=0} \\
&\quad+\sum\limits_{\rho_{1}+\rho_{2}+\rho_{3}=3} \wedge_{\alpha} D^{\rho_{1}}w D^{\rho_{2}}w D^{\rho_{3}}w |_{y=0}+\sum\limits_{j_x=0}^{2 } \sum\limits_{j_y=0}^{2 }\wedge_{\alpha}\partial _{x}  \partial _{y}^{j_x}\partial _{y}^{j_y} w \partial _{y}^{2-j_x} \partial _{y}^{2-j_y} w|_{y=0}.
 \label{5.2.24}
\end{aligned}
\end{eqnarray}
We justify the formula (\ref{5.2.24}) for $k=3$.

Now, using the same algorithm, we are going to prove formula (\ref{5.2.005}) by induction on $k$. For notational convenience, we denote
\begin{eqnarray}
\mathcal{A}_k:=\bigg\{\sum\limits_{l=0}^{k-1}\epsilon^{2l}\sum\limits_{j=2}^{\max\{2,k-l\}}\sum\limits_{\rho \in A^j_{k,l}}C_{\wedge_{\alpha},k,l,\rho^1,\rho^2,...,\rho^j}\prod^{j}_{i=1}D^{\rho^i}w|_{y=0}\bigg\}.
\end{eqnarray}
Using this notation, we will prove $\partial_y^{2k+1}w|_{y=0}-(\partial_t-\epsilon^2\partial_x^2)^k(\partial_x  P-U) \in \mathcal{A}_k$. Assuming that formula (\ref{5.2.005}) holds for $k=n$, we will show that it also holds for $k=n+1$ as follows. Then we differentiate the vorticity equation $2n+1$ times with respect to $y$  to obtain
\begin{eqnarray}
\partial _{y}^{2n+3}w|_{y=0}&=& \left\{(\partial _{t}-\epsilon^2\partial_x^2)\partial _{y}^{2n+1}w+\partial _{y}^{2n+1}(u\partial _{x}w+v\partial _{y}w)+\partial _{y}^{2n+1}w\right\}|_{y=0}\nonumber\\
&=&\left\{(\partial _{t}-\epsilon^2\partial_x^2)\partial _{y}^{2n+1}w+\sum\limits_{j=0}^{2n } \wedge_{\alpha}\partial _{x}  \partial _{y}^{j} w  \partial _{y}^{2n-j} w|_{y=0}+\partial _{y}^{2n+1}w\right\}|_{y=0}.
\label{5.2.006}
\end{eqnarray}
By routine checking, one may show that the last three terms of (\ref{5.2.006}) belong to
$\mathcal{A}_{k+1}$, so it only remains to deal with the term $(\partial _{t}-\epsilon^2\partial_x^2)\partial _{y}^{2n+1}w$.
Thanks to the induction hypothesis, there exist constants $C_{\wedge_{\alpha},n,l,\rho^1,\rho^2,...,\rho^j}$ such that
\begin{eqnarray}
\begin{aligned}
\partial_y^{2n+1}w|_{y=0}&=C_n\sum\limits_{s=0}^{n}(\partial_t-\epsilon^2\partial^2_x)^n(\partial_x P -U)\\
&\quad+\sum\limits_{l=0}^{k-1}\epsilon^{2l}\sum\limits_{j=2}^{\max\{2,n-l\}}\sum\limits_{\rho \in A^j_{n,l}}C_{\wedge_{\alpha},n,l,\rho^1,\rho^2,...,\rho^j}\prod^{j}_{i=1}D^{\rho^i}w|_{y=0},
\end{aligned}
\end{eqnarray}
so we have, up to a relabeling of the indices $\rho^i$,
\begin{eqnarray}
\begin{aligned}
(\partial_t-\epsilon^2\partial^2_x)\partial_y^{2n+1}w|_{y=0}&=C_n\sum\limits_{s=1}^{n+1}(\partial_t-\epsilon^2\partial^2_x)^{s}(\partial_x P -U)\\
&\quad+\sum\limits_{l=0}^{n-1}\epsilon^{2l}\sum\limits_{j=2}^{\max\{2,n-l\}}\sum\limits_{\rho \in A^j_{n,l}} \widetilde{C}_{\wedge_{\alpha},n,l,\rho^1,\rho^2,...,\rho^j}(\partial_t-\epsilon^2\partial^2_x)D^{\rho^1}w\prod^{j}_{i=2}D^{\rho^i}w|_{y=0}\\
&\quad-\sum\limits_{l=0}^{n-1}\epsilon^{2l+2}\sum\limits_{j=2}^{\max\{2,n-l\}}\sum\limits_{\rho \in A^j_{n,l}}\widetilde{\widetilde{C}}_{\wedge_{\alpha},n,l,\rho^1,\rho^2,...,\rho^j}\partial_x D^{\rho^1}w\partial_x D^{\rho^2}w\prod^{j}_{i=3}D^{\rho^i}w|_{y=0},
\label{5.2.007}
\end{aligned}
\end{eqnarray}
where $\widetilde{C}_{\wedge_{\alpha},n,l,\rho^1,\rho^2,...,\rho^j}$ and $\widetilde{\widetilde{C}}_{\wedge_{\alpha},n,l,\rho^1,\rho^2,...,\rho^j}$ are some new constants depending on $C_{\wedge_{\alpha},n,l,\rho^1,\rho^2,...,\rho^j}$. It is worth noting that the last term on the right-hand side of (\ref{5.2.007})
belongs to $\mathcal{A}_{n+1}$, so it remains to check whether the second term on right-hand
side of (\ref{5.2.007}) also belongs to $\mathcal{A}_{n+1}$.
This completes the proof.
\hfill $\Box$

\begin{Lemma}\label{y5.2.2}
Let $s \geq 4$ be an  integer, $\gamma\geq 1$, $\sigma>\gamma+\frac{1}{2}$ and $\delta\in (0,1)$. If  $w\in H^{s,\gamma}_{\sigma,\delta}$  solves $(\ref{5.2.0001})$  and $(\ref{5.2.0002})$, we have the following  estimates:\\
(i) When $|\alpha|\leq s-1$,
\begin{eqnarray}
\bigg|\int_{\mathbb{T} } D^{\alpha}w\partial _{y}D^{\alpha}wdx | _{y=0}\bigg|\leq \frac{1}{12}\|(1+y)^{\gamma+\alpha_2+1}\partial^2_y D^\alpha w\|_{L^2}^2+C\|w\|^2_{H^{s,\gamma}_g},
\end{eqnarray}
(ii) When $|\alpha|=s$, $\alpha_2$ is even,
\begin{eqnarray}
\begin{aligned}
\bigg|\int_{\mathbb{T} } D^{\alpha}w\partial _{y}D^{\alpha}wdx | _{y=0}\bigg|&\leq \frac{1}{12}\|(1+y)^{\gamma+\alpha_2}\partial_y D^\alpha w\|^2_{L^2}+C_{s,\gamma,\sigma,\delta}(1+\|w\|^2_{H^{s,\gamma}_g})^{s-2}\|w\|^2_{H^{s,\gamma}_g}\\
&\quad+C_s \sum\limits_{l=0}^{\frac{s}{2}}\|\partial^{l}_t(\partial_xP-U)\|_{H^{s-2l}(\mathbb{T})}^{2},
\end{aligned}
\end{eqnarray}
(iii) When $|\alpha|=s$, $\alpha_2$ is odd,
\begin{eqnarray}
\begin{aligned}
\bigg|\int_{\mathbb{T} } D^{\alpha}w\partial _{y}D^{\alpha}wdx | _{y=0}\bigg|&\leq \frac{1}{12}\|(1+y)^{\gamma+\alpha_2+1}\partial_x^{\alpha_1-1}\partial_y^{\alpha_2+2} w\|^2_{L^2}+C_{s,\gamma,\sigma,\delta}(1+\|w\|^2_{H^{s,\gamma}_g})^{s-2}\|w\|^2_{H^{s,\gamma}_g}\\
&\quad+C_s \sum\limits_{l=0}^{\frac{s}{2}}\|\partial^{l}_t(\partial_xP-U)\|_{H^{s-2l}(\mathbb{T})}^{2}.
\end{aligned}
\end{eqnarray}
\end{Lemma}
\textbf{Proof}.
\emph{Case 1.} When $|\alpha|\leq s-1$, by using the following trace estimate
\begin{eqnarray}
\int_{\mathbb{T} } |f|dx \bigg |_{y=0} \leq C\left(\int^1_0\int_{\mathbb{T}}|f|dxdy+\int^1_0\int_{\mathbb{T}}|\partial_y f|dxdy\right),
\end{eqnarray}
we have
\begin{eqnarray}
\bigg|\int_{\mathbb{T} } D^{\alpha}w\partial _{y}D^{\alpha}wdx | _{y=0}\bigg|\leq \frac{1}{12}\|(1+y)^{\gamma+\alpha_2+1}\partial^2_y D^\alpha w\|_{L^2}^2+C\|w\|^2_{H^{s,\gamma}_g}.
\end{eqnarray}
\emph{Case 2.} When $|\alpha|=s$, $\alpha_2=2k$ for some $k\in \mathbb{N}$, we can apply boundary reduction Lemma \ref{y5.2.1} to $\partial_y D^\alpha  | _{y=0}$ and obtain
\begin{eqnarray}
\begin{aligned}
\int_{\mathbb{T} } D^{\alpha}w\partial _{y}D^{\alpha}wdx | _{y=0}&=C_k\sum\limits_{s=0}^{k}\int_{\mathbb{T} } D^{\alpha}w(\partial_t-\epsilon^2\partial_x^2)^k\partial_x^{\alpha_1}(\partial_xP-U)dx | _{y=0}\\
&\quad+\sum\limits_{l=0}^{k-1}\epsilon^{2l}\sum\limits_{j=2}^{\max\{2,k-l\}}\sum\limits_{\rho \in A^j_{k,l}}C_{\wedge_{\alpha},k,l,\rho^1,\rho^2,...,\rho^j}\int_{\mathbb{T} }D^\alpha w \partial_x^{\alpha_1}(\prod^{j}_{i=1}D^{\rho^i}w)dx|_{y=0},
\end{aligned}
\end{eqnarray}
then we again apply the simple trace estimate to control the boundary integral as follows:
\begin{eqnarray}
\begin{aligned}
\bigg|\int_{\mathbb{T} } D^{\alpha}w\partial _{y}D^{\alpha}wdx | _{y=0}\bigg|&\leq \frac{1}{12}\|(1+y)^{\gamma+\alpha_2}\partial_y D^\alpha w\|^2_{L^2}+C_{s,\gamma,\sigma,\delta}(1+\|w\|_{H_g^{s,\gamma}})^{s-2}\|w\|^2_{H^{s,\gamma}_g}\\
&\quad+C_s \sum\limits_{l=0}^{\frac{s}{2}}\|\partial^{l}_t(\partial_xP-U)\|_{H^{s-2l}(\mathbb{T})}^{2}.
\end{aligned}
\end{eqnarray}
\emph{Case 3.} When $|\alpha|=s$, $\alpha_2=2k+1$ for some $k\in \mathbb{N}$, using integration by parts in the $x$-variable, we have
\begin{eqnarray}
\int_{\mathbb{T} } D^{\alpha}w\partial _{y}D^{\alpha}wdx | _{y=0}=-\int_{\mathbb{T}}\partial_x D^\alpha w \partial_x^{\alpha_1-1}\partial_y^{\alpha_2+1}wdx |_{y=0}.
\end{eqnarray}
The term $\partial_x D^\alpha w|_{y=0}=\partial_x^{\alpha_1+1}\partial_y^{2k+1}w | _{y=0}$ has an odd number of $y$ derivatives, then we have
\begin{eqnarray}
\begin{aligned}
\bigg|\int_{\mathbb{T} } D^{\alpha}w\partial _{y}D^{\alpha}wdx | _{y=0}\bigg|&\leq \frac{1}{12}\|(1+y)^{\gamma+\alpha_2+1}\partial_x^{\alpha_1-1}\partial_y^{\alpha_2+2} w\|^2_{L^2}+C_{s,\gamma,\sigma,\delta}(1+\|w\|_{H_g^{s,\gamma}})^{s-2}\|w\|^2_{H^{s,\gamma}_g}\\
&\quad+C_s \sum\limits_{l=0}^{\frac{s}{2}}\|\partial^{l}_t(\partial_xP-U)\|_{H^{s-2l}(\mathbb{T})}^{2}.
\end{aligned}
\end{eqnarray}
\hfill $\Box$

\begin{Proposition}\label{p5.2.1}
Let $s \geq 4$ be an   integer, $\gamma\geq 1$, $\sigma>\gamma+\frac{1}{2}$ and $\delta\in (0,1)$. If  $w\in H^{s,\gamma}_{\sigma,\delta}$  solves $(\ref{5.2.0001})$  and $(\ref{5.2.0002})$, we have the following uniform (in $\epsilon$) estimate:
\begin{eqnarray}
\begin{aligned}
&\frac{1}{2}\frac{d}{dt}\sum\limits_{\substack{ |\alpha|\leq s\\ \alpha_{1}\leq s-1}}\|(1+y)^{ \gamma+ \alpha_{2}}D^{\alpha}w\|_{L^{2}}^{2}+\sum\limits_{\substack{ |\alpha|\leq s\\ \alpha_{1}\leq s-1}}\|(1+y)^{ \gamma+ \alpha_{2}}D^{\alpha}w\|_{L^{2}}^{2}
 \\
&\leq -\frac{1}{2}\sum\limits_{\substack{ |\alpha|\leq s\\ \alpha_{1}\leq s-1}}\|(1+y)^{ \gamma+ \alpha_{2}}\partial _{y}D^{\alpha}w\|_{L^{2}} ^{2}-\epsilon^2\sum\limits_{\substack{ |\alpha|\leq s\\ \alpha_{1}\leq s-1}} \|(1+y)^{\gamma+ \alpha_{2}}\partial _{x}D^{\alpha}w\|_{L^{2}}^2\\
&\quad \quad+C_{s, \gamma,\sigma,\delta  }\|w\|_{H^{s,\gamma}_{g}}^2(\|\partial_{x}^{s}U\|_{{L^\infty}(\mathbb{T})}^2+\|w\|_{H^{s,\gamma}_{g}}^2)+C_{s,\gamma,\sigma,\delta}(1+\|w\|_{H_g^{s,\gamma}})^{s-2}\|w\|^2_{H^{s,\gamma}_g}\\
&\quad \quad+C_s \sum\limits_{l=0}^{\frac{s}{2}}\|\partial^{l}_t(\partial_xP-U)\|_{H^{s-2l}(\mathbb{T})}^{2}.
\label{5.2.002}
\end{aligned}
\end{eqnarray}
\end{Proposition}
\textbf{Proof}.
Applying the operator $D^{\alpha}$ on $(\ref{5.2.0001})_{1}$ with $\alpha=(\alpha_{1},\alpha_{2}),~|\alpha|\leq s, ~\alpha_{1}\leq s-1$,
\begin{eqnarray}
(\partial _{t} +u\partial _{x} +v\partial _{y}+1-\epsilon^2\partial _{x}^{2}-\partial _{y}^{2})D^{\alpha}w
=- \sum\limits_{0<\beta\leq \alpha}\binom{ \alpha}{\beta} \left\{D^{\beta} u\partial _{x}D^{\alpha-\beta}w + D^{\beta}v\partial _{y}D^{\alpha-\beta}w\right\} .
 \label{5.2.1}
\end{eqnarray}
Multiplying $(\ref{5.2.1})$ by $(1+y)^{2\gamma+2\alpha_{2}}D^{\alpha}w$ and integrating it over $\mathbb{T}\times \mathbb{R}_{+}$ yield
\begin{eqnarray}
\begin{aligned}
&\frac{1}{2}\frac{d}{dt}\|(1+y)^{ \gamma+ \alpha_{2}}D^{\alpha}w\|_{L^{2}}^{2}+\|(1+y)^{ \gamma+ \alpha_{2}}D^{\alpha}w\|_{L^{2}}^{2}
+\|(1+y)^{ \gamma+ \alpha_{2}}\partial _{y}D^{\alpha}w\|_{L^{2}}^{2}+ \epsilon^2\|(1+y)^{\gamma+ \alpha_{2}}\partial _{x}D^{\alpha}w\|_{L^{2}}^2   \\
&=-\int_{\mathbb{T} } D^{\alpha}w\partial _{y}D^{\alpha}wdx | _{y=0}-(2\gamma+2\alpha_{2}) \iint(1+y)^{2\gamma+2\alpha_{2}-1} D^{\alpha}w\partial _{y}   D^{\alpha}w\\
&\quad+(\gamma+\alpha_{2}) \iint(1+y)^{2\gamma+2\alpha_{2}-1} v  | D^{\alpha}w|^{2}  \\
&\quad - \sum\limits_{0<\beta\leq \alpha}\binom{ \alpha}{\beta} \iint(1+y)^{2\gamma+2\alpha_{2}}D^{\alpha}w  \left\{D^{\beta} u\partial _{x}D^{\alpha-\beta}w + D^{\beta}v\partial _{y}D^{\alpha-\beta}w\right\} . \label{5.2.2}
\end{aligned}
\end{eqnarray}
Indeed, $(\ref{5.2.2})$ can be obtained by integration by parts and the boundary condition.
We need to estimate the integral equation $(\ref{5.2.2})$ term by term. Obviously, the first term on the right-hand side of $(\ref{5.2.2})$ follows from Lemmas \ref{y5.2.1}-\ref{y5.2.2}.

Secondly,  using H$\ddot{o}$lder's inequality, we have
\begin{eqnarray}
&&\left|(2\gamma+2\alpha_{2}) \iint(1+y)^{2\gamma+2\alpha_{2}-1} D^{\alpha}w\partial _{y}   D^{\alpha}w  \right|\nonumber \\
 &&\leq (2\gamma+2\alpha_{2}) \|\frac{1}{1+y }  \|_{L^{\infty}} \|(1+y)^{\gamma+ \alpha_{2}}D^{\alpha}w \|_{L^{2}}
  \|(1+y)^{\gamma+ \alpha_{2}}\partial_{y} D^{\alpha}w \|_{L^{2}} \nonumber \\
&&\leq  C_{ s,\gamma } \|w\|_{H^{s,\gamma}_{g}}^2+\frac{1}{4}\|(1+y)^{\gamma+ \alpha_{2}}\partial_{y} D^{\alpha}w \|_{L^{2}}^{2} ,  \label{5.2.5}
\end{eqnarray}
and using Lemma \ref{y5.4.4},
\begin{eqnarray}
\left|   (\gamma+\alpha_{2}) \iint(1+y)^{2\gamma+2\alpha_{2}-1} v  | D^{\alpha}w|^{2}  \right|
 &\leq&  \|\frac{v}{1+y}\|_{L^{\infty}}  \|(1+y)^{\gamma+ \alpha_{2}}D^{\alpha}w \|_{L^{2}}^{2} \nonumber \\
& \leq &  C_{s, \gamma,\sigma,\delta  }(\|w\|_{H^{s,\gamma}_{g}}+\|\partial_{x}^{s}U\|_{{L^2}(\mathbb{T})})\|w\|_{H^{s,\gamma}_{g}}^2.  \label{5.2.6}
\end{eqnarray}
Lastly, noticing the fact $\partial _{y}v=-\partial _{x}u, ~\partial _{y}u=w$, it follows that the last term on the right-hand side of equation $(\ref{5.2.2})$ has the following three cases, \\ (i) the first case:
$$ J_{1}= \iint(1+y)^{2\gamma+2\alpha_{2}}D^{\alpha}w  \partial _{x}^{\eta}vD^{k}w ,\ \   \eta\in [1,s-1],    $$
(a) for $\eta=s-1$, we get
\begin{eqnarray}
|J_{1}|& =& \left| \iint(1+y)^{2\gamma+2\alpha_{2}}D^{\alpha}w  \partial _{x}^{s-1}v\partial _{y}^{2}w  \right|\nonumber \\
 &\leq &  \|(1+y)^{\gamma+ \alpha_{2}}D^{\alpha}w \|_{L^{2}}   \|\frac{\partial _{x}^{s-1}v +y\partial_x^sU} {1+y}\|_{L^{2}}
  \|(1+y)^{\gamma+k_2}D ^{k }w\|_{L^{\infty}}   \nonumber \\
  &\quad& + \|(1+y)^{\gamma+ \alpha_{2}}D^{\alpha}w \|_{L^{2}}   \|\partial_x^s U\|_{{L^\infty}(\mathbb{T})}
  \|(1+y)^{\gamma+k_2} D ^{k }w\|_{L^{\infty}}   \nonumber \\
&\leq & C_{s, \gamma,\sigma,\delta  }\|w\|_{H^{s,\gamma}_{g}} ^2(\|w\|_{H^{s,\gamma}_{g}}+\|\partial_{x}^{s}U\|_{{L^\infty}(\mathbb{T})}),
\label{5.2.7}
\end{eqnarray}
(b) for $\eta=1,2,\cdots,s-2$, we derive
\begin{eqnarray}
|J_{1}|
&\leq &   \|(1+y)^{\gamma+ \alpha_{2}}D^{\alpha}w \|_{L^{2}}   \|\frac{\partial _{x}^{\eta}v } {1+y}\|_{L^{\infty}}
  \|(1+y)^{\gamma+k_{2}}D^{k}w\|_{L^{2}}   \nonumber \\
&\leq & C_{s, \gamma,\sigma,\delta  }\|w\|_{H^{s,\gamma}_{g}}^2 (\|w\|_{H^{s,\gamma}_{g}}+\|\partial_{x}^{s}U\|_{{L^2}(\mathbb{T})}).
\label{5.2.8}
\end{eqnarray}
(ii) the second case:
$$ J_{2}= \iint(1+y)^{2\gamma+2\alpha_{2}}D^{\alpha}w  \partial _{x}^{\eta}uD^{k}w ,\ \ \eta\in [1,s ],    $$
(a) for $\eta=s $, we get
\begin{eqnarray}
|J_{2}|& =& \left| \iint(1+y)^{2\gamma+2\alpha_{2}}D^{\alpha}w  \partial _{x}^{s }u\partial _{y} w  \right|
 \leq   \|(1+y)^{\gamma+ \alpha_{2}}D^{\alpha}w \|_{L^{2}}   \|\partial _{x}^{s }(u-U)\|_{L^{2}}
  \|(1+y)^{\gamma+k_2} D^k w\|_{L^{\infty}}   \nonumber \\
  &\quad&+\|(1+y)^{\gamma+ \alpha_{2}}D^{\alpha}w \|_{L^{2}}   \|\partial _{x}^{s }U \|_{{L^\infty}(\mathbb{T})}
  \|(1+y)^{\gamma+k_2} D^k w\|_{L^{2}}   \nonumber \\
&\leq & C_{s, \gamma,\sigma,\delta  }\|w\|_{H^{s,\gamma}_{g}}^2 (\|w\|_{H^{s,\gamma}_{g}}+\|\partial_{x}^{s}U\|_{{L^\infty}(\mathbb{T})}),
\label{5.2.9}
\end{eqnarray}
(b) for $\eta=1,2,\cdots,s-1$, we derive
\begin{eqnarray}
|J_{2}|
&\leq &   \|(1+y)^{\gamma+ \alpha_{2}}D^{\alpha}w \|_{L^{2}}   \| \partial _{x}^{\eta}u   \|_{L^{\infty}}
  \|(1+y)^{\gamma+k_{2}}D^{k}w\|_{L^{2}}   \nonumber \\
&\leq & C_{s, \gamma,\sigma,\delta  }\|w\|_{H^{s,\gamma}_{g}}^2 (\|w\|_{H^{s,\gamma}_{g}}+\|\partial_{x}^{s}U\|_{{L^2}(\mathbb{T})}).
\label{5.2.10}
\end{eqnarray}
(iii) the last case:
$$ J_{3}= \iint(1+y)^{2\gamma+2\alpha_{2}}D^{\alpha}w  D^{\theta}w D^{k}w ,\ \ \theta\in [0,s-1 ],  \theta+k=s,  $$
we obtain
\begin{eqnarray}
|J_{3}|
&\leq & C_{s, \gamma,\sigma,\delta  }\|w\|_{H^{s,\gamma}_{g}} (\|w\|_{H^{s,\gamma}_{g}}+\|\partial_{x}^{s}U\|_{{L^2}(\mathbb{T})})\|(1+y)^{\gamma+ \alpha_{2}}D^{\alpha}w \|_{L^{2}}\nonumber \\
 &\leq& C_{s, \gamma,\sigma,\delta  }\|w\|_{H^{s,\gamma}_{g}}^{2}(\|w\|_{H^{s,\gamma}_{g}}+\|\partial_{x}^{s}U\|_{{L^2}(\mathbb{T})})^2+\frac{1}{12} \|(1+y)^{\gamma+ \alpha_{2}}D^{\alpha}w \|_{L^{2}}  ^{2}.
\label{5.2.12}
\end{eqnarray}
Hence, combining $(\ref{5.2.5})$-$(\ref{5.2.12})$, $(\ref{5.2.2})$ reduces to inequality (\ref{5.2.002}).
\hfill $\Box$

\subsection{Estimates on $g_{s}$ }
In this subsection, we will estimate the norm of $D^{\alpha}w$ with weight  $(1+y)^{\gamma+\alpha_{2}}$, here $\alpha=(s, 0)$, i.e., $(1+y)^{\gamma}g_{s}$, because  we have estimated the norm of $(1+y)^{\gamma+\alpha_{2}} D^{\alpha}w$,    $\alpha=(\alpha_{1}, \alpha_{2})$, $\alpha_{1}\leq s-1$ in subsection 2.1.
The evolution equations for $w$ and $u-U$ as follows
\begin{equation}\left\{
\begin{array}{ll}
(\partial _{t} +u\partial _{x} +v\partial _{y}+1-\partial _{y}^{2}-\epsilon^2\partial _{x}^{2}) w
=0,\\
(\partial _{t} +u\partial _{x} +v\partial _{y}+1-\partial _{y}^{2}-\epsilon^2\partial _{x}^{2}) (u-U)
=-  (u-U)\partial _{x} U  .
\end{array}
 \label{5.3.01}         \right.\end{equation}

\begin{Proposition}\label{p5.2.2}
Under the same assumption of Proposition \ref{p5.2.1}, we have the following estimate:
\begin{eqnarray}
\begin{aligned}
&\frac{d}{dt}\|(1+y)^{ \gamma}g_{s}\|_{L^{2}}^{2}+\|(1+y)^{ \gamma}g_{s}\|_{L^{2}}^{2} \\
&\leq -\frac{1}{2}\epsilon^2\|(1+y)^{\gamma}\partial _{x}g_{s}\|_{L^{2}} ^{2} -\frac{1}{2}\|(1+y)^{\gamma}\partial _{y}g_{s}\|_{L^{2}} ^{2}+\|w\|_{H^{s,\gamma}_{g}}^2 \\
&\quad+ C\|\partial_{x}^{s}(\partial_xP-U)\|_{{L^2}(\mathbb{T})}^2+C_{ \gamma,\delta  }\|\partial_{x}^{s}U\|_{{L^{2}}(\mathbb{T})}^2 \|w\|_{H^{s,\gamma}_{g}}^2\\
&\quad+ C_{s, \gamma,\sigma,\delta  }\left(1+\|w\|_{H^{s,\gamma}_{g}}+\|\partial_{x}^{s}U\|_{{L^\infty}(\mathbb{T})}\right)
\left(\|w\|_{H^{s,\gamma}_{g}}+\|\partial_{x}^{s+1}U\|_{{L^\infty}(\mathbb{T})}\right)\|w\|_{H^{s,\gamma}_{g}}.
\label{5.3.00003}
\end{aligned}
\end{eqnarray}
\end{Proposition}
\textbf{Proof}.
Applying the operator $\partial_{x}^{s}$ on $(\ref{5.3.01})_{1}$ and $(\ref{5.3.01})_{2}$ respectively, we obtain the following  equations
\begin{equation}\left\{
\begin{array}{ll}
(\partial _{t} +u\partial _{x} +v\partial _{y}+1-\partial _{y}^{2}-\epsilon^2\partial _{x}^{2})\partial_{x}^{s}w+\partial_{x}^{s}v\partial _{y}w
=- \sum\limits_{0\leq j<s}\binom{s}{j} \partial _{x}^{s-j}u\partial _{x}^{j+1}w - \sum\limits_{1\leq j<s}\binom{s}{j} \partial _{x}^{s-j}v\partial _{y}\partial _{x}^{j }w  ,\\
(\partial _{t} +u\partial _{x} +v\partial _{y}+1-\partial _{y}^{2}-\epsilon^2\partial _{x}^{2})\partial_{x}^{s}(u-U)+\partial_{x}^{s}v\partial _{y}u
=- \sum\limits_{0\leq j<s}\binom{s}{j} \partial _{x}^{s-j}u\partial _{x}^{j+1}(u-U) \\
- \sum\limits_{1\leq j<s}\binom{s}{j} \partial _{x}^{s-j}v\partial _{y}\partial _{x}^{j }u- \sum\limits_{0\leq j\leq s}\binom{s}{j} \partial _{x}^{j}(u-U)\partial _{x}^{s-j+1}U  .
\end{array}
 \label{5.3.1}         \right.\end{equation}
To overcome the difficult term  $\partial_{x}^{s}v$, we use the cancellation property as in \cite{mw}. Indeed, eliminating $\frac{\partial _{y}w}{w}\times (\ref{5.3.1})_{2}$ from $ (\ref{5.3.1})_{1}$, and letting  $a(t,x,y)=\frac{\partial _{y}w}{w}$, it follows
\begin{eqnarray}
\begin{aligned}
&
 (\partial _{t} +u\partial _{x} +v\partial _{y}+1-\partial _{y}^{2})g_{s} +\partial_{x}^{s}(u-U) (\partial _{t} +u\partial _{x} +v\partial _{y}+1-\partial _{y}^{2})a    \\
&= 2\epsilon^2\partial_x^{s+1}(u-U)\partial_x a+2\partial_{x}^{s}w\partial _{y}a - \sum\limits_{j=0}^{s-1}\binom{s}{j} g_{j+1}\partial _{x}^{s-j}u
- \sum\limits_{j=1}^{s-1}\binom{s}{j}\partial _{x}^{s-j}v \{\partial _{y}\partial _{x}^{j }w-a \partial _{y}\partial _{x}^{j }u \}\\
&\quad+a\sum\limits_{j=0}^s\binom{s}{j} \partial _{x}^{j}(u-U)\partial _{x}^{s-j+1}U . \label{5.3.2}
\end{aligned}
\end{eqnarray}
Applying the operator $\partial_{y}$ on $(\ref{5.2.0002})_{1}$ yields
\begin{eqnarray}
(\partial _{t} +u\partial _{x} +v\partial _{y}+1 )\partial _{y}w=2\epsilon^2\partial_x\partial_y w+\partial _{y}^{3}w-w\partial _{x}w+\partial _{x}u\partial _{y}w.  \label{5.3.3}
\end{eqnarray}
Using $(\ref{5.3.3})$ and equation $(\ref{5.2.0002})_{1}$, we get
\begin{eqnarray}
 (\partial _{t} +u\partial _{x} +v\partial _{y}+1 )a&=&\frac{(\partial _{t} +u\partial _{x} +v\partial _{y}+1 )\partial _{y}w}{w} -
\frac{\partial _{y}w(\partial _{t} +u\partial _{x} +v\partial _{y}+1 ) w}{w^{2}}   \nonumber\\
&=&\epsilon^2\frac{\partial_x^2\partial_y w}{w}-\epsilon^2a\frac{\partial_x^2w}{w}+\frac{\partial _{y}^{3}w}{w}-a\frac{\partial _{y}^{2}w}{w}-  g_1+a \partial _{x}U,
\label{5.3.4}
\end{eqnarray}
and
\begin{equation}\left\{
\begin{array}{ll}
  \partial _{y}a=\frac{ \partial _{y}^{2}w}{w} -\frac{\partial _{y}w \partial _{y} w}{w^{2}},   \\
 \partial _{y}^{2}a=\frac{\partial _{y}^{3}w}{w}-a\frac{\partial _{y}^{2}w}{w}- 2a \partial _{y}a .
\end{array}
 \label{5.3.5}         \right.\end{equation}
Now inserting $(\ref{5.3.4})$ and $(\ref{5.3.5})_{2}$ into $(\ref{5.3.2})$ yields
\begin{eqnarray}
\begin{aligned}
&
 (\partial _{t} +u\partial _{x} +v\partial _{y}+1-\partial _{y}^{2}-\epsilon^2 \partial_x^2)g_{s}   \\
&= 2\epsilon^2\big\{\partial_x^{x+1}(u-U)-\frac{\partial_xw}{w}\partial_x^2(u-U)\big\}\partial_x a+2g_{s}\partial _{y}a -g_1\partial _{x}^{s}U- \sum\limits_{j=1}^{s-1}\binom{s}{j} g_{j+1}\partial _{x}^{s-j}u
\\
&\quad- \sum\limits_{j=1}^{s-1}\binom{s}{j}\partial _{x}^{s-j}v \{\partial _{y}\partial _{x}^{j }w-a \partial _{y}\partial _{x}^{j }u \}+a\sum\limits_{j=0}^{s-1}\binom{s}{j} \partial _{x}^{j}(u-U)\partial _{x}^{s-j+1}U . \label{5.3.6}
\end{aligned}
\end{eqnarray}
 Multiplying $(\ref{5.3.6})$ by $(1+y)^{2\gamma}g_{s}$, and integrating it over $\mathbb{T}\times \mathbb{R}_{+}$,
\begin{eqnarray}
\begin{aligned}
&\frac{1}{2}\frac{d}{dt}\|(1+y)^{ \gamma}g_{s}\|_{L^{2}}^{2}+\|(1+y)^{ \gamma}g_{s}\|_{L^{2}}^{2}+\|(1+y)^{\gamma}\partial _{y}g_{s}\|_{L^{2}} ^{2}+\epsilon^2\|(1+y)^{\gamma}\partial _{x}g_{s}\|_{L^{2}} ^{2}  \\
&=2\epsilon^2\iint (1+y)^{2\gamma}g_{s}\big\{\partial_x^{x+1}(u-U)-\frac{\partial_xw}{w}\partial_x^2(u-U)\big\}\partial_x a\\
&\quad+\int_{\mathbb{T} } g_{s}\partial _{y}g_{s}dx | _{y=0}-2\gamma \iint(1+y)^{2\gamma-1}g_{s}\partial _{y}  g_{s}
+\gamma \iint(1+y)^{2\gamma-1} v  |g_{s}|^{2}  + 2\iint(1+y)^{2\gamma}|g_{s}|^{2}\partial _{y}a\\
&\quad - \sum\limits_{j=1}^{s-1}\binom{s}{j} \iint(1+y)^{2\gamma}g_{s} g_{j+1}\partial _{x}^{s-j}u
- \sum\limits_{j=1}^{s-1}\binom{s}{j}\iint(1+y)^{2\gamma}g_{s}\partial _{x}^{s-j}v \{\partial _{y}\partial _{x}^{j }w-a \partial _{x}^{j }w \}\\
&\quad-\iint(1+y)^{2\gamma}g_{s}g_1\partial _{x}^{s}U+\sum\limits_{j=0}^{s-1}\binom{s}{j}\iint(1+y)^{2\gamma}g_{s}a \partial _{x}^{j}(u-U)\partial _{x}^{s-j+1}U, \label{5.3.7}
\end{aligned}
\end{eqnarray}
which can be obtained by integration by parts and the boundary condition.
We estimate equation $(\ref{5.3.7})$ by terms. Firstly $w \in C([0,T];H^{s+4,\gamma}_{\sigma,\delta})$, it follows from the  definition of $H^{s+4,\gamma}_{\sigma,\delta}$ that $(1+y)^\sigma w \geq \delta$ and $|(1+y)^{\sigma+\alpha}D^\alpha w|\leq \delta^{-1}$ for all $|\alpha| \leq 2$. Thus, we have $\|(1+y)\partial_x a\|_{L^\infty} \leq \delta^{-2}+\delta^{-4}$ and $\|\frac{\partial_x w}{w}\|_{L^\infty} \leq \delta^{-2}$, and hence
\begin{eqnarray}
\begin{aligned}
&2\epsilon^2\iint (1+y)^{2\gamma}g_{s}\big\{\partial_x^{x+1}(u-U)-\frac{\partial_xw}{w}\partial_x^2(u-U)\big\}\partial_x a\\
&\leq2\epsilon^2 C_{\delta}\|(1+y)^{\gamma}g_s\|_{L^2}(\|(1+y)^{\gamma-1}\partial_x^{s+1}(u-U)\|_{L^2}+\|(1+y)^{\gamma-1}\partial_x^s(u-U)\|_{L^2}).\\
\end{aligned}
\end{eqnarray}
In addition,   using Lemma \ref{y5.4.1} and the following estimate
\begin{equation}\left\{
\begin{array}{ll}
 |w|^{-1}\leq \delta^{-1}(1+y)^{\sigma},   \\
 |\partial _{y}  w| \leq \delta^{-1}(1+y)^{-\sigma-1},  \ \
| \partial _{y}^{2}w| \leq \delta^{-1}(1+y)^{-\sigma-2},
\end{array}
 \label{5.2.60}         \right.\end{equation}
we have
\begin{eqnarray}
\begin{aligned}
\|(1+y)^{\gamma-1}\partial_x^{s+1}(u-U)\|_{L^2} &\leq C_{\gamma, \sigma, \delta}\|(1+y)^{\gamma-\sigma-1}\frac{\partial_x^{s+1}(u-U)}{w}\|_{L^2}\\
&\leq C_{\gamma, \sigma, \delta}\left(\|\partial_x^{s+1}U\|_{L^{2}(\mathbb{T})}+\|(1+y)^{\gamma}\partial_x g_s\|_{L^2}+\|(1+y)^{\gamma-1}\partial_x^{s} (u-U)\|_{L^2}\right),
\end{aligned}
\end{eqnarray}
and
\begin{eqnarray}
\begin{aligned}
\|(1+y)^{\gamma-1}\frac{\partial_x w}{w}\partial_x^{s}(u-U)\|_{L^2} \leq C_{\gamma, \sigma, \delta}\|(1+y)^{\gamma-1}\partial_x^{s} (u-U)\|_{L^2}.
\end{aligned}
\end{eqnarray}
Then we can obtain
\begin{eqnarray}
\begin{aligned}
&\bigg|2\epsilon^2\iint (1+y)^{2\gamma}g_{s}\big\{\partial_x^{x+1}(u-U)-\frac{\partial_xw}{w}\partial_x^2(u-U)\big\}\partial_x a\bigg|\\
& \leq\frac{1}{2} \epsilon^2 \|(1+y)^\gamma \partial_x g_s\|_{L^2}^2 + \epsilon^2 C_{s,\gamma, \sigma, \delta}\left(\|w\|_{H^{s,\gamma}_g}+\|\partial_x^{s+1}U\|_{L^2(\mathbb{T})}\right)\|w\|_{H^{s,\gamma}_g}^2.
\end{aligned}
\end{eqnarray}
Secondly,  we have the fact that
\begin{eqnarray}
\partial _{y}g_{s} | _{y=0}&=&\left(\partial_{x}^{s}\partial _{y}w-\frac{\partial _{y}w}{w} \partial_{x}^{s}w -\frac{\partial^{2} _{y}w}{w} \partial_{x}^{s}(u-U) +\frac{\partial_{y}w\partial_{y}w}{w^{2}} \partial_{x}^{s}(u-U) \right) | _{y=0} \nonumber\\
&=&\left(\partial_{x}^{s}(\partial _{x}p-U)-ag_{s} +\frac{\partial^{2} _{y}w}{w} \partial_{x}^{s}U\right) | _{y=0}.  \label{3.10}
\end{eqnarray}
Then according to the trace theorem, we get
\begin{eqnarray}
&&\left|\int_{\mathbb{T} } g_{s}\partial _{y}g_{s}dx | _{y=0}  \right|\nonumber \\
 &&\leq \iint \left| ag_{s}^{2} \right|dxdy+\iint \left| \partial _{y}ag_{s}^{2} \right|dxdy+ 2\iint \left| ag_{s}\cdot \partial _{y} g_{s}  \right|dxdy
  \nonumber \\
   &&\quad+\iint (\partial _{y}g_{s}+g_{s})\partial_{x}^{s}(\partial _{x}p-U)dxdy+\iint g_s(\frac{\partial^{2} _{y}w}{w}+\frac{\partial^{3} _{y}w}{w}-\frac{\partial^{2} _{y}w\partial _{y}w}{w^2})\partial_{x}^{s}Udxdy
  \nonumber \\
  &&\quad+\iint \partial_y g_s\frac{\partial^{2} _{y}w}{w}\partial_{x}^{s}Udxdy
  \nonumber \\
&& \leq  C_{ \gamma,\sigma } \|w\|_{H^{s,\gamma}_{g}}^2(1+\|\partial_{x}^{s}U\|_{{L^{\infty}}(\mathbb{T})}^2)
+C\|\partial_{x}^{s}(\partial _{x}p-U)\|_{{L^{2}}(\mathbb{T})}^2
+\frac{1}{4}\|(1+y)^{\gamma}\partial_{y} g_s\|_{L^{2}}^{2} ,  \label{5.3.11}
\end{eqnarray}
which with $(\ref{5.2.60})$ gives the facts that $\|a\|_{L^{\infty}}\leq \delta^{-2}$ and $\|\partial _{y}a\|_{L^{\infty}}\leq \delta^{-2}+\delta^{-4}$.
Next, using H$\ddot{o}$lder inequality,  we obtain
\begin{eqnarray}
\left| 2\gamma \iint(1+y)^{2\gamma-1}g_{s}\partial _{y}  g_{s}  \right|
 &\leq& 2\gamma \|\frac{1}{1+y }  \|_{L^{\infty}} \|(1+y)^{\gamma} g_{s}\|_{L^{2}} \|(1+y)^{\gamma}\partial _{y}g_{s}\|_{L^{2}} \nonumber \\
& \leq & C_{ \gamma }  \|w\|_{H^{s,\gamma}_g}^2+\frac{1}{4}\|(1+y)^{\gamma}\partial_{y} g_{s}\|_{L^{2}}^{2} ,  \label{3.12}
\end{eqnarray}
and  using Lemma \ref{y5.4.4}
\begin{eqnarray}
\begin{aligned}
\left|   \gamma \iint(1+y)^{2\gamma-1} v  |g_{s}|^{2} \right|
& \leq \|\frac{v}{1+y}\|_{L^{\infty}}\|(1+y)^{ \gamma}g_{s}\|_{L^{2}}^{2}\\
&\leq  C_{s, \gamma,\sigma,\delta  }(\|w\|_{H^{s,\gamma}_{g}}+\|\partial_{x}^{s}U\|_{{L^2}(\mathbb{T})}) \|w\|_{H^{s,\gamma}_g}^{2}.  \label{3.13}
\end{aligned}
\end{eqnarray}
It follows from to $(\ref{5.2.60})$ and $(\ref{5.3.5})$ that
\begin{eqnarray}
\left| 2\iint(1+y)^{2\gamma}|g_{s}|^{2}\partial _{y}a  \right|
  \leq  2  \|\partial _{y}a  \|_{L^{\infty}} \|(1+y)^{\gamma} g_{s}\|_{L^{2}} ^{2}\leq C_{ \delta }  \|w\|_{H^{s,\gamma}_g}^{2},
\label{5.3.14}
\end{eqnarray}
and
\begin{eqnarray}
\left|    \sum\limits_{j=1}^{s-1}\binom{s}{j} \iint(1+y)^{2\gamma}g_{s} g_{j+1}\partial _{x}^{s-j}u  \right|
& \leq &  \sum\limits_{j=1}^{s-1}\binom{s}{j}     \|\partial _{x}^{s-j}u  \|_{L^{\infty}} \|(1+y)^{\gamma} g_{s}\|_{L^{2}}
   \|(1+y)^{\gamma} g_{j+1}\|_{L^{2}}\nonumber \\
& \leq &  C_{s, \gamma,\sigma,\delta  }(\|w\|_{H^{s,\gamma}_{g}}+\|\partial_{x}^{s}U\|_{{L^2}(\mathbb{T})})^2\|w\|_{H^{s,\gamma}_{g}}.
\label{5.3.15}
\end{eqnarray}
Lastly, for $j=2,3,\cdots, s-1$,
\begin{eqnarray}
&& \left|   \sum\limits_{j=1}^{s-1}\binom{s}{j}\iint(1+y)^{2\gamma}g_{s}\partial _{x}^{s-j}v \{\partial _{y}\partial _{x}^{j }w-a \partial _{x}^{j }w \}   \right|
 \nonumber \\
&& \leq    \sum\limits_{j=1}^{s-1}\binom{s}{j}  \|(1+y)^{\gamma} g_{s}\|_{L^{2}} \|\frac{\partial _{x}^{s-j}v } {1+y}\|_{L^{\infty}}
 \left(  \|(1+y)^{\gamma+1} \partial _{y}\partial _{x}^{j }w\|_{L^{2}}+ \|(1+y)a  \|_{L^{\infty}} \|(1+y)^{\gamma} \partial _{x}^{j }w\|_{L^{2}} \right)  \nonumber \\
&& \leq C_{s, \gamma,\sigma,\delta  }(\|w\|_{H^{s,\gamma}_{g}}+\|\partial_{x}^{s}U\|_{{L^2}(\mathbb{T})})\|w\|_{H^{s,\gamma}_{g}}^2,
\label{5.3.16}
\end{eqnarray}
which with $(\ref{5.2.60})$ gives the fact that $\|(1+y) a \|_{L^{\infty}} \leq \delta^{-2}$.  And for $j=1$,
\begin{eqnarray}
&& \left|  \iint(1+y)^{2\gamma}g_{s}\partial _{x}^{s-1}v \{\partial _{y}\partial _{x} w-a \partial _{x} w \}   \right|
 \nonumber \\
&& \leq    \|\frac{\partial _{x}^{s-1}v +y\partial_x^s U} {1+y}\|_{L^{2}} \left(  \|(1+y)^{\gamma+1} \partial _{y}\partial _{x}w\|_{L^{\infty}}+ \|(1+y)a  \|_{L^{\infty}} \|(1+y)^{\gamma} \partial _{x} w\|_{L^{\infty}} \right)\|(1+y)^{\gamma} g_{s}\|_{L^{2}}
   \nonumber \\
&& \quad+   \|\partial_x^s U\|_{L^{\infty}(\mathbb T)}\left(  \|(1+y)^{\gamma+1} \partial _{y}\partial _{x}w\|_{L^{2}}+ \|(1+y)a  \|_{L^{\infty}} \|(1+y)^{\gamma} \partial _{x} w\|_{L^{2}} \right)\|(1+y)^{\gamma} g_{s}\|_{L^{2}}\nonumber \\
&& \leq C_{s, \gamma,\sigma,\delta  }\|w\|_{H^{s,\gamma}_{g}}^2(\|w\|_{H^{s,\gamma}_{g}}+\|\partial_{x}^{s}U\|_{{L^\infty}(\mathbb{T})}).
\label{5.3.17}
\end{eqnarray}
For the term $| \iint(1+y)^{2\gamma}g_{s}g_1\partial _{x}^{s}U| $, we have
\begin{eqnarray}
&&\left| \iint(1+y)^{2\gamma}g_{s}g_1\partial _{x}^{s}U   \right|\nonumber \\
&&\leq \|(1+y)^{\gamma} g_{s}\|_{L^{2}}\|(1+y)^{\gamma} g_{1}\|_{L^{2}}\|\partial_{x}^{s}U\|_{{L^\infty}(\mathbb{T})}\nonumber \\
&& \leq C_{s, \gamma,\sigma,\delta  }\|\partial_{x}^{s}U\|_{{L^\infty}(\mathbb{T})}(\|w\|_{H^{s,\gamma}_{g}}+\|\partial_{x}^{s}U\|_{{L^2}(\mathbb{T})})\|w\|_{H^{s,\gamma}_{g}} .
\end{eqnarray}
For the last term, we have
\begin{eqnarray}
&&\left|\iint(1+y)^{2\gamma}g_{s}a \partial _{x}^{j}(u-U)\partial _{x}^{s-j+1}U  \right|\nonumber \\
&&\leq \|(1+y)^{\gamma} g_{s}\|_{L^{2}}\|(1+\gamma)a\|_{L^{\infty}}\|(1+y)^{\gamma-1} \partial _{x}^{j}(u-U)\|_{L^{2}}\|\partial_{x}^{s-j+1}U\|_{{L^\infty}(\mathbb{T})}\nonumber \\
&& \leq C_{s, \gamma,\sigma,\delta  }\|\partial_{x}^{s+1}U\|_{{L^\infty}(\mathbb{T})}(\|w\|_{H^{s,\gamma}_{g}}+\|\partial_{x}^{s}U\|_{{L^2}(\mathbb{T})})\|w\|_{H^{s,\gamma}_{g}}.
\end{eqnarray}
Hence combining $(\ref{5.3.11})$-$(\ref{5.3.17})$ and $(\ref{5.3.7})$ leads to (\ref{5.3.00003}).
\hfill $\Box$

\subsection{Weighted $H^s$ estimate on $w$ }
Now we can derive the weighted $H^s$ estimate on $w$ by employing Proposition \ref{p5.2.1} and Proposition \ref{p5.2.2}.
\begin{Proposition}\label{p5.2.3}
Under the same assumption of Proposition \ref{p5.2.1}, we have the following estimate:
\begin{eqnarray}
\begin{aligned}
\|w\|_{H^{s,\gamma}_{g}}^2 &\leq \left\{\|w_0\|_{H^{s,\gamma}_{g}}^2+\int_0^t F(\tau)d \tau\right\}\\
&\quad \times\left\{1-C(\frac{s}{2}-1)\left(\|w_0\|_{H^{s,\gamma}_{g}}^2+\int_0^t F(\tau) d\tau\right)^{\frac{s-2}{2}}t\right\}^{-\frac{2}{s-2}},
\label{5.2.011}
\end{aligned}
\end{eqnarray}
where $C>0$ is a constant independent of $\epsilon$ and $t$. The function $F(t)$ is expressed by
\begin{eqnarray}
\begin{aligned}
F(t)=\mathcal{P}(\|\partial_{x}^{s+1}U\|_{{L^\infty}(\mathbb{T})})+C_s \sum\limits_{l=0}^{\frac{s}{2}}\|\partial^{l}_t(\partial_xP-U)\|_{H^{s-2l}(\mathbb{T})}^{2},
\end{aligned}
\end{eqnarray}
and $\mathcal{P}(\cdot)$ denotes a polynomial.
\end{Proposition}
\textbf{Proof}.
According to Proposition \ref{p5.2.1} and Proposition \ref{p5.2.2}, we know from the definition of $\|\cdot\|_{H^{s,\gamma}_g}$ that
\begin{eqnarray}
\frac{d}{dt} \|w\|_{H^{s,\gamma}_{g}}^2 &\leq& C_{s, \gamma,\sigma,\delta  }\left(1+\|w\|_{H^{s,\gamma}_{g}}+\|\partial_{x}^{s}U\|_{{L^\infty}(\mathbb{T})}\right)
\left(\|w\|_{H^{s,\gamma}_{g}}+\|\partial_{x}^{s+1}U\|_{{L^\infty}(\mathbb{T})}\right)\|w\|_{H^{s,\gamma}_{g}}\nonumber\\
&\quad&+C_{ \gamma,\delta  }\|\partial_{x}^{s}U\|_{{L^{\infty}}(\mathbb{T})}^2 \|w\|_{H^{s,\gamma}_{g}}^2+C_{s,\gamma,\sigma,\delta}(1+\|w\|_{H_g^{s,\gamma}})^{s-2}\|w\|^2_{H^{s,\gamma}_g}\nonumber\\
&\quad&+C_s \sum\limits_{l=0}^{\frac{s}{2}}\|\partial^{l}_t(\partial_xP-U)\|_{H^{s-2l}(\mathbb{T})}^{2}\nonumber\\
&\leq& C_{s,\gamma,\sigma,\delta}\|w\|_{H_g^{s,\gamma}}^{s}+\mathcal{P}(\|\partial_{x}^{s+1}U\|_{{L^\infty}(\mathbb{T})})+C_s \sum\limits_{l=0}^{\frac{s}{2}}\|\partial^{l}_t(\partial_xP-U)\|_{H^{s-2l}(\mathbb{T})}^{2},
\end{eqnarray}
and hence, it follows by using the comparison principle of ordinary differential equations that
\begin{eqnarray}
\begin{aligned}
\|w\|_{H^{s,\gamma}_{g}}^2 &\leq \left\{\|w_0\|_{H^{s,\gamma}_{g}}^2+\int_0^t F(\tau)d \tau\right\}\\
&\quad \times \left\{1-C(\frac{s}{2}-1)\left(\|w_0\|_{H^{s,\gamma}_{g}}^2+\int_0^t F(\tau) d\tau\right)^{\frac{s-2}{2}}t\right\}^{-\frac{2}{s-2}},
\end{aligned}
\end{eqnarray}
provided that
\begin{eqnarray}
\begin{aligned}
1-C(\frac{s}{2}-1)\left(\|w_0\|_{H^{s,\gamma}_{g}}^2+\int_0^t F(\tau) d\tau\right)^{\frac{s-2}{2}}t > 0,
\end{aligned}
\end{eqnarray}
this proves Proposition \ref{p5.2.3}.
\hfill $\Box$
\subsection{Weighted $L^{\infty}$ estimates on lower order terms  }\label{su5.4}
In this subsection, we will estimate the weighted $L^{\infty}$ on $D^{\alpha} w$ for $|\alpha|\leq 2$ by using the classical maximum principles. More precisely, we will derive two parts: an $L^\infty$-estimate on  $I:=\sum \limits_{|\alpha|\leq 2} |(1+y)^{\gamma+\alpha_{2}}D^{\alpha}w|_{L^{2}}$ and a lower bound estimate on $B_{(0,0)}:=(1+y)^{\sigma} w$.
\begin{Lemma}\label{y5.2.3}
Under the same assumption of Proposition \ref{p5.2.1}, we have the following estimate: \\
For any $s \geq 4$,
\begin{eqnarray}
\begin{aligned}
\|I(t)\|_{L^\infty(\mathbb{T}\times\mathbb{R}^+)}\leq \max \left\{           \|I(0)\|_{L^\infty(\mathbb{T}\times\mathbb{R}^+)},6C^2W(t)^2
\right\} e^{C\left(1+G(t)\right)t},
\label{5.3.31}
\end{aligned}
\end{eqnarray}
and for any $s\geq 6$,
\begin{eqnarray}
\begin{aligned}
\|I(t)\|_{L^\infty(\mathbb{T}\times\mathbb{R}^+)}\leq  \left\{           \|I(0)\|_{L^\infty(\mathbb{T}\times\mathbb{R}^+)}+C(1+W(t))W(t)^2t
\right\} e^{C\left(1+G(t)\right)t}.
\label{5.3.32}
\end{aligned}
\end{eqnarray}
In addition, if $s \geq 4$, we also have,
\begin{eqnarray}
\begin{aligned}
\mathop{\min}\limits_{\mathbb{T}\times\mathbb{R}^+}(1+y)^{\sigma}w(t)\geq \max \left\{         1-C\left(1+G(t)\right)te^{C\left(1+G(t)\right)t}
\right\}\cdot\min \left\{\mathop{\min}\limits_{\mathbb{T}\times\mathbb{R}^+}(1+y)^{\sigma}w_0-CW(t)t
\right\},
\label{5.3.33}
\end{aligned}
\end{eqnarray}
where positive constant $C$ depends on $s, \gamma,\sigma$, and $\delta$ only. The functions W and $G : [0,T] \rightarrow \mathbb{R}^+$ are respectively defined by
\begin{eqnarray}
G(t):=\mathop{\sup}\limits_{ [0,t]}\|w(\tau)\|_{H^{s,\gamma}_{g}}+\mathop{\sup}\limits_{ [0,t]}\|\partial_{x}^{s}U(\tau)\|_{{L^2}(\mathbb{T})} \quad and \quad W(t):=\mathop{\sup}\limits_{ [0,t]}\|w(\tau)\|_{H^{s,\gamma}_{g}}.
\label{5.2.012}
\end{eqnarray}
\end{Lemma}
\textbf{Proof}.
For simplicity, we denote,
\begin{eqnarray*}
\begin{aligned}
I:=\sum \limits_{|\alpha|\leq 2} \left|(1+y)^{\sigma+\alpha_{2}}D^{\alpha}w\right|^2
\end{aligned}
\end{eqnarray*}
and
\begin{eqnarray*}
\begin{aligned}
B_{\alpha}:=(1+y)^{\sigma}D^{\alpha}w,
\end{aligned}
\end{eqnarray*}
then $B_{(0,0)}=(1+y)^{\gamma}w$.
By direct computations, we obtain
\begin{eqnarray}
\begin{aligned}
\left(\partial_t+u\partial_x+v\partial_y-\partial^2_y-\epsilon^2\partial^2_x\right)B_{\alpha}=\sum\limits_{i=1}^{3} S_{i},
\label{5.3.34}
\end{aligned}
\end{eqnarray}
where
\begin{eqnarray*}
\begin{aligned}
S_1=\left(\frac{\sigma+\alpha_2}{1+y}v+\frac{(\sigma+\alpha_2)(\sigma+\alpha_2-1)}{(1+y)^2}-1\right)B_{\alpha},\ \ S_2=-\frac{2(\sigma+\alpha_2)}{1+y}\partial_y B_{\alpha},
\end{aligned}
\end{eqnarray*}
and
\begin{eqnarray*}
\begin{aligned}
S_3=-\sum\limits_{0< \beta \leq \alpha}\binom{\alpha}{\beta}\left\{(1+y)^{\beta_2}\left(D^{\beta}u B_{\alpha-\beta+e_1}+\frac{D^{\beta} v B_{\alpha-\beta+e_2}}{1+y}\right)\right\}.
\end{aligned}
\end{eqnarray*}
Multiplying the equation (\ref{5.3.34}) by $2B_{\alpha}$, we conclude
\begin{eqnarray}
&&\left(\partial_t+u\partial_x+v\partial_y-\epsilon^2\partial^2_x-\partial^2_y\right)I\nonumber\\
&&\leq2\sum \limits_{|\alpha|\leq 2} \left(|S_1 B_{\alpha}|+|S_2 B_{\alpha}|+|S_3 B_{\alpha}|-\epsilon^2|\partial_x B_{\alpha}|^2-|\partial_y B_{\alpha}|^2\right)\nonumber\\
&&\leq  C_{s, \gamma,\sigma,\delta  }( 1+\|w(s)\|_{H^{s,\gamma}_{g}}+\|\partial_{x}^{s}U\|_{{L^2}(\mathbb{T})})I+C_{\delta}I+\sum \limits_{|\alpha|\leq 2}|\partial_y B_{\alpha}|^2-2\epsilon^2\sum \limits_{|\alpha|\leq 2}|\partial_x B_{\alpha}|^2\nonumber\\
&&\quad -2\sum \limits_{|\alpha|\leq 2}|\partial_y B_{\alpha}|^2+2\|(1+y)^{\beta_2}(D^{\beta}u +\frac{D^{\beta} v }{1+y})\|_{L^{\infty}}\sum \limits_{|\alpha|\leq 2}\left(\sum\limits_{0\leq \beta < \alpha}|B_{\alpha-\beta+e_1}+B_{\alpha-\beta+e_2}|\right)B_{\alpha}\nonumber\\
&&\leq  C_{s, \gamma,\sigma,\delta  }( 1+\|w(s)\|_{H^{s,\gamma}_{g}}+\|\partial_{x}^{s}U\|_{{L^2}(\mathbb{T})})I,
\end{eqnarray}
where we have used Lemma \ref{y5.4.4} and Young inequality.
Applying the classical maximum principle (see Lemma \ref{y5.1.6} for instance) for parabolic equations, we have
\begin{eqnarray}
\begin{aligned}
&\|I(t)\|_{L^\infty(T\times\mathbb{R}^+)}\\
&\leq \max \left\{            e^{C\left(1+G(t)\right)t}\|I(0)\|_{L^\infty(T\times\mathbb{R}^+)},\mathop{\max}\limits_{\tau \in [0,t]}\left(e^{C\left(1+G(t)\right)(t-\tau)}\|I(\tau)|_{y=0}\|_{L^\infty(\mathbb{T})}\right)
\right\}.
\label{5.3.36}
\end{aligned}
\end{eqnarray}
 To derive a lower bound estimate on $B_{(0,0)}$, we have
 \begin{eqnarray}
\begin{aligned}
\left(\partial_t+u\partial_x+(v+\frac{2\sigma}{1+y})\partial_y-\partial^2_y-\epsilon^2\partial^2_x\right)B_{(0,0)}=\left(\frac{\sigma}{1+y}v+\frac{\sigma(\sigma-1)}{(1+y)^2}-1\right)B_{(0,0)}.
\label{5.3.37}
\end{aligned}
\end{eqnarray}
Applying the classical minimum principle  (see Lemma \ref{y5.1.7} for instance) for (\ref{5.3.37}), we get
\begin{eqnarray}
\begin{aligned}
&\mathop{\min}\limits_{\mathbb{T}\times\mathbb{R}^+}(1+y)^{\sigma}w(t)\\
&\geq \max \left\{         1-C\left(1+G(t)\right)te^{C\left(1+G(t)\right)t}
\right\}\cdot\min \left\{\mathop{\min}\limits_{\mathbb{T}\times\mathbb{R}^+}(1+y)^{\sigma}w_0,        \mathop{\min}\limits_{[0,t]\times\mathbb{T}}w|_{y=0}
\right\}.
\label{5.3.38}
\end{aligned}
\end{eqnarray}
Next we start estimating the boundary values, using  Lemma \ref{y5.4.3}, we obtain
\begin{eqnarray}
\begin{aligned}
\|I(\tau)|_{y=0}\|_{L^\infty(\mathbb{T})}&\leq 3C^2\sum \limits_{|\alpha|\leq 2}\left(\|D^{\alpha}w\|_{L^2}^2+\|\partial_xD^{\alpha}w\|_{L^2}^2+\|\partial_y^2 D^{\alpha}w\|_{L^2}^2\right)\\
&\leq 6C^2\|w\|_{H^{s,\gamma}_{g}}^{2},
\end{aligned}
\end{eqnarray}
which, together with (\ref{5.3.36}), gives (\ref{5.3.31}).

According to (\ref{5.1.5}) and boundary condition $\partial^{\alpha_1}_x u=\partial^{\alpha_1}_x v=0$, we have
\begin{eqnarray}
\begin{aligned}
\partial_t D^{\alpha}w|_{y=0}=D^{\alpha}(\epsilon^2\partial_x^2+\partial_y^2 w-w-u\partial_x w-v\partial_y w)|_{y=0}.
\end{aligned}
\end{eqnarray}
When $\alpha_2=0$ and $\alpha_1 \leq 2$,
\begin{eqnarray}
\begin{aligned}
\partial_t D^{\alpha}w|_{y=0}=(\epsilon^2\partial^{\alpha_1+2}_{x} w+\partial^{\alpha_1}_{x}\partial_y^2 w-\partial^{\alpha_1}_{x}w)|_{y=0}.
\end{aligned}
\end{eqnarray}
When $1 \leq \alpha_2 \leq 2$ and $\alpha_1=0$,
\begin{eqnarray}
\begin{aligned}
\partial_t D^{\alpha}w|_{y=0}
&=\left(\epsilon^2\partial^{\alpha_2}_{y}\partial^{2}_{x} w+\partial^{\alpha_2} _{y}\partial _{y}^{2}w-\partial^{\alpha_2} _{y}w-\partial^{\alpha_2} _{y}(u\partial_x w+v\partial_y w)\right)|_{y=0}\\
&=\left(\epsilon^2\partial^{\alpha_2}_{y}\partial^{2}_{x} w+\partial^{\alpha_2} _{y}\partial _{y}^{2}w-\partial^{\alpha_2} _{y}w-\alpha_2 w \partial_x \partial^{\alpha_2-1} _{y}w\right)|_{y=0}.
\end{aligned}
\end{eqnarray}
When $\alpha_2=1$ and $\alpha_1 =1$,
\begin{eqnarray}
\begin{aligned}
\partial_t D^{\alpha}w|_{y=0}=(\epsilon^2\partial_{y}\partial^{3}_{x} w+\partial _{x}\partial _{y}^{3}w-\partial _{x}\partial _{y}w-\partial _{x}w\partial _{x}w-w\partial _{x}^2w)|_{y=0}.
\end{aligned}
\end{eqnarray}
For $s\geq 6$, using Lemma \ref{y5.4.4}, we have
 \begin{eqnarray}
\begin{aligned}
\|\partial_t I|_{y=0}\|_{{L^\infty}(\mathbb{T})}&\leq C_{s, \gamma  } \|D^{\alpha}w\partial_t D^{\alpha}w|_{y=0}\|_{{L^\infty}(\mathbb{T})}\\
&\leq  C_{s, \gamma  }\|w(\tau)\|_{H^{s,\gamma}_{g}}(\|w(\tau)\|_{H^{s,\gamma}_{g}}+\|w(\tau)\|_{H^{s,\gamma}_{g}}^2).
\end{aligned}
\end{eqnarray}
Hence, a direct integration yields
 \begin{eqnarray}
\begin{aligned}
\|I(t)|_{y=0}\|_{{L^\infty}(\mathbb{T})}\leq  \|I(0)|_{y=0}\|_{{L^\infty}(\mathbb{T})}+C_{s, \gamma  } (1+W(t))W(t)^2t,
\end{aligned}
\end{eqnarray}
which, together with (\ref{5.3.36}), gives (\ref{5.3.32}). For $s\geq4$,
 \begin{eqnarray}
\begin{aligned}
\|\partial_t w|_{y=0}\|_{{L^\infty}(\mathbb{T})}&=\|(\partial_y^2 w-w)|_{y=0}\|_{{L^\infty}(\mathbb{T})}\\
&\leq  CW(t).
\end{aligned}
\end{eqnarray}
Hence, a direct integration gives
 \begin{eqnarray}
\begin{aligned}
\mathop{\min}\limits_{\mathbb{T}}{w(t)}|_{y=0}\geq \mathop{\min}\limits_{\mathbb{T}}{w_0}|_{y=0}-CW(t)t,
\end{aligned}
\end{eqnarray}
which together with (\ref{5.3.38}) implies (\ref{5.3.33}). The proof of Lemma \ref{y5.2.3} is thus complete.
\hfill $\Box$
\section{Local-in-time existence and uniqueness}
\subsection{Local-in-time existence }
In this subsection, we go back to use the symbol $(u^\epsilon, v^\epsilon, w^\epsilon)$ instead of $(u, v, w)$ from Subsections \ref{su5.1}-
\ref{su5.4} to denote the solution to the regularized system (\ref{5.2.0001}). To obtain   the local-in-time solution of the initial-boundary value problem (\ref{5.1.1})-(\ref{5.1.2}), we will construct the solution to the Prandtl equations (\ref{5.1.1}) by passing to the limit $\epsilon \rightarrow 0^+$ in the regularized Prandtl equations (\ref{5.2.0001}). We only sketch the proof into five steps and more details can be found in \cite{mw}.\\
$\textbf{Step 1.}$ According to the definition of $F$, assumption (\ref{5.1.8}), and the regularized Bernoulli's law,
\begin{eqnarray}
\|F\|_{L^\infty} \leq M < +\infty.
\end{eqnarray}
Thus we derive the uniform estimate for any $\epsilon \in [0,1]$ and any $t \in [0,T_1]$
\begin{eqnarray}
\|w^\epsilon\|_{H^{s,\gamma}_{g}} \leq 4 \|w_0^\epsilon\|_{H^{s,\gamma}_{g}},
\end{eqnarray}
provided that $T_1$ is chosen by (\ref{5.2.011}) such that
\begin{eqnarray*}
T_1 :=\min \left\{\frac{3\|w_0\|_{H^{s,\gamma}_{g}}^2}{C_{s, \gamma,\sigma,\delta  }M},\frac{1-2^{2-s}}{2^{s-2}C_{s, \gamma,\sigma,\delta  }\|w_0\|_{H^{s,\gamma}_{g}}^{s-2}}\right\}.
\end{eqnarray*}
$\textbf{Step 2.}$ When $s \geq 6$, we know from definition (\ref{5.2.012}) of $\Omega$ and $G$ that for any $t \in [0,T_1]$,
\begin{eqnarray}
\Omega(t) \leq 4 \|w_0^\epsilon\|_{H^{s,\gamma}_{g}} \quad and \quad G(t) \leq 4 \|w_0^\epsilon\|_{H^{s,\gamma}_{g}}+M.
\label{5.2.015}
\end{eqnarray}
Thus, if we choose
\begin{eqnarray*}
T_2 :=\min \left\{T_1,\frac{1}{64\delta^2 C_{s, \gamma }(1+4 \|w_0\|_{H^{s,\gamma}_{g}}) \|w_0\|_{H^{s,\gamma}_{g}}^2},\frac{\ln 2}{C_{s, \gamma,\sigma,\delta  }(1+4 \|w_0\|_{H^{s,\gamma}_{g}}+M)}\right\},
\end{eqnarray*}
then using inequality (\ref{5.3.32}) and initial assumption
\begin{eqnarray}
\sum \limits_{|\alpha|\leq 2}  |(1+y)^{\sigma+\alpha_{2}}D^{\alpha}w | ^{2}\leq\frac{1}{4\delta^{2}},
\end{eqnarray}
we have the upper bound
\begin{eqnarray}
\|\sum \limits_{|\alpha|\leq 2}  |(1+y)^{\sigma+\alpha_{2}}D^{\alpha}w^\epsilon (t) | ^{2}\|_{L^{\infty}(\mathbb{T}\times \mathbb{R}^+)}\leq\frac{1}{\delta^{2}}
\label{5.2.019}
\end{eqnarray}
for all $t \in [0,T_2]$. When $s \geq 4$, from the initial hypothesis $\|\omega_0\| \leq C\delta^{-1}$, we have the same estimate (\ref{5.2.019}) for all $t\in [0,T_2]$.\\
$\textbf{Step 3.}$
Let us choose
\begin{eqnarray*}
T_3 :=\min \left\{T_1,\frac{\delta}{8 C_{s, \gamma }\|w_0\|_{H^{s,\gamma}_{g}}},\frac{1}{6C_{s, \gamma,\sigma,\delta  }N},\frac{\ln 2}{C_{s, \gamma,\sigma,\delta  }N}\right\},
\end{eqnarray*}
where $N :=1+4\|w_0\|_{H^{s,\gamma}_{g}}+M$. Then by using (\ref{5.2.012}) and (\ref{5.2.015}), we derive the uniform estimate for any $\epsilon \in [0,1]$ and any $t \in [0,T_3]$,
\begin{eqnarray}
\min \limits_{\mathbb{T}\times \mathbb{R}^+}(1+y)^{\sigma} w^\epsilon(t) \geq \delta.
\end{eqnarray}
$\textbf{Step 4.}$
In summary, the above uniform estimates  hold for any $t \in [0,T]$, if $T$ is chosen to satisfy $T=\min\{T_1,T_2,T_3\}$. Further using almost equivalence relation  (\ref{5.4.1}), we have
\begin{eqnarray}
\sup  \limits_{0 \leq t \leq T}( \|w^\epsilon\|_{H^{s,\gamma}}+\|u^\epsilon-U\|_{H^{s,\gamma-1}}) \leq C\left( 4\|w_0\|_{H^{s,\gamma}_{g}}+\sup \limits_{0 \leq t \leq T}\|\partial_x^s U\|_{L^2}\right) < + \infty.
\label{5.2.016}
\end{eqnarray}
From the equation (\ref{5.2.0001}), system (\ref{5.2.0002}), (\ref{5.2.016}) and Lemma \ref{y5.4.3}, we also have $\partial_t w^{\epsilon}$	 and $\partial_t (u^{\epsilon}-U)$ are uniformly
bounded in $L^{\infty}([0,T];H^{s-2,\gamma})$ and $L^{\infty}([0,T];H^{s-2,\gamma-1})$ respectively. By Lions-Aubin lemma and the compact embedding of $H^{s,\gamma}$ in $H^{s^{\prime}}_{loc}$
, we have after taking a subsequence, as $\epsilon_k \rightarrow 0^+$,
\begin{equation}\left\{
\begin{array}{ll}
\omega^{\epsilon_k} \stackrel{*}{\rightharpoonup} \omega, &{\rm in} \quad L^{\infty}([0,T];H^{s,\gamma}),\\
\omega^{\epsilon_k} \rightarrow \omega, &{\rm in} \quad C([0,T];H^{s^{\prime}}_{loc}),\\
u^{\epsilon_k}-U \stackrel{*}{\rightharpoonup} u-U,  &{\rm in} \quad L^{\infty}([0,T];H^{s,\gamma-1}),\\
u^{\epsilon_k} \rightarrow u, &{\rm in} \quad C([0,T];H^{s^{\prime}}_{loc}),
\end{array}
         \right.\end{equation}
for all $s^{\prime} < s$, where
\begin{equation}\left\{
\begin{array}{ll}
w=\partial_y u \in L^{\infty} ([0,T];H^{s,\gamma})\cap \bigcap_{s^\prime <s}C([0,T];H^{s^{\prime}}_{loc}),\\
 u-U \in L^{\infty} ([0,T];H^{s,\gamma-1})\cap \bigcap_{s^\prime <s}C([0,T];H^{s^{\prime}}_{loc}).
\end{array}
         \right.\end{equation}
$\textbf{Step 5.}$
Using the local uniform
convergence of $\partial_x u^{\epsilon_k}$, we also have the pointwise convergence of $v^{\epsilon_k}:$ as $\epsilon \rightarrow 0^+$,
\begin{eqnarray}
v^{\epsilon_k}=-\int^y_0 \partial_x u ^{\epsilon_k}~ dy \rightarrow -\int^y_0 \partial_x u ~dy=: v.
\end{eqnarray}
Thus, passing the limit $\epsilon_k \rightarrow 0^+$ in the initial-boundary value problem (\ref{5.2.0001}), we get that the limit
$(u, v)$ solves the initial-boundary value problem (\ref{5.1.4}),(\ref{5.1.2}) in the classical sense. Furthermore, we obtain
$(u, v, b)$ solves the initial-boundary value problem (\ref{5.1.1})-(\ref{5.1.2}) in the classical sense. This completes the proof
of the existence.
\subsection{Uniqueness of solutions}
In this subsection, we are going to prove the uniqueness to 2D magnetic Prandtl model.  Let us denote $(\bar u -\bar v)=(u_1,v_1)-(u_2,v_2)$, $\bar w =w_1-w_2$, $\bar b_1 =b_{11}-b_{12}$ and $a_2=\frac{\partial_y w_2}{w_2}$. It is easy to check that $\bar g =\bar w -a_2 \bar u=w_2 \partial_y(\frac{\bar u}{w_2})$ and the evolution equation on $\bar g$ is as follows
\begin{eqnarray}
\begin{aligned}
(\partial_t+u_1\partial_x+v_1\partial_y-\partial^2_y+1)\bar g
&=(\partial_t+u_1\partial_x+v_1\partial_y-\partial^2_y+1)\bar w-a_2(\partial_t+u_1\partial_x+v_1\partial_y-\partial^2_y+1)\bar u\\
&\quad -\bar u(\partial_t+u_1\partial_x+v_1\partial_y-\partial^2_y+1)a_2-2\partial_y a_2 \bar w.
 \label{5.33.40}
\end{aligned}
\end{eqnarray}
Next, we calculate the values of the first three terms on the right-hand side of the equality respectively, recalling the vorticity system (\ref{5.1.5}), we have
\begin{eqnarray}
\begin{aligned}
(\partial_t+u_i\partial_x+v_i\partial_y+1)\partial_y w_i=\partial_y^3 w_i+\partial_x u_i \partial_y w_i -w_i \partial_x w_i,
\end{aligned}
\end{eqnarray}
then, according to the definition of $a_i$, we get
\begin{eqnarray}
(\partial_t+u_i\partial_x+v_i\partial_y+1)a_i&=&\frac{(\partial_t+u_i\partial_x+v_i\partial_y+1)\partial_y w_i}{w_i}-\frac{\partial_y w_i(\partial_t+u_i\partial_x+v_i\partial_y+1) w_i}{w_i^2}\nonumber\\
&=&\frac{\partial^3_y w_i}{w_i}+a_i\partial_x u_i-\partial_x w_i -a_i \frac{\partial_y^2 a_i}{w_i},
\end{eqnarray}
which, combined with the fact
\begin{eqnarray*}
\frac{\partial^3_y w_i}{w_i}-a_i \frac{\partial_y^2 a_i}{w_i}=\partial_y^2 a_i+2a_i \partial_y a_i,
\end{eqnarray*}
yields,
\begin{eqnarray}
\begin{aligned}
(\partial_t+u_i\partial_x+v_i\partial_y-\partial_y^2+1)a_i=a_i\partial_x u_i-\partial_x w_i +2a_i \partial_y a_i.
\end{aligned}
\end{eqnarray}
Furthermore, we conclude that
\begin{eqnarray}
\begin{aligned}
(\partial_t+u_1\partial_x+v_1\partial_y-\partial_y^2+1)a_2
=a_2\partial_x u_2-\partial_x w_2 +2a_2 \partial_y a_2+(\bar u \partial_x+\bar v \partial_y)a_2.
 \label{5.33.41}
\end{aligned}
\end{eqnarray}
For the estimate of $(\partial_t+u_1\partial_x+v_1\partial_y-\partial^2_y+1)\bar u$ and $(\partial_t+u_1\partial_x+v_1\partial_y-\partial^2_y+1)\bar w$, from the equation (\ref{5.1.4}) and (\ref{5.1.5}), we can derive that
\begin{eqnarray}
\begin{aligned}
(\partial_t+u_1\partial_x+v_1\partial_y-\partial^2_y+1)\bar u=-\bar u \partial_x u_2-\bar v \partial_y  u_2,
\end{aligned}
\end{eqnarray}
and
\begin{eqnarray}
\begin{aligned}
(\partial_t+u_1\partial_x+v_1\partial_y-\partial^2_y+1)\bar w=-\bar u \partial_x w_2-\bar v \partial_y  w_2.
 \label{5.33.42}
\end{aligned}
\end{eqnarray}
Using (\ref{5.33.40}) and (\ref{5.33.41})-(\ref{5.33.42}), we have
\begin{eqnarray}
\begin{aligned}
(\partial_t+u_1\partial_x+v_1\partial_y-\partial^2_y+1)\bar g=-2\bar w \partial_y a_2 -\bar u(\bar u \partial_x a_2+\bar v \partial_y a_2+2a \partial_y a_2).
\end{aligned}
\end{eqnarray}
Now we derive $L^2$ estimate on $\bar g$. For any $t \in (0,T]$, multiplying by $2 \bar g$ and then integrating by parts over $\mathbb{T}\times\mathbb{R}_{+}$, we obtain
\begin{eqnarray}
\begin{aligned}
&\frac{d}{dt}\|\bar g\|_{L^2}^2+2\|\bar g\|_{L^2}^2+2\|\partial_y \bar g\|_{L^2}^2\\
&=\int_{\mathbb{T}}\bar g \partial_y \bar g |_{y=0}dx-2\iint \bar g \bar u(\bar u \partial_x a_2+\bar v \partial_y a_2+2a \partial_y a_2)-4\iint \bar g \bar w \partial_y a_2  \\
&\quad -2\iint \bar g(u_1\partial_x \bar g+v_1\partial_y \bar g).
 \label{5.33.43}
\end{aligned}
\end{eqnarray}
We now need to estimate the integral equation (\ref{5.33.43}) term by term. Indeed, applying the simple trace theorem and Young's inequality,
\begin{eqnarray}
|\int_{\mathbb{T}}\bar g \partial_y \bar g |_{y=0}dx|
&\leq&|\iint a_2 |\bar g|^2  dxdy|+|\iint \partial_y a_2 |\bar g|^2  dxdy|+2|\iint  a_2 \bar g \partial_y \bar g  dxdy|\nonumber\\
&\leq& \frac{1}{2}\|\partial_y \bar g\|_{L^2}^2+C_{\sigma,\delta  }\|\bar g\|_{L^2}^2,
\end{eqnarray}
where we have used the fact $ \partial_y \bar g |_{y=0}=-a_2 \bar g  |_{y=0}$ ($\partial_y \bar w=0$). We claim $\|\frac{\bar u}{1+y}\|_{L^2} \leq C_{\sigma,\delta} \|\bar g\|_{L^2}$, so by Lemma \ref{y5.4.4},
\begin{eqnarray}
-2\iint \bar g \bar u(\bar u \partial_x a_2+\bar v \partial_y a_2+2a \partial_y a_2) &\leq&  2\|(1+y)(\bar u \partial_x a_2+\bar v \partial_y a_2+2a \partial_y a_2)\|_{L^\infty}\|\frac{\bar u}{1+y}\|_{L^2}\|\bar g \|_{L^2}\nonumber\\
&\leq& C_{ \gamma,\sigma,\delta  }(1+\|w_i\|_{H^{4,\gamma}_{g}}+\|\partial_{x}^{4}U\|_{{L^2}(\mathbb{T})})\|\bar g\|_{L^2}^2.
\end{eqnarray}
Below we give the fact that $\|\frac{\bar u}{1+y}\|_{L^2}$ can be controlled by $\|\bar g\|_{L^2}$, since $\delta \leq (1+y)^{\delta} w_2\leq \delta^{-1}$,
\begin{eqnarray}
\begin{aligned}
\|\frac{\bar u}{1+y}\|_{L^2} \leq \delta^{-1} \|(1+y)^{-\sigma-1}\frac{\bar u}{w^2}\|_{L^2} \leq C_{\sigma,\delta}\|(1+y)^{-\sigma} \partial(\frac{\bar u}{w_2})\|_{L_2} \leq C_{\sigma,\delta} \|\bar g\|_{L^2}.
\end{aligned}
\end{eqnarray}
in addition, we also have
\begin{eqnarray}
\begin{aligned}
\| \bar w\|_{L^2} \leq \|\bar g \|_{L^2}+\delta^{-2}\|\frac{\bar u}{1+y}\|_{L^2}\leq C_{\sigma,\delta} \|\bar g\|_{L^2},
\end{aligned}
\end{eqnarray}
and
\begin{eqnarray}
\begin{aligned}
-4\iint  \bar g \bar w \partial_y a_2 \leq C_{\sigma,\delta} \|\bar g\|_{L^2}^2.
\end{aligned}
\end{eqnarray}
For the last term in (\ref{5.33.43}), using the integration by parts, boundary condition $(u_1, v_1)|_{y=0}$ and $\partial_x u_1+\partial_y v_1=0$, we easily show
\begin{eqnarray}
\begin{aligned}
-2\iint \bar g(u_1\partial_x \bar g+v_1\partial_y \bar g)=0.
\end{aligned}
\end{eqnarray}
Combining all of the above estimates, we have
\begin{eqnarray}
\begin{aligned}
\frac{d}{dt}\|\bar g\|_{L^2}^2 \leq C_{ \gamma,\sigma,\delta  }(1+\|w_i\|_{H^{4,\gamma}_{g}}+\|\partial_{x}^{4}U\|_{{L^2}(\mathbb{T})})\|\bar g\|_{L^2}^2.
\end{aligned}
\end{eqnarray}
 which,  by Gronwall's inequality, gives
\begin{eqnarray}
\begin{aligned}
\|\bar g(t)\|_{L^2}^2 \leq \|\bar g(0)\|_{L^2}^2e^{Ct},
\end{aligned}
\end{eqnarray}
here, $C=C_{ \gamma,\sigma,\delta  }(1+\|w_i\|_{H^{4,\gamma}_{g}}+\|\partial_{x}^{4}U\|_{{L^2}(\mathbb{T})})$, and this implies $\bar g=0$ due to $u_1|_{t=0}=u_2|_{t=0}$. Since $w_2\partial_y (\frac{u_1-u_2}{w_2})=\bar g =0$, we have
\begin{eqnarray}
\begin{aligned}
u_1-u_2=q w_2
\end{aligned}
\end{eqnarray}
for some function $q=q(t,x)$. By using the  Oleinik's monotonicity assumption $w_2 >0 $ and boundary condition $u_1|_{y=0}=u_2|_{y=0}=0$, we can get $q=0$, and hence $u_1=u_2$. Further, using $\partial_x u+\partial_y v=0$ and $\partial_y u+\partial_y^2 b_1=0$, then $v_i$  and $b_{1i}$ can be uniquely determined (i.e., $v_1=v_2$ and $b_{11}=b_{12}$).

\section{Bibliographic Comments}

Now let us briefly review the background and corresponding results about the boundary layer. Actually, Prandtl proposed the basic rules for describing a phenomenon at the Heidelberg international mathematical conference in 1904. He pointed out that there are two regions in the flow of a solid: a thin layer near an object, viscous friction plays an important role; outside this thin layer, friction is negligible.  Prandtl called this thin layer as the boundary layer, which can be described by the so-called Prandtl equations. Later on, the researches on the boundary layer were appeared, such as \cite{CH,EE}. But until 1999, Oleinik and Samokhin \cite{OS} gave a systematic study of the Prandtl equations from a mathematical point of view, which is the fundamental research on boundary layer system.  The Prandtl system is obtained as  a simplification of the Navier-Stokes system and describes the motion of a fluid with small viscosity about a solid body in a thin layer which is formed near its surface owing to the adhesion of the viscous fluid to the solid surface.
The Prandtl boundary layer theory has been extensively studied in different references, see \cite{AWXY,CS,CLS,cwz,cwz1,gvm,qindong}.

The first well-known result was developed by Oleinik and Samokhin in \cite{OS}, where under the
monotonicity condition on tangential velocity with respect to the normal variable to the boundary, the local (in time) well-posedness of Prandtl equations was obtained by using the Crocco transformation and von Mises transformation.  Since then, Crocco transformation and von Mises transformation have been used in boundary layer problems. For example, Xin and Zhang \cite{XZ} established a global existence of weak solutions to the two-dimensional Prandtl system for the pressure is favourable, which had generalized the local well-posedness results of Oleinik \cite{OS}.
 Gong, Guo and Wang \cite{GGW5} investigated the local spatial existence of solutions for the compressible Navier-Stokes boundary layer equations by using von Mises transformation.
The coordinate transformation was also used to MHD boundary layer in \cite{hly}.
Xie and Yang in \cite{xy} also obtained the local existence of solutions to the MHD boundary layer system as a general shear flow.
 Thereafter, Xie and Yang \cite{lxy1,LXY}, Gao, Huang and Yao \cite{ghy} studied the local well-posedness of MHD boundary layer problems by the energy method.

There are not many literature  about the Prandtl-Hartmann regime, only \cite{CRWZ,dongqin,gly,GGW,LXY0,LWXY,xy1}, and even fewer literature  that apply coordinate transformation to the equation of Prandtl-Hartmann regime except for \cite{GGW} where Gong and  Wang obtained the global existence of a weak solution to the mixed Prandtl-Hartmann boundary layer problem by  introducing the Crocco transformation for the case that the pressure is favourable ($\partial_{x}P<0$).
The authors in \cite{CRWZ,dongqin,gly,LXY0,LWXY,xy1} studied the well-posedness of the Prandtl-Hartmann regime with the constant outer flow $U$ which means $\partial_{x}P=0$.
For example, Dong and Qin \cite{dongqin}(see also Chapter 2) established the global well-posedness of solutions in $H^s(1\leq s\leq 4)$ with $ U= \text{constant}$ and without monotonicity condition and a lower bound, and Qin and Dong \cite{qindong1} (see also Chapter 4) also proved the local well-posedness of solutions to 2D mixed Prandtl equations in Sobolev  space with $ U=\text{constant}$.

From above statement, a natural question  is that whether or not the existence of local solutions can be achieved for the Prandtl-Hartmann regime equations in both favourable case ($\partial_{x}P<0$) and infavourable case ($\partial_{x}P>0$). However, the answer to this question is still open before the paper (\cite{qinwang2}) has been finished. This means that the paper  \cite{qinwang2} (see also this chapter) has positively answered this question. In this direction, we would like to mention some works. Masmoudi and Wong  \cite{mw} proved the local existence and uniqueness of solutions for the classical Prandtl system by the energy method in a polynomial weighted  Sobolev space,  and Fan, Ruan, Yang \cite{FRY} proved the  local well-posedness for the compressible Prandtl boundary layer equations by the same method. Motivated by \cite{mw,FRY}, our purpose of this chapter is to prove the local existence of solutions  to the Prandtl-Hartmann regime model under the outer flow  $U\neq\text{constant}$. Although the method used here is the same as that in \cite{mw,FRY},  the Prandtl-Hartmann regime model has one more tangential magnetic components than Prandtl equations have, which leads to the boundary data $y=0$ being much more complicated.  To overcome this difficulty, some more delicated estimates and mathematical induction method are employed.

\chapter{ Local existence of solutions to 3D Prandtl equations with a special structure}
In this chapter, we shall consider the 3D Prandtl equation in a periodic domain and prove the local existence and uniqueness of solutions by the energy method in a polynomial weighted  Sobolev space. Compared to the existence and uniqueness of solutions to the classical Prandtl equations
where the Crocco transform has always been used with the general outer flow $U\neq\text{constant}$, this
Crocco transform is not needed here for 3D  Prandtl equations.
 We use the skill of cancellation mechanism and construct a new unknown function to show that   the
 existence and uniqueness of solutions to 3D Prandtl equations (cf.  Masmoudi and Wong, Comm. Pure Appl. Math.,
    68(10)(2015), 1683-1741) which extends from the two dimensional case in  \cite{qinwang2}
to the present three dimensional case  with a special structure. The content of this chapter is adapted from \cite{qinwang1}.
\section{Introduction}
Experimental data and theoretical analysis show that in many important practical cases, for fluids whose viscosity is small, there exists a thin transition layer near the boundary, in which the behavior of flow changes dramatically, this phenomenon is called boundary layer theory in \cite{OS}.
In this chapter, we shall consider the 3D Prandtl equations in the periodic domain  $\{(t,x,y,z)\big|t>0, (x,y)\in \mathbb{T}^{2}, z\in\mathbb{R^{+}}\}$:
\begin{equation}\left\{
\begin{array}{ll}
\partial _{t}u+(u\partial _{x} +v\partial _{y} +w\partial_{z})u+\partial _{x}P=\partial _{z}^{2}u,\\
\partial _{t}v+(u\partial _{x} +v\partial _{y} +w\partial_{z})v+\partial _{y}P=\partial _{z}^{2}v,\\
\partial _{x}u+\partial _{y}v+\partial_{z}w=0,\\
(u,v)|_{t=0}=(u_{0}(x,y,z),v_{0}(x,y,z)),\\
(u,v,w)|_{z=0}=0,\\
  \lim\limits_{z\rightarrow+\infty}(u,v)=(U(t,x,y), V(t,x,y)).
\end{array}
 \label{6.1.1}         \right.\end{equation}
Here $\mathbb{T}^{2}\subseteq \mathbb{R}^{2}$,
$(u,v,w)=(u(t,x,y,z),v(t,x,y,z),w(t,x,y,z))$ denote the velocity field, $(U(t,x,y), V(t,x,y))$ and $P$ are the traces of the tangential velocity field and the pressure of the the Euler flow respectively which satisfy Bernoulli's law
\begin{equation}\left\{
\begin{array}{ll}
\partial _{t}U+(U\partial _{x} +V\partial _{y} )U+\partial _{x}P=0,\\
\partial _{t}V+(U\partial _{x} +V\partial _{y}  )V+\partial _{y}P=0.
\end{array}
 \label{6.1.2}         \right.\end{equation}

Despite of its importance in physics, there are very few mathematical results
on the Prandtl equations in three space variables, and the three-dimensional results are based on special structures, such as \cite{lwy1,lwy}, where authors constructed a solution of the three-dimensional Prandtl equations $(\ref{6.1.1})$ with the structure
\begin{eqnarray}
 (u(t,x,y,z),K(t, x,y)u(t,x,y,z),w(t,x,y,z)),
 \label{6.1.3}
\end{eqnarray}
and the outer Euler flow takes the following form on the boundary $\{z=0\}$,
$$(U(t,x,y), K(t, x,y)U(t,x,y),0 ). $$
But in what follows, we shall consider the following the equivalent system of system (\ref{6.1.1}) (see Proposition \ref{p6.1.2}  for the specific proof),
\begin{equation}\left\{
\begin{array}{ll}
\partial _{t}u+(u\partial _{x} +K u\partial _{y} +w\partial_{z})u =\partial _{z}^{2}u^\epsilon-\partial_x P,\\
\partial _{x}u+\partial _{y}(K u)+\partial_{z}w=0,\\
 u  \big|_{t=0}= u_{0}(x,y,z) ,\\
(u,w)\big|_{z=0}=0,\quad  \lim\limits_{z\rightarrow+\infty}u=U(t, x, y),
\end{array}   \label{6.1.222}      \right.\end{equation}
where  $U$ satisfies
\begin{eqnarray}
\partial_t U+ U\partial_x U+ K U\partial_y U+\partial_x P=0.
\end{eqnarray}

Hence, the regularized equivalent system of  equations (\ref{6.1.222}) can be expressed as follows,
\begin{equation}\left\{
\begin{array}{ll}
\partial _{t}u^\epsilon+(u^\epsilon\partial _{x} +K^\epsilon u^\epsilon\partial _{y} +w^\epsilon\partial_{z})u^\epsilon -\epsilon^2\partial _{x}^{2}u^\epsilon-\epsilon^2\partial _{y}^{2}u^\epsilon=\partial _{z}^{2}u^\epsilon-\partial_x P^\epsilon,\\
\partial _{x}u^\epsilon+\partial _{y}(K^\epsilon u^\varepsilon)+\partial_{z}w^\epsilon=0,\\
 u ^\epsilon \big|_{t=0}= u_{0}(x,y,z) ,\\
(u^\epsilon,K^{^\epsilon}u^\epsilon,w^\epsilon)\big|_{z=0}=0,\quad  \lim\limits_{z\rightarrow+\infty}u^\epsilon=U(t, x, y)
\end{array}
 \label{6.1.4}         \right.\end{equation}
with  regularized Bernoulli's law
\begin{eqnarray}
\partial_t U+ U\partial_x U+ K^\epsilon U\partial_y U-\epsilon^2 \partial_x^2 U -\epsilon^2 \partial_y^2 U +\partial_x P^\epsilon=0.
\label{6.1.44}
\end{eqnarray}
Let the vorticity $\varphi^\epsilon=\partial_{z}u^\epsilon$, then equations $(\ref{6.1.4})$ satisfy the following vorticity system: for any $\epsilon >0$,
\begin{equation}\left\{
\begin{array}{ll}
\partial _{t}\varphi^\epsilon+(u^\epsilon\partial _{x} +K^\epsilon u^\epsilon\partial _{y}+w^\epsilon\partial_{z})\varphi^\epsilon -\varepsilon\partial _{x}^{2}\varphi^\epsilon-\varepsilon\partial _{y}^{2}\varphi^\epsilon= \frac{1}{2}\partial_y K^\epsilon \partial_z |u^\epsilon|^2 +\partial _{z}^{2}\varphi^\epsilon,\\
\partial _{x}u^\epsilon+\partial _{y}(K^\epsilon u^\epsilon)+\partial_{z}w^\epsilon=0,\\
\varphi^\epsilon(0,x,y,z)=\partial_{z}u_{0},\\
\partial _{z}\varphi^\epsilon|_{z=0}=\partial_x P^\epsilon,
\end{array}
 \label{6.1.5}         \right.\end{equation}
where the velocity field $(u^\epsilon, w^\epsilon)$ is given by
\begin{eqnarray}
u^\epsilon (t, x, y, z) = U-\int_{z}^{+\infty} \varphi^\epsilon (t, x, y, \tilde{z})d\tilde{z},
\end{eqnarray}
and
\begin{eqnarray}
w^\epsilon (t, x, y, z) = -\int_{0}^{z} \partial_x u^\epsilon (t, x, y, \tilde{z})d\tilde{z}-\int_{0}^{z} \partial_y (K^\epsilon u^\epsilon (t, x, y, \tilde{z}))d\tilde{z}.
\end{eqnarray}

The three-dimensional Prandtl equation describes the relationship between the velocity and boundary shear stress in the boundary layer and is used to study the behavior and characteristics of fluids within the boundary layer. However, it is challenging to obtain solutions for the boundary layer due to the presence of complex three-dimensional boundary layer structures and secondary flows. Therefore, a specific structure $v(t,x,y,z)=K(t,x,y)u(t,x,y,z)$ is employed to study the three-dimensional Prandtl equation. This structure simplifies the problem by reducing it to a set of equations involving only $u$ and $w$. In Proposition  \ref{p6.1.2}, the Burgers equation $\partial_t K+2 u  (\partial_x+K\partial_y)K=0$ is derived, where the coefficient $K$  is shown to be independent of
$z$. The classical Burgers equation, a nonlinear  partial differential equation, describes the velocity field of a fluid with both viscous and nonlinear transport terms. It finds wide applications in studying shock waves, shock wave interactions, and nonlinear phenomena. The introduction of the Burgers equation simplifies the problem and provides a better mathematical formulation to describe the behavior of the boundary layer. This simplification aids in understanding the characteristics of the boundary layer and offers a potential approach to explore solutions of the three-dimensional Prandtl equation.

The main idea behind this approximation  scheme $(\ref{6.1.222})\rightarrow (\ref{6.1.4})$ or $(\ref{6.1.5})$ is to introduce viscous terms, namely $\epsilon^2\partial_{x}^2u^\epsilon$, $\epsilon^2\partial_{y}^2u^\epsilon$, $\epsilon^2\partial_{x}^2\varphi^\epsilon$ and
$\epsilon^2\partial_{y}^2 \varphi^\epsilon$, in order to prevent the loss of $xy$-derivatives. This regularization offers the advantage that our new weighted
$H^s$ and $L^\infty$ a priori estimates also remain valid for $\varphi^\epsilon$. It is this key aspect that allows us to derive the uniform (in $\epsilon$) estimates in the subsequent section. However, a consequence of this regularization is the emergence of additional terms, such as $\frac{\partial_{xy}\varphi}{\varphi}$, $\frac{\partial_{xy}^2\varphi}{\varphi}$, $\frac{\partial_{xy}\partial_{z}\varphi}{\varphi}$, during the estimation process. Nevertheless, these terms can be controlled within the function space $C([0,T]; H^{s,\gamma}_{\sigma, \delta})$. In addition, it is worth emphasizing that the replacement of Bernoulli's law with the regularized Bernoulli's law is crucial. Without this regularization, it would not be possible to simultaneously satisfy the conditions conditions $u|_{z=0} = 0$ and $\lim_{z  \rightarrow +\infty} u = U$.

Next, we introduce the weighted Sobolev space and define the space $H^{s,\gamma}_{\sigma,\delta}$ by
$$H^{s,\gamma}_{\sigma,\delta}:=\left\{\varphi:\mathbb{T}^{2}\times\mathbb{R^{+}}\rightarrow\mathbb{R},
\|\varphi\|_{H^{s,\gamma}}< +\infty, (1+z)^{\sigma}|\varphi| \geq\delta,
\sum \limits_{|\alpha|\leq 2}  |(1+z)^{\sigma+\alpha_{3}}D^{\alpha}\varphi | ^{2}\leq\frac{1}{\delta^{2}} \right\},$$
where $D^{\alpha}:=\partial_{x}^{\alpha_{1}}\partial_{y}^{\alpha_{2}}\partial_{z}^{\alpha_{3}}$, $\alpha_{1}+\alpha_{2}+\alpha_{3}=s$,  $s\geq 5$, $\gamma\geq 1$, $\sigma>\gamma+\frac{1}{2}$ and $\delta\in (0,1)$. We define the  weighted norm as
  \begin{eqnarray}
  \|\varphi\|_{H^{s,\gamma}}^{2}:=\sum \limits_{|\alpha|\leq s} \|(1+z)^{\gamma+\alpha_{3}}D^{\alpha}\varphi\|_{L^{2}}^{2}
  \label{norm1}
   \end{eqnarray}
and
\begin{eqnarray}
\|\varphi\|_{H^{s,\gamma}_{g}}^{2}:=\|(1+z)^{\gamma}g_{s}\|^{2}_{L^{2}}+
\sum\limits_{\substack{ |\alpha|\leq s\\ \alpha_{1}+\alpha_{2}\leq s-1}}  \|(1+z)^{\gamma+\alpha_{3}}D^{\alpha}\varphi\|_{L^{2}}^{2}.
  \label{norm2}
   \end{eqnarray}
Here
$$g_{s}:=\partial_{xy}^{s}\varphi-\frac{\partial_{z}\varphi}{\varphi}\partial_{xy}^{s}(u-U)
=\sum\limits_{i=0}^{ s}\left( \partial_{x}^{i}\partial_{y}^{s-i}\varphi- \frac{\partial_{z}\varphi}{\varphi}  \partial_{x}^{i}\partial_{y}^{s-i} (u-U) \right),$$
and
$$g_{s}=(g_{s})_{x}+(g_{s})_{y}=\partial_x\partial_{xy}^{s-1}\varphi-\frac{\partial_{z}\varphi}{\varphi}\partial_x\partial_{xy}^{s-1}(u-U)+\partial_y\partial_{xy}^{s-1}\varphi
-\frac{\partial_{z}\varphi}{\varphi}\partial_y\partial_{xy}^{s-1}(u-U),$$
provided that $\varphi = \partial_z u >0$.
Hence, we can calculate that
$$g_{j+1}=\partial_{x}g_{j}+\partial_{y}g_{j}+\partial_{xy}^{j}u\partial_{x}a+\partial_{xy}^{j}u\partial_{y}a,\quad
a=\frac{\partial_{z}\varphi}{\varphi}.$$

As pointed out that, the difference between norm (\ref{norm1}) and norm (\ref{norm2}) are that the weighted $L^2$-norm of $\partial_{xy}\varphi$ is replaced by that of $g_s$, which avoids the loss of $xy$-derivative by the nonlinear cancellation. On the other hand, these two weighted norms are the almost equivalent, that is, there exists a positive constant
$C$ such that the following inequality holds for any $\varphi \in H^{s,\gamma}_{\sigma,\delta}(\mathbb{T}^{2}\times\mathbb{R}^{+})$, $\gamma\geq 1$, $\sigma>\gamma+\frac{1}{2}$ and $\delta\in (0,1)$:
 \begin{eqnarray}
C_{\delta }\|\varphi\|_{H^{s,\gamma}_{g}}\leq  \|\varphi\|_{H^{s,\gamma} }+ \|u-U\|_{H^{s,\gamma-1} }
\leq C_{s, \gamma,\sigma,\delta  }\left(\|\varphi\|_{H^{s,\gamma}_{g}}+\|\partial_{xy}^s U\|_{L^2(\mathbb{T}^2)}\right)
  ,
\end{eqnarray}

In this chapter, for convenience, we simply write
$$\partial_{xy}^{s}:= \sum\limits_{i=0}^{ s} \partial_{x}^{i}\partial_{y}^{s-i} , \quad  \iiint\cdot:=\iint_{\mathbb{T}^{2}}\int_{\mathbb{R^{+}}}\cdot dxdydz.$$

Firstly, we state our main result as follows.
\begin{Theorem}\label{t6.1.1}
Given any   integer  $s \geq 5$, and real numbers $\gamma, \sigma, \delta$ satisfying $\gamma\geq 1$, $\sigma>\gamma+\frac{1}{2}$ and $\delta\in (0,1)$.   Assume the following conditions on the initial data,  the outer flow $U$ and $K(x,y)$:

(i)~Suppose that the initial data $u_0-U(0,x) \in H^{s, \gamma-1}$ and $\partial_z u_0 \in H^{s,\gamma}_{\sigma, 2\delta}$ satisfy the compatibility conditions $u_0|_{z=0}$ and $\lim \limits_{z \rightarrow +\infty} u_{0}=U|_{t=0}$. In addition, when $s = 5$, it is further assumed that $\delta \geq 0$ is chosen small enough such that
 \begin{eqnarray}
 \|\omega_0\|_{H^{s,\gamma} } \leq C \delta^{-1}
 \label{6.1.888}
 \end{eqnarray}
with a generic constant $C$.

(ii)~The outer flow $U$ is supposed to satisfy
 \begin{eqnarray}
\sup \limits_{t} \sum\limits_{l=0}^{\frac{s}{2}+1}\|\partial_t^l U\|_{H^{s-2l+2}(\mathbb{T}^2)} < + \infty.
\label{6.1.8}
\end{eqnarray}

(iii)~The $K(x,y)$ is supposed to satisfy
 \begin{eqnarray}
\|\partial^{s}_{xy} K\|_{L^{\infty}(\mathbb{T}^2)} < + \infty.
\label{6.1.88}
\end{eqnarray}

Then there exist a time $T := T(K, s, \gamma, \sigma, \delta,  \|\varphi_0\|_{H^{s,\gamma} },U)$ such that the initial boundary value problem
(\ref{6.1.1})-(\ref{6.1.3}) has a unique classical solution $(u, v, w)$ (v=Ku) satisfying
$$u-U \in L^{\infty}([0,T];H^{s,\gamma-1}) \cap C([0,T]; H^{s}-\varphi)$$
and
$$ \partial_z u \in L^{\infty}([0,T];H^{s,\gamma}_{\sigma, \delta}) \cap C([0,T]; H^{s}-\varphi),$$
where $H^s-\varphi$ is the space $H^s$ endowed with its weak topology.
\end{Theorem}

As we have known,   Crocco transformation  and the outer flow  $(U,V)=\text{constant}$  have been used in three-dimensional prandtl system, such as \cite{LMD,Lx,lin,lwy2,lwy1,lwy,PX,px}.
But, we need neither the Crocco transformation nor the the outer flow  $(U,V)=\text{constant}$ in this chapter. In fact, we  first use the energy method to prove the local existence of three-dimensional Prandtl system based on a cancellation property. Our result is an improvement of the   existence and  uniqueness of solutions   in Masmoudi and Wong \cite{mw}  from the two dimensional case to the three dimensional with special structure case.

The chapter is arranged as follows. Section 6.1 is an introduction, which contains many inequalities used in this paper. In Section 6.2, we shall give the uniform estimates with the weighted norm. The difficulty of this paper is the estimates on the boundary at $z=0$, and we give the regular pattern at $z=0$ for $s\geq 5$ eventually through a series of precise calculations, see Lemmas \ref{y6.2.1}-\ref{y6.2.2}. In Section 6.3, we shall prove local existence and uniqueness of solutions to the Prandtl system to problem (\ref{6.1.1}).

\subsection{Preliminaries}
\begin{Proposition} \label{p6.1.2}
Under the assumption $\partial_z u >0$ and special structures $v(t,x,y,z)=K(t, x, y)u(t,x,y,z)$, to study
the problem (\ref{6.1.1}) is equivalent to studying the following reduced problem:
\begin{equation}\left\{
\begin{array}{ll}
\partial _{t}u+(u\partial _{x} +K u\partial _{y} +w\partial_{z})u =\partial _{z}^{2}u^\epsilon-\partial_x P,\\
\partial _{x}u+\partial _{y}(K u)+\partial_{z}w=0,\\
 u  \big|_{t=0}= u_{0}(x,y,z) ,\\
(u,w)\big|_{z=0}=0,\quad  \lim\limits_{z\rightarrow+\infty}u=U(t, x, y),
\end{array}
 \label{6.4.40}         \right.\end{equation}
 where  $K$ is independent of $t$, and $U$ satisfies
 \begin{equation}\left\{
\begin{array}{ll}
\partial_t U+ U\partial_x U+ K U\partial_y U+\partial_x P=0,\\
\partial_y P-K\partial_x P=0.
\end{array}
        \right.\end{equation}
\end{Proposition}
\textbf{Proof}.
Using special structures $v(k,x,y,z)=K(t, x, y)u(k,x,y,z)$ for the equation $(\ref{6.1.1})_2$, we have
 \begin{eqnarray}
\partial_t (Ku)+(u \partial_x +Ku \partial_y +w \partial_z)(Ku)+\partial_y P - K\partial_z^2 u=0,
\label{6.4.41}
\end{eqnarray}
which, by the equation $(\ref{6.1.1})_1$, yields
 \begin{eqnarray}
u\left[\partial_t K+u( \partial_x + \partial_y)K\right]-K\partial_x P + \partial_y P=0.
\label{6.4.42}
\end{eqnarray}
Noting that $K(t,x,y)$ is independent of $z$, by differentiating (\ref{6.4.42}) with respect to $z$, it follows that
 \begin{eqnarray}
\partial_z u\partial_t K+2 u \partial_z u (\partial_x+K\partial_y)K=0.
\label{6.4.43}
\end{eqnarray}
Due to the monotonicity assumption, we can further get
 \begin{eqnarray}
\partial_t K+2 u  (\partial_x+K\partial_y)K=0.
\label{6.4.44}
\end{eqnarray}
Differentiating (\ref{6.4.44}) with respect to $z$ gives the $Burgers$ equation
 \begin{eqnarray}
  (\partial_x+K\partial_y)K=0.
\label{6.4.45}
\end{eqnarray}
Plugging (\ref{6.4.45}) into (\ref{6.4.44}), it follows
 \begin{eqnarray}
 \partial_t K=0,
\label{6.4.46}
\end{eqnarray}
which, combined (\ref{6.4.42}) with (\ref{6.4.45}), gives
 \begin{eqnarray}
-K\partial_x P + \partial_y P=0.
\label{6.4.47}
\end{eqnarray}
The above equation means that $(\partial_x P, \partial_y P)$ is parallel to both the velocity field of the outer Euler flow and the tangential velocity field in the boundary layer. Now, we assume that  the classical solution $(u, v, w)$  to the problem (\ref{6.1.1}) satisfies $v(t,x,y,z)=K(x,y)u(t,x,y,z)$, then  by using (\ref{6.4.45}) and (\ref{6.4.47}), $W(t, x, y, z):\equiv v(t,x,y,z)-K(x,y)u(t,x,y,z)$ satisfies the following problem:
\begin{equation}\left\{
\begin{array}{ll}
\partial _{t}W+(u\partial _{x} +v\partial _{y} +w\partial_{z})W+(\partial_y v - K\partial_y u)W=\partial _{z}^{2}W,\\
W|_{t=0}=0,\\
W|_{z=0}=0,\\
  \lim\limits_{z\rightarrow+\infty}W=0.
\end{array}
 \label{6.4.48}         \right.\end{equation}
The equation (\ref{6.4.48}) has only trivial solution $W\equiv0$ based on the energy argument. Therefore,   studying
the problem (\ref{6.1.1}) is equivalent to studying  the following reduced problem for only two
unknown functions $u$ and $w$,
\begin{equation}\left\{
\begin{array}{ll}
\partial _{t}u+(u\partial _{x} +K u\partial _{y} +w\partial_{z})u =\partial _{z}^{2}u^\epsilon-\partial_x P,\\
\partial _{x}u+\partial _{y}(K u)+\partial_{z}w=0,\\
 u  \big|_{t=0}= u_{0}(x,y,z) ,\\
(u,w)\big|_{z=0}=0,\quad  \lim\limits_{z\rightarrow+\infty}u=U(t, x, y) .
\end{array}
 \label{6.4.49}         \right.\end{equation}
Hence, the proof of   Proposition \ref{p6.1.2} is complete.
\hfill $\Box$

\begin{Proposition}
Let $s\geq 5$ be an integer, $\gamma\geq 1$, $\sigma>\gamma+\frac{1}{2}$ and $\delta\in (0,1)$. Then for any $\varphi \in H^{s,\gamma}_{\sigma,\delta}(\mathbb{T}^{2}\times\mathbb{R}^{+})$,
we have the following inequality
 \begin{eqnarray}
C_{\delta }\|\varphi\|_{H^{s,\gamma}_{g}}\leq  \|\varphi\|_{H^{s,\gamma} }+ \|u-U\|_{H^{s,\gamma-1} }
\leq C_{s, \gamma,\sigma,\delta  }\left(\|\varphi\|_{H^{s,\gamma}_{g}}+\|\partial_{xy}^s U\|_{L^2(\mathbb{T}^2)}\right)
  ,
\label{6.4.1}
\end{eqnarray}
where $C_{s, \gamma,\sigma,\delta  }>0$ is a constant and only depends on $s, \gamma,\sigma,\delta  $.
\end{Proposition}
\textbf{Proof}.
Firstly, by the definition of $\|\varphi\|_{H^{s,\gamma}}$ and $\|u-U\|_{H^{s,\gamma-1}}$ that
\begin{align}
    \|\varphi\|_{H^{s,\gamma}} + \sum^s_{k=0} \|(1+z)^{\gamma-1} \partial_{xy}^k (u-U)\|_{L^2}
   & \leq \|\varphi\|_{H^{s,\gamma}} + \|u-U\|_{H^{s,\gamma-1}} \nonumber\\
    &\leq 2 \left(\|\varphi\|_{H^{s,\gamma}} + \sum^s_{k=0} \|(1+z)^{\gamma-1}
    \partial_{xy}^k (u-U) \|_{L^2}\right).
    \label{6.4.0001}
\end{align}
Furthermore, applying Wirtinger's inequality  and Lemma \ref{y6.4.1}, we get
\[\left\{\begin{aligned}
    \sum^s_{k=1} \|(1+z)^{\gamma-1} \partial_{xy}^k (u-U)\|_{L^2} & \le \frac{1+\pi^{-2s}}{1-\pi^{-2}}
    \| (1+z)^{\gamma-1} \partial_{xy}^s (u-U)\|_{L^2} \\
    \|(1+z)^{\gamma-1} (u-U) \|_{L^2} & \le \frac{2}{2\gamma-1} \|(1+z)^\gamma \varphi\|_{L^2} \le \frac{2}{2\gamma-1} \|\varphi\|_{H^{s,\gamma}}, \label{6.4.0002}
\end{aligned}\right.\]
and hence, there exists a constant $C_{s,\gamma} > 0$ such that
\begin{align}
\|\varphi\|_{H^{s,\gamma}} +  \|(1+z)^{\gamma-1} \partial_{xy}^s (u-U)\|_{L^2}
    &\leq \|\varphi\|_{H^{s,\gamma}} + \sum^s_{k=0} \|(1+z)^{\gamma-1} \partial_{xy}^k
    (u-U)\|_{L^2}\nonumber\\
    &\leq C_{s,\gamma} \left(\|\varphi\|_{H^{s,\gamma}} +
    \|(1+z)^{\gamma-1} \partial_{xy}^s (u-U) \|_{L^2}\right). \label{6.4.0003}
\end{align}
Therefore, according to Lemma \ref{y6.4.2} and inequalities (\ref{6.4.0001})-(\ref{6.4.0003}), it suffices
to prove
\begin{align}
    C_\delta \|\varphi\|_{H^{s,\gamma}_g} \le \|\varphi\|_{H^{s,\gamma}} + \|(1+z)^{\gamma-1} \partial_{xy}^s
    (u-U)\|_{L^2}
    \leq C_{s,\gamma,\sigma,\delta} \left(\|\varphi\|_{H^{s,\gamma}_g} + \|\partial_{xy}^s U\|_{L^2} \right)
\end{align}
for some constants $C_\delta$ and $C_{s,\gamma, \sigma, \delta} > 0$.
\hfill $\Box$

\begin{Lemma} \label{y6.4.1}
Let $f:\mathbb{T}^{2} \times\mathbb{R}^{+}\rightarrow \mathbb{R}$, \\
 (i)  if $\lambda>-\frac{1}{2}$ and $\lim\limits_{z\rightarrow+\infty}f(x,y ,z)=0$, then
  \begin{eqnarray}
  \|(1+z)^{\lambda}f\|_{L^{2}(\mathbb{T}^{2} \times\mathbb{R}_{+})} \leq \frac{2}{2\lambda+1} \|(1+z)^{\lambda+1}\partial_{z}f\|_{L^{2}(\mathbb{T}^{2} \times\mathbb{R}_{+})};
 \label{6.4.10}
 \end{eqnarray}
 (ii)  if $\lambda<-\frac{1}{2}$ and $ f(x,y,z )|_{z=0}=0$, then
 \begin{eqnarray}
 \|(1+z)^{\lambda}f\|_{L^{2}(\mathbb{T}^{2} \times\mathbb{R}^{+})} \leq
 -\frac{2}{2\lambda+1} \|(1+z)^{\lambda+1}\partial_{z}f\|_{L^{2}(\mathbb{T}^{2} \times\mathbb{R}^{+})}.
\label{6.4.11}
\end{eqnarray}
\end{Lemma}
\textbf{Proof}.
The proof is elementary. Noting that the term
\begin{eqnarray*}
&&\|(1+z)^{\lambda}f\|_{L^{2}(\mathbb{T}^2 \times\mathbb{R}_{+})} ^{2}
=\frac{1}{2\lambda+1}\int_{\mathbb{T}^2}\int_{\mathbb{R}_{+}}f^{2}( x,y,z)d(1+z)^{2\lambda+1} dxdy  \nonumber \\
&&=\frac{1}{2\lambda+1}\int_{\mathbb{T}^2}(1+z)^{2\lambda+1}f^{2}( x,y,z) \big|_{z=0}^{z=+\infty}dxdy
-\frac{2}{2\lambda+1} \int_{\mathbb{T}^2}\int_{\mathbb{R}_{+}} (1+z)^{2\lambda+1}f ( x,y,z) \partial_{z}f(x,y,z)dxdydz \nonumber \\
&& \leq \left|-\frac{2}{2\lambda+1}\right| \|(1+z)^{\lambda }  f\|_{L^{2}(\mathbb{T}^2 \times\mathbb{R}_{+})}  \|(1+z)^{\lambda+1}\partial_{z}f\|_{L^{2}(\mathbb{T}^2 \times\mathbb{R}_{+})},
 \end{eqnarray*}
thus we prove  $(\ref{6.4.10})$-$(\ref{6.4.11})$.
\hfill $\Box$

\begin{Lemma} \label{y6.4.2}
Let $s\geq 5$ be an integer, $\gamma\geq 1$, $\sigma>\gamma+\frac{1}{2}$ and $\delta\in (0,1)$. If $\varphi \in H^{s,\gamma}_{\sigma,\delta}(\mathbb{T}^{2}\times\mathbb{R}^{+})$, then  for any $k=0,1,2,\cdots, s$,
 \begin{eqnarray}
 \|(1+z)^{\gamma} g_{k}\|_{L^{2}}\leq\|(1+z)^{\gamma}\partial _{xy}^{k}\varphi\|_{L^{2}}+\delta^{-2}\|(1+z)^{\gamma-1}\partial _{xy}^{k}(u-U)\|_{L^{2}}  ,
 \label{6.4.5}
\end{eqnarray}
where $g_{k}:=\partial_{xy}^{k}\varphi-\frac{\partial_{z}\varphi}{\varphi}\partial_{xy}^{k}(u-U)$. In addition, if $u \big|_{z=0} = 0$, then for any $k=1, 2, \ldots, s$,
 \begin{eqnarray}
 \|(1+z)^{\gamma}\partial _{xy}^{k}\varphi\|_{L^{2}}+\|(1+z)^{\gamma-1}\partial _{xy}^{k}(u-U)\|_{L^{2}}
\leq C_{s, \gamma,\sigma,\delta  }\left(\|(1+z)^{\gamma} g_k\|_{L^2}+\|\partial_{xy}^k U\|_{L^2(\mathbb{T}^2)}\right) , \label{6.4.4}
\end{eqnarray}
where $C_{s, \gamma,\sigma,\delta  }>0$ is a constant and only depends on $s, \gamma,\sigma,\delta  $.
\end{Lemma}

\textbf{Proof}.
The proof is similar to that given in \cite{mw}, it is only related to the defined weighted Sobolev space and norms, not to the condition $(\ref{6.1.4})_2$, so we outline it here for brevity.\\
Firstly, according to the definition  of $H^{s,\gamma}_{\sigma,\delta}$, we have
\begin{eqnarray}
\delta(1+z)^{-\sigma} \leq |\varphi|\leq \delta^{-1}(1+z)^{-\sigma}, \ \ |\partial_{z}\varphi| \leq\delta^{-1}(1+z)^{-\sigma-1},   \label{6.4.6}
\end{eqnarray}
so $\frac{|\partial_z \varphi|}{\varphi} \leq \delta^{-2}(1+z)^{-1}$, for any $k=0,1,2,\cdots, s$,
 \begin{align*}
 \|(1+z)^{\gamma} g_{k}\|_{L^{2}}&\leq \|(1+z)^{\gamma }\partial _{xy}^{k}\varphi\|_{L^{2}}+ \|(1+z)^{\gamma }\frac{\partial _{z}\varphi}{\varphi}\partial _{xy}^{k}(u-U)\|_{L^{2}}    \nonumber\\
 &\leq \|(1+z)^{\gamma }\partial _{xy}^{k}\varphi\|_{L^{2}}+ \|(1+z)\frac{\partial _{z}\varphi}{\varphi}\|_{L^{\infty}}\|(1+z)^{\gamma-1}\partial _{xy}^{k}(u-U)\|_{L^{2}}\nonumber\\
 &\leq \|(1+z)^{\gamma}\partial _{xy}^{k}\varphi\|_{L^{2}}+\delta^{-2}\|(1+z)^{\gamma-1}\partial _{xy}^{k}(u-U)\|_{L^{2}}.
 \end{align*}
Next,
using $(\ref{6.4.11})$ in Lemma \ref{y6.4.1} with the index $\gamma-\sigma-1<-\frac{3}{2}<-\frac{1}{2}$, we get
\begin{align}
\|(1+z)^{\gamma-1}\partial _{xy}^{k}(u-U)\|_{L^{2}}&\leq \delta^{-1}\|(1+z)^{\gamma-\sigma-1}\frac{\partial _{xy}^{k}(u-U)}{\varphi}\|_{L^{2}}    \nonumber\\
&\leq C_{\gamma,\sigma} \delta^{-1} \left(\left\|\frac{\partial _{xy}^{k}(u-U)}{\varphi}|_{z=0}\right\|_{L^{2}}+\left\|(1+z)^{\gamma-\sigma}\partial_{z}\left(\frac{\partial _{xy}^{k}(u-U)}{\varphi}\right)\right\|_{L^{2}}\right) \nonumber\\
&\leq C_{\gamma,\sigma} \delta^{-2}\left(\|\partial_{xy}^{k}U\|_{{L^2}(\mathbb{T}^2)}+\|(1+z)^{\gamma }g_{k}\|_{L^{2}}\right).
\label{6.4.8}
\end{align}
Finally, for $k=0,1,2,\cdots, s$, according to $(\ref{6.4.8})$ and the triangle inequality, we have
\begin{align}
\|(1+z)^{\gamma}\partial _{xy}^{k}\varphi\|_{L^{2}}&\leq
\|(1+z)^{\gamma }g_{k}\|_{L^{2}}+ \|(1+z)^{\gamma }\frac{\partial _{z}\varphi}{\varphi}\partial _{xy}^{k}(u-U)\|_{L^{2}}    \nonumber\\
&\leq C_{\gamma,\sigma,\delta} \left(\|\partial_{xy}^{k}U\|_{{L^2}(\mathbb{T}^2)}+\|(1+z)^{\gamma }g_{k}\|_{L^{2}}\right).
\label{6.4.9}
\end{align}
The proof is thus complete.
\hfill $\Box$

\begin{Lemma}(\cite{ghy}) \label{y6.4.3}
Let $f:\mathbb{T}^{2} \times\mathbb{R}^{+}\rightarrow \mathbb{R}$ and $ f(x,y,z )|_{z=0}=0$. Then there exists a constant $C>0$ such that
\begin{eqnarray}
\begin{aligned}
&  \|f(x,y,z )\|_{L^{\infty}(\mathbb{T}^{2} \times \mathbb{R}^{+})}  \\
&   \leq C
 (\|f(x,y,z)\|_{L^{2}(\mathbb{T}^{2}\times \mathbb{R}^{+} )}   +\|\partial_{z} f(x,y,z)\|_{L^{2}(\mathbb{T}^{2}\times \mathbb{R}^{+} )}
  +\partial_{xy}^{1}\|f(x,y,z)\|_{L^{2}(\mathbb{T}^{2}\times \mathbb{R}^{+} )}\\
 &\quad +\|\partial_{z}\partial_{xy}^{1}f(x,y,z)\|_{L^{2}(\mathbb{T}^{2}\times \mathbb{R}^{+} )}
  +\|\partial_{xy}^{2}f(x,y,z)\|_{L^{2}(\mathbb{T}^{2}\times \mathbb{R}^{+} )}+\|\partial_{z}\partial_{xy}^{2}f(x,y,z)\|_{L^{2}(\mathbb{T}^{2}\times \mathbb{R}^{+} )}  ).
\label{6.4.12}
\end{aligned}
\end{eqnarray}

\end{Lemma}
\textbf{Proof}.
Firstly, we have
\begin{eqnarray}
&& f^{2}(x,y,z)=\int_{0}^{z}\partial_{\tau}[f^{2}(x,y,\tau)]d\tau =2\int_{0}^{z}f (x,y,\tau)\partial_{\tau} f (x,y,\tau) d\tau   \nonumber \\
&& \leq C\left(\int_{0}^{z} f^{2}(x,y,\tau) d\tau \right)^{\frac{1}{2}} \left(\int_{0}^{z}| \partial_{\tau}f (x,y,\tau)|^{2} d\tau \right)^{\frac{1}{2}} \nonumber \\
&&  \leq C\left(\int_{0}^{+\infty}| f(x,y,z)|^{2} dz+\int_{0}^{+\infty} |\partial_{z}f (x,y,z)|^{2} dz \right).
\label{6.4.29}
\end{eqnarray}
Secondly, we apply the two-dimensional Sobolev inequality for $x, y$ variables to get
\begin{eqnarray}
  |f(x,y,z)|^{2} \leq C(\|f( z)\|_{L^{2}(\mathbb{T}^{2})}^{2}+ \|\partial_{xy}^{1}f( z)\|_{L^{2}(\mathbb{T}^{2})}^{2}
+\|\partial_{xy}^{2}f( z)\|_{L^{2}(\mathbb{T}^{2})}^{2} )
\label{6.4.30}
\end{eqnarray}
and
\begin{eqnarray}
 |\partial_{z}f (x,y,z)|^{2} \leq C(\|\partial_{z}f( z)\|_{L^{2}(\mathbb{T}^{2})}^{2}+ \|\partial_{z}\partial_{xy}^{1}f( z)\|_{L^{2}(\mathbb{T}^{2})}^{2}
+\|\partial_{z}\partial_{xy}^{2}f( z)\|_{L^{2}(\mathbb{T}^{2})}^{2} ).
\label{6.4.31}
\end{eqnarray}
At last, substituting (\ref{6.4.30}) and (\ref{6.4.31}) into (\ref{6.4.29}),  we complete the proof of estimate  (\ref{6.4.12}).

\hfill $\Box$

\begin{Lemma} \label{y6.4.4}
Let the vector field $(u, Ku, w)$ defined on $\mathbb{T}^2 \times \mathbb{R}^+$ satisfy the condition $\partial_x u+ \partial_y (Ku)+ \partial_z w=0$.  If $s\geq 5$ is an integer, $\gamma\geq 1$, $\sigma>\gamma+\frac{1}{2}$ and $\delta\in (0,1)$, then for any $\varphi \in H^{s,\gamma}_{\sigma,\delta}(\mathbb{T}^{2}\times\mathbb{R}^{+})$,
we have the following inequalities: :\\
$(i)$ For $k=0, 1, 2,  \cdots, s $,
 \begin{eqnarray}
\|(1+z)^{\gamma-1} \partial_{xy}^k(u-U)\|_{L^2}
\leq C_{s, \gamma,\sigma,\delta  }\left(\|\varphi\|_{H^{s,\gamma}_{g}}+\|\partial_{xy}^s U\|_{L^2(\mathbb{T}^2)}\right).
\label{6.4.13}
\end{eqnarray}
(ii) For $k=0, 1, 2,  \cdots, s-1 $, if $K \in C^{s+1}(\mathbb{T}^2)$, then
\begin{eqnarray}
\left \|\frac{\partial_{xy}^{k} w +z\tilde{U}}{1+z}  \right\|_{L^{2}}
 \leq C_{s, \gamma,\sigma,\delta  }\left(\|\varphi\|_{H^{s,\gamma}_{g}}+\|\partial_{xy}^s U\|_{L^2(\mathbb{T}^2)}\right).
\label{6.4.14}
\end{eqnarray}
where $\tilde{U}= \partial_x\partial^{k}_{xy} U + \partial_y\partial^{k}_{xy} (KU)$\\
(iii) For all $| \alpha  | \leq s$,
\begin{align}
&\|(1+z)^{\gamma+\alpha_3} D^\alpha \varphi\|_{L^2} \leq \left\{
\begin{aligned}
&C_{s, \gamma,\sigma,\delta  }\left(\|\varphi\|_{H^{s,\gamma}_{g}}+\|\partial_{xy}^s U\|_{L^2(\mathbb{T}^2)}\right) \ \ \mathrm{if}~\alpha =(\alpha_1, \alpha_2, 0) ,\\
&\|\varphi\|_{H^{s,\gamma}_{g}} \ \ \mathrm{if}~\alpha \neq(\alpha_1, \alpha_2, 0).\\
\end{aligned} \label{6.4.15}
\right.
\end{align}
(iv) For $k=0, 1, 2,  \cdots, s $,
\begin{align}
&\|(1+z)^{\gamma} g_k\|_{{L^2}(\mathbb{T}^2)} \leq \left\{
\begin{aligned}
&C_{s, \gamma,\sigma,\delta  }\left(\|\varphi\|_{H^{s,\gamma}_{g}}+\|\partial_{xy}^s U\|_{L^2(\mathbb{T}^2)}\right) \ \ \mathrm{if}~k \neq s ,\\
&\|\varphi\|_{H^{s,\gamma}_{g}} \ \ \mathrm{if}~k=s.\\
\end{aligned} \label{6.4.16}
\right.
\end{align}
(v) For $k=0, 1, 2,  \cdots, s-2$,
\begin{eqnarray}
 \|  \partial_{xy}^{k} u  \|_{L^{\infty}} \leq C_{s, \gamma,\sigma,\delta  }\left(\|\varphi\|_{H^{s,\gamma}_{g}}+\|\partial_{xy}^s U\|_{L^2(\mathbb{T}^2)}\right).
\label{6.4.17}
\end{eqnarray}
(vi) For $k=0, 1, 2,  \cdots, s-3$, if $K \in C^{s+1}(\mathbb{T}^2)$, then
\begin{eqnarray}
 \left\|  \frac{\partial_{xy}^{k} w}{1+z} \right \|_{L^{\infty}} \leq C_{s, \gamma,\sigma,\delta  }\left(\|\varphi\|_{H^{s,\gamma}_{g}}+\|\partial_{xy}^s U\|_{L^2(\mathbb{T}^2)}+1\right).
\label{6.4.18}
\end{eqnarray}
(vii) For all $| \alpha  | \leq s-3$,
\begin{eqnarray}
\|(1+z)^{\gamma+\alpha_3} D^\alpha \varphi\|_{L^\infty} \leq C_{s, \gamma,\sigma,\delta  }\|\varphi\|_{H^{s,\gamma}_{g}}
\label{6.4.19}
\end{eqnarray}
where $C_{s, \gamma,\sigma,\delta  }>0$ is a constant and only depends on $s, \gamma,\sigma,\delta  $.
\end{Lemma}
\textbf{Proof}.
First, by the definition of $\| \cdot \|_{H^{s, \gamma-1}}$ and Lemma \ref{y6.4.2},  we can directly derive the inequality (\ref{6.4.13}). Using Lemma \ref{y6.4.1}, (\ref{6.4.13}) and the fact that $\partial_x u+ \partial_y (Ku)+ \partial_z w=0$, we have
\begin{align}
 \left \|\frac{\partial_{xy}^{k} w +z \tilde{U}}{1+z}  \right\|_{L^{2}}
 &\leq 2 \| \partial_{xy}^{k }(-\partial_x u- \partial_y (Ku))+\tilde{U}\|_{L^{2}  }  \nonumber \\
&\leq 2 \| \partial_{xy}^{k +1}(u-U)\|_{L^{2}  } +C \sum \limits^{k} \limits_{i=0} \|\partial_y\{\partial_{xy}^{i} K \partial_{xy}^{k-i}(u-U)\}\|_{L^2} \nonumber \\
&  \leq C_{s, \gamma,\sigma,\delta  }\left(\|\varphi\|_{H^{s,\gamma}_{g}}+\|\partial_{xy}^s U\|_{L^2(\mathbb{T}^2)}\right).
\end{align}
The proofs of  (\ref{6.4.15}) and (\ref{6.4.16}) are obvious, so we omit it here for brevity. Second, for $k=0, 1, 2,  \cdots, s-2$, applying Lemma \ref{y6.4.3} and (\ref{6.4.13}), and (\ref{6.4.15}), we have
\begin{align}
 \|  \partial_{xy}^{k} (u-U)  \|_{L^{\infty}} &\leq C( \|  \partial_{xy}^{k} (u-U)  \|_{L^{2}}+\|  \partial_{xy}^{k+1} (u-U) \|_{L^{2}}+\|  \partial_{xy}^{k+2 } (u-U) \|_{L^{2}}+\|  \partial_{xy}^{k} \varphi \|_{L^{2}}\nonumber\\
&\quad+\|   \partial_{xy}^{k+1}\varphi\|_{L^{2}}+\|  \partial_{xy}^{k+2 } \varphi \|_{L^{2}}) \nonumber  \\
&\leq C_{s, \gamma,\sigma,\delta  }\left(\|\varphi\|_{H^{s,\gamma}_{g}}+\|\partial_{xy}^s U\|_{L^2(\mathbb{T}^2)}\right).
\end{align}
Hence, by the triangle inequality and $\|\partial_{xy}^k U\|_{L^2(\mathbb{T}^2)} \leq \|\partial_{xy}^s U\|_{L^2(\mathbb{T}^2)}$, we can justify (\ref{6.4.17}).\\
Next, for $k=0, 1, 2,  \cdots, s-3$, using   (\ref{6.4.13})-(\ref{6.4.15}) and Lemma \ref{y6.4.3}, we conclude that
\begin{align}
 \left\|  \frac{\partial_{xy}^{k} w }{1+z}  \right\|_{L^{\infty}}& \leq \left\|  \frac{z\tilde{U} }{1+z}  \right\|_{L^{\infty}}+\left\|  \frac{\partial_{xy}^{k} w+z\tilde{U} }{1+z}  \right\|_{L^{\infty}} \nonumber \\
 &  \leq C \Big(\|   \tilde{U}   \|_{L^{2}(\mathbb{T}^2)}+\|  \partial_{xy} \tilde{U}   \|_{L^{2}(\mathbb{T}^2)}+\|  \partial_{xy}^{2} \tilde{U}   \|_{L^{2}(\mathbb{T}^2)}+\|(1+z)^{-1}(\partial_{xy}^{k} w+z\tilde{U})\|_{L^2} \nonumber \\
&\quad+\|(1+z)^{-1}(\partial_{xy}^{k} \partial_{z}w+\tilde{U})\|_{L^2}+\|(1+z)^{-2}(\partial_{xy}^{k} w+z\tilde{U})\|_{L^2} \nonumber  \\
&\quad+\|(1+z)^{-1}(\partial_{xy}^{k+1} w+z\partial_{xy}\tilde{U})\|_{L^2} +\|(1+z)^{-1}(\partial_{xy}^{k+2} w+z\partial_{xy}^2\tilde{U})\|_{L^2} \nonumber  \\
&\quad+\|(1+z)^{-1}(\partial_{xy}^{k+1} \partial_{z}w+\partial_{xy}\tilde{U})\|_{L^2}+\|(1+z)^{-2}(\partial_{xy}^{k+1} w+z\partial_{xy}\tilde{U})\|_{L^2} \nonumber  \\
&\quad+\|(1+z)^{-1}(\partial_{xy}^{k+2} \partial_{z}w+\partial_{xy}^2\tilde{U})\|_{L^2}+\|(1+z)^{-2}(\partial_{xy}^{k+2} w+z\partial_{xy}^2\tilde{U})\|_{L^2} \Big )
 \nonumber \\
& \leq C_{s, \gamma,\sigma,\delta  } \left(\|\varphi\|_{H^{s,\gamma}_{g}}+\|\partial_{xy}^s U\|_{L^2(\mathbb{T}^2)}+1\right).
\label{6.4.20}
\end{align}
Last, inequality (\ref{6.4.19}) follows directly from  inequality (\ref{6.4.15}) in Lemma \ref{y6.4.3}.
\hfill $\Box$

\begin{Lemma}(\cite{[01]} p.369 Theorem 9)\label{pmax}
Denote the operator
\begin{eqnarray*}
\mathcal{L}u=-\sum^{n}_{i,j=1}a^{ij}u_{x_i x_j}+\sum^{n}_{i=1}b^{i}u_{x_i }+cu,
\end{eqnarray*}
where $a^{ij}$, $b^{i}$ and  $c$  are continuous.  Assume $u \in C_{1}^{2}([0,T]\times\mathbb{T}^2\times \mathbb{R}^{+}) \cap C([0,T]\times\overline{\mathbb{T}^2\times \mathbb{R}^{+}}) $, and $c\geq 0$ in $[0,T]\times\mathbb{T}^2\times \mathbb{R}^{+}$.

(i)If
\begin{eqnarray*}
u_{t}+\mathcal{L}u \leq 0
\end{eqnarray*}
in $[0,T]\times\mathbb{T}^2\times \mathbb{R}^{+}$, then
\begin{eqnarray*}
\max \limits_{[0,T]\times\overline{\mathbb{T}^2\times \mathbb{R}^{+}}} \leq \max \limits_{\partial([0,T]\times\mathbb{T}^2\times \mathbb{R}^{+})}u^{+}.
\end{eqnarray*}

(ii)If
\begin{eqnarray*}
u_{t}+\mathcal{L}u \geq 0
\end{eqnarray*}
in $[0,T]\times\mathbb{T}^2\times \mathbb{R}^{+}$, then
\begin{eqnarray*}
\min\limits_{[0,T]\times\overline{\mathbb{T}^2\times \mathbb{R}^{+}}} \geq -\max \limits_{\partial([0,T]\times\mathbb{T}^2\times \mathbb{R}^{+})}u^{-},
\end{eqnarray*}
where $u^{+}=\max\{u,0\}$ and $u^{-}=-\min\{u,0\}$.
\end{Lemma}

  Lemma \ref{y6.4.5} is the maximum principle for bounded solutions to parabolic equations.
\begin{Lemma}\label{y6.4.5}(Maximum principle for parabolic equations)
Let $\epsilon \geq 0$. If $H \in C([0,T]; C^{2}(\mathbb{T}^2\times \mathbb{R}^{+}) \cap C^{1}([0,T]; C^{0}(\mathbb{T}^2\times \mathbb{R}^{+}))$ is a bounded function that satisfies the differential inequality
\begin{eqnarray*}
\left\{\partial_{t}+b_{1}\partial_{x}+b_2\partial_{y}+b_3\partial_{z}-\epsilon^{2}\partial_{x}^{2}-\epsilon^{2}\partial_{y}^{2}-\partial_{z}^{2} \right\}H \leq fH \ \ \ in \ [0,T]  \times \mathbb{T}^2 \times \mathbb{R}^{+},
\end{eqnarray*}
where the coefficients $b_{1}$, $b_{2}$, $b_3$ and $f$ are continuous and satisfy
\begin{eqnarray}
\left\|\frac{b_{3}}{1+z}\right\|_{L^{\infty}([0,T] \times \mathbb{T}^2\times \mathbb{R}^{+})} < +\infty,
\end{eqnarray}
and
\begin{eqnarray}
 \|f\|_{L^{\infty}([0,T] \times \mathbb{T}^2 \times \mathbb{R}^{+})} \leq \lambda,
\end{eqnarray}
then for any $t \in [0,T]$,
\begin{eqnarray}
\sup \limits_{\mathbb{T}^2 \times \mathbb{R}^{+}} H(t) \leq \max\left\{e^{\lambda t}\|H(0)\|_{L^{\infty}(\mathbb{T}^2 \times \mathbb{R}^{+})}, \max \limits_{\tau \in [0,T]}\left\{e^{\lambda(t-\tau)}\|H(\tau)\big{|}_{z=0}\|_{L^{\infty}(\mathbb{T}^2)}\right\}\right\}.
\label{max}
\end{eqnarray}
\end{Lemma}
\textbf{Proof}.
For any $\mu > 0$, let us define
\begin{eqnarray}
\mathcal{H} := e^{-\lambda t} H - \mu \| \frac{b_3}{1+z} \|_{L^\infty}\cdot t - \mu \ln (1+z).
\end{eqnarray}
Then for any $\tilde{t} \in (0,T]$, we can get
\begin{align}
& \left\{\partial_{t}+b_{1}\partial_{x}+b_2\partial_{y}+b_3\partial_{z}-\epsilon^{2}\partial_{x}^{2}-\epsilon^{2}\partial_{y}^{2}-\partial_{z}^{2}+(\lambda - f) \right\}\mathcal{H} \nonumber\\
 &=\left\{\partial_{t}+b_{1}\partial_{x}+b_2\partial_{y}+b_3\partial_{z}-\epsilon^{2}\partial_{x}^{2}-\epsilon^{2}\partial_{y}^{2}-\partial_{z}^{2}+(\lambda - f) \right\} e^{-\lambda t} H\nonumber\\
 &\quad-\left\{\partial_{t}+b_{1}\partial_{x}+b_2\partial_{y}+b_3\partial_{z}-\epsilon^{2}\partial_{x}^{2}-\epsilon^{2}\partial_{y}^{2}-\partial_{z}^{2}+(\lambda - f) \right\}\mu \| \frac{b_3}{1+z} \|_{L^\infty} \cdot t\nonumber\\
 &\quad-\left\{\partial_{t}+b_{1}\partial_{x}+b_2\partial_{y}+b_3\partial_{z}-\epsilon^{2}\partial_{x}^{2}-\epsilon^{2}\partial_{y}^{2}-\partial_{z}^{2}+(\lambda - f) \right\}\mu \ln (1+z)
 \nonumber\\
 &=e^{-\lambda t}\left\{\partial_{t}+b_{1}\partial_{x}+b_2\partial_{y}+b_3\partial_{z}-\epsilon^{2}\partial_{x}^{2}-\epsilon^{2}\partial_{y}^{2}-\partial_{z}^{2} - f \right\}  H-\mu \| \frac{b_3}{1+z} \|_{L^\infty}\nonumber\\
 &\quad-(\lambda - f) \mu \| \frac{b_3}{1+z} \|_{L^\infty}\cdot t-\frac{\mu b_3}{1+z} -\frac{\mu }{(1+z)^2}\nonumber\\
 &< 0, \quad    \mathrm{in} \quad [0, \tilde{t}]   \times \mathbb{T}^2 \times \mathbb{R}^{+},
\end{align}
where we have used $-\mu\left\|\frac{b_{3}}{1+z}\right\|_{L^{\infty}([0,T] \times \mathbb{T}^2\times \mathbb{R}^{+})} -\frac{\mu b_3}{1+z} <0$, $f-\lambda \leq 0$, and the following fact:
\begin{eqnarray*}
\left\{\partial_{t}+b_{1}\partial_{x}+b_2\partial_{y}+b_3\partial_{z}-\epsilon^{2}\partial_{x}^{2}-\epsilon^{2}\partial_{y}^{2}-\partial_{z}^{2} \right\}H \leq fH \ \ \ in \ [0,T]  \times \mathbb{T}^2 \times \mathbb{R}^{+}.
\end{eqnarray*}
Now, we want to use the classical maximum principle (Lemma \ref{pmax}) for parabolic equations on $\mathcal{H}$. To apply this maximum principle, we need the function $\mathcal{H}$ to be defined and continuous on the closure of the domain $[0, \tilde{t}]   \times \mathbb{T}^2 \times \mathbb{R}^{+}$. Since $\ln(1+z)$ can blow up to infinity as $z$  approaches $-1$, we need to make sure that $R$ is chosen to be large enough such that $\ln(1+z)$ is well-defined and finite for all $(x,y,z) \in \mathbb{T}^2 \times \mathbb{R}^{+}$. And we also need to ensure that the maximum inequality itself is established.

So by the classical maximum principle  for parabolic equations, we have
\begin{align}
    \max_{[0,\tilde{t}] \times \mathbb{T}^2 \times [0,R]} \mathcal{H}  & \leq \max \left\{\|\max\mathcal{H}(0)\|_{L^\infty(\mathbb{T}^2 \times \mathbb{R}^{+})}, \max_{\tau \in [0, \tilde{t}]}
   \|\mathcal{H}(\tau) |_{z=0}\|_{L^\infty (\mathbb{T}^2)}  \right\}\nonumber\\
    & \leq \max \left\{\|H(0)\|_{L^\infty (\mathbb{T}^2 \times \mathbb{R}^{+})}, \max_{\tau \in [0, \tilde{t}]}
    \left\{e^{-\lambda \tau}
    \|H(\tau) |_{z=0} \|_{L^\infty (\mathbb{T}^2) } \right\} \right\}
\end{align}
provided that $R \geq \exp \left(\frac{1}{\mu} \|H\|_{L^\infty(\mathbb{T}^2 \times \mathbb{R}^{+})} \right) -1$.
Therefore, for any $(x,y,z) \in \mathbb{T}^2 \times \mathbb{R}^{+}$, we have
\begin{align}
    &  H(\tilde{t},x,y) - \mu \left\| \frac{b_3}{1+z}\right\|_{L^\infty} e^{\lambda \tilde{t}}  \tilde{t} - \mu e^{\lambda \tilde{t}} \ln (1+z)\nonumber\\
   & \leq  \max \left\{e^{\lambda \tilde{t}} \|H(0)\|_{L^\infty (\mathbb{T}^2 \times \mathbb{R}^{+})}, \max_{\tau \in [0,\tilde{t}]}
   \left \{ e^{\lambda (\tilde{t} -\tau)} \|H(\tau) |_{z=0} \|_{L^\infty (\mathbb{T}^2 )} \right\} \right\}
\end{align}
which is just inequality (\ref{max}), after taking the limit $\mu \rightarrow 0^+$ and replace the arbitrary time $\tilde{t}$ by $t$. Hence, this completes the proof of the lemma.

\hfill $\Box$

Lemma \ref{y6.4.6} is the lower bound estimate on bounded solutions for parabolic
equations.

\begin{Lemma}\label{y6.4.6}(Minimum principle for parabolic equations)
Let $\epsilon \geq 0$. If $H \in C([0,T]; C^{2}(\mathbb{T}^2\times \mathbb{R}^{+}) \cap C^{1}([0,T]; C^{0}(\mathbb{T}^2\times \mathbb{R}^{+}))$ is a bounded function with
 \begin{eqnarray*}
\kappa(t) :=\min \left\{ \min \limits_{\mathbb{T}^2 \times \mathbb{R}^{+}}H(0), \min \limits_{[0,T] \times \mathbb{T}^2} H \big{|}_{z=0} \right\} \geq 0
\end{eqnarray*}
and satisfies
\begin{eqnarray*}
\left\{\partial_{t}+b_{1}\partial_{x}+b_2\partial_{y}+b_3\partial_{z}-\epsilon^{2}\partial_{x}^{2}-\epsilon^{2}\partial_{y}^{2}-\partial_{z}^{2} \right\}H = fH
\end{eqnarray*}
where the coefficients $b_{1}$, $b_{2}$, $b_3$ and $f$ are continuous and satisfy
\begin{eqnarray}
\left\|\frac{b_{3}}{1+z}\right\|_{L^{\infty}([0,T] \times \mathbb{T}^2\times \mathbb{R}^{+})} <  +\infty,
\end{eqnarray}
and
\begin{eqnarray}
 \|f\|_{L^{\infty}([0,T] \times \mathbb{T}^2 \times \mathbb{R}^{+})} \leq \lambda,
\end{eqnarray}
then for any $t \in [0,T]$,
\begin{eqnarray}
\sup \limits_{\mathbb{T}^2 \times \mathbb{R}^{+}} H(t) \geq (1-\lambda t e^{\lambda t}) \kappa(t).
\label{min}
\end{eqnarray}

\end{Lemma}
\textbf{Proof}.
For any fixed $\tilde{t} \in [0,T]$ and $\mu > 0$, let us define
\begin{eqnarray}
 h := e^{- \lambda t}\left(H - \kappa (\tilde{t})\right) + \left(\lambda \kappa (\tilde{t}) + \mu \left\|
\frac{b_3}{1+z} \right\|_{L^\infty} \right)\cdot t + \mu \ln (1+z).
\end{eqnarray}

Then for any $\tilde{t} \in (0,T]$, we have
 \begin{align}
& \left\{\partial_{t}+b_{1}\partial_{x}+b_2\partial_{y}+b_3\partial_{z}-\epsilon^{2}\partial_{x}^{2}-\epsilon^{2}\partial_{y}^{2}-\partial_{z}^{2}+(\lambda - f) \right\}h \nonumber\\
 &=\left\{\partial_{t}+b_{1}\partial_{x}+b_2\partial_{y}+b_3\partial_{z}-\epsilon^{2}\partial_{x}^{2}-\epsilon^{2}\partial_{y}^{2}-\partial_{z}^{2}+(\lambda - f) \right\}e^{-\lambda t}  H\nonumber\\
  &\quad-\left\{\partial_{t}+b_{1}\partial_{x}+b_2\partial_{y}+b_3\partial_{z}-\epsilon^{2}\partial_{x}^{2}-\epsilon^{2}\partial_{y}^{2}-\partial_{z}^{2}+(\lambda - f) \right\}e^{- \lambda t} \kappa (\tilde{t})\nonumber\\
 &\quad+\left\{\partial_{t}+b_{1}\partial_{x}+b_2\partial_{y}+b_3\partial_{z}-\epsilon^{2}\partial_{x}^{2}-\epsilon^{2}\partial_{y}^{2}-\partial_{z}^{2}+(\lambda - f) \right\}\left(\lambda \kappa (\tilde{t}) + \mu \left\|
\frac{b_3}{1+z} \right\|_{L^\infty} \right) \cdot t\nonumber\\
 &\quad+\left\{\partial_{t}+b_{1}\partial_{x}+b_2\partial_{y}+b_3\partial_{z}-\epsilon^{2}\partial_{x}^{2}-\epsilon^{2}\partial_{y}^{2}-\partial_{z}^{2}+(\lambda - f) \right\}\mu \ln (1+z)
 \nonumber\\
 &=0-\left\{\partial_{t}+(\lambda - f) \right\}e^{- \lambda t} \kappa (\tilde{t})+\left\{\partial_{t}+(\lambda - f) \right\}\lambda \kappa (\tilde{t}) \cdot t+\mu \| \frac{b_3}{1+z} \|_{L^\infty}\nonumber\\
 &\quad+(\lambda - f) \mu \| \frac{b_3}{1+z} \|_{L^\infty}\cdot t+\frac{\mu b_3}{1+z} +\frac{\mu }{(1+z)^2}\nonumber\\
 &=\lambda  \kappa (\tilde{t})(e^{- \lambda t}+1)+(\lambda t-e^{- \lambda t})(\lambda - f)  \kappa (\tilde{t})+(\lambda - f) \mu \| \frac{b_3}{1+z} \|_{L^\infty}\cdot t+\frac{\mu b_3}{1+z} +\frac{\mu }{(1+z)^2}\nonumber\\
 &> 0, \quad    \mathrm{in} \quad [0, \tilde{t}]   \times \mathbb{T}^2 \times \mathbb{R}^{+},
\end{align}
where we have used $\lambda t-e^{- \lambda t} > 0$ for $t>0$, $\mu\left\|\frac{b_{3}}{1+z}\right\|_{L^{\infty}([0,T] \times \mathbb{T}^2\times \mathbb{R}^{+})}+\frac{\mu b_3}{1+z} \geq0$, $\lambda-f \geq 0$, and the following fact:
\begin{eqnarray*}
\left\{\partial_{t}+b_{1}\partial_{x}+b_2\partial_{y}+b_3\partial_{z}-\epsilon^{2}\partial_{x}^{2}-\epsilon^{2}\partial_{y}^{2}-\partial_{z}^{2} \right\}H = fH \ \ \ in \ [0,T]  \times \mathbb{T}^2 \times \mathbb{R}^{+},
\end{eqnarray*}
so by the classical maximum principle (Lemma \ref{pmax}) for parabolic equations, we have
 \begin{align}
    \min_{[0,\tilde{t}]\times \mathbb{T}^2 \times [0,R]} h \geq -\max \limits_{\partial([0,T]\times\mathbb{T}^2\times \mathbb{R}^{+})} \left(-\min\{0, h\}\right)\geq 0,
 \end{align}
provided that
\begin{eqnarray}
R \geq \exp \left(\frac{1}{\mu} \big( \|H\|_{L^\infty(\mathbb{T}^2 \times \mathbb{R}^{+})} + \kappa (\tilde{t})\big)\right) - 1.
\end{eqnarray}

Here, we have ensured that $R$ is chosen to be large enough such that $\ln(1+z)$ is well-defined and finite for all $(x,y,z) \in \mathbb{T}^2 \times \mathbb{R}^{+}$.  And we also ensured that the above inequality itself is established, that is
\begin{eqnarray}
\min \left\{e^{- \lambda t}\left(H - \kappa (\tilde{t})\right) + \left(\lambda \kappa (\tilde{t}) + \mu \left\|
\frac{b_3}{1+z} \right\|_{L^\infty} \right) t + \mu \ln (1+z)\right\} \geq 0.
\end{eqnarray}

Taking the limit $R \rightarrow + \infty$, and then $\mu \rightarrow 0^+$, we obtain
 \begin{eqnarray}
H(\tilde{t}) \ge (1-\lambda \tilde{t} e^{\lambda \tilde{t}} )\kappa(\tilde{t})
 \end{eqnarray}
which implies inequality (\ref{min}) if we replace the arbitrary time $\tilde{t}$ by $t$. Hence, this completes the proof of the lemma.

\hfill $\Box$

\section{Uniform estimates }
In this section, we will estimate the norm $\|\varphi\|_{H^{s,\gamma}_{g}}^{2}$ to   problem $(\ref{6.1.5})$ by the energy method, whose idea  comes from \cite{mw}.
\subsection{Estimates on $D^{\alpha}\varphi$ }
We will estimate the norm with weight $(1+z)^{ \gamma+ \alpha_{3}}$ on $D^{\alpha}\varphi$, $\alpha=(\alpha_{1},\alpha_{2}, \alpha_{3})$, $|\alpha|\leq s$, and $\alpha_{1}+\alpha_{2}$ should be restricted as $\alpha_{1}+\alpha_{2}\leq s-1$. Otherwise, it is impossible to directly estimate the norm with $\partial _{xy}^{s}\varphi$.

\begin{Lemma}(Reduction of boundary data)\label{y6.2.1}
If $\varphi$ solves $(\ref{6.1.4})$ and $(\ref{6.1.5})$, then on the boundary at $z=0$, we have,
\begin{equation}\left\{
\begin{array}{ll}
\partial _{z}\varphi|_{z=0}=\partial _{x}P,\\
\partial _{z}^{3}\varphi|_{z=0}=(\partial _{t}-\epsilon^2\partial_x^2-\epsilon^2\partial_y^2)\partial _{x}P+(\varphi \partial_x \varphi+K\varphi \partial_y \varphi)\big|_{z=0}-\partial_y K |\varphi|^2 \big|_{z=0}. \label{6.3.21}
\end{array}
         \right.\end{equation}
For any $5 \leq 2k+1 \leq s$, there are some constants  $C_{K,k,l,\rho^1,\rho^2,...,\rho^j}$,   not depending on $\epsilon$ (in this section, the superscript  $\epsilon$ is removed for simplification) or $(u,Ku, w)$, such that
\begin{eqnarray}
\begin{aligned}
\partial_z^{2k+1} \varphi|_{z=0}&=(\partial_t-\epsilon^2\partial^2_x-\epsilon^2\partial^2_y)^k\partial_x P \\
&\quad+\sum\limits_{h=0}^{2k-1}\partial_{xy}^{h} K\sum\limits_{l=0}^{k-1}\epsilon^{2l}\sum\limits_{j=2}^{\max\{2,k-l\}}\sum\limits_{\rho \in A^j_{k,l}}C_{K,k,l,\rho^1,\rho^2,...,\rho^j}\prod^{j}_{i=1}D^{\rho^i}\varphi \big{|}_{z=0},
\label{6.3.22}
\end{aligned}
\end{eqnarray}
where $A^j_{k,l}:=\{\rho:=(\rho^1,\rho^2,...,\rho^j)\in \mathbb{N}^{2j}; \sum\limits_{i=1}^{j}\rho^i_1 \leq k+2l-1, \sum\limits_{i=1}^{j}\rho^i_2 \leq k+2l-1, \sum\limits_{i=1}^{j}\rho^i_{12} \leq k+2l-1, \sum\limits_{i=1}^{j}\rho_3^i\leq 2k-2l-2, ~and~ |\rho^i| \leq 2k-l-1 ~for~ all~ i=1,2,...,j\}$.
\end{Lemma}
\textbf{Proof}.
According to the equation $(\ref{6.1.5})$ and the boundary condition $\partial _{z}\varphi|_{z=0}=\partial _{x}P$, we get
\begin{align}
\partial _{z}^{3}\varphi\big{|}_{z=0}&= \partial_z\left\{\partial _{t}\varphi+(u\partial _{x} +Ku\partial _{y}+w\partial_{z})\varphi -\epsilon^2\partial _{x}^{2}\varphi-\epsilon^2\partial _{y}^{2}\varphi- \frac{1}{2}\partial_y K \partial_z |u|^2 \right\}\big|_{z=0}\nonumber\\
&=(\partial _{t}-\epsilon^2\partial_x^2-\epsilon^2\partial_y^2)\partial _{x}P+(\varphi \partial_x \varphi+K\varphi \partial_y \varphi)\big|_{z=0}-\partial_y K |\varphi|^2 \big|_{z=0},
\label{6.3.23}
\end{align}
and
\begin{eqnarray}
\begin{aligned}
\varphi _{z}^{2n+1}\varphi \big{|}_{z=0}&= \left\{\left(\partial _{t}-\epsilon^2\partial_x^2\varphi-\epsilon^2\partial_y^2\varphi\right)\partial _{z}^{2n-1}\varphi+\partial _{z}^{2n-1}(u\partial _{x} \varphi +Ku\partial _{y}\varphi+w\partial_{z}\varphi-\frac{1}{2}\partial_y K \partial_z |u|^2)\right\}\big|_{z=0}.
\label{6.3.24}
\end{aligned}
\end{eqnarray}
Hence, for $n=2$, we arrive at
\begin{align}
\partial _{z}^{5}\varphi \big|_{z=0}
&= \left\{(\partial _{t}-\epsilon^2\partial_x^2-\epsilon^2\partial_y^2)\partial_z^{3}\varphi+\partial _{z}^{3}(u\partial _{x} +Ku\partial _{y}+w\partial_{z})\varphi-\partial _{z}^{3}(\frac{1}{2}\partial_y k \partial_z |u|^2) \right\}\Big|_{z=0}\nonumber\\
&=(\partial _{t}-\epsilon^2\partial_x^2-\epsilon^2\partial_y^2)^2\partial_x P+(\partial _{t}-\epsilon^2\partial_x^2-\epsilon^2\partial_y^2)\left\{(\varphi \partial_x \varphi+K\varphi \partial_y \varphi)-\partial_y K |\varphi|^2\right\}\big|_{z=0}\nonumber\\
&\quad+(3\varphi\partial_x\partial_z^2 \varphi+2\partial_z \varphi \partial_x \partial_z \varphi-2\partial_x \varphi \partial^2_z \varphi)|_{z=0}+K(3\varphi\partial_y\partial_z^2 \varphi+2\partial_z \varphi \partial_y \partial_z \varphi-2\partial_y \varphi \partial^2_z \varphi)|_{z=0}\nonumber\\
&\quad-\partial_y K(4\varphi\partial_z^2 \varphi+3\partial_z\varphi\partial_z \varphi)\big|_{z=0}.
\label{6.3.25}
\end{align}
Since the right-hand side of (\ref{6.3.25}) except for the second term in the formula is in the desired form, we only need to deal with the terms $(\partial _{t}-\epsilon^2\partial_x^2-\epsilon^2\partial_y^2)\left\{(\varphi \partial_x \varphi+K\varphi \partial_y \varphi)-\partial_y K |\varphi|^2\right\}\big|_{z=0}$. Using the evolution equations for $\varphi$, $\partial_x\varphi$ and $\partial_y\varphi$ as well as $(u, w) \big|_{z=0}=0$, we may check that
\begin{align}
&(\partial _{t}-\epsilon^2\partial_x^2-\epsilon^2\partial_y^2)\left\{(\varphi \partial_x \varphi+K\varphi \partial_y \varphi)-\partial_y K |\varphi|^2\right\}\big|_{z=0}\nonumber \\
&=\left\{\partial_t \varphi (\partial_x \varphi +K\partial_y \varphi -2\partial_y  K \varphi)+\varphi(\partial _{t}\partial_x \varphi+K\partial _{t}\partial_y \varphi)\right\}\big|_{z=0}\nonumber \\
&\quad-(\epsilon^2\partial_x^2+\epsilon^2\partial_y^2)\left\{(\varphi \partial_x \varphi+K\varphi \partial_y \varphi)-\partial_y K |\varphi|^2\right\}\big|_{z=0} \nonumber \\
&= \left( \epsilon^2\partial _{x}^{2}\varphi+\epsilon^2\partial _{y}^{2}\varphi+ \partial _{z}^{2}\varphi \right)(\partial_x \varphi +K\partial_y \varphi -2\partial_y  K \varphi)\big|_{z=0} \nonumber \\
&\quad+\varphi \left\{\left( \epsilon^2\partial _{x}^{3}\varphi+\epsilon^2\partial _x\partial _{y}^{2}\varphi+ \partial _x\partial _{z}^{2}\varphi \right)+K\left( \epsilon^2\partial _{y}\partial _{x}^{2}\varphi+\epsilon^2\partial _{y}^{3}\varphi+ \partial _y\partial _{z}^{2}\varphi \right)\right\}\big|_{z=0} \nonumber \\
& \quad-(\epsilon^2\partial_x^2+\epsilon^2\partial_y^2)\left\{(\varphi \partial_x \varphi+K\varphi \partial_y \varphi)-\partial_y K |\varphi|^2\right\}\big|_{z=0},
\label{6.3.26}
\end{align}
where  we have used the following facts,
\begin{align*}
\partial _{t} \varphi|_{z=0}&= \left\{-(u\partial _{x} +Ku\partial _{y}+w\partial_{z})\varphi +\epsilon^2\partial _{x}^{2}\varphi+\epsilon^2\partial _{y}^{2}\varphi+ \partial _{z}^{2}\varphi+\frac{1}{2}\partial_y K \partial_z |u|^2 \right\}\big{|}_{z=0} \nonumber\\
&=\left( \epsilon^2\partial _{x}^{2}\varphi+\epsilon^2\partial _{y}^{2}\varphi+ \partial _{z}^{2}\varphi \right)\big{|}_{z=0}, \nonumber\\
\partial _{t}\partial_x \varphi|_{z=0}&= \partial_x\left\{-(u\partial _{x} +Ku\partial _{y}+w\partial_{z})\varphi +\epsilon^2\partial _{x}^{2}\varphi+\epsilon^2\partial _{y}^{2}\varphi+ \partial _{z}^{2}\varphi+\frac{1}{2}\partial_y K \partial_z |u|^2 \right\}\big{|}_{z=0} \nonumber\\
&=\left( \epsilon^2\partial _{x}^{3}\varphi+\epsilon^2\partial _x\partial _{y}^{2}\varphi+ \partial _x\partial _{z}^{2}\varphi \right)\big{|}_{z=0}, \nonumber\\
\partial _{t}\partial_y \varphi|_{z=0}&=\left( \epsilon^2\partial _{y}\partial _{x}^{2}\varphi+\epsilon^2\partial _{y}^{3}\varphi+ \partial _y\partial _{z}^{2}\varphi \right)\big{|}_{z=0}.
\end{align*}
Combining   (\ref{6.3.26}) and (\ref{6.3.25}), we have
\begin{align}
\partial _{z}^{5}\varphi \big|_{z=0}
&=(\partial _{t}-\epsilon^2\partial_x^2-\epsilon^2\partial_y^2)^2\partial_x P+(3\varphi\partial_x\partial_z^2 \varphi+2\partial_z \varphi \partial_x \partial_z \varphi-2\partial_x \varphi \partial^2_z \varphi)|_{z=0}\nonumber\\
&\quad+K(3\varphi\partial_y\partial_z^2 \varphi+2\partial_z \varphi \partial_y \partial_z \varphi-2\partial_y \varphi \partial^2_z \varphi)|_{z=0}-\partial_y K(4\varphi\partial_z^2 \varphi+3\partial_z\varphi\partial_z \varphi)\big|_{z=0}\nonumber \\
&\quad+ \left( \epsilon^2\partial _{x}^{2}\varphi+\epsilon^2\partial _{y}^{2}\varphi+ \partial _{z}^{2}\varphi \right)(\partial_x \varphi +K\partial_y \varphi -2\partial_y  K \varphi)\big|_{z=0} \nonumber \\
&\quad+\varphi \left\{\left( \epsilon^2\partial _{x}^{3}\varphi+\epsilon^2\partial _x\partial _{y}^{2}\varphi+ \partial _x\partial _{z}^{2}\varphi \right)+K\left( \epsilon^2\partial _{y}\partial _{x}^{2}\varphi+\epsilon^2\partial _{y}^{3}\varphi+ \partial _y\partial _{z}^{2}\varphi \right)\right\}\big|_{z=0} \nonumber \\
& \quad-(\epsilon^2\partial_x^2+\epsilon^2\partial_y^2)\left\{(\varphi \partial_x \varphi+K\varphi \partial_y \varphi)-\partial_y K |\varphi|^2\right\}\big|_{z=0},
\end{align}
which implies that
\begin{align}
\partial _{z}^{5}\varphi \big|_{z=0}
&=(\partial _{t}-\epsilon^2\partial_x^2-\epsilon^2\partial_y^2)^2\partial_x P\nonumber\\
&\quad+(4\varphi\partial_x\partial_z^2 \varphi+2\partial_z \varphi \partial_x \partial_z \varphi-\partial_x \varphi \partial^2_z \varphi)|_{z=0}+K(4\varphi\partial_y\partial_z^2 \varphi+2\partial_z \varphi \partial_y \partial_z \varphi-\partial_y \varphi \partial^2_z \varphi)|_{z=0}\nonumber \\
&\quad-\partial_y K(6\varphi\partial_z^2 \varphi+3\partial_z\varphi\partial_z \varphi)\big|_{z=0}\nonumber \\
&\quad+ \left( \epsilon^2\partial _{x}^{2}\varphi+\epsilon^2\partial _{y}^{2}\varphi \right)(\partial_x \varphi +K\partial_y \varphi -2\partial_y  K \varphi)\big|_{z=0} +\varphi \left( \epsilon^2\partial _{x}^{3}\varphi+\epsilon^2\partial _x\partial _{y}^{2}\varphi \right)\big|_{z=0} \nonumber \\
&\quad+K\left( \epsilon^2\partial _{y}\partial _{x}^{2}\varphi+\epsilon^2\partial _{y}^{3}\varphi \right)\big|_{z=0} \nonumber \\
& \quad-(\epsilon^2\partial_x^2+\epsilon^2\partial_y^2)\left\{(\varphi \partial_x \varphi+K\varphi \partial_y \varphi)-\partial_y K |\varphi|^2\right\}\big|_{z=0}.
\label{r3.25}
\end{align}
We start completing the verification of formula  (\ref{6.3.22}) for $k=2$. When $k=2$, formula  (\ref{6.3.22}) becomes the following the  equality
\begin{eqnarray}
\begin{aligned}
\partial_z^{5} \varphi|_{z=0}&=(\partial_t-\epsilon^2\partial^2_x-\epsilon^2\partial^2_y)^2\partial_x P +\sum\limits_{h=0}^{3}\partial_{xy}^{h} K\sum\limits_{l=0}^{1}\epsilon^{2l}\sum\limits_{\rho \in A^2_{2,l}}C_{K,k,l,\rho^1,\rho^2}\prod^{2}_{i=1}D^{\rho^i}\varphi \big{|}_{z=0},
\label{r3.26}
\end{aligned}
\end{eqnarray}
where  $\sum\limits_{i=1}^{2}\rho^i_1 \leq 1+2l, \sum\limits_{i=1}^{2}\rho^i_2 \leq 1+2l, \sum\limits_{i=1}^{2}\rho^i_{12} \leq 1+2l, \sum\limits_{i=1}^{2}\rho_3^j\leq 2-2l, ~and~ |\rho^i| \leq 3-l ~for~ all~ i=1,2$.

When $l=0$, then $\epsilon^{2l}$=1, by  the second, third and fourth terms on the right-hand side of the equation (\ref{r3.25}),  it holds  that
\begin{eqnarray*}
\sum\limits_{i=1}^{2}\rho^i_1 \leq 1, \ \sum\limits_{i=1}^{2}\rho^i_2 \leq 1, \ \sum\limits_{i=1}^{2}\rho^i_{12} \leq 1,  \ \sum\limits_{i=1}^{2}\rho_3^j\leq 2, \ and \ |\rho^i| \leq 3 \ for \  all \ i=1,2.
\end{eqnarray*}
When $l=1$, then $\epsilon^{2l}$=$\epsilon^{2}$, by  the fifth, sixth and seventh terms on the right-hand side of the equation (\ref{r3.25}),  it holds  that
\begin{eqnarray*}
\sum\limits_{i=1}^{2}\rho^i_1 \leq 3, \ \sum\limits_{i=1}^{2}\rho^i_2 \leq 3, \ \sum\limits_{i=1}^{2}\rho^i_{12} \leq 3,  \ \sum\limits_{i=1}^{2}\rho_3^j\leq 0, \ and \ |\rho^i| \leq 2 \ for \  all \ i=1,2.
\end{eqnarray*}
Note that although $\epsilon^2\partial _{x}^{3}\varphi$ and $\epsilon^2\partial _{y}^{3}\varphi$ appears in the seventh term of the right-hand side of the equation (\ref{r3.25}), it can be canceled out by $-(\epsilon^2\partial_x^2+\epsilon^2\partial_y^2)\left(\varphi \partial_x \varphi+K\varphi \partial_y \varphi\right)$ in the eighth term. Therefore, $\ |\rho^i| \leq 2$ is indeed true.

Moreover,  $\partial_{xy}^{h} K$ is not only related to $k$ in (\ref{6.3.22}) and (\ref{r3.26}) in fact, but also related to other quantities. However, due to assumption (\ref{6.1.88}), we do not need to accurately give the relationship between $\partial_{xy}^{h} K$ and other quantities.

According to the above argument, the  equality (\ref{r3.25}) satisfies (\ref{r3.26}), we have completed the verification of formula  (\ref{6.3.22}) for $k=2$.

Now, using the same algorithm, we are going to prove formula (\ref{6.3.22}) by induction on $k$. For notational convenience, we denote
\begin{eqnarray}
\mathcal{A}_k:=\bigg\{\sum\limits_{h=0}^{2k-1}\partial_{xy}^{h} K\sum\limits_{l=0}^{k-1}\epsilon^{2l}\sum\limits_{j=2}^{\max\{2,k-l\}}\sum\limits_{\rho \in A^j_{k,l}}C_{K,k,l,\rho^1,\rho^2,...,\rho^j}\prod^{j}_{i=1}D^{\rho^i}\varphi\big{|}_{z=0}\bigg\}.
\label{6.3.27}
\end{eqnarray}
Using this notation, we will prove $\partial_{z}^{2k+1}\varphi\big{|}_{z=0}-(\partial_t-\epsilon^2\partial_x^2-\epsilon^2\partial_y^2)^k\partial_x  P \in \mathcal{A}_k$. Assuming that formula (\ref{6.3.22}) holds for $k=n$, we will show that it also holds for $k=n+1$ as follows. Then we differentiate the vorticity equation $2n+1$ times  with respect to $z$ to obtain
\begin{align}
\partial _{z}^{2n+3}\varphi\big{|}_{z=0}&= \left\{\left(\partial _{t}-\epsilon^2\partial_x^2\varphi-\epsilon^2\partial_y^2\varphi\right)\partial _{z}^{2n+1}\varphi+\partial _{z}^{2n+1}(u\partial _{x} \varphi +Ku\partial _{y}\varphi+w\partial_{z}\varphi-\partial_y K \varphi u)\right\}\big|_{z=0}\nonumber\\
&=\left(\partial _{t}-\epsilon^2\partial_x^2\varphi-\epsilon^2\partial_y^2\varphi\right)\partial _{z}^{2n+1}\varphi \big|_{z=0} \nonumber\\
&\quad  + \sum^{2n+1}_{m=1}\binom{2n+1}{m} \partial^{m-1}_z \varphi (\partial_x+K\partial_y-\partial_y K) \partial_z^{2n-m+1}\varphi\big|_{z=0} \nonumber\\
&\quad - \sum^{2n+1}_{m=2}\binom{2n+1}{m} (\partial_x+K\partial_y) \partial^{m-2}_z \varphi  \partial_z^{2n-m+2}\varphi\big|_{z=0}.
\label{6.3.28}
\end{align}
The last two terms of (\ref{6.3.28}) satisfy the situation where $k=n+1$, $j=2$ and $l=0$ in (\ref{6.3.22}). So $\sum\limits_{i=1}^{2}\rho_3^i$ and $|\rho^i|$  need to satisfy $\sum\limits_{i=1}^{j}\rho_3^i\leq 2n$ and  $|\rho^i| \leq 2n+1$ for all $i=1,2,...,j$ respectively. By routine checking, it is easy to obtain $\sum\limits_{i=1}^{j}\rho_3^i=(m-1)+2n-m+1~\mathrm{or} ~(m-2)+2n-m+2=2n$ and $\max|\rho^i|  \leq 2n+1$ (i.e.,  $m=1$ in $ (\partial_x+K\partial_y) \partial_z^{2n-m+1}\varphi$),
one may show that the last two terms of (\ref{6.3.28}) belong to
$\mathcal{A}_{k+1}$, so it only remains to deal with the term $(\partial _{t}-\epsilon^2\partial_x^2-\epsilon^2\partial_y^2)\partial _{z}^{2n+1}\varphi\big{|}_{z=0}$.

Thanks to the induction hypothesis, there exist constants $C_{K,n,l,\rho^1,\rho^2,...,\rho^j}$ such that
\begin{align}
\partial_z^{2n+1}\varphi\big{|}_{z=0}&=(\partial_t-\epsilon^2\partial^2_x-\epsilon^2\partial^2_y)^n\partial_x P\nonumber\\
&\quad+\sum\limits_{h=0}^{2k-1}\partial_{xy}^{h}K\sum\limits_{l=0}^{k-1}\epsilon^{2l}\sum\limits_{j=2}^{\max\{2,n-l\}}\sum\limits_{\rho \in A^j_{n,l}}C_{K,n,l,\rho^1,\rho^2,...,\rho^j}\prod^{j}_{i=1}D^{\rho^i}\varphi\big{|}_{z=0}.
\label{6.3.29}
\end{align}
Thus we have, up to a relabeling of the indices $\rho^i$,
\begin{align}
&(\partial_t-\epsilon^2\partial^2_x-\epsilon^2\partial^2_y)\partial_z^{2n+1}\varphi\big{|}_{z=0} \nonumber\\
&=(\partial_t-\epsilon^2\partial^2_x-\epsilon^2\partial^2_y)^{n+1}\partial_x P \nonumber\\
&\quad+\sum\limits_{h=0}^{2n-1}\partial_{xy}^{h}K\sum\limits_{l=0}^{n-1}\epsilon^{2l}\sum\limits_{j=2}^{\max\{2,n-l\}}\sum\limits_{\rho \in A^j_{n,l}} \widetilde{C}_{K,n,l,\rho^1,\rho^2,...,\rho^j}(\partial_t-\epsilon^2\partial^2_x-\epsilon^2\partial^2_y)D^{\rho^1}\varphi\prod^{j}_{i=2}D^{\rho^i}\varphi\big{|}_{z=0}\nonumber\\
&\quad-\sum\limits_{h=0}^{2n-1}\partial_{xy}^{h}K\sum\limits_{l=0}^{n-1}\epsilon^{2l+2}\sum\limits_{j=2}^{\max\{2,n-l\}}\sum\limits_{\rho \in A^j_{n,l}}\widetilde{\widetilde{C}}_{K,n,l,\rho^1,\rho^2,...,\rho^j}(\partial_x+\partial_y) D^{\rho^1}\varphi(\partial_x+\partial_y) D^{\rho^2}\varphi\prod^{j}_{i=3}D^{\rho^i}\varphi\big{|}_{z=0},
\label{6.3.30}
\end{align}
where $\widetilde{C}_{K,n,l,\rho^1,\rho^2,...,\rho^j}$ and $\widetilde{\widetilde{C}}_{K,n,l,\rho^1,\rho^2,...,\rho^j}$ are some new constants depending on $C_{K,n,l,\rho^1,\rho^2,...,\rho^j}$.

For $\mathcal{A}_{n}$, $\rho^{i}\leq 2n-l-1$; for $\mathcal{A}_{n+1}$, $\rho^{i}\leq 2n-l+1$. It can be seen that every time the number of $n$ increases by 1, the maximum value of the order $\rho^{i}$ will increase by 2.  Due to

\begin{align}
-\sum\limits_{h=0}^{2n-1}\partial_{xy}^{h}K\sum\limits_{l=0}^{n-1}\epsilon^{2l}\sum\limits_{j=2}^{\max\{2,n-l\}}\sum\limits_{\rho \in A^j_{n,l}}\widetilde{\widetilde{C}}_{K,n,l,\rho^1,\rho^2,...,\rho^j}\prod^{j}_{i=2}D^{\rho^i}\varphi\big{|}_{z=0} \in \mathcal{A}_{n+1},
\end{align}
the last term on the right-hand side of (\ref{6.3.30}) belongs to $\mathcal{A}_{n+1}$. So it remains to check whether the second term on right-hand
side of (\ref{6.3.30}) also belongs to $\mathcal{A}_{n+1}$.

Indeed, differentiating the vorticity equation $(\ref{6.1.5})_1$ $\rho_1^1$ times,   $\rho_2^1$ times, $\rho_3^1$ times  with respect to $x,y,z$ respectively,  and then calculating  at $z=0$, we have, by denoting $e_3=(0,0,1)$,
\begin{align}
&(\partial_t-\epsilon^2\partial_x^2-\epsilon^2\partial_y^2)D^{\rho^1}\varphi \big{|}_{z=0} \nonumber \\
&=-\sum_{\substack{\beta \le \rho^1 \\
\beta_3 \geq 1}}\binom{\rho^1}{\beta} D^{\beta-e_3} \varphi \partial_x D^{\rho_1-\beta} \varphi \big{|}_{z=0}
-\sum_{\substack{\beta \le \rho^1 \\
\beta_3 \geq 1}}\binom{\rho^1}{\beta} D^{\beta-e_3} (K\varphi) \partial_y D^{\rho_1-\beta} \varphi \big{|}_{z=0}\nonumber\\
&\quad +\sum_{\substack{\beta \le \rho^1 \\
\beta_3 \geq 2}}\binom{\rho^1}{\beta} D^{\beta-2e_3} (\partial_x \varphi +\partial_y(K\varphi)) \partial_z D^{\rho_1-\beta} w \big{|}_{z=0}+\partial_z^2 D^{\rho^1}\varphi \big{|}_{z=0}\nonumber\\
&\quad +\sum_{\substack{\beta \le \rho^1 \\
\beta_3 \geq 1}}\binom{\rho^1}{\beta} D^{\beta-e_3} (\partial_{y}K\varphi)  D^{\rho_1-\beta} \varphi \big{|}_{z=0}.
\label{r3.31}
\end{align}
Using (\ref{r3.31}), one may justify by a routine counting of indices that the second term
on the right-hand side of (\ref{6.3.30}) belongs to $\mathcal{A}_{n+1}$.
This thus completes the proof.
\hfill $\Box$

\begin{Lemma}\label{y6.2.2}
Let $s \geq 5$ be an   integer, $\gamma\geq 1$, $\sigma>\gamma+\frac{1}{2}$ and $\delta\in (0,1)$. If  $\varphi \in H^{s,\gamma}_{\sigma,\delta}$  solves $(\ref{6.1.4})$  and $(\ref{6.1.5})$, we have the following  estimates:\\
(i)\;when $|\alpha|\leq s-1$,
\begin{eqnarray}
\bigg|\iint_{\mathbb{T}^2 } D^{\alpha}\varphi\partial _{z}D^{\alpha} \varphi dxdy \big{|} _{z=0}\bigg|\leq \frac{1}{12}\|(1+z)^{\gamma+\alpha_3+1}\partial^2_z D^\alpha \varphi\|_{L^2}^2+C\|\varphi\|^2_{H^{s,\gamma}_g},
\end{eqnarray}
(ii)\;when $|\alpha|=s$, $\alpha_3$ is even,
\begin{align}
\bigg|\iint_{\mathbb{T} ^2} D^{\alpha}\varphi\partial _{z}D^{\alpha}\varphi dxdy \big{|} _{z=0}\bigg|\leq &\frac{1}{12}\|(1+z)^{\gamma+\alpha_2}\partial_z D^\alpha \varphi\|^2_{L^2}+ \sum\limits_{l=0}^{\frac{s}{2}}\|\partial^{l}_t\partial_xP\|_{H^{s-2l}(\mathbb{T}^2)}^{2}\nonumber\\
&+C_{s,\gamma,\sigma,\delta}\|\partial_{xy}^s K\|_{L^{\infty}(\mathbb{T}^2)}(1+\|\varphi\|^2_{H^{s,\gamma}_g})^{s-2}\|\varphi\|^2_{H^{s,\gamma}_g},
\end{align}
(iii)\;when $|\alpha|=s$, $\alpha_3$ is odd,
\begin{align}
\bigg|\iint_{\mathbb{T} ^2} D^{\alpha}\varphi\partial _{z}D^{\alpha}\varphi dxdy \big{|} _{z=0}\bigg|\leq& \frac{1}{12}\|(1+z)^{\gamma+\alpha_3+1}\partial_x^{\alpha_1-1}\partial_y^{\alpha_2}\partial_z^{\alpha_3+2} \varphi\|^2_{L^2}+\sum\limits_{l=0}^{\frac{s}{2}}\|\partial^{l}_t\partial_xP\|_{H^{s-2l}(\mathbb{T})^2}^{2}\nonumber\\
&+C_{s,\gamma,\sigma,\delta}\|\partial_{xy}^s K\|_{L^{\infty}(\mathbb{T}^2)}(1+\|\varphi\|^2_{H^{s,\gamma}_g})^{s-2}\|\varphi \|^2_{H^{s,\gamma}_g}.
\end{align}
\end{Lemma}
\textbf{Proof}.
\emph{Case 1.} When $|\alpha|\leq s-1$,  using the following trace estimate
\begin{eqnarray}
\iint_{\mathbb{T}^2 } |f|dxdy \big{|}_{z=0} \leq C\left(\int^1_0\int_{\mathbb{T}^2}|f|dxdydz+\int^1_0\int_{\mathbb{T}^2}|\partial_z f|dxdydz\right),
\end{eqnarray}
we have
\begin{eqnarray}
\bigg|\int_{\mathbb{T}^2 } D^{\alpha}\varphi \partial _{z}D^{\alpha}\varphi dxdy \big{|} _{z=0}\bigg|\leq \frac{1}{12}\|(1+z)^{\gamma+\alpha_3+1}\partial^2_z D^\alpha \varphi\|_{L^2}^2+C\|\varphi\|^2_{H^{s,\gamma}_g}.
\end{eqnarray}
\emph{Case 2.} When $|\alpha|=s$, $\alpha_3=2k$ for some $k\in \mathbb{N}$, we can apply   Lemma \ref{y6.2.1} to $\partial_z D^\alpha \varphi | _{z=0}$ to obtain
\begin{align}
&\int_{\mathbb{T} } D^{\alpha}\varphi\partial _{z}D^{\alpha}\varphi dxdy \big{|} _{z=0} \nonumber \\
&=\int_{\mathbb{T}^2 } D^{\alpha}\varphi(\partial_t-\epsilon^2\partial_x^2-\epsilon^2\partial_y^2)^k\partial_x^{\alpha_1}\partial_y^{\alpha_2}\partial_x P dx dy\big{|} _{z=0}\nonumber\\
&\quad+\sum\limits_{l=0}^{2k-1}\partial_{xy}^{l} K\sum\limits_{l=0}^{k-1}\epsilon^{2l}\sum\limits_{j=2}^{\max\{2,k-l\}}\sum\limits_{\rho \in A^j_{k,l}}C_{K,k,l,\rho^1,\rho^2,...,\rho^j}\int_{\mathbb{T}^2 }D^\alpha \varphi \partial_x^{\alpha_1}\partial_y^{\alpha_2}(\prod^{j}_{i=1}D^{\rho^i}\varphi)dx\big{|}_{z=0}.
\end{align}
Then we again apply the simple trace estimate to control the boundary integral as follows
\begin{align}
\bigg|\int_{\mathbb{T} ^2} D^{\alpha}\varphi \partial _{z}D^{\alpha}\varphi dxdy \big{|} _{z=0}\bigg|\leq& \frac{1}{12}\|(1+z)^{\gamma+\alpha_3}\partial_z D^\alpha \varphi\|^2_{L^2}+ \sum\limits_{l=0}^{\frac{s}{2}}\|\partial^{l}_t\partial_xP\|_{H^{s-2l}(\mathbb{T}^2)}^{2}\nonumber\\
&+C_{s,\gamma,\sigma,\delta}\|\partial_{xy}^s K\|_{L^{\infty}(\mathbb{T}^2)}(1+\|\varphi\|_{H_g^{s,\gamma}})^{s-2}\|\varphi\|^2_{H^{s,\gamma}_g}.
\end{align}
\emph{Case 3.} When $|\alpha|=s$, $\alpha_3=2k+1$ for some $k\in \mathbb{N}$, using integration by parts in the $x,y$ variables, we have
\begin{eqnarray}
\int_{\mathbb{T}^2 } D^{\alpha}\varphi \partial _{z}D^{\alpha}\varphi dxdy \big{|} _{z=0}=-\int_{\mathbb{T}^2}\partial_x D^\alpha \varphi \partial_x^{\alpha_1-1}\partial_y^{\alpha_2}\partial_z^{\alpha_3+1}\varphi dxdy \big{|}_{z=0}.
\end{eqnarray}
The term $\partial_x D^\alpha \varphi \big{|}_{z=0}=\partial_x^{\alpha_1-1}\partial_y^{\alpha_2}\partial_z^{\alpha_3+1}\varphi \big{|} _{z=0}$ has an odd number of derivatives in $z$, then we get
\begin{align}
\bigg|\int_{\mathbb{T} ^2} D^{\alpha}\varphi\partial _{y}D^{\alpha}\varphi dxdy \big{|} _{z=0}\bigg|\leq& \frac{1}{12}\|(1+z)^{\gamma+\alpha_3+1}\partial_x^{\alpha_1-1}\partial_y^{\alpha_2}\partial_z^{\alpha_3+2}\varphi \|^2_{L^2}+ \sum\limits_{l=0}^{\frac{s}{2}}\|\partial^{l}_t\partial_xP\|_{H^{s-2l}(\mathbb{T}^2)}^{2}\nonumber\\
&+C_{s,\gamma,\sigma,\delta}\|\partial_{xy}^s K\|_{L^{\infty}(\mathbb{T}^2)}(1+\|\varphi\|_{H_g^{s,\gamma}})^{s-2}\|\varphi\|^2_{H^{s,\gamma}_g}.
\end{align}
 \hfill $\Box$

\begin{Proposition}\label{p6.2.1}
Let $s \geq 5$ be an   integer, $\gamma\geq 1$, $\sigma>\gamma+\frac{1}{2}$ and $\delta\in (0,1)$. If  $\varphi \in H^{s,\gamma}_{\sigma,\delta}$  solves $(\ref{6.1.4})$  and $(\ref{6.1.5})$, then we have the following uniform (in $\epsilon$) estimate:
\begin{eqnarray}
\begin{aligned}
\frac{d}{dt}\|(1+z)^{ \gamma+ \alpha_{3}}D^{\alpha}\varphi\|_{L^{2}}^{2} &\leq C_{s, \gamma,\sigma,\delta  }\left(1+\|\partial_{xy}^{s} K\|_{L^\infty(\mathbb{T}^2)}\right)\left(\|\varphi\|_{H^{s,\gamma}_{g}}+\|\partial_{xy}^s U\|_{L^2(\mathbb{T}^2)}+1\right)\|\varphi \|_{H^{s,\gamma}_{g}}^{2} \\
& \quad +C_{s,\gamma,\sigma,\delta}\|\partial_{xy}^s K\|_{L^{\infty}(\mathbb{T}^2)}(1+\|w\|_{H_g^{s,\gamma}})^{s-2}\|\varphi\|^2_{H^{s,\gamma}_g}\\
&\quad +\sum\limits_{l=0}^{\frac{s}{2}}\|\partial^{l}_t\partial_xP\|_{H^{s-2l}(\mathbb{T}^2)}^{2}. \label{6.2.13}
\end{aligned}
\end{eqnarray}
\end{Proposition}
\textbf{Proof}.
Applying the operator $D^{\alpha}$ on $(\ref{6.1.5})_{1}$ with $\alpha=(\alpha_{1},\alpha_{2},\alpha_{3}),~|\alpha|\leq s, ~\alpha_{1}+\alpha_{2}\leq s-1$,
\begin{align}
&(\partial _{t} +u\partial _{x} +Ku\partial _{y}+w\partial _{z}-\varepsilon^2\partial _{x}^{2}\varphi-\varepsilon^2\partial _{y}^{2}\varphi-\partial _{z}^{2})D^{\alpha}\varphi \nonumber\\
&=- \sum\limits_{0<\beta\leq \alpha}\binom{ \alpha}{\beta} \left\{D^{\beta} u\partial _{x}D^{\alpha-\beta}\varphi  + D^{\beta}(Ku)\partial _{y}D^{\alpha-\beta}\varphi + D^{\beta} w\partial _{z}D^{\alpha-\beta}\varphi \right\}\nonumber\\
&\quad+ \sum\limits_{0 \leq \beta\leq \alpha}\binom{ \alpha}{\beta} D^{\beta} (u\partial_y K )D^{\alpha-\beta}\varphi .
 \label{6.2.1}
 \end{align}
Multiplying (\ref{6.2.1}) by $(1+z)^{2\gamma+2\alpha_{3}}D^{\alpha}\varphi$ and integrating it over $\mathbb{T}^{2}\times \mathbb{R}^{+}$ yield
\begin{eqnarray}
\begin{aligned}
&\frac{1}{2}\frac{d}{dt}\|(1+z)^{ \gamma+ \alpha_{3}}D^{\alpha}\varphi\|_{L^{2}}^{2}
+\|(1+z)^{ \gamma+ \alpha_{3}}\partial _{z}D^{\alpha}\varphi\|_{L^{2}}^{2}+\varepsilon^2\|(1+z)^{ \gamma+ \alpha_{3}}\partial _{z}D^{\alpha}\varphi\|_{L^{2}}^{2}+\varepsilon^2\|(1+z)^{ \gamma+ \alpha_{3}}\partial _{z}D^{\alpha}\varphi\|_{L^{2}}^{2}   \\
&=-\iint_{\mathbb{T} ^{2}} D^{\alpha}\varphi\partial _{z}D^{\alpha}\varphi dxdy | _{z=0}-(2\gamma+2\alpha_{3}) \iiint(1+z)^{2\gamma+2\alpha_{3}-1} D^{\alpha}\varphi\partial _{z}   D^{\alpha}\varphi \\
&\quad +(\gamma+\alpha_{3}) \iiint(1+z)^{2\gamma+2\alpha_{3}-1} w| D^{\alpha}\varphi|^{2}  \\
&\quad  - \sum\limits_{0<\beta\leq \alpha}\binom{ \alpha}{\beta} \iiint(1+z)^{2\gamma+2\alpha_{3}}D^{\alpha}\varphi  \left\{D^{\beta} u\partial _{x}D^{\alpha-\beta}\varphi + D^{\beta}(Ku)\partial _{y}D^{\alpha-\beta}\varphi+ D^{\beta} w\partial _{z}D^{\alpha-\beta}\varphi \right\} \\
&\quad  + \sum\limits_{0<\beta\leq \alpha}\binom{ \alpha}{\beta} \iiint(1+z)^{2\gamma+2\alpha_{3}}D^{\alpha}\varphi  \left\{D^{\beta} (u \partial_y K )D^{\alpha-\beta}\varphi  \right\}.
\label{6.2.2}
\end{aligned}
\end{eqnarray}
Indeed, $(\ref{6.2.2})$ can be obtained by integration by parts and by using the boundary condition, i.e.,
\begin{eqnarray}
\begin{aligned}
& \iiint(1+z)^{2\gamma+2\alpha_{3}}D^{\alpha} \varphi(u\partial _{x} +ku\partial _{y} +w\partial _{z} )D^{\alpha}\varphi \\
&=-(\gamma +\alpha_{3})\iiint(1+z)^{2\gamma+2\alpha_{3}-1} w |D^{\alpha}\varphi |^{2},
 \label{6.2.3}
\end{aligned}
\end{eqnarray}
and
\begin{eqnarray}
\begin{aligned}
 -\iiint(1+z)^{2\gamma+2\alpha_{3}}D^{\alpha}\varphi \partial _{z}^{2}D^{\alpha}\varphi & = \|(1+z)^{\gamma+ \alpha_{3}}\partial _{z}D^{\alpha}\varphi\|_{L^{2}}+\iint_{\mathbb{T}^{2} } D^{\alpha}\varphi\partial _{z}D^{\alpha}\varphi dxdy | _{z=0}\\
& \quad +  (2\gamma+2\alpha_{3})  \iiint(1+z)^{2\gamma+2\alpha_{3}-1}D^{\alpha}\varphi\partial _{z} D^{\alpha}\varphi.
 \label{6.2.4}
\end{aligned}
\end{eqnarray}
We need to estimate   equation $(\ref{6.2.2})$ term by term. Obviously, the first term on the right-hand side of $(\ref{6.2.2})$ follows from Lemmas \ref{y6.2.1}-\ref{y6.2.2}.

Secondly, using the H$\ddot{o}$lder inequality, we have
\begin{eqnarray}
&&\left|(2\gamma+2\alpha_{3}) \iiint(1+z)^{2\gamma+2\alpha_{3}-1} D^{\alpha}\varphi\partial _{z}   D^{\alpha}\varphi \right|\nonumber \\
 &\leq& (2\gamma+2\alpha_{3}) \|\frac{1}{1+z }  \|_{L^{\infty}} \|(1+z)^{\gamma+ \alpha_{3}}D^{\alpha}\varphi\|_{L^{2}}
  \|(1+z)^{\gamma+ \alpha_{3}}\partial_{z} D^{\alpha}\varphi \|_{L^{2}} \nonumber \\
  & \leq & C_{ \gamma,\alpha } \|(1+z)^{\gamma+ \alpha_{3}}D^{\alpha}\varphi \|_{L^{2}}^{2}+\frac{1}{4}\|(1+z)^{\gamma+ \alpha_{3}}\partial_{z} D^{\alpha}\varphi \|_{L^{2}}^{2} \nonumber \\
& \leq & C_{ \gamma,\alpha } \|\varphi \|_{H^{s, \gamma}_g}^{2}+\frac{1}{4}\|(1+z)^{\gamma+ \alpha_{3}}\partial_{z} D^{\alpha}\varphi \|_{L^{2}}^{2} ,  \label{6.2.5}
\end{eqnarray}
which, together with Lemma \ref{y6.4.4}, yields
\begin{eqnarray}
\left|   (\gamma+\alpha_{3}) \iiint(1+z)^{2\gamma+2\alpha_{3}-1} w  | D^{\alpha}\varphi|^{2}  \right|
 &\leq&  \|\frac{w}{1+z}\|_{L^{\infty}}  \|(1+z)^{\gamma+ \alpha_{3}}D^{\alpha}\varphi\|_{L^{2}}^{2} \nonumber \\
& \leq &  C_{s, \gamma,\sigma,\delta  }\left(\|\varphi\|_{H^{s,\gamma}_{g}}+\|\partial_{xy}^s U\|_{L^2(\mathbb{T}^2)}+1\right)\|\varphi\|_{H^{s,\gamma}_{g}}^2.  \label{6.2.6}
\end{eqnarray}
Last,  it follows that the first term of the third line on the right-hand side of equation $(\ref{6.2.2})$ has the following three cases,
\begin{eqnarray}
J_1=- \sum\limits_{0<\beta\leq \alpha}\binom{ \alpha}{\beta} \iiint(1+z)^{2\gamma+2\alpha_{3}}D^{\alpha}\varphi  D^{\beta} u\partial _{x}D^{\alpha-\beta}\varphi,
\end{eqnarray}
\begin{eqnarray}
J_2=- \sum\limits_{0<\beta\leq \alpha}\binom{ \alpha}{\beta} \iiint(1+z)^{2\gamma+2\alpha_{3}}D^{\alpha}\varphi  D^{\beta}(Ku)\partial _{y}D^{\alpha-\beta}\varphi,
\end{eqnarray}
and
\begin{eqnarray}
J_3=- \sum\limits_{0<\beta\leq \alpha}\binom{ \alpha}{\beta} \iiint(1+z)^{2\gamma+2\alpha_{3}}D^{\alpha}\varphi   D^{\beta} w\partial _{z}D^{\alpha-\beta}\varphi.
\end{eqnarray}
Here $J_1,J_2$ and $J_3$ can be estimated in the following.

For $J_{1}$, when $\beta_3 > 0$, we have
\begin{align}
|J_{1}|&\leq   \left|\sum\limits_{0<\beta\leq \alpha}\binom{ \alpha}{\beta} \iiint(1+z)^{2\gamma+2\alpha_{3}}D^{\alpha}\varphi  D^{\beta-e_3} \varphi D^{\alpha-\beta+e_1}\varphi\right| \nonumber \\
& \leq C \| (1+z)^{\gamma+\beta_{3}-1}  D^{\beta-e_3} \varphi\cdot (1+z)^{\gamma+\alpha_{3}-\beta_{3}}  D^{\alpha-\beta+e_1}\varphi\|_{L^2} \|(1+z)^{\gamma+\alpha_{3}}D^{\alpha}\|_{L^2}\nonumber \\
&\leq C\|(1+z)^{\gamma+\alpha_{3}}D^{\alpha}\|_{L^2}\nonumber\\
&\quad \times\left\{
\begin{aligned}
& \| (1+z)^{\gamma+\beta_{3}-1}  D^{\beta-e_3} \varphi\|_{L^\infty}\| (1+z)^{\gamma+\alpha_{3}-\beta_{3}}  D^{\alpha-\beta+e_1}\varphi\|_{L^2} , ~\beta \leq s-2,\\
& \| (1+z)^{\gamma+\beta_{3}-1}  D^{\beta-e_3} \varphi\|_{L^2}\| (1+z)^{\gamma+\alpha_{3}-\beta_{3}}  D^{\alpha-\beta+e_1}\varphi\|_{L^\infty} , ~\beta \geq s-3,\\
\end{aligned}
\right.\nonumber \\
& \leq C_{s, \gamma,\sigma,\delta  }\|\varphi \|_{H^{s,\gamma}_{g}}^{3}.
\label{6.2.7}
\end{align}
When $\beta_3 = 0$,  $\beta_1+ \beta_2 = s-1$, we get
\begin{align}
|J_{1}|&\leq   \left|\sum\limits_{0<\beta\leq \alpha}\binom{ \alpha}{\beta} \iiint(1+z)^{2\gamma+2\alpha_{3}}D^{\alpha}\varphi  \partial_{xy}^{s-1} u \partial_xD^{\alpha+1-s}\varphi\right| \nonumber \\
& \leq C \left(\|  \partial_{xy}^{s-1} (u-U)\|_{L^2}+ \| \partial_{xy}^{s-1} U \|_{L^2}\right) \|(1+z)^{\gamma+\alpha_{3}} \partial_xD^{\alpha+1-s}\varphi\|_{L^\infty}\|(1+z)^{\gamma+\alpha_{3}}D^{\alpha}\|_{L^2}\nonumber \\
& \leq C\left(\left\|\frac{1}{(1+z)^{\gamma-1}}\right\|_{L^\infty}\| (1+z)^{\gamma-1}  \partial_{xy}^{s-1} (u-U)\|_{L^2}+ \| \partial_{xy}^{s-1} U \|_{L^2}\right)\|\varphi \|_{H^{s,\gamma}_{g}}^{2} \nonumber \\
& \leq C_{s, \gamma,\sigma,\delta  }\left(\|\varphi\|_{H^{s,\gamma}_{g}}+\|\partial_{xy}^s U\|_{L^2(\mathbb{T}^2)}\right)\|\varphi \|_{H^{s,\gamma}_{g}}^{2}.
\end{align}
When $\beta_3 = 0$,  $\beta_1+ \beta_2 \leq s-2$, we derive
\begin{align}
|J_{1}|&\leq \left|\sum\limits_{0<\beta\leq \alpha}\binom{ \alpha}{\beta} \iiint(1+z)^{2\gamma+2\alpha_{3}}D^{\alpha}\varphi  \partial_{x}^{\beta_1}\partial_{y}^{\beta_2} u\partial_{x}^{\alpha_1-\beta_1+1}\partial_{y}^{\alpha_2-\beta_2}\partial_{z}^{\alpha_3}\varphi\right|\nonumber\\
& \leq C \|  \partial_{x}^{\beta_1}\partial_{y}^{\beta_2} u\|_{L^\infty} \|(1+z)^{\gamma+\alpha_{3}} \partial_{x}^{\alpha_1-\beta_1+1}\partial_{y}^{\alpha_2-\beta_2}\partial_{z}^{\alpha_3}\varphi\|_{L^2}\|(1+z)^{\gamma+\alpha_{3}}D^{\alpha}\|_{L^2}\nonumber \\
& \leq C_{s, \gamma,\sigma,\delta  }\left(\|\varphi\|_{H^{s,\gamma}_{g}}+\|\partial_{xy}^s U\|_{L^2(\mathbb{T}^2)}\right)\|\varphi \|_{H^{s,\gamma}_{g}}^{2}.
\end{align}

For $J_2$,
when $\beta_3 > 0$, we obtain
\begin{align}
|J_{2}|&\leq   \left|\sum\limits_{0<\beta\leq \alpha}\binom{ \alpha}{\beta} \iiint(1+z)^{2\gamma+2\alpha_{3}}D^{\alpha}\varphi  D^{\beta}(Ku)\partial _{y}D^{\alpha-\beta}\varphi\right| \nonumber \\
& \leq C \| (1+z)^{\gamma+\beta_{3}-1}  D^{\beta-e_3}(K\varphi) \cdot (1+z)^{\gamma+\alpha_{3}-\beta_{3}}  D^{\alpha-\beta+e_2}\varphi\|_{L^2} \|(1+z)^{\gamma+\alpha_{3}}D^{\alpha}\|_{L^2}\nonumber \\
&\leq C\|(1+z)^{\gamma+\alpha_{3}}D^{\alpha}\|_{L^2}\nonumber\\
&\quad\times\left\{
\begin{aligned}
& \| (1+z)^{\gamma+\beta_{3}-1}  D^{\beta-e_3}(K\varphi)\|_{L^\infty}\| (1+z)^{\gamma+\alpha_{3}-\beta_{3}}  D^{\alpha-\beta+e_2}\varphi\|_{L^2} , ~\beta \leq s-2,\\
& \| (1+z)^{\gamma+\beta_{3}-1}  D^{\beta-e_3}(K\varphi)\|_{L^2}\| (1+z)^{\gamma+\alpha_{3}-\beta_{3}}  D^{\alpha-\beta+e_2}\varphi\|_{L^\infty} , ~\beta \geq s-3,\\
\end{aligned}
\right.\nonumber \\
& \leq C_{s, \gamma,\sigma,\delta  }\|\partial_{xy}^s K\|_{L^{\infty}(\mathbb{T}^2)}\|\varphi \|_{H^{s,\gamma}_{g}}^{3}.
\end{align}
When $\beta_3 = 0$,  $\beta_1+ \beta_2 = s-1$,  we conclude
\begin{align}
|J_{2}|&\leq   \left|\sum\limits_{0<\beta\leq \alpha}\binom{ \alpha}{\beta} \iiint(1+z)^{2\gamma+2\alpha_{3}}D^{\alpha}\varphi  \partial^{s-1}(Ku)\partial _{y}D^{\alpha+1-s}\varphi\right| \nonumber \\
& \leq C \left(\|  \partial_{xy}^{s-1} K(u-U)\|_{L^2}+ \| \partial_{xy}^{s-1} KU \|_{L^2}\right) \|(1+z)^{\gamma+\alpha_{3}} \partial_yD^{\alpha+1-s}\varphi\|_{L^\infty}\|(1+z)^{\gamma+\alpha_{3}}D^{\alpha}\|_{L^2}\nonumber \\
& \leq C\left(\left\|\frac{1}{(1+z)^{\gamma-1}}\right\|_{L^\infty}\| (1+z)^{\gamma-1}  \partial_{xy}^{s-1} K(u-U)\|_{L^2}+ \| \partial_{xy}^{s-1} KU \|_{L^2}\right)\|\varphi \|_{H^{s,\gamma}_{g}}^{2} \nonumber \\
& \leq C_{s, \gamma,\sigma,\delta  }\|\partial_{xy}^s K\|_{L^{\infty}(\mathbb{T}^2)}\left(\|\varphi\|_{H^{s,\gamma}_{g}}+\|\partial_{xy}^s U\|_{L^2(\mathbb{T}^2)}\right)\|\varphi \|_{H^{s,\gamma}_{g}}^{2}.
\end{align}
When $\beta_3 = 0$,  $\beta_1+ \beta_2 \leq s-2$, we get
\begin{align}
|J_{2}|&\leq \left|\sum\limits_{0<\beta\leq \alpha}\binom{ \alpha}{\beta} \iiint(1+z)^{2\gamma+2\alpha_{3}}D^{\alpha}\varphi  \partial_{x}^{\beta_1}\partial_{y}^{\beta_2} (Ku)\partial_{x}^{\alpha_1-\beta_1}\partial_{y}^{\alpha_2-\beta_2+1}\partial_{z}^{\alpha_3}\varphi\right|\nonumber\\
& \leq C \|  \partial_{x}^{\beta_1}\partial_{y}^{\beta_2} (Ku)\|_{L^\infty} \|(1+z)^{\gamma+\alpha_{3}} \partial_{x}^{\alpha_1-\beta_1}\partial_{y}^{\alpha_2-\beta_2+1}\partial_{z}^{\alpha_3}\varphi\|_{L^2}\|(1+z)^{\gamma+\alpha_{3}}D^{\alpha}\|_{L^2}\nonumber \\
& \leq C_{s, \gamma,\sigma,\delta  }\|\partial_{xy}^s K\|_{L^{\infty}(\mathbb{T}^2)}\left(\|\varphi\|_{H^{s,\gamma}_{g}}+\|\partial_{xy}^s U\|_{L^2(\mathbb{T}^2)}\right)\|\varphi \|_{H^{s,\gamma}_{g}}^{2}.
\end{align}

For $J_3$, when $\beta_3 =0$, $\beta_1+ \beta_2 =s-2$ or  $s-1$, we achieve
\begin{align}
|J_{3}|&\leq\left|\sum\limits_{0<\beta\leq \alpha}\binom{ \alpha}{\beta} \iiint(1+z)^{2\gamma+2\alpha_{3}}D^{\alpha}\varphi   \partial_x^{\beta_1}\partial_y^{\beta_2} w\partial_x^{\alpha_1-\beta_1}\partial_{y}^{\alpha_2 - \beta_2}\partial_z^{\alpha_3+1}\varphi\right| \nonumber \\
& \leq C\|(1+z)^{\gamma+\alpha_{3}}D^{\alpha}\varphi\|_{L^2}\|(1+z)^{-1}(\partial_{x}^{\beta_1}\partial_{y}^{\beta_2} w + z\tilde{U})\|_{L^2}   \|(1+z)^{\gamma+\alpha_3+1} \partial_x^{\alpha_1-\beta_1}\partial_y^{\alpha_2 - \beta_2} \partial_z^{\alpha_3+1} \varphi\|_{L^\infty} \nonumber \\
& \quad +C\|(1+z)^{\gamma+\alpha_{3}}D^{\alpha}\varphi\|_{L^2}\|(1+z)^{-1} z\tilde{U}\|_{L^\infty}
\|(1+z)^{\gamma+\alpha_3+1} \partial_x^{\alpha_1-\beta_1}\partial_y^{\alpha_2 - \beta_2} \partial_z^{\alpha_3+1} \varphi\|_{L^2} \nonumber \\
&  \leq C_{s, \gamma,\sigma,\delta  }\left(\|\varphi\|_{H^{s,\gamma}_{g}}+\|\partial_{xy}^s U\|_{L^2(\mathbb{T}^2)}\right)\|\varphi \|_{H^{s,\gamma}_{g}}^{2}.
\end{align}
When $\beta_3 =0$, $\beta_1+ \beta_2 \leq s-3$, we deduce
\begin{align}
|J_{3}|&\leq\left|\sum\limits_{0<\beta\leq \alpha}\binom{ \alpha}{\beta} \iiint(1+z)^{2\gamma+2\alpha_{3}}D^{\alpha}\varphi   \partial_x^{\beta_1}\partial_y^{\beta_2} w\partial_x^{\alpha_1-\beta_1}\partial_{y}^{\alpha_2 - \beta_2}\partial_z^{\alpha_3+1}\varphi\right| \nonumber \\
& \leq C\|(1+z)^{\gamma+\alpha_{3}}D^{\alpha}\varphi\|_{L^2}\|(1+z)^{-1}\partial_x^{\beta_1}\partial_y^{\beta_2} w\|_{L^\infty}\|(1+z)^{\gamma+\alpha_3+1} \partial_x^{\alpha_1-\beta_1}\partial_y^{\alpha_2 - \beta_2} \partial_z^{\alpha_3+1} \varphi\|_{L^2}\nonumber \\
& \leq C_{s, \gamma,\sigma,\delta  }\left(\|\varphi\|_{H^{s,\gamma}_{g}}+\|\partial_{xy}^s U\|_{L^2(\mathbb{T}^2)}+1\right)\|\varphi \|_{H^{s,\gamma}_{g}}^{2}.
\end{align}
When $\beta_3 =1$, we obtain
\begin{align}
|J_{3}|&\leq\left|\sum\limits_{0<\beta\leq \alpha}\binom{ \alpha}{\beta} \iiint(1+z)^{2\gamma+2\alpha_{3}}D^{\alpha}\varphi   \partial_x^{\beta_1}\partial_y^{\beta_2} [\partial_x u +\partial_y(Ku)]\partial_x^{\alpha_1-\beta_1}\partial_{y}^{\alpha_2 - \beta_2}\partial_z^{\alpha_3}\varphi\right| \nonumber \\
& \leq C\|(1+z)^{\gamma+\alpha_{3}}D^{\alpha}\varphi\|_{L^2}\nonumber\\
&\quad \times\|(1+z)^{\gamma}\partial_x^{\beta_1}\partial_y^{\beta_2} [\partial_x u +\partial_y(Ku)]\cdot (1+z)^{\gamma+\alpha_3}\partial_x^{\alpha_1-\beta_1}\partial_{y}^{\alpha_2 - \beta_2}\partial_z^{\alpha_3}\varphi\|_{L^2}\nonumber \\
&\leq C\|(1+z)^{\gamma+\alpha_{3}}D^{\alpha}\|_{L^2}\nonumber\\
&\quad\times\left\{
\begin{aligned}
& \| (1+z)^{\gamma}\partial_x^{\beta_1}\partial_y^{\beta_2} [\partial_x u +\partial_y(Ku)]\|_{L^\infty}\| (1+z)^{\gamma+\alpha_3}\partial_x^{\alpha_1-\beta_1}\partial_{y}^{\alpha_2 - \beta_2}\partial_z^{\alpha_3}\varphi\|_{L^2} , ~\beta_{1}+\beta_{2} \leq s-4,\\
& \| (1+z)^{\gamma}\partial_x^{\beta_1}\partial_y^{\beta_2} [\partial_x u +\partial_y(Ku)]\|_{L^2}\| (1+z)^{\gamma+\alpha_3}\partial_x^{\alpha_1-\beta_1}\partial_{y}^{\alpha_2 - \beta_2}\partial_z^{\alpha_3}\varphi\|_{L^\infty} , ~s-3\leq \beta_{1}+\beta_{2} \leq s-1,\\
\end{aligned}
\right.\nonumber \\
&  \leq C_{s, \gamma,\sigma,\delta  }(1+\|\partial_{y}\partial_{xy}^{s-1} K\|_{L^\infty(\mathbb{T}^2)})\left(\|\varphi\|_{H^{s,\gamma}_{g}}+\|\partial_{xy}^s U\|_{L^2(\mathbb{T}^2)}+1\right)\|\varphi \|_{H^{s,\gamma}_{g}}^{2}.
\end{align}
When $\beta_3  \geq 2$, we get
\begin{align}
|J_{3}|&\leq\left|\sum\limits_{0<\beta\leq \alpha}\binom{ \alpha}{\beta} \iiint(1+z)^{2\gamma+2\alpha_{3}}D^{\alpha}\varphi   \partial_x^{\beta_1}\partial_y^{\beta_2} [\partial_x \varphi +\partial_yK\varphi +K \partial_y \varphi]\partial_x^{\alpha_1-\beta_1}\partial_{y}^{\alpha_2 - \beta_2}\partial_z^{\alpha_3-1}\varphi\right| \nonumber \\
& \leq C \|(1+z)^{\gamma+\alpha_{3}}D^{\alpha}\|_{L^2}\nonumber \\
&\quad \times \| (1+z)^{\gamma}    \partial_x^{\beta_1}\partial_y^{\beta_2} [\partial_x \varphi +\partial_yK\varphi +K \partial_y \varphi] \cdot (1+z)^{\gamma+\alpha_3-1}\partial_x^{\alpha_1-\beta_1}\partial_{y}^{\alpha_2 - \beta_2}\partial_z^{\alpha_3-1}\varphi\|_{L^2}\nonumber \\
&\leq C\|(1+z)^{\gamma+\alpha_{3}}D^{\alpha}\|_{L^2}\nonumber\\
&\quad\times\left\{
\begin{aligned}
& \| (1+z)^{\gamma}    \partial_x^{\beta_1}\partial_y^{\beta_2} [\partial_x \varphi +\partial_yK\varphi +K \partial_y \varphi]\|_{L^2}\| (1+z)^{\gamma+\alpha_3-1}\partial_x^{\alpha_1-\beta_1}\partial_{y}^{\alpha_2 - \beta_2}\partial_z^{\alpha_3-1}\varphi\|_{L^\infty} , ~\beta_1+\beta_2 \geq 2,\\
& \| (1+z)^{\gamma}    \partial_x^{\beta_1}\partial_y^{\beta_2} [\partial_x \varphi +\partial_yK\varphi +K \partial_y \varphi]\|_{L^\infty}\| (1+z)^{\gamma+\alpha_3-1}\partial_x^{\alpha_1-\beta_1}\partial_{y}^{\alpha_2 - \beta_2}\partial_z^{\alpha_3-1}\varphi\|_{L^2} , ~\beta_{1}+\beta_{2}=1,\\
\end{aligned}
\right.\nonumber \\
&  \leq C_{s, \gamma,\sigma,\delta  }(1+\|\partial_{y}\partial_{xy}^{s-1} K\|_{L^\infty(\mathbb{T}^2)})\left(\|\varphi\|_{H^{s,\gamma}_{g}}+\|\partial_{xy}^s U\|_{L^2(\mathbb{T}^2)}+1\right)\|\varphi \|_{H^{s,\gamma}_{g}}^{2}.
\end{align}
When $\beta_3  \geq 3$, we get
\begin{align}
|J_{3}|&\leq\left|\sum\limits_{0<\beta\leq \alpha}\binom{ \alpha}{\beta} \iiint(1+z)^{2\gamma+2\alpha_{3}}D^{\alpha}\varphi   \partial_x^{\beta_1}\partial_y^{\beta_2}\partial_z^{\beta_3-2} [\partial_x \varphi +\partial_yK\varphi +K \partial_y \varphi]\partial_x^{\alpha_1-\beta_1}\partial_{y}^{\alpha_2 - \beta_2}\partial_z^{\alpha_3-\beta_3+1}\varphi\right| \nonumber \\
& \leq C \|(1+z)^{\gamma+\alpha_{3}}D^{\alpha}\|_{L^2}\nonumber \\
&\quad \times \| (1+z)^{\gamma+\beta_{3}-2}    \partial_x^{\beta_1}\partial_y^{\beta_2}\partial_z^{\beta_3-2} [\partial_x \varphi +\partial_yK\varphi +K \partial_y \varphi] \cdot (1+z)^{\gamma+\alpha_3-\beta_{3}+1}\partial_x^{\alpha_1-\beta_1}\partial_{y}^{\alpha_2 - \beta_2}\partial_z^{\alpha_3-\beta_3+1}\varphi\|_{L^2}\nonumber \\
&\leq C\|(1+z)^{\gamma+\alpha_{3}}D^{\alpha}\|_{L^2}\nonumber\\
&\quad\times\left\{
\begin{aligned}
& \| (1+z)^{\gamma+\beta_{3}-2}    \partial_x^{\beta_1}\partial_y^{\beta_2}\partial_z^{\beta_3-2} [\partial_x \varphi +\partial_yK\varphi +K \partial_y \varphi]\|_{L^\infty} \nonumber\\
&\quad\cdot \| (1+z)^{\gamma+\alpha_3-\beta_{3}+1}\partial_x^{\alpha_1-\beta_1}\partial_{y}^{\alpha_2 - \beta_2}\partial_z^{\alpha_3-\beta_3+1}\varphi\|_{L^2} , ~\beta_{1}+\beta_{2} \leq s-5,\\
& \| (1+z)^{\gamma+\beta_{3}-2}    \partial_x^{\beta_1}\partial_y^{\beta_2}\partial_z^{\beta_3-2} [\partial_x \varphi +\partial_yK\varphi +K \partial_y \varphi]\|_{L^2} \nonumber\\
&\quad\cdot \| (1+z)^{\gamma+\alpha_3-\beta_{3}+1}\partial_x^{\alpha_1-\beta_1}\partial_{y}^{\alpha_2 - \beta_2}\partial_z^{\alpha_3-\beta_3+1}\varphi\|_{L^\infty} , ~\beta_{1}+\beta_{2}= s-4, s-3,\\
\end{aligned}
\right.\nonumber \\
&  \leq C_{s, \gamma,\sigma,\delta  }(1+\|\partial_{y}\partial_{xy}^{s-1} K\|_{L^\infty(\mathbb{T}^2)})\|\varphi \|_{H^{s,\gamma}_{g}}^{3}.
\end{align}
The estimate of $J_{4}$ is similar to that of $J_1$, then, when $\beta > 0$, we conclude
\begin{eqnarray}
|J_4| \leq C_{s, \gamma,\sigma,\delta  }\|\partial_{y} K\|_{L^\infty(\mathbb{T}^2)})\|\varphi \|_{H^{s,\gamma}_{g}}^{3}.
\end{eqnarray}
When $\beta_3 = 0$,  $\beta_1+ \beta_2 = s-1$, we get
\begin{eqnarray}
|J_4| \leq C_{s, \gamma,\sigma,\delta  }\|\partial_{y}\partial_{xy}^{s-1} K\|_{L^\infty(\mathbb{T}^2)}\left(\|\varphi\|_{H^{s,\gamma}_{g}}+\|\partial_{xy}^s U\|_{L^2(\mathbb{T}^2)}+1\right)\|\varphi \|_{H^{s,\gamma}_{g}}^{2}.
\end{eqnarray}
When $\beta_3 = 0$,  $\beta_1+ \beta_2 \leq s-1$, we achieve
\begin{eqnarray}
|J_4| \leq C_{s, \gamma,\sigma,\delta  }\|\partial_{y}\partial_{xy}^{s-1} K\|_{L^\infty(\mathbb{T}^2)}\left(\|\varphi\|_{H^{s,\gamma}_{g}}+\|\partial_{xy}^s U\|_{L^2(\mathbb{T}^2)}+1\right)\|\varphi \|_{H^{s,\gamma}_{g}}^{2}.
\label{6.2.12}
\end{eqnarray}
Hence, combining $(\ref{6.2.2})$, $(\ref{6.2.5})$-$(\ref{6.2.6})$ with $(\ref{6.2.7})$-$(\ref{6.2.12})$, we conclude
\begin{eqnarray}
\begin{aligned}
\frac{d}{dt}\|(1+z)^{ \gamma+ \alpha_{3}}D^{\alpha}\varphi\|_{L^{2}}^{2} &\leq C_{s, \gamma,\sigma,\delta  }\left(1+\|\partial_{xy}^{s} K\|_{L^\infty(\mathbb{T}^2)}\right)\left(\|\varphi\|_{H^{s,\gamma}_{g}}+\|\partial_{xy}^s U\|_{L^2(\mathbb{T}^2)}+1\right)\|\varphi \|_{H^{s,\gamma}_{g}}^{2} \\
& \quad +C_{s,\gamma,\sigma,\delta}\|\partial_{xy}^s K\|_{L^{\infty}(\mathbb{T}^2)}(1+\|w\|_{H_g^{s,\gamma}})^{s-2}\|\varphi\|^2_{H^{s,\gamma}_g}\\
&\quad +\sum\limits_{l=0}^{\frac{s}{2}}\|\partial^{l}_t\partial_xP\|_{H^{s-2l}(\mathbb{T}^2)}^{2}, \label{2.13}
\end{aligned}
\end{eqnarray}
which completes the proof.
 \hfill $\Box$

\subsection{Estimates on $g_{s}$ }
In this subsection, we will estimate the norm of $\partial_{x}^{i}\partial_{y}^{s-i}\varphi$ with weight  $(1+z)^{\gamma }$. To overcome the difficult term  $\partial_{z}^{s}w$, we use the usual cancellation property, i.e., we estimate the norm of $(1+z)^{\gamma}g_{s}$.
\begin{Proposition}\label{p6.2.2}
Under the same assumption of Proposition \ref{p6.2.1}, we have the following estimate:
\begin{eqnarray}
\begin{aligned}
\frac{d}{dt}\|(1+z)^{ \gamma}g_{s}\|_{L^{2}}^{2}& \leq  C_{s,\gamma, \sigma, \delta} \left(1+\|\partial_{xy}^s K \|_{L^\infty}\right)\left(1+\|\varphi\|_{H^{s,\gamma}_{g}}+\|\partial_{xy}^s U\|_{L^\infty(\mathbb{T}^2)}\right)\\
&\quad \times\left(\|\varphi\|_{H^{s,\gamma}_{g}}+\|\partial_{xy}^{s+1} U\|_{L^\infty(\mathbb{T}^2)}\right)\|\varphi\|_{H^{s,\gamma}_{g}}\\
&\quad+C_{\gamma, \delta }\left\{1+\|\partial_{xy}U\|^2_{L^\infty(\mathbb{T}^2)}\right\}\|\varphi\|^2_{H^{s,\gamma}_g}+C\|\partial_{xy}^{s}\partial_xP\|_{L^2(\mathbb{T}^2)}
.
\end{aligned}
\end{eqnarray}
\end{Proposition}
\textbf{Proof}.
Applying the operator $\partial_{xy}^{s}$ on $(\ref{6.1.5})_{1}$ and $(\ref{6.1.4})_{1}$ respectively, we obtain the following  equations
\begin{equation}\left\{
\begin{array}{ll}
(\partial _{t} +u\partial _{x} +Ku\partial _{y}+w\partial _{z}-\epsilon^2\partial_x^2 -\epsilon^2\partial^2_y -\partial _{z}^{2})\partial_{xy}^{s}\varphi +\partial_{xy}^{s}w\partial _{z}\varphi\\
=- \sum\limits_{0\leq j<s}\binom{s}{j} \partial _{xy}^{s-j}u\partial _{x}\partial _{xy}^{j }\varphi- \sum\limits_{0\leq j<s}\binom{s}{j} \partial _{xy}^{s-j}(Ku) \partial _{y}\partial _{xy}^{j }\varphi\\
\quad- \sum\limits_{1\leq j<s}\binom{s}{j} \partial _{x}^{s-j}w\partial _{z}\partial _{xy}^{j }\varphi + \sum\limits_{0\leq j \leq s}\binom{s}{j} \partial _{x}^{s-j}(\partial_y K u)\partial _{xy}^{j }\varphi   ,\\
(\partial _{t} +u\partial _{x} +Ku\partial _{y}+w\partial _{z}-\epsilon^2\partial _{x}^{2}-\epsilon^2\partial _{y}^{2}-\partial _{z}^{2})\partial_{xy}^{s}(u-U)+\partial_{xy}^{s}w\varphi \\
=- \sum\limits_{0\leq j<s}\binom{s}{j} \partial _{xy}^{s-j}u\partial _{x}\partial _{xy}^{j }(u-U)- \sum\limits_{0\leq j<s}\binom{s}{j} \partial _{xy}^{s-j}(Ku) \partial _{y}\partial _{xy}^{j }(u-U) \\
\quad- \sum\limits_{1\leq j<s}\binom{s}{j} \partial _{xy}^{s-j}w\partial _{xy}^{j }\varphi - \sum\limits_{0\leq j \leq s}\binom{s}{j} \partial _{xy}^{j}(u-U)\partial _{xy}^{s-j }\partial_x U\\
 \quad- \sum\limits_{0\leq j \leq s}\binom{s}{j} \partial _{xy}^{j}(u-U)\partial _{xy}^{s-j }(K\partial_y U ).
\end{array}
 \label{6.3.1}         \right.\end{equation}
Eliminating $\frac{\partial _{z}\varphi}{\varphi}\times (\ref{6.3.1})_{2}$ from $ (\ref{6.3.1})_{1}$, and letting
$a(t,x,y,z)=\frac{\partial _{z}\varphi}{\varphi}$, it follows
\begin{eqnarray}
\begin{aligned}
&
 (\partial _{t} +u\partial _{x} +Ku\partial _{y}+w\partial _{z}-\epsilon^2\partial_x^2-\epsilon^2\partial_y^2-\partial _{z}^{2})g_{s}\\
&  \quad+\partial_{xy}^{s}(u-U) (\partial _{t} +u\partial _{x} +Ku\partial _{y}+w\partial _{z}-\epsilon^2\partial_x^2-\epsilon^2\partial_y^2-\partial _{z}^{2})a    \\
&= 2\partial_{xy}^{s}\varphi\partial _{z}a+2\epsilon^2 \partial_x\partial^{s}_{xy} (u-U) \partial_x a+2\epsilon^2 \partial_y\partial^{s}_{xy} (u-U) \partial_y a \\
&\quad
- \sum\limits_{j=0}^{s-1}\binom{s}{j} (\partial _{xy}^{s-j}u(\partial _{x}\partial _{xy}^{j }\varphi-a\partial _{x}\partial _{xy}^{j }(u-U))+\partial _{xy}^{s-j}(Ku)(\partial _{y}\partial _{xy}^{j }\varphi-a\partial _{y}\partial _{xy}^{j }(u-U) ))  \\
&\quad - \sum\limits_{j=1}^{s-1}\binom{s}{j}\partial _{xy}^{s-j}w (\partial _{z}\partial _{xy}^{j }\varphi-a \partial _{xy}^{j }\varphi )+\sum\limits_{0\leq j \leq s}\binom{s}{j} \partial _{xy}^{j}(u-U)\partial _{xy}^{s-j }\partial_x U\\
&\quad+a\sum\limits_{0\leq j \leq s}\binom{s}{j} \partial _{xy}^{j}(u-U)\partial _{xy}^{s-j }(K\partial_y U)+ \sum\limits_{0\leq j \leq s}\binom{s}{j} \partial _{x}^{s-j}(\partial_y K u)\partial _{xy}^{j }\varphi . \label{6.3.2}
\end{aligned}
\end{eqnarray}
Applying the operator $\partial_{z}$ on $(\ref{6.1.5})_{1}$ yields
\begin{eqnarray}
\begin{aligned}
&(\partial _{t} +u\partial _{x} +Ku\partial _{y} +w\partial _{z})\partial _{z}\varphi \\
&=\epsilon^2(\partial_x^2+\partial_y^2)\partial_z \varphi+\partial _{z}^{3}\varphi-\varphi\partial _{x}\varphi-K\varphi\partial_{y}\varphi+\partial _{x}u\partial _{z}\varphi+\partial _{y}(Ku)\partial _{z}\varphi+\partial_y K(\varphi^2+u\partial_z \varphi).  \label{6.3.3}
\end{aligned}
\end{eqnarray}
Using $(\ref{6.3.3})$ and equation $(\ref{6.1.5})_{1}$, we get
\begin{align}
 &(\partial _{t} +u\partial _{x} +Ku\partial _{y}+w\partial _{z} )a \nonumber \\
 &=\frac{(\partial _{t} +u\partial _{x} +Ku\partial _{y} +w\partial _{z})\partial _{z}\varphi}{\varphi} -\frac{\partial_z \varphi(\partial _{t} +u\partial _{x} +Ku\partial _{y} +w\partial _{z})\varphi}{\varphi^2}\nonumber \\
&=\epsilon^2 \left\{\frac{(\partial_x^2+\partial_y^2)\partial_z \varphi}{\varphi}-a\frac{(\partial_x^2+\partial_y^2)\varphi}{\varphi}\right\}+\left\{\frac{\partial_z^3 \varphi}{\varphi}-a\frac{\partial_z^2 \varphi}{\varphi}\right\}-\partial_x \varphi-K\partial_y \varphi+a\partial_x u+ a\partial_y(Ku)\nonumber \\
&\quad+\partial_y K \varphi,  \label{6.3.4}
\end{align}
and
\begin{equation}\left\{
\begin{array}{ll}
  \partial _{z}a=\frac{ \partial _{z}^{2}\varphi}{\varphi} -\frac{\partial _{z}w \partial _{z}\varphi}{\varphi^{2}},   \\
 \partial _{z}^{2}a=\frac{\partial _{z}^{3}\varphi}{\varphi}-a\frac{\partial _{z}^{2}\varphi}{\varphi}- 2a \partial _{z}a, \\
 \partial _{x}^{2}a=\frac{\partial _{x}^{2}\partial_z\varphi}{\varphi}-a\frac{\partial _{x}^{2}\varphi}{\varphi}- 2\frac{\partial_x \varphi}{\varphi} \partial _{x}a,\\
 \partial _{y}^{2}a=\frac{\partial _{y}^{2}\partial_z\varphi}{\varphi}-a\frac{\partial _{y}^{2}\varphi}{\varphi}- 2\frac{\partial_y \varphi}{\varphi} \partial _{y}a.
\end{array}
 \label{6.3.5}         \right.\end{equation}
 Substituting (\ref{6.3.5}) into (\ref{6.3.4}), we obtain an equation for $a$:
 \begin{eqnarray}
\begin{aligned}
&(\partial _{t} +u\partial _{x} +Ku\partial _{y}+w\partial _{z}-\epsilon^2\partial_x^2-\epsilon^2\partial_y^2-\partial _{z}^{2} )a\\
&=2\epsilon^2(\frac{\partial_x \varphi}{\varphi} \partial _{x}a+\frac{\partial_y \varphi}{\varphi} \partial _{y}a)+2a \partial_z a-g_{1x}+a\partial_x U -Kg_{1y}+K a\partial_y U\\
&\quad+a\partial_y K u+\partial_y K \varphi.
\label{6.3.444}
\end{aligned}
\end{eqnarray}

Now inserting $(\ref{6.3.444})$ and $(\ref{6.3.5})_{2}$ into $(\ref{6.3.2})$ yields
\begin{eqnarray}
\begin{aligned}
&
 (\partial _{t} +u\partial _{x} +Ku\partial _{y}+w\partial _{z}-\epsilon^2\partial_x^2-\epsilon^2\partial_y^2-\partial _{z}^{2})g_{s}    \\
&= 2g_{s}\partial _{z}a+2\epsilon^2\partial_x a\left( \partial_x\partial^{s}_{xy} (u-U)-\frac{\partial_x \varphi}{\varphi}\partial^{s}_{xy} (u-U)\right) +2\epsilon^2\partial_y a\left( \partial_y\partial^{s}_{xy} (u-U)-\frac{\partial_y \varphi}{\varphi}\partial^{s}_{xy} (u-U)\right) \\
&\quad
- \sum\limits_{j=1}^{s-1}\binom{s}{j} g_{j+1}\partial _{xy}^{s-j}u-g_{x1}\partial_{xy}^s U- \sum\limits_{j=1}^{s-1}\binom{s}{j} g_{j+1}\partial _{xy}^{s-j}(Ku)-  g_{1y}\partial _{xy}^{s}(Ku)+g_{1y}(K\partial_{xy}^{s}(u-U))   \\
&\quad - \sum\limits_{j=1}^{s-1}\binom{s}{j}\partial _{xy}^{s-j}w (\partial _{z}\partial _{xy}^{j }\varphi-a \partial _{xy}^{j }\varphi )+a\sum\limits_{0\leq j \leq s-1}\binom{s}{j} \partial _{xy}^{j}(u-U)\partial _{xy}^{s-j }\partial_x U\\
&\quad+a\sum\limits_{0\leq j \leq s-1}\binom{s}{j} \partial _{xy}^{j}(u-U)\partial _{xy}^{s-j }(K\partial_y U)+ \sum\limits_{0\leq j \leq s}\binom{s}{j} \partial _{x}^{s-j}(\partial_y K u)\partial _{xy}^{j }\varphi \\
&\quad -( a\partial_y K u+\partial_y K \varphi)\partial_{xy}^{s}(u-U). \label{6.3.6}
\end{aligned}
\end{eqnarray}
 Multiplying $(\ref{6.3.6})$ by $(1+z)^{2\gamma}g_{s}$, and integrating it over $\mathbb{T}^{2}\times \mathbb{R}^{+}$, we arrive at
\begin{align}
&\frac{1}{2}\frac{d}{dt}\|(1+z)^{ \gamma}g_{s}\|_{L^{2}}^{2}  +\|(1+z)^{\gamma}\partial _{z}g_{s}\|_{L^{2}} ^{2}+\epsilon^2\|(1+z)^{\gamma}\partial _{x}g_{s}\|_{L^{2}} ^{2}+\epsilon^2\|(1+z)^{\gamma}\partial _{y}g_{s}\|_{L^{2}} ^{2} \nonumber \\
&=\iint_{\mathbb{T}^{2} } g_{s}\partial _{z}g_{s}dxdy | _{z=0}-2\gamma \iiint(1+z)^{2\gamma-1}g_{s}\partial _{z}  g_{s}
+\gamma \iiint(1+z)^{2\gamma-1} w  |g_{s}|^{2}  \nonumber\\
&\quad +2\iiint(1+z)^{2\gamma}|g_{s}|^{2}\partial _{z}a+2\epsilon^2\iiint(1+z)^{2\gamma}g_{s}\partial_x a\left( \partial_x\partial^{s}_{xy} (u-U)-\frac{\partial_x \varphi}{\varphi}\partial^{s}_{xy} (u-U)\right)\nonumber\\
&\quad +2\epsilon^2\iiint(1+z)^{2\gamma}g_s\partial_y a\left( \partial_y\partial^{s}_{xy} (u-U)-\frac{\partial_y \varphi}{\varphi}\partial^{s}_{xy} (u-U)\right) \nonumber\\
&\quad - \sum\limits_{j=1}^{s-1}\binom{s}{j} \iiint(1+z)^{2\gamma}g_s g_{j+1}\partial _{xy}^{s-j}u-\iiint(1+z)^{2\gamma}g_s g_{x1}\partial_{xy}^s U\nonumber\\
&\quad- \sum\limits_{j=1}^{s-1}\binom{s}{j} \iiint(1+z)^{2\gamma}g_sg_{j+1}\partial _{xy}^{s-j}(Ku)- \iiint(1+z)^{2\gamma}g_s g_{1y}\partial _{xy}^{s}(Ku)\nonumber\\
&\quad+\iiint(1+z)^{2\gamma}g_s g_{1y}(K\partial_{xy}^{s}(u-U))- \sum\limits_{j=1}^{s-1}\binom{s}{j}\iiint(1+z)^{2\gamma}g_s\partial _{xy}^{s-j}w (\partial _{z}\partial _{xy}^{j }\varphi-a \partial _{xy}^{j }\varphi )\nonumber\\
&\quad +\sum\limits_{0\leq j \leq s-1}\binom{s}{j} \iiint(1+z)^{2\gamma}g_s a\partial _{xy}^{j}(u-U)\partial _{xy}^{s-j }\partial_x U\nonumber\\
&\quad+\sum\limits_{0\leq j \leq s-1}\binom{s}{j}\iiint(1+z)^{2\gamma}g_s a \partial _{xy}^{j}(u-U)\partial _{xy}^{s-j }(K\partial_y U)\nonumber\\
&\quad+ \sum\limits_{0\leq j \leq s}\binom{s}{j}\iiint(1+z)^{2\gamma}g_s \partial _{x}^{s-j}(\partial_y K u)\partial _{xy}^{j }\varphi -\iiint(1+z)^{2\gamma}g_s( a\partial_y K u+\partial_y K \varphi)\partial_{xy}^{s}(u-U)
, \label{6.3.7}
\end{align}
which can be obtained by integration and by using the boundary condition, i.e.,
\begin{eqnarray}
\begin{aligned}
& \iiint(1+z)^{2\gamma}g_{s}(u\partial _{x} +Ku\partial _{y}+w\partial _{z} )g_{s} \\
&=\frac{1}{2} \iiint(1+z)^{2\gamma}\left\{u(|g_{s}|^{2})_{x}+ u(|g_{s}|^{2}) _{y}+w (|g_{s}|^{2}) _{z}\right\}
 =-\gamma \iiint(1+z)^{2\gamma-1}w |g_{s}|^{2},
 \label{6.3.8}
\end{aligned}
\end{eqnarray}
and
\begin{align}
 \iiint(1+z)^{2\gamma}g_{s}\partial _{z}^{2} g_{s}   &=-\|(1+z)^{\gamma}\partial _{z}g_{s}\|_{L^{2}}-\iint_{\mathbb{T}^{2}} g_{s}\partial _{z}g_{s}dxdy | _{z=0}
 -2\gamma \iiint(1+z)^{2\gamma-1}g_{s}\partial _{z}  g_{s}.
 \label{6.3.9}
\end{align}
By definition of $H^{s,\gamma}_{\sigma,\delta}$, we have
\begin{equation}\left\{
\begin{array}{ll}
|\varphi|^{-1}\leq \delta^{-1}(1+z)^{\sigma},   \\
| \partial _{x} \varphi|, |\partial _{y} \varphi| \leq \delta^{-1}(1+z)^{-\sigma },  \\
|\partial _{z}\varphi|, |\partial _{x} \partial _{z}\varphi|, |\partial _{y} \partial _{z}\varphi| \leq \delta^{-1}(1+z)^{-\sigma-1},\\
|\partial _{z}^{2}\varphi|  \leq \delta^{-1}(1+z)^{-\sigma-2}.
\end{array}
 \label{6.2.60}         \right.\end{equation}
We estimate equation $(\ref{6.3.7})$ by terms. Firstly,  we have the fact that
\begin{align}
\partial _{z}g_{s} \big{|} _{z=0}&=\left(\partial_{xy}^{s}\partial _{z}\varphi-\frac{\partial _{z}\varphi}{\varphi} \partial_{xy}^{s}\varphi -\frac{\partial^{2} _{z}\varphi}{\varphi} \partial_{xy}^{s}(u-U) +\frac{\partial_{z}\varphi\partial_{z}\varphi}{\varphi^{2}} \partial_{xy}^{s}(u-U) \right) \big{|} _{z=0} \nonumber\\
&=\partial_{xy}^{s}\partial _{x}P-ag_{s} \big{|} _{z=0} +\frac{\partial^{2} _{z}\varphi}{\varphi} \partial_{xy}^{s}U \big{|} _{z=0}.  \label{6.3.10}
\end{align}
Then according to the trace theorem, together with the facts that $\|a\|_{L^{\infty}}\leq \delta^{-2}$ and $\|\partial _{z}a\|_{L^{\infty}}\leq \delta^{-2}+\delta^{-4}$, we get
\begin{align}
&\left|\iint_{\mathbb{T}^{2} } g_{s}\partial _{z}g_{s}dxdy | _{z=0}  \right|\nonumber\\
&=\left|\iint_{\mathbb{T}^{2} } g_{s}| _{z=0}\partial_{xy}^{s}\partial _{x}Pdxdy   \right|+\left|\iint_{\mathbb{T}^{2} } a g_{s}^2| _{z=0}dxdy   \right|+\left|\iint_{\mathbb{T}^{2} } \frac{\partial^{2} _{z}\varphi}{\varphi} \partial_{xy}^{s}U g_{s}| _{z=0}dxdy   \right|\nonumber\\
&\leq  C\iint_{\mathbb{T}^{2} } (g_{s}+\partial_{z}g_{s})\partial_{xy}^{s}\partial _{x}Pdxdydz  +C\iint_{\mathbb{T}^{2} } (g_{s}+\partial_{z}g_{s})g_{s}dxdydz  +C\iint_{\mathbb{T}^{2} } \frac{\partial^{2} _{z}\varphi}{\varphi} \partial_{xy}^{s}U g_{s}dxdydz            \nonumber\\
&+C\iint_{\mathbb{T}^{2} } \partial_{z}\left(\frac{\partial^{2} _{z}\varphi}{\varphi}\right) \partial_{xy}^{s}U g_{s}dxdydz  +C\iint_{\mathbb{T}^{2} } \frac{\partial^{2} _{z}\varphi}{\varphi} \partial_{xy}^{s}U  \partial_{z}g_{s}dxdydz\nonumber\\
&\leq \frac{1}{4}\|(1+z)^{\gamma}\partial_{z} g_{s}\|_{L^{2}}^{2}+C_{\delta}\left\{1+\|\partial_{xy}U\|^2_{L^\infty(\mathbb{T}^2)}\right\}\|\varphi\|^2_{H^{s,\gamma}_g}+C\|\partial_{xy}^{s}\partial_xP\|_{L^2(\mathbb{T}^2)}^2.   \label{6.3.11}
\end{align}

Using the  H$\ddot{o}$lder inequality,  we obtain
\begin{eqnarray}
\left| 2\gamma \iiint(1+z)^{2\gamma-1}g_{s}\partial _{z}  g_{s}  \right|
 &\leq& 2\gamma \|\frac{1}{1+z }  \|_{L^{\infty}} \|(1+z)^{\gamma} g_{s}\|_{L^{2}} \|(1+z)^{\gamma}\partial _{z}g_{s}\|_{L^{2}} \nonumber \\
& \leq & C_{ \gamma } \|\varphi\|^2_{H^{s,\gamma}_g}+\frac{1}{4}\|(1+z)^{\gamma}\partial_{z} g_{s}\|_{L^{2}}^{2} ,  \label{6.3.12}
\end{eqnarray}
and using Lemma \ref{y6.4.4},
\begin{align}
\left|   \gamma \iiint(1+z)^{2\gamma-1} w  |g_{s}|^{2} \right| & \leq \|\frac{w}{1+z}\|_{L^{\infty}}\|(1+z)^{ \gamma}g_{s}\|_{L^{2}}^{2}\nonumber\\
&\leq  C_{s, \gamma,\sigma,\delta  }\left(1+\|\varphi \|_{H^{s,\gamma}_{g}}+\|\partial_{xy}U\|^2_{L^\infty(\mathbb{T}^2)}\right)\|\varphi \|_{H^{s,\gamma}_{g}} ^2.  \label{6.3.13}
\end{align}

Thus, it follows from to $(\ref{6.2.60})$ and $(\ref{6.3.5})$ that
\begin{eqnarray}
\left| 2\iiint(1+z)^{2\gamma}|g_{s}|^{2}\partial _{z}a  \right|
  \leq  2  \|\partial _{z}a  \|_{L^{\infty}} \|(1+z)^{\gamma} g_{s}\|_{L^{2}} ^{2}\leq C_{ \delta }  \|\varphi \|_{H^{s,\gamma}_{g}} ^2,
\label{6.3.15}
\end{eqnarray}
and by Lemma \ref{y6.4.2},
\begin{eqnarray}
&&\left| 2\epsilon^2\iiint(1+z)^{2\gamma}g_{s} \partial_x a\left( \partial_x\partial^{s}_{xy} (u-U)-\frac{\partial_x \varphi}{\varphi}\partial^{s}_{xy} (u-U)\right)   \right| \nonumber \\
&&\leq  2\epsilon^2C_{\delta} \|(1+z)^{\gamma} g_{s}\|_{L^{2}}\left(\|(1+z)^{\gamma-1}\partial_{xy}^s\partial_x(u-U)\|_{L^2}+\|(1+z)^{\gamma-1}\partial_{x}^s(u-U)\|_{L^2}\right),
\end{eqnarray}
which, together with $(\ref{6.2.60})$, implies the fact that $\|(1+z)(\partial _{x}a+\partial _{y}a)\|_{L^{\infty}} \leq \delta^{-2}+\delta^{-4}$.
In addition, by using Lemma \ref{y6.4.1} and  the following relation
\begin{eqnarray}
\varphi\partial_z\left(\frac{\partial_x\partial_{xy}^{s}(u-U)}{\varphi}\right)=g_{(s+1)x}=\partial_{x} g_s +\partial_x a\partial_{xy}^s(u-U),
\end{eqnarray}
we get
\begin{align}
&\|(1+z)^{\gamma-1}\ \partial_x\partial^{s}_{xy}(u-U)\|_{L^2} \nonumber\\
&\leq  \delta^{-1}\|(1+y)^{\gamma-\sigma-1}\frac{\partial_x\partial_{xy}^{s}(u-U)}{\varphi}\|_{L^2}\nonumber\\
&\leq C_{\gamma, \sigma, \delta}\left(\|\partial_{x}\partial_{xy}^{s}U\|_{L^{2}(\mathbb{T})}+\left\|(1+z)^{\gamma}\varphi\partial_z\left(\frac{\partial_x\partial_{xy}^{s}(u-U)}{\varphi}\right)\right\|_{L^2}\right)\nonumber\\
&\leq C_{\gamma, \sigma, \delta}\left(\|\partial_{x}\partial_{xy}^{s}U\|_{L^{2}(\mathbb{T}^2)}+\|(1+z)^{\gamma} \partial_x g_{s}\|_{L^2}+\|(1+z)\partial_x a\|_{L^{\infty}}\|(1+z)^{\gamma-1}\partial_{xy}^s(u-U)\|_{L^2}\right) .
\end{align}
Then, we can obtain
\begin{eqnarray}
&&\left| 2\epsilon^2\iiint(1+z)^{2\gamma}g_{s} \partial_x a\left( \partial_x\partial^{s}_{xy} (u-U)-\frac{\partial_x \varphi}{\varphi}\partial^{s}_{xy} (u-U)\right)   \right| \nonumber \\
&&\leq  \epsilon^2C_{\gamma, \sigma, \delta} \left(\|\varphi\|_{H^{s,\gamma}_{g}}+\|\partial_{xy}^s U\|_{L^2(\mathbb{T}^2)}+\|\partial_{x}\partial_{xy}^{s}U\|_{L^{2}(\mathbb{T}^2)}+\|(1+z)^{\gamma} \partial_x g_{s}\|_{L^2}\right)\|\varphi\|_{H^{s,\gamma}_{g}}.
\label{6.3.1555}
\end{eqnarray}
Similarly, we also have
\begin{eqnarray}
&&\left| 2\epsilon^2\iiint(1+z)^{2\gamma}g_{s} \partial_y a\left( \partial_y\partial^{s}_{xy} (u-U)-\frac{\partial_y \varphi}{\varphi}\partial^{s}_{xy} (u-U)\right)   \right| \nonumber \\
&&\leq  \epsilon^2C_{\gamma, \sigma, \delta} \left(\|\varphi\|_{H^{s,\gamma}_{g}}+\|\partial_{xy}^s U\|_{L^2(\mathbb{T}^2)}+\|\partial_{y}\partial_{xy}^{s}U\|_{L^{2}(\mathbb{T}^2)}+\|(1+z)^{\gamma} \partial_y g_{s}\|_{L^2}\right)\|\varphi\|_{H^{s,\gamma}_{g}}.
\label{6.3.15555}
\end{eqnarray}
Now we estimate some of the remaining complex terms in (\ref{6.3.7}). Thus we need to discuss the classification to complete the estimation. For $i \leq s-3$, we use the usual inequality to derive
\begin{align}
\|(1+z)^{\gamma}g_i\|_{L^\infty} &\leq \|(1+z)^{\gamma}\partial^i_{xy}\varphi\|_{L^\infty}+\|(1+z)^{\gamma-1}a\|_{L^\infty}\|\partial_{xy}^i (u-U)\|_{L^\infty} \nonumber \\
&\leq C_{\delta} \left(\|(1+z)^{\gamma+\alpha_3}\partial^i_{xy}\varphi\|_{L^\infty}+\|\partial_{xy}^i u\|_{L^\infty}+\|\partial_{xy}^i U\|_{L^\infty}\right)\nonumber \\
&\leq C_{s,\gamma, \sigma, \delta}\left(\|\varphi\|_{H^{s,\gamma}_{g}}+\|\partial_{xy}^s U\|_{L^2(\mathbb{T}^2)}\right).
\end{align}
By Lemma \ref{y6.4.4}, we obtain
\begin{align}
\left| \iiint(1+z)^{2\gamma}g_s g_{j+1}\partial _{xy}^{s-j}u\right|
&\leq \left\{
\begin{aligned}
&\|(1+z)^{\gamma}  g_{s}\|_{L^2}\|(1+z)^{\gamma}  g_{j+1}\|_{L^2}\|\partial _{xy}^{s-j}u\|_{L^\infty}, ~j=[2, s-1],\\
&\|(1+z)^{\gamma}  g_{s}\|_{L^2}\|(1+z)^{\gamma}  g_2\|_{L^\infty}\|\partial _{xy}^{s-1}(u-U)+\partial _{xy}^{s-1}U\|_{L^2}, ~j=1.\\
\end{aligned}
\right.\nonumber \\
&\leq C_{s,\gamma, \sigma, \delta}\left(\|\varphi\|_{H^{s,\gamma}_{g}}+\|\partial_{xy}^s U\|_{L^2(\mathbb{T}^2)}\right)^2\|\varphi\|_{H^{s,\gamma}_{g}},
\label{6.3.155}
\end{align}
for $j=1,2,\cdots, s-1$, and
\begin{eqnarray}
\left|\iiint(1+z)^{2\gamma}g_s g_{x1}\partial_{xy}^s U\right| \leq C_{s,\gamma, \sigma, \delta} \left(\|\varphi\|_{H^{s,\gamma}_{g}}+\|\partial_{xy}^s U\|_{L^2(\mathbb{T}^2)}\right)\|\partial_{xy}^s U\|_{L^2(\mathbb{T}^\infty)}\|\varphi\|_{H^{s,\gamma}_{g}}.
\end{eqnarray}
Similarly,
\begin{eqnarray}
\begin{aligned}
&\left| \sum\limits_{j=1}^{s-1}\binom{s}{j} \iiint(1+z)^{2\gamma}g_sg_{j+1}\partial _{xy}^{s-j}(Ku)\right|\leq C_{s,\gamma, \sigma, \delta}\|\partial_{xy}^{s-1} K\|_{L^\infty}\left(\|\varphi\|_{H^{s,\gamma}_{g}}+\|\partial_{xy}^s U\|_{L^2(\mathbb{T}^2)}\right)^2\|\varphi\|_{H^{s,\gamma}_{g}}
\end{aligned}
\end{eqnarray}
and
\begin{align}
&\left|- \iiint(1+z)^{2\gamma}g_s g_{1y}\partial _{xy}^{s}(Ku)+\iiint(1+z)^{2\gamma}g_s g_{1y}(K\partial_{xy}^{s}(u-U))\right|\nonumber\\
&\leq C_{s,\gamma, \sigma, \delta}\|\partial_{xy}^{s-1} K\|_{L^\infty}\|(1+z)^{\gamma}  g_{s}\|_{L^2}\|(1+z)^{\gamma}  g_{1y}\|_{L^\infty}\|-\partial _{xy}^{s}(Ku)+(K\partial_{xy}^{s}(u-U))\|_{L^2}\nonumber\\
& \leq C_{s,\gamma, \sigma, \delta}\|\partial_{xy}^{s-1} K\|_{L^\infty}\left(\|\varphi\|_{H^{s,\gamma}_{g}}+\|\partial_{xy}^s U\|_{L^2(\mathbb{T}^2)}\right)^2\|\varphi\|_{H^{s,\gamma}_{g}}.
\label{6.3.017}
\end{align}
It should be noted that  the highest order of $u$ will not reach $s$, due to cancelling each other out,  in the   inequality  (\ref{6.3.017}).\\
Last, for $j=3,\cdots, s-1$,
\begin{align}
&\left|   \iiint(1+z)^{2\gamma}g_s\partial _{xy}^{s-j}w (\partial _{z}\partial _{xy}^{j }\varphi-a \partial _{xy}^{j }\varphi )  \right| \nonumber\\
&  \leq\left\|\frac{\partial_{xy}^{s-j}w}{1+z}\right\|_{L^\infty} \left(\|(1+z)^{\gamma+1}\partial_{xy}^{j}\partial_{z}\varphi\|_{L^2}+\|(1+z)a\|_{L^2}\|(1+z)^{\gamma}\partial_{xy}^{j}\varphi
\|_{L^2}   \right)\|(1+z)^{\gamma}g_{s}\|_{L^2}                                                                                                            \nonumber\\
&\leq C_{s,\gamma, \sigma, \delta}\left(\|\varphi\|_{H^{s,\gamma}_{g}}+\|\partial_{xy}^s U\|_{L^2(\mathbb{T}^2)}\right)\|\varphi\|_{H^{s,\gamma}_{g}}^2
\label{6.3.17}
\end{align}
and for $j=1, 2$,
\begin{align}
&\left|   \iiint(1+z)^{2\gamma}g_s\partial _{xy}^{s-j}w (\partial _{z}\partial _{xy}^{j}\varphi-a \partial _{xy}^{j}\varphi )  \right| \nonumber\\
&\leq\left\|\frac{\partial^{s-j}_{xy}w+z\tilde{U}}{1+z}\right\|_{L^2}
\left(\|(1+z)^{\gamma+1}\partial_{xy}\partial_{z}\varphi\|_{L^\infty}+\|(1+z)a\|_{L^\infty}\|(1+z)^{\gamma}\partial_{xy}^{j}\varphi
\|_{L^\infty}   \right)\|(1+z)^{\gamma}g_{s}\|_{L^2}\nonumber\\
&\quad+\|\tilde{U}\|_{L^\infty(\mathbb{T}^2)}
\left(\|(1+z)^{\gamma+1}\partial_{xy}\partial_{z}\varphi\|_{L^2}+\|(1+z)a\|_{L^\infty}\|(1+z)^{\gamma}\partial_{xy}^{j}\varphi
\|_{L^2}   \right)\|(1+z)^{\gamma}g_{s}\|_{L^2}\nonumber\\
&\leq C_{s,\gamma, \sigma, \delta}\left(\|\varphi\|_{H^{s,\gamma}_{g}}+\|\partial_{xy}^s U\|_{L^2(\mathbb{T}^2)}\right)\|\varphi\|_{H^{s,\gamma}_{g}}^2.
\label{6.3.18}
\end{align}
For $j=1,2,\cdots, s-1$, by Lemma \ref{y6.4.4} and using the fact  $\|(1+z) a \|_{L^{\infty}} \leq \delta^{-2}$, it follows
\begin{eqnarray}
&& \left|    \iiint(1+z)^{2\gamma}g_s a\partial _{xy}^{j}(u-U)\partial _{xy}^{s-j }\partial_x U  \right|
 \nonumber \\
&& \leq    C \|(1+z)^{\gamma}g_s\|_{L^2}\|(1+z)a\|_{L^\infty}\|(1+z)^{\gamma-1}\partial_{xy}^{j}(u-U)\|_{L^2}\|\partial _{x}\partial _{xy}^{s-j}U\|_{L^\infty(\mathbb{T}^2)}  \nonumber \\
&& \leq C_{s,\gamma, \sigma, \delta}\|\partial _{x}\partial _{xy}^{s}U\|_{L^\infty(\mathbb{T}^2)} \left(\|\varphi\|_{H^{s,\gamma}_{g}}+\|\partial_{xy}^s U\|_{L^2(\mathbb{T}^2)}\right)\|\varphi\|_{H^{s,\gamma}_{g}},
\label{6.3.19}
\end{eqnarray}
and
\begin{eqnarray}
&& \left|    \iiint(1+z)^{2\gamma}g_s a \partial _{xy}^{j}(u-U)\partial _{xy}^{s-j }(K\partial_y U)  \right|
 \nonumber \\
&& \leq C_{s,\gamma, \sigma, \delta}\|\partial _{xy}^{s}(K\partial_yU)\|_{L^\infty(\mathbb{T}^2)} \left(\|\varphi\|_{H^{s,\gamma}_{g}}+\|\partial_{xy}^s U\|_{L^2(\mathbb{T}^2)}\right)\|\varphi\|_{H^{s,\gamma}_{g}}.
\label{6.3.199}
\end{eqnarray}
For $j=1,2,\cdots, s$, we achieve
\begin{align}
\left| \iiint(1+z)^{2\gamma}g_s \partial _{x}^{s-j}(\partial_y K u)\partial _{xy}^{j }\varphi\right|
&\leq C_{\delta}\|\varphi\|_{H^{s,\gamma}_{g}}\left\{
\begin{aligned}
&\|(1+z)^{\gamma} \partial _{xy}^{j }\varphi\|_{L^2}\|\partial _{x}^{s-j}(\partial_y K u)\|_{L^\infty}, ~j \geq 2,\\
&\|(1+z)^{\gamma} \partial _{xy}^{j }\varphi\|_{L^\infty}\|\partial _{x}^{s-j}(\partial_y K u)\|_{L^2}, ~j=0,1.\\
\end{aligned}
\right.\nonumber \\
&\leq C_{s,\gamma, \sigma, \delta}\|\partial _{xy}^{s-2}K\|_{L^\infty(\mathbb{T}^2)}\left(\|\varphi\|_{H^{s,\gamma}_{g}}+\|\partial_{xy}^s U\|_{L^2(\mathbb{T}^2)}\right)^2\|\varphi\|_{H^{s,\gamma}_{g}},\label{6.3.1999}
\end{align}
and
\begin{align}
-\iiint(1+z)^{2\gamma}g_s( a\partial_y K u+\partial_y K \varphi)\partial_{xy}^{s}(u-U) &\leq  C\|\partial_y K \|_{L^\infty}\|(1+z)^{\gamma}g_s\|_{L^2}\|(1+z)^{\gamma-1}\partial_{xy}^{s}(u-U)\|_{L^2}
\nonumber \\
&\quad \times (\|(1+z)a\|_{L^\infty}\|u\|_{L^\infty}+\|(1+z)^{1-\gamma}\|_{L^\infty}\|(1+z)^{\gamma}\varphi\|_{L^\infty})\nonumber\\
&\leq C_{s,\gamma, \sigma, \delta} \|\partial_y K \|_{L^\infty}\left(\|\varphi\|_{H^{s,\gamma}_{g}}+\|\partial_{xy}^s U\|_{L^2(\mathbb{T}^2)}\right)^2\|\varphi\|_{H^{s,\gamma}_{g}}.
\label{6.3.20}
\end{align}

Hence combining (\ref{6.3.11})-(\ref{6.3.15}), (\ref{6.3.1555})-(\ref{6.3.15555}), (\ref{6.3.155})-(\ref{6.3.20}) with (\ref{6.3.7}) leads to
\begin{eqnarray}
\begin{aligned}
\frac{d}{dt}\|(1+z)^{ \gamma}g_{s}\|_{L^{2}}^{2}& \leq  C_{s,\gamma, \sigma, \delta}\left(1+ \|\partial_{xy}^s K \|_{L^\infty}\right)\left(1+\|\varphi\|_{H^{s,\gamma}_{g}}+\|\partial_{xy}^s U\|_{L^\infty(\mathbb{T}^2)}\right)\\
&\quad \times\left(\|\varphi\|_{H^{s,\gamma}_{g}}+\|\partial_{xy}^{s+1} U\|_{L^\infty(\mathbb{T}^2)}\right)\|\varphi\|_{H^{s,\gamma}_{g}}\\
&\quad+C_{\gamma, \delta }\left\{1+\|\partial_{xy}U\|^2_{L^\infty(\mathbb{T}^2)}\right\}\|\varphi\|^2_{H^{s,\gamma}_g}+C\|\partial_{xy}^{s}\partial_xP\|_{L^2(\mathbb{T}^2)}
. \label{6.3.21}
\end{aligned}
\end{eqnarray}
 \hfill $\Box$

\subsection{Weighted $H^s$ estimate on $\varphi$ }
Now we can derive the weighted $H^s$ estimate on $\varphi$ by employing Propositions \ref{p6.2.1}-\ref{p6.2.2}.
\begin{Proposition}\label{p6.2.3}
Under the same assumption of Proposition \ref{p6.2.1}, we have the following estimate
\begin{eqnarray}
\begin{aligned}
\|\varphi\|_{H^{s,\gamma}_{g}}^2 &\leq \left\{\|\varphi_{0}\|_{H^{s,\gamma}_{g}}^2+\int_0^t F(\tau)d \tau\right\}\\
&\quad \times \left\{1-C(s-1)\left(\|\varphi_{0}\|_{H^{s,\gamma}_{g}}^2+\int_0^t F(\tau) d\tau\right)^{s-1}t\right\}^{-\frac{1}{s-1}},
\label{6.3.2222}
\end{aligned}
\end{eqnarray}
where $C>0$ is a constant independent of $\epsilon$ and $t$. The function $F(t)$ is expressed by
\begin{eqnarray}
\begin{aligned}
F(t)=\mathcal{P}(\|\partial_{xy}^{s+1}U\|_{{L^\infty}(\mathbb{T}^2)},
\|\partial_{xy}^{s}K\|_{{L^\infty}(\mathbb{T}^2)})+C \sum\limits_{l=0}^{\frac{s}{2}}\|\partial^{l}_t\partial_xP\|_{H^{s-2l}(\mathbb{T}^2)}^{2},
\label{MK}
\end{aligned}
\end{eqnarray}
and $\mathcal{P}(\cdot, \cdot)$ denotes a polynomial.
\end{Proposition}
\textbf{Proof}.
According to Propositions \ref{p6.2.1}-\ref{p6.2.2}, we deduce from the definition of $\|\cdot\|_{H^{s,\gamma}_g}$ that
\begin{align}
\frac{d}{dt} \|\varphi\|_{H^{s,\gamma}_{g}}^2 &\leq C_{s,\gamma, \sigma, \delta} \|\varphi\|_{H^{s,\gamma}_{g}}^{2s}+\sum\limits_{l=0}^{\frac{s}{2}}\|\partial^{l}_t\partial_xP\|_{H^{s-2l}(\mathbb{T}^2)}^{2}+\mathcal{P}(\|\partial_{xy}^{s+1}U\|_{{L^\infty}(\mathbb{T}^2)},
\|\partial_{xy}^{s}K\|_{{L^\infty}(\mathbb{T}^2)})
\end{align}
and hence, it follows by using the comparison principle of ordinary differential equations that
\begin{eqnarray}
\begin{aligned}
\|\varphi\|_{H^{s,\gamma}_{g}}^2 &\leq \left\{\|\varphi_0\|_{H^{s,\gamma}_{g}}^2+\int_0^t F(\tau)d \tau\right\}\\
&\quad \times \left\{1-C(s-1)\left(\|\varphi_0\|_{H^{s,\gamma}_{g}}^2+\int_0^t F(\tau) d\tau\right)^{s-1}t\right\}^{-\frac{1}{s-1}},
\end{aligned}
\end{eqnarray}
provided that
\begin{eqnarray}
\begin{aligned}
1-C(s-1)\left(\|\varphi_0\|_{H^{s,\gamma}_{g}}^2+\int_0^t F(\tau) d\tau\right)^{s-1}t > 0.
\end{aligned}
\end{eqnarray}
This proves Proposition \ref{p6.2.3}.
\hfill $\Box$

\subsection{Weighted $L^{\infty}$ estimates on lower order terms  }
In this subsection, we will estimate the weighted $L^{\infty}$ on $D^{\alpha} \varphi$ for $|\alpha|\leq 2$ by using the classical maximum principles. More precisely, we will derive two parts: an $L^\infty$-estimate on  $I:=\sum \limits_{|\alpha|\leq 2} \left|(1+z)^{\sigma+\alpha_{3}}D^{\alpha}\varphi\right|^2$ and a lower bound estimate on $B_{(0,0,0)}:=(1+z)^{\sigma} \varphi$. Here
\begin{eqnarray*}
\begin{aligned}
I:=\sum \limits_{|\alpha|\leq 2} \left|(1+z)^{\sigma+\alpha_{3}}D^{\alpha}\varphi\right|^2
\end{aligned},\;\;B_{\alpha}:=(1+z)^{\sigma+\alpha_{3}}D^{\alpha}\varphi,
\end{eqnarray*}
and $B_{(0,0, 0)}=(1+z)^{\sigma}\varphi$.
\begin{Lemma}\label{y6.2.3}
Under the same assumption of Proposition \ref{p6.2.1}, we have the following estimates. \\
For any $s \geq 5$,
\begin{eqnarray}
\begin{aligned}
\|I(t)\|_{L^\infty(\mathbb{T}^2\times\mathbb{R}^+)}\leq \max \left\{           \|I(0)\|_{L^\infty(\mathbb{T}^2\times\mathbb{R}^+)},6C^2W(t)^2
\right\} e^{C\left(1+G_1(t)\right)t},
\label{6.3.31}
\end{aligned}
\end{eqnarray}
and for any $s\geq 7$,
\begin{eqnarray}
\begin{aligned}
\|I(t)\|_{L^\infty(\mathbb{T}^2\times\mathbb{R}^+)}\leq  \left\{           \|I(0)\|_{L^\infty(\mathbb{T}^2\times\mathbb{R}^+)}+C[1+Y(t)W(t)]W(t)^2t
\right\} e^{C\left(1+G_1(t)\right)t}.
\label{6.3.32}
\end{aligned}
\end{eqnarray}
In addition, if $s \geq 5$, we also have,
\begin{eqnarray}
\begin{aligned}
\mathop{\min}\limits_{\mathbb{T}^2\times\mathbb{R}^+}(1+z)^{\sigma}\varphi(t)\geq \max \left\{         1-C\left(1+G_2(t)\right)te^{C\left(1+G_2(t)\right)t}
\right\}\cdot\min \left\{\mathop{\min}\limits_{\mathbb{T}^2\times\mathbb{R}^+}(1+z)^{\sigma}\varphi_0-CW(t)t
\right\},
\label{6.3.33}
\end{aligned}
\end{eqnarray}
where constant $C>0$ depends on $s, \gamma,\sigma$, and $\delta$ only. The functions $G_1$, $G_2$, W and Y : $[0,T] \rightarrow \mathbb{R}^+$ are respectively defined by
\begin{eqnarray}
\begin{aligned}
G_1(t):&=\mathop{\sup}\limits_{ [0,t]}\|\partial_{xy}^{2}K(\tau)\|_{{L^{\infty}}(\mathbb{T}^2)}+\mathop{\sup}\limits_{ [0,t]}\|w(\tau)\|_{H^{s,\gamma}_{g}}+\mathop{\sup}\limits_{ [0,t]}\|\partial_{xy}^{s}U(\tau)\|_{{L^2}(\mathbb{T}^2)}\\
&\quad+\mathop{\sup}\limits_{ [0,t]}\|\partial_{xy}^{2}K(\tau)\|_{{L^{\infty}}(\mathbb{T}^2)}\left(\mathop{\sup}\limits_{ [0,t]}\|w(\tau)\|_{H^{s,\gamma}_{g}}+\mathop{\sup}\limits_{ [0,t]}\|\partial_{xy}^{s}U(\tau)\|_{{L^2}(\mathbb{T}^2)}\right) ,\\
G_2(t):&=\mathop{\sup}\limits_{ [0,t]}\|w(\tau)\|_{H^{s,\gamma}_{g}}+\mathop{\sup}\limits_{ [0,t]}\|\partial_{xy}^{s}U(\tau)\|_{{L^2}(\mathbb{T}^2)},
\label{6.2.012}
\end{aligned}
\end{eqnarray}
and
\begin{eqnarray}
  \ W(t):=\mathop{\sup}\limits_{ [0,t]}\|w(\tau)\|_{H^{s,\gamma}_{g}}, \ \ \ \ Y(t):=\mathop{\sup}\limits_{ [0,t]}\|\partial_{xy}K(\tau)\|_{{L^{\infty}}(\mathbb{T}^2)}.
\label{6.2.0122}
\end{eqnarray}
\end{Lemma}
\textbf{Proof}.
This proof will be presented in two distinct steps, making use of the classical maximum principle for parabolic equations. In the first step, we will establish weighted $L^\infty$ controls on  $D^\alpha \varphi$, employing the maximum principle outlined in Lemma \ref{y6.4.5}. These controls will depend on the boundary values of  $D^\alpha \varphi$ at $z=0$. In the second step, we will obtain estimates on the boundary values of $D^\alpha \varphi$, utilizing either Sobolev embedding or growth rate control argument. \\
$\textbf{Step 1.}$

By direct computations, we obtain
\begin{eqnarray}
\begin{aligned}
\left(\partial_t+u\partial_x+Ku\partial_y+w\partial_z-\partial^2_z-\epsilon^2\partial^2_x-\epsilon^2\partial^2_y-\partial_y K u\right)B_{\alpha}=\sum\limits_{i=1}^{3} S_{i},
\label{6.3.34}
\end{aligned}
\end{eqnarray}
where
\begin{eqnarray*}
\begin{aligned}
S_1=\left(\frac{\sigma+\alpha_3}{1+z}w+\frac{(\sigma+\alpha_3)(\sigma+\alpha_3-1)}{(1+z)^2}\right)B_{\alpha},\ \ S_2=-\frac{2(\sigma+\alpha_3)}{1+z}\partial_z B_{\alpha},
\end{aligned}
\end{eqnarray*}
and
\begin{eqnarray*}
\begin{aligned}
S_3=-\sum\limits_{0< \beta \leq \alpha}\binom{\alpha}{\beta}\left\{(1+z)^{\beta_3}\left(D^{\beta}u B_{\alpha-\beta+e_1}+D^{\beta}(Ku) B_{\alpha-\beta+e_2}+\frac{D^{\beta} w B_{\alpha-\beta+e_3}}{1+z}+D^{\beta}(\partial_{y} Ku) B_{\alpha-\beta}\right)\right\}.
\end{aligned}
\end{eqnarray*}
Here, $e_1 := (1,0 , 0)$, $e_2 := (0,1, 0)$ and $e_3 := (0,0, 1)$. Using Lemma \ref{y6.4.4},
we have the following pointwise controls on $Q_\alpha,
R_\alpha$ and $S_\alpha$: for $|\alpha| \leq 2$,
\begin{equation}
\left\{\begin{aligned}
|S_{1}|  &\leq C_{s, \gamma,\sigma,\delta  } \left(1 + \|\varphi\|_{H_{g}^{s,\gamma}} + \|\partial_{xy}^s U\|_{L^2 (\mathbb{T}^2)} \right),
 \quad   |S_{2}|  \le C_\sigma, \\
    |S_3| & \leq C_{s, \gamma,\sigma,\delta  }(1+\|\partial_{xy}^{2}K\|_{{L^{\infty}}(\mathbb{T}^2)}) \left( \|\varphi\|_{H_{g}^{s,\gamma}} + \|\partial_{xy}^s U\|_{L^2 (\mathbb{T}^2)}\right)
    \sum_{0 < \beta \le \alpha}  |B_{\alpha-\beta+e_i}|, \ \ i=0,1,2,3,
\end{aligned}\right. \label{6.11.1}
\end{equation}
where $C_\sigma$ and $C_{s, \gamma,\sigma,\delta  }$ are some universal constants which are
independent of the solution $\varphi$.

By (\ref{6.11.1}), multiplying the equation (\ref{6.3.34}) by $2B_{\alpha}$  gives
\begin{eqnarray}
&&\left(\partial_t+u\partial_x+Ku\partial_y+w\partial_z-\partial^2_z-\epsilon^2\partial^2_x-\epsilon^2\partial^2_y-2\partial_y K u\right)I\nonumber\\
&&=2\sum \limits_{|\alpha|\leq 2} \left(|S_1 B_{\alpha}|+|S_2 B_{\alpha}|+|S_3 B_{\alpha}|-\epsilon^2|\partial_x B_{\alpha}|^2-\epsilon^2|\partial_y B_{\alpha}|^2-|\partial_z B_{\alpha}|^2\right)\nonumber\\
&&\leq  C_{s, \gamma,\sigma,\delta  }(1+\|\partial_{xy}^{2}K\|_{{L^{\infty}}(\mathbb{T}^2)})( 1+\|\varphi(s)\|_{H^{s,\gamma}_{g}}+\|\partial_{xy}^{s}U\|_{{L^2}(\mathbb{T}^2)})I,
\label{parabolic1.1}
\end{eqnarray}
where we have used Lemma \ref{y6.4.4}.

Applying the classical maximum principle Lemma \ref{y6.4.5} for parabolic equations (\ref{parabolic1.1}), we have
\begin{eqnarray}
\begin{aligned}
&\|I(t)\|_{L^\infty(\mathbb{T}^2\times\mathbb{R}^+)}\\
&\leq \max \left\{            e^{C\left(1+G_1(t)\right)t}\|I(0)\|_{L^\infty(\mathbb{T}^2\times\mathbb{R}^+)},\mathop{\max}\limits_{\tau \in [0,t]}\left(e^{C\left(1+G_1(t)\right)(t-\tau)}\|I(\tau)\big{|}_{z=0}\|_{L^\infty(\mathbb{T}^2)}\right)
\right\}.
\label{6.3.36}
\end{aligned}
\end{eqnarray}

 To derive a lower bound estimate on $B_{(0,0,0)}$, we have
 \begin{eqnarray}
\begin{aligned}
&\left(\partial_t+u\partial_x+Ku\partial_{y}+(w+\frac{2\sigma}{1+z})\partial_z-\partial^2_z-\epsilon^2\partial^2_x-\epsilon^2\partial^2_y-\partial_{y}Ku\right)B_{(0,0,0)}\\
&=\left(\frac{\sigma}{1+z}w+\frac{\sigma(\sigma-1)}{(1+z)^2}\right)B_{(0,0,0)}.
\label{6.3.37}
\end{aligned}
\end{eqnarray}
Applying the classical maximum principle Lemma \ref{y6.4.6} for (\ref{6.3.37}), we get
\begin{eqnarray}
\begin{aligned}
&\mathop{\min}\limits_{\mathbb{T}^2\times\mathbb{R}^+}(1+z)^{\sigma}\varphi(t)\\
&\geq \max \left\{         1-C\left(1+G_2(t)\right)te^{C\left(1+G_2(t)\right)t}
\right\}\times\min \left\{\mathop{\min}\limits_{\mathbb{T}^2\times\mathbb{R}^+}(1+z)^{\sigma}\varphi_0,        \mathop{\min}\limits_{[0,t]\times\mathbb{T}^2}\varphi\big{|}_{z=0}
\right\}.
\label{6.3.38}
\end{aligned}
\end{eqnarray}

$\textbf{Step 2.}$

Based on inequalities (\ref{6.3.38}) and (\ref{6.3.36}), we have successfully managed to control the underlying  quantities $I$ and $B_{(0,0,0)}$ by considering their initial and boundary values. Nonetheless, the problem does not provide specific boundary values for $I$ and $B_{(0,0,0)}$, thus it becomes necessary to estimate them in this step.

In order to control $I|_{z=0}$, we will make use of the Sobolev embedding argument when $s\geq5$, and the growth rate control argument when $s\geq7$. By combining these boundary estimates on $I$ with (\ref{6.3.38}), we can eventually derive the desired inequalities (\ref{6.3.31}) and (\ref{6.3.32}).

By Lemma \ref{y6.4.3}, we obtain
\begin{eqnarray}
\begin{aligned}
\|I(\tau)|_{z=0}\|_{L^\infty(\mathbb{T})}&\leq C^2\sum \limits_{|\alpha|\leq 2}\big(\|D^{\alpha}\varphi\|_{L^2}^2+\|\partial_{z}D^{\alpha}\varphi\|_{L^2}^2+\|\partial_{xy}^{1}D^{\alpha}\varphi\|_{L^2}^2+\|\partial_{z}\partial_{xy}^{1}D^{\alpha}\varphi\|_{L^2}^2\\
&\quad +\|\partial_{xy}^{2}D^{\alpha}\varphi\|_{L^2}^2+\|\partial_{z}\partial_{xy}^{2}D^{\alpha}\varphi\|_{L^2}^2\big)\\
&\leq 6C^2\|\varphi\|_{H^{s,\gamma}_{g}}^{2},
\end{aligned}
\end{eqnarray}
which,  together with (\ref{6.3.36}),  gives (\ref{6.3.31}).

Now, we derive inequality (\ref{6.3.32}), according to (\ref{6.1.5}) and boundary condition $\partial^{\alpha_1}_x\partial^{\alpha_2}_y u=\partial^{\alpha_1}_x \partial^{\alpha_2}_y w=0$, we have
\begin{eqnarray}
\begin{aligned}
\partial_t D^{\alpha}\varphi\big{|}_{z=0}=D^{\alpha}(\epsilon^2\partial_x^2+\epsilon^2\partial_y^2+\partial_y^2 +\partial_{y}Ku-u\partial_x -Ku\partial_y-w\partial_z )\varphi\big{|}_{z=0}.
\end{aligned}
\end{eqnarray}
When $\alpha_3=0$ and $\alpha_1+\alpha_2 \leq 2$,
\begin{eqnarray}
\begin{aligned}
\partial_t D^{\alpha}\varphi\big{|}_{z=0}=(\epsilon^2\partial^{\alpha_1+2}_{x}\partial^{\alpha_2}_{y}+\epsilon^2\partial^{\alpha_1}_{x}\partial^{\alpha_2+2}_{y}
+\partial^{\alpha_1}_{x}\partial^{\alpha_2}_{y}\partial_z^2 )\varphi\big{|}_{z=0}.
\end{aligned}
\end{eqnarray}
When $1 \leq \alpha_3 \leq 2$ and $\alpha_1=\alpha_2=0$,
\begin{eqnarray}
\begin{aligned}
\partial_t D^{\alpha}\varphi\big{|}_{z=0}
&=\big(\epsilon^2\partial^{\alpha_3}_{z}\partial^{2}_{x}+\epsilon^2\partial^{\alpha_3}_{z}\partial^{2}_{y}+\partial^{\alpha_3+2}_{z}+(3\alpha_3-2)\partial_{y}K\varphi\partial_{z}^{\alpha-1}\\
&\quad-\alpha_3\varphi\partial_{z}^{\alpha_3-1}\partial_{x}-K\alpha_3\varphi\partial_{z}^{\alpha_3-1}\partial_{y}\big)\varphi\big{|}_{z=0}.
\end{aligned}
\end{eqnarray}
When $\alpha_3=1$ and $\alpha_1+\alpha_2 =1$,
\begin{eqnarray}
\begin{aligned}
\partial_t D^{\alpha}\varphi|_{z=0}&=(\epsilon^2\partial_{z}\partial_{xy}^{1}\partial^{2}_{x} +\epsilon^2\partial_{z}\partial_{xy}^{1}\partial^{2}_{y} +
\partial _{xy}^{1}\partial _{z}^{3}+\partial_{xy}^{1}\partial_{y}K\varphi+2\varphi\partial_{xy}+\varphi\partial_{xy}^{1}\partial_x\\
&\quad+\partial_{xy}^{1}\varphi\partial_{x}+K\varphi\partial_{xy}^{1}\partial_y+\partial_{xy}^{1}(K\varphi)\partial_{y}
)\varphi\big{|}_{z=0}.
\end{aligned}
\end{eqnarray}
For $s\geq 7$, using Lemma \ref{y6.4.4}, we have
 \begin{eqnarray}
\begin{aligned}
\|\partial_t I\big{|}_{z=0}\|_{{L^\infty}(\mathbb{T}^2)}&\leq C_{s, \gamma  } \|D^{\alpha}\varphi\partial_t D^{\alpha}\varphi \big{|}_{z=0}\|_{{L^\infty}(\mathbb{T}^2)}\\
&\leq  C_{s, \gamma  }\|\varphi(\tau)\|_{H^{s,\gamma}_{g}}(\|\varphi(\tau)\|_{H^{s,\gamma}_{g}}+\|\partial_{xy}K(\tau)\|_{{L^{\infty}}(\mathbb{T}^2)}\|\varphi(\tau)\|_{H^{s,\gamma}_{g}}^2).
\end{aligned}
\end{eqnarray}
Hence, a direct integration yields
 \begin{eqnarray}
\begin{aligned}
\|I(t)\big{|}_{z=0}\|_{{L^\infty}(\mathbb{T}^2)}\leq  \|I(0)\big{|}_{z=0}\|_{{L^\infty}(\mathbb{T}^2)}+C_{s, \gamma  } [1+Y(t)W(t)]W(t)^2t,
\end{aligned}
\end{eqnarray}
which, together with (\ref{6.3.36}), gives (\ref{6.3.32}). For $s\geq 5$,
 \begin{eqnarray}
\begin{aligned}
\|\partial_t \varphi|_{z=0}\|_{{L^\infty}(\mathbb{T}^2)}\leq  CW(t).
\end{aligned}
\end{eqnarray}
Thus, a direct integration leads to
 \begin{eqnarray}
\begin{aligned}
\mathop{\min}\limits_{\mathbb{T}^2}{\varphi(t)}|_{z=0}\geq \mathop{\min}\limits_{\mathbb{T}^2}{\varphi_0}\big{|}_{z=0}-CW(t)t,
\end{aligned}
\end{eqnarray}
which, together with (\ref{6.3.38}), gives (\ref{6.3.33}). The proof of Lemma \ref{y6.2.3} is thus complete.
 \hfill $\Box$

\section{Local-in-time existence and uniqueness}
\subsection{Local-in-time existence }
In this subsection, we go back to use the symbol $(u^\epsilon, K^\epsilon u^{\epsilon}, w^\epsilon, \varphi^{\epsilon})$ instead of $(u, Ku, w, \varphi)$ from Subsections 2.1-2.4 to denote the solution to the regularized system (\ref{6.1.5}). To obtain the local-in-time solution of the initial-boundary value problem (\ref{6.1.1})-(\ref{6.1.2}) or (\ref{6.1.4})-(\ref{6.1.44}), we will construct the solution to the Prandtl equations (\ref{6.1.4}) by passing to the limit $\epsilon \rightarrow 0^+$ in the regularized Prandtl equations (\ref{6.1.5}). We only sketch the proof into five steps and more details can be found in \cite{mw}.\\

\begin{Lemma}
Let $s\ge 5$ be an integer, $\gamma \ge 1, \sigma > \gamma + \frac{1}{2},
\delta \in (0, \frac{1}{2})$ and $\epsilon \in (0,1]$. If
$\varphi_0 \in H^{s,\gamma}_{\sigma, 2\delta}, U$ and $P^\epsilon$ are given and satisfy the
regularized Bernoulli's law (\ref{6.1.44}) and the regularity assumption (\ref{6.1.8}), when $s=5$, we further assume that $\delta > 0$ is chosen small
enough such that the initial hypothesis (\ref{6.1.888}) holds. Then there exists a uniform
life-span $T := T (s,\gamma, \sigma, \delta, \|\varphi_0\|_{H^{s,\gamma}}, U) > 0$, which is independent of $\epsilon$, such that the
regularized vorticity system (\ref{6.1.5}) has a solution
$\varphi^{\epsilon} \in C([0,T]; H^{s,\gamma}_{\sigma, \delta}) \cap C^1 ([0,T]; H^{s-2,\gamma})$ with the following uniform
(in $\epsilon$) estimates:

(i) (Uniform  Weighted $H^s$ Estimate) For any $\epsilon \in [0,1]$ and any $t \in [0,T]$,

\begin{eqnarray}
    \|\varphi^\epsilon\|_{H^{s,\gamma}_{g}} \leq 4 \|\varphi_0\|_{H^{s,\gamma}_{g}}.
\end{eqnarray}

(ii) (Uniform Weighted $L^\infty$ Bound) For any $\epsilon \in [0,1]$ and $t \in [0,T]$,
\begin{eqnarray}
    \left\|\sum \limits_{|\alpha|\leq 2}  |(1+z)^{\sigma+\alpha_{3}}D^{\alpha}\varphi^\epsilon (t) | ^{2}\right\|_{L^{\infty}(\mathbb{T}^2\times \mathbb{R}^+)}\leq\frac{1}{\delta^{2}}.
\end{eqnarray}
(iii) (Uniform Weighted $L^\infty$ Lower Bound) For any $\epsilon \in [0,1]$ and $t \in
[0,T]$,
\begin{eqnarray}
\min \limits_{\mathbb{T}^2\times \mathbb{R}^+}(1+z)^{\sigma} \varphi^\epsilon(t) \geq \delta.
\end{eqnarray}
\end{Lemma}
\textbf{Proof}.
$\textbf{Step 1.}$ According to the definition of $F$, assumption (\ref{6.1.8}) and (\ref{6.1.88}), and the regularized Bernoulli's law,
\begin{eqnarray}
\|F\|_{L^\infty} \leq C_{s, \gamma,\sigma,\delta  }\left(M_{U}+M_{K}\right) < +\infty,
\end{eqnarray}
which, together with (\ref{6.3.2222}), implies
\begin{eqnarray}
\begin{aligned}
\|\varphi^\epsilon\|_{H^{s,\gamma}_{g}}^2 &\leq \left\{\|\varphi_{0}\|_{H^{s,\gamma}_{g}}^2+C_{s, \gamma,\sigma,\delta  }\left(M_{U}+M_{K}\right)t_{1}\right\}\\
&\quad \times \left\{1-C_{s, \gamma,\sigma,\delta  }\left(\|\varphi_{0}\|_{H^{s,\gamma}_{g}}^2+C_{s, \gamma,\sigma,\delta  }\left(M_{U}+M_{K}\right)t_{1}\right)^{s-1}t_{2}\right\}^{-\frac{1}{s-1}},
\label{6.2.01555}
\end{aligned}
\end{eqnarray}
for $t_{1}$ and $t_{2}$ to be chosen later. If choosing
\begin{eqnarray}
t_{1} \leq \frac{3\|\varphi_0\|_{H^{s,\gamma}_{g}}^2}{C_{s, \gamma,\sigma,\delta  }(M_{U}+M_{K})}
\end{eqnarray}
 in (\ref{6.2.01555}), then we can  get
\begin{eqnarray}
\begin{aligned}
\|\varphi^\epsilon\|_{H^{s,\gamma}_{g}}^2 &\leq 4\|\varphi_{0}\|_{H^{s,\gamma}_{g}}^{2}\times \left\{1-C_{s, \gamma,\sigma,\delta  }\left(4\|\varphi_{0}\|_{H^{s,\gamma}_{g}}^{2}\right)^{s-1}t_{2}\right\}^{-\frac{1}{s-1}}.
\label{6.2.015555}
\end{aligned}
\end{eqnarray}
Further, choosing
\begin{eqnarray}
t_{2} \leq \frac{1-2^{2-2s}}{2^{2s-2}C_{s, \gamma,\sigma,\delta  }\|\varphi_0\|_{H^{s,\gamma}_{g}}^{2s-2}}
\end{eqnarray}
 in (\ref{6.2.015555}) gives
\begin{eqnarray}
\begin{aligned}
\|\varphi^\epsilon\|_{H^{s,\gamma}_{g}}^2 &\leq 4\|\varphi_{0}\|_{H^{s,\gamma}_{g}}^{2}\times \left(2^{2-2s}\right)^{-\frac{1}{s-1}} \leq 16\|\varphi_{0}\|_{H^{s,\gamma}_{g}}^{2}.
\end{aligned}
\end{eqnarray}
Note that the functions in both brackets in  (\ref{6.2.01555}) are monotonically increasing with respect to time $t_{1}$ and $t_{2}$, thus we derive the uniform estimate for any $\epsilon \in [0,1]$ and any $t \in [0,T_1]$,
\begin{eqnarray}
\|\varphi^\epsilon\|_{H^{s,\gamma}_{g}} \leq 4 \|\varphi_0\|_{H^{s,\gamma}_{g}},
\label{6.2.0155555}
\end{eqnarray}
provided that $T_1$ is chosen  such that
\begin{eqnarray*}
T_1 :=\min \left\{\frac{3\|\varphi_0\|_{H^{s,\gamma}_{g}}^2}{C_{s, \gamma,\sigma,\delta  }(M_{U}+M_{K})},\frac{1-2^{2-2s}}{2^{2s-2}C_{s, \gamma,\sigma,\delta  }\|\varphi_0\|_{H^{s,\gamma}_{g}}^{2s-2}}\right\}.
\end{eqnarray*}
$\textbf{Step 2.}$ When $s \geq 7$, we know from definition (\ref{MK}) of $F$, (\ref{6.2.0122}) of $W$ and (\ref{6.2.012}) of $G_1$ that for any $t \in [0,T_1]$,
\begin{eqnarray}
F(t) \leq M_{K}, \quad W(t) \leq 4 \|\varphi_0\|_{H^{s,\gamma}_{g}} \quad and \quad G_1(t) \leq M_{K}+(4 \|\varphi_0\|_{H^{s,\gamma}_{g}}+M_{U})(1+M_{K}),
\label{6.2.015}
\end{eqnarray}
which, together with (\ref{6.3.32}), gives
\begin{eqnarray}
\begin{aligned}
&\|I(t)\|_{L^\infty(\mathbb{T}^2\times\mathbb{R}^+)}\\
&\leq  \left\{           \|I(0)\|_{L^\infty(\mathbb{T}^2\times\mathbb{R}^+)}+C[1+ 4 M_{K}\|\varphi_0\|_{H^{s,\gamma}_{g}}]16\|\varphi_0\|_{H^{s,\gamma}_{g}}^2t_3
\right\} e^{C\left(1+M_{K}+(4 \|\varphi_0\|_{H^{s,\gamma}_{g}}+M_{U})(1+M_{K})\right)t_4}
\label{step3.32}
\end{aligned}
\end{eqnarray}
for $t_{3}$ and $t_{4}$ to be chosen later. If choosing
\begin{eqnarray}
t_{3} \leq \frac{1}{64\delta^2 C(1+4 M_{K}\|\varphi_0\|_{H^{s,\gamma}_{g}}) \|\varphi_0\|_{H^{s,\gamma}_{g}}^2},
\end{eqnarray}
 in (\ref{step3.32}), then we can  get
\begin{eqnarray}
\begin{aligned}
\|I(t)\|_{L^\infty(\mathbb{T}^2\times\mathbb{R}^+)}\leq  \left\{           \|I(0)\|_{L^\infty(\mathbb{T}^2\times\mathbb{R}^+)}+\frac{1}{4\delta^{2}}
\right\} e^{C\left(1+M_{K}+(4 \|\varphi_0\|_{H^{s,\gamma}_{g}}+M_{U})(1+M_{K})\right)t_4}.
\end{aligned}
\end{eqnarray}
Further, choosing
\begin{eqnarray}
t_{4} \leq \frac{\ln 2}{C[1+M_{K}+(4 \|\varphi_0\|_{H^{s,\gamma}_{g}}+M_{U})(1+M_{K})]}
\end{eqnarray}
 in  the above inequality, gives
\begin{eqnarray}
\begin{aligned}
\|I(t)\|_{L^\infty(\mathbb{T}^2\times\mathbb{R}^+)}\leq 2 \left\{           \|I(0)\|_{L^\infty(\mathbb{T}^2\times\mathbb{R}^+)}+\frac{1}{4\delta^{2}}
\right\}.
\end{aligned}
\end{eqnarray}
Then using  initial assumption
\begin{eqnarray}
\sum \limits_{|\alpha|\leq 2}  |(1+z)^{\sigma+\alpha_{3}}D^{\alpha}\varphi_0 | ^{2}\leq\frac{1}{4\delta^{2}},
\end{eqnarray}
we have the upper bound
\begin{eqnarray}
\Big\|\sum \limits_{|\alpha|\leq 2}  |(1+z)^{\sigma+\alpha_{3}}D^{\alpha}\varphi^\epsilon (t) | ^{2}\Big\|_{L^{\infty}(\mathbb{T}^2\times \mathbb{R}^+)}\leq\frac{1}{\delta^{2}}
\label{6.2.019}
\end{eqnarray}
provided that $T_2$ is chosen  such that
\begin{eqnarray*}
T_2 :=\min \left\{T_1,\frac{1}{64\delta^2 C(1+4 M_{K}\|\varphi_0\|_{H^{s,\gamma}_{g}) \|\varphi_0\|_{H^{s,\gamma}_{g}}^2}},\frac{\ln 2}{C[1+M_{K}+(4 \|\varphi_0\|_{H^{s,\gamma}_{g}}+M_{U})(1+M_{K})]}\right\},
\end{eqnarray*}
for all $t \in [0,T_2]$. When $s \geq 5$, from the initial hypothesis $\|\varphi_0\| \leq C\delta^{-1}$, we have the same estimate (\ref{6.2.019}) for all $t\in [0,T_2]$.\\
$\textbf{Step 3.}$
From (\ref{6.3.33}), we have
\begin{eqnarray}
\begin{aligned}
\mathop{\min}\limits_{\mathbb{T}^2\times\mathbb{R}^+}(1+z)^{\sigma}\varphi(t)\geq \max \left\{         1-CNte^{CNt}
\right\}\cdot\min \left\{\mathop{\min}\limits_{\mathbb{T}^2\times\mathbb{R}^+}(1+z)^{\sigma}\varphi_0-4C|\varphi_0\|_{H^{s,\gamma}_{g}}t
\right\},
\end{aligned}
\end{eqnarray}
where $N :=1+G_2(t)\leq1+M_{U}+4 \|\varphi_0\|_{H^{s,\gamma}_{g}}$.

Similarly, let us choose
\begin{eqnarray*}
T_3 :=\min \left\{T_1,\frac{\delta}{8 C\|\varphi_0\|_{H^{s,\gamma}_{g}}},\frac{1}{6CN},\frac{\ln 2}{CN}\right\}.
\end{eqnarray*}
Hence, we derive the uniform estimate for any $\epsilon \in [0,1]$ and any $t \in [0,T_3]$,
\begin{eqnarray}
\min \limits_{\mathbb{T}^2\times \mathbb{R}^+}(1+z)^{\sigma} \varphi^\epsilon(t) \geq \delta.
\end{eqnarray}
$\textbf{Step 4.}$
In summary, the above uniform estimates  hold for any $t \in [0,T]$, if $T$ is chosen to satisfy $T=\min\{T_1,T_2,T_3\}$. Further using almost equivalence relation  (\ref{6.4.1}), Lemma \ref{y6.4.4} and (\ref{6.2.0155555}), we have
\begin{align}
&\sup  \limits_{0 \leq t \leq T}( \|\varphi^\epsilon\|_{H^{s,\gamma}}+\|u^\epsilon-U\|_{H^{s,\gamma-1}}) \nonumber\\
&\leq 4C\|\varphi_0\|_{H^{s,\gamma}_{g}}+ C\sup \limits_{0 \leq t \leq T}\left(\sum_{k=0}^{s}\|(1+z)^{\gamma-1}\partial_{xy}^{k}(u^\epsilon-U)\|_{L^{2}}\right)\nonumber\\
&\leq C\|\varphi_0\|_{H^{s,\gamma}_{g}}+ C_{s, \gamma,\sigma,\delta  }\sup \limits_{0 \leq t \leq T}\left(\|\varphi\|_{H^{s,\gamma}_{g}}+\|\partial_{xy}^s U\|_{L^2(\mathbb{T}^2)}\right)\nonumber\\
& \leq C_{s, \gamma,\sigma,\delta  }\left(\|\varphi_0\|_{H^{s,\gamma}_{g}}+\sup \limits_{0 \leq t \leq T}\|\partial_{xy}^s U\|_{L^2(\mathbb{T}^2)}\right)\nonumber\\
&< + \infty.
\label{6.2.016}
\end{align}
From the equation (\ref{6.1.4}),  (\ref{6.1.5}), (\ref{6.2.016}) and Lemma \ref{y6.4.4}, we also have $\partial_t \varphi^{\epsilon}$	 and $\partial_t (u^{\epsilon}-U)$ are uniformly
bounded in $L^{\infty}([0,T];H^{s-2,\gamma})$ and $L^{\infty}([0,T];H^{s-2,\gamma-1})$ respectively. By Lions-Aubin lemma and the compact embedding of $H^{s,\gamma}$ in $H^{s^{\prime}}_{loc}$
, we have after taking a subsequence, as $\epsilon_k \rightarrow 0^+$,
\begin{equation}\left\{
\begin{array}{ll}
\varphi^{\epsilon_k} \stackrel{*}{\rightharpoonup} \varphi, &{\rm in} \quad L^{\infty}([0,T];H^{s,\gamma}),\\
\varphi^{\epsilon_k} \rightarrow \varphi, &{\rm in} \quad C([0,T];H^{s^{\prime}}_{loc}),\\
u^{\epsilon_k}-U \stackrel{*}{\rightharpoonup} u-U,  &{\rm in} \quad L^{\infty}([0,T];H^{s,\gamma-1}),\\
u^{\epsilon_k} \rightarrow u, &{\rm in} \quad C([0,T];H^{s^{\prime}}_{loc}),
\end{array}\label{step4.1}
         \right.\end{equation}
for all $s^{\prime} < s$, where
\begin{equation}\left\{
\begin{array}{ll}
\varphi=\partial_z u \in L^{\infty} ([0,T];H^{s,\gamma})\cap \bigcap_{s^\prime <s}C([0,T];H^{s^{\prime}}_{loc}),\\
 u-U \in L^{\infty} ([0,T];H^{s,\gamma-1})\cap \bigcap_{s^\prime <s}C([0,T];H^{s^{\prime}}_{loc}).
\end{array}\label{step4.2}
         \right.\end{equation}
$\textbf{Step 5.}$
Using the local uniform
convergence of $\partial_x u^{\epsilon_k}$ and $\partial_y (K^{\epsilon_k}u ^{\epsilon_k})$, we also have the pointwise convergence of $w^{\epsilon_k}:$ as $\epsilon \rightarrow 0^+$,
\begin{eqnarray}
w^{\epsilon_k}=-\int^{z}_{0} \partial_x u ^{\epsilon_k}~ dz-\int^{z}_{0} \partial_y (K^{\epsilon_k}u ^{\epsilon_k})~ dz \rightarrow -\int^{z}_{0} \partial_x u ~ dz-\int^{z}_{0} \partial_y (Ku )~ dz=: w.
\end{eqnarray}
Combining (\ref{step4.1}) and (\ref{step4.2}), one may justify the pointwise convergences of all terms in the
regularized Prandtl equations (\ref{6.1.4}). Thus, passing the limit $\epsilon_k \rightarrow 0^+$ in the initial-boundary value problem (\ref{6.1.4}) and the regularized Bernoulli's law
(\ref{6.1.44}), we get that the limit $(u, Ku, w)$ solves the initial-boundary value problem (\ref{6.1.222}) in the classical sense. By the equivalent system, we get that the limit $(u, Ku, w)$ solves the initial-boundary value problem (\ref{6.1.1}) in the classical sense.

Lastly, in order to complete the proof, it remains to justify that $\varphi \in L^{\infty}([0,T];H^{s,\gamma}_{\sigma, \delta})$ and matching condition $(\ref{6.1.1})_6$. Since $D^\alpha   \varphi^{\epsilon_k}$ converges to $D^\alpha  \varphi$ pointwisely for all $|\alpha| \leq 2$ and that  $\varphi^{\epsilon_k} $ satisfies

\begin{eqnarray}
    \sum_{|\alpha| \leq 2} |(1+z)^{\sigma + \alpha_3} D^\alpha  \varphi^{\epsilon_k} |^2  \le \frac{1}{\delta^2}
    \label{step5.1}
\end{eqnarray}
and
\begin{eqnarray}
    \min \limits_{\mathbb{T}^2\times \mathbb{R}^+}(1+z)^\sigma \varphi^{\epsilon_k} \geq \delta,
    \label{step5.2}
\end{eqnarray}
we deduce that (\ref{step5.1}) and  (\ref{step5.2}) still holds for $\varphi$, and hence, $\varphi\in  L^\infty(0,T;H^s) $.

Also, by the Lebesgue's dominated convergence theorem,
\begin{eqnarray}
    \int^{+\infty}_0 \varphi  dz = \lim_{\epsilon_k \rightarrow 0^+} \int^{+\infty}_0 \varphi^{\epsilon_k}  dz = U
\end{eqnarray}
and
\begin{eqnarray}
    \int^{+\infty}_0 K\varphi  dz = \lim_{\epsilon_k \rightarrow 0^+} \int^{+\infty}_0 K\varphi^{\epsilon_k}  dz = KU=V
\end{eqnarray}
which are equivalent to  condition $(\ref{6.1.1})_6$ because $\varphi = \partial_z  u
> 0$. This completes the proof
of the existence.
 \hfill $\Box$

\subsection{Uniqueness of solutions}
The aim of this section is to prove the uniqueness of $H^{s,\gamma}_{\sigma, \delta}$ solutions to the 3D  Prandtl model. To show the uniqueness, we will generalize the nonlinear
cancelation  to the $L^2$ comparison. This motivates us to consider the quantity $\bar{g}$ below.

Specifically, the uniqueness of $H^{s,\gamma}_{\sigma, \delta}$ solutions to the Prandtl equations (\ref{6.1.1}) is a direct consequence of the following $L^2$ comparison principle.

\begin{Lemma}($L^2$ Comparison Principle) \label{y.Uniqueness}
For any $s \ge 5, \gamma \geq 1, \sigma > \gamma + \frac{1}{2}$ and $\delta \in (0,1)$,
let $(u_i, v_i, w_{i})$ solve the Prandtl equations (\ref{6.1.222}) with  the vorticity $\varphi_i
:= \partial_z u_i \in C([0,T]; H^{s,\gamma}_{\sigma, \delta}) \cap C^1 ([0,T]; H^{s-2,\gamma})$ for $i = 1, 2$. Define
$\bar{g} := \varphi_1 - \varphi_2 + \frac{\partial_z\varphi_2}{\varphi_2} (u_1 - u_2)$. Then
we have
\begin{eqnarray}
\begin{aligned}
\|\bar g\|_{L^2}^2 +\int^t_0 \|\partial_z \bar{g} \|^2_{L^2}\leq \|\bar g_{0}\|_{L^2}^2 +C_{\gamma, \sigma, \delta, \varphi, K, U} \int^t_0 \|\bar{g} \|^2_{L^2},
\end{aligned}
\end{eqnarray}
where the constant $ C_{\gamma, \sigma, \delta, \varphi, K, U}$ depends on $\gamma, \sigma,
\delta, \|w_i\|_{H^{5,\gamma}_{g}}, \|\partial_{y}K\|_{{L^{\infty}}(\mathbb{T}^2)}$ and $\|\partial_{xy}^{5}U\|_{{L^2}(\mathbb{T}^2)}$
 only.

\end{Lemma}
\textbf{Proof}.
Let us denote $(\bar u, K\bar{u}, \bar{w})=(u_1, Ku_1, w_1)-(u_2, Ku_2, w_2)$, $\bar {\varphi} ={\varphi}_{1}-{\varphi}_{2}$, $a_2=\frac{\partial_z \varphi_2}{\varphi_{2}}$ and $\bar g =\bar {\varphi} -a_2 \bar u$. It is easy to check that $\bar g =\bar {\varphi} -a_2 \bar u=\varphi_2 \partial_z(\frac{\bar u}{\varphi_2})$ and the evolution equation on $\bar g$ is as follows
\begin{eqnarray}
\begin{aligned}
&(\partial_t+u_1\partial_x+Ku_1\partial_y+w_1\partial_z-\partial^{2}_{z})\bar g\\
&=(\partial_t+u_1\partial_x+Ku_1\partial_y+w_1\partial_z-\partial^2_z)\bar \varphi-a_2(\partial_t+u_1\partial_x+Ku_1\partial_y+w_1\partial_z-\partial^2_z)\bar u\\
&\quad -\bar u(\partial_t+u_1\partial_x+Ku_1\partial_y+w_1\partial_z-\partial^2_z)a_2-2\partial_z a_2 \bar {\varphi}.
 \label{6.33.40}
\end{aligned}
\end{eqnarray}
Next, we calculate the values of the first three terms on the right-hand  side of the equality (\ref{6.33.40}) respectively. Recalling the vorticity system (\ref{6.1.5}), we have
\begin{eqnarray}
\begin{aligned}
(\partial_t+u_i\partial_x+Ku_i\partial_y+w_i\partial_z)\partial_z \varphi_i&=\partial_z^3 \varphi_i+[\partial_x u_i+\partial_y (Ku_i)] \partial_{z} \varphi_{i} -\varphi_i \partial_x \varphi_i-K\varphi_i \partial_y \varphi_i\\
&\quad+\partial_{y}K(\partial_{z} \varphi u+\varphi^{2}),
\end{aligned}
\end{eqnarray}
then, according to the definition of $a_i$, we get
\begin{align}
(\partial_t+u_i\partial_x+Ku_i\partial_y+w_i\partial_z)a_i&=\frac{(\partial_t+u_i\partial_x+Ku_i\partial_y+w_i\partial_z)\partial_z \varphi_i}{\varphi_i}-\frac{\partial_z \varphi_i(\partial_t+u_i\partial_x+Ku_i\partial_y+w_i\partial_z) \varphi_i}{\varphi_i^2}\nonumber\\
&=\frac{\partial^3_z \varphi_i}{\varphi_i}+a_i\partial_x u_i+a_i\partial_y (Ku_i)-\partial_x \varphi_i-K\partial_y \varphi_i -a_i \frac{\partial_z^2 \varphi_i}{\varphi_i}+\partial_{y}K\varphi_{i},
\end{align}
which, combined with the fact
\begin{eqnarray*}
\frac{\partial^3_z \varphi_i}{\varphi_i}-a_i \frac{\partial_y^2 a_i}{\varphi_i}=\partial_z^2 a_i+2a_i \partial_z a_i,
\end{eqnarray*}
implies
\begin{eqnarray}
\begin{aligned}
(\partial_t+u_i\partial_x+Ku_i\partial_y+w_i\partial_z-\partial_{z}^{2})a_i=a_i\partial_x u_i+a_i\partial_y (Ku_i)-\partial_x \varphi_i-K\partial_y \varphi_i +\partial_{y}K\varphi +2a_i \partial_z a_i.
\end{aligned}
\end{eqnarray}
Furthermore, we conclude that
\begin{eqnarray}
\begin{aligned}
&(\partial_t+u_1\partial_x+Ku_1\partial_y+w_1\partial_z-\partial_{z}^{2})a_2\\
&=a_2\partial_x u_2+a_2\partial_y (Ku_2)-\partial_x \varphi_2-K\partial_y \varphi_2 +\partial_{y}K\varphi_{2} +2a_2 \partial_z a_2+(\bar u \partial_x+K\bar {u} \partial_y+\bar{w}\partial_{z})a_2.
 \label{6.33.41}
\end{aligned}
\end{eqnarray}
For the estimate of $(\partial_t+u_1\partial_x+Ku_1\partial_y+w_1\partial_z-\partial_{z}^{2})\bar u$ and $(\partial_t+u_1\partial_x+Ku_1\partial_y+w_1\partial_z-\partial_{z}^{2})\bar \varphi$, from the evolution equation on $u$ and on $\varphi$, we can derive that
\begin{eqnarray}
\begin{aligned}
(\partial_t+u_1\partial_x+Ku_1\partial_y+w_1\partial_z-\partial_{z}^{2})\bar u&=-\bar u \partial_x u_2-K\bar {u} \partial_y  u_2-\bar {w} \partial_z  u_2,
\end{aligned}
\end{eqnarray}
and
\begin{eqnarray}
\begin{aligned}
(\partial_t+u_1\partial_x+Ku_1\partial_y+w_1\partial_z-\partial_{z}^{2})\bar {\varphi}&=-\bar u \partial_x \varphi_2-K\bar {u} \partial_y  \varphi_2-\bar {w} \partial_z  \varphi_2+\partial_{y}K(\bar{u}\varphi_{2}+u_1\bar{\varphi}).
 \label{6.33.42}
\end{aligned}
\end{eqnarray}
Using (\ref{6.33.40}) and (\ref{6.33.41})-(\ref{6.33.42}), we have
\begin{eqnarray}
\begin{aligned}
(\partial_t+u_1\partial_x+Ku_1\partial_y+w_1\partial_z-\partial_{z}^{2})\bar g&= -\bar u(\bar u \partial_xa_{2}+K\bar {u} \partial_ya_{2}+\bar{w}\partial_{z}a_{2}+2a_{2}\partial_{z}a_{2})-2\bar {\varphi} \partial_z a_2\\
&\quad+\partial_y{K}(u_1\bar{\varphi}-\bar{u}a_{2}u_{2}).
\end{aligned}
\end{eqnarray}
Now we derive $L^2$ estimate on $\bar g$. For any $t \in (0,T]$, multiplying by $2 \bar g$ and then integrating by parts over $\mathbb{T}^2\times\mathbb{R}_{+}$, we obtain
\begin{eqnarray}
\begin{aligned}
&\frac{d}{dt}\|\bar g\|_{L^2}^2+2\|\partial_z \bar g\|_{L^2}^2\\
&=\iint_{\mathbb{T}^{2}}\bar g \partial_z \bar g \big{|}_{z=0}dxdy-2\iiint \bar g \bar u(\bar u \partial_xa_{2}+K\bar {u} \partial_ya_{2}+\bar{w}\partial_{z}a_{2}+2a_{2}\partial_{z}a_{2})-4\iiint \bar g \bar {\varphi} \partial_z a_2  \\
&\quad+2\iiint \bar{g}\partial_y{K}(u_1\bar{\varphi}-\bar{u}a_{2}u_{2}) -2\iiint \bar g(u_1\partial_x\bar{g}+Ku_1\partial_y\bar{g}+w_1\partial_z\bar{g}).
 \label{6.33.43}
\end{aligned}
\end{eqnarray}
We need to estimate all terms on the right-hand side of equation (\ref{6.33.43}). Applying the simple trace theorem and Young's inequality,
\begin{align}
\left|\int_{\mathbb{T}^2}\bar {g} \partial_z \bar {g} |_{z=0}dxdy\right|
&\leq|\iint a_2 |\bar g|^2  dxdydz|+|\iint \partial_z a_2 |\bar g|^2  dxdydz|+2|\iint  a_2 \bar g \partial_z \bar g  dxdydz|\nonumber\\
&\leq \frac{1}{2}\|\partial_z \bar g\|_{L^2}^2+C_{\sigma,\delta  }\|\bar g\|_{L^2}^2,
\label{w6.3.20}
\end{align}
where we have used the fact $ \partial_z \bar g \big{|}_{z=0}=-a_2 \bar g  |_{z=0}$ . We claim $\|\frac{\bar u}{1+z}\|_{L^2} \leq C_{\sigma,\delta} \|\bar g\|_{L^2}$, so by Lemma \ref{y6.4.4},
\begin{align}
&-2\iint \bar g \bar u(\bar u \partial_x a_2+\bar v \partial_y a_2+2a \partial_y a_2) \nonumber\\
&\leq  C\|(1+z)(\bar u \partial_xa_{2}+K\bar {u} \partial_ya_{2}+\bar{w}\partial_{z}a_{2}+2a_{2}\partial_{z}a_{2})\|_{L^\infty}\left\|\frac{\bar u}{1+z}\right\|_{L^2}\|\bar g \|_{L^2}\nonumber\\
&\leq C_{ \gamma,\sigma,\delta  }\left(1+\|w_i\|_{H^{5,\gamma}_{g}}+\|\partial_{xy}^{5}U\|_{{L^2}(\mathbb{T}^2)}+\|K\|_{{L^{\infty}}(\mathbb{T}^2)}\|w_i\|_{H^{5,\gamma}_{g}}
+\|K\|_{{L^{\infty}}(\mathbb{T}^2)}\|\partial_{xy}^{5}U\|_{{L^2}(\mathbb{T}^2)}\right)\|\bar g\|_{L^2}^2.
\end{align}
Now we shall prove  that $\|\frac{\bar u}{1+z}\|_{L^2}$ can be controlled by $\|\bar g\|_{L^2}$. Indeed, since $\delta \leq (1+z)^{\delta} \varphi_2\leq \delta^{-1}$ and Lemma \ref{y6.4.1},
\begin{eqnarray}
\begin{aligned}
\left\|\frac{\bar u}{1+z}\right\|_{L^2} \leq \delta^{-1}\left \|(1+z)^{-\sigma-1}\frac{\bar u}{\varphi_2}\right\|_{L^2} \leq C_{\sigma,\delta}\left\|(1+z)^{-\sigma} \partial_{z}\left(\frac{\bar u}{\varphi_2}\right)\right\|_{L_2} \leq C_{\sigma,\delta} \|\bar g\|_{L^2}.
\end{aligned}
\end{eqnarray}
In addition, we also have
\begin{eqnarray}
\begin{aligned}
\| \bar \varphi\|_{L^2} \leq \|\bar g \|_{L^2}+\delta^{-2}\left\|\frac{\bar u}{1+z}\right\|_{L^2}\leq C_{\sigma,\delta} \|\bar g\|_{L^2},
\end{aligned}
\end{eqnarray}
then we can obtain that
\begin{eqnarray}
\begin{aligned}
-4\iint  \bar g \bar \varphi \partial_z a_2 \leq C_{\sigma,\delta} \|\bar g\|_{L^2}^2
\end{aligned}
\end{eqnarray}
and
\begin{eqnarray}
\begin{aligned}
2\iiint \bar{g}\partial_y{K}(u_1\bar{\varphi}-\bar{u}a_{2}u_{2})\leq C_{ \gamma,\sigma,\delta  }\|\partial_{y}K\|_{{L^{\infty}}(\mathbb{T}^2)}\left(\|w_i\|_{H^{5,\gamma}_{g}}
+\|\partial_{xy}^{5}U\|_{{L^2}(\mathbb{T}^2)}\right)\|\bar g\|_{L^2}^2.
\end{aligned}
\end{eqnarray}
For the last term in (\ref{6.33.43}), using the integration by parts, boundary condition $(u_1, w_1)|_{z=0}=0$ and $\partial_x u_1+\partial_y (Ku_1)+\partial_z w_1=0$, we readily  show
\begin{eqnarray}
\begin{aligned}
-2\iiint \bar g(u_1\partial_x\bar{g}+Ku_1\partial_y\bar{g}+w_1\partial_z\bar{g})=0.
\label{w6.3.26}
\end{aligned}
\end{eqnarray}
Combining all of the above estimates (\ref{w6.3.20})-(\ref{w6.3.26}), we can derive from (\ref{6.33.43})
\begin{eqnarray}
\begin{aligned}
\frac{d}{dt}\|\bar g\|_{L^2}^2 \leq C_{ \gamma,\sigma,\delta  }\left[1+\left(1+\|\partial_{y}K\|_{{L^{\infty}}(\mathbb{T}^2)}\right)\left(\|w_i\|_{H^{5,\gamma}_{g}}+\|\partial_{xy}^{5}U\|_{{L^2}(\mathbb{T}^2)}\right)\right]\|\bar g\|_{L^2}^2.
\end{aligned}
\end{eqnarray}
Integrating the above inequality over  (0, t), then the proof is complete.
\hfill $\Box$

Now applying  Gronwall's inequality for Lemma \ref{y.Uniqueness}, yields
\begin{eqnarray}
\begin{aligned}
\|\bar g(t)\|_{L^2}^2 \leq \|\bar g(0)\|_{L^2}^2e^{Ct},
\end{aligned}
\end{eqnarray}
here  $C=C_{ \gamma,\sigma,\delta  }\left[1+\left(1+\|\partial_{y}K\|_{{L^{\infty}}(\mathbb{T}^2)}\right)\left(\|w_i\|_{H^{5,\gamma}_{g}}+\|\partial_{xy}^{5}U\|_{{L^2}(\mathbb{T}^2)}\right)\right]$, and this implies $\bar g=0$ due to $u_1|_{t=0}=u_2|_{t=0}$. Since $\varphi_{2}\partial_z (\frac{u_1-u_2}{\varphi_{2}})=\bar g =0$, we have
\begin{eqnarray}
\begin{aligned}
u_1-u_2=q \varphi_{2}
\end{aligned}
\end{eqnarray}
for some function $q=q(t,x)$. By using the  Oleinik's monotonicity assumption $\varphi_{2} >0 $ and boundary condition $u_1|_{z=0}=u_2|_{z=0}=0$, we can get $q=0$, and hence $u_1=u_2$. Furthermore, using $\partial_{x} u+\partial_{y}(K u)+\partial_{z} w=0$, then $w_i$ $(i=1,2)$   can be uniquely determined (i.e., $w_1=w_2$ ). This completes the proof of the uniqueness of solutions.


\section{Bibliographic Comments}

Let us now briefly review the background and corresponding results about the boundary layer. The mathematical studies on the Prandtl boundary layer have a very long history. Moore \cite{M} gave an analytical framework of the three-dimensional boundary layer in 1956.
In the following decades, the 2-dimensional boundary layer theory developed rapidly.
The first well-known result was developed by Oleinik and Samokhin in \cite{OS}, where under the
monotonicity condition on tangential velocity with respect to the normal variable to the boundary, the local (in time) well-posedness of Prandtl equations was obtained by using the Crocco transformation and von Mises transformation under the outer flow $U\neq\text{constant}$.  Since then, Crocco transformation and von Mises transformation have been used in boundary layer problems. For example, Xin and Zhang \cite{XZ} established a global existence of weak solutions to the two-dimensional Prandtl system for the pressure is favourable (i.e., $\partial_{x}P<0$), which had generalized the local well-posedness results of Oleinik \cite{OS}.
But, Masmoudi and Wong  \cite{mw} proved the local existence and uniqueness of solutions for  classical Prandtl  system by energy methods in a polynomial weighted  Sobolev space under the outer flow $U\neq\text{constant}$,  and Fan, Ruan, Yang \cite{FRY} proved the  local well-posedness for the compressible Prandtl boundary layer equations by the same method.
The authors Qin, Wang and Liu \cite{qinwang2} proved the local existence and uniqueness of the Prandtl-Hartmann regime with  the general outer flow $U\neq\text{constant}$, and Dong and Qin \cite{dongqin}(see also Chapter 2) established the global well-posedness of solutions in $H^s(1\leq s\leq 4)$ with $ U= \text{constant}$ and without monotonicity condition and a lower bound.

Compared to the 2D case, the results of the three-dimensional boundary layer equations were very few.
Until recent years, Liu, Yang and Wang \cite{lwy1,lwy} obtained the local and global existence of weak solutions to the three-dimensional Prandtl equations with a special structure by Crocco transformation.
Liu, Yang and Wang \cite{lwy2} gave an ill-posedness criterion for the Prandtl equations in three-dimensional space  under the outer flow  $(U,V)=\text{constant}$, which shows that the monotonicity condition on tangential velocity fields is not sufficient for the well-posedness of the three-dimensional Prandtl equations. Based on this result,
they proved the existence and uniqueness of local solutions of three-dimensional boundary layer equations without any special structure in the Gevrey function space, see \cite{LMD,Lx}.
Lin and Zhang \cite{lin}  proved the almost global existence of classical solutions to the
3D Prandtl system with the initial data which lie within $\varepsilon$ of a stable shear flow under the outer flow  $(U,V)=\text{constant}$.
Pan and Xu \cite{PX} proved the global existence of solutions to the three dimensional axially
symmetric Prandtl boundary layer equations with small initial data under $(U,V)=0$.

\printindex

\end{document}